\documentclass{memo-l}

\usepackage[english]{babel}
\usepackage[utf8]{inputenc}
\usepackage[toc,page]{appendix}
\usepackage{mathrsfs}
\usepackage{mathabx}
\usepackage{amsmath}
\usepackage{amsfonts}
\usepackage{amssymb}
\usepackage{amscd}
\usepackage{amsthm}
\usepackage{graphicx}
\graphicspath{ {./images/} }
\usepackage{stmaryrd}
\usepackage{mathrsfs, hyperref}
\usepackage{tikz-cd}
\DeclareFontFamily{U}{MnSymbolA}{}
\DeclareSymbolFont{mnsymbols}{U}{MnSymbolA}{m}{n}
\DeclareFontShape{U}{MnSymbolA}{m}{n}{
<-6> MnSymbolA5
<6-7> MnSymbolA6
<7-8> MnSymbolA7
<8-9> MnSymbolA8
<9-10> MnSymbolA9
<10-12> MnSymbolA10
<12-> MnSymbolA12}{}
\DeclareMathSymbol{\rcirclearrowleft}{\mathop}{mnsymbols}{250}
\DeclareMathOperator{\colim}{colim}
\usetikzlibrary{matrix,arrows,decorations.pathmorphing}
\usepackage[all,cmtip]{xy}
\usepackage[normalem]{ulem}

\newcommand{\pE}{\pi^* \! E}

\makeatletter
\def\addressname{Department}

\makeatletter
\def\addressname{\itshape Place of preparation}

\def\@maketitlepage{%
    \thispagestyle{empty}%
    \begingroup
        \topskip\z@skip
        \null
        \vfil
        \begingroup
            \LARGE\bfseries \centering
            \openup\medskipamount
            \@title\par
            \vspace{24pt}%
            \def\and{\par\medskip}%
            \centering
            \mdseries\authors\par
        \endgroup
        \ifx\@empty\addresses\else
            \par\vspace{12pt}
            \begingroup
            \centering\normalsize
            \def\author##1{} 
            \def\\{\unskip, \ignorespaces}%
            \def\address##1##2{\addressname: \scshape\ignorespaces##2\par}%
            \def\curraddr##1##2{\curraddrname: \ignorespaces##2\par}%
            \def\email##1##2{\emailaddrname: \ttfamily\ignorespaces##2\par}%
            \def\urladdr##1##2{\urladdrname: \ttfamily\ignorespaces##2\par}%
            \addresses
            \endgroup
        \fi
        \bigskip
        \vfil
        \ifx\@empty\@translators\else
            \vfill
            \begin{center}
            \@settranslators
            \end{center}
        \fi
        \vfil\vfil
    \endgroup
}
\makeatother

\newtheorem{theorem}{Theorem}[chapter]
\newtheorem{lemma}[theorem]{Lemma}
\newtheorem{proposition}[theorem]{Proposition}
\newtheorem{corollary}[theorem]{Corollary}
\newtheorem{conjecture}[theorem]{Conjecture} 

\theoremstyle{definition}
\newtheorem{definition}[theorem]{Definition}
\newtheorem{example}[theorem]{Example}

\theoremstyle{remark}
\newtheorem{remark}[theorem]{Remark}

\numberwithin{section}{chapter}
\numberwithin{equation}{chapter}

\makeindex

\begin{document}

\frontmatter

\title{Equivariant Instanton Homology}

\author{Mike Miller Eismeier}
\address{\mbox{Department of Mathematics, University of California, Los Angeles}}
\curraddr{\mbox{Department of Mathematics {\&} Statistics, University of Vermont}}
\email{Mike.Miller-Eismeier@uvm.edu}
\thanks{The author's work on this monograph was partially supported by NSF grant number DMS-1708320.}
\thanks{This is the accepted version of a manuscript which has been accepted for publication in \textit{Memoirs of the American Mathematical Society}. The final version of record will be available from the AMS. \copyright~2025 American Mathematical Society.}

\date{August 18, 2025}



\dedicatory{Dedicated twice to Nancy Miller.}

\begin{abstract}
We define four versions of equivariant instanton Floer homology ($I^+, I^-, I^\infty$ and $\widetilde I$) for a class of 3-manifolds and $SO(3)$-bundles over them including all rational homology spheres. These versions are analogous to the four flavors of monopole and Heegaard Floer homology theories. This construction is functorial for a large class of 4-manifold cobordisms, and agrees with Donaldson's definition of equivariant instanton homology for integer homology spheres. Furthermore, one of our invariants is isomorphic to Floer's instanton homology for admissible bundles, and we calculate $I^\infty$ in all cases it is defined, away from characteristic 2.

The appendix, possibly of independent interest, defines an algebraic construction of three equivariant homology theories for dg-modules over a dg-algebra, the equivariant homology $H^+_A(M)$, the coBorel homology $H^-_A(M)$, and the Tate homology $H^\infty_A(M)$. The constructions of the appendix are used to define our invariants.
\end{abstract}

\maketitle

\tableofcontents

\chapter{Introduction}
In \cite{Fl1}, Andreas Floer introduced the instanton homology groups $I(Y)$, $\mathbb Z/8$-graded abelian groups associated to integer homology 3-spheres. These form a sort of TQFT in which oriented cobordisms $W: Y_0 \to Y_1$ induce homomorphisms on the corresponding instanton homology groups. Since then, similar TQFT-style invariants have found themselves a powerful tool in 3- and 4-dimensional topology, especially the related monopole Floer homology of \cite{KMSW} and the Heegaard Floer homology of \cite{OS}.

In ideal circumstances, the instanton homology groups are defined by a chain complex generated by \emph{irreducible flat $SU(2)$ connections} up to isomorphism (equivalently, representations $\pi_1(Y) \to SU(2)$ whose image is non-abelian, modulo conjugacy by elements of $SU(2)$); the component of the differential between two flat connections $\alpha_-, \alpha_+$ is given by an algebraic count of solutions on the cylinder $\mathbb R \times Y$ to the \emph{ASD equation} $$F_{\mathbf{A}}^+ = 0,$$ where $\mathbf{A}$ is a connection on the trivial $SU(2)$-bundle over the cylinder which is asymptotically equal to the $\alpha_\pm$. We can think of this as the ``Morse chain complex" of the Chern-Simons functional on the space of irreducible connections modulo gauge equivalence, $\mathcal B_Y^* = \mathcal A_Y^*/\mathcal G$, whose critical points are the flat connections and gradient flow equation is (formally) the ASD equation. While these equations depend on a choice of metric on the 3-manifold $Y$, the homology groups are an invariant of the $Y$ itself.

Floer's theory is constrained to homology 3-spheres because of the presence of \emph{reducible connections}. While the instanton chain complex above is still a chain complex for rational homology spheres, the proof that its homology is independent of the choice of metric fails in the presence of \emph{$U(1)$-reducible connections} (corresponding to representations $\pi_1(Y) \to SU(2)$ with image lying inside a circle subgroup). One would need to take these reducible connections into account in the definition of the chain complex, but this cannot be done naively: while the instanton chain complex is a Morse complex for the space of irreducible connections, which is an infinite-dimensional manifold, the gauge group $\mathcal G$ does not act freely on the entire space of connections, and so the configuration space of \emph{all} connections modulo gauge, $\mathcal B_Y$, is not a manifold.

Austin and Braam in \cite{AB} resolve this difficulty for a class of 3-manifolds (including all rational homology spheres) by defining an invariant called the \emph{equivariant instanton homology} of $Y$, a $\mathbb Z/8$-graded $\mathbb R$-vector space $I^G_*(Y)$ with an action of $\mathbb R[U] = H^*(BSO(3);\mathbb R)$, to be a form of $SO(3)$-equivariant Morse theory on an infinite-dimensional $SO(3)$-manifold $\widetilde{\mathcal B}_Y$ with $$\widetilde{\mathcal B}_Y/SO(3) = \mathcal B_Y.$$ The manifold $\widetilde{\mathcal B}_Y$ might be called the configuration space of \emph{framed connections} on (the trivial $SO(3)$-bundle over) $Y$. Their invariant is defined using the equivariant de Rham complex as a model form the equivariant (co)homology of a smooth $G$-manifold, and thus inherently uses real coefficients.

Floer also defined, in \cite{Fl2}, instanton homology groups for $SO(3)$-bundles $E$ over 3-manifolds $Y$ satisfying the admissibility criterion that $w_2(E) \in H^2(Y;\mathbb Z/2)$ lifts to a non-torsion class in $H^2(Y;\mathbb Z)$; in particular, $b_1(Y) > 0$. In this case, there are no reducible connections, and the homology of the Floer complex is a well-defined invariant of the pair $(Y,E)$. This case is important for his work on surgery triangles in instanton homology.

Using this, Kronheimer and Mrowka introduce \emph{framed instanton homology groups} $I^\#(Y,E)$ for an arbitrary $SO(3)$-bundle over a 3-manifold in \cite{KM1} by studying the instanton homology of $(Y \# T^3, E \# Q)$ for a certain admissible bundle $Q$ over $T^3$. This is meant to be a version of the (non-equivariant) Morse homology of the space of framed connections $\widetilde{\mathcal B}_E$.

In this paper, we jointly generalize Floer's work on admissible bundles and Austin-Braam's work for rational homology spheres; to speak of both in the same breath, we say that an $SO(3)$-bundle $E$ over a 3-manifold $Y$ is \textit{weakly admissible} if either $w_2(E)$ has no lifts to a torsion class in $H^2(Y;\mathbb Z)$, or if $b_1(Y) = 0$.

We take an alternate approach to Kronheimer and Mrowka's to the framed instanton homology groups: instead of taking a connected sum with $T^3$, we work on the space of framed connections $\widetilde{\mathcal B}_E$ itself. We do this with a sort of Morse-Bott complex for a smooth $G$-manifold equipped with an equivariant Morse function; our definition is partly inspired by the Morse-Bott complex introduced for monopole Floer homology in \cite{Lin}. This uses Lipyanskiy's notion of the \textit{geometric chain complex} $C_*^{\text{gm}}(X;R)$ of a smooth manifold $X$, introduced in \cite{Lip}, whose homology gives the usual singular homology of $X$.

While there are technical obstructions to carrying this out for all 3-manifolds, this has the advantage of providing more structure: for $(Y,E)$ a weakly admissible bundle and $R$ a commutative ground ring, we can define a $\mathbb Z/8$-graded chain complex of $R$-modules, $\widetilde{CI}(Y,E,\pi;R)$, which carries the action of the differential graded algebra $C_*(SO(3);R)$. (This is what we find to be the cleanest notion of a chain complex with an action of the Lie group $SO(3)$.) This chain complex depends on further data $\pi$, including a metric on the 3-manifold itself and a perturbation of the functional defining the Morse complex, but this turns out to be \textit{mostly} inessential: associated to a perturbation $\pi$ on a pair $(Y, E)$ is an element of a finite set $\sigma(Y,E)$ of \textit{signature data}; this will be defined in Chapter \ref{sec:3d-sigdata}. For concreteness, we remark that if $Y$ is a rational homology sphere whose universal abelian cover $\tilde Y$ has $H^1(\tilde Y; \mathbb C) = 0$, then $\sigma(Y, E)$ consists of a single element for all $E$. The set $\sigma(Y,E)$ corresponds precisely to the ``natural classes of perturbations'' stated in the main theorem of \cite{AB}.

In fact, the TQFT structure of the usual instanton Floer homology groups can be lifted to the level of the homology groups of $\widetilde{CI}$. This is the main theorem of this paper.

\begin{theorem}\label{main}There is a category $\mathsf{Cob}_{3,b}^{U(2),w}$ of based connected 3-manifolds $(Y,E,\sigma,b)$ equipped with weakly admissible $U(2)$-bundles and signature data, whose morphisms are certain `weakly admissible' oriented connected cobordisms $(W,\widetilde{\mathbf{E}})$ equipped with a path between the basepoints on the ends. There is also a category of relatively $\mathbb Z/8$-graded $R$-modules with an action of $H_*(SO(3);R)$, and a functor $$\widetilde{I}: \mathsf{Cob}_{3,b}^{U(2),w} \to \mathsf{Mod}^{r,\mathbb Z/8}_{H_*(SO(3);R)}.$$ The relatively graded group $\widetilde I(Y,E,\sigma;R)$ is called the \textit{framed instanton homology} of the triple $(Y,E,\sigma)$. When $Y$ is a rational homology sphere equipped with the trivial bundle, the relative grading may be lifted to an absolute $\mathbb Z/8$-grading.
\end{theorem}

Note the sudden shift to $U(2)$ bundles in this statement. This is to pin down signs in the definition of the differential and cobordism maps, and nothing else: the majority of this text is written in the context of $SO(3)$-bundles, and we only make the passage to $U(2)$-bundles in Chapter \ref{sec:or-canonical}, in our discussion of orientations on the moduli spaces. In particular, if we had chosen to work over a coefficient ring where $1 = -1$, we may omit discussion of $U(2)$-bundles entirely. Furthermore, if we only care about the underlying isomorphism class of group $\tilde I(Y, E, \pi)$, this only depends on the underlying $SO(3)$-bundle: see Corollary \ref{no-u2}.

In the most important cases (including trivial bundles), the weak admissibility condition includes a condition called $\rho$-monotonicity, which is defined in terms of the Atiyah-Patodi-Singer $\rho$ invariant for flat connections and the signature data $\sigma$. Roughly, their sum should always increase across the cobordism. This condition has to do with achieving transversality normal to the reducible locus; when it fails, it is not clear how to try to define the cobordism maps.

For a more precise statement, see Chapter \ref{sec:3d-sigdata} for Definition \ref{sigdata} of signature data; the $\rho$-invariant and $\rho$-monotonicity condition are introduced at the end of Chapter \ref{sec:4d-index}, and Definition \ref{admiss-cob} gives the definition of weakly (and fully) admissible bundles. Finally, Chapter \ref{sec:Floer-eqinst} contains Definition \ref{cobcat} of the weakly admissible cobordism category (and below it, two relatives) as well as Theorem \ref{framed-functor} defining the framed instanton functor. 

It is our expectation that the notion of signature data above is inconsequential:

\begin{conjecture}If $Y$ is a rational homology sphere, $\widetilde{CI}_*(Y,E,\sigma;R)$ is independent of the choice of signature data $\sigma$ up to graded $C_*(SO(3);R)$-equivariant quasi-isomorphism.
\end{conjecture}

\begin{remark}
The above conjecture has been established since this article first appeared as a preprint: see \cite[Theorem 2.9]{DME}. In particular, while the equivariant quasi-isomorphism type of $\widetilde{CI}(Y, E, \sigma)$ is independent of $\sigma$, the equivariant homotopy type is not. 
\end{remark}

From here, if $R$ is a principal ideal domain, we may construct a complex $CI^+_*(Y,E,\sigma;R)$ using the \emph{bar construction} (topologically, the Borel construction) on $\widetilde{CI}$. Its homology groups are denoted $I^+_*(Y,E,\sigma;R)$ and form a relatively graded module under the action of $H^*(BSO(3);R)$. This construction is standard in algebraic topology, and reviewed in the appendix. 

One appealing feature of the monopole and Heegaard Floer theories (which are in some sense $S^1$-equivariant homology theories) is the existence of two other variants of the Floer homology groups that fit into an exact triangle. (Largely, we choose our notation to fit with that of Heegaard Floer theory.) 

In the appendix, we describe algebraic constructions $C^\bullet(A,M)$, where $\bullet = +, -, \infty$, where $A$ is a dg-algebra and $M$ a dg-module. $C^-(A,M)$ is a \emph{cobar} construction, dual to the above bar construction, and $C^\infty(A,M)$ is the \emph{Tate complex}, constructed as a comparison between a variant $C^+$ and $C^-$. Thus up to a sort of twist (often only a grading shift), we have an exact sequence $$0 \to C^{+,\text{tw}}(A,M) \to C^-(A,M) \to C^\infty(A,M) \to 0,$$ leading to an exact triangle of homology groups. The Tate homology groups $H^\infty(A,M)$ satisfy a short list of axioms, including that $H^\infty(A,A) = 0$: they are essentially the homology theory which sets the irreducibles to zero.

Applying a version of this for $\mathbb Z/8$-graded complexes to the framed instanton complex $\widetilde{CI}$ as a module over $C_*(SO(3);R)$, and taking homology, we arrive at the equivariant instanton homology groups $$I^+(Y,E,\sigma;R),\;\; I^-(Y,E,\sigma;R),\;\; I^\infty(Y,E,\sigma;R).$$ As desired, these are all modules over $H^{-*}(BSO(3);R)$, and fit into an exact triangle. 

\begin{theorem}There are functors from the cobordism category of pointed 3-manifolds equipped with weakly admissible $U(2)$ bundles and signature data to the category of relatively $\mathbb Z/8$-graded $H^{-*}(BSO(3);R)$-modules, $$I^\bullet_*: \mathsf{Cob}_{3,b}^{U(2),w} \to \mathsf{Mod}^{r,\mathbb Z/8}_{H^{-*}(BSO(3);R)},$$ $\bullet = +, -, \infty$. When $Y$ is a rational homology sphere, $I^\bullet(Y,E,\sigma;R)$'s relative grading lifts to an absolute $\mathbb Z/8$-grading. Furtheremore, there is a long exact sequence of $R$-modules $$\cdots \to I^+ \xrightarrow{[3]} I^- \to I^\infty \xrightarrow{[-4]} I^+ \to \cdots$$ where $[n]$ denotes that the map increases grading by $n$.\footnote{Note that the grading shift is only meaningful when $Y$ is a rational homology sphere and $E$ is trivial, for which the grading is absolute.} The connecting map $I^-(Y;R) \to I^\infty(Y;R)$ is always an $H^{-*}(BSO(3);R)$-module homomorphism. When $\frac 12 \in R$ or $2 = 0 \in R$, the other two connecting maps are also $H^{-*}(BSO(3);R)$-module homomorphisms.\end{theorem}

It seems plausible that the connecting maps are $H^{-*}(BSO(3);R)$-equivariant for all $R$. Establishing this relates to resolving a technical question about the algebras $C_*(G;R)$; see Chapter \ref{sec:gpalg} for details. \\

The main tool we use to establish invariance properties and perform calculations and compare with existing theories is a collection of spectral sequences that calculate these equivariant homology groups.

\begin{theorem}Let $(Y,E,\sigma)$ be a closed oriented 3-manifold equipped with weakly admissible $U(2)$- bundle and signature datum. Supposing $\pi$ is a regular perturbation for $(Y,E,\sigma)$, we denote by $\mathfrak C_\pi$ the set of critical $SO(3)$-orbits.

Then there are $\mathbb Z/8 \times \mathbb Z$-graded spectral sequences $$\bigoplus_{\alpha \in \mathfrak C_\pi} H^\bullet_{SO(3)}(\alpha;R) \to I^\bullet(Y,E,\sigma;R).$$ These spectral sequences can be used to compute $\widetilde I$ and $I^-$, and if a map induces an isomorphism on the spectral sequences for any $I^\bullet$ on some finite page $E^r$, it induces an isomorphism on $I^\bullet$ itself.
\end{theorem}

This is proved in two parts, Theorems \ref{IndexSS} and \ref{EqIndexSS}, first for nonequivariant homology and then for equivariant homology.

Using the fact that Tate homology of a free $A$-module vanishes, the groups $I^\infty(Y,E)$ are especially computable. For an admissible bundle, they vanish, and we can compare $I^+$ to Floer's invariant $I(Y,E;R)$ for admissible bundles. 

\begin{theorem}\label{intro-admiss}If $E$ is an admissible bundle over a 3-manifold $Y$, then $$I^\infty(Y,E;R) = 0, \;\;\;\; \text{and} \;\;\; I^+(Y,E;R) \cong I(Y,E;R).$$ If we also have $\frac 12 \in R$, the action of $H^*(BSO(3);R) = R[U]$ on $I^+(Y,E;R)$ is taken by the isomorphism to the standard Floer-theroetic $U$-map.
\end{theorem}

The isomorphisms are given in Theorem \ref{admissible-iso} and Corollary \ref{tate-iso}.

Let $\frac 12 \in R$. In \cite{Don}, under favorable conditions on $Y$ Donaldson introduced three chain complexes for $Y$ an integer homology sphere equipped with the trivial $SU(2)$-bundle: the framed complex $\widetilde{CF}(Y;R)$, the equivariant homology complex $\overline{\overline{CF}}(Y;R)$, and the equivariant cohomology complex $\uline{\uline{CF}}(Y;R)$. The first complex has an action of the exterior algebra $\Lambda(u)$, with $|u| = 3$, and the second two complexes have an action of $R[U]$, where $|U| = -4$. In Chapter \ref{secDCI}, we prove the following.

\begin{theorem}For $(Y,E)$ a 3-manifold equipped with a weakly admissible bundle, there is a finite dimensional $\Lambda(u)$-module $DCI(Y,E;R)$ and finite type $R[U]$-modules $\overline{DCI}^\pm(Y,E;R)$, so that there are equivariant quasi-isomorphisms 

\begin{align*}DCI(Y,E;R) &\simeq \widetilde{CI}(Y,E;R) \\
\overline{DCI}^+(Y,E;R) &\simeq CI^+(Y,E;R) \\
\overline{DCI}^-(Y,E;R) &\simeq CI^-(Y,E;R)
\end{align*}

Furthermore, when $Y$ is an integer homology sphere, we have equalities
\begin{align*}DCI(Y,E;R) &= \widetilde{CF}(Y,E;R) \\
\overline{DCI}^+ &= \overline{\overline{CF}}(Y,E;R) \\
\overline{DCI}^-(Y,E;R) &= \uline{\uline{CF}}(Y,E;R),
\end{align*}
up to a rescaling of basis.
\end{theorem}

Here finite type means that it is a direct sum of finitely many simple pieces: for $\overline{DCI}^+$, they are $R$ with trivial $U$-action and $R\llbracket U^*\rrbracket$, where $U$ contracts against $U^*$; for $\overline{DCI}^-$, they are $R$ with trivial $U$-action and $R[U]$ with canonical $U$-action. 

It is by passing through this isomorphism that we show that the isomorphisms Theorem \ref{intro-admiss} above preserve the $U$-action; this is Corollary \ref{admiss-U}.

While more complicated than the Tate calculations for admissible bundles, we are able to exploit the isomorphisms above to calculate instanton Tate homology for an arbitrary rational homology 3-sphere in Chapter \ref{sec:ex-tate}.

\begin{theorem}\label{intro-tate}Let $(Y,E,\sigma)$ be a rational homology 3-sphere equipped with $U(2)$-bundle and signature datum, and suppose $R$ is a PID in which $2$ is invertible.

We write $R[H^2(Y)]$ to mean the group algebra of the finite group $H^2(Y;\mathbb Z) \cong H_1(Y;\mathbb Z)$. If $c = c_1 E \in H^2(Y)$, there is an action of $\mathbb Z/2$ on $H^2(Y)$ given by the involution $x \mapsto c-x$. Further define an action of $\mathbb Z/2$ on the ring $R[U^{1/2}, U^{-1/2}]$, acting on the basis by $U^{n/2} \mapsto (-1)^n U^{n/2}$. Here $|U^{1/2}| = -2$.

Then there is a canonical isomorphism of $\mathbb Z/8$-graded modules over $R[U]$,  $$I^\infty(Y,E,\sigma;R) \cong R[U^{1/2},U^{-1/2}\rrbracket \otimes_{R[\mathbb Z/2]} R[H^2(Y)].$$
\end{theorem}

\begin{remark}In \cite{AB}, an equivariant instanton Floer homology is associated to any 3-manifold with $b_1(Y) = 0$ or $H_1(Y;\mathbb Z) \cong \mathbb Z^a \oplus (\mathbb Z/2)^b$ equipped with the trivial $SO(3)$ bundle, though their analysis is unchanged in the case of $b_1(Y) = 0$ and nontrivial $SO(3)$ bundle. 

In the case that $(Y,E)$ supports an $SO(2)$-reducible flat connection --- that is, if 
$$\langle w_2(E), H_2(Y;\mathbb Z) \bmod 2\rangle, \text{  and  } 2\text{Tors}\left(H_1(Y;\mathbb Z\right) \neq 0,$$
there is a topological obstruction to equivariant transversality. 

When this is the case, the unperturbed moduli space of $SO(2)$-reducible flat connections is a disjoint union of tori $T^{b_1(Y)}$. After applying an $SO(3)$-invariant perturbation, the space of $SO(2)$-invarinat flowlines between various reducible critical points lying in the same $T^{b_1}$ are smooth manifolds of dimension at most $b_1(Y) - 1$; they have a map to a space of Fredholm operators, given by sending a reducible ASD connection splitting as $\mathbf{A} \cong \theta \oplus A$ to the ASD operator corresponding to $A$; in this specific case, this operator has index $0$. This map must be made disjoint from the loci of operators with cokernel of dimension at least 1; this is a subspace of codimension 2. Thus when $b_1(Y) \geq 3$, we may have some non-trivial interChapter with this locus which cannot be removed by a small perturbation. A similar argument shows that there is a further obstruction to achieving invariance than those we have already identified. 

Austin and Braam's invariant in the case that $b_1(Y) > 0$ is defined by careful (non-generic) choices of equivariant perturbation. We expect that the instanton complex defined in this paper can be extended to that level of generality, as well as the proof of its invariance up to signature data. However, the analysis is somewhat more delicate, and the choice of perturbation is necessarily non-generic. The cases with $b_1(Y) > 0$ also add some difficulty in defining a cobordism category that instanton homology is functorial on, as many cobordisms for which the cobordism maps are naturally and essentially uniquely defined do not necessarily compose to cobordisms in the same class.
\end{remark}
%
%
%
\section{Survey of the homology theories}
As a number of different chain complexes are introduced in this paper and exist in the literature, we survey here the definitions and relationships. 

Here $(Y,E)$ is a closed oriented 3-manifold with weakly admissible $U(2)$-bundle, $\pi$ is a regular perturbation and $R$ is a PID. This means, in particular, that there are finitely many critical orbits of the perturbed Chern-Simons functional on $\widetilde{\mathcal B}^e_E$. The finite set whose elements are connected components of the critical set is written $\mathfrak C_\pi$; an element $\alpha \in \mathfrak C_\pi$ is an $SO(3)$-space, either a point, $S^2$, or $SO(3)$ itself. If $E$ is equipepd with a trivialization, each $\alpha \in \mathfrak C_\pi$ has an associated grading $i(\alpha) \in \mathbb Z/8$; otherwise, we instead have relative gradings $i(\alpha, \beta) \in \mathbb Z/8$. For uniformity of notation we often write the relative grading as if it were given by an absolute grading.

The first chain complex introduced is $\widetilde{CI}(Y,E,\pi;R)$. As an $R$-module, this is given by $$\bigoplus_{\alpha \in \mathfrak C_\pi} C_*^{\text{gm}}(\alpha;R)[i(\alpha)].$$ The individual terms $C_*^{\text{gm}}(\alpha;R)$ are the \emph{geometric chain complexes} of the orbits $\alpha$ described in Chapter \ref{sec:Floer-gchain}. A generic basis element of $C_*^{\text{gm}}(\alpha)$ is given by a ``strong $\delta$-chain", but it is helpful to imagine that a generator of $C_*^{\text{gm}}(\alpha;R)$ is given by a smooth map $\sigma: P \to \alpha$, where $P$ is a compact smooth manifold with corners, and $\sigma$ is considered up to diffeomorphism of the domain. (This is essentially a special case of the more general notion of $\delta$-chain.) 

The boundary operator on $\widetilde{CI}(Y,E,\pi;R)$ is given as the sum of the geometric boundary operator (sending $\sigma: P\to \alpha$ to $\partial P \to \alpha$) and a fiber product map with the moduli spaces whose properties are discussed in the first half of this paper, taking a map $\sigma: P \to \alpha$ to a map 
\begin{equation}\label{eqn:fib-moduli}
  P \times_{\alpha} \overline{\mathcal M}(\alpha, \beta) \to \beta. 
\tag{*}
\end{equation}
This chain complex carries the action of a $C_*^{\text{gm}}(SO(3);R)$-module.

After this, one applies the results of the appendix to construct chain complexes $$CI^+(Y,E,\pi;R),\;\;\;\; CI^-(Y,E,\pi;R),\;\;\;\; CI^\infty(Y,E,\pi;R).$$ These are constructed algebraically, and we will not give detailed descriptions here of a generic element of these chain complexes. Writing $\bullet$ for one of $\{+, -, \infty\}$, suffice it to say that if $G$ is a connected Lie group and $H$ a connected subgroup, then for $\alpha = G/H$, there are chain complexes $\hat C^\bullet_G(\alpha;R)$ (the completed group homology complexes) so that $H^+_G(\alpha)$ computes group homology of $H$, while $H^-_G(\alpha)$ computes group cohomology of $H$ in negative degrees and $H^\infty_G(H)$ computes Tate homology of $H$. These all having some degree shift and an action of $H^{-*}(BG;R)$ induced by restriction to cohomology of $BH$. Further, $$CI^\bullet(Y,E,\pi;R) = \bigoplus_{\alpha \in \mathfrak C_\pi} \hat C^\bullet_{SO(3)}(\alpha;R).$$ 

Next when $\frac 12 \in R$ we define a finite-dimensional chain complex $DCI(Y,E,\pi;R)$, given as an $R$-module by $$DCI(Y,E,\pi;R) = \bigoplus_{\alpha \in \mathfrak C_\pi} H_*(\alpha;R)[i(\alpha)],$$ with an action of $H_*(SO(3);R)$; this is now an exterior algebra on a single degree-3 generator, the fundamental class of $SO(3)$. One may think of this as contributing, for each irreducible orbit $\alpha$, a copy of $R$ in degree $i(\alpha)$ and $i(\alpha)+3$ with the action of $H_*(SO(3);R)$ taking the first to the second, for each $SO(2)$-reducible a copy of $R$ in degrees $i(\alpha)$ and $i(\alpha)+2$, and for each full reducible a copy of $R$ in degree $i(\alpha)$. The differential counts $0$-dimensional moduli spaces between irreducibles, as well as the degrees of maps of moduli spaces between different orbits. One of these maps is often called the $U$-map, which we write $U_{\text{Fl}}$ to distinguish from later algebraic terms, also written $U$.

There is a $C_*(SO(3);R)$-equivariant quasi-isomorphism $$DCI(Y,E,\pi;R) \simeq \widetilde{CI}(Y,E,\pi;R).$$ 

When $Y$ is an integer homology sphere and $E$ is the trivial bundle, $DCI(Y,\pi;R)$ is isomorphic to the complex Donaldson writes as $\widetilde{HF}(Y,\pi;R)$ in \cite[Chapter~7.3.3]{Don}; though Donaldson only writes his for rational coefficients, there is no difficulty extending the definition to the broader case $\frac 12 \in R$.

One may immediately apply the constructions of the appendix to $DCI(Y,E,\pi;R)$. Applying the plus-homology construction we arrive at $DCI^+(Y,E,\pi;R)$, and applying the minus-homology construction we arrive at $DCI^-(Y,E,\pi;R)$. If $U$ is a degree $-4$ element and $U^*$ is a degree $4$ element, the underlying $R$-modules are given by $$DCI^+(Y,E,\pi;R) = DCI(Y,E,\pi;R)\llbracket U^*\rrbracket$$ and $$DCI^-(Y,E,\pi;R) = DCI(Y,E,\pi;R)[U].$$ One thinks of each irreducible here as contributing what looks like a copy of $ESO(3)$. 

Applying a trick of Seidel and Smith from \cite{SS-Invol}, we pass from these to complexes $\overline{DCI}^+$ and $\overline{DCI}^-$. The first has underlying $R$-module given by $$\bigoplus_{\alpha \in \mathfrak C_\pi} H^{SO(3)}_*(\alpha; R)[i(\alpha)]$$ and the second by $$\bigoplus_{\alpha \in \mathfrak C_\pi} H_{SO(3)}^{\dim \alpha-*}(\alpha; R)[i(\alpha)].$$ That is, each irreducible contributes a copy of $R$ (in degree $i(\alpha)$ or $i(\alpha)+3$), while each $SO(2)$-reducible contributes a tower $R\llbracket U^{*/2}\rrbracket$ or $R[U^{1/2}][2]$ respectively, where $|U^{*/2}| = 2$ and $|U^{1/2}| = -2$. Lastly, each full reducible contributes a copy of $R\llbracket U^*\rrbracket$ or $R[U]$, respectively. The differentials involve large powers of $U_{\text{Fl}}$.

We have $C^{-*}(BSO(3);R)$-equivariant quasi-isomorphisms $$CI^+(Y,E,\pi; R) \simeq DCI^+(Y,E,\pi; R) \simeq \overline{DCI}^+(Y,E,\pi;R)$$ and $$CI^-(Y,E,\pi; R) \simeq DCI^-(Y,E,\pi; R) \simeq \overline{DCI}^-(Y,E,\pi;R).$$

For an integer homology sphere $Y$ equipped with the trivial bundle $E$, the complex $\overline{DCI}^+(Y,E,\pi;R)$ is isomorphic as a $U$-module to the complex Donaldson writes as $\overline{\overline{CF}}(Y;R)$, and similarly $\overline{DCI}^-(Y,E,\pi;R) \cong \uline{\uline{CF}}(Y;R).$ 

Donaldson also defines complexes $\overline{CF}(Y;R)$ and $\uline{CF}(Y;R)$. These are best understood as the quotient and fixed points of the chain complex $DCI(Y,\pi;R) \cong \widetilde{CF}(Y,\pi;R)$ under the $H_*(SO(3);R)$-action. We do not use these, but do notice them as appearing in certain spectral sequences for integer homology spheres. In that context, we write the resulting homology groups as $\overline I$ and $\uline I$.

There are also Froyshov's reduced instanton Floer homology groups, which he writes $\widetilde{HF}$ and Donaldson writes $\widehat{HF}$. We identify these in Chapter \ref{examples} as being the image of $I^+(Y;R)$ inside $I^-(Y;R)$ (or more precisely, we see this at the level of an $E^\infty$ page of a spectral sequence calculating those), but otherwise do not use them. If we were to give these groups a notation, we could call them $\widehat I$. 

There is also Floer’s original version of instanton homology, $I(Y)$, which is defined only for integer homology spheres and does not use the reducibles. This, too, may be seen in terms of the spectral sequence for $I^+$ (ignoring the reducible piece). 

Finally, Kronheimer and Mrowka define instanton homology groups $I^\#(Y,E;R)$ for all pairs of 3-manifolds and $SO(3)$-bundles, by taking the connected sum with a pair $(T^3, E)$ where $w_2(E)$ is Poincar\'e dual to $T^2 \times \{*\}$, also called `framed instanton homology'. These are not the same as the groups $\widetilde I(Y,E;R)$. A calculation of \cite{scaduto2015instantons} using Fukaya's connected sum theorem shows that when $Y$ is an integer homology sphere, the group $I^\#(Y;R)$ may be calculated using a chain complex very much like that defining $DCI(Y,\pi;R)$, but they differ in one component of the matrix defining the differential: the term $U_{\text{Fl}}$ in $\partial_{DCI}$ is instead given by $U_{\text{Fl}} - 8$ in the complex defining $I^\#(Y;R)$. So one finds instead that $I^\#(Y,E;R)$ is a sort of deformation of $\widetilde I(Y,E;R)$. We do not discuss this relationship further here. 

In summary, there are nine versions of instanton homology for integer homology spheres, $$\widetilde{I}, I^+, I^-, I^{\infty}, \overline I, \uline I, \widehat I, I^{\#}, I;$$ all but the last of these are also defined for pairs of a rational homology sphere and $U(2)$-bundle. In this paper we will work most extensively with the first four of these. We will briefly mention the next three `reduced’ homology theories $\overline I, \uline I,$ and $\hat I$ when discussing spectral sequence calculations for integer homology spheres. We will not mention $I^{\#}$ and $I$ any further.

In addition to the complexes defining each of these homology theories, we will also use five important complexes for calculation, the Donaldson models $$DCI, \; DCI^+, \; \overline{DCI}^+, \; DCI^-, \; \overline{DCI}^-.$$ The first of these has homology naturally isomorphic to $\widetilde I$, the second two of these have homology naturally isomorphic to $I^+$, and the last two of these have homology naturally isomorphic to $I^-$.  
\subsection*{Convention} 
All $3$- and $4$- manifolds are assumed connected and oriented unless stated otherwise. 

\subsection*{Organization}
Chapter \ref{chap:2} introduces the framed configuration spaces on which we attempt to do Morse theory, and discusses the action of $SO(3)$ on them, and then explains how to complete these topological spaces using Sobolev spaces and obtain the structure of a Banach manifold.

The technical heart of the paper is in Chapters \ref{chap:3} and \ref{chap:4}, where we use the standard holonomy perturbations in instanton Floer theory to show that we can achieve equivariant transversality: a generic perturbation gives rise to a finite set of critical $SO(3)$-orbits, and the space of trajectories between them forms a smooth manifold. As long as the dimension is sufficiently small, this has a compactification to a topological manifold with corners, and we use these compactifications to define the complex $\widetilde{CI}$. This includes a calculation of the reducible perturbed instantons in Chapter \ref{sec:4d-red} and their indices in Chapter \ref{sec:4d-index}, which is then used to give a condition which guarantees we can achieve transversality normal to the reducible locus. This condition uses the notion of signature data introduced in Chapter \ref{sec:3d-sigdata}.

Part \ref{part1} provides us with the technical machine (the moduli spaces and their properties) we need. Part \ref{part2} is devoted to defining and calculating equivariant instanton homology given this machine, and may be read independently of the first part of the paper.

In Chapter \ref{chap:6}, we define the invariants $CI^\bullet$. First we review Lipyanskiy's geometric homology, which is a crucial technical tool in our definition of $\widetilde{CI}$: it allows us to define the chain complex without triangulating the moduli spaces appearing in \eqref{eqn:fib-moduli}, and only requires the use of moduli spaces of small dimension; in particular, small enough that the Uhlenbeck bubbling phenomenon does not arise. The chain complex $\widetilde{CI}(Y,E,\pi;R)$, which depends on a choice of metric and perturbation, is defined in Chapter \ref{sec:Floer-eqinst}. This construction comes with cobordism maps, and these provide us with the usual invariance properties, giving Theorem \ref{main}.

In Chapter \ref{sec:Floer-filt} we explain the notion of index filtration for $\widetilde{CI}$, which gives rise to a $\mathbb Z/8 \times \mathbb Z$-graded spectral sequence. This is not a filtration on $\widetilde{CI}$ in the standard sense, but rather what we call a \emph{periodic filtration}. In Chapter \ref{sec:fourflavors} we use the idea of periodic filtrations to define the equivariant instanton homology complexes $CI^\bullet$, following the construction of Chapter \ref{PeriodicMachine}. 

In Chapter \ref{chap:7}, we use the index spectral sequences of Chapters \ref{sec:Floer-filt} and \ref{sec:flourflavors} to carry out some calculations and comparisons.

In Chapter \ref{sec:ex-admissible} we warm up with a calculation of $I^\bullet(Y,E)$ for admissible bundles, giving Theorem \ref{intro-admiss} above. 

After this, in Chapter \ref{secDCI} we compare $\widetilde{CI}$ to Donaldson's complex $DCI$ (written in \cite{Don} as $\widetilde{CF}$), and in particular show that their homologies are isomorphic, justifying the notation $H(\widetilde{CI}) = \widetilde I$. We are then able to show $I^+(Y;R) \cong \overline{\overline{HF}}(Y;R)$ for rings containing $1/2$, as well as $I^-(Y;R) \cong \uline{\uline{HF}}(Y;R)$. In the same chapter, we extend the definition of Donaldson's chain complexes to rational homology spheres. 

Chapter \ref{sec:ex-tate} uses the Donaldson model, and a localization/periodicity theorem for Tate homology, to give the explicit formula for $I^\infty$ mentioned as Theorem \ref{intro-tate} above. Chapter \ref{examples} contains a calculation of $I^\bullet$ for the 3-sphere, lens spaces, and the Poincar\'e homology sphere, as well as a description of the spectral sequences for integer homology spheres and a brief comment on the connection to Fr\o yshov's reduced groups. We conclude in Chapter \ref{sec:ex-orrev} with a discussion of orientation-reversal and duality in $CI^\bullet(Y)$. 

The appendix describes the algebraic constructions $C^\bullet(A,M)$ and their invariance properties, as well as providing calculational tools. In particular, if $A$ is the group algebra $C_*(G;R)$, where $G$ is a compact Lie group, then $C^\bullet(G, G/H)$ is calculated for all $\bullet$ and any closed subgroup $H \subset G$. Chapter \ref{PeriodicMachine} concludes with an extension of this machinery to the $\mathbb Z/8$-graded case, when our complexes come equipped with well-behaved periodic filtrations.
%
%
%
\subsection*{Acknowledgements}I thank Ciprian Manolescu for his constant advice and support. I also thank Ali Daemi, Kim Fr\o yshov, Tom Mrowka, Chris Scaduto, and Matt Stoffregen for useful discussions about gauge theory, and Sucharit Sarkar for conversations on equivariant homology. I also thank Kevin Carlson, Tyler Lawson, and Mike Hill, for helpful conversations during the preparation of the appendix, and especially Aaron Royer for many patient discussions. 

I would also like to thank Mariano Eccheveria, Robert Ladu, and several anonymous referees for a careful reading and discussion of parts of this article. Finally, Gard Olav Helle deserves special thanks for his careful read of the appendix. In particular, the final version of this manuscript uses the clean approach to completion presented in \cite{GOH}, greatly simplifying the discussion. 

\newpage

\mainmatter

\part{Properties of moduli spaces}\label{part1}


\chapter{Configuration spaces and their reducibles}\label{chap:2}
Let $Y$ be an oriented closed 3-manifold equipped with a Riemannian metric, and $E \to Y$ an $SO(3)$-bundle over $Y$; it is worth mentioning that $SO(3)$-bundles over a 3-complex are determined up to isomorphism by their second Stiefel-Whitney class $w_2(E) \in H^2(Y;\mathbb Z/2)$. Associated to $E$ by the adjoint representation $SO(3) \curvearrowright \mathfrak{so}(3)$ is the adjoint bundle $\mathfrak g_E \subset \text{End}(E)$, the subbundle given by skew-adjoint endomorphisms. The space $\mathcal A_E$ of orthogonal connections on $E$ is affine over $\Omega^1(\mathfrak g_E)$. There is a gauge group, $\mathcal G_E$, the set of smooth bundle automorphisms of $E$ that cover the identity; equivalently, this is the set of smooth sections of the non-principal $SO(3)$ bundle $\text{Aut}(E)$, the associated bundle to $E$ under the conjugation action $SO(3) \curvearrowright SO(3)$. Using the isomorphism $\text{Inn}(SU(2)) \cong SO(3)$, we may form the bundle $$\widetilde{\text{Aut}}(E) = \text{Aut}(E) \times_{SO(3)} SU(2).$$ We say that $\sigma \in \mathcal G_E$ is an \textit{even gauge transformation} if $\sigma$ lifts to a section of $\widetilde{\text{Aut}}(E)$, and denote the group of even gauge transformations by $\mathcal G^e_E$. Obstruction theory applied to sections of $\text{Aut}(E)$ provides a short exact sequence of groups $$1 \to \mathcal G^e_E \to \mathcal G_E \to H^1(Y;\mathbb Z/2) \to 0.$$

There is a map, the \textit{Chern-Simons functional}, $$\text{cs}: \mathcal A_E \to \mathbb R,$$ defined as follows. Pick a compact oriented 4-manifold $X$, equipped with an $SO(3)$ bundle $E_X$, such that $\partial(X,E_X) = (Y,E)$. Given a connection $A$ on $E$, extend it to a connection $A_X$ on $E_X$. Then define $$\text{cs}(A) = \int_X \text{Tr}(F_{A_X}^2).$$ This gives a function that depends only on the choice $(X,E_X)$, not on the extension $A_X$. If $(X',E_{X'},A_{X'})$ is another extension, we may define $M = X \cup_Y \overline{X'}$ with the obvious choice of $SO(3)$-bundle and connection over it and invoke the Chern-Weil formula $$-2 \pi^2 p_1(E') = \int_M \text{Tr}(F_{A'}^2).$$ Thus $\text{cs}$ is defined in general up to an $8\pi^2 \mathbb Z$ ambiguity (using that $4 \mid p_1(E')$). If $(X', E_{X'}) = (X, E_X)$, then $p_1(E_X \cup E_{X'}) = 0$, giving us a well-defined functional $\text{cs}$ on $\mathcal A_E$ conditional on that choice of bounding 4-manifold $(X, E_X)$.

The gauge group $\mathcal G_E$ acts on $\mathcal A_E$ by $\sigma(A) = A - (\nabla_A \sigma) \sigma^{-1}$, where we take the covariant derivative of $\sigma$ by considering it as a section of $\text{End}(E) = E \otimes E^*$, where $E$ is considered now as an oriented vector bundle equipped with metric. We denote the quotient by this group action $$\mathcal B_E = \mathcal A_E/\mathcal G_E,$$ the configuration space of connections on $E$, and also write the \textit{even configuration space} as the quotient $$\mathcal B^e_E := \mathcal A_E/\mathcal G^e_E.$$ Immediately from the definition we see that the stabilizer at $A$ by the action of $\mathcal G_E$ or $\mathcal G^e_E$ is precisely the subset of (even) $A$-parallel gauge transformations; thus elements of the stabilizer are determined by their value at a single point, and evaluating at a point $b \in Y$ gives an isomorphism to a subgroup of $SO(3) = \text{Aut}(E_b)$.

The Chern-Simons functional $\text{cs}: \mathcal A_E \to \mathbb R$ does not descend naively to $\mathcal B_E$ --- its value may change after applying an element of gauge group $\mathcal G_E$. However, if $u$ is a gauge transformation, $$\text{cs}(u(A)) = \text{cs}(A) + 8\pi^2 k$$ for some integer $k$, so $\text{cs}$ descends to a continuous map $\mathcal B_E \to \mathbb R/8\pi^2 \mathbb Z$. This statement is little more than saying that there is always a 4-manifold with $SO(3)$-bundle $(X,E)$ over which the gauge transformation $u$ extends, for then one sees that $$\text{cs}(u(A)) = \int_X \text{tr}\left(F_{u(A)}^2\right) = \int_X \text{tr}\left(F_A^2\right);$$ that the integrals are equal is Stokes' theorem. It is this circle-valued functional, on a slightly modified space, that we hope to do Morse theory with.

\section{Framed connections and the framed configuration space}\label{sec:config-framed}
The configuration space $\mathcal B_E$ is in no sense a manifold, because $\mathcal G_E$ does not act freely on the space of connections. To free up the action of the gauge group on the space of connections so that the quotient by gauge is a manifold (at least heuristically, at this point), we pick a basepoint $b \in Y$ and consider the space of \emph{framed connections} 
$$\widetilde{\mathcal A}_E = \mathcal A_E \times E_b.$$ 
We call a point $p \in E_b$ a framing because it determines an isomorphism $SO(3) \cong E_b$ sending the identity to the point $p$; here we are thinking of $E$ as a principal bundle, not a vector bundle. A gauge transformation $\sigma$ evaluates to $\sigma(b) \in \text{Aut}(E_b)$, which acts on $E_b$ on the left by applying the automorphism. (By definition, to say $f\in\text{Aut}(E_b)$ means $f(pg) = f(p)g$.) This gives us an action of $\mathcal G_E$ on the space of framed connections $\widetilde{\mathcal A}_E$. There is further a natural \emph{right} $SO(3)$ action on $\widetilde{\mathcal A}_E$ by acting on $E_b$ by translation. If $\sigma \in \mathcal G^e_E$, $(A,p) \in \widetilde{\mathcal A}_E$, and $g \in SO(3)$, we have that 
$$\sigma \cdot ((A,p) \cdot g) = (\sigma \cdot (A,p)) \cdot g,$$ 
because the gauge group acts by automorphisms of the right $G$-set $E_b$.

Because any gauge transformation in the stabilizer of $A$ is $A$-parallel, it is trivial if its value at any point is the identity. As a consequence, $\mathcal G^e_E$ acts \textit{freely} on $\widetilde{\mathcal A}_E$. Its quotient under the gauge group action, denoted $\widetilde{\mathcal B}^e_E$, retains the right $SO(3)$-action. We call this the \textit{framed configuration space} of connections on $E$. The stabilizer of $[A,p] \in \widetilde{\mathcal B}^e_E$ in $SO(3)$ is the isomorphic image of the stabilizer of $A \in \mathcal A_E$ under the action of the even gauge group following evaluation $\mathcal G^e_E \to \text{Aut}(E_b)$ and the natural isomorphism $\text{Aut}(E_b) \cong SO(3)$. There is a map $$\widetilde{\mathcal B}^e_E \to \mathcal B^e_E,$$ given by quotienting by the leftover $SO(3)$ action or equivalently by forgetting the framing. Said another way, orbits of the $SO(3)$-action on $\widetilde{\mathcal B}_E$ or $\widetilde{\mathcal B}^e_E$ correspond to gauge equivalence classes of connections, and a point on an orbit $\alpha \subset \widetilde{\mathcal B}_E$ corresponding to $[A]$ is an equivalence class of framings, with $\alpha \cong E_b/\text{Stab}([A])$.

We will soon define Hilbert manifold completions of these spaces, in the context of which we will find that $\widetilde{\mathcal B}^e_E$ is a smooth Hilbert manifold with smooth $SO(3)$-action. These remarks also apply to the full gauge group, with quotient $\widetilde{\mathcal B}_E = \widetilde{\mathcal A}_E/\mathcal G_E$.

\begin{remark}One could also define $\widetilde{\mathcal B}^e_E = \mathcal A_E/\widetilde{\mathcal G}_E^{e,b}$, quotienting by the group of \textit{based} gauge transformations (gauge transformations with $\sigma(b) = \text{Id}$). This works just as well here, and the right $SO(3)$ action corresponds to the inverse of the left $SO(3)$ action given by $\mathcal G_E^e/\mathcal G_E^{e,b}$. However, we find later discussions of the 4-manifold configuration spaces and the restriction maps to their ends clearer in the language of framed connections, so this is our preferred model.
\end{remark}

\section{The equivalent $U(2)$-model}\label{sec:config-U2}
We would like to understand the reducible subspaces of $\widetilde{\mathcal B}^e_E$ under the $SO(3)$ action. The $SO(2)$ fixed points are more easily understood if we introduce an auxiliary construction. Pick a principal $U(2)$-bundle $\tilde E$ and an isomorphism $$\tilde E \times_{U(2)} SO(3) \cong E.$$ This is possible because $U(2)$-bundles on 3-manifolds are classified by their first Chern class (and $SO(3)$ bundles by their second Steifel-Whitney class), so we only need to know that we can pick a lift of $w_2$ to an integral cohomology class; that this is possible follows from the Bockstein long exact sequence and the fact that $H^3(Y;\mathbb Z) \cong \mathbb Z$ has no 2-torsion. Thinking of $\tilde E$ instead as a complex vector bundle and explicitly identifying the quotient homomorphism $U(2) \to SO(3)$, the construction $\tilde E \times_{U(2)} SO(3) \cong E$ produces the oriented 3-plane bundle $\mathfrak{su}(\tilde E) \subset \text{End}(\tilde E)$ of skew-Hermitian endomorphisms of $\tilde E$ as the associated $SO(3)$-bundle.

Fix a connection $A_0$ on the determinant complex line bundle $\text{det}(\tilde E) := \lambda$ (if $c_1(\tilde E)$ is finite order, we can choose this to be the flat connection, unique up to gauge transformation). We consider the space of connections with fixed determinant connection $$\mathcal A_{\tilde E}^{\det} = \{A \in \mathcal A_{\tilde E} \mid \text{tr}(A) = A_0 \in \Omega^1(Y;i\mathbb R)\};$$ this is also affine over $\Omega^1(\mathfrak g_E)$. To every connection on $\tilde E$ there is an associated connection on $E$, giving us a map $\mathcal A_{\tilde E} \to \mathcal A_E$ which is a bijection when restricted to $\mathcal A_{\tilde E}^\text{det}$. To enhance this to a bijection of spaces of framed connections, let 
$$\widetilde{\mathcal A}_{\tilde E}^\text{det} = \mathcal A_{\tilde E}^{\text{det}} \times (\tilde E_b \times_{U(2)} SO(3));$$ 
this carries a natural right action by $SO(3)$. Furthermore, the natural gauge group acting on $\mathcal A_{\tilde E}^\text{det}$ is the set of gauge transformations whose (pointwise) determinant is $1$, denoted 
$$\mathcal G_{\tilde E}^{\text{det}} = \Gamma(\text{Aut}(\tilde E) \mid \det \gamma = 1).$$ 
We can further identify this latter bundle of groups as isomorphic to $\widetilde{\text{Aut}}(E)$ and thus there is a surjective homomorphism $$\mathcal G_{\tilde E}^{\text{det}} \to \mathcal G^e_E.$$ It is not a bijection: recall that the latter group is defined as a subset of $\mathcal G_E = \Gamma(\text{Aut}(E))$. Its kernel is the 2-element set of gauge transformations whose pointwise values are $\pm 1 \in SU(2)$. This subgroup acts trivially on $\widetilde{\mathcal A}_{\tilde E}^\text{det}$.

Now it is easy to verify that the bijection $\widetilde{\mathcal A}_{\tilde E}^\text{det} \to \widetilde{\mathcal A}_E$ is equivariant under the actions of the gauge groups, and thus after quotienting we have an equivariant diffeomorphism \begin{equation}\label{U(2)}\widetilde{\mathcal B}_{\tilde E}^\text{det} \xrightarrow{\cong} \widetilde{\mathcal B}^e_E.\end{equation} The group acting on the former space is $SU(2)$, but $\pm 1$ act trivially, and so passing to the quotient $PSU(2) \cong SO(3)$ we identify these two configuration spaces as $SO(3)$-spaces. Furthermore, the natural definition of Chern-Simons functional on this first space is sent to our Chern-Simons functional on the second space, and there is no difference in the resulting gauge theory. We only consider this construction auxiliary because it depends on the unnecessary input data of $\tilde E$ and $A_0$. When working on the $U(2)$-bundle, we usually prefer to speak of the $SU(2)$ action, as $SU(2)$ naturally sits inside the $U(2)$ gauge group, even though $\pm 1$ act trivially on the framed configuration space.

In the simple case that $\tilde E = \eta_1 \oplus \eta_2$, the associated $SO(3)$-bundle is $i\mathbb R \oplus (\eta_1 \otimes \eta_2^{-1})$, and a connection on the former respecting the splitting is taken under the bijection $\mathcal A_{\tilde E}^{\text{det}} \to \mathcal A_E$ to a connection which respects the latter splitting.

\section{Reducibles on $3$-manifolds}\label{sec:config-red-3d}
We can explicitly describe the stabilizers that arise under the action of the even gauge group on $\mathcal A_E$, and hence the orbit types in $\widetilde{\mathcal B}^e_E$. In understanding these, it's also convenient to carry out the same analysis for the full gauge group $\mathcal G_E$ and its quotient $\widetilde{\mathcal B}_E$. Recall that $$\mathcal G_E/\mathcal G^e_E \cong H^1(Y;\mathbb Z/2),$$ identifying the latter with the group of obstructions to lifting a section of $\text{Aut}(E)$ to a section of $\widetilde{\text{Aut}}(E)$. This leaves an $H^1(Y;\mathbb Z/2)$ action on $\widetilde{\mathcal B}^e_E$, whose quotient is $\widetilde{\mathcal B}_E$, the quotient by all gauge transformations. This is a useful gadget to keep track of, especially in light of the following lemma (a version of \cite[Lemma~4.28]{DK}), calculating the stabilizers of connections under both $\mathcal G^e_E$ and $\mathcal G_E$.

\begin{lemma}\label{stab1}Let $A$ be an $SO(3)$-connection on $E$. Then the stabilizer $\Gamma_A$ under the action of the full gauge group $\mathcal G_E$ is $$C(H_A) \subset SO(3),$$ the centralizer of the holonomy group of $A$ at some choice of basepoint $b$. Let $\pi: SU(2) \to SO(3)$ be the projection. $\Gamma_A \cap \mathcal G^e_E$ is $$C_{SU(2)}(\pi^{-1}H_A) \subset SO(3),$$ the set of elements of $SO(3)$ that fix $\pi^{-1}H_A$ under conjugation, considering $SO(3) = \text{Inn}(SU(2))$.
\end{lemma}

\begin{proof}
If $SO(3)$ acts smoothly on some manifold $M$, consider the $M$-bundle $E \times_{SO(3)} M$. If $A$ is an $SO(3)$-connection on $E$, there is a natural connection induced on the associated $M$-bundle. Then $m \in M_x$ extends to a parallel section of this $M$-bundle if and only if it is fixed under the action of $H_A$ on $M$. Taking $M$ to be $SO(3)$ and $SU(2)$ equipped with the conjugation action gives the desired result.\end{proof}
A group-theoretic calculation shows that the only subgroups of $SU(2)$ that arise as centralizers are the center $\mathbb Z/2$, the circle subgroups $U(1)$, and the full group $SU(2)$. So the only possible stabilizers of the action of $\mathcal G^e_E$ on $\mathcal A_E$ are 
$$\{I\}, SO(2), SO(3).$$
A similar calculation shows that the subgroups of $SO(3)$ that arise as centralizers are additionally $O(2)$, the diagonal subgroup $V_4$, and $\mathbb Z/2$. (Calculate centralizers of elements first, then calculate the possible intersections of these.) Comparing stabilizers, we see that the action of $H^1(Y;\mathbb Z/2)$ on $\mathcal B^e_E$ giving rise to the quotient by the full gauge group $\mathcal B_E$ is free except at reducibles with full stabilizer $O(2), V_4$, or $\mathbb Z/2$, where the action of $H^1(Y;\mathbb Z/2)$ has stabilizer isomorphic to $\mathbb Z/2, V_4$, and $\mathbb Z/2$, respectively.

When the stabilizer of a connection $A$ on $E$ is $SO(3) \subset \mathcal G^e_E$, the bundle must be trivial and the connection gauge equivalent to the trivial connection; however, it needn't be equivalent by an \textit{even} gauge transformation. In fact, because $H^1(Y;\mathbb Z/2)$ acts freely on points in $\mathcal B^e_E$ with full stabilizer and there is only one such in the full quotient, there are $H^1(Y;\mathbb Z/2)$ different elements of $\widetilde{\mathcal B}^e_E$ with stabilizer $SO(3)$. 

Inside the $U(2)$-model, \cite[Section~3]{SS} identifies the action of $H^1(Y;\mathbb Z/2)$ as sending a connection on $\tilde E$ to the corresponding connection on $\tilde E \otimes \xi_{\mathbb C}$, where $\xi_{\mathbb C}$ is the complexification of a real line bundle (equivalently, a complex line bundle equipped with a unitary connection with holonomy in $\pm 1$.)

We will ultimately be interested in reducible \textit{orbits} of the $SO(3)$-action. To describe the reducible orbits in a $G$-space $X$ with orbit type isomorphic to $G/H$, it suffices to describe the set of $H$-fixed points $X^H$, and the action of the \textit{Weyl group} 
$$W(H) = N_G(H)/H$$ 
on $X^H$: every orbit $Gx$ of a point $x$ with stabilizer conjugate to $H$ intersects $X^H$ nontrivially; and in fact, if $x \in X^H$, then 
$$Gx \cap X^H = W(H)x.$$ 
Thus, for instance, a $G$-invariant function on the subspace of points whose stabilizer contains a conjugate of $H$ is determined uniquely by a $W(H)$-invariant function on $X^H$. The description in terms of fixed subspaces and Weyl groups tends to be easier to state and prove, so we largely prefer that language as long as possible. The only case we actually use is that of $SO(2) \subset SO(3)$, which has Weyl group isomorphic to $\mathbb Z/2$.

\begin{proposition}\label{action1}Suppose $E$ is an $SO(3)$-bundle over a closed 3-manifold $Y$. Denote $\tilde E$ a $U(2)$ bundle with $c_1(\tilde E) \mod 2 = w_2(E)$ and write $\lambda = \det{\tilde{E}}$, as in Chapter \ref{sec:config-U2}.\begin{enumerate}
\item The set of points of $\widetilde{\mathcal B}^e_E$ fixed by $SO(2) \subset SO(3)$ is identified with $$\bigsqcup_\eta \mathcal B_\eta := \bigsqcup_\eta \mathcal A_\eta/\mathcal G_\eta$$ where $\eta$ varies over isomorphism classes of complex line bundles, $\mathcal A_\eta$ is the configuration space of unitary connections on $\eta$, affine over $\Omega^1(Y;i\mathbb R)$, and the gauge group $\mathcal G_\eta$ is the space of sections of $\textup{Aut}(\eta)$, which is the same as the space of maps $\textup{Map}(Y, S^1)$. The $SO(3)$ orbit of the $SO(2)$-fixed point space is the set of all reducible connections. The action of the Weyl group sends a class of connection $[A]$ on $\eta$ to the class of $[A_0 - A]$ on $\lambda \otimes \eta^{-1}$.
\item If $E$ is trivial, $\widetilde{\mathcal B}^e_E$ has $SO(3)$ fixed point set in bijection with $H^1(Y;\mathbb Z/2)$; if $E$ is nontrivial, the even configuration space has no $SO(3)$ fixed points.
\end{enumerate}
\end{proposition}

\begin{proof}Following equation (\ref{U(2)}), it suffices to find the fixed subspaces in $\widetilde{\mathcal B}_{\tilde E}^{\text{det}}$. Consider the subset of framed connections $(\tilde A, p)$ in $\widetilde{\mathcal A}_{\tilde E}^{\det}$ whose connection term $\tilde A$ has stabilizer consisting of gauge transformationgs with $\sigma(b)$ in the diagonal subgroup $$S(U(1) \times U(1)) \subset SU(2),$$ using the framing $p$ to specify the isomorphism $$\text{Aut}(\tilde E_b) \cong U(2).$$ All other circle stabilizers in $\mathcal G_{\tilde E}^{\text{det}}$ are conjugate to this subgroup, and this fixed subspace of the space of framed connections projects to the subset of $SO(2)$-fixed connections in $\widetilde{\mathcal B}_{\tilde E}^{\text{det}}$. Our assumption implies that $\tilde A$ has holonomy contained in the diagonal subgroup 
$$U(1) \times U(1) \subset U(2).$$ This leaves invariant a splitting $$\tilde E_b \cong \mathbb C^2 \cong \mathbb C \oplus \mathbb C,$$ the first isomorphism given by the framing $p$, and hence gives a parallel splitting $\tilde E \cong \eta \oplus \eta'$. Sending $(\tilde A, p)$ to the corresponding connection $A'$ on $\eta$ gives us the map from this subset to $\mathcal A_\eta$, and because any gauge transformation of $\eta$ extends to a determinant-1 gauge transformation of $\tilde E$ (act by the inverse in the $\eta'$ coordinate), this descends to a well-defined map $$\left(\widetilde{\mathcal B}_{\tilde E}^{\text{det}}\right)^{U(1)} \to \bigsqcup_\eta \mathcal B_\eta$$ modulo gauge. Conversely, any connection on $\eta$ induces a connection on $\tilde E$ of the specified form, and any gauge transformation of $\eta$ is induced by a unique determinant-1 gauge transformation of $\tilde E$. Any framing on $\eta$ then induces a framing on $\tilde E$, but there is a determinant-1 gauge transformation changing any framing on $\eta$ to another, so after modding out by gauge the choice of framing on $\eta$ didn't matter. This gives the stated bijection for those connections whose induced splitting gives a line bundle of topological type $\eta$ in the first coordinate; but note that on a 3-manifold, a $U(2)$-bundle is determined by its first Chern class, and so we have an isomorphism $$\tilde E \cong \eta \oplus (\lambda \otimes \eta^{-1}).$$ Thus every complex line bundle $\eta$ arises in such a splitting.

The action of the Weyl group is to swap the coordinates of the framing; this then swaps the components $[A' \oplus (A_0 - A')]$ of the connection, as stated in the lemma. In particular, it only fixes $A'$ if there is an isomorphism $$\eta \cong \lambda \otimes \eta^{-1}$$ and this isomorphism sends $A'$ to $A_0 - A'$. This is only possible if $\lambda$ is twice an integral class, and so $w_2(E) = 0$; but if one such choice is made, all others are affine over $H^1(Y;\mathbb Z/2)$, tensoring the whole bundle (and hence each component) with $\xi_{\mathbb C}$, the complexification of a real line bundle, or equivalently a complex line bundle with holonomy in $\pm 1$.\end{proof}

\begin{corollary}\label{orbits1}Let $E$ and $\tilde E$ be as in the previous proposition; write $c = c_1(\lambda) = c_1(\tilde E)$. Then the reducible subspace of $\widetilde{\mathcal B}^e_E$, consisting of framed connections with nontrivial stabilizer, is a disjoint union over connected components labeled by pairs $\{z_1, z_2\} \subset H^2(Y;\mathbb Z)$ with $z_1 + z_2 = c$, where the $z_i$ are cohomology classes corresponding to complex line bundles $\eta_i$.

If $z_1 \neq z_2$, this connected component is a fiber bundle over $\mathcal B_{\eta_1}\cong \mathcal B_{\eta_2}$ with fiber $S^2$, where $SO(3)$ acts trivially on the base and via the standard action on the fiber. 
If $z_1 = z_2$, this connected component contains a unique fully reducible connection.
\end{corollary}

\begin{proof}Consider the subspace of framed reducibles corresponding to the splitting $\tilde E \cong \eta_1 \oplus \eta_2$; we write this space as $\widetilde{\mathcal B}^{\text{red}}_{\eta_1, \eta_2}$. The quotient of this space by the $SO(3)$-action is the same as the quotient of its $SO(2)$-fixed subspace by the action of $\mathbb Z/2$. In the case that $z_1 \neq z_2$, the $SO(2)$-fixed subspace is $\mathcal B_{\eta_1} \sqcup \mathcal B_{\eta_2}$, and the action of the Weyl group identifies these. Thus the desired fiber bundle is the quotient map $$\widetilde{\mathcal B}^{\text{red}}_{\eta_1, \eta_2} \to \mathcal B_{\eta_1} \cong \mathcal B_{\eta_2}.$$ 

When $z_1 = z_2 =: z$, corresponding to a complex line bundle $\eta$, the $SO(2)$-fixed subspace is $\mathcal B_{\eta}$. Its quotient by the action of the Weyl group is connected, so $\widetilde{\mathcal B}^{\text{red}}_{\eta_1, \eta_2}$ is connected. In this case, $E$ is trivial, so we may choose $\tilde E$ trivial for convenience of discussion; then the action of the Weyl group may be described as complex conjugation, using the isomorphism $\eta \cong \xi_{\mathbb C}$ for a unique real line bundle $\xi$; the induced connection is the fully reducible connection.
\end{proof}

\begin{remark}\label{H1action} The full group of gauge transformations preserves the reducible set, so $\mathcal G_E/\mathcal G^e_E \cong H^1(Y;\mathbb Z/2)$ acts on $\text{Red}(Y,E)$. If $\beta: H^1(Y;\mathbb Z/2) \to H^2(Y;\mathbb Z)$ is the Bockstein homomorphism, then in the notation of the above corollary, the action on $\text{Red}(Y,E)$ is given by sending $$x \cdot \{z_1, z_2\} \mapsto \{z_1 + \beta x, z_2 + \beta x\}.$$ In particular, this action is free and transitive on $\text{Red}_*(Y,E)$. 

If $x$ stabilizes $\{z_1, z_2\}$, then $z_2 = z_1 + \beta x$, and we see that $c = 2z_1 + \beta x$, and we see that we may rechoose $c$ to be $\beta x$; then the above corresponds to the reduction $$\mathbf{E} \cong \xi_{\mathbb C} \oplus \mathbb R \cong \xi \oplus \xi \oplus \mathbb R,$$ where $\xi$ is the real line bundle corresponding to $x$. 

The discussion after Lemma \ref{stab1} shows that the $H^1(Y;\mathbb Z/2)$ action has stabilizer equal to $\mathbb Z/2$ at an $O(2)$ connection, and otherwise acts freely on the reducible set. So in fact we see that $H^1(Y;\mathbb Z/2)$ acts freely on the set of reducible components that do not contain an $O(2)$-connection, and has stabilizer $\mathbb Z/2$ on those that do.
\end{remark}

We summarize the content of this section as notation: 

\begin{definition}\label{red-def}Let $Y$ be a closed oriented 3-manifold and $E$ an $SO(3)$-bundle over $Y$. We write $\textup{Red}(Y,E)$ for the set of connected components of the reducible subspace of $\widetilde{\mathcal B}^e_E$. This may be written as $$\textup{Red}(Y,E) = \textup{Red}_{*}(Y,E) \sqcup \textup{Red}_{SO(2)}(Y, E),$$ where the first term refers to those components containing a fully reducible orbit, and the second refers to those components entirely consisting of framed connections whose stabilizers are conjugate to $SO(2)$. 

If we fix a choice of $c \in H^2(Y;\mathbb Z)$ that reduces mod $2$ to $w_2(E)$, we are furnished with a bijection between $\textup{Red}(Y,E)$ and the set of unordered pairs $\{z_1, z_2\} \subset H^2(Y;\mathbb Z)$ with $z_1 + z_2 = c$. The set $\textup{Red}_*(Y,E)$ is sent to the 1-element sets\footnote{that is, those pairs with $z_1 = z_2$; equivalently the set of $z$ with $2z = c$}, and the set $\textup{Red}_{SO(2)}(Y,E)$ is sent to the 2-element sets $\{z_1, z_2\}$ with $z_1 \neq z_2$.

We write elements of $\textup{Red}(Y,E)$ as $\{z_1, z_2\}$, making the choice of $c$ implicit.
\end{definition}

In the case that $b_1(Y) = 0$, we determine the cardinality of these sets as a simple exercise in algebraic topology.

\begin{proposition}\label{red-enumerate}Let $Y$ be a closed oriented 3-manifold with $b_1(Y) = 0$, and $E$ an $SO(3)$-bundle over it. We may determine the number of reducible components as follows.\begin{itemize}

\item If $E$ is trivial, then \begin{align*}|\textup{Red}_*(Y,E)| &= |H^1(Y;\mathbb Z/2)| \\ \big|\textup{Red}_{SO(2)}(Y,E)\big| &= \left(|H^2(Y;\mathbb Z)| - |H^1(Y;\mathbb Z/2)|\right)/2.\end{align*}

\item If $E$ is nontrivial, $\textup{Red}_*(Y,E)$ is empty, and $$\big|\textup{Red}_{SO(2)}(Y,E)\big| = |H^2(Y;\mathbb Z)|/2.$$
\end{itemize}\end{proposition}

\begin{proof}If $E$ is trivial, we choose the integral lift $\lambda$ of $w_2(E) = 0$ to be $0$. Then $\text{Red}_*(Y,E)$ is in bijection with 2-torsion elements of $H^2(Y;\mathbb Z)$. The Bockstein exact sequence $$0 = H^1(Y;\mathbb Z) \xrightarrow{\mod 2} H^1(Y;\mathbb Z/2) \xrightarrow{\beta} H^2(Y;\mathbb Z) \xrightarrow{\times 2} H^2(Y;\mathbb Z)$$ shows that the set of 2-torsion elements is in bijection with $H^1(Y;\mathbb Z/2)$, the bijection sending $x \in H^1(Y;\mathbb Z/2)$ to $\beta x \in H^2(Y;\mathbb Z)$. The set $\text{Red}_{SO(2)}(Y,E)$ is in bijection with the set of pairs $$\{z, -z\} \in H^2(Y;\mathbb Z)$$ for $z$ not 2-torsion. This proves both equalities.

If $E$ is nontrivial, and $c$ is an integral lift of $w_2 E$, then a set $\{z_1, z_2\}$ with $z_1 + z_2 = c$ cannot be a singleton: if $c = 2z_1$, then $c$ reduces to $0$ mod $2$, and thus $w_2 E = 0$, so $\text{Red}_*(Y,E)$ is empty. Thus every pair $(z_1, z_2)$ with $z_1 + z_2 = c$ consists of distinct elements, and the $\mathbb Z/2$ action swapping $z_1$ and $z_2$ has no fixed points. Because a pair $(z_1, z_2)$ is determined completely by $z_1 \in H^2(Y;\mathbb Z)$, the desired equality is proven.

As a sanity check, observe that $|H^2(Y;\mathbb Z)|$ is divisible by $|H^1(Y;\mathbb Z/2)|$, which is a power of 2, and we assumed $H^1(Y;\mathbb Z/2)$ to be nontrivial in assuming that $w_2(E)$ is nontrivial, so $|H^2(Y;\mathbb Z)|/2$ is indeed an integer.
\end{proof}

\section{Configurations on the cylinder}\label{sec:config-cyl}
There are two more kinds of configuration spaces we should consider. First, we should have a configuration space of connections on the bundle $\pE$ over $\mathbb R \times Y$ pulled back from $Y$ by projection to the first factor. We write $\mathcal A^{(4)}_{\pE}$ as the space of connections which, when restricted sufficiently far on each end, are pullbacks of connections on $Y$. (We will denote any gauge group or space of connections on a 4-manifold with a (4) unless there is no risk of confusion.) When we later give a Hilbert manifold modification of this construction, it will instead be replaced by an appropriate space of connections which have exponential decay to certain constant connections on the end. The space of framed connections on $\pE$ is $$\widetilde{\mathcal A}^{(4)}_{\pE} := \mathcal A^{(4)}_{\pE} \times (\pE)_{(0,b)},$$ with the framing chosen at $(0,b) \in \mathbb R \times Y$. The gauge group acting on this is $\mathcal G^{(4),e}_{\pE}$ --- even gauge transformations which are constant on the ends. As usual, the framed configuration space is the space $$\widetilde{\mathcal B}^{(4),e}_{\pE} := \widetilde{\mathcal A}^{(4)}_{\pE}/\mathcal G^{(4),e}_{\pE},$$ with the right $SO(3)$ action on $\pE_{(0,b)}$. The `restriction to the ends' map is slightly more complicated than what one might do naively: there is a map $$\text{ev}_-: \widetilde{\mathcal A}^{(4)}_{\pE} \to \widetilde{\mathcal A}_{E}$$ given by sending $(\mathbf{A}, p)$ to $(\mathbf{A}_{-\infty}, \gamma_{-\infty}(\mathbf{A}, p))$, where in the latter term we send $p$ to its parallel transport under $\mathbf{A}$ along $\gamma$ traversed backwards from $(0,b)$. When the constant framing of the path (assigning the same $p \in \pE_{(t,b)}$ for all $t$) is $\mathbf{A}$-parallel --- if $\mathbf{A}$ is in temporal gauge, for instance --- this map is particularly simple: $\gamma_{-\infty}(\mathbf{A}, p) = p$. This map is equivariant on the left under the restriction-to-$(-\infty)$ map $$\mathcal G^{(4),e}_{\pE} \to \mathcal G^e_E,$$ and on the right under the $SO(3)$-action, because parallel transport is an isomorphism $$\pE_{(0,b)} \to \pE_{(-\infty,b)}$$ of right $SO(3)$-sets. There is also a corresponding map $\text{ev}_+$.

Thus the evaluation maps descend to a map of right $SO(3)$-spaces $$[\text{ev}_-] \times [\text{ev}_+]: \widetilde{\mathcal B}^{(4),e}_{\pE} \to \widetilde{\mathcal B}^e_{E} \times \widetilde{\mathcal B}^e_{E}.$$

\noindent We are most frequently interested in the spaces of connections on $\mathbb R \times Y$ with specified limits $[\text{ev}_{\pm}] \in \alpha_\pm$ inside previously specified $SO(3)$-orbits. Furthermore, we should pick a relative homotopy class $z \in \pi_1(\mathcal B^e_E, [\alpha_-], [\alpha_+])$. Because $\mathcal A_E$ is contractible and $\mathcal G^e_E$ acts with connected stabilizers, these are in bijection with $$\pi_0 \mathcal G^e_E \cong H^3(Y;\pi_3 \widetilde{\text{Aut}}(E_b)) \cong \mathbb Z$$ by an obstruction theory argument. Given limiting framed connections $[A_\pm,p_\pm] \in \alpha_\pm$, we choose lifts (denoted $A_-$ and $A_+$) so that a path from $A_-$ to $A_+$ in $\mathcal A_E$ projects to the relative homotopy class $z$ in $\widetilde{\mathcal B}^e_E$. We then define $$\widetilde{\mathcal A}_{\pE,z}(A_-, A_+) = \{(\mathbf{A},p) \in \widetilde{\mathcal A}^{(4)}_{\pE} \mid \mathbf{A}_{\pm \infty} = A_{\pm \infty}\}.$$ (Note that the endpoint conditions only include the connection and nothing about the framing.) The natural gauge group here is the group of even gauge transformations that are harmonic on the ends: $$\mathcal G^{e,h}_{\pE}(A_-, A_+) = \{\sigma \in \mathcal G^{(4),e}_{\pE} \mid d_{\pi^*A_{\pm}} \sigma = 0 \text{ on the corresponding end}\}.$$ The following is straightforward, and justifies this simplification. 

\begin{lemma}The natural map $$\widetilde{\mathcal A}_{\pE, z}(A_-, A_+)/\mathcal G^{e,h}_{\pE}(A_-, A_+) \to [\textup{ev}]^{-1}(\alpha_- \times \alpha_+)$$ is an isomorphism of $SO(3)$-spaces, where the final space denotes the corresponding subset of $\widetilde{\mathcal B}^{(4),e}_{\pE}.$ \end{lemma} 

While the latter space is clearly the object of interest, the former is easier to describe, so we prefer to use it as our working definition of the configuration space of framed connections between two orbits. We record this as a definition:

\begin{definition}\label{naiveconf}If $\alpha_\pm$ are $SO(3)$-orbits in $\widetilde{\mathcal B}^e_E$, the \textit{naive configuration space}\footnote{We use the term naive to contrast this space of smooth connections which are constant near $\infty$ to the later \textit{Hilbert} spaces of connections which decay exponentially at $\infty$.} of smooth framed connections from $\alpha_-$ to $\alpha_+$ is 
$$\widetilde{\mathcal B}^e_{\pE,z}(\alpha_-, \alpha_+) := \widetilde{\mathcal A}_{\pE,z}(A_-, A_+) / \mathcal G^{e,h}_{\pE}(A_-, A_+),$$ 
where $A_-$ and $A_+$ are representatives for the connections in the $\alpha_\pm$, and a path from $A_-$ to $A_+$ projects to the homotopy class $z$. This carries a continuous $SO(3)$ action, acting by translation on $\pE_{(0,b)}$ on the right. Identifying 
$$\alpha_\pm = [A_\pm, p] \cong E_b/\Gamma_{[A_\pm]}$$ 
with $p$ varying over all possible framings at the basepoint, the endpoint maps 
$$[\textup{ev}_{\pm}]: \widetilde{\mathcal B}^e_{\pE,z}(\alpha_-, \alpha_+) \to \alpha_\pm$$ 
are induced by $\widetilde{\mathcal A}_{\pE}(A_-, A_+) \to E_b,$ sending $(\mathbf{A}, p)$ to $\gamma_{\pm \infty}(\mathbf{A}, p)$, respectively. There is a left translation action of $\mathbb R$, induced by the family of diffeomorphisms $\tau_t: \mathbb R \to \mathbb R$ with $\tau_t(s) = s+t$, given by sending $$\tau_t^*(\mathbf{A}, p) = (\tau_t^* \mathbf{A}, \gamma_t(\tau_t^*\mathbf{A}, \tau_t^*p)).$$ In the framing coordinate, we are parallel transporting $\tau_t^*p$ forward by $t$ along $\gamma$, so that it is sent back to the basepoint.
\end{definition}

\begin{remark}\label{rmk:t-action}
Suppose $\mathbf{A}$ in temporal gauge has mass $m_{\mathbf A}(x) = |F_{\mathbf{A}(x)}|$ as a function of $\Bbb R$. Then the action of $t \in \Bbb R$ sends $m_{\tau_t^* \mathbf A}(x) = m_{\mathbf A}(x - t)$. That is, as $t$ goes to infinity, the mass of $\tau_t^* \mathbf{A}$ goes to $-\infty$, and vice versa.
\end{remark}

We have a corresponding version of Proposition \ref{action1} for the configuration spaces on cylinders. We do not repeat the setup of the determinant-1 $U(2)$-model, but refer to it freely in the proof below. Lemma \ref{stab1} still gives that the only possible stabilizers of elements in the 4-dimensional configuration space are the trivial group, the full group $SO(3)$, and circle subgroups conjugate to $SO(2)$.

\begin{proposition}\label{action2}Let $E$ be an $SO(3)$-bundle over a 3-manifold $Y$, and $\alpha_\pm$ are $SO(3)$-orbits in $\widetilde{\mathcal B}^e_E$. We have the following descriptions of the reducible subspaces in the configuration space $\widetilde{\mathcal B}^e_{\pE,z}(\alpha_-, \alpha_+)$ under the $SO(3)$ action: \begin{enumerate}
\item If the $\alpha_\pm$ are not both reducible orbits belonging to the same component $\textup{Red}(Y,E)$, or after concatenating a path from $\alpha_+$ to $\alpha_-$ in the (simply connected) reducible component they lie in, the homotopy class $z$ is nontrivial, then $SO(3)$ acts freely on the configuration space.

\item If $\alpha_\pm \in \textup{Red}_{SO(2)}(Y,E)$ are $SO(2)$-reducible orbits lying in the same connected component of reducibles which contains no fully reducible point, labelled by $\{z_1, z_2\} \subset H^2(Y)$, and the homotopy class $z$ is trivial, then the reducible subspace is a fiber bundle with base $\mathcal B_{\eta_1}(\alpha_-, \alpha_+)$ and fiber $S^2$, and in particular consists of one connected component. There are no fully reducible points.

\item If $\alpha_- = \alpha_+$ are the same unique fully reducible orbit in a fixed component $\{z_1, z_1\} \in \textup{Red}_*(Y,E)$, and $z$ is the trivial homotopy class, there is a unique fully reducible point in the configuration space on the cylinder, lying in the unique connected component of reducible orbits.
\end{enumerate}
\end{proposition}

\begin{proof}It's easier to find the reducibles inside the larger space $\widetilde{\mathcal B}^{(4),e}_{\pE}$ and afterwards take the intersection with $[\text{ev}]^{-1}(\alpha_- \times \alpha_+)$. As in Proposition \ref{action1}, working with connections on $\pi^*\tilde E$ with fixed trace, a connection representing an $SO(2)$-fixed point has holonomy inside $U(1) \times U(1)$ and hence induces a global splitting $\pi^* \tilde E \cong \eta_1 \oplus \eta_2$. As before, this gives a natural correspondence between the $SO(2)$-fixed points and $\bigsqcup_\eta \mathcal B^{(4)}_{\eta_1} \sqcup \mathcal B^{(4)}_{\eta_2}$ when $\eta_1 \neq \eta_2$ and with the single space $\mathcal B^{(4)}_\eta$ when $\eta_1 = \eta_2$. In the first case the involution given by the Weyl group action swaps the two spaces and in the second case it acts with a single fixed point. 

Write $\widetilde{\mathcal B}^{(4),\text{red}}_\eta$ for the subspace of $\widetilde{\mathcal B}^{(4),e}_{\pE}$ consisting of reducible connections for which the corresponding line bundle is of topological type $\eta$; in particular, we do not specify the stabilizer. If $\eta$ is an $SO(2)$-reducible component, the above shows that the map $$\widetilde{\mathcal B}^{(4),\text{red}}_\eta \to \widetilde{\mathcal B}^{(4),\text{red}}_\eta/SO(3)$$ is a constant rank submersion, with fiber $S^2$. The base may be identified with $\mathcal B^{(4)}_\eta$: it is identified as the $SO(2)$-fixed point space modulo the Weyl group action.

To calculate the intersection with $[\text{ev}]^{-1}(\alpha_- \times \alpha_+)$, observe that the line bundles $\eta_i$ restrict to isomorphic line bundles on both ends (because the inclusion $Y \hookrightarrow \mathbb R \times Y$ is a homotopy equivalence), and so the pair of line bundles associated to $\alpha_-$ and $\alpha_+$ must be the same. Furthermore, after being put in temporal gauge, the relative homotopy class $z$ traced out by one of these induced connections lies inside the simply connected space $\widetilde{\mathcal B}^{e}_{\eta}(\alpha_-, \alpha_+)$. This gives us the restrictions in (2)-(3).
\end{proof}

\begin{remark}If the $SO(3)$-bundle $E$ is nontrivial, there are no full reducibles, and so only cases (1) and (2) arise.
\end{remark}

When the limiting connections are flat, we can distinguish the relative homotopy classes $z$ by the values of the integral calculating the ``relative Pontryagin number" of the connection $\bf{A}$ by $$p_1(\mathbf{A}) := \frac{1}{8\pi^2} \int_{\mathbb R \times Y} F_{\bf A}^2;$$ this takes on a well-defined value modulo $\mathbb Z$, and the different discrete values it can take in $\mathbb R$ parameterize the set of connected components $\pi_0\widetilde{\mathcal B}^e_{\pE}(\alpha_-, \alpha_+)$. When the limits are not flat, one could pick a base connection $A_0$ on the cylinder interpolating between them and compute $\int (F_A^2 - F_{A_0}^2)$; alternatively, later the limits will be critical points of a perturbed Chern-Simons functional and we will define the action of a connection in terms of this perturbed functional.

\section{Configurations on a cobordism}\label{sec:config-4d}
Suppose we are given a complete Riemannian manifold $W$ with two cylindrical ends (i.e., specified isometries to $(-\infty, 0] \times Y_1$ and $[0, \infty) \times Y_2$) and an $SO(3)$ bundle $\bf E$ with specified isomorphisms over the ends to bundles $E_i \to Y_i$). We think of this as a cobordism $Y_1 \to Y_2$. Here $E_i \to Y_i$ is an $SO(3)$ bundle, framed over a basepoint $b_i \in Y_i$. Furthermore, we pick a smooth embedded path $\gamma: \mathbb R \to W$ with $\gamma(t) = (b_1, t)$ for $t < 0$ sufficiently small and $\gamma(t) = (b_2, t)$ for $t$ sufficiently large, and a trivialization of $\gamma^*E$ restricting to the given trivializations on the ends. This path $\gamma$ will serve the role $\mathbb R \times \{b\}$ did for the cylinder. Given orbits $\alpha_i \subset \widetilde{\mathcal B}^e_{E_i}$, we construct the configuration space $\widetilde{\mathcal B}^e_{\mathbf{E}, z}(\alpha_1, \alpha_2)$ much as we did before. The basepoint is instead $\gamma(0)$, and the framing portion of the endpoint maps given by parallel transport to $\pm \infty$ along $\gamma$. There is still a decomposition of the configuration space into disconnected components, which we still label by $z$, but this no longer has a description in terms of relative homotopy classes. Our above discussion on relative Pontryagin numbers, however, still does, as will the definition of an action using the perturbed Chern-Simons functional.

As in the previous sections, we analyze the reducible points in terms of the connected components they lie in. 

\begin{definition}\label{red-def-W}Let $W$ be a 4-manifold with two cylindrical ends, as above, equipped with an $SO(3)$-bundle $\mathbf{E}$. We write $\textup{Red}(W,\mathbf{E})$ to denote the set of connected components of the reducible subspace of $\widetilde{\mathcal B}^e_{\mathbf{E}}$. This may be written as the disjoint union $$\textup{Red}_*(W,\mathbf{E}) \sqcup \textup{Red}_{SO(2)}(W,\mathbf{E}),$$ where the former denotes components that include some fully reducible orbit. If $\alpha_-$ and $\alpha_+$ label orbits in the configuration spaces $\widetilde{\mathcal B}^e_{E_i}$ of the ends, we write $\textup{Red}(W,\mathbf{E})(\alpha_-, \alpha_+)$ for the set of connected components of the reducible subspace $\widetilde{\mathcal B}^{e}_{\mathbf{E}}(\alpha_-, \alpha_+)$.
\end{definition}

We enumerate the reducible components in the following.

\begin{proposition}\label{action3}Suppose $W$ is a Riemannian 4-manifold equipped with $SO(3)$ bundle $\bf E$, with one incoming cylindrical end $(Y_1, E_1)$ and one outgoing cylindrical end $(Y_2,E_2)$. The reducible subspaces of the $SO(3)$-action on the configuration space $\widetilde{\mathcal B}^{e}_{\mathbf{E}}(\alpha_-, \alpha_+)$ is as follows.\begin{enumerate}
    \item If either of the $\alpha_\pm$ are irreducible, or $\beta w_2 \mathbf{E} \in H^3(W;\mathbb Z)$ is nonzero, then there are no reducible points in the configuration space.

    \item If $\beta w_2 \mathbf{E} = 0$, then fix an integral lift $c$ of $w_2 \mathbf{E}$, and use the restriction of $c$ to the ends to determine the bijection between $\textup{Red}(Y_i, E_i)$ and 2-element sets of cohomology classes in $H^2(Y_i; \mathbb Z)$ as in Definition \ref{red-def}. Then the connected components $\textup{Red}_{SO(2)}(W,\mathbf{E})$ are in bijection with pairs $\{z_1, z_2\} \subset H^2(W;\mathbb Z)$ with $z_1 + z_2 = c$ and $z_1 \neq z_2$. The set $\textup{Red}_{SO(2)}(W,\mathbf{E})(\alpha_-, \alpha_+)$ is the subset of $\textup{Red}_{SO(2)}(W,\mathbf{E})$ consisting of pairs that restrict to the $\{x_1, x_2\}$ corresponding to each $\alpha_\pm$.\footnote{Note in particular that taking the intersection of the reducible subspace with $[\text{ev}]^{-1}(\alpha_- \times \alpha_+)$ --- that is, specifying the limits --- does not increase the number of connected components. Furthermore, either of $\alpha_\pm$ may be $SO(2)$-reducible or fully reducible.}

    \item If $\mathbf{E}$ is nontrivial there are no fully reducible points. Otherwise, fix a trivialization; this produces a bijection of $\textup{Red}_{*}(Y_i, E_i)$ with 2-torsion cohomology classes in $H^2(Y_i;\mathbb Z)$. Then the components $\textup{Red}_*(W,\mathbf{E})$ are in bijection with 2-torsion cohomology classes in $H^2(W;\mathbb Z)$. If $\alpha_\pm$ are both fully reducible orbits, then $\textup{Red}_*(W,\mathbf{E})(\alpha_-, \alpha_+)$ is the subset of those 2-torsion classes in $H^2(W;\mathbb Z)$ that restrict to the 2-torsion classes in $H^2(Y_i;\mathbb Z)$ labelled by the $\alpha_\pm$. In each component with a fully reducible point, there is a unique such point.\footnote{It may be that $W$ has more full reducibles than the $Y_i$; this is true if, for instance, each $Y_i$ is an integer homology sphere but $H_1(W;\mathbb Z)$ has 2-torsion.}
\end{enumerate}
\end{proposition}

\begin{proof}If $\beta w_2 \mathbf{E} \neq 0$, then $\mathbf{E}$ cannot be written as the direct sum of a trivial line bundle and an oriented 2-plane bundle: $\beta w_2(\mathbf{E}) = e(\mathbf{E})$, the Euler class, which is the obstruction to finding a nonvanishing section; but an $SO(2)$-reducible point induces such a splitting (there is a parallel section). Note that this is also the obstruction to finding an integral lift of $w_2(\mathbf{E})$. So suppose now that $c \in H^2(W;\mathbb Z)$ is an integral lift of $w_2 \mathbf{E}$.

Almost identically to Proposition \ref{action2}, the $SO(2)$-fixed point space of $\widetilde{\mathcal B}^{e}_{\bf E}$ is identified with the disjoint union over configuration spaces of connections on line bundles over $W$: $\sqcup_{\eta} \mathcal B^{W}_{\eta}$. If $\eta_1 \oplus \eta_2$ has first Chern class $c$, then the Weyl group acts on $\mathcal B_{\eta_i}$ by sending a connection $[A]$ on one to the connection $[A - A_0]$ on the other, where $A_0$ is a fixed connection on the complex line bundle $\lambda$. Quotienting by the Weyl group we are only left with components labelled by pairs $\{z_1, z_2\}$ with $z_1 + z_2 = c$, and if $z_1 = z_2$, the Weyl group fixes a unique point in $\mathcal B^W_\eta$, the corresponding fully reducible point.
\end{proof}

\section{Tangent spaces}\label{sec:config-tangent}
Before we introduce the Hilbert manifold versions of these constructions, we investigate the tangent spaces and differentials involved in the naive constructions (where all connections and gauge transformations are smooth). 

The space $\mathcal A_E$ of all smooth connections on $E$ is an affine space over $\Omega^1(\mathfrak g_E)$; the natural Riemannian metric on $T_A \mathcal A_E = \Omega^1(\mathfrak g_E)$ is the (incomplete!) $L^2$ inner product on forms, using a bi-invariant inner product on the Lie algebra bundle $\mathfrak g_E$ to define the Hodge star. With this, we can define $\nabla \text{cs}$ as a vector field on $\mathcal A_E$; we may identify $(\nabla \text{cs})(A) = *F_A$. The action of the gauge group $\mathcal G^e_E$ at $A \in \mathcal A_E$ has differential 

$$d_A: \Omega^0(\mathfrak g_E) \cong T_e \mathcal G^e(E) \to T_A \mathcal A_E \cong \Omega^1(\mathfrak g_E).$$

The kernel of this map is the space of $A$-harmonic sections of $\mathfrak g_E$. The differential of the action of $\mathcal G^e_E$ on the fiber $E_b$ (whose tangent spaces are all canonically identified with $\mathfrak g_b$, the fiber of $\mathfrak g_E$ at $b$) is given by evaluation at the basepoint $\text{ev}_b: \Omega^0(\mathfrak g_E) \to \mathfrak g_b$. We write 
$$\tilde d_A: \Omega^0(\mathfrak g_E) \to \Omega^1(\mathfrak g_E) \oplus \mathfrak g_b$$ 
for the differential of the action of $\mathcal G^e_E$ on $\widetilde{\mathcal A}_E$ at $(A,p)$. Unfortunately, there is no natural identification of the cokernel of $\tilde d_A$ --- which we hope to eventually identify with $T_{[A,p]} \widetilde{\mathcal B}^e_E$ --- with the $L^2$ orthogonal complement of its image, since the image is not $L^2$ closed. Its closure is $\text{Im}(d_A) \oplus \mathfrak g_b$, the essential difficulty here being that point-evaluation is not $L^2$-continuous. Rather, this orthogonal complement is precisely the same as $\text{ker}(d_A^*)$. However, we can identify $$(\text{Im}(d_A) \oplus \mathfrak g_b)/\text{Im}(\tilde d_A) \cong \mathfrak g_A^\perp,$$ the orthogonal complement inside $\mathfrak g_b$ of the subspace $\mathfrak g_A$ of elements that extend to $A$-parallel sections. Equivalently, $\mathfrak g_A$ is the tangent space to the image of $\Gamma_A$ in $\text{Aut}(E_b)$. 

Picking a choice of (necessarily non-orthogonal) complement of $\text{Im}(\tilde d_A)$ inside $\text{Im}(d_A) \oplus \mathfrak g_b$ to identify with $\mathfrak g_A^\perp$, we may still decompose 
$$T_{(A,p)} \widetilde{\mathcal A}_E = \text{ker}(d_A^*) \oplus \text{Im}(\tilde d_A) \oplus \mathfrak g_A^\perp.$$ 
We will henceforth write this as 
$$\mathcal T_{A,p} = \mathcal K_A \oplus \mathcal G_A \oplus \mathfrak g_A^\perp,$$ 
giving us the isomorphism 
$$T_{[A,p]}\widetilde{\mathcal B}^e_E \cong \mathcal K_A \oplus \mathfrak g_A^\perp.$$ 
This can not be a decomposition into a direct sum of locally trivial bundles (the last factor has dimension varying between $0, 2,$ and $3$), but it is when restricted to the submanifold of $\widetilde{\mathcal B}^e_E$ consisting of framed connections with stabilizer of \emph{fixed} conjugacy class (irreducible, conjugate to $SO(2)$, or fully reducible), where $\mathfrak g_A^\perp$ actually defines a locally trivial vector bundle.

Because $\text{cs}$ is, up to a constant, fixed under the action of $\mathcal G^e_E$ and not just $\mathcal G^{e,b}_E$, its $L^2$ gradient on $\mathcal A_E$ must take values in $\mathcal K_A$. It is clear from $$\nabla \text{cs}(A) = *F_A$$ that the critical points of $\text{cs}$ are precisely the flat connections.

If $A_-$ and $A_+$ are connections on the ends of $\pE$ over $\mathbb R \times Y$ (and $\alpha_\pm$ the corresponding orbits in the configuration space), we constructed the 4-dimensional configuration space in Definition \ref{naiveconf} as a quotient of $\widetilde{\mathcal A}_{\pE,z}(A_-, A_+)$. If we pick a base connection $A_0$ in this space, we get a description 
$$\widetilde{\mathcal A}_{\pE,z}(A_-, A_+) = \{A_0 + \Omega^1_c(\pi^*\mathfrak g_E)\} \times (\pE)_{0,b},$$
and in particular of the space of unframed connections as an affine space. The Lie algebra of the gauge group $\mathcal G^{e,h}_{\pE}(A_-, A_+)$ is $\Omega^0_h(\mathfrak g_E)$, the space of sections of $\mathfrak g_E$ that are $A_\pm$-harmonic on the corresponding ends. The differential of the gauge group action at some $(\mathbf{A}, p)$ is given by $\sigma \mapsto (d_{\mathbf{A}} \sigma, \sigma(b)).$ Picking a complement $\mathfrak g_A^\perp$ of $\text{Im}(\tilde d_{\mathbf{A}})$, we again have our decomposition 
$$T_{(\mathbf{A},p)} \widetilde{\mathcal A}_E = \text{ker}(d_{\mathbf A}^*) \oplus \text{Im}(\tilde d_{\mathbf A}) \oplus \mathfrak g_{\mathbf A}^\perp.$$ 

We also have the endpoint maps to $(A_\pm, E_b) \cong SO(3)$. Before modding out by gauge, the positive endpoint map sends $\text{ev}_+(\mathbf{A},p) = \gamma_\infty(\mathbf{A},p)$, the parallel transport of $p$ to $\infty$ using $\mathbf{A}$. The differential of this map at $(\mathbf{A},p)$ is a linear map 
$$\Omega^1_c(\mathfrak g_E) \oplus \mathfrak g_{0,b} \to \mathfrak g_b,$$ 
equal to the integral 
$$\int_0^\infty \gamma^*_{\mathbf A}\sigma(t,b)dt.$$ 

Here we use parallel transport backwards along $\gamma$ using $\mathbf{A}$ so that $\sigma(t,b)$ is taken to the single vector space $\mathfrak g_{0,b}$, where it makes sense to integrate. Observe that if $\mathbf{A}$ is in temporal gauge, this is just $\int_0^\infty \sigma(t,b)dt$.

Nothing changes in the above description when passing to a general Riemannian 4-manifold $W$ with cylindrical ends. 

\section{Hilbert manifold completions}\label{sec:config-Hilbert-3d}
To have a useful transversality theory, we must replace the naive ``infinite-dimensional manifolds" $\mathcal A,\; \mathcal G$, and so on with certain completions which are manifolds modeled on Hilbert spaces. If $E$ is an $SO(3)$-bundle over a 3-manifold, define $\mathcal A_{E,k}$ to be the set of $L^2_k$ connections on $E$; more precisely, if $\Omega^1_k$ denotes the Hilbert space of $L^2_k$ 1-forms, then 
$$\mathcal A_{E,k} = A_0 + \Omega^1_k(\mathfrak g_E)$$
for any choice of smooth base connection $A_0$. The resulting set does not depend on the choice of smooth connection. Similarly, we may define 
$$\mathcal G_{E,k+1} = \{\sigma \in L^2_{k+1}(\text{End}(E)) \mid \sigma(x) \text{ is an isomorphism for all } x\},$$ 
given the constraint that $k+1 \geq 2$ so that gauge transformations are automatically continuous and we may take point evaluations as in the definition. Then $\mathcal G^e_{E,k+1}$ is the subset of sections which lift continuously to $\widetilde{\text{Aut}}(E)$. As in \cite[Section~4.2]{DK}, the space $\mathcal G^e_{E,k+1}$ becomes a Banach Lie group, and we have a smooth action of $\mathcal G^{e}_{E,k+1}$ on $\widetilde{\mathcal A}_{E,k} = \mathcal A_{E,k} \times E_b$; the same argument as the case of smooth gauge transformations shows that this action is free. The following is \cite[Proposition~3.1.5.10]{Bou}. To facilitate comparison to the version stated there, note that Bourbaki's definition of immersion includes the assumption that the differential has closed image with a closed complement.

\begin{proposition}Suppose $G$ is a Banach Lie group and $X$ a Banach manifold on which $G$ acts. Suppose that $G$ acts freely and properly on $X$, and such that if $\rho(x): T_eG \to T_xX$ is the differential of the action at the point $x$, the image of $\rho(x)$ is closed and has some closed complement for all $x$. Then the quotient topology on $X/G$ is Hausdorff, and there is a unique smooth structure on $X/G$ such that $\pi: X \to X/G$ is a submersion. Furthermore, $X \to X/G$ is a principal $G$-bundle.
\end{proposition}

We should see that all these assumptions apply to the action of $\mathcal G^e_{E,k+1}$ on $\widetilde{\mathcal A}_{E,k}$, giving us the quotient manifold $\widetilde{\mathcal B}^e_{E,k}$. We record this as a lemma.

\begin{lemma}$\widetilde{\mathcal B}^e_{E,k}$ is a Hausdorff Hilbert manifold if $k \geq 1$.
\end{lemma}

\begin{proof}Properness of this action would follow from properness of the action on the unframed space of connections, because the map $\widetilde{\mathcal A}_{E,k} \to \mathcal A_{E,k}$ is proper (because $E_b$ is compact). That the action of $\mathcal G_{E,k+1}$ on $\mathcal A_{E,k}$ is proper is proved in \cite[Section~2.3.7]{DK}, by a bootstrapping argument. 

Because $\widetilde{\mathcal A}_{E,k}$ is a Hilbert manifold, we only need to verify that the image of the differential is closed; write this explicitly as the sum 
$$(d_A, \text{ev}_b): \Omega^0_{k+1}(\mathfrak g_E) \to \Omega^1_k(\mathfrak g_E) \oplus \mathfrak g_b.$$ 
First note that $\text{Im}(d_A)$ is closed: there is an $L^2$-orthogonal decomposition 
$$\Omega^1_k(\mathfrak g_E) = \text{ker}(d_A^*) \oplus \text{Im}(d_A);$$ 
that this is true in $\Omega^1_k$ and not just $\Omega^1_0$ follows from elliptic regularity (if $d_A \sigma = \eta$ where $\eta \in L^2_k$, then $\sigma \in L^2_{k+1}$).

Then suppose $(\eta_n, g_n) \to (\eta, g)$ is a convergent sequence, for which there is a sequence of $\sigma_n \in \Omega^0_{k+1}$ with $d_A \sigma_n = \eta_n$ and $\sigma_n(b) = g_n$. Because $\text{Im}(d_A)$ is closed, there is some $\sigma$ with $\eta = d_A \sigma$. Thus there is some sequence $\psi_n \in \text{ker}(d_A)$ so that $\sigma_n + \psi_n \to \sigma$. As long as $k+1 \geq 2$, point-evaluation is continuous, and so $g_n + \psi_n(b) \to g$; because $g_n \to g$, we see that $\psi_n(b) \to 0$. Thus $\sigma(b) = g$ and $d_A \sigma = \eta$ as desired.
\end{proof}

This provides somewhat less inspiring charts than the case of irreducible connections and unbased gauge group, where one can explicitly present slices for the action of the gauge group as those connections with 
$$d^*_{A_0}(A-A_0) = 0.$$ 
One can at least use the decomposition 
$$\mathcal T_{A,k} = \mathcal K_{A,k} \oplus \mathcal G_{A,k} \oplus \mathfrak g_A^\perp$$
to provide slices for the based gauge group as exponentials of $\mathcal K_{A,k} \oplus \mathfrak g_A^\perp$. In addition to the natural tangent bundle $\mathcal T\mathcal A_{E,k} = \Omega^1_k(\mathfrak g_E)$, we may also define completed tangent bundles $\mathcal T_j \widetilde{\mathcal A}_{E,k} = \Omega^1_j(\mathfrak g_E)$ for any $j \leq k$.

After passing to this completion, the analogue of Proposition \ref{action1} remains true. Whenever we construct a gauge transformation via parallel transport, it has one greater order of differentiability than the connection we used to define parallel transport.

\section{The 4-dimensional case}\label{sec:config-Hilbert-4d}
In contrast to the case of homology 3-spheres, where to define Floer's instanton homology we only needed to study instantons on the cylinder with irreducible limits on the ends, the same flavor of Sobolev completion will not suffice to define the configuration space of 4-manifolds with general limits. Instead, we will need to use \textit{weighted} Sobolev spaces. Recall $\pE$ is an $SO(3)$-bundle on $\mathbb R \times Y$, given as the pullback of some bundle $E$ on $Y$ under the projection map. Let $f_\delta$ be a smooth positive function on $\mathbb R$ such that $f_\delta(t) = e^{\delta|t|}$ for $|t| \geq 1$. 

Suppose we have fixed limiting connections $A_-$ and $A_+$ on $E$. Here we are fixing \emph{connections} and not gauge equivalence classes; we write the corresponding orbits in $\widetilde{\mathcal B}^e_E$ as $\alpha_\pm$. Pick a smooth base connection $A_0$ on $\pE$ over $\mathbb R \times Y$ that agrees with the pullback connections of $A_-$ and $A_+$ near $\pm \infty$ respectively; it traces out a relative homotopy class $z \in \pi_1(\widetilde{\mathcal B}^e_E, A_-, A_+)$. Then one may define the space of $L^2_{k,\delta}$ sections of $\pE$ as those sections for which 
$$\int_{\mathbb R \times Y} f_\delta^2\left(|\sigma|^2 + |\nabla_{A_0}\sigma|^2 + \dots |\nabla_{A_0}^k \sigma|^2\right) < \infty;$$ 
this is the $L^2_{k,\delta}$ norm with respect to $A_0$. Then we say $$\mathcal A_{\pE,z,k,\delta}(A_-, A_+) := \{A_0 + \Omega^1_{k,\delta}(\mathfrak g_E)\}.$$ This only depends on $A_0$ through the choice of $z$, but it clearly depends on the limits $A_\pm$; if two pairs $(A_\pm^i)$ of limits are gauge equivalent by gauge transformations in the same homotopy class, then there is a gauge transformation taking 
$$\mathcal A_{\pE,z,k,\delta}(A^1_-, A^1_+) \cong \mathcal A_{\pE,z,k,\delta}(A^2_-, A^2_+).$$ 
This set depends only on $\delta$, not the particular function $f_\delta$ chosen above. 

We need a Banach Lie group $\mathcal G^{e,h}_{\pE,k+1,\delta}(A_-, A_+)$. Set 
$$\Gamma_\pm = \Gamma(\alpha_\pm) \subset \mathcal G^e_{E,k+1};$$ 
this is precisely the set of $\alpha_\pm$-harmonic gauge transformations. As before, we make sense of the $L^2_{k+1,\delta}$ norm by considering gauge transformations as certain sections of $\pi^*(\text{End}(E))$. If we only cared about constant limits (as is the case when the $\alpha_\pm$ are irreducible), we could take the Lie group 
$$\mathcal G_{\pE, k+1, \delta} = \{\sigma \in \Gamma(\text{Aut}(\pE)) \subset \Gamma(\text{End}(\pE)) \mid (\sigma - 1) \in L^2_{k+1,\delta}(\text{End}(\pE)\}.$$ 
That is, roughly, gauge transformations in this group are sections of $\pi^*(E \otimes E^*)$ that are pointwise automorphisms and decay exponentially quickly to the identity. This is again a Banach Lie group that acts smoothly on $\mathcal A_{\pE,z,k,\delta}(A_-, A_+)$; this is essentially the content of \cite[Section~4.3]{Don}. We write 
\begin{align*}
\mathcal G^{h}_{\pE, k+1, \delta}&(\alpha_\pm) = \\
&\{\sigma \in \Gamma(\text{Aut}(E)) \mid \exists \psi \text{ with parallel ends and } \sigma - \psi \in L^2_{k+1,\delta}(\text{End}(E))\}.
\end{align*}
This is again a Banach Lie group, but now we have the surjective endpoint evaluation map 
$$\mathcal G^{h}_{\pE, k+1, \delta}(A_-, A_+) \to \Gamma_- \times \Gamma_+.$$ 
There is also the natural subgroup $\mathcal G^{e,h}_{\pE, k+1, \delta}(A_-, A_+)$, those transformations that lift to sections of $\widetilde{\text{Aut}}(\pE)$. Following Definition \ref{naiveconf}, we use this to define the configuration space as follows. 

\begin{definition}\label{conf}The configuration space of framed connections on $\pE$ between $\alpha_-$ and $\alpha_+$ is defined to be $$\widetilde{\mathcal B}^{e}_{\pE, z, k, \delta}(A_-, A_+) = \widetilde{\mathcal A}_{\pE, z, k, \delta}(A_-, A_+)/\mathcal G^{e,h}_{\pE, k+1, \delta}(A_-, A_+).$$ This carries an action of $SO(3)$ by acting on the fiber above the basepoint on the right. The endpoint maps given by parallel transport along $\mathbb R \times \{b\}$ are right $SO(3)$ maps to $[\alpha_\pm] \cong E_b/\Gamma_\pm$; these maps are smooth and furthermore submersions. Equivariant with respect to these is a (left) $\mathbb R$ action by pullback of the connection and parallel transport of the framing.
\end{definition}

Again, we should check that this group action is well-behaved.

\begin{lemma}$\widetilde{\mathcal B}^{e}_{\pE, z, k, \delta}(\alpha_-, \alpha_+)$ carries the natural structure of a smooth Hilbert manifold when $k \geq 2$. 
\end{lemma}
\begin{proof}
We need to check that the action on $\widetilde{\mathcal A}_{\pE,z,k,\delta}(A_-, A_+)$ is proper, and that the tangent spaces to orbits are closed. 

Let $(\sigma_n, A_n)$ be a sequence of $L^2_{k+1,\delta}$ gauge transformations and $L^2_{k,\delta}$ connections so that $\sigma_n^*A_n \to B$ in the $L^2_{k,\delta}$ topology, and so that $A_n \to A$ in the $L^2_{k,\delta}$ topology. Write $B_n = \sigma_n^* A_n$. We want to show that $\sigma_n$ then has a convergent subsequence in the $L^2_{k+1,\delta}$ topology.

On any chart with compact closure, choose an arbitrary trivialization of $\pE$, so that we represent $A_n$ and $B_n$ as matrices of 1-forms and $\sigma_n$ as a map to $SO(3)$. Then as in \cite[Section~2.3.7]{DK} we have the formula
$$d\sigma_n = \sigma_n A_n - B_n \sigma_n;$$ 
from this it follows that an $L^2_j$ bound on $\sigma_n$ implies an $L^2_{j+1}$ bound as long as $j \leq k$, and an $L^2$ bound follows from compactness of $SO(3)$ and of the chart itself. Therefore $\|\sigma_n\|_{L^2_{k+1}}$ is uniformly bounded; by compactness, we may choose a convergent subsequence $\sigma_n \to \sigma$ in $L^2_k$. But then we may identify $\sigma^*A = B$, and therefore we have $$d\sigma - d\sigma_n = \sigma A - \sigma_n A_n - B \sigma + B_n \sigma_n,$$ which may be rewritten as $$(\sigma - \sigma_n)(A) + \sigma_n(A-A_n) - B(\sigma - \sigma_n) - (B - B_n)(\sigma_n).$$ Then because all of $\sigma - \sigma_n, A - A_n$, and $B - B_n$  go to $0$ in $L^2_k$, we see that $d\sigma - d\sigma_n \to 0$ in $L^2_k$, and therefore $\sigma_n \to \sigma$ in $L^2_{k+1}$.

This shows that $\sigma_n \to \sigma$ in $L^2_{k+1,\text{loc}}$. We should check convergence on the ends. We may write $A_n$ and $B_n$ as the flat connection $\alpha$ plus an $L^2_{k,\delta}$ 1-form valued in $\mathfrak g_E$. Choose a trivialization of $E$ on a chart $U$ of $Y$, and extend that to a trivialization of $\pE$ over $\mathbb R \times U$. Consider this on one end at a time; for convenience we say $(-\infty, 0] \times U$. We write the flat connection $\alpha$ as $d + \omega_\alpha$ in this chart, where $\omega_\alpha$ is a particular matrix of 1-forms, and then the connections in this trivialization take the form $d + \omega_\alpha + A_n$ for an $L^2_{k,\delta}$ 1-form $A_n$.

The above formula now gives us 
$$d\sigma_n = \sigma_n (\omega_\alpha + A_n) - (\omega_\alpha + B_n) \sigma_n = \sigma_n A_n - B_n \sigma_n + (\sigma_n \omega_\alpha - \omega_\alpha \sigma_n);$$ 
more simply stated, this is 
$$d_\alpha \sigma_n = \sigma_n A_n - B_n \sigma_n.$$ 
For the limit $\sigma \in L^2_{k+1,\text{loc}}$, we then have $(d_\alpha \sigma) = \sigma A  - B\sigma$. Because the $L^2_{0,\delta}$ norm of $\sigma A$ agrees with that of $A$, we see that $(d_\alpha \sigma)$ has an $L^2_{0,\delta}$ bound; because $\sigma$ is the limit of gauge transformations which are asymptotic to $1$ on the ends, $\sigma$ is asymptotic to $1$ on the ends, and so $\sigma$ is $L^2_{0,\delta}$.

We also have $$d_\alpha\sigma - d_\alpha\sigma_n = \sigma(A-A_n) + (\sigma - \sigma_n)A_n - (B - B_n)(\sigma) - B_n(\sigma - \sigma_n).$$ Therefore $$\|d_\alpha \sigma - d_\alpha \sigma_n\|_{L^2_{0,\delta}} \leq \|A - A_n\| + 2\|A_n\| + \|B - B_n\| +2\|B_n\|,$$ and in particular $\|\sigma_n\|_{L^2_{1,\delta}}$ is uniformly bounded. By the compactness theorem for weighted Sobolev spaces, $\sigma_n \to \sigma$ in $L^2_{1,\delta'}$ for $\delta' < \delta$. Inducting with the previous formula, we may see that $\|\sigma_n\|_{L^2_{k+1,\delta}}$ is uniformly bounded. Therefore $\sigma_n \to \sigma$ in $L^2_{k+1,\delta'}$ for $\delta' < \delta$. The above formula again shows using Sobolev multiplication that $\|\sigma - \sigma_n\|_{L^2_{k+1,\delta}} \to 0$ on this chart, and by running this argument a few times $\sigma_n \to \sigma$ in $L^2_{k+1,\delta}$ on all of $\mathbb R \times Y$, as desired.

The extension to the entire gauge group (with gauge transformations that are exponentially decaying to $\alpha$-harmonic gauge transformations) is formal, as the quotient by the subgroup used above is compact, equal to the space of harmonic gauge transformations on the two ends.

To see that tangent spaces are closed, it suffices to show that $\Omega^1_{k,\delta}(\mathfrak g_E)$ admits a closed splitting $\text{Im}(d_{\mathbf{A}}) \oplus \text{ker}(d_{\mathbf{A}}^*)$; we want to show that every element may be written uniquely as a sum of that form.

That these are closed subspaces follows from the closed range theorem for densely defined operators (thinking of this as a densely defined operator in $L^2$): the range of $d_{\mathbf{A}}$ is closed if and only if its image is the $L^2$ orthogonal complement of $\text{ker}(d_{\mathbf{A}}^*) \cap L^2_{-k,-\delta}$. 

If $\psi \in \text{ker}(d_{\mathbf{A}}^*) \cap L^2_{-k,-\delta}$, then by separation of variables on the end, using that the signature operator on $Y$ has no nonzero eigenvalues of magnitude less than $\delta$, we see that $\psi \in L^2_{-k,\delta,\text{ext}}$: it is the sum of a section which is constant and $\alpha$- or $\beta$-harmonic on the ends and an $L^2_{-k,\delta}$ section; similarly elliptic regularity and the assumption that $\mathbf{A}$ is $L^2_{k,\delta}$ implies that $\psi \in L^2_{k,\delta,\text{ext}}$. Now we may apply the usual integration by parts trick to show that $\langle d_{\mathbf{A}} \sigma, \psi \rangle = \langle \sigma, d_{\mathbf{A}}^* \psi \rangle = 0$ for all $\sigma$.

What is left is to see that the equation $\Delta_{\mathbf{A}} \sigma = -d_{\mathbf{A}}^* \psi$ has a solution for any $\psi \in \Omega^1_{k,\delta}$; to see this, choose $\psi$ so that $d_{\mathbf{A}}^* \psi$ is in the orthogonal complement of $\text{Im} \Delta_{\mathbf{A}}$. But as before we may see that if 
$$0 = \langle \Delta_{\mathbf{A}} \sigma, d_{\mathbf{A}}^* \psi \rangle = \langle d_{\mathbf{A}} \sigma, d_{\mathbf{A}} d_{\mathbf{A}}^* \psi \rangle,$$ 
then in particular $d_{\mathbf{A}} d_{\mathbf{A}}^* \psi = 0$, and as such $d_{\mathbf{A}}^* \psi$ is parallel. Because $\psi$ is asymptotic to $0$ on the ends, we see that $d_{\mathbf{A}}^* \psi = 0$, and hence indeed
$$\Delta_{\mathbf{A}} \sigma = -d_{\mathbf{A}}^* \psi$$ 
is solvable for all $\psi$. Uniqueness follows because if $\Delta_{\mathbf{A}} \sigma = 0$, then because $\sigma$ is asymptotic to $0$ on the ends, by the usual integration by parts trick $d_{\mathbf{A}} \sigma = 0$; again because $\sigma$ is asymptotically zero, we see that $\sigma$ is globally zero. 
\end{proof}

The discussion of tangent bundles and reducible configurations (Proposition \ref{action2}) from the naive case carry over without change.

Finally, suppose we are given a complete Riemannian manifold $W$ with two cylindrical ends (i.e., specified isometries to $(-\infty, 0] \times Y_1$ and $[0,\infty) \times \overline{Y_2}$), an $SO(3)$ bundle $\bf E$ with specified isomorphisms over the ends to bundles $E_i \to Y_i$, and specified $\alpha_i$ on the $E_i$. Further suppose $W$ is equipped with an embedded path $\gamma$ in $W$, cylindrical and agreeing with specified basepoints $b_i$ on the ends. Then the construction of configuration spaces $\widetilde{\mathcal B}^{e}_{\mathbf{E},k,\delta}(\alpha_1, \alpha_2)$ with equivariant endpoint maps to $(E_i)_{b_i}/\Gamma_{\alpha_i}$ carries through with essentially no change from above. When defining cobordism maps in instanton homology, the data of the embedded path will only ultimately matter up to homotopy of such paths, which reduces to the data of a relative homotopy class $(D^1, \{0\}, \{1\}) \to (\overline W, Y_1, Y_2)$.


\chapter{Critical points and perturbations}\label{chap:3}
\section{Holonomy perturbations}\label{sec:3d-pert}
In order to carry out the construction of a Morse-like homology theory on our configuration space $\widetilde{\mathcal B}^{e}_{E,k}$, the Chern-Simons functional will usually need to be perturbed so that the critical sets are isolated orbits of $SO(3)$, and the moduli spaces of trajectories are smooth manifolds of the appropriate dimension. Here we recall the standard perturbations used in this situation, as well as their basic properties, from \cite[Section~3.2]{KM1}; they work relative to a knot $K$, but we may simply take $K = \varnothing$. These perturbations originate in Floer's definition of instanton homology and have been used consistently in the development of the subject.

Suppose we are given a collection of immersions 
$$q_i: S^1 \times D^2 \looparrowright Y, \;\;\;\; i = 1, \dots, N,$$ 
such that the $q_i$ all agree in a small neighborhood of $\{1\} \times D^2$. Further pick a conjugation-invariant smooth function $h: SO(3)^N \to \mathbb R$ and a nonnegative 2-form $\mu$ on $D^2$ that integrates to 1 and vanishes near the boundary. Taking the holonomy of a connection $A$ along the family of curves $q_i(e^{2\pi i t},z)$ parameterized by $i$ and $z$ gives a map 
$$\text{Hol}_q(A): D^2 \to SO(3)^N;$$
we are interested in the functions (called $SO(3)$-cylinder functions) 
$$f_{q,h,\mu}(A) = \int_{D^2} h(\text{Hol}_q(A)) \mu.$$ 
These functions $f_{q,h,\mu}$ have smooth extensions to the Hilbert manifolds $\mathcal A_{E,k}$, where (as on $\mathcal A_E$) they are invariant under the action of the gauge group. By ignoring the framing, we can extend this trivially to a map $f_{q,h,\mu}: \widetilde{\mathcal A}_{E,k} \to \mathbb R$ invariant under both the left action of the \textit{full} gauge group and the right action of $SO(3)$. Thus they descend to $SO(3)$-invariant smooth functions $f_{q,h,\mu}: \widetilde{\mathcal B}^e_{E,k} \to \mathbb R$. Note that they are furthermore invariant under the remaining $H^1(Y;\mathbb Z/2)$ action. 

We will construct a version of cylinder functions which breaks the $H^1(Y; \mathbb Z/2)$ symmetry. In the context of the $U(2)$-model, we choose an integral cohomology class $c \in H^2(Y;\mathbb Z)$ reducing mod 2 to $w_2 \mathbf{E}$ and hence a lift of $\mathbf{E}$ to the $U(2)$-bundle $\lambda \oplus \mathbb C$, where $c_1 \lambda = c$. Fixing a base connection $A_0$ on $\lambda$, the configuration space is the set of connections for which $\text{tr}(A) = A_0$. Fixing further a framing of the $U(2)$-bundle at $b$, associated to any connection with fixed trace is the holonomy map 
$$\text{Hol}_A: \Omega_b Y \to U(2),$$
which is invariant under the action of the group of determinant-$1$ gauge transformations based at $b$. Quotienting by this action, we obtain an $SO(3) \cong PU(2)$-equivariant map \[\text{Hol}: \widetilde B^e_{E,k} \times \Omega Y \to U(2),\] with $PU(2)$ acting on $U(2)$ by conjugation. We also have, associated to $A_0$, a map 
$$\text{Hol}_{A_0}: \Omega_b Y \to U(1).$$ 
Because $\text{tr}(A) = A_0$, we have 
$$\det \text{Hol}_A = \text{Hol}_{A_0}.$$ 
To define a relative holonomy map valued in $SU(2)$ and equivariant under the action of $SO(3) \cong PSU(2)$, we use a square root of $\text{Hol}_{A_0}$, defined on a covering space 
\[\begin{tikzcd}
	{\Omega'_b Y} & {\Omega_b Y} \\
	{U(1)} & {U(1)}
	\arrow["p", from=1-1, to=1-2]
	\arrow["{\sqrt{\text{Hol}_{A_0}}}"', from=1-1, to=2-1]
	\arrow["{\text{Hol}_{A_0}}", from=1-2, to=2-2]
	\arrow["{z^2}"', from=2-1, to=2-2]
\end{tikzcd}\]
Including $U(1)$ as the diagonal subgroup of $U(2)$, we may define an $SO(3) \cong PU(2)$-equivariant relative holonomy map 
$$\text{Hol}' := \text{Hol} \sqrt{\text{Hol}_{A_0}}^{-1}: \widetilde{\mathcal B}_Y \times \Omega'_b Y \to SU(2).$$ 
If one further acts on connections by the operation of tensoring with $\chi_{\mathbb C}$, where $\chi$ is a real line bundle, this describes the action of 
$$\mathcal G_E/\mathcal G^e_E = H^1(Y;\mathbb Z/2)$$ 
on this space. Abusing notation to write $\chi: \pi_1 Y \to \pm 1$ for the holonomy of the unique flat connection on this real line bundle, the corresponding map $\chi_{\mathbb C}: \pi_1 Y \to U(1)$ has determinant $1$; so 
$$\text{Hol}'_{A \otimes \chi}(\gamma) = \chi(\gamma) \text{Hol}'_A.$$
That is, $\text{Hol}'$ \emph{is not invariant under the $H^1(Y;\mathbb Z/2)$ action}.\\

To define cylinder functions using this relative holonomy map, the collection of immersions $q_i$ define a map $q: D^2 \to (\Omega_b Y)^N$. Because the disc is contractible, one may choose a lift $q': D^2 \to (\Omega'_b Y)^N$. Choosing an $SO(3)$-invariant function $h: SU(2)^N \to \mathbb R$, we set \[f_{q,h,\mu}(A) = \int_{D^2} h(\text{Hol}'_{q'}(A)) \mu.\] 

We call these $SU(2)$-cylinder functions. Note that $SO(3)$-cylinder functions are the special case of $H^1(Y;\mathbb Z/2)$-invariant $SU(2)$-cylinder functions.\\

We're interested in the gradient flow equation of the Chern-Simons functional on $\widetilde{\mathcal B}$; the relevant perturbations (called ``holonomy perturbations") will be the formal gradients of the cylinder functions described above, so we will need to know that the formal gradients are well-behaved. The most convenient way to do so is to calculate them explicitly. We write 
$$\nabla h: SO(3)^N \to \mathfrak{so}(3)^N$$ 
for the gradient of $h$ using the Lie group trivialization of the tangent bundle and the standard inner product on $\mathfrak{so}(3)$; write 
$$\nabla^j h: SO(3)^N \to \mathfrak{so}(3)$$ 
for the $j$th component. Kronheimer and Mrowka give a formula for this in \cite[Equation~(74)]{KM1} along the base $[-\varepsilon, \varepsilon] \times D^2 \hookrightarrow Y$ of the embeddings $q_i$ as $$\nabla_f := (\nabla f_{q,h,\mu})(A) = \ast \left(\sum_{i=1}^N (q_i)_* ((\nabla^i h)(\text{Hol}_{\mathbf{q}}(A)) \mu)\right),$$ where $(q_i)_*$ is the pushforward of differential forms on each tangent space. This is defined on the rest of $\text{Im}(\mathbf{q})$ via parallel transport, and then extended by 0 to the rest of $Y$. The following properties of the formal gradient $\nabla_f$ are enumerated as in \cite[Proposition~3.5]{KM1}.

\begin{theorem}\label{pert}Let $f: \mathcal A_{E,k}\to \mathbb R$ be a fixed $SU(2)$-cylinder function as above. Then the formal gradient $\nabla_f$ defines a smooth section of the tangent bundle $\mathcal T_k \mathcal A_{E,k}$, and its first derivative $D\nabla_f$, considered as a section of the bundle $\textup{Hom}(\mathcal T_k, \mathcal T_k)$, extends to $\textup{Hom}(\mathcal T_j, \mathcal T_j)$ for all $j \leq k$.\\
Further, fixing a reference connection $A_0$, we have the following pointwise norm bounds for constants $C, C_k, C'$:
\begin{align*}\|\nabla_f(A)\|_{L^\infty} &\leq C \\
\|\nabla_f(A)\|_{L^2_k} &\leq C_k(1+\|A-A_0\|_{L^2_k})^k\\
\forall \infty \geq p\geq 2, \|\nabla_f(A) - \nabla_f(A')\|_{L^p} &\leq C' \|A - A'\|_{L^p}.
\end{align*}
More generally, for any $n$ we may find continuous increasing functions $k_n$ (depending on $f$) such that $$\|(D^n \nabla_f)(A)\|_{L^2_k} \leq k_n(\|A - A_0\|_{L^2_k}),$$ where on the left side we're taking an operator norm. Furthermore, $f$ (and thus $\nabla_f$) are invariant under the action of the full gauge group.
\end{theorem}

Suppose now we have a collection of $SU(2)$-cylinder functions $f_i$; we may associate an increasing sequence of constants $(K_i)$ so that when $\sum_i K_i |a_i|$ converges, $\sum a_i \nabla_{f_i}$ converges in $C^k$ on compact sets. One may take $K_i$ to be, for instance, 
$$K_{i-1} + \sum_{n=1}^i k_n(f_i)(i).$$ Then we may define the $K$-weighted $\ell^1$ Banach space as the space of sequences $(a_i)$ with $\|a\|_{\mathcal P} = \sum K_i |a_i| < \infty$. There is a map $\mathcal P \to C^\infty(\widetilde{\mathcal B}^e_{E,k}, \mathbb R)$ given by sending 
$$\pi = \sum_i a_i \mapsto \sum_i a_i f_i =: f_\pi$$ 
and a map to smooth sections of the tangent bundle by sending $\pi$ to its formal gradient, which we still call $\nabla_\pi$. We call $\mathcal P$ a Banach space of perturbations and the induced functions $f_\pi$ (which are now $L^1$ sums of $SU(2)$ cylinder functions, but not necessarily cylinder functions themselves) \textit{holonomy perturbations}. 

The results of the previous theorem extend to smoothness of the $\mathcal P$-parameterized section of the tangent bundle, and parameterized inequalities for $\nabla_\pi(A)$ where now $\pi$ may vary within $\mathcal P$ (with an extra factor of $\|\pi\|_{\mathcal P}$ in the right-hand sides); this is \cite[Proposition~3.7]{KM1}. 

The following is the relevant geometric fact about $SO(3)$ and $SU(2)$-cylinder functions --- and therefore also holonomy perturbations. It is well-known, though usually stated nonequivariantly and for irreducible connections. Note that in the first statement we have quotiented by the entire (as opposed to even) based gauge group.

\begin{lemma}\label{dense}Given any compact $SO(3)$-invariant submanifold $S \subset \widetilde{\mathcal B}_{E,k}$, the restriction of $SO(3)$-cylinder functions to $S$ are dense in the space of smooth invariant functions $C^\infty(S)^{SO(3)}$. This is furthermore true even if we demand the cylinder functions vanish in a small invariant open neighborhood $\mathcal O$ of a finite set of orbits.

Similarly, given a compact $SO(3)$-invariant submanifold $S \subset \widetilde{\mathcal B}^{\textup{det}}_{\tilde E} \cong \widetilde{\mathcal B}^e_E$, the restriction of $SU(2)$-cylinder functions to $S$ are dense in the space of smooth invariant functions on $S$.
\end{lemma}

\begin{proof}
If $y \in Y$ is the basepoint, consider the maps $\mathcal A_{E,k} \to SO(3)$ given by $\text{Hol}(\gamma_i)$ where $\gamma_i$ are some family of immersed loops, based at $y$, all with the same germ there. (`Immersed' is implicit for all loops in the rest of this argument.) These maps are equivariant with respect to the action of the gauge group on $SO(3)$ (conjugating by the value at the basepoint), and hence descend to $SO(3)$-equivariant maps $\widetilde{\mathcal B}_{E,k} \to SO(3)$. The claim is that these maps separate points and tangent vectors in $\widetilde{\mathcal B}_{E,k}$; then, by compactness of $S$ there is some finite collection $\gamma_i$, $1 \leq i \leq N$, such that $\text{Hol}_A(\gamma_i): S \to SO(3)^N$ is an embedding.

If two connections have the same holonomy along every immersed loop $\gamma$ based at $y$ with specified germ, they have the same holonomy along every piecewise smooth immersed loop: we may write any piecewise smooth loop $\gamma$ based at $y$ as being the composition of an arbitrarily small loop and a loop with the specified germ, and so $\text{Hol}_A \gamma = \text{Hol}_{A'} \gamma$ for any piecewise smooth based loop $\gamma$. For any point $p \in Y$, pick a path $\alpha$ from $y$ to $p$, and define 
$$\sigma(p) = \text{Hol}_{A}(\alpha) \text{Hol}_{A'}^{-1}(\alpha);$$ 
that this is well-defined follows from the assumption that the (based) holonomies always agree along closed loops. Taking derivatives along any smooth curve based at $y$ we see that $A - d_A \sigma = A'$, and clearly $\sigma(y) = 1$. So based holonomy separates connections modulo based gauge. Equivalently, holonomy separates framed connections modulo gauge.

At the level of tangent spaces, we run the same argument with any $\omega$ such that $d\text{Hol}_A (\gamma)(\omega) = 0$ for all $\gamma$ to see that $\omega = d_A \sigma$ for some $\sigma \in \Omega^0(\mathfrak g_E)$. Pick a path $\alpha$ from $y$ to an arbitrary point $p$; then we have a natural isomorphism 
$$\mathfrak g_p \cong T\text{Isom}(E_y, E_p)$$
coming from the framing on $E_y$, and so we may define 
$$\sigma(p) = d\text{Hol}_A(\alpha)(\omega);$$ 
as before this is independent of the choice of $\alpha$ and it is easy to verify from this definition that $d_A \sigma = \omega$. Thus based holonomy seperates tangent vectors in $\widetilde{\mathcal B}_{E,k}$.

Now choose some finite set of embedded based curves $\gamma_i$ with the same germ so that $\text{Hol}_A(\gamma)$ embeds $S$ into $SO(3)^N$. Given $h: S \to \mathbb R$ an $SO(3)$-invariant function, extend it to $\tilde h: SO(3)^N \to \mathbb R$ arbitrarily using the equivariant tubular neighborhood theorem. We may then approximate $h$ by holonomy perturbations by picking a sequence of embeddings $q_i^j$ of solid tori of small radius around the curves so that the radii go to zero as $j \to \infty$; then the cylinder functions 
$$f_{q^j,\tilde h, \mu}: S \to \mathbb R$$ 
approach $h$ in $C^\infty(S)$ for any choice of 2-form $\mu$, as desired.

The final claim about $SO(3)$-cylinder functions follows by applying the previous paragraphs to embed the submanifold $S$ and the disjoint finite set of orbits $C$ into $SO(3)^N$; if $\mathcal O$ is, for instance, the inverse image of a small neighborhood of the image of $C$ in $SO(3)^N$ (whose closure is disjoint from the image of $S$). Then choose the extended function $\tilde h$ to be zero on this small neighborhood.

The corresponding fact for $SU(2)$-cylinder functions is a straightforward modifications of the same proof.
\end{proof}

Following the lemma, we choose our Banach space of perturbations $\mathcal P_E$ to be generated by a choice of a countable set of holonomy perturbations $\{\pi_i\}$, generated by $(h, \gamma, \mu)$. Here 
\begin{itemize}
\item $h$ varies over a countable dense set of $SO(3)$-invariant functions on $SU(2)^N$ which vanish in a neighborhood of $(\pm 1)^N$,
\item $\gamma$ varies over a countable dense subset of the space of $N$ immersions of $S^1 \times D^2$ that share a germ around $\{1\} \times D^2$, 
\item $\mu$ is specified, and 
\item $N$ varies over all positive integers. 
\end{itemize} 

We say that a Banach space of perturbations $\mathcal P_E$ arising in this way is \textit{sufficiently large}. The first condition will be useful in technical arguments later, so that the Hessian at fully reducible connections is the same no matter the perturbation.

The particular countable family of holonomy perturbations is inessential (other than these denseness conditions). While the particular choice of countable dense set to choose is noncanonical, the union of any two gives another; any theorem on independence of perturbation in $\mathcal P_E$ extends to show that the particular choice of $\mathcal P_E$ is irrelevant. 

In practice we will need to make further perturbations, but we want to do so without changing the existing functional at a finite set $C$ of `acceptable' critical orbits. We say that holonomy perturbations vanishing on a small open neighborhood $\mathcal O$ of this set $C$ \textit{adapted} to $\mathcal O$; the lemma shows these are in large supply. If $\mathcal P_E$ is a Banach space of perturbations and $\pi_0$ is some fixed perturbation, we denote $$\mathcal P_{E,\mathcal O} := \{\pi \in \mathcal P_E \mid f_\pi\big|_{\mathcal O} = f_{\pi_0}\big|_{\mathcal O}\}.$$

It is well-known that the space of flat connections is compact. The set of critical points of the unperturbed Chern-Simons functional framed flat connections, modulo even gauge; it's equivalently described as a sort of framed projective representation variety of $E$. We recall the following compactness principle for critical points from \cite[Lemma~3.8]{KM1}; their argument is unchanged for the framed moduli space, modulo even gauge.

\begin{lemma}\label{prop1}Let $\mathcal P$ be a Banach space of perturbations. The map $$F: \mathcal P \times \widetilde{\mathcal B}^e_{E,k} \to \mathcal P \times \mathcal T_{k-1} \widetilde{\mathcal B}^e_{E,k},$$  $(\pi, [A]) \mapsto (\pi, \nabla(\textup{cs}+f_\pi)([A]))$, is proper. In particular, if $\mathfrak C_\bullet \subset \mathcal P \times \widetilde{\mathcal B}^e_{E,k}$ denotes the subset whose fiber over $\pi \in \mathcal P$ is the set of critical points of $\textup{cs} + f_\pi$ --- that is, $\mathfrak C_\bullet = F^{-1}(\mathcal P \times \{0\})$ --- then the projection $\mathfrak C_\bullet \to \mathcal P$ is proper.
\end{lemma}

\section{Linear analysis}\label{sec:3d-linear}
It follows immediately from \cite[Equation~2.18]{Don}, modified to fit our normalization of $\text{cs}$, that $$\text{cs}(A+\omega) = \text{cs}(A) + \int_Y \text{tr}(d_A \omega \wedge \omega + \frac 23 \omega \wedge \omega \wedge \omega)$$ that the Hessian of the Chern-Simons functional on $\mathcal A_E$ is precisely 
$$*d_A: \Omega^1_k(\mathfrak g_E) \to \Omega^1_{k-1}(\mathfrak g_E).$$ 
(When working with Sobolev completions, the Hessian is a smooth map between tangent bundles of different regularity.) Because the tangent space at a framed connection $[A,p] \in \widetilde{\mathcal B}^e_{E,k}$ can be decomposed as $\mathcal K_{A,k} \oplus \mathfrak g_A^\perp$ (recall that $\mathcal K_{A,k}$ is defined to be the Coulomb slice $\text{ker}(d_A^*) \subset \Omega^1_k(\mathfrak g_E)$), this identifies the Hessian of the Chern-Simons functional \emph{at a flat connection} as 
$$\text{Hess}_A = *d_A \oplus 0: \mathcal K_{A,k} \oplus \mathfrak g_A^\perp \to \mathcal K_{A,k-1} \oplus \mathfrak g_A^\perp.$$ 

We call the summand $*d_A: \mathcal K_{A,k} \to \mathcal K_{A,k-1}$ the normal Hessian, $\text{Hess}^\nu_{A}$, which is the Hessian restricted to the normal space to the orbit. The Hessian is a summand of a larger elliptic operator, the extended Hessian 
$$\widehat{\text{Hess}}_A: \Omega^0_k(\mathfrak g_E) \oplus \Omega^1_k(\mathfrak g_E) \to \Omega^0_{k-1}(\mathfrak g_E) \oplus \Omega^1_{k-1}(\mathfrak g_E)$$ 
by the matrix 
\[\begin{pmatrix}0 & -d_A^* \\ - d_A & * d_A\end{pmatrix}.\] 
It is clear that $\widehat{\text{Hess}}_A^2 = \Delta_A$ on $\Omega^0 \oplus \Omega^1$ and thus $\widehat{\text{Hess}}_A$ is a self-adjoint elliptic operator. At a flat connection, where $(*d_A)^2 = 0$, the extended Hessian has the further decomposition under 
$$\Omega^1(\mathfrak g_E) = \mathcal K_{A,k} \oplus \mathcal G_{A,k}$$ 
as the orthogonal direct sum of $\text{Hess}^\nu_{A}$ on $\mathcal K_{A,k}$ and the signature operator $S_A$ on $\Omega^0 \oplus \mathcal G_{A,k}$. The kernel of the Hessian is identified with the kernel of the Laplacian on 0- and 1-forms (and is thus invertible if the Laplacian is).

Adding a perturbing term, the Hessian of the perturbed Chern-Simons functional $\text{Hess}_{A, \pi}$ differs only by the gradient of the vector field $\nabla_\pi$ at $[A]$, written $D_A \nabla_\pi$. We define the perturbed extended Hessian the same way. The perturbation of $\widehat{\text{Hess}}_{A}$ is a compact self-adjoint operator, and so the resulting perturbed extended Hessian operator is still self-adjoint first order elliptic. At critical points of $\text{cs}+f_\pi$, using the decomposition $\Omega^0 = \mathcal H^0_A \oplus \text{Im}(d_A^*)$, we see that 
$$\widehat{\text{Hess}}_{A,\pi} = 0 \oplus \text{Hess}_{A,\pi}^\nu \oplus S_A$$ 
on 
$$\mathcal H^0_A \oplus (\text{Im}(d_A^*) \oplus \mathcal G_{A,k}) \oplus \mathcal K_{A,k}.$$ 

Thus at critical points $A$ of the perturbed Chern-Simons functional, $\text{Hess}^\nu_{A,\pi}$ is a summand of a self-adjoint elliptic operator, and so it has discrete eigenvalues and the space it acts on enjoys a direct sum decomposition over these (finite-dimensional) eigenspaces.

\begin{definition}A critical orbit $[A,p] \in \widetilde{\mathcal B}^e_{E,k}$ of the perturbed function $\textup{cs} + f_\pi$ is called $\text{\em{nondegenerate}}$ if $\text{\em{Hess}}^\nu_{[A],\pi}$ is invertible, or, equivalently, the kernel of $\widehat{\text{\em{Hess}}}_{A,\pi}$ on $\Omega^0_k(\mathfrak g_E) \oplus \Omega^1_k(\mathfrak g_E)$ consists only of harmonic 0-forms.\end{definition}
This terminology agrees with \cite{Don}; if the connection $A$ is furthermore irreducible, and so its extended Hessian operator is invertible, Donaldson calls $A$ \textit{acyclic}.

\section{Reducible critical points}\label{sec:3d-red}
Because of the linear nature of the $SO(2)$-fixed subspace of $\widetilde{\mathcal B}^e_E$, $$\left(\widetilde{\mathcal B}^e_E\right)^{SO(2)} = \bigsqcup_\eta \mathcal B_{\eta},$$ it is especially straightforward to determine the critical points of the Chern-Simons functional. When $b_1(Y) = 0$, there is a unique flat connection in each component (corresponding to the calculation $\text{Hom}(\pi_1(Y), \mathbb Z) = 0$), at which the Hessian inside the $SO(2)$-fixed subspace is simply the restriction of $*d: \Omega^1_k(Y;i\mathbb R) \to \Omega^1_{k-1}(Y;i\mathbb R)$ to $$\text{Hess}_\theta = *d: \text{ker}(d^*)_k \to \text{ker}(d^*)_{k-1}.$$ This depends only on the underlying metric on $Y$, not the connection $A$. Because $b_1(Y) = 0$, Hodge theory guarantees that this is an isomorphism, and so these critical points are nondegenerate in the $SO(2)$-fixed locus. 

\begin{proposition}\label{red1}Suppose $Y$ is a Riemannian 3-manifold with $b_1(Y) = 0$, equipped with an $SO(3)$-bundle $E$. There is a positive constant $\varepsilon$ so that for any perturbation $\pi \in \mathcal P$ with $\|\pi\| \leq \varepsilon$, there is precisely one critical point of $\textup{cs}+f_\pi$ in each component $\mathcal B_{\eta,k}$ of the $SO(2)$-fixed point space, and the Hessian of each inside $\mathcal B_{\eta, k}$ is nondegenerate. 
\end{proposition}

\begin{proof}The compactness principle Lemma \ref{prop1} also applies to $\mathcal B_{\eta, k}$. The set of nondegenerate critical points in $\mathfrak C_\bullet$ (which we call $\mathfrak C_\bullet^*$) is open, and hence the subset of degenerate points is closed; because proper maps are closed, the set of perturbations $\pi$ for which some perturbed critical point $[A]$ is cut out \textit{non-transversely} is closed. Conversely, regular perturbations (written $\mathcal P^*$) form an open set. Thus, because every reducible is cut out transversely in the fixed locus for $\pi = 0$, the same is true for $\|\pi\| \leq \varepsilon$ for some small $\varepsilon$. Finally, if $p: \mathfrak C_\bullet \to \mathcal P$ is the projection, the map 
$$p: p^{-1}(\mathcal P^*) \to \mathcal P^*$$
is still a proper local diffeomorphism and hence the inverse function theorem guarantees that (for sufficiently small $\|\pi\|$) there is precisely one critical point of $\text{cs}+f_\pi$ in each component $\mathcal B_{\eta, k}$ of the $SO(2)$-fixed point set.
\end{proof}

There is a class of $SO(3)$-bundles over 3-manifolds with $b_1(Y) > 0$ for which we can avoid reducibles entirely, called admissible bundles (sometimes `non-trivial admissible' bundles). These were introduced in Floer's work on Dehn surgery in \cite{floer1995instanton}, and were extensively used in Kronheimer and Mrowka's study of instanton knot homology, beginning with \cite{KM1}. In defining instanton Floer homology we will restrict to these two somewhat orthogonal cases: $b_1 = 0$ or $E$ admissible.

\begin{definition}\label{admiss-3m}An $SO(3)$-bundle over a 3-manifold $Y$ is \textit{admissible} if every lift of $w_2(E)$ to a class in $H^2(Y;\mathbb Z)$ is non-torsion. We say that an $SO(3)$-bundle $E$ over $Y$ is \textit{weakly admissible} if either $E$ is admissible or $b_1(Y) = 0$.\end{definition}

Ultimately, the instanton chain complex is described in terms of critical \textit{orbits}, not fixed points, so we give a description of these. The class of weakly admissible bundles allows for a succinct description of the class of reducible critical orbits.

\begin{proposition}\label{redcrit}Let $Y$ be a Riemannian 3-manifold equipped with weakly admissible $SO(3)$-bundle $E$, equipped with a perturbation $\pi$ small enough that Proposition \ref{red1} applies. Then we have the following description of the reducible critical $SO(3)$-orbits of $\textup{cs}+f_\pi$ in $\widetilde{\mathcal B}^e_{E,k}$.
\begin{itemize}
    \item If $b_1(Y) = 0$, then there is a unique critical orbit lying in each reducible component $\{z_1, z_2\} \in \textup{Red}(Y,E)$ of the configuration space. In the components $\textup{Red}_*(Y,E)$ that contain a fully reducible orbit, the full reducible is the critical orbit.
    \item If $E$ is admissible, there are no reducible critical orbits.
\end{itemize}
\end{proposition}

\begin{proof}First suppose $b_1(Y) = 0$. By Proposition \ref{red1}, there is a unique critical point in each connected component of the $SO(2)$-fixed subspace. The components of the $SO(2)$-fixed subspaces were enumerated in Corollary \ref{orbits1}.

In general, a critical orbit of $\text{cs}$ in the component $\{z_1, z_2\} \in \text{Red}(Y,E)$ is a gauge equivalence class of flat connection respecting the isomorphism 
$$E \cong i\mathbb R \oplus \eta_1 \otimes \eta_2^{-1}.$$ 
In particular, $\eta_1 \otimes \eta_2^{-1}$ would support a flat connection, and thus have torsion first Chern class. Because $c = z_1 + z_2$, this implies that $c - 2z_2$ is a torsion cohomology class; as this also represents $w_2(E)$, this implies that $E$ is not admissible. Because $E$ supports no reducible critical points for the unperturbed Chern-Simons functional, the same is true of all sufficiently small perturbations.
\end{proof}

\section{Transversality for critical points}\label{sec:3d-trans}
Before discussing the 4-dimensional case, where we will assume that the limiting connections are nondegenerate, we should verify that this situation is achievable! This leads us to the first transversality theorem of this paper; its proof is a model for the rest of the transversality results we will need. It essentially follows the corresponding proof in \cite{AB}.

\begin{theorem}\label{trans1}Let $Y$ be a Riemannian 3-manifold, equipped with an $SO(3)$-bundle $E$. Suppose either that $b_1(Y) = 0$ or that $E$ is admissible. Given a sufficiently large Banach space of perturbations $\mathcal P$, there is a dense open set of $\pi \in \mathcal P$ for which there are finitely many critical orbits of the perturbed Chern-Simons functional $\widetilde{\mathcal B}^{e}_{E,k}$, on each of which the normal Hessian is invertible.
\end{theorem}

The proof proceeds by achieving transversality in each locus inductively. In the course of achieving transversality over the reducible locus, we need a simple functional analysis lemma.

\begin{lemma}\label{cpert}Let $H$ be a Hilbert space, equipped with a densely defined self-adjoint Fredholm operator $T: H \to H$; write $H_0 = \textup{ker}(T)$ and $H_1 = \textup{im}(T)$ for the associated orthogonal splitting of $H$. Let $K: H \to H$ be a bounded self-adjoint operator with $K_{00}: H_0 \to H_0$ injective. Then for sufficiently small positive $\varepsilon$, the map $T + \varepsilon K$ is injective; in particular, an isomorphism. 
\end{lemma}
\begin{proof}If $D = H_0 \oplus D_1$ is the domain of $T$, then because $T_{11}: D_1 \to H_1$ is an isomorphism, we have a bound $\|T_{11} x_1\| \geq C_T \|x_1\|$ for some $C_T > 0$ (where the norms are defined on the appropriate domains). Similarly we have $\|K_{00} x_0\| \geq C_K \|x_0\|$ for some $C_K > 0$. Now suppose $(T + \varepsilon K)(x_0 + x_1) = 0$ for some $x = (x_0, x_1) \in H$. This means precisely that 
\begin{align*}K_{00} x_0 + K_{10} &x_1 = 0,\\
T_{11} x_1 + \varepsilon (K_{01} x_0 &+ K_{11} x_1) = 0.
\end{align*}

These two equations give us the bounds 
\begin{align*}
C_K \|x_0\| \leq &\|K_{00} x_0\| = \|K_{10} x_1\| \leq \|K\| \|x_1\|,\\
C_T \|x_1\| \leq \|T_{11} x_1\| &= \varepsilon \|K_{01} x_0 + K_{11} x_1\| \leq \varepsilon\|K\|(\|x_0\| + \|x_1\|).
\end{align*}

Combining these two inequalities we then have $$C_T \|x_1\| \leq \varepsilon \|K\| (1 + \|K\|/C_K) \|x_1\|.$$ As soon as $$\varepsilon < \frac{C_T}{\|K\| (1 + \|K\|/C_K)},$$ this implies that $x_1 = 0$, and hence that $x_0 = 0$. 
\end{proof}

\begin{proof}[Proof of Theorem \ref{trans1}]By the discussion in the proof of Proposition \ref{red1}, there is an open subset $U$ of $\mathcal P$ so that all reducibles are cut out transversely in their own locus. For $\pi \in U$, there are thus a finite number of reducible $\pi$-critical orbits. Recall that $\text{cs}$ and $SU(2)$-holonomy perturbations are invariant under the action of $SO(3)$, but not of $H^1(Y;\mathbb Z/2)$. 

The normal Hessian to an $SO(2)$-reducible fixed point might have nontrivial kernel, a finite-dimensional subspace of the tangent space. 

Enumerate an element of each $SO(2)$-reducible critical orbit that is not fully reducible as $[A_i]$, and let $D_i$ denote a small disc inside the kernel of $\text{Hess}^\nu_{A_i} = Z_i$. The manifold $D_i$ is an $SO(2)$-manifold; we may identify the orbit through $D_i$ as a vector bundle over $S^2$ with generic fiber $D_i$; call this neighborhood $V_i$. An $SO(3)$-invariant function on $V_i$ is the same thing as an $SO(2)$-invariant function on $D_i$. At an $SO(2)$-reducible corresponding to a reduction $E \cong \mathbb R \oplus \zeta$, the kernel of this operator is $H^1(Y;\zeta)$. The $SO(2)$-action is by scalar multiplication on this complex vector space. The restriction of some $SO(2)$-invariant nondegenerate quadratic function $q$ on $Z_i$ gives a smooth function on $D_i$ with Hessian at the origin equal to the identity. Note that this implies that there are an even number of positive and negative eigenvalues. 

We will use holonomy perturbations approximating this smooth function on $D_i$ to correct the fact that the Hessian has kernel. More precisely, if $\gamma_j$ is a sequence of curves for which $\text{Hol}_A(\gamma)$ gives an equivariant embedding of the neighborhoods $V_i$ into $SU(2)^N$, we choose a function $\tilde h$ on $SU(2)^N$ which agrees in a neighborhood of each $A_i$ with $v \mapsto q(v)$ (identifying a neighborhood of each $A_i$ with the tangent space at that point). Now pick a sequence of holonomy perturbations $\pi_n$ arising from approximations of $\tilde h$ and $\gamma_j$. For sufficiently large $n$, the restriction of $D_{A_i}\nabla_{\pi_n}$ to $\text{ker}(\text{Hess}_{A_i})$ is very close to $q$.

If we decompose $\mathcal K_{A_i,k}$ as the direct sum of $Z_i$ and its orthogonal complement, then one may write $\text{Hess}_{A_i, \pi}$ in block-matrix form with respect to this splitting; the component $Z_i \to Z_i$ is the Hessian restricted to the submanifold $D_i$ at $A_i$. Therefore, we can apply Lemma \ref{cpert} to $T = \text{Hess}_{A_i,\pi}$ and $K = D_{A_i}\nabla \pi_n$ for large $n$ to see that $\pi+\varepsilon \pi_n$ is a regular perturbation at reducible connections for sufficiently large $n$ and sufficiently small $\varepsilon$. A similar argument works for the fully reducible connections, now using a nondegenerate $SO(3)$-invariant quadratic form.

Similarly to Proposition \ref{red1}, nondegeneracy for reducible critical points is an open condition\footnote{Note that \ref{red1} shows that nondegeneracy \textit{in the reducible locus} is an open condition, and here we want nondegeneracy in $\widetilde{\mathcal B}^e_E$; however, the proof requires only trivial modifications.} in perturbations $\mathcal P$, and we have just seen that it is also a dense condition. 

Now we have reduced ourselves to furthermore showing that perturbations for which the irreducible critical points are nondegenerate are dense in $\mathcal P_{\leq \epsilon}$. This is an application of the Sard-Smale theorem, using that $\widehat{\nabla}_\pi(\mathbf{A})$ separates points and tangent vectors as $\pi$ varies over $\mathcal P$. (For the irreducible case, see \cite{Fl1} and \cite[Section~5.5.1]{Don}.) The open-ness is slightly more delicate, and uses that nondegenerate critical points are isolated; therefore, for $\pi \in U$, reducible critical points are isolated. This implies that the projection $p^{-1}(U)^* \to U$ is proper, where the asterisk denotes that we restrict to irreducible connections. This is enough to conclude, as in Proposition \ref{red1}.
\end{proof}

\begin{remark}There is not much additional difficulty in proving this theorem for $SO(3)$-holonomy perturbations, which are invariant under the $H^1(Y;\mathbb Z/2)$ action. However, one needs to do more bookkeeping at the other kinds of reducibles; those with stabilizer $\mathbb Z/2$ are the most difficult, being labelled by their Euler classes in twisted cohomology. 
\end{remark}

\begin{remark}In the above proof, we only used the assumption that $E$ is weakly admissible --- that is, either $b_1(Y) = 0$ or $E$ is admissible --- to ensure that the reducible flat connections are cut out transversely in the reducible locus. This is not a serious issue: one may easily first choose a perturbation for which the reducible $\pi$-flat connections \textit{are} cut out transversely in the reducible locus. However, transversality for the moduli spaces of flowlines is much more restrictive, and forces the weak admissibility assumption on us --- if $b_1(Y) > 0$ and $E$ admits reducible connections, then there will be trajectories which have higher expected dimension in the reducible locus than the irreducible locus. That in mind, we see no harm in making the assumption a little early.
\end{remark}

Note that the critical set of $\text{cs}+f_\pi$ is canonically identified with that in any Hilbert manifold completion of lower regularity by an elliptic regularity argument; we are not being remiss in leaving the regularity $k$ out from our notation $\mathfrak C_\pi$ for the critical sets.

\section{Signature data of reducible critical points}\label{sec:3d-sigdata}
Suppose $Y$ is a rational homology sphere. Fix a lift $c$ of $w_2(E)$ to an integral cohomology class; then the connected components $\text{Red}(Y,E)$ of the configurations of reducible connections on $E$ are in bijection with pairs $\{z_1, z_2\} \subset H^2(Y;\mathbb Z)$ with $z_1 + z_2 = c$. (See Definition \ref{red-def} for the definition of this set and a discussion of its enumeration.) 

The component labeled by $\{z_1, z_2\}$ corresponds to equivalence classes of framed connections (with fixed determinant connection) on $\tilde E = \eta_1 \oplus \eta_2$ respecting the direct sum decomposition, where $c_1(\eta_i) = z_i$ and $\det \tilde E = \lambda$, with $c_1(\lambda) = c$. In 
$$E \cong \mathbb R \oplus (\eta_1 \otimes \eta_2^{-1})$$ 
these again correspond to connections respecting the direct sum decomposition. In particular, if $A$ represents the unique class of flat connection on $E$ that preserves the splitting, we may identify $\widehat{\text{Hess}}_A(\mathfrak g_E)$ with the direct sum 
$$\widehat{\text{Hess}}_\theta(\mathbb R) \oplus \widehat{\text{Hess}}_\rho(\eta_1 \otimes \eta_2^{-1}),$$ 
where $\rho$ is the flat connection given by restricting $A$ to $\eta_1 \otimes \eta_2^{-1}$. (and thus, the unique equivalence class of flat connection on this oriented 2-plane bundle). In particular, because the space of harmonic forms for a flat connection is isomorphic to the cohomology groups for the local system of coefficients defined by its holonomy representation, we have $$\text{ker } \widehat{\text{Hess}}^\nu_A \cong H^1(Y;\rho).$$ Our next goal is to obtain a more accessible calculation of this space.

\begin{lemma}Let $Y$ be a rational homology sphere and fix an abelian representation $\rho: \pi_1 Y \to U(1)$. If $\tilde Y$ is the universal abelian cover of $Y$ --- the closed 3-manifold arising as the covering space of $Y$ corresponding to $[\pi_1 Y, \pi_1 Y]$ --- then we may identify $H^1(Y; \rho)$ with the subspace of $H^1(\tilde Y; \mathbb C)$ on which $gx = \rho(g) x$ for all $g \in H_1(Y)$ (the `$\rho$-eigenspace'). In particular, if $H^1(\tilde Y; \mathbb C) = 0$, then $H^1(Y; \rho) = 0$ for all local coefficient systems $\rho$.
\end{lemma}

\begin{proof}Let $C_*(Y; \mathbb C)$ be the cellular chain complex of some finite CW-decomposition of $Y$, and $C_*(\tilde Y; \mathbb C)$ the cellular chain complex with $H_1(Y)$-action. The cohomology groups $H^*(Y; \rho)$ are (essentially by definition) isomorphic to the homology groups of the cochain complex 
$$\text{Hom}_{H_1 Y}\left(C_*(\tilde Y; \mathbb C), \mathbb C\right),$$
where $H_1(Y)$ acts on $\mathbb C$ via $\rho$. By definition, this cochain complex is isomorphic to the $\rho$-eigenspace of $C^*(\tilde Y; \mathbb C)$. Because the chain complex $C^*(\tilde Y; \mathbb C)$ splits as a direct sum of its eigenspaces, labelled by homomorphisms $H_1(Y) \to U(1)$, we see that the cohomology $H^*(\tilde Y; \mathbb C)$ splits as a direct sum of the direct sum of its eigenspaces $H^*_\rho(\tilde Y; \mathbb C)$, each of which is isomorphic to the homology of the eigenspace $C^*_\rho(\tilde Y; \mathbb C)$, which is then isomorphic to $H^*(Y; \rho)$.
\end{proof}

\begin{proposition}\label{Eigencounting}Let $Y$ be a rational homology sphere equipped with $SO(3)$-bundle $E$. For $\varepsilon > 0$ sufficiently small, the set of perturbations $\pi$ with nondegenerate reducible critical points and $\|\pi\| < \varepsilon$ is an open set, which decomposes as the disjoint union of nonempty open subspaces, each labelled by a function $$N_\pi: \textup{Red}_{SO(2)}(Y,E) \to \mathbb Z_{\geq 0}$$ with $N_\pi(\{z_1, z_2\}) \leq \textup{dim}_{\mathbb C } H^1\!\left(Y; \eta_1 \otimes \eta_2^{-1}\right)$.
\end{proposition}

\begin{proof}That the set of such perturbations $\pi$ is open is a mild extension of Proposition \ref{red1}.

Given any $\varepsilon$-small perturbation $\pi$, let $\beta: [0,1] \to [0,1]$ be a function equal to $1$ near $0$ and equal to $0$ near $1$, and consider the path $\pi(t) = \beta(t) \pi$ from $\pi$ to $0$; we have $\|\pi(t)\| < \varepsilon$ for all $t$. If 
$$r = \{z_1, z_2\} \in \text{Red}_{SO(2)}(Y,E)$$ 
denotes a connected component of the reducible subset of the configuration space (here $z_1 + z_2$ is a fixed integral lift of $w_2 E$ and $z_1 \neq z_2$), by Proposition \ref{red1}, there is a unique $\pi(t)$-critical orbit in the component labeled by $r$, and the proof shows it varies smoothly in $t$. Choose a particular path $A(t)$ of framed connections so that $[A(t)]$ lies in the $\pi(t)$-critical orbit in the reducible component labelled by $r$. We obtain an associated continuous path of self-adjoint Fredholm operators $$\widehat{\text{Hess}}_{A(t),\pi(t)}: \Omega^0_k \oplus \Omega^1_k \to \Omega^0_{k-1} \oplus \Omega^1_{k-1}.$$

There is an associated spectral flow, $\text{sf}\left(\widehat{\text{Hess}}_{A(t),\pi(t)}\right) \in \mathbb Z$, defined as the intersection number of the eigenvalues with the line $\lambda = - c$ for sufficiently small positive $c$; it is essentially the number of eigenvalues that go from negative to nonnegative, counted with sign. (This definition allows for the operators at $t = 0,1$ to have kernel, as opposed to intersecting with the line $\lambda = 0$.) In fact, $\widehat{\text{Hess}}_{A(t),\pi(t)}$ splits as a direct sum corresponding to the reductions $\mathbb R \oplus (\eta_1 \otimes \eta_2^{-1})$ (here $c_1(\eta_i) = z_i$); the operator in the first component is constant, and in the second component is complex linear, so this spectral flow is actually an even integer. The spectral flow of this path depends only on the ending perturbation $\pi$, so we define 
$$N_\pi(r) = \frac{1}{2}\text{sf}\left(\widehat{\text{Hess}}_{A(t),\pi(t)}\right).$$ 

The kernel of $\widehat{\text{Hess}}_{A(t),\pi(t)}$ splits as the direct sum of harmonic $0$-forms and a space of $1$-forms. Note that $A(t)$ is $SO(2)$-reducible for all $t$ (and never fully reducible), which implies that the dimension of the space of harmonic $0$-forms is constant in $t$, and so does not contribute to the spectral flow.

The spectral flow depends on $\pi$ continuously (that is, sufficiently close perturbations give the same integer) when $\widehat{\text{Hess}}_{A(0),\pi}$ has no kernel other than the harmonic 0-forms. To say that $\widehat{\text{Hess}}_{A(0),\pi}$ has no kernel other than harmonic 0-forms is precisely to say that $\pi$ is regular at the critical point $A(0)$. Therefore, the spectral flow to $0$ is a continuous function on the open set of interest, the space of regular perturbations.

As long as $\varepsilon$ is chosen small enough, the only eigenvalues that are close enough to $0$ to cross above or below it are the zero eigenvalues in $\text{ker}\;\widehat{\text{Hess}}_{A(1)}$, the unperturbed flat connection in this reducible component. The eigenvalues corresponding to harmonic $0$-forms stay constant, so we only need to think about the other zero eigenvalues.

As we are counting the number of eigenvalues that go from negative to nonnegative as we pass from $\widehat{\text{Hess}}_{A(0),\pi}$ to the unperturbed extended Hessian, the spectral flow must be nonnegative (said another way, the only relevant eigenvalues are the zero eigenvalues at $t = 1$, and those can only contribute to spectral flow by going down; reversing the direction, the only way to contribute to the spectral flow is to have an eigenvalue go upwards to $0$ as we approach $t \to 1$). Because there are 
$$\dim_{\mathbb R} H^1\left(Y; \eta_1 \otimes \eta_2^{-1}\right)$$ 
of these eigenvalues, this spectral flow is bounded above by that number.

Thus we have a continuous map from the space of reducible-regular perturbations to the discrete space $$\text{Map}(\text{Red}_{SO(2)}(Y,E), \mathbb Z_{\geq 0})$$ 
satisfying the stated bound; this decomposes the open set of reducible-regular perturbations as the disjoint union of open subsets labelled by these functions 
$$N_\pi: \text{Red}(Y, E) \to \mathbb Z_{\geq 0}$$ 
with upper bounds as in the statement of the proposition.

Given any function $N$ satisfying these bounds, we may construct a small perturbation $\pi$ so that $N$ arises as the spectral-flow function $N_\pi$. Choose an $S^1$-invariant quadratic function with $2N_\pi(r)$ positive eigenvalues on the space $\text{ker } \text{Hess}^\nu_A$, where $A$ is the unique critical point in the component labelled by $r$; note that on any given complex line, this function takes the form $cr^2$ where $r$ is the radius, for some constant $c$. 

Pick a holonomy perturbation that well approximates that quadratic function, thinking of $\text{ker } \text{Hess}^\nu_A$ as an $S^1$-submanifold of $\widetilde{\mathcal B}^e_{E,k}$ by exponentiating. Because holonomy perturbations are invariant under the full gauge group, the same quadratic function is introduced at the reducible labelled $x \cdot r$. Doing this as we vary over $\textup{Red}_{SO(2)}(Y,E)$ gives us a perturbation which restricts to a nondegenerate quadratic function, with the appropriate number of eigenvalues at each reducible critical point. Then if we choose $c$ sufficiently small, the perturbation $c \cdot \pi$ will be nondegenerate at each reducible and have $c \cdot \|\pi\| < \epsilon$.

Note that there is no analogous discussion in the components $\text{Red}_*(Y,E)$ containing a fully reducible point, because this point is the unique critical point in that component and is nondegenerate for all small $\pi$: the spectral flow is trivial.
\end{proof}

\begin{remark}It is likely the case that the open set labelled by $N_\pi: \textup{Red}_{SO(2)}(Y,E) \to \mathbb Z_{\geq 0}$ is \textit{connected}, but we will not find this necessary, so do not prove it here. In a sense similar to Remark \ref{infty-stuff}, the `space of regular perturbations with fixed signature data' is contractible. This will not be true in the 4-dimensional setting.
\end{remark}

\begin{remark}If one works with $SO(3)$-holonomy perturbations, which are invariant under $H^1(Y;\mathbb Z/2)$, one needs to do this bookkeping for additional types of reducibles, each with different demands on the signature data.
\end{remark}

The essential count turns out not to be the number of negative eigenvalues, but rather the \emph{signature} of the corresponding real vector space. As long as we remember the dimension of the underlying vector space, this is equivalent information.

\begin{definition}\label{sigdata}Let $(Y, E)$ be an $SO(3)$-bundle over a rational homology sphere, and fix a complex line bundle $\lambda$ with $w_2(E) = c_1(\lambda) \pmod 2$, which induces a bijection of $\textup{Red}_{SO(2)}(Y,E)$ with pairs $\{z_1, z_2 \}\subset H^2(Y; \mathbb Z)$ with $z_1 + z_2 = c_1(\lambda)$; the choice of $\lambda$ is essentially irrelevant. We write $\eta_1$ and $\eta_2$ for the corresponding complex line bundles.

We say that a \textit{signature datum} on $(Y, E)$ is a choice of function 
$$\sigma: \textup{Red}_{SO(2)}(Y,E) \to 2\mathbb Z$$ 
with 
$$|\sigma_\pi(\{z_1, z_2\})| \leq \textup{dim}_{\mathbb R} H^1\left(Y; \eta_1 \otimes \eta_2^{-1}\right)$$ 
and 
$$\sigma_\pi(\{z_1, z_2\}) \equiv \textup{dim}_{\mathbb R} H^1\left(Y; \eta_1 \otimes \eta_2^{-1}\right) \pmod 4.$$

Given a small regular perturbation $\pi$, the associated signature datum $\sigma_\pi$ is 
$$\sigma_\pi(\{z_1, z_2\}) = \textup{dim}_{\mathbb R} H^1\left(Y; \eta_1 \otimes \eta_2^{-1}\right) - 4N_\pi(\{z_1, z_2\}).$$

The set of \textit{signature data} is denoted $\sigma(Y, E)$.
\end{definition}

The name `signature datum' refers to the fact that $\sigma$ essentially chooses the symmetric bilinear form on each complex vector space $\text{ker } \text{Hess}^\nu_A$ that we perturb in the direction of (which are determined by their signature). The constraints are precisely those that the signature of an $S^1$-invariant bilinear form on a complex vector space must satisfy.

\begin{remark}If the universal abelian cover $\tilde Y$ of $Y$ (corresponding to the subgroup $[\pi_1 Y, \pi_1 Y] \subset \pi_1 Y$) is a rational homology sphere, then $|\sigma(Y, E)| = 1$ for all $E$. If the cover $\tilde Y'$ corresponding to the subgroup $[\pi_1 Y, \pi_1 Y]^2$ is a rational homology sphere, then $|\sigma(Y, \text{triv})| = 1$. 
\end{remark}

\chapter{Moduli spaces of instantons}\label{chap:4}
\section{Moduli spaces for cylinders and cobordisms}\label{sec:4d-moduli}
In Definition \ref{conf}, we defined the configuration space of framed connections on a bundle $\pE$ over a cylinder with specified limits as a quotient $$\widetilde{\mathcal B}^e_{k,\delta}(\alpha_-, \alpha_+) := \widetilde{\mathcal A}_{k,\delta}(A_-, A_+)\big/ \mathcal G_{k+1,\delta}^{e, h}(A_-, A_+).$$ The trivial vector bundle over this configuration space with fiber $\Omega^{2,+}_{k-1,\delta}(\text{End}(\pE))$ carries a compatible linear action of $\mathcal G_{\pE,k+1,\delta}^{e,h}(A_-, A_+)$, given by the action of the gauge group on $\pi^* \! \mathfrak g_E$; it therefore descends to a vector bundle $\mathcal S_{k-1,\delta}$ over the quotient ($\mathcal S$ for `self-dual') with fibers isomorphic to $\Omega^{2,+}_{k-1,\delta}$. The smoothly varying section $\mathbf{A} \mapsto F_{\mathbf{A}}^+$, which is equivariant under the gauge group action on the left and invariant under the $SO(3)$ action on the right, descends to an $SO(3)$-invariant smooth section of $\mathcal S_{k-1,\delta}$. 

\begin{definition}\label{moduli}The (unperturbed) (even) \emph{moduli space of framed instantons} on the bundle $\pE \to \mathbb R \times Y$ with limits $\alpha_-, \alpha_+$ is the $SO(3)$-invariant subspace 
$$\widetilde{\mathcal M}_E(\alpha_-, \alpha_+) \subset \widetilde{\mathcal B}^e_{k,\delta}(\alpha_-, \alpha_+),$$ 
consisting of gauge equivalence classes of framed connections $(\textbf{A},p) \in \widetilde{\mathcal A}_{k,\delta}(A_-, A_+)$ with $F_{\textbf{A}}^+ = 0,$ where $A_\pm$ are connections whose gauge equivalence class is $\alpha_\pm$. Equivalently, it is the zero set of the section $F_{\mathbf{A}}^+$ of $\mathcal S_{k-1,\delta}$. There are equivariant endpoint maps $\textup{ev}_{\pm}: \widetilde{\mathcal M}_E(\alpha_-, \alpha_+) \to \alpha_\pm$. We denote the quotient as $$\mathcal M_E(\alpha_-, \alpha_+) := \widetilde{\mathcal M}_E(\alpha_-, \alpha_+)\big/SO(3).$$
\end{definition}

In our notation for the moduli space, we drop the symbols $k,\delta,$ and $e$. The first two are unnecessary, because as long as the perturbation $\pi$ is of regularity $L^2_j$, the moduli spaces for all $1 < k \leq j$ are canonically identified by an elliptic regularity result. The requirement $k > 1$ is necessary to define the moduli spaces, as we require gauge transformations in $L^2_{k+1}$ to have continuous point-evaluation maps, which requires $k+1 > 2$. Similarly $\delta > 0$ is irrelevant so long as it is taken less than the absolute value of any nonzero eigenvalue of any of the eigenvalues of the limiting operator $\widehat{\text{Hess}}_{\alpha}$. The even moduli space is not the same as its quotient by the full gauge group, but we will never work with non-even moduli spaces, so we drop the superscript anyway.

Note that if the $A_\pm$ are not flat connections, there are no unperturbed instantons with those limits and $L^2$ curvature. 

If $\pi$ is a regular perturbation, so that all critical points of $\text{cs}+f_\pi$ are nondegenerate, we instead work with perturbed moduli spaces of instantons; we will need to do this to guarantee these moduli spaces are smooth manifolds.

Write $\mathcal P_E$ for the Banach space of perturbations of Chapter \ref{sec:3d-pert}. Given $\pi \in \mathcal P_E$ corresponding to $f_\pi: \widetilde{\mathcal B}^e_E \to \mathbb R$, the natural perturbation on the connections $\mathbf{A}$ over $\mathbb R \times Y$ is $$(\widehat \nabla_\pi)(\mathbf{A}) = \left(dt \wedge (\nabla f_\pi)(i_t^*\mathbf{A})\right)^+,$$ where $i_t^*$ is pullback to the slice $\{t\} \times Y$.

The operator $\widehat \nabla_\pi$ gives a well-defined, smoothly varying section of the trivial bundle $$\widetilde{\mathcal A}_{k,\delta}(A_-, A_+) \times \Omega^{2,+}_{k-1,\delta}(\text{End}(\pE)) \to \widetilde{\mathcal A}_{k,\delta}(A_-, A_+),$$ 
equivariant under the action of $\mathcal G^{e,h}_{k+1,\delta}(A_-, A_+)$, and the rest of Theorem \ref{pert} carries over with the obvious modifications (including a version with dependence on $\pi$). 

In particular, because the above perturbations are invariant under the action of the full 4-dimensional gauge group $\mathcal G^h_{k+1,\delta}(A_-, A_+)$, as well as the $SO(3)$-action, $F_{\mathbf{A}}^+ + \widehat\nabla_\pi(\mathbf{A})$ defines a smooth section 
\begin{equation}\label{ASDsec}\mathcal P \times \widetilde{\mathcal B}^e_{k,\delta}(\alpha_-, \alpha_+) \to \mathcal P \times \mathcal S_{k-1,\delta};\end{equation}
the fact that $F_{\mathbf{A}}^+ + \widehat \nabla_\pi(\mathbf{A}) \in \Omega^{2,+}_{k-1,\delta}\left(\text{End}(\pE)\right)$ follows from the assumption that $\mathbf{A}_t$ decays exponentially on the ends to connections with $*F_{A_\pm} = - \widehat \nabla_\pi(A_\pm)$. 

As in Definition \ref{moduli}, we say that the \textit{perturbed moduli space} of framed instantons for the perturbation $\pi$, written $\widetilde{\mathcal M}_{E,\pi}(\alpha_-, \alpha_+)$, is the set of gauge equivalence classes of framed connections with $F_{\mathbf{A}}^+ = - \widehat \nabla_\pi(\mathbf{A})$. As in the unperturbed case, $\mathbf{A}$ can only have $L^2$ curvature if $*F_{A_\pm} = - \nabla_\pi(A_{\pm})$, and so the limiting connections are critical points for $\text{cs}+f_\pi$. We call these \emph{$\pi$-flat connections}.

Using the gradient equation, we can define two useful notions of energy for an instanton on the cylinder.

\begin{proposition}\label{energy}For a perturbed connection $\mathbf{A}$ limiting to $\pi$-flat connections on the ends, we call the expression $$\int \|F_{\mathbf{A}} + *_3\nabla_\pi(\mathbf{A}(t))\|^2$$ the \emph{analytic energy} of $\mathbf{A}$. If $\mathbf{A}$ is a $\pi$-perturbed instanton, the analytic energy of $\mathbf{A}$ depends only on the endpoints $\alpha_-, \alpha_+$ of $\mathbf{A}$, and the homotopy class that $\mathbf{A}$ traces out, written $z \in \pi_1(\widetilde{\mathcal B}_E, \alpha_-, \alpha_+).$ If $A_+$ and $A_-$ are the limiting connections of $\mathbf{A}$, it is equal to $2\left(\textup{cs}_\pi(\alpha_-) - \textup{cs}_\pi(\alpha_+)\right).$ 

Correspondingly, given a pair of $\pi$-flat connections $\alpha_-, \alpha_+$ and a homotopy class $z$ between them, we call $$\mathcal E^\pi_z(\alpha_-, \alpha_+)=2\left(\textup{cs}_\pi(\alpha_-) - \textup{cs}_\pi(\alpha_+)\right)$$ the \textit{topological energy} of $(\alpha_-, \alpha_+, z)$. 

Whenever there is an instanton going from $\alpha_-$ to $\alpha_+$ in the homotopy class $z$, we have $\mathcal E^\pi_z(\alpha_-, \alpha_+) \geq 0$  with equality if and only if $\mathbf{A}$ is constant. If $z_i$ is a homotopy class from $\alpha_i$ to $\alpha_{i+1}$, then $$\mathcal E^\pi_{z_1 \ast z_2}(\alpha_1, \alpha_3) = \mathcal E^\pi_{z_1}(\alpha_1, \alpha_2) + \mathcal E^\pi_{z_2}(\alpha_2,\alpha_3).$$ If $1 \in \pi_1(\widetilde{\mathcal B}_E)$ is the positive generator, then $\mathcal E^\pi_{z+1}(\alpha_-, \alpha_+) = 64\pi^2 + \mathcal E^\pi_z(\alpha_-, \alpha_+)$. 
\end{proposition}
\begin{proof}The given expression is gauge invariant, so we may assume $\mathbf{A}$ is in temporal gauge and write $\mathbf{A} = d/dt + \mathbf{A}(t)$. The curvature is $F_{\mathbf{A}} = F_{\mathbf{A}(t)} + dt \wedge A'(t)$, and so we may write this integral as 
$$\int \|A'(t)\|^2 + \|F_{\mathbf{A}(t)} + \nabla_\pi\left(\mathbf{A}(t)\right)\|^2 dt.$$ 
Because 
$$F_{\mathbf{A}} + dt \wedge \nabla_\pi(\mathbf{A}(t)) = F_{\mathbf{A}} + *(*_3 \nabla_\pi(\mathbf{A}(t)))$$ 
is anti-self-dual, and $*(dt \wedge \omega) = *_3 \omega$ for a 1-form $\omega$ pulled back from $Y$, we see that 
$$A'(t) = - *_3 F_{\mathbf{A}(t)} - \nabla_\pi(\mathbf{A}(t)) = - \nabla(\text{cs} + f_\pi)(\mathbf{A}(t)).$$ 
Using that $\mathbf{A}(t)$ is a gradient flowline, we may rewrite the integral as 

\begin{align*}\int \|A'(t) + \nabla(\text{cs}&+f_\pi)(\mathbf{A}(t))\|^2 - 2\langle A'(t), \nabla(\text{cs}+f_\pi)(\mathbf{A}(t))\rangle dt \\
&= -2 \int \langle A'(t), \nabla(\text{cs}+f_\pi)(\mathbf{A}(t)) \rangle dt \\
&= -2 \int d(\text{cs}+f_{\pi})(A'(t))\\
&= -2\left((\text{cs}+f_\pi)(A_+ - (\text{cs}+f_\pi)(A_-)\right) \in \mathbb R,\end{align*}

the desired result. Now it is clear that this only depends on the connected component of the connection $\mathbf{A}$. (Here we are using the path to pin down the real lift of this difference of Chern-Simons values.) 

That $\mathcal E^\pi_0(\alpha, \alpha) = 0$ is clear from the existence of the constant trajectory. Additivity is clear by picking a connections $\mathbf{A}_i$ for each that are constant sufficiently far down the cylindrical ends, and gluing those together; the difference in $\text{cs}+f_\pi$ is additive.

Lastly, we need to determine $\mathcal E^\pi_1(\alpha, \alpha)$; by definition, this is the same as determining 
$$2(\text{cs}+f_\pi)(g(A)) - 2(\text{cs}+f_\pi)(A)$$ 
for a gauge transformation $g$ generating $\pi_0 \mathcal G^e_E$. Because $f_\pi$ is invariant under the full gauge group, we are only asking to determine the difference in the Chern-Simons functional. Pick a trajectory $\mathbf{A}$ going from $A$ to $g(A)$, constant near the ends. 
By definition, this is given by computing $2\int_{I \times Y} \text{Tr}(F_{\mathbf{A}}^2)$. Because $A$ differs from $g(A)$ by a gauge transformation, we may glue these together to get a connection on an $SO(3)$-bundle $\mathbf{E}$ over $S^1 \times Y$, the bundle $\mathbf{E}$ being the mapping torus of the automorphism $g$; noting that $F_{\mathbf{A}}^2$ is zero near the ends, this implies we may compute the curvature integral just as well over this closed 4-manifold, where it is equal to $-2 \cdot 8\pi^2 p_1(\mathbf{E})$. We choose the positive generator of $\pi_1(\widetilde{\mathcal B}^e_E)$ to be the one that makes this integral positive; it remains to compute it. This is computed more generally in \cite[Equation~(24)]{KM1}, and in this case the minimal $p_1(\mathbf{E})$ is $4$, as desired. Because of the slight difference that we use specifically even gauge transformations, we provide a short proof, given in the lemma that follows. 
\end{proof}

\begin{lemma}Let $E$ be an $SO(3)$-bundle over a 3-manifold $Y$. Given an even gauge transformation $g \in \Gamma(\textup{Aut}(E))$, the minimal first Pontryagin class of the mapping torus of $g$ on $S^1 \times Y$ is $4$.\end{lemma}
\begin{proof}
The bundle $\mathbf{E}$ has 
$$w_2 \mathbf{E} \in H^2(S^1 \times Y;\mathbb Z/2) \cong H^2(Y;\mathbb Z/2) \oplus H^1(Y;\mathbb Z/2).$$ 
The first factor of $w_2 \mathbf{E}$ is precisely the second Stiefel-Whitney class of its restriction to $Y$, which is $w_2 E$. The second term of $w_2 \mathbf{E}$ is precisely the obstruction class $o(g) \in H^1(Y;\mathbb Z/2)$ for lifting a section of $\text{Aut}(E)$ to a section of $\widetilde{\text{Aut}}(E)$. 

To see this, pick a loop $\gamma$ in $Y$ and consider the corresponding torus in $S^1 \times Y$. The map $\gamma \mapsto w_2(\mathbf{E})\big|_{S^1 \times \gamma}$ defines a homomorphism $H_1(Y;\mathbb Z) \to \mathbb Z/2$. If $o(g) = 0$, then the restriction of $g$ to $\gamma$ lifts to $\widetilde{\text{Aut}}(E)$, and so is homotopic to the identity (because $\pi_1 \widetilde{\text{Aut}}(E_x) = 0$), so the corresponding bundle on the torus $S^1 \times \gamma$ is trivial. 

Conversely, every bundle over $S^1 \times Y$ may be constructed as the mapping torus of a bundle over $Y$, and both $w_2(E)$ and $o(g)$ of the bundle and mapping are determined by $w_2 \mathbf{E}$. Given fixed $w_2(\mathbf{E}) = w_2(E) \oplus o(g)$, the Dold-Whitney theorem identifies the possible values of $p_1(\mathbf{E})$ as those elements $p \in H^4(S^1 \times Y;\mathbb Z)$ with 
$$p \equiv (w_2 \mathbf{E})^2 \pmod 4,$$ 
where the mod $4$ indicates we use the \emph{Pontryagin square} on even-dimensional cohomology $H^2(X;\mathbb Z/2) \to H^4(X;\mathbb Z/4)$. Expanding this, we obtain 
$$p_1(\mathbf{E}) = [S^1] \cdot 2o(g) \cdot w_2(E) \in H^4(S^1 \times Y;\mathbb Z/4).$$
In particular, as long as $o = 0$, the only condition is $p_1(\mathbf{E}) \equiv 0 \pmod 4$. 
\end{proof}

\begin{corollary}\label{triv-component}Given a 3-manifold $Y$ equipped with $SO(3)$-bundle $E$ and a regular perturbation $\pi$, if $\alpha$ is a $\pi$-critical orbit, the constant connection is the unique $\pi$-perturbed instanton $\mathbf{A}$ in the component labelled by the trivial homotopy class $0 \in \pi_1(\widetilde{\mathcal B}^e_E,\alpha)$.
\end{corollary}
\begin{proof}
The constant connection at a $\pi$-flat connection is always a solution to the perturbed ASD equations. Because $\mathcal E^\pi_0(\alpha, \alpha) = 0$, and the only connections with analytic energy equal to zero are constant, so is any solution of the perturbed ASD equations in this component of the space of connections.
\end{proof}

For a Riemannian manifold $W$, with two cylindrical ends, orientably isometric to $(-\infty, 0] \times Y_1$ and $[0,\infty) \times Y_2$, written as $W: Y_1 \to Y_2$, we can no longer consider constant perturbations (pulled back from $Y$). Further, we need to be able to interpolate between perturbations, and have enough available to prove transversality results. As before, the precise formulation of the perturbations here is similar to that in \cite[Section~3.2]{KM1}. 

Fix once and for all a smooth function $\beta_0: [0, \infty) \to [0,1]$ with $\beta_0(x) = 0$ for $x \leq 1$ and $\beta_0(x) = 1$ for $x \geq 2$.

Given a perturbation $\pi \in \mathcal P_{E_2}$, we may define the associated 4-dimensional perturbation $\widehat{\nabla}_{\pi}$ on $[0,\infty) \times Y_2$ as $$(\widehat \nabla_{\pi})(\mathbf{A}) = \beta_0(t)\big(dt \wedge \nabla_{\pi}(\mathbf{A}(t))\big)^+,$$ and identically zero elsewhere, and similarly for a perturbation in $\mathcal P_{E_1}$. 

An \emph{end perturbation} on $W$ is labelled by is $$(\pi_-, \pi_+) \in \mathcal P_{E_1} \times \mathcal P_{E_2} =: \mathcal P^{(4)}_{\text{end}}.$$ The terms $\pi_\pm$ denote the constant perturbations on the ends, and the $\pi_i$ are perturbing terms on a compact piece of the end, as above. (While $\mathcal P^{(4)}_{\text{end}}$ depends on $W$ and on the choice of isometry between the ends and $[0,\infty) \times Y$, this dependence is implicit and should be clear from context.) 

These end perturbations are the same as the perturbations used in \cite{KM1}; they introduce further perturbations supported in $(1,2) \times Y$, but as those are insufficient for later purposes, we use an altogether different source for our holonomy perturbations on the interior. 

We will need one further type of perturbation, which is mainly useful to achieve transversality at flat connections and reducible connections. These were introduced in the setting of closed 4-manifolds in \cite{Kronheimer}, and are similar to the notion of holonomy perturbations along thickened loops in \cite{donaldson1987orientation} and \cite{Fr3}.

\begin{definition}\label{ThickLoop}Let $W$ be a compact Riemannian manifold with boundary, equipped with an $SO(3)$-bundle $\mathbf{E}$, and $W^\circ$ its interior. A \emph{collection of thickened loops} in $W$ consists of a choice of closed ball $B \subset W^\circ$, a finite collection $q_1, \cdots, q_n$ of smooth submersions $q_i: S^1 \times B \to W^\circ$ so that $q_i(1, b) = b$ and $q_i(-, b)$ is an immersion for all $1 \leq i \leq n$ and $b \in B$; we will write $\vec q$ for the entire collection. 

If $\gamma: S^1 \to X^\circ$ is a curve, we write $\textup{Hol}_\gamma(\mathbf{A}) \in \widetilde{\textup{Aut}}(\mathbf{E}_{\gamma(1)})$ for the parallel transport around the loop, thought of as based at the image of $1 \in S^1$; the lift to $\widetilde{\textup{Aut}}$ is made canonical by the choice that $\textup{Hol}_\gamma(\mathbf{A}) = 1$ if $\mathbf{A}$ is fully reducible above $\gamma$. This is preserved by \emph{even} based gauge transformations. Similarly we write $\textup{Hol}_{\vec \gamma}(\mathbf{A}) \in \widetilde{\textup{Aut}}(\mathbf{E}_{\vec \gamma(1)})^n$ for an $n$-tuple of curves with a common basepoint.

An \emph{interior holonomy perturbation} is induced by the data of a collection of thickened loops $q_1, \cdots, q_n$, as well as a choice of smooth, (even) gauge-equivariant bundle map 
$$r: \widetilde{\textup{Aut}}(\mathbf{E})^N\big|_B \to \mathfrak g_{\mathbf{E}}\big|_B$$ 
and a smooth self-dual 2-form $\omega$ on $W$ with support in the interior of $B$. 

The interior holonomy perturbation corresponding to $(\vec q, r, \omega)$ is the map 

\begin{align*}\widehat{\nabla}_\pi: \mathcal A_{\mathbf{E}} &\to \Omega^{2,+}(W; \mathfrak g_{\mathbf{E}})\\
\widehat{\nabla}_\pi(\mathbf{A})(p) &= r\left(\textup{Hol}_{\vec q_i(p, -)}(\mathbf{A})\right) \otimes \omega(p).\end{align*}
Note that $\omega(p) = 0$ when $p \not\in B$, making this formula well-defined.
\end{definition}

We will apply these perturbations on the submanifold $$W' = W \setminus \big([1,\infty) \times Y\big).$$

Being gauge-equivariant, these perturbations descend to sections of the bundle $\mathcal V^+ \to \widetilde{\mathcal B}$. We will shortly discuss extentions to the configuration spaces modeled on Sobolev spaces.

Soon, we will choose a countable dense set of the possible choices of data above.

These perturbations satisfy the following properties, analogous to those of \cite[Proposition~3.7]{KM1}.

\begin{proposition}\label{pertfacts}Suppose $\pi_\pm$ are 3-manifold perturbations on $(Y_\pm)$, with fixed critical orbits $\alpha_\pm$. 
Suppose $\pi$ is a 4-dimensional perturbation on the Riemannian manifold $W$ with cylindrical ends defined above, including both interior and end perturbations. The map $\widehat{\nabla}_\pi$ extends to smooth maps $$\widehat{\nabla}_\pi: \mathcal A_{\mathbf{E},k,\delta}(\alpha_-, \alpha_+) \to \Omega^{2,+}_{k,\delta}(W;\mathfrak g_E)$$ for any $k \geq 2$ and $\delta > 0$ sufficiently small. They satisfy the following properties. 
\begin{enumerate}
\item The map $(\pi, \mathbf A) \mapsto (D\widehat{\nabla}_\pi)(\mathbf A)$ from $\mathcal P \times \mathcal A_{k,\delta}$ to $\text{Hom}\left(\Omega^1_{k,\delta}(W), \Omega^{2,+}_{k,\delta}(W)\right)$ extends smoothly to a map with codomain $\text{Hom}(\Omega^1_{j,\delta'}, \Omega^{2,+}_{j,\delta'})$ for all $j \le k$ and $\delta' \le \delta$.
\item Fixing a base connection $\mathbf{A}_0$, for each $n \geq 0$ and $k \geq 2$, there is an increasing continuous function $c_n$ so that $$\|D^n\widehat \nabla_\pi(\mathbf{A})\|_{L^2_{k,\delta}} \leq c_n\left(\|\mathbf{A}-\mathbf{A}_0\|_{L^2_{k,\delta}}\right);$$ when $n > 0$ this should be interpreted as the operator norm.
\item $\|\widehat{\nabla}_\pi(\mathbf{A})\|_{L^\infty} \leq K$, independent of $\mathbf{A}$. 
\item There is a constant $C$ so that for all $1 \leq p \leq \infty$, and all $L^2_k$ connections $\mathbf{A}, \mathbf{A}'$, we have 
$$\|\widehat{\nabla}_\pi(\mathbf{A}) - \widehat{\nabla}_\pi(\mathbf{A})\|_{L^p} \leq C\|\mathbf{A} - \mathbf{A}_0\|_{L^p}.$$
\item Suppose $\sigma$ is a gauge transformation of $\mathbf{E}$, defined away from some finite set of points; and suppose $\mathbf{A}_i$ are globally defined connections with $\sigma^*\mathbf{A}_1 = \mathbf{A}_2$ where defined. Then $$\widehat{\nabla}_\pi((\mathbf{A}_1)) = \widehat{\nabla}_\pi(\mathbf{A}_2).$$
\item Suppose we have a finite set $\mathbf{x} \subset W$, a sequence of gauge transformations $\sigma_n$ defined on $W \setminus \mathbf{x}$ which are $L^2_{k+1,\delta}$ on the ends and locally $L^2_k$ on the complement of $\mathbf{x}$, and a sequence of connections $\mathbf{A}_n$ on $\mathbf{E}$, so that over any compact set $K$ for any $p \geq 2$, the connection $\sigma_n(\mathbf{A}_n)$ is bounded in $L^p_1$ over $K$. Then $(\sigma_n)_* \widehat{\nabla}_\pi(\mathbf{A}_n)$ has a subsequence which is Cauchy in $L^p(W)$. If there is an $L^p$ connection $\mathbf{A}$ on $\mathbf{E}$ so that $\sigma_n(\mathbf{A}_n) \to \mathbf{A}$ in $L^p$ on compact sets away from $\mathbf{x}$, then the limit of this Cauchy subsequence is $\widehat{\nabla}_\pi(\mathbf{A})$.
\end{enumerate}
\end{proposition}

The proof must be divided into two parts, depending on the contribution from the end perturbations and the interior perturbations. In the former case, everything except the last point is given in \cite[Proposition~3.7]{KM1}. These calculations for interior perturbations are given in \cite[Section 3]{Kronheimer}; the last point is argued in \cite[Lemma~3.4]{Kronheimer}, and that argument works with only slight modifications for end perturbations as well.

As with the 3-dimensional perturbations and the end perturbations, then, we see that for any countable set of perturbations $\pi_i$ there is an increasing sequence of constants $C_i > 0$ so that for any sequence $d_i \in \mathbb R$ with $$\sum |d_i| C_i < \infty,$$ the limit $\sum d_i \widehat{\nabla}_{\pi_i}$ has the same properties as above; one defines the Banach space of perturbations to be the $L^1$ space weighted by $C_i$ (that is, the space of such sequences $d_i$, with norm $\|(d_i)\| = \sum |d_i| C_i$). All inequalities in the above may then be made dependent on $\pi$, adding a factor of $\|\pi\|$ on the right-hand side. 

\begin{definition}\label{pertdef}In this section, and in what follows, we used the following notation for Banach spaces of perturbations. \begin{itemize}
\item On a 3-manifold $Y$ equipped with $SO(3)$-bundle $E$, the space $\mathcal P_E$ given by weighted $L^1$ sums of a countable dense set of data defining cylinder functions (with derivatives vanishing to all orders at the fully reducible connections), \\
\item On a 4-manifold $(W,\mathbf{E})$, with cylindrical ends $(Y_1, E_1)$ and $(Y_2, E_2)$, the space $$\mathcal P^{(4)}_{\textup{end}} = \mathcal P_{E_1} \oplus \mathcal P_{E_2},$$ where the spaces $\mathcal P_{E_i}$ contribute a perturbation of the form $$\beta_0(t) \left(dt \wedge \nabla_{\pi}(\mathbf{A}(t))\right)^+$$ on the end, $\beta_0$ being a fixed cutoff function supported in $(1 \infty)$ (with $W$ left implicit in the notation $\mathcal P^{(4)}$), \\
\item On the same 4-manifold $(W, \mathbf{E})$, the space $\mathcal P^{(4)}_{\textup{int}}$ is given by $L^1$ sums of a countable dense set of interior holonomy perturbations, and we write $$\mathcal P^{(4)} = \mathcal P^{(4)}_{\textup{end}} \oplus \mathcal P^{(4)}_{\textup{int}}.$$
\item The affine subspace $\mathcal P^{(4)}_c \subset \mathcal P^{(4)}$ of perturbations with $(\pi_-, \pi_+) \in \mathcal P_{E_1} \oplus \mathcal P_{E_2}$ fixed (but left implicit in notation); this subspace is affine over $\mathcal P_{\textup{int}}.$
\end{itemize}

Fix $\epsilon > 0$ so that for any reducible flat connection $\alpha$, the only eigenvalues of $\widehat{\textup{Hess}}^\nu_{\alpha,0}$ with absolute value at most $\epsilon$ are the elements of the kernel; one may fix a connected neighborhood $0 \in U_\epsilon \subset \mathcal P_E$ so that this remains true for all $\pi \in U_\epsilon$ and any reducible $\pi$-flat connection $\alpha$. 

This implies that $\dim \textup{Eig}_{\pm \epsilon}(\widehat{\textup{Hess}}^\nu_{\alpha,\pi}) = 0$ for all $\pi \in U_\epsilon$ and $\alpha$ a $\pi$-flat connection. In particular, this means that $$\dim \textup{Eig}_{-\epsilon \leq \lambda \leq \epsilon}(\widehat{\textup{Hess}}^\nu_{\alpha,\pi}) = \dim \ker \widehat{\textup{Hess}}^\nu_{\alpha,0}$$ for all such $(\pi, \alpha)$. 

If we let $\delta$ be a constant $0 < \delta \ll \epsilon$, then we denote by $\mathcal P_{E,\delta} \subset U_\epsilon$ the subset of those $\pi \in U_\epsilon$ so that $\widehat{\textup{Hess}}^\nu_{\alpha,\pi}$ has no eigenvalues of absolute value at most $\delta$, for any $\pi$-flat connection $\alpha$ (not necessarily reducible). This subset does not include $0$, but by Theorem \ref{trans1}, $\bigcup_{\delta \to 0} \mathcal P_{E, \delta}$ is dense in $U_\epsilon$. There are corresponding open subspaces $\mathcal P^{(4)}_\delta$ and $\mathcal P^{(4)}_{\textup{end}, \delta}$ for which the perturbations on the 3-manifold at $\infty$ lie in $\mathcal P_{E_i, \delta}$. These spaces will be important for the later weighted Sobolev theory. 
\end{definition}

So long as $\delta$ is sufficiently small, there is a union of connected components of $\mathcal P_{E,\delta}$ so that the enumeration of reducible flat connections of Proposition \ref{red1} holds. This follows because modifying the count of reducible flat connections would require one of them to become critical, and so the count is the same at every point in the component; further, there is some $\pi \in \mathcal P_{E,\delta}$ arbitrarily close to zero, so long as $\delta$ is chosen sufficiently small relative to the distance to zero. 

Additionally, we may decompose $\mathcal P_{E, \delta}$ into connected open sets depending on the spectral flow of $\widehat{\text{Hess}}^\nu_{\alpha_t, \pi_t}$ for a generic path $0 \to \pi \in \mathcal P_{E, \delta}$ at the reducible $\pi_t$-flat connections $\alpha$, precisely as in Proposition \ref{Eigencounting}; this will be relevant later when calculating certain operator indices.
We still have a notion of analytic energy of a configuration on $W$, but it is no longer a topological invariant, and will mostly be useful in keeping track of what happens during compactification.

\begin{definition}Let $W$ be a Riemannian manifold with cylindrical ends, equipped with an $SO(3)$-bundle and a perturbation $\pi = (\pi_-, \pi_1, \pi_2, \pi_+)$ as above. Let $\alpha_\pm$ be critical orbits with respect to $\pi_\pm$. The \emph{analytic energy} of a connection $\mathbf{A} \in \widetilde{\mathcal B}^{e}_{E,k,\delta}(\alpha_-,\alpha_+)$ is defined to be $$\mathcal E^\pi_{\textup{an}}(\mathbf{A}) = \int_W \|F^+_{\mathbf{A}} + \widehat{\nabla}_\pi(\mathbf{A})\|^2.$$
\end{definition}

This gives a real-valued lift of $\text{cs}_{\pi_+}(\alpha_+) - \text{cs}_{\pi_-}(\alpha_-)$ for $\pi$-ASD connections which only depends on their homotopy class $z$, and agrees with
$$\int \|F_{\mathbf{A}} + *_3\nabla_\pi(\mathbf{A}(t))\|^2$$
for connections on the cylinder with constant perturbation.

\section{Linear analysis and index theory}\label{sec:4d-linear}
Throughout this section, a perturbation $\pi$ on $W$ is fixed so that the limiting perturbations $\pi_i$ on $(Y_i,E_i)$ have finitely many critical orbits in $\mathfrak C_{\pi_i}$, all nondegenerate and such that each $\text{Hess}^\nu_{\alpha, \pi}$ has all eigenvalues of absolute value larger than a fixed constant $\delta$ (independent of $\alpha$).

We fix two such critical orbits, $\alpha_\pm = [A_\pm, E_b]$. Here recall that the elements of $\alpha_\pm$ vary over gauge equivalence classes of pairs $(A_\pm, p)$, where $p \in E_b$ is a point in the fiber of $E$ above the basepoint $b$, thought of as a framing of $E$ above that point.

To analyze the local structure of the moduli spaces around a framed instanton $(\mathbf{A},p)$ on $\mathbb R \times Y$ framed above $(0,b)$, we should restrict to the \textit{framed Coulomb slice} $\left(\mathbf{A}+\text{ker}(d_\mathbf{A}^*)\right) \times (\pE)_b \subset \widetilde{\mathcal A}_{k,\delta}(A_-, A_+)$: this is just the usual Coulomb slice with an additional $SO(3)$ coordinate for the framing. Every framed connection sufficiently close to $(\mathbf{A},p)$ is gauge equivalent to one in the framed Coulomb slice, and the representation is unique up to the action of $\Gamma_{\mathbf{A}} \hookrightarrow \text{Aut}(E_b) \cong SO(3)$ (the last isomorphism depending on the choice of framing at the basepoint). This acts trivially on the connection coordinate, and by translation in the framing coordinate.

We also find it convenient to pass to a smaller $\Gamma_A$-invariant open set in the framing coordinate, so that our subset is identified with 
\begin{equation}\label{eCoulomb}\text{ker}(d_{\mathbf{A}}^*) \times N(\Gamma_{\mathbf{A}}) \cong \mathcal K^{(4)}_{\mathbf{A},k} \times \Gamma_{\mathbf{A}} \times \mathfrak g_\mathbf{A}^\perp;\end{equation} 
we embed the normal bundle in $SO(3)$ by exponentiating. We refer to this as the \textit{extended Coulomb slice}. The action of $\Gamma_{\mathbf{A}}$ is again identified with left translation in the framing coordinate. Here $\mathfrak g_{\mathbf{A}}^\perp \subset \mathfrak g_b$ is the orthogonal complement to $\mathfrak g_{\mathbf{A}}$, the subset of $\mathfrak g_b$ which extends to $\mathbf{A}$-parallel sections; equivalently, this is the tangent space to $\Gamma_{\mathbf{A}}$. 

Inside the extended Coulomb slice around $(\mathbf{A},p)$, (say, $\widetilde {\mathcal K}^{(4)}_{k,\delta}$), there is a $\Gamma_\mathbf{A}$-equivariant trivialization of the bundle $\Omega^{2,+}_{k-1,\delta}(\text{End}(E))$, descending to a local trivialization of $\mathcal S_{k-1,\delta}$ over $\widetilde{\mathcal K}_{k,\delta}/\Gamma_\mathbf{A} \subset \widetilde{\mathcal B}^e_{k,\delta}(\alpha_-, \alpha_+)$. Using the second expression in (\ref{eCoulomb}), an element $q \in N(\Gamma_{\mathbf{A}})$ can be expressed uniquely as $u \cdot (e^{\xi}p)$ for $\xi \in \mathfrak g_\mathbf{A}^\perp$ and $u \in \Gamma_{\mathbf{A}}$; the equivariant trivialization sends 
$$(\mathbf{A}+a, u e^{\xi} p,\omega) \mapsto ((\mathbf{A}+a, u e^{\xi} p),  u \omega).$$
In this equivariant trivialization, the section defined by (\ref{ASDsec}) is given as 
$$(\pi, \mathbf{A}+a, u e^{\xi} p) \mapsto u\left(d_{\mathbf{A}}a + a \wedge a\right)^+ + u\left(\widehat \nabla_\pi(\mathbf{A}+a) - \widehat \nabla_\pi(\mathbf{A})\right).$$

We can thus identify the derivative of the perturbed map (for fixed $\pi$ and at $(\mathbf{A},p)$ a $\pi$-perturbed framed instanton) $\widetilde{\mathcal K}^{(4)}_{k,\delta} \to \Omega^2_{k-1,\delta}$ as the \textit{perturbed ASD operator}, $$d_{\mathbf{A}}^+a + (D_{\mathbf{A}}\nabla_{\pi})(a) \oplus 0: \mathcal K^{(4)}_{k,\delta} \oplus \mathfrak g_{\mathbf{A}}^\perp \to \Omega^{2,+}_{k-1,\delta};$$ we write this as $D_{\mathbf{A},\pi}$, and the corresponding perturbed \textit{normal ASD operator} as $$D^\nu_{\mathbf{A},\pi}: \mathcal K^{(4)}_{k,\delta} \to \Omega^{2,+}_{k-1,\delta}.$$ 

Just as in Chapter \ref{sec:3d-linear}, we gain control over the normal ASD operator with domain $\mathcal K^{(4)}_{k,\delta}$ by expressing it as a summand of a larger elliptic operator and including a gauge fixing condition. Define $Q_{\mathbf{A},\pi}$, the perturbed \textit{extended ASD operator}, by $$Q_{\mathbf{A},\pi}: \Omega^1_{k,\delta}(\pi^*\mathfrak g_E) \to \Omega^0_{k-1,\delta}(\pi^*\mathfrak g_E) \oplus \Omega^{2,+}_{k-1,\delta}(\pi^*\mathfrak g_E),$$ given by $$Q_{\mathbf{A},\pi}(a,\xi) = (d_{\mathbf{A}}^* a, d_{\mathbf{A}}^+ a + D_{\mathbf{A}}\widehat{\nabla}_\pi(a)).$$ This is the operator usually used in the instanton theory for integer homology spheres, as applied to weighted Sobolev spaces.

Now recall that we have an $L^2$-orthogonal splitting 
$$\Omega^1_{k,\delta}(\mathfrak g_{\mathbf{E}}) = \text{Im}(d_{\mathbf{A}}) \oplus \text{ker}(d_{\mathbf{A}}^*);$$
that every element may be written uniquely as a sum in this way is essentially the statement that the map $$\Delta_{\mathbf{A}}: \Omega^0_{k+1,\delta}(\mathfrak g_{\mathbf{E}})/\mathfrak g_{\mathbf{A}} \to \Omega^0_{k-1,\delta}(\mathfrak g_{\mathbf{E}})/\mathfrak g_{\mathbf{A}}$$ 
is an isomorphism, assuming $\delta$ is not an eigenvalue of the Laplacian of the limiting connections, which is true as long as $\delta$ is chosen sufficiently small.

Now in this splitting $\Omega^1 = \text{Im}(d_{\mathbf{A}}) \oplus \text{ker}(d_{\mathbf{A}}^*)$, and rewriting the first term as $\Omega^0_{k+1,\delta}/\mathfrak g_{\mathbf{A}}$ under the isomorphism $d_{\mathbf{A}}^{-1}$, we may write 
$$Q_{\mathbf{A},\pi} = \begin{pmatrix}\Delta_{\mathbf{A}} & 0 \\ 0 & D^\nu_{\mathbf{A},\pi}\end{pmatrix}.$$ 
The fact that $D^\nu_{\mathbf{A},\pi}$ takes values in $\text{ker}(d_{\mathbf{A}}^*)$ is the linearization of the statement that $F_{\mathbf{A}}^+ + \widehat{\nabla}_\pi(\mathbf{A})$ is a gauge invariant quantity. In particular, the normal ASD operator is Fredholm, and because the index of the top left operator is $- \dim \mathfrak g_{\mathbf{A}}$, the index of the normal ASD operator $D^\nu_{\mathbf{A},\pi}$ is the index of $Q_{\mathbf{A},\pi}$ plus the dimension of $\mathfrak g_{\mathbf{A}}$.

We recall the basic Fredholm property from \cite[Sections~3.2~and~3.3.1]{Don}.

\begin{proposition}\label{Fred1}Suppose $W$ is a Riemannian 4-manifold with cylindrical ends, isometric to $(-\infty, 0] \times Y_1$ and $[0,\infty) \times \overline{Y_2}$, equipped with an $SO(3)$-bundle $\mathbf{E}$ restricting to the pullbacks of fixed $SO(3)$-bundles $E_i$ on the ends. In this situation, we say that $(W,\mathbf{E})$ is a cobordism from $(Y_1, E_1)$ to $(Y_2, E_2)$. Suppose $\pi$ is a fixed perturbation on $\mathbf{E}$, restricting to regular perturbations $\pi_i$ on the ends, and $A_i$ are fixed nondegenerate critical points of $\textup{cs}_{Y_i} + f_{\pi_i}$. 

Let $\mathbf{A}$ be a choice of connection in $\mathcal A_{\mathbf{E},k,\delta}(A_1,A_2)$. 

If $\delta > 0$ is less than than the absolute value of any eigenvalue of $\widehat{\text{\em{Hess}}}^\nu_{A_i, \pi}$, then $Q_{\mathbf{A},\pi}$ is Fredholm, and has index independent of such $\delta$. If the $A_i$ are irreducible, we may even take $\delta = 0$.
\end{proposition}

There are two essential points. The first is that when a connection $\mathbf{A} = \mathbf{A}(t)$ on the cylinder is in temporal gauge, the operator $Q_{\mathbf{A},\pi}$ can be written in the form $$Q_{\mathbf{A}, \pi} = \frac{d}{dt} + \widehat{\text{Hess}}_{\mathbf{A}(t), \pi}.$$ (See \cite[Section~2.5]{Don}.) Secondly, if $e^{\sigma(t)} = f_\delta(t)$ is the function used to define the weighted Sobolev spaces, multiplication by $f_\delta$ is an isometry $L^2_{k,\delta} \to L^2_k$, and a first order linear differential operator $D$ of the form $\frac{d}{dt} + L_t$ is taken under this isometry to $\frac{d}{dt} + L_t - \sigma'(t)$. So to study our operator $Q^\nu_{\mathbf{A},\pi}$ on weighted Sobolev spaces, we should equivalently study 
$$\frac{d}{dt} + \widehat{\text{Hess}}_{\mathbf{A}(t),\pi} - \sigma'(t)$$ 
where $\sigma(t) = -\delta t$ for $t \ll 0$ and $\sigma(t) = \delta t$ for $t \gg 0$. (This is well-explained in \cite[Section~3.3.1]{Don} and the beginning of \cite[Section~3.3]{Lin}.) Once we have done this, the above result for the ASD operator is a consequence of general theory; it is proven for $L_t$ a family of almost self-adjoint first order differential operators as \cite[Proposition~14.2.1]{KMSW}.

An important consequence of this description is that for $\mathbf{A} = \mathbf{A}(t)$ in temporal gauge, we can describe the index of $Q_{\mathbf{A},\pi}$ as the \textit{spectral flow} of the family of operators $\widehat{\text{Hess}}_{\mathbf{A}(t), \pi} - \sigma'(t)$ between $\widehat{\text{Hess}}_{A_1, \pi} + \delta I$ and $\widehat{\text{Hess}}_{A_2, \pi} - \delta I$, the (algebraic) intersection number of the paths the eigenvalues take with $0 \in \mathbb R$. 

Because each $\mathcal A_{\mathbf{E},k,\delta,z}(A_1,A_2)$ is connected (even contractible) and index is a homotopy invariant of self-adjoint Fredholm operators, the index of $Q_{\mathbf{A},\pi}$ only depends on the homotopy class $z$, not the actual choice of connection $\mathbf{A}$. For the same reason, this also agrees with the index of $Q^\nu_{\mathbf{A},\pi'}$ for any other perturbation which restricts to the same perturbation sufficiently far on the ends. 

It is perhaps worth observing that the index of $D^\nu_{\mathbf{A},\pi}$, defined to be the derivative of the section operator $\widetilde{\mathcal B}^e_{\mathbf{E},k,\delta} \to \mathcal S_{k-1,\delta}$ normal to an orbit, \emph{does} have its index jump as we pass from irreducibles to reducibles. This makes sense, thinking of $D^\nu$ as measuring the expected codimension of the orbit through $\mathbf{A}$ in the entire moduli space: the codimension is larger at smaller-dimensional orbits.

Following the definition after \cite[Lemma~3.13]{KM1}, we use this to give the following definition.

\begin{definition}\label{grading}In the situation of Proposition \ref{Fred1}, let $z$ denote 
a connected component of $\widetilde{\mathcal A}_{\mathbf{E},k, \delta}(\alpha, \beta)$. We write the \textit{unframed grading} $\overline{\text{\em{gr}}}^W_z(\alpha, \beta) = \text{\em{ind}}(Q^\nu_{\mathbf{A},\pi})$ for any choice of $\mathbf{A}$ in the component $z$, and the \textit{relative grading} between the orbits $\alpha$, $\beta$, with respect to the path $z$, is $\text{\em{gr}}_z(\alpha, \beta) = \overline{\text{\em{gr}}}_z(\alpha, \beta) + 3 - \dim \alpha$.

When $W$ is the cylinder with the constant perturbation, we drop the superscript $W$.
\end{definition}

The relative grading here will be quite natural in the definition of the framed instanton Floer complex, whose differential is defined in terms of fiber products with moduli spaces; $\text{gr}^W_z(\alpha, \beta)$ is the expected dimension of a fiber of the map $$\text{ev}_-: \widetilde{\mathcal M}_{\mathbf{E},\pi,z}(\alpha, \beta) \to \alpha.$$ It is worth noting here that $3 - \dim \alpha = \dim \mathfrak g_{\alpha} = \dim \Gamma_{\alpha}$.

The more immediate point of this definition is that the relative grading is \textit{additive}. In the particular case that $W$ is isometric to the cylinder $\mathbb R \times Y$, $z$ corresponds to a relative homotopy class between $[\alpha]$ and $[\beta]$ in $\mathcal B^e_{E,k}$. If $w$ is a path from $[\beta]$ to $[\gamma]$, then $$\text{gr}_z(\alpha, \beta) + \text{gr}_w(\beta, \gamma) = \text{gr}_{z \cdot w}(\alpha, \gamma)$$ ($z \cdot w$ the concatenated path). This follows by computing these as a spectral flow: if $\mathbf{A}_z(t) = \mathbf{A}_z \in \widetilde{\mathcal A}_{k, \delta,z}(\alpha, \beta)$ is in temporal gauge and $\mathbf{A}_z(t) = \alpha$ (resp. $\beta$) for $t \ll 0$ (resp $t \gg 0$), the index of $Q_{\mathbf{A}_z,\pi}$ is the spectral flow of the path
$$\widehat{\text{Hess}}_{A_z(t),\pi} - \sigma'(t),$$ 
where $\sigma'(t) = - \delta$ for $t \ll 0$ and $\sigma'(t) = \delta$ for $t \gg 0$. Making a similar choice of $\mathbf{A}_w$, if we glue together $\mathbf{A}_z$ and $\mathbf{A}_w$ sufficiently far out on the ends, we can find a connection $\mathbf{A}_{z \cdot w}$ in the component corresponding to $z \cdot w$ as the concatenation of the paths $\mathbf{A}_z(t)$ and $\mathbf{A}_w(t)$. However, we cannot concatenate the corresponding paths of self-adjoint operators yet; for large $t$, the first ends at $\widehat{\text{Hess}}_{\beta,\pi} - \delta$ and the second begins at $\widehat{\text{Hess}}_{\beta,\pi} + \delta$. To actually concatenate them, we must traverse the path

$$\widehat{\text{Hess}}_{\beta,\pi} - (1-2t)\delta;$$ 

doing so changes $(3-\dim \beta) = \dim \text{ker}\left(\widehat{\text{Hess}}_{\beta, \pi}\right)$ negative eigenvalues to positive, and so

$$\overline{\text{gr}}_z(\alpha, \beta)  +(3-\dim \beta) + \overline{\text{gr}}_w(\beta, \gamma) = \overline{\text{gr}}_{z \cdot w}(\alpha, \gamma).$$

Additivity of this grading for general cobordisms is also true. We record this as a proposition.

\begin{proposition}\label{gr-add}Suppose we have cobordisms $(W_1, \mathbf{E}_1)$ from $(Y_1, E_1)$ to $(Y_2, E_2)$ and $(W_2, \mathbf{E}_2)$ from $(Y_2, E_2)$ to $(Y_3, E_3)$, equipped with paths $\gamma_1: \mathbb R \to W_1$ and $\gamma_2: \mathbb R \to W_2$ between the basepoints of the $Y_i$. Suppose the $W_i$ are equipped with perturbations $\hat \pi_i$, which restrict to regular perturbations $\pi_j$ on each $Y_j$, and furthermore suppose each $(Y_j, E_j, \pi_j)$ is equipped with some $\pi_j$-flat connection $A_j$. 

We can define the composed cobordism $(W^T_{12}, \mathbf{E}_{12})$ from $(Y_1, E_1)$ to $(Y_3, E_3)$, identifying $(T, \infty) \times \overline Y_2$ with $(-\infty, -T) \times Y_2$ on the cylindrical ends; there is a corresponding perturbation $\hat \pi$ interpolating between $\pi_1$ and $\pi_3$. We denote by $z_i$ a component of $\widetilde{\mathcal A}_{\mathbf{E}_i, k, \delta}(A_i, A_{i+1})$; there is a component $z_1 \ast z_2$ corresponding to gluing representative connections of these components along the ends. In this situation, we have $$\textup{gr}^{W_1}_{z_1}(\alpha_1, \alpha_2) + \textup{gr}^{W_2}_{z_2}(\alpha_2, \alpha_3) = \textup{gr}^{W^T_{12}}_{z_1\ast z_2}(\alpha_1, \alpha_3).$$ 
\end{proposition}

This follows from the additivity theorem of the index when the limiting operators over the ends have no kernel (\cite[Proposition~3.9]{Don}) and the relation of operators on weighted Sobolev spaces to unweighted spaces, given by conjugating by the weighting function: the operator $Q_{\mathbf{A}_1, \pi_1}$ on its positive end, after conjugating by the weighting function $e^\sigma$, takes the form $\frac{d}{dt} + \widehat{\text{Hess}}_{A_2,\pi_2} - \delta$ and on the negative end of $W_2$ takes the form $\frac{d}{dt} + \widehat{\text{Hess}}_{A_2, \pi_2} + \delta$. To glue these we first need to interpolate between $-\delta$ and $+\delta$, moving the 
$$(3-\dim \alpha_2) = \dim \text{ker}\left(\widehat{\text{Hess}}_{A_2,\pi_2}\right)$$ 
negative eigenvalues across $0$; this is observed as \cite[Proposition~3.10]{Don}, identifying the index on weighted spaces with what he denotes $\text{ind}^+(P)$.

\section{Uhlenbeck compactness for framed instantons}\label{sec:4d-compact}
The following definition is precisely \cite[Definition~16.1.1]{KMSW}. Note that there are no framings involved yet.

\begin{definition}Let $Y$ be a Riemannian 3-manifold, equipped with $SO(3)$-bundle $E$ and regular perturbation $\pi$; then $\pi$ has finitely many critical points in $\mathcal B^e_{E,k}$, which we write a generic point of as $\alpha$. We say that a \emph{trajectory} from $\alpha_-$ to $\alpha_+$ is an equivalence class of nonconstant $\pi$-perturbed instanton $\mathbf A \in \mathcal M_{E}(\alpha_-, \alpha_+) \subset \mathcal B^e_{\pE,k,\delta}(\alpha_-,\alpha_+)$ on $\mathbb R \times Y$ under the translation action. 

The \emph{homotopy class} of a trajectory is the element of $\pi_1(\mathcal B^e_{E,k}, \alpha_-, \alpha_+)$ it traces out; these are in noncanonical bijection with $\mathbb Z$. The topological energy of a trajectory was defined above as 
$$2((\textup{cs}+f_\pi)(\textup{ev}_+  \mathbf{A}) -  (\textup{cs} + f_\pi)(\textup{ev}_-\mathbf{A}));$$ 
even though the Chern-Simons functional of an individual connection is only defined in $\mathbb R/8\pi^2\mathbb Z$, this difference of boundary components of a connection on a cylinder is defined in $\mathbb R$.

A \emph{broken trajectory} from $\alpha_-$ to $\alpha_+$ consists of a finite sequence of $\pi$-perturbed instantons $\mathbf{A}_i$ on $\pE$ over $\mathbb R \times Y$ (say $1 \leq i \leq n$), with $$\textup{ev}_-\mathbf{A}_1 = \alpha_-, \; \; \textup{ev}_+ \mathbf{A}_n = \alpha_+, \;\; \text{ and } \textup{ev}_+ \mathbf{A}_i = \textup{ev}_- \mathbf{A}_{i+1} \text{  for 0 }< \text{ i } < \text{ n}.$$ The homotopy class of a broken trajectory is the composite of the homotopy classes given by the individual trajectories; the energy of a broken trajectory is defined as $$\mathcal E^\pi(\mathbf{A}_1, \cdots, \mathbf{A}_n) = \sum_{i=1}^n \mathcal E^\pi(\mathbf{A}_i).$$
\end{definition}

We topologize this exactly as in \cite[Page~276]{KMSW}. As per the author's taste, we present this in terms of sequences, following \cite[Page~116]{Don}: give the space of broken trajectories the final topology so that the following sequences of unbroken instantons, and their natural generalizations to sequences of broken instantons, converge to their stated limits: $\mathbf{A}_i$, where $i \in \mathbb N$, converges to the broken trajectory $(\mathbf{B}_1, \cdots, \mathbf{B}_n)$, if there is a sequence $(T_i^1, \cdots, T_i^n)$ of real numbers with $T_i^j \leq T_i^{j+1}$ for $0 < j < n$, and the successive differences $T^{j+1}_i - T^j_i \to \infty$ as $i \to \infty$, so that the pullbacks $\tau_{T^i_j}^* \mathbf{A}$ converge as $i \to \infty$ to $\mathbf{B}_j$ in the $L^2_{k,\delta}$ topology. Note that the energy of a trajectory is continuous with respect to this `chain convergence'. 

This form of noncompactness, trajectories breaking into a composite of lower-index trajectories, is familiar in Morse theory. There is another kind of noncompactness familiar in the instanton theory: Uhlenbeck bubbling.

\begin{definition}Let $(Y,E,\pi)$ be as above. An \emph{ideal instanton} is a solution $\mathbf{A}$ to the $\pi$-ASD equations on $(\mathbb R \times Y, \pE)$, along with a finite (possibly empty) collection of points $x_i \in \mathbb R \times Y$ and integer weights $k_i \geq 1$ at each point. An \emph{ideal trajectory} is an $\mathbb R$-equivalence class of nontrivial ideal instantons, where nontrivial means either that $\mathbf{A}$ is nonconstant or that $\{x_i\}$ is nonempty; equivalently, a nontrivial ideal instanton is one so that $\tau_t^*(\mathbf{A},x,k)$ is not gauge equivalent to $(\mathbf{A},x,k)$ for any $t$. Similarly, a \emph{broken ideal trajectory} is a sequence of ideal trajectories $(\mathbf{A}_i)_{i=1}^n$ with $\textup{ev}_+ \mathbf{A}_i = \textup{ev}_- \mathbf{A}_{i+1}$.

The energy of an ideal instanton $(\mathbf{A},x_i,k_i)$, where $1 \leq i \leq n$, is 
$$\mathcal E^\pi(\mathbf{A}) + 16\pi^2 \sum_i k_i.$$ 
If the homotopy class of the instanton $\mathbf{A}$ is $z$, the homotopy class of the ideal instanton $(\mathbf{A},x_i,k_i)$ is $z + \sum k_i$. The energy and homotopy class of a broken ideal instanton are defined to be additive under concatenation.
\end{definition}

We say that a sequence of $\pi$-perturbed instantons $\mathbf{A}_n$ converges to an ideal instanton $(\mathbf{A},x_i,k_i)$ if there is a sequence of gauge transformations $\sigma_n$ defined on $(\mathbb R \times Y) \setminus \{x_i\}$ such that $\sigma_n^* \mathbf{A}_n \to \mathbf{A}$ in the $L^p_1$ topology on compact subsets of $(\mathbb R \times Y) \setminus \{x_i\}$ for all $\infty > p \geq 2$, and such that the density measures converge: 
$$2|F_{\mathbf{A}_n}|^2 \to 2|F_{\mathbf{A}}|^2 + \sum_{i=1}^n 64\pi^2 k_i \delta_{x_i}.$$

There is then a natural extension of this to a definition of convergence to (and of) broken ideal trajectories. The space of broken ideal trajectories from $\alpha$ to $\beta$ in the homotopy class $z$ is written $\check{\mathcal M}_{E,z,\pi}(\alpha, \beta)$. Note that this does \emph{not} include constant trajectories, or broken `trajectories' for which one of the components is constant.

As observed in \cite{Kronheimer}, the non-local nature of holonomy perturbations means that we cannot expect better convergence than $L^p_1$ under bubble-limits.

The Uhlenbeck compactness theorem for the ASD equations on the cylinder is the following.

\begin{proposition}The subspace of broken ideal trajectories in $\check{\mathcal M}_{E,z,\pi}(\alpha, \beta)$ with a fixed energy bound $\mathcal E^\pi(\mathbf{A}) \leq C$ is compact.
\end{proposition}

We do not repeat the proof, which can be seen in \cite[Section~5.1]{Don}; the corresponding fact for compact cylinders is \cite[Proposition~3.20]{KM1}.

A somewhat stronger statement, Proposition \ref{properflowlines} below, is true; we will only use it briefly, but find it to be somewhat interesting. 

As in Definition \ref{pertdef}, for any real numbers $\epsilon \gg \delta > 0$, we let $\mathcal P_{E,\delta} \subset \mathcal P_E$ be the open subset of perturbations on $(Y,E)$ so that, for each $\pi$-critical point $\alpha$, all eigenvalues of $\widehat{\text{Hess}}^\nu_{\alpha,\pi}$ have absolute value larger than $\delta$, and so that for every pair $(\pi, \alpha)$ of perturbation and $\pi$-flat \emph{reducible} connection, the operator $\widehat{\text{Hess}}^\nu_{\alpha,\pi}$ has no eigenvalues of absolute value $\epsilon$. (In particular, we assume $\pi$ is a regular perturbation.) 

Because the projection $\mathfrak C_\pi \to \mathcal P_{E,\delta}$ of the parameterized critical set (in $\widetilde{\mathcal B}^e_E$) to the space of perturbations is a proper submersion, it is in particular a locally trivial fiber bundle. Therefore, for some small open set $U$ around any $\pi_0 \in \mathcal P_{E,\delta}$, we have a canonical bijection $\mathfrak C_{\pi_0} \cong \mathfrak C_{\pi}$ for any $\pi \in U$; in fact, if we fix $\alpha \in \mathfrak C_{\pi_0}$, we may choose a smooth map $s_\alpha: U \to \mathcal A_{E,k}$, so that $s(\pi)$ is a $\pi$-flat connection which is identified under the above bijection with $\alpha$. Choose once and for all, for each homotopy class $z$, a smooth map 
$$r_z: \mathcal A_{E,k} \times \mathcal A_{E,k} \to \mathcal A^{(4)}_{\pE,z,k,\delta},$$ 
sending $(A_-, A_+)$ to a connection which is constant at $A_-$ for $t \leq -1$ and constant at $A_+$ for $t \geq 1$ and in the homotopy class $z$. 

We may use the $s_\alpha$ to define the parameterized configuration space $\mathcal P_\delta \mathcal A^{(4)}_{\pE,z,k,\delta}$, whose elements $(\pi, \alpha_\pm, \mathbf{A})$ consist of a perturbation $\pi \in \mathcal P_{E,\delta}$, a choice of two $\pi_\pm$-flat connections $\alpha_\pm$, and a connection $\mathbf{A} \in \mathcal A^{(4)}_{\pE,z,k,\delta}(\alpha_-, \alpha_+)$.

This set is given the structure of a smooth Banach manifold by patching together charts of the form $$U \times \mathcal A_{\pE,z,k,\delta}(\alpha_-,\alpha_+) \cong r_z(s_{\alpha_-} \pi, s_{\alpha_+} \pi) + \Omega^1_{k,\delta}(\mathfrak g_E).$$ Write $$\mathcal P_\delta \widetilde{\mathcal A}^{(4)}_{\pE,z,k,\delta} = \mathcal P_\delta \mathcal A^{(4)}_{\pE,z,k,\delta} \times \pE_{(0,b)}$$ for the parameterized space of framed connections, which inherits a smooth structure and a smooth right action of $SO(3)$. It carries a smooth projection map 
$$\mathcal P_{\delta} \widetilde{\mathcal A}^{(4)}_{\pE,z,k,\delta} \to \mathcal P_{E,\delta}.$$
It carries the action of a bundle of Banach Lie groups over $\mathcal P_{E,\delta}$; the fiberwise quotient gives a topological space $\mathcal P_{\delta} \widetilde{\mathcal B}^e_{\pE,z,k,\delta}$. We may define the parameterized moduli space 
$$\mathcal P_{\delta} \widetilde{\mathcal M} \subset \mathcal P_{\delta} \mathcal B^e_{\pE,z,k,\delta}$$ 
as the equivalence classes of triples $(\pi, \mathbf{A},p)$, where $\pi \in \mathcal P_{E,\delta}$ is a perturbation, $\mathbf{A}$ is a $\pi$-perturbed instanton, and $p$ is a framing. 

Because the perturbed ASD equations do not depend on the framing, this set inherits the right $SO(3)$-action. The quotient of $\mathcal P_{\delta} \widetilde{\mathcal M}$ by this $SO(3)$ action is the parameterized moduli space of instantons, which we denote $\mathcal P_\delta \mathcal M$. If we take the quotient by the $\mathbb R$ action, throw out the constant trajectories, and incorporate bubble-limits into the topology, we may extend this to a space $\mathcal P_\delta \check{\mathcal M}$, the parameterized moduli space of ideal broken trajectories. (We will soon discuss the version of this appropriate to the framed setting, which is slightly more subtle.)

\begin{proposition}\label{properflowlines}Let $\mathcal P_{\delta} \check{\mathcal M}^{\leq C}$ be the subspace of $\mathcal P_{\delta} \check{\mathcal M}$ consisting of those pairs $(\pi, \mathbf{A})$ so that $\mathcal E^\pi(\mathbf{A}) \leq C$. Then the projection map $\mathcal P_{\delta} \check{\mathcal M}^{\leq C} \to \mathcal P_{\delta}$ is proper.
\end{proposition}

\begin{proof}Suppose we have a sequence $(\pi_n, \mathbf{A}_n)$ of perturbations and unbroken instantons so that $\pi_n \to \pi$. We want to show that there is a subsequence of $\mathbf{A}_n$ which converges to a broken ideal $\pi$-trajectory. (The general case where $\mathbf{A}_n$ is itself a broken ideal trajectory provides no further difficulty.) First, because the possible $\text{ev}_\pm$ take values in a finite set, choose a subset of $\mathbf{A}_n$ so that, for $n$ large, $\text{ev}_\pm \mathbf{A}_n$ correspond to $\alpha_\pm \in \mathfrak C_\pi \cong \mathfrak C_{\pi_n}$. We may now apply \cite[Lemma~4.3]{Don}, which establishes that any $\pi$-instanton with sufficiently small energy and $\text{ev}_- \mathbf{A} = \alpha$ is gauge equivalent on $(-\infty, 0) \times Y$ to $\alpha + a$ for some $a$ with a uniform bound on $|a(t)| e^{-\delta t}$, as well as the derivatives $|\nabla^{(\ell)}(a)(t)| e^{\delta t}$ for $\ell \leq k$. The constants in these uniform bounds are bounded for the convergent sequence $\pi_n \to \pi$, and $\delta$ is fixed. In particular, for our solutions $\mathbf{A}_n$, gauge equivalent to $\alpha_n + a_n$, this is enough for the Arzela-Ascoli and dominated convergence argument in \cite[Lemma~5.1]{Don} to imply that $a_n \to a$ for some function $a$ satisfying the same bounds; in fact we must have 
$$F_{\alpha + a}^+ + \widehat{\nabla}_\pi(\alpha+a) = 0$$ 
as this is the pointwise limit of the corresponding equations 
$$F_{\mathbf{A}_n}^+ + \widehat{\nabla}_\pi(\mathbf{A}_n) = 0.$$

Now that we have control over the ends, everything else is standard: a uniform bound on $\mathcal E^{\pi_n}(\mathbf{A}_n)$ implies a uniform bound on the $L^2$ norm of $F_{\mathbf{A}}$ on compact sets, so one has a limit on compact sets after accounting for bubbling, and then a limit on the whole line to a broken ideal trajectory; this uses that there is a uniform positive lower bound on the minimal energy of a nontrivial $\pi_n$-instanton. To see this, recall that the energy of an instanton may be written as $2\int |\mathbf{A}'(t)|^2$ as in the proof of Proposition \ref{energy}. If $\mathbf{A}_n$ is a sequence of nontrivial $\pi_n$-instantons with $\mathcal E^{\pi_n}(\mathbf{A}_n) \to 0$, this implies that the distance between the endpoints $\alpha_\pm^n$ goes to zero. But $\alpha_{\pm}^n \to \alpha_\pm$, so $\alpha_- = \alpha_+$; but then because $\alpha^n_\pm$ is sent to $\alpha_\pm$ under the bijection $\mathfrak C_{\pi_n} \cong \mathfrak C_{\pi}$, we see that $\alpha_-^n = \alpha_+^n$ for large $n$. However, this would imply that $\mathcal E^{\pi_n}(\mathbf{A}_n)$ is a multiple of $64\pi^2$ for large $n$, and also may be made arbitrarily close to zero by taking $n$ large; therefore $\mathcal E^{\pi_n}(\mathbf{A}_n) = 0$ for large $n$. However, this implies $\mathbf{A}_n'(t) = 0$, but we assumed $\mathbf{A}_n$ was nontrivial.

The reason we restrict to perturbations in $\mathcal P_{E,\delta}$ is for the definition of the parameterized moduli space as $L^2_{k,\delta}$ connections; we always want to take $\delta$ less than the eigenvalues of the extended Hessian, and otherwise would need to choose $\delta$ depending on the perturbation.
\end{proof}

We conclude the discussion of compactness for cylinders by defining the object of interest to us: the compactification of the \emph{framed} moduli space of trajectories.

\begin{definition}Let $Y$ be a Riemannian 3-manifold, equipped with $SO(3)$-bundle $E$ and regular perturbation $\pi$. We say that a \emph{framed ideal trajectory} from $\alpha_-$ to $\alpha_+$ is an equivalence class of $\pi$-perturbed framed instanton $(\mathbf A,p) \in \widetilde{\mathcal B}^e_{\pE,k,\delta}(\alpha_-,\alpha_+)$ equipped with a (possibly empty) points of points $x \in \mathbb R \times Y$ and positive integer weights $k_x$, where none of the $x$ are the basepoint $(0,b)$. We demand the ideal instanton is nontrivial, in the sense that either the set of points $x_i$ is nonempty or the trajectory $\mathbf{A}$ is nonconstant.

A \emph{deframed ideal trajectory} from $\alpha_-$ to $\alpha_+$ is an ideal trajectory from $\alpha_-$ to $\alpha_+$ so that $(0,b)$ is a weight-point $x$ with $k_{(0,b)} > 0$.\footnote{Here deframed is meant to indicate that the framing has been removed via the placement of a $\delta$-mass at the basepoint; not all the weight-sets are allowed, and in particular a non-ideal instanton without a framing is not a `deframed ideal trajectory'.}

A \emph{framed broken trajectory} from $\alpha_-$ to $\alpha_+$ is a finite sequence whose elements are either framed ideal trajectories $(\mathbf{A}_i,p_i,x_i, k_i)$ or deframed ideal trajectories $(\mathbf{A},x_i,k_i)$ with
\begin{align*}&\textup{ev}_-(\mathbf{A}_1,p_1) \in \alpha_-, \; \; \textup{ev}_+ (\mathbf{A}_n,p_n) \in \alpha_+, \;\; \textup{and}\\
\textup{ev}_+ (\mathbf{A}_i,p_i) &= \textup{ev}_- (\mathbf{A}_{i+1},p_{i+1}) \textup{  for }0 < i < n,\textup{ when } \mathbf{A}_i, \mathbf{A}_{i+1} \textup{ are both framed}.\end{align*}

The set of framed broken trajectories in the homotopy class $z$ from $\alpha_-$ to $\alpha_+$ is written $\overline{\mathcal M}_{E,z,\pi}(\alpha_-, \alpha_+)$. We write a generic element as $(\mathbf{A},p)$, even though $\mathbf{A}$ may be broken and not every piece must be framed.
\end{definition}

This set is topologized as follows. A sequence $(\mathbf{A}_n,p_n)$ of framed instantons converges to a framed ideal instanton $(\mathbf{A},p,x_i,k)$ if there is a sequence of gauge transformations $\sigma_n$, defined on $(\mathbb R \times Y) \setminus \{x_i\}$, and in particular defined on the basepoint $(0,b)$, so that $\sigma_n^*\mathbf{A}_n \to \mathbf{A}$ converge in the sense of ideal instantons above, and $\sigma_n^*p_n \to p$.

However, if the underlying trajectores $\mathbf{A}_n$ converge to an ideal trajectory with nontrivial weight at $(0,b)$, so that the $\sigma_n$ are undefined at $(0,b)$, then it doesn't make sense to compare $\sigma_n^*p_n$ and $p$. In this situation \emph{we lose the framing via bubbling at the basepoint}. In this case, the sequence $(\mathbf{A}_n,p_n)$ of framed instantons converges to a deframed ideal trajectory. Incorporating this into the topology on broken trajectories is straightforward.

\begin{proposition}The space $\overline{\mathcal M}_{E,z,\pi}(\alpha_-, \alpha_+)$ is compact.
\end{proposition}

\begin{proof}We begin by recalling that there is a surjective map $$\overline{\mathcal M}_{E,z,\pi}(\alpha_-, \alpha_+) \to \check{\mathcal M}_{E,z,\pi}(\alpha_-,\alpha_+).$$ This map is proper: if a sequence $(\mathbf{A}_n,p_n)$ of framed broken trajectories has underlying sequence of broken ideal trajectories converge to $\mathbf{A}$, the framings are either incomparable (and that limit component of $\mathbf{A}_n$ is a deframed ideal trajectory), or (an appropriate sequence of translations of) $\sigma_n^*p_n$ are all defined, and live in the compact space $E_b \cong SO(3)$; so some subsequence converges, as desired. Doing this for the finitely many components of the limit $\mathbf{A}$ constructs an element of $\overline{\mathcal M}$ that a subsequence of $(\mathbf{A}_n, p_n)$ converges to.
Because we have constructed a proper map to a compact space, the total space $\overline{\mathcal M}$ is compact.
\end{proof}

Similarly, the projection from the parameterized moduli space $\mathcal P_\delta \overline{\mathcal M}^{\leq C}$ to $\mathcal P_{E, \delta}$ is proper.

There are straightforward extensions of these to moduli spaces on cobordisms, which we now state. In what follows, we write $\mathcal P^{(4)}_\delta$ for the subspace $$\mathcal P_{E_1, \delta} \oplus \mathcal P^{(4)}_{\text{int}} \oplus \mathcal P_{E_2, \delta}.$$  
(As usual, the manifold $W$ is left implicit in the notation.) 

\begin{definition}
Let $W$ be an oriented Riemannian 4-manifold, with cylindrical ends oriented isometric to $(-\infty, 0] \times Y_1$ and $[0,\infty) \times Y_2$, equipped with an $SO(3)$-bundle $\mathbf{E}$ and specified isomorphisms to the pullback of bundles on $Y_i$ over the ends. Suppose $W$ is equipped with an embedding $\gamma:\mathbb R \hookrightarrow W$ which agrees for $|t|$ sufficiently large with $(t,b_1)$ or $(t,b_2)$, depending on the sign of $t$, and write $b = \gamma(0)$ as the basepoint of $W$. Suppose $W$ is equipped with a perturbation $\pi$, restricting to fixed regular perturbations $\pi_i$ on the ends. Connected components of the space of connections on $W$ limiting to $\alpha_\pm$, critical points of $\pi_\pm$, are in bijection with $\mathbb Z$, and labeled by $z$.

Then a \emph{broken ideal $W$-trajectory} from $\alpha_-$ to $\alpha_+$ in the homotopy class $z$ is a triple of a broken ideal $\pi_-$-trajectory on $Y_1$ from $\alpha_-$ to some $\beta$ (possibly constant, if $\alpha_- = \beta$), an ideal $\pi$-instanton on $W$ from $\beta$ to some $\gamma$, and a broken ideal $\pi_+$-trajectory on $Y_2$ from $\gamma$ to $\alpha_+$ (possibly constant, if $\gamma = \alpha_+$). The homotopy class is the composite of the corresponding homotopy classes of non-ideal instantons, then summing the weights. We denote the set of broken ideal $W$-trajectories from $\alpha_-$ to $\alpha_+$ in the homotopy class $z$ as $\check{\mathcal M}^{W,\pi}_{\mathbf{E},z}(\alpha_-, \alpha_+)$.

A \emph{framed broken ideal $W$-trajectory} is the same, but with in addition a framing on each component where there is no weight at the basepoint of $Y_i$ or $W$, and so that the neighboring evaluations agree: $$\textup{Hol}^{0 \to \infty}_{\mathbf{A}_i,\gamma_i}(p_i) = \textup{Hol}^{0 \to -\infty}_{\mathbf{A}_{i+1}, \gamma_{i+1}}(p_{i+1}).$$ We denote the set of framed broken ideal $W$-trajectories from $\alpha_-$ to $\alpha_+$ in the homotopy class $z$ as $\overline{\mathcal M}_{\mathbf{E},z,\pi}(\alpha_-, \alpha_+)$. We similarly denote the parameterized space by $$\mathcal P_{L,\delta} \overline{\mathcal M}_{\mathbf{E},z}.$$
\end{definition}

Note in particular that we include the `deframed trajectories' here as in the cylindrical case, and at these, we have lost the framing at the basepoint. Here $(W,\mathbf{E})$ is collapsed in the notation to simply the bundle $\mathbf{E}$.

A sequence of instantons $\mathbf{A}_n$ on $W$ from $\alpha_-$ to $\alpha_+$ converges to a broken trajectory $(\mathbf{A}_-, \mathbf{A}, \mathbf{A}_+)$ if the following hold. First, for each end, we demand that there is a sequence of real numbers $T_n^j$ with 
\begin{align*}
\lim_{j \to \infty} T_n^{j+1} - T_n^j &= \pm \infty \\
\lim_{n \to \infty} T_n^0 &= \pm \infty;
\end{align*}

the limit should be negative on the negative end and and positive on the positive end. These real numbers should be such that the translates $\tau_{T_n^j}^* \mathbf{A}_n$ converge to $(\mathbf{A}_-)_j$ in 
$$\mathcal B_{\pE,k,\delta}\left((-\infty, N] \times Y_1\right)$$ 
for all $N$, and similarly for the positive end. Note that the assumption that $T_n^0 \to \pm \infty$ implies that these translates are eventually defined on all of $(-\infty, N] \times Y_1$ for any fixed $N$, and similarly $[-N, \infty) \times Y_2$. 

Second, we demand that $\mathbf{A}_n \to \mathbf{A}$ in $\mathcal B_{\mathbf{E},k,\delta}$. 

One similarly accounts for the bubbling phenomenon to define limits to broken \textit{ideal} trajectories: there should be a sequence of gauge transformations defined away from the bubble points so that convergence of $\sigma_n(\mathbf{A}_n)$ is $L^p_1$ on compact sets away from the bubble-points for all $\infty > p \geq 2$. 

We have the following analogue of Proposition \ref{properflowlines} for $W$-trajectories, which will be useful to us. The proof is a verbatim combination of what was already said in the case of the ends (that is, in the cylindrical case), as well as a discussion of what happens on compact subsets of $W$.

\begin{proposition}\label{propertrajectories}Let $\mathcal P_{\delta} \check{\mathcal M}_{\mathbf{E}}^{\leq C}$ be the subspace of $\mathcal P_{\delta} \check{\mathcal M}_{\mathbf{E}}$ consisting of those pairs $(\pi, \mathbf{A})$ so that $\mathcal E^\pi(\mathbf{A}) \leq C$. Then the projection map $\mathcal P_{\delta} \check{\mathcal M}_{\mathbf{E}}^{\leq C} \to \mathcal P^{(4)}_{\delta}$ is proper.\end{proposition}

\begin{proof}Fix $(\mathbf{A}_n, \pi_n)$ with $\pi_n \to \pi$. We have already dealt with the ends in the discussion of the cylinder, and so we will focus on $W$ itself. We need to show that there is some sequence of gauge transformations, defined on the complement of some finite set, so that there is a connection $\mathbf{A}$ and a subsequence (still written $\mathbf{A}_n$) so that $\sigma_n^*\mathbf{A}_n \to \mathbf{A}$. What we need is to show that that a uniform energy bound gives, over any given compact set, a uniform bound on the $L^2$ norm of $F_{\mathbf{A}_n}$. Because $$F^+_{\mathbf{A}_n} = \widehat{\nabla}_{\pi_n}(\mathbf{A}_n)$$ is bounded in $L^\infty$ independent of $n$, Uhlenbeck's compactness theorem will then guarantee that there is a sequence of gauge transformations $\sigma_n$ defined on the complement of a finite set $\mathbf{x}$, and an $L^p_1$ connection $\mathbf{A}$ in Coulomb gauge with respect to some reference smooth connection, so that $\sigma_n(\mathbf{A}_n) \to \mathbf{A}$ in $L^p_1$ on compact sets of $W \setminus \mathbf{x}$. Proposition \ref{pertfacts} (6) guarantees that 
$$\widehat{\nabla}_{\pi_n}(\sigma_n)_*(\mathbf{A}_n) \to_{L^p} \widehat{\nabla}_\pi(\mathbf{A}),$$ 
and so 
$$0 = (\sigma_n)_*\left(F^+_{\mathbf{A}_n} + \widehat{\nabla}_{\pi_n}(\mathbf{A}_n)\right) \to_{L^p} F^+_{\mathbf{A}} + \widehat{\nabla}_\pi(\mathbf{A});$$ 
therefore, $\mathbf{A}$ is indeed a $\pi$-instanton. Because $\mathbf{A}$ is in Coulomb gauge with respect to a reference smooth connection, the usual bootstrapping arguments imply $\mathbf{A}$ is smooth.

We will further need to ensure that if $\mathbf{A}_n$ converges to an ideal perturbed instanton, the weight of the $\delta$-mass at the ideal point $x$ is a multiple of $64\pi^2$.

First, we will resolve the question of curvature. Consider the restriction of $\mathbf{A}$ to the complement of the ends $[1,\infty) \times Y$, a compact manifold $W'$ with boundary $Y$. The perturbations vanish on the boundary of $W'$ (and on the interior are the previously discussed interior holonomy perturbations). On the ends the perturbations $\pi_n$ are of the form (eg, for $[0,\infty) \times Y_1$) $$\beta_0(t) (\pi_-)_n,$$ where $\beta_0: [0,\infty) \to [0,1]$ is zero for $t \leq 1$ and $1$ for $t \geq 2$.

Remember that for a $\pi$-instanton $\mathbf{A}$ (these chosen arbitrarily), its analtyic energy $\mathcal E^{\text{an}}_\pi(\mathbf{A})$ was defined to be equal to $$\mathcal E^{\text{top}}_\pi(\mathbf{A}) = 2\left((\text{cs}+f_{\pi_-})(\alpha_-) - (\text{cs}+f_{\pi_+})(\alpha_+)\right).$$ Cutting this off at the boundary of the ends (we say $\mathbf{A}$ restricted to $\{t\} \times Y_1$ is $\mathbf{A}_-(t)$, and restricted to $\{t\} \times Y_2$ is $\mathbf{A}_+(t)$), where the perturbations are all zero, this decomposes into three pieces:

\begin{align*}2\left((\text{cs}+f_{\pi_-})(\alpha_-) - (\text{cs})(\mathbf{A}_-(0))\right)\\
+ 2\left((\text{cs})(\mathbf{A}_-(0)) - (\text{cs})(\mathbf{A}_+(0))\right)\\
+ 2\left((\text{cs})(\mathbf{A}_+(0)) - (\text{cs}+f_{\pi_+})(\alpha_+)\right)\end{align*}

First we study how the value of $\text{cs}+f_{\pi(t)}$ changes along the ends: consider the function $(\text{cs}+f_{\pi(t)})(\mathbf{A}_-(t))$. Here, $\pi(t) = \beta_0(t) \pi_-$. Then 
$$\int d\left((\text{cs}+f_{\pi(t)})(\mathbf{A}_-(t))\right) = \text{cs}(\mathbf{A}_-(0)) - (\text{cs}+f_{\pi_-})(\alpha_-).$$
We may expand the integral as 
$$\int (\text{cs}+f_{\pi(t)})'(\mathbf{A}_-(t)) dt + \int d(\text{cs}+f_{\pi(t)})\left(\frac{d}{dt} \mathbf{A}_-(t)\right).$$ 
Note that the second integral is the same as 
$$\int \big\langle *F_{\mathbf{A}_-(t)} + \nabla_{\pi(t)}(\mathbf{A}_-(t)), \frac{d}{dt} \mathbf{A}_-(t)\big\rangle;$$ 
by assumption that $\mathbf{A}_-(t)$ satisfies the time-dependent gradient flow equations, this integral is $\int \|\frac{d}{dt} \mathbf{A}_-(t)\|^2 dt$. Because $f_{\pi(t)} = \beta_0(t) f_{\pi_0}$, the first integral is the same as 
$$\int \beta_0'(t) f_{\pi_-}(\mathbf{A}_-(t)).$$ 
This is uniformly bounded above and below because the support of $\beta_0'$ is compact, and the function $f_{\pi_-}$ is bounded. Thus up to an amount bounded by a continuous function of $\pi$, the first part of the Chern-Simons difference above is $2\int \|\frac{d}{dt} \mathbf{A}_-(t)\|^2$; a similar discussion implies the same of the the third part of the Chern-Simons difference and $2\int \|\frac{d}{dt} \mathbf{A}_+(t)\|^2$.

Because $\pi_n \to \pi$, we obtain a bound on the energy of any $\pi_n$-instanton $\mathbf{A}_n$ on the ends, uniform in $n$. 

Altogether we see that an energy bound on $\mathbf{A}$ gives a bound on the $L^2$ norm of the $F_{\mathbf{A}_n}$ on $W'$, as well as a bound on the $L^2$ norms of $\mathbf{A}_\pm(t)$. Because $\mathbf{A}_-(t)$ satisfies the equations 
$$\frac{d}{dt} \mathbf{A}_-(t) = -*F_{\mathbf{A}_-(t)} - \nabla_{\pi(t)}(\mathbf{A}_-(t)),$$ 
and $\|\nabla_{\pi(t)} A\|_{L^\infty}$ is uniformly bounded, we get a uniform $L^2$ bound on curvature on any compact piece of the negative end of $W$; similarly with the positive end. Because this bound holds for all $\mathbf{A}_n$ and is uniform in $n$, as in Proposition \ref{properflowlines} we may conclude the existence of the appropriate gauge transformations. 

What remains is to check is that the energy lost at the points is as expected.

Now suppose $\mathbf{A}_n$ is a sequence of $\pi_n$-perturbed instantons on $W$ and suppose there is a sequence of gauge transformations $\sigma_n^*$, defined on the complement of a finite set of points, so that $\sigma_n(\mathbf{A}_n) \to \mathbf{A}$ in $L^p_1$ on compact sets in the complement of those points, for $p \geq 2$. Call one of these points $x$; we will show the energy lost at $x$ is a non-negative multiple of $64\pi^2$; call this energy $c$. Abuse notation and rewrite $\mathbf{A}_n = \sigma_n(\mathbf{A}_n)$ for convenience, write $S(r)$ for the geodesic sphere of radius $r$ around $x$, and write $\mathbf{A}_n(r)$ for its restriction to $S(r)$. 

Because each $\mathbf{A}_n(r) \to \mathbf{A}(r)$ in the $L^2_{1/2}$ topology, we see that 
$$\text{lim}_n \text{cs}(\mathbf{A}_n(r)) \equiv \text{cs}(\mathbf{A}(r)) \bmod 64\pi^2 \mathbb Z,$$ 
We also have, by definition, 
$$\text{cs}(\mathbf{A}_n(r)) = \int_{B(r)} \text{Tr}\left(F_{\mathbf{A}_n}^2\right);$$ 
because $\mathbf{A}_n$ is a $\pi$-instanton, we have 
$$\text{cs}(\mathbf{A}_n(r)) = \int_{B(r)} \|F_{\mathbf{A}_n}\|^2 -2\int_{B(r)} \|\widehat{\nabla}_{\pi_n}(\mathbf{A}_n)\|.$$ 
Because $\|\widehat{\nabla}_\pi(\mathbf{A})\|_{L^\infty} \leq C\|\pi\|$, the second term goes to zero as $r \to 0$, and the first term limits to $c$, the energy lost at $x$. 

Putting this together, we find 
$$c = \lim_{r \to 0} \lim_{n \to \infty} \text{cs}(\mathbf{A}_n(r)) \equiv \lim_{r \to 0} \text{cs}(\mathbf{A}(r)) = 0 \bmod 64\pi^2 \mathbb Z,$$
as desired. 
\end{proof}

\begin{corollary}The space $\overline{\mathcal M}^{\leq C}_{\mathbf{E},\pi}(\alpha_-, \alpha_+)$ of broken ideal trajectories on $(W,\mathbf{E})$ with uniform energy bound is compact.
\end{corollary}

The following theorem summarizes the content of this section. 

\begin{theorem}\label{comp}Let $(Y,E)$ be a 3-manifold equipped with an $SO(3)$-bundle and regular perturbation $\pi$. There is a natural compactification of the space of unparameterized trajectories $$\widetilde{\mathcal M}^0_{E,z,\pi}(\alpha,\beta) \subset \overline{\mathcal M}_{E,z,\pi}(\alpha,\beta)$$ as a compact $SO(3)$-space equipped with equivariant endpoint maps.\footnote{Here recall that $\widetilde{\mathcal M}^0$ is $\widetilde{\mathcal M}/\mathbb R$, less the constant trajectories.} The added strata $\overline{\mathcal M} \setminus \widetilde{\mathcal M}$ are given as a union of fiber products of moduli spaces of lower dimension (in codimension equal to the number of intermediary orbits between $\alpha$ and $\beta$) and strata corresponding to Uhlenbeck bubbling.

The same is true for a cobordism $(W,\mathbf{E})$ equipped with a perturbation $\pi$ which is regular on the bounding manifolds $(Y_i, E_i)$: the space of framed instantons $\widetilde{\mathcal M}_{\mathbf{E},z,\pi}(\alpha, \beta)$ has a natural compactification to the compact $SO(3)$-space $\overline{\mathcal M}_{\mathbf{E},z,\pi}(\alpha, \beta)$.
\end{theorem}

We conclude with some discussions on energy bounds. On a cylinder $\mathbb R \times Y$, equipped with a regular perturbation $\pi$, there are finitely many critical orbits $\alpha$. For a fixed pair $\alpha, \beta$, the space of connections from $\alpha$ to $\beta$ has
$$\pi_1 \widetilde B^e_E(\alpha, \beta) \cong \mathbb Z$$
components. We may define two functions of $z \in \pi_1 \widetilde{\mathcal B}^e_E(\alpha, \beta)$. The first is the topological energy $\mathcal E^\pi(z)$, defined by taking a connection $A$ in the component labelled by $z$, and evaluating $$2\left((\text{cs}+f_\pi)(\beta)-(\text{cs}+f_\pi)(\alpha)\right) \in \mathbb R;$$ we use the path $A$ and the homotopy lifting property of the covering $\mathbb R \to \mathbb R/8\pi^2$ to pin down a real lift of this difference. The second is $\text{gr}_z(\alpha, \beta) \in \mathbb Z$. These are both affine functions; writing $z \in \mathbb Z$, we have $\mathcal E^\pi(z) = 64\pi^2 z + c$ and $\text{gr}_z(\alpha, \beta) = 8z + c'$. In particular, we see that
$$\mathcal E^\pi(z) - 8\pi^2\text{gr}_z(\alpha, \beta)$$
is constant. Therefore an energy bound is equivalent to an upper bound on $\text{gr}_z(\alpha, \beta)$.

Suppose that the cobordism
$$(W,\mathbf{E}): (Y_1, E_1) \to (Y_2, E_2)$$
is equipped with a perturbation $\pi$ so that its restriction $\pi_i$ to the ends is regular, and all $\pi$-perturbed instanton moduli spaces up to some energy bound $C$ are cut out transversely. Unfortunately, we do not have a notion of perturbed topological energy for a cobordism, so we must be slightly more precise in our phrasing.

The set $\pi_0\left(\widetilde{\mathcal B}^e_{\mathbf{E}}(\alpha, \beta)\right)$ of components of the space of trajectories from $\alpha$ to $\beta$ is affine over both $\pi_1(\widetilde{\mathcal B}^e_{E_1}, \alpha)$ and $\pi_1(\widetilde{\mathcal B}^e_{E_2},\beta)$, and hence we may write the set of components as affine over $\mathbb Z$. Then $\text{gr}_z(\alpha, \beta)$ is affine with distortion $8$. Now, on this same set, we have a \emph{partially defined} function to $\mathbb R$, given by the analyic energy
$$\mathcal E^\pi_{\text{an}}(\mathbf{A}) = \int_W \|F_{\mathbf{A}} + \widehat{\nabla}_\pi(\mathbf{A})\|^2$$
for a $\pi$-instanton $\mathbf{A}$ in the homotopy class $z$. Of course, this is only defined when such $\pi$-instantons exist, but it is independent of the choice when they do. Furthermore, when defined, this function is affine with distortion $64\pi^2$. This is enough to say that an energy bound $\mathcal E^\pi_{\text{an}}(\mathbf{A}) \leq C$ is equivalent to a grading bound $\text{gr}_z(\alpha, \beta) \leq C'$ for $\pi$-instantons $\mathbf{A}$.

In what follows we will often only care about the case $\text{gr}_z(\alpha, \beta) \leq 10 - \dim \alpha$; there is a corresponding energy bound for $\pi$-instantons, for $\pi$ fixed (or for $\pi$ chosen from some compact set).

Because components of the moduli space with $\text{gr}_z(\alpha, \beta) < 3 - \dim \alpha$ are empty when cut out transversely, and bubbling decreases $\text{gr}_z(\alpha,\beta)$ by multiples of $8$, we have the following.

\begin{corollary}\label{nobubbling}When $\pi$ is a regular perturbation and $\textup{gr}_z(\alpha, \beta) \leq 10 - \dim \alpha$, the compactification $\overline{\mathcal M}_{\mathbf{E},z,\pi}(\alpha,\beta)$ has no strata corresponding to Uhlenbeck bubbling.
\end{corollary}

\section{Reducible instantons on the cylinder and cobordisms}\label{sec:4d-red}
The goal of this section is to characterize the reducible $\pi$-ASD connections on a cobordism $W$. We begin with the especially simple case of the cylinder.

\begin{proposition}\label{red-flowlines}Let $E$ be an $SO(3)$-bundle over a rational homology sphere $Y$ equipped with a Riemannian metric. Let $\mathcal P_E$ denote a Banach space of perturbations on $(Y,E)$. By Proposition \ref{redcrit}, if $\pi \in \mathcal P_E$ has $\|\pi\|$ sufficiently small, then there is a unique critical orbit $\alpha$ in every reducible component $\textup{Red}(Y,E)$, which is the unique fully reducible orbit if the reducible component has a fully reducible point.

Consider moduli spaces of reducible instantons on the cylinder $\mathbb R \times Y$, equipped with the constant perturbation $\pi$ on the cylinder. By Proposition \ref{action2}, the only moduli spaces which can possibly be nonempty are $\widetilde{\mathcal M}_{E,0,\pi}(\alpha, \alpha)$, the spaces of trajectories from $\alpha$ to $\alpha$ in the trivial homotopy class. There is exactly one instanton in each of these: the constant trajectory at $\alpha$. Therefore, the moduli space $\overline{\mathcal M}_{E,0,\pi}(\alpha, \alpha) = \varnothing$, by definition.
\end{proposition}

\begin{proof}In other words, the constant trajectory (which is always a solution) remains the only solution. Because each component in $\textup{Red}(Y,E)$ lies inside $\widetilde{\mathcal B}^e_{E,0,k,\delta}(\alpha, \alpha)$ (the zero denoting the trivial homotopy class), and an instanton in this component has energy equal to zero, the trajectory must be constant, as in Corollary \ref{triv-component}.
\end{proof}

In fact, even in the case of cobordisms (or cylinders with nonconstant perturbations), the enumeration will be in terms of topological information.

\begin{proposition}\label{red-solutions-nobplus}Let $W$ be an oriented Riemannian 4-manifold equipped with an $SO(3)$-bundle $\mathbf{E}$ and with one incoming cylindrical end $(Y_1,E_1)$ and one outgoing cylindrical end $(Y_2, E_2)$. Suppose $b^+(W) = 0$.

Consider the space of perturbations on $W$ so that the perturbations $\pi_\pm$ on the ends lie in $\mathcal P^{\textup{red}}_{E,\delta}$, the subspace of perturbations so that the only eigenvalues of $\widehat{\textup{Hess}}^\nu_{\alpha, \pi_\pm}$ acting on $\Omega^1(Y;i\mathbb R)$ have absolute value at least $\delta$; abusing notation\footnote{This space is larger than the usual $\mathcal P^{(4)}_\delta$ defined in Definition \ref{pertdef}, as the requirement is easier to satisfy; as the results will apply to the usual subset $\mathcal P^{(4)}_\delta$, and this larger set will not be used outside this section, we do not feel this abuse of notation is harmful.}, we still call this space $\mathcal P^{(4)}_\delta$.

If $\beta(w_2(\mathbf{E})) \neq 0$ or one of $(Y_i, E_i)$ is admissible, then $(W,\mathbf{E})$ admits no $\pi$-perturbed reducible instantons in $\widetilde{\mathcal M}_{\mathbf{E},k,\delta}$ whatsoever.

If $\mathbf{E}$ is trivial and $b_1(W) = 0$, there is a neighborhood $\mathcal P^{(4)}_{\textup{end}} \subset \in U \subset \mathcal P^{(4)}$ so that, for all $\pi \in U$, all fully reducible connections are cut out transversely. We write $U_\delta = U \cap \mathcal P^{(4)}_\delta$; this set contains $0$ so long as $\delta$ is sufficiently small. 

Fix $(\pi_-, \pi_+) \in \mathcal P^{(4)}_{\textup{end}}$; write $U_{\delta, c} = U_\delta \cap \mathcal P^{(4)}_c$ for the portion contained inside the affine space through $(\pi_-, \pi_+)$ over $\mathcal P^{(4)}_{\text{int}}$. 

For all $\pi$ in a residual subset of $U_{\delta,c}$ and for every component in $\textup{Red}(W,\mathbf{E})$, the reducible $\pi$-instantons in this component comprise a finite set of orbits, which are cut out transversely inside the reducible locus. 

If $\pi\in \mathcal P^{(4)}$ is sufficiently small (with respect to some constant $\epsilon_{W,C}$), then for every component in $\textup{Red}(W,\mathbf{E})$ with $\mathcal E^\pi(\mathbf{A}) \leq C$, there is further a unique orbit of $\pi$-perturbed reducible instantons in that component; for those components $\textup{Red}_*(W,\mathbf{E})$ containing a fully reducible orbit, that orbit is the unique reducible instanton in that component.
\end{proposition}

For bundles with $\beta w_2 \neq 0$ or with admissible ends, Proposition \ref{action3} shows that there are not even reducible \emph{connections} with $L^2$ curvature. Thus there is only something interesting to say for cobordisms with rational homology sphere ends.

There is also a version when $b^+(W) > 0$. Our goal is to avoid reducibles in this case, as they cannot be cut out transversely.

\begin{proposition}\label{red-bplus-pos}Let $W$ be an oriented Riemannian 4-manifold equipped with an $SO(3)$-bundle $\mathbf{E}$, with one incoming cylindrical end $(Y_1,E_1)$ and one outgoing cylindrical end $(Y_2, E_2)$. Suppose $b^+(W) > 0$.

If $\beta(w_2(\mathbf{E})) \neq 0$ or one of $(Y_i, E_i)$ is admissible, then $(W,\mathbf{E})$ admits no $\pi$-perturbed reducible instantons in $\widetilde{\mathcal M}_{\mathbf{E},\pi}$ whatsoever.

If $b_1(W) < b^+(W)$, and $\mathbf{E}$ is nontrivial but neither $E_i$ is admissible, for any fixed $\pi_\pm \in \mathcal P^{(4)}_{\textup{end},\delta}$, there exists a residual set of the affine slice $\pi \in \mathcal P^{(4)}_{c}$ for which there are no reducible $\pi$-ASD connections in any component $z$.

If $\mathbf{E}$ is trivial, each fully reducible component always has the full reducible as a solution, no matter the perturbation; if $b^+(W) > b_1(W)$, these fully reducible solutions are never cut out transversely in the $SO(2)$-fixed locus, for any $\pi \in \mathcal P^{(4)}_{\delta}$.
\end{proposition}

We prove these simultaneously, much like Proposition \ref{redcrit}: we define a Banach manifold of reducible solutions to the ASD equations, equipped with a map to $U_\delta$ (or $\mathcal P^{(4)}_\delta$, in the case that there are no fully reducible connections), which is proper; when $b_1(W) = b^+(W) = 0$, the projection is a local diffeomorphism above points near $0$. Recall from Proposition \ref{action2} that the $SO(2)$-fixed subspace is a disjoint union over copies of $\mathcal B_{\eta}(\alpha_-, \alpha_+)$, as $\eta$ varies over certain complex line bundles. Note that this configuration space is unframed, despite being a subset of the framed configuration space of $SO(3)$-connections on $E$. Setting this up takes some small amount of work because the Hilbert manifold the equation is \emph{defined} on, $\mathcal{B}_{\eta}(\alpha_-, \alpha_+)$, depends on the limiting orbits, and hence depends on $\pi$ (which determines the reducible critical orbits). While we discussed such a configuration space in the previous section, we did not show that this (unframed) quotient was a Banach manifold: that is usually not true, but is in the special case of $SO(2)$-bundles.

\begin{lemma}\label{pi-Banach}Let $W$ be an oriented Riemannian 4-manifold equipped with an $SO(3)$-bundle $\mathbf{E}$ and with one incoming cylindrical end $(Y_1,E_1)$ and one outgoing cylindrical end $(Y_2, E_2)$, both rational homology spheres. Fix a lift of $\mathbf{E}$ to a $U(2)$-bundle $\widetilde{\mathbf{E}}$; if none exist, then $W$ admits no reducible connections whatsoever by Proposition \ref{action3}. Write $\lambda$ for the complex line bundle $\det(\widetilde{\mathbf{E}})$ and fix a connection $\mathbf{A}_0$ on $\lambda$.

There is a Banach manifold $\mathcal P_{\delta} \mathcal B_{\eta, k, \delta}$, whose objects are pairs $(\pi, \mathbf{A})$ defined as follows. First, $\pi \in \mathcal P^{(4)}_\delta$ is a perturbation as defined above; second, $\mathbf{A}$ is an $L^2_{k,\delta}$ connection on $\widetilde{\mathbf{E}}$ which respects a fixed splitting
$$\widetilde{\mathbf{E}} \cong \eta \oplus (\lambda \otimes \eta^{-1})$$
and has $\pi_\pm$-flat limits on the ends, considered up to equivalence by the group of $L^2_{k+1,\delta,\textup{ext}}$ maps $W \to S^1$, meaning that these maps exponentially decay to constant maps on the ends.

This Banach manifold comes equipped with a smooth submersion $\mathcal P_{\delta} \mathcal B_{\eta, k, \delta} \to \mathcal P^{(4)}_{\delta}$. There is a smooth vector bundle $\mathcal S_{k-1,\delta} \to \mathcal P_{\delta} \mathcal B_{\eta, k, \delta}$ with fiber isomorphic to $\Omega^{2,+}_{k-1,\delta}(i\mathbb R)$, and a smooth Fredholm section $s(\pi, \mathbf{A}) = F_{\mathbf{A}}^+ + \widehat{\nabla}_\pi(\mathbf{A})$.
\end{lemma}

\begin{proof}Denote by $\eta_i$ the restrictions of $\eta$ to the ends. Write $\mathcal P_\delta \mathfrak C_\eta$ for the set whose elements are triples $(\pi, \alpha_1, \alpha_2)$, where $\pi \in \mathcal P^{(4)}_\delta$ is a perturbation restricting to $\pi_i$ on the ends, and $\alpha_i \in \mathcal B_{\eta_i}$ are equivalence classes of $\pi_i$-flat connections on $\eta_i$. With minimal change, Lemma \ref{prop1} asserts that the map $\mathcal P_{\delta} \mathfrak C_\eta \to \mathcal P^{(4)}_{\delta}$, sending the parameterized critical set of gauge equivalence classes of connections on $\eta$ to the space of perturbations, is proper. That $\mathcal P_{\delta} \mathfrak C_\eta$ is a manifold and that the projection above is a covering map follow quickly from the assumption that every $\pi \in \mathcal P^{(4)}_{\delta}$ is such that $\widehat{\textup{Hess}}^\nu_{A_i, \pi_i}$ never has eigenvalues of absolute value at most $\delta$. Thus, $\mathcal P_{\delta} \mathfrak C_\eta$ has charts given by local sections $\psi$ of the projection to $\mathcal P^{(4)}_{\delta}$. Fix such a chart, given by $V \subset \mathcal P^{(4)}_{\delta}$.

We may lift the above section $\psi: V \to \mathcal P_{\delta} \mathfrak C_\eta$ to a smooth map
$$p: V \to \mathcal A_{\eta_1, k} \times \mathcal A_{\eta_2, k},$$
so that $p(\pi_1, \cdot, \cdot, \pi_2)$ is a pair of connections in the same gauge equivalence as $\psi(\pi_1, \cdot, \cdot, \pi_2)$. That is, this is a smooth map picking out critical points for the end-perturbations $\pi_i$. This entire discussion depends only on the connections over the ends.

For this chart $V$, there are two possible situations: it could be that there are no reducible connections over the cobordism restricting to the configuration $\psi(\pi_1, 0, 0, \pi_2)$ on the ends, in which case $\mathcal P_{\delta} \mathcal B_{\eta, k, \delta}$ is empty above this chart and is tautologically a manifold.

More interestingly, if there is some reducible connection extending that configuration, then may choose a smooth map $e: \mathcal A_{\eta_1, k} \times \mathcal A_{\eta_2, k} \to \mathcal A^{(4)}_{\eta, k,\delta}$ so that $e(A_1, A_2)$ is a connection $\mathbf{A}$ on $\eta$ so that $\mathbf{A}$ is constant and equal to $A_i$ on the ends. Composing these, $e \circ p: \mathcal P^{(4)}_{\delta} \to \mathcal A^{(4)}_{\eta, k, \delta}$ is a smooth map choosing a connection $(e\circ p)(\pi)$ which is constant and equal to the $\pi_\pm$-critical connections $p(\pi)$ on the ends.

We may use this to define $$V \mathcal A_{\eta,k,\delta}^{(4)} = \{(\pi, \mathbf{A}) \mid \pi \in V,  \mathbf{A} - (e\circ p)(\pi) \in \Omega^1_{k,\delta}(W;i\mathbb R)\}.$$ In particular, this Banach manifold is diffeomorphic to $V \times \Omega^1_{k,\delta}(W;i\mathbb R).$ Clearly it comes equipped with a smooth projection to $V$. There is the trivial bundle $\mathcal S_{k-1,\delta}^+ = \Omega^{2,+}_{k-1,\delta}(W;i\mathbb R)$ over $V \times \mathcal A_{\eta,k,\delta}^{(4)}.$ Given a connection $\mathbf{A}$ on $\eta$, the induced connection on $\eta \oplus (\lambda \otimes \eta^{-1})$ has curvature equal to $2F_{\mathbf{A}} - F_{\mathbf{A_0}}$, the non-central part of the curvature of $\eta \oplus (\lambda \otimes \eta^{-1})$. Then the section $s: V \mathcal A^{(4)}_{\eta, k, \delta} \to \mathcal S_{k-1, \delta}$ is given by $$(\pi, \mathbf{A}) \mapsto 2F_{\mathbf{A}}^+ - F_{\mathbf{A}_0}^+ + \widehat{\nabla}_\pi(\mathbf{A}).$$

This space a smooth action by the Banach Lie group $\mathcal G_{\eta, k+1,\delta,\text{ext}}$, preserving the projection map to $\mathcal P^{(4)}_{\delta}$ and the section $s$. While the action is not free, the stabilizer is the same at every point: it is the group of constant maps $W \to S^1$. The action of $\mathcal G_{\eta,k+1,\delta}$ factors through its quotient by the subgroup of constant maps, and the action of \emph{this} quotient group is free. The fact that `$\mathbf{A}$-harmonic gauge transformations' are the same for all $\mathbf{A}$ is a convenient and unique aspect of the case of $SO(2)$-bundles.

Now it is easy to verify that this action is proper, such that each orbit has closed complemented tangent space, giving us a quotient manifold $V \mathcal B_{\eta, k, \delta}$ with a projection map and a section of a vector bundle $\mathcal S_{k-1,\delta}$ (still the trivial bundle). That the section on $V \mathcal B_{\eta, k, \delta}$ is Fredholm follows from the same fact fiberwise, on a fixed Banach manifold $\mathcal B_{\eta, k, \delta}(\alpha_1, \alpha_2)$. Patching these together for different charts $V$ of $\mathcal P_{\delta} \mathfrak C_\eta$ gives us a smooth Banach manifold structure on $\mathcal P_{\delta} \mathcal B_{\eta, k, \delta}$.
\end{proof}

Note that this Banach manifold is precisely the same as a component of the $SO(2)$-fixed point space of $\mathcal P_{\delta} \widetilde{\mathcal B}_{\mathbf{E},k,\delta}$, expressed in terms of the $U(2)$ model.

In the case of reducible trajectories, the compactness properties are especially strong. Note that there are no broken or ideal trajectories in the following statement.

\begin{lemma}\label{C-proper}The map $\mathcal P_\delta \mathcal M_{k,\delta} \to \mathcal P^{(4)}_{\delta}$ is proper.
\end{lemma}

\begin{proof}This is an application of Proposition \ref{propertrajectories}. Start by choosing a sequence $(\pi_n, \mathbf{A}_n)$ of perturbed reducible instantons so that $\pi_n \to \pi$.

First, observe that because every fiber of $\mathcal P_\delta \mathcal B_{\eta, k, \delta} \to \mathcal P^{(4)}_\delta$ is \emph{connected}, the analytic energy $\mathcal E^\pi(\mathbf{A})$, for a reducible $\pi$-instanton $\mathbf{A}$, depends only on $\pi$, and depends continuously on $\pi$ at that. In particular, we have $\mathcal E^{\pi}(\mathbf{A}_n) \leq C$ for some constant $C$ (without already assuming an energy bound on the $\mathbf{A}_n$).

In particular, by the compactness theorem there is a subsequence of $\mathbf{A}_n$ and a broken ideal trajectory $\mathbf{A}$ such that a subsequence of $\mathbf{A}_n$ conerges to $\mathbf{A}$ as a broken bubble-limit. Because all of the $\mathbf{A}_n$ are reducible, so is $\mathbf{A}$. Bubbling cannot change the topological type of $\eta$ (a line bundle over a manifold is determined by its restriction to the complement of any finite set, so long as the manifold is dimension at least $2$), and so for large $n$ we have that $\mathbf{A}_n$ lies in the same component as $\mathbf{A}$. Thus, we know that $\mathcal E^{\pi_n}(\mathbf{A}_n) \to \mathcal E^{\pi}(\mathbf{A})$, implying no energy is lost at bubble-points.

So the only noncompactness can arise from breaking of trajectories. We already know by Proposition \ref{red-flowlines} that on the cylinder $\mathbb R \times Y$, where $Y$ is equipped with a weakly admissible bundle and regular perturbation, reducible trajectories are constant. So after passing to a subsequence and applying suitable gauge transformations, the configuration $\mathbf{A}_n$ converges in $L^2_{k,\delta}$ an honest reducible connection on $\eta$ with the same limits at $\pm \infty$, as desired, with no energy lost at the ends.
\end{proof}

The following unique continuation lemma is quoted essentially verbatim from \cite[Lemma~7.1.3]{KMSW}; we will use it frequently.

\begin{lemma}\label{Unique}Let $H$ be a real Hilbert space and $I = [t_1, t_2]$ a closed interval, equipped with a family $L(t)$ of unbounded operators $D \to H$ with common dense domain $D\subset H$, so that $L(t) = L_+(t) + L_-(t)$ for $L_+$ symmetric on $D$ and $L_-(t)$ skew-adjoint and bounded on the whole of $H$. Further suppose that the time-derivative $L'(t)$ is a well-defined operator $D \to H$, defined pointwise, which has a bound $$\|L'(t)x\| \leq C_1(\|L(t)x\| + \|x\|)$$ uniform in $t\in I$ and $x \in D$.

Let $f: I \to H$ be a continuous map and $z: I \to D$ be a solution of the equation $z' + Lz = f$, where we have the bound on the inhomogeneous term $\|f(t)\| \leq C_2\|z(t)\|$, uniform in $t$. Then if $z(t) = 0$ for some $t \in I$, then $z$ is identically zero.
\end{lemma}

We will first justify the open set $U_\delta$ used in the statement of Proposition \ref{red-solutions-nobplus}.

\begin{lemma}\label{Trivial-transverse}Suppose $b_1(W) = 0$. If $b^+(W) = 0$ or $\mathbf{E}$ is nontrivial, then all fully reducible connections are cut out transversely for any end perturbation; therefore, there is a connected open set 
$$\mathcal P^{(4)}_{\textup{end}} \subset U\subset \mathcal P^{(4)}$$
for which every fully reducible connection is cut out transversely. When $\delta > 0$ is sufficiently small, the intersection $U_\delta = U \cap \mathcal P^{(4)}_\delta$ contains $0$.
\end{lemma}

\begin{proof}If $\mathbf{E}$ is nontrivial, there are no fully reducible connections; if $b^+(W) = 0$, then for any fully reducible connection $\mathbf{A}$, we have
$$\text{coker}(Q^\nu_{\mathbf{A},\pi}) = \mathcal H^+_{\mathbf{A}} = \mathcal H^+ = 0$$
by the Hodge theorem. Therefore every fully reducible connection is cut out transversely above $0$; furthermore, because we assume the 3-dimensional holonomy perturbations $f_\pi$ vanish in a neighborhood of fully reducible connections, we necessarily have $D_{A}\nabla_\pi = 0$ for any fully reducible connection $A$. 
In particular, because the end perturbations take the form 
$$\widehat{\nabla}_\pi \mathbf{A} = \beta_0(t) \big(dt \wedge \nabla_\pi \mathbf{A}(t)\big),$$
we see that $D_{\mathbf{A}} \widehat{\nabla}_\pi = 0$ for any full reducible $\mathbf{A}$ on $W$ and any end perturbation $\pi$; so $Q_{\mathbf{A},\pi}$ is surjective iff $Q_{\mathbf{A},0}$ is.

Because the projection from the parameterized configuration space of fully reducible connections to $\mathcal P^{(4)}_\delta$ is a proper local diffeomorphism, we see that every fully reducible connection is cut out transversely in a neighborhood of $\mathcal P^{(4)}_{\text{end}}$, and we call this neighborhood $U$.

To say that $0 \in \mathcal P^{(4)}_\delta$ just means that for the reducible flat connections $\alpha$, all eigenvalues of $\widehat{\text{Hess}}^\nu_{\alpha}$, acting on $\Omega^1(W;i\mathbb R)$, have absolute value larger than $\delta$. Because the kernel of this operator is $\mathcal H^1(W) = 0$, and the eigenvalues form a discrete closed set, we see that for the finitely many reducible flat connections $\alpha$, all eigenvalues have absolute value larger than some $\delta > 0$, as desired.
\end{proof}

When $b_1(W) = 0$ and $b^+(W) > 0$ but $\mathbf{E}$ is trivial, the set $U_{\delta}$ is empty: full reducibles are always $\pi$-perturbed ASD connections, regardless of $\pi$, but the index of $D_{\mathbf{A},\pi}$ is always negative, so they cannot be cut out transversely.

\begin{lemma}\label{zero-Banach}Let $(W,\mathbf{E})$ be a Riemannian manifold with two cylindrical ends, both of which are modelled on rational homology spheres, with $b_1(W) = 0$. Suppose either that $b^+(W) = 0$ or that $\mathbf{E}$ is nontrivial, and fix a component of the $SO(2)$-fixed point space, labeled by a line bundle $\eta$. We write $U_\delta \subset \mathcal P^{(4)}_\delta$ for the open set in which all fully reducible connections are cut out transversely; in the latter case, when $\mathbf{E}$ is non-trivial, this containment is an equality.

Then the zero set $$U_\delta \mathcal M_{\eta} \subset U_\delta \mathcal B_{\eta,k,\delta}$$ is a smooth Banach submanifold.

Similarly, for any fixed $(\pi_-, \pi_+) \in \mathcal P^{(4)}_{\textup{end},\delta}$, the zero set $U_{\delta, c} \mathcal M_\eta$ is a Banach submanifold of $U_\delta \mathcal B_{\eta, k, \delta}$. 
\end{lemma}

\begin{proof}By assumption, we need only consider the $SO(2)$-reducible ASD connections which are not fully reducible.

We want to show that
$$Q'_{\mathbf{A},\pi} = \big(Q_{\mathbf{A},\pi}, \widehat{\nabla}_{\pi'}(\mathbf{A})\big): \Omega^1_{k,\delta}(W;i\mathbb R) \oplus \mathcal P^{(4)} \to \Omega^{2,+}_{k-1,\delta}(W;i\mathbb R)$$ is surjective for all reducible $\pi$-instantons $\mathbf{A}$ in the component corresponding to $\eta$ (of class $L^2_{k,\delta}$).

To see this, consider $\psi \in \text{coker}(Q'_{\mathbf{A},\pi})$; necessarily, $\psi \in \text{coker}(Q_{\mathbf{A},\pi})$, and in particular $Q^*_{\mathbf{A},\pi}\psi = 0$. We may write this as

$$Q^*_{\mathbf{A},0}\psi = -(D_{\mathbf{A}}\widehat{\nabla}_\pi)^*\psi;$$ because $D_{\mathbf{A}}\widehat{\nabla}_\pi$ extends to $L^2_{j,\delta'}$ for $j \le k$ and $\delta' \le \delta$, the usual elliptic bootstrapping arguments imply that $\psi \in L^2_{k,\delta}$ by Proposition \ref{pertfacts} (1).

Recall that the interior holonomy perturbations are supported inside a compact submanifold $W' \subset W$, whose complement is the portion $[1,\infty) \times Y$ of the cylindrical ends. Now, the argument in \cite[Lemma~3.7]{Kronheimer} shows that so long as $\pi$ is not fully reducible,\footnote{If $\pi$ is fully reducible, these perturbations are identically zero. The hypothesis of $SO(2)$-reducibility is used to ensure that, at any point in $W'$, one may find a collection of loops so that $\text{Hol}_{\mathbf{A}} \not \in (\pm 1)^N$; at these points of $U(1)^N \subset SU(2)^N$, any equivariant map $r: SU(2)^N \to \mathfrak{su}(2)$ must necessarily send $r((\pm 1)^N) = 0$; points in $U(1)^N \setminus (\pm 1)^N$ are sent into $\mathfrak u(1) \subset \mathfrak{su}(2)$, but are not otherwise constrained.} the image of $\widehat{\nabla}_{\pi'}(\mathbf{A})$ is dense in continuous sections of $\Lambda^{2,+}(i\mathbb R)$ over $W'$ which vanish on the boundary; in particular, because $\psi \in \text{coker}(Q')$, it should be $L^2$ orthogonal to all of these, and hence vanish on the whole of $W'$.

Now $\psi$ satisfies the perturbed ASD equations on the end, written $[0, \infty) \times Y$, with $\psi \in L^2_k$ and $\psi\big|_{\{0\} \times Y} = 0$. So long as $k \geq 2$, the map $\psi$ defines an element of
$$C^1\left([0, \infty), \Omega^1_0(Y)\right) \cap C^0\left([0, \infty), \Omega^1_1(Y)\right);$$
here the subscripts $k,\delta$ indicate we are taking 1-forms of Sobolev class $L^2_{k}$. Further, the equation $Q_{\mathbf A, \pi}^* \psi = 0$ can be expressed as the pair of ODE
\begin{align*}
d^*\psi(t) &= 0 \\
\psi'(t) - *d\psi(t) &= D_{\mathbf{A}(t)} \nabla_{\pi(t)} \psi(t).
\end{align*}

This ODE satisfies the conditions of Lemma \ref{Unique} so long as $\pi(t)$ and $\mathbf{A}(t)$ are $C^1$ paths; the fact that $\mathbf{A}(t)$ is a $C^1$ path in $\Omega^1_0(Y)$ follows because $\mathbf{A}$ is $L^2_2$, and $\pi(t) = \beta_0(t) \pi_\infty$ is smooth.

Therefore $\psi(0) = 0$ implies that $\psi = 0$ on this entire end. Repeating for both ends, we find that any $\psi \in \text{coker}(Q'_{\mathbf{A},\pi})$ is necessarily zero, as desired.
\end{proof}

\begin{proof}[Proof of Propositions \ref{red-solutions-nobplus} and \ref{red-bplus-pos}]
Until the final paragraph of this proof, we assume that $b_1(W) = 0$, and that $\mathbf{E}$ is nontrivial or $b^+(W) = 0$.

Fix $(\pi_-, \pi_+) \in \mathcal P^{(4)}_{\text{end}, \delta}$.

The set $U_{\delta,c} \mathcal M_{\eta,k,\delta}$ is a smooth submanifold of $U_{\delta,c} \mathcal B_{\eta, k, \delta}$. For fixed $\pi$, the set of all $\pi$-instantons $\mathbf{A}$ are cut out transversely for all in $\mathcal B_{\eta, k, \delta}$ if and only if the projection
$$p_\eta: U_{\delta,c} \mathcal M_{\eta} \to U_{\delta,c}$$
has $\pi$ as a regular value. By the Sard-Smale theorem, such $\pi$ are in large supply: they form a dense set. Furthermore, the projection $p_\eta$ is proper by Lemma \ref{C-proper}. Therefore, the regular values form an open set of $U_{\delta,c}$.

In the case that $b^+(W) > 0$, this means that $\mathcal M_{\eta, \pi} = \varnothing$ for $\pi$ a regular value, as the index of the operator $D_{\mathbf{A},\pi}$ is negative. In the case that  $b^+(W) = 0$, the moduli space $\mathcal M_{\eta, \pi}$ for regular $\pi$ is a finite set, and $0$ is a regular value.

What remains is to see what the solutions are when $\pi = 0$; then by properness, we immediately understand the solutions for all sufficiently small $\pi$. But for a cobordism $(W,\mathbf{E})$ with $b^+(W) = 0$ but $b_1(W) = n$, the space of reducible unperturbed ASD connections in each reducible component is the torus $T^n$; this is \cite[Lemma~1.6]{Dae}, and again follows from Hodge theory, as the ASD equation is affine in the $SO(2)$-reducible case. See also the related \cite[Lemma~2]{Fr1}. In the case we most care about, $b_1(W) = b^+(W) = 0$, this means that there is a unique $\pi$-perturbed reducible ASD connection in each reducible component, so long as $\pi$ is sufficiently small.

Intersecting over all line bundles $\eta$, these open dense sets become residual sets. For the statement of Proposition \ref{red-solutions-nobplus}, we want the enumeration above to hold for all line bundles up to a certain energy level; the assumption that $b^+(W) = 0$ implies that there are only finitely many line bundles whose topological energy is bounded above by a given constant $C$, and so intersecting over the finitely many relevant open sets, we find an open set containing $0$ so that the enumeration of reducibles up to energy $C$ is correct for all $\pi$ in this open set.

If $\mathbf{E}$ is trivial, the trivial connection $\mathbf{A}$ is always an ASD connection, no matter the perturbation. For $\pi \in \mathcal P^{(4)}_{\delta}$, we have $\text{ind}(D_{\mathbf{A},\pi}) = \text{ind}(D_{\mathbf{A},0})$. When $b^+(W) > b_1(W)$, this index is negative, and so $\mathbf{A}$ cannot be cut out transversely.
\end{proof}

\section{Index calculations}\label{sec:4d-index}
With a strong grasp on the regularity properties of reducible instantons (internal to the reducible locus), we move on to understanding the index of the ASD operator \emph{normal} to the reducible locus. 

Let $(W,\mathbf{E})$ be a cobordism from $(Y_1, E_1)$ to $(Y_2, E_2)$, thought of as a manifold with cylindrical ends, and let $\pi$ be a perturbation so that the ends $\pi_\pm$ are regular perturbations on $Y_i$. Suppose $\delta$ is larger than the least nonzero eigenvalue of the extended Hessian operators of the $\pi_\pm$-critical points. Choose $\pi_\pm$ sufficiently small that Proposition \ref{red1} holds; this is guaranteed by choosing $\pi \in \mathcal P^{(4)}_{\delta}$. 

If $\mathbf{A}$ is a $\pi$-perturbed instanton, we would like to compute the index of $$Q_{\mathbf{A},\pi}: \Omega^1_{k,\delta}(W;\mathfrak g_E) \to \Omega^{2,+}_{k-1,\delta} \oplus \Omega^0_{k-1,\delta}.$$ In fact, by invariance properties of the index, this only depends on the component $\mathbf{A}$ sits inside and on the value of $\pi$ on the ends. Correspondingly, we assume $\pi$ has no component corresponding to an interior perturbation: we suppose that $\pi$ is supported on the ends.

The most useful tool for computing the index of a differential operator on a compact manifold with boundary is the Atiyah-Patodi-Singer index theorem. Let $X$ be a compact manifold with boundary, with metric of product type near the boundary. Let $D$ be an elliptic differential operator acting on sections of bundles $V_1, V_2$ over $W'$ (thought of as a compact manifold with boundary) with specified isomorphisms near the boundary components $[0,\infty) \times Y$ from the given operator to $\frac{d}{dt} + A$, where $A$ is a self-adjoint elliptic operator whose \textit{nonzero} eigenvalues $c$ have $|c| > \delta$. Then one may consider $D$ as an operator 
$$D_{APS}: L^2_1(X,V_1;P) \to L^2(X,V_2),$$ 
where here $P$ means that we demand that \emph{the restriction of $\sigma \in L^2_1(W,V_1;P)$ to the boundary lies in the subspace spanned by the negative eigenvalues}: that is, if $P$ is the projection operator corresponding to the eigenspaces with $\lambda \geq 0$, we demand that $P\sigma = 0$. Then even if $A$, the operator at the boundary, has kernel, $D_{APS}$ is a Fredholm operator with a well-defined index. If 
$$\hat X := X \cup_{\partial X} [0,\infty) \times \partial X$$ 
with the appropriate product metric on the ends, we may just as well consider $D$ as an operator on $\hat X$. We write $I(D)$ to mean the index of $D$, computed as a map between weighted Sobolev spaces of weight $\delta > 0$, and $I_{APS}(D)$ to mean the index computed as $D_{APS}$.

\begin{lemma}As long as $\delta$ is less than any nonzero eigenvalue of the boundary operators $A_i$, we have $I_{APS}(D) = I(D)$.
\end{lemma}

\begin{proof}This is essentially \cite[Proposition~3.11]{atiyah1975spectral}. When thought of as a map on weighted Sobolev spaces, $\text{ker}(D)$ clearly consists of $L^2$ solutions of $Df$ on $\hat X$. Conversely, every $L^2$ solution has as much smoothness as the operator itself, and in the asymptotic expansion $\sum a_\lambda e^{-t\lambda} \phi_\lambda$ on the ends, $a_\lambda$ must be zero when $\lambda \leq 0$ if the solution is to be $L^2$. So a solution decays exponentially with weight at least $\delta > 0$, where $\delta$ is less than any positive eigenvalue of $A$.

The kernel of $D^*$ is computed in weighted Sobolev spaces as a subspace of $L^2_{k,-\delta}$. To decay exponentially slower than $e^{-\delta t}$, we must demand that $a_\lambda = 0$ for $\lambda \geq \delta$. Because positive eigenvalues are at least as large as $\delta$, we see that a solution is, on the ends, a sum of an exponentially decaying solution and a solution which is constant in time. This is precisely the kernel of $D^*$ on extended $L^2$ sections, as in the computation of $I(D_{APS})$ given in  \cite[Proposition~3.11]{atiyah1975spectral}. Because the two kernels and cokernels agree, we have $I_{APS}(D) = I(D)$.
\end{proof}

So it suffices to compute $I_{APS}(Q_{\mathbf{A},\pi})$, and for this, we have the Atiyah-Singer index theorem. We will need some simple lemmas first before we can calculate this index.

\begin{lemma}\label{APSadd}Suppose $W$ is a compact manifold with boundary of product type, equipped with an elliptic operator $D$ that is of product type near the boundary. Suppose further that there is an oriented closed submanifold $Y \subset W$ with a neighborhood of product type so that the operator may be written as $d/dt + A$ on this neighborhood. Write $D_1$ for the operator on the compact manifold whose positive boundary contains $Y$, and $D_2$ for the operator on the compact manifold whose negative boundary contains $Y$; say $A$ is the operator at $Y$. Then $$I_{APS}(D) = I_{APS}(D_1) + I_{APS}(D_2) + \dim \textup{ker}(A).$$
\end{lemma}
\begin{proof}This follows immediately from the index theorem itself, \cite[Theorem~3.10]{atiyah1975spectral}. The only term which is not additive is $-h/2$, where $h$ is the kernel of the boundary operator. Because we are contributing two extra copies of $-\dim \text{ker}(A)/2$ on the right, we counterbalance that by adding $\dim \text{ker}(A)$.
\end{proof}

Finally, we will need to know how the index of spectral flow is computed.

\begin{lemma}\label{sfindex}Let $D = d/dt + A_t$ be an operator on $[0,1] \times Y$, where $A_t$ is a time-dependent self-adjoint elliptic operator, possibly with kernel. We define the spectral flow $\textup{sf}(A_t)$ to be the intersection number of the graph of the spectra with the line $\lambda(A_t) = -\delta$, which $A_0$ and $A_1$ do not intersect. This is the aggregate number of eigenvalues that go from $\leq - \delta$ to $\geq 0$, counted with sign.

Then 
$$I_{APS}(D) = \textup{sf}(A_t) - \dim \textup{ker}(A_1).$$
\end{lemma}
\begin{proof}This may be proved using separation of variables. A similar formula is stated below \cite[Theorem~7.4]{atiyah1996spectral}, only giving the spectral flow term; their argument uses the periodic boundary conditions on $[0,1] \times Y$, which corresponds to the projection to nonnegative eigenvalues at $t = 0$ and positive eigenvalues at $t = 1$, whereas the APS boundary conditions stated above use the spectral projection to \emph{nonnegative} eigenvalues at $t=1$. The second operator has smaller domain, of codimension $\dim \text{ker}(A_1)$, having included the demand that an element $f$ of the domain projects nontrivially to $\text{ker}(A_1)$. Because index is addititive under composition and the inclusion of this subspace has index $-\dim\text{ker}(A_1)$, the theorem follows.
\end{proof}

Now $W$ may be decomposed as the union of a compact manifold $X$ and the two cylindrical ends. This decomposition will provide the desired calculation. Before stating the result, we recall the definition of one of the terms that will appear.

\begin{definition}Let $\alpha$ be an orthogonal (resp. unitary) flat connection on a real (complex) vector bundle $E$ of dimension $n$ over a closed manifold $Y$. The Atiyah-Patodi-Singer $\rho$-invariant of $\alpha$ is defined to be 
$$\eta_\alpha(0) - n \cdot \eta_\theta(0).$$ 
Here $\eta_\alpha$ is the Atiyah-Patodi-Singer $\eta$ invariant associated to the $\alpha$-twisted signature operator 
\begin{align*}
B: \Omega^{\textup{ev}}(W,E) &\to \Omega^{\textup{ev}}(W,E)\\
B\omega = (-1)^{\textup{deg } \omega}&\;(*d_\alpha - d_\alpha*)\omega,
\end{align*}
as is studied in \cite{atiyah1975spectral2}. This constant is denoted $\rho(\alpha)$.

If $\mathbf{E}$ has a reduction labelled by $\{\zeta_1, \zeta_2\}$, where $\zeta_i$ are complex line bundles and
$$c_1(\zeta_1) + c_1(\zeta_2) = \lambda,$$
we define $\rho(\{\zeta_1, \zeta_2\})$ to be the $\rho$-invariant of the induced $SO(3)$-flat connection on $\mathbf{E}$, given as $\mathbb R \oplus (\zeta_1 \otimes \zeta_2^{-1})$.
\end{definition}

The following equivalent computations of the $\rho$-invariant given here are listed in the discussion in \cite[Section~2.9]{hedden2011chern}. We state them without proof.

\begin{lemma}Suppose $E \to Y$ is an $SO(3)$-bundle, equipped with a reduction $E \cong \mathbb R \oplus (\zeta_1 \otimes \zeta_2^{-1})$. Write $\kappa = \zeta_1 \otimes \zeta_2^{-1}$, thought of as a flat complex line bundle.
If $z_i = c_1(\zeta_i)$, then $\rho(\{z_1, z_2\})$ coincides with $\rho(\kappa) + \rho(\kappa^{-1})$. Further, if 
$$(W,\mathbf{E}): (Y_1,E_1) \to (Y_2,E_2)$$ 
is a cobordism eqiupped with a compatible reduction  
$$\mathbf{E} \cong \mathbb R \oplus \left(\zeta_1^W \otimes (\zeta_2^{W})^{-1}\right).$$ 
Again, write $\kappa = \zeta_1^W \otimes (\zeta_2^W)^{-1}$, write $\kappa_i$ for the restrictions to each component, and write $z_i^W = c_1(\zeta_i^W)$. Then 
$$\rho(\kappa_2) - \rho(\kappa_1) = \textup{Sign}(W) - \textup{Sign}_\kappa(W),$$ 
where $\textup{Sign}_\kappa(W)$ is the index of the twisted signature operator on $W$.
\end{lemma}

Now we calculate the Atiyah-Patodi-Singer index of the operator $Q_{\mathbf{A},\pi}$, where $\mathbf{A}$ is a reducible trajectory on the 4-manifold $W$. By a homotopy, put $\mathbf{A}$ in a form so that it is constant sufficiently far on the ends, and constant at the boundary $\{0\} \times Y$ of the ends, constant at the unique unperturbed flat connection on $Y$ in the corresponding component of reducibles. For the statement of the following theorem, recall Definition \ref{sigdata} of signature data on a pair $(Y,E)$, and in particular signature data associated to a perturbation.

\begin{proposition}\label{index-calculation}
Suppose $(W,\mathbf{E}): (Y_1, E_1,\pi_1) \to (Y_2, E_2,\pi_2)$ is a cobordism (with cylindrical ends) between rational homology spheres, equipped with a perturbation $\pi \in \mathcal P^{(4)}_{\delta}$ which is regular on the ends. Then by Definition \ref{pertdef}, the perturbations on the ends are sufficiently small that Proposition \ref{red1} applies. 

Let $\mathbf{A}$ be a reducible connection on $\mathbf{E}$ which, sufficiently far on the ends, is constant and equal to a $\pi$-critical point; suppose $\mathbf{A}$ is in the component labeled by $r \in \textup{Red}(W,\mathbf{E})$; let $r_i \in \textup{Red}(Y_i,E_i)$ be the restrictions to the ends. Suppose $r$ corresponds to the pair of complex line bundles $\{\eta_1, \eta_2\}$, corresponding to cohomology classes $\{z_1, z_2\}$. Because the ends are rational homology spheres, we may write $z_i \in H^2_c(W;\mathbb Q)$; the compactly supported cohomology ring has a cup-product with values in $\mathbb Q$.

Let $S(r) = 1$ if $r$ is a component of $SO(2)$-reducibles and $S(r) = 3$ if $r$ contains a full reducible.

Then 
\begin{align*}I_{APS}(Q_{\mathbf{A},\pi}) &= -2(z_1 - z_2)^2 +3(b^1-b^+) \\
&+ \frac{\rho(r_2) - \rho(r_1)}{2} + \frac{\sigma_{\pi_2}(r_2) - \sigma_{\pi_1}(r_1)}{2} - \frac{S(r_1) + S(r_2)}{2}.
\end{align*}
\end{proposition}

\begin{proof}Let $W_N$ be the compact submanifold of $W$ given by including the first $[0,N]$ of each end; for $N$ sufficiently large, the operator $Q_{\mathbf{A}}^\nu$ is of product type near the boundary.

By splitting $W_N$ into three pieces, $[-N, 0] \times Y_1 \cup X \cup [0, N] \times Y_2$, we may decompose the operator $Q_{\mathbf{A},\pi}$ into its pieces on these three corresponding ends. Write $Q_\pm$ for the pieces of $Q_{\mathbf{A},\pi}$ on the corresponding ends, $Q$ for the piece on $W'$, and $A_i$ the connections $\mathbf{A}$ restricts to on $\{0\} \times Y_i$.

First we calculate $I_{APS}(Q_-)$. Write $\mathbf{A}(t)$ for the restriction of $\mathbf{A}$ to the negative end, where $t \in [-N, 0]$. The path $\mathbf{A}(t)$ is homotopic to a path $\mathbf{A}_f(t)$ so that $\mathbf{A}_f(t)$ is in the unique gauge equivalence class of $\pi(t)$-flat connection in its reducible component: This is the path used in the proof of Proposition \ref{Eigencounting}, and our assumption that $\pi_i$ on the ends lie in the sets $\mathcal P_{E_i, \delta}$ is to ensure this doesn't go awry. (See the remark immediately after Definition \ref{pertdef}.) In particular, because $Q_{\mathbf{A},\pi}$ is the operator $\frac{d}{dt} + \widehat{\text{Hess}}_{\mathbf{A}(t),\pi(t)}$, this spectral flow is by definition equal to the function $2N_{\pi_1}(r_1)$ defined in that proposition. Thus by Lemma \ref{sfindex}, $$I_{APS}(Q_-) = 2N_{\pi_1}(r_1) - \dim \text{ker} \widehat{\text{Hess}}_{A_1}.$$ A similar discussion gives $$I_{APS}(Q_+) = -S(r_2) - 2N_{\pi_2}(r_2).$$

What remains is to apply the index theorem to $D$ on $X$. If we write $Q_\theta$ to mean the corresponding ASD operator for the trivial connection, this can be read off from \cite[Proposition~8.4.1]{morgan1994l2} as giving $$I_{APS}(D) - 3I_{APS}(D_\theta) = -2p_1(\mathbf{A}) + \frac{\rho(r_2) - \rho(r_1)}{2} + 3 - \frac{h_1 + h_2}{2}.$$ 

Here $h_i = \dim \text{ker} \widehat{\text{Hess}}_{A_i}$. A detailed computation is provided in \cite[Proposition~2.6]{hedden2011chern}, but note that our sign conventions on boundary orientations and the definition of $p_1$ are the negative of theirs. Hodge theory provides the equality $$I_{APS}(Q_\theta) = -(1-b^1 + b^+),$$ so we obtain $$I_{APS}(Q)= -2p_1(\mathbf{A}) + 3(b^1-b^+) + \frac{\rho(r_2) - \rho(r_1)}{2} - \frac{h_1 + h_2}{2}.$$

Summing over these and including boundary kernel terms as in Lemma \ref{APSadd}, we obtain $$-2p_1 \mathbf{E} + 3(b^1 - b^+) + \frac{\rho(r_2) - \rho(r_1)}{2} + 2N_{\pi_1}(r_1) - \frac{h_1}{2} + \frac{h_2}{2} - 2N_{\pi_2}(r_2) - S(r_2).$$ 

Now $\mathbf{A}$ induces a reduction $\mathbf{E} \cong \mathbb R \oplus (\eta_1 \otimes \eta_2^{-1})$. Pontryagin classes are preserved under stabilization, so we want to compute $p_1(\eta_1 \otimes \eta_2^{-1})$. For a complex line bundle $\zeta$, we have $p_1 \zeta = c_1(\zeta)^2$, and considering the classes $z_i$ in $H^2_c(W;\mathbb Q)$ corresponding to $\eta_i$, we obtain $p_1 \mathbf{E} = (z_1 - z_2)^2$. 

We focus now on the last few terms. If we write 
$$D(r_i) = \dim_{\mathbb R} H^1(Y_i;\eta_1 \otimes \eta_2^{-1}),$$ 
then we have $h_1 = D(r_1) + S(r_1)$ and $h_2 = D(r_2) + S(r_2)$. Because the dimension of a vector space equipped with a nondegenerate symmetric bilinear form is the number of positive eigenvalues plus the number of negative eigenvalues, we see that $D(r)-4N_{\pi}(r)$ is the number of positive eigenvalues less the number of negative eigenvalues, and hence $$2N_{\pi_1}(r_1) - \frac{h_1}{2} + \frac{h_2}{2} - 2N_{\pi_2}(r_2) - S(r_2) = \frac{\sigma_{\pi_2}(r_2) - \sigma_{\pi_1}(r_1)}{2} - \frac{S(r_1) + S(r_2)}{2}.$$ 
\end{proof}

The above calculation did not at all depend on the fact that the connection $\mathbf{A}$ was reducible. In general, the same formula holds, where if the restriction of $\mathbf{A}$ to one of the ends is irreducible, we write $S(r_1) = 0$ for that end; and the first term should be read $-2 p_1(\mathbf{E})$, defined as a curvature integral for a connection $\mathbf{A}$ on the compact manifold $W$ which restricts to the relevant flat connections on a neighborhood of the boundary. 

Note that if we choose a different connection $\mathbf{A}'$ that restricts to the same flat connections on the boundary, the only thing that can possibly change in the index formula is $-2 p_1(\mathbf{E})$. Fix a base connection $\mathbf{A}_0$, equal to the desired flat connections near the boundary. Consider the double of $W$, with $\mathbf{A}_0$ on one half and, on the other half, an arbitrary connection $\mathbf{A}$ restricting to the desired flat connections on $\partial W$. Clearly $-2p_1$ of this new connection on a closed manifold is $-2p_1(\mathbf{A_0}) -2p_1(\mathbf{A})$. Note that $-2 p_1 \mathbf{E}$ is constant mod $8$ on a closed manifold: it reduces to 
$$2(w_2 \mathbf{E})^2 \in H^4(X;2\mathbb Z/8) = 2\mathbb Z/8,$$
where here we take the Pontryagin square to write $w_2^2 \in \mathbb Z/4$. Then we immediately have the following corollary.

\begin{corollary}\label{gr-mod-8}The relative grading $\textup{gr}_z(\alpha, \beta) \in \mathbb Z$ is independent, modulo $8$, of the choice of $z$. Therefore, we may unambiguously write $\textup{gr}(\alpha, \beta) \in \mathbb Z/8$.
\end{corollary}

The above index calculation in mind for reducible connections, we combine the $\rho$ and signature terms into a single function.

\begin{definition}Let $(Y,E)$ be a rational homology sphere equipped with an $SO(3)$-bundle and small regular perturbation $\pi$. We define the perturbed $\rho$-invariant to be $\rho_\pi(r) = \rho(r) + \sigma_\pi(r)$, where $\sigma_\pi$ is the signature datum associated to $\pi$ as in Definition \ref{sigdata}.
\end{definition}

Now if we write $D^\nu_{\mathbf{A},\pi}$ for the normal ASD operator, the linearization of the section defining the moduli spaces in $\widetilde{\mathcal B}^{e}_{\mathbf{E},z,k,\delta}(\alpha_1, \alpha_2)$ restricted to the normal space to an orbit, then we see by the discussion at the beginning of Chapter \ref{sec:4d-linear} that 
$$I(D^\nu_{\mathbf{A},\pi}) = I_{APS}(Q_{\mathbf{A}}^\nu) + S(r),$$ 
where $S(r) = \dim \mathfrak{g}_{\mathbf{A}}$ is the dimension of the space of $\mathbf{A}$-parallel gauge transformations. Further, at an $SO(2)$-reducible connection $\mathbf{A}$, the normal ASD operator splits as a sum of a `reducible part' and an `irreducible part', corresponding to the splitting of $\mathfrak g_E \cong E \cong \mathbb R \oplus \lambda$; write 
$$D^\nu_{\mathbf{A},\pi} = D_{\mathbf{A},\pi}^{\text{red}} \oplus D^{\text{irred}}_{\mathbf{A},\pi};$$ 
one has $I(D^\text{red}) = b^1 - b^+$. This is the component that $S(r) = 1$ contributes to, as the operator $Q^{\text{irred}}$ is the same as the operator $Q_\theta$ written in the above proof, and hence has index $-1+b^1 - b^+$. 

Correspondingly, we see that at an $SO(2)$-reducible $\mathbf{A}$ we have $$I(D^{\text{irred}}_{\mathbf{A},\pi}) = -2(z_1 - z_2)^2 +2(b^1-b^+) + \frac{\rho_{\pi_2}(r_2) - \rho_{\pi_1}(r_1)}{2} + 1 - \frac{S(r_1) + S(r_2)}{2}.$$ 

\begin{definition}\label{badred}Let $(W,\mathbf{E})$ be a cobordism $(Y_1, E_1) \to (Y_2, E_2)$ equipped with some small perturbation $\pi$, regular at the ends. We say that an $SO(2)$-reducible $r$ on $(W,\mathbf{E})$ is \emph{good} if $I(D^{\textup{irred}}_{\mathbf{A},\pi}) \geq 0$, and \emph{bad} otherwise.
\end{definition}

\begin{remark}In fact, $D^{\text{irred}}_{\mathbf{A},\pi}$ is a complex linear operator, so its index is \emph{even}.
\end{remark}

There is a natural condition on a cobordism-with-perturbation $(W,\mathbf{E},\pi)$ that reduces the class of bad reducibles to a simple, sometimes avoidable, set.

\begin{definition}Let $(W,\mathbf{E},\pi): (Y_1, E_1, \pi_1) \to (Y_2, E_2, \pi_2)$ be a cobordism. For an $SO(2)$-reducible component $r$ on $(W,\mathbf{E})$, we write its restriction to the two ends as $r_i$. We say that $(W,\mathbf{E},\pi)$ is \emph{$\rho$-monotonic} if for every $SO(2)$-reducible component $r$, we have $$\rho_{\pi_1}(r_1) \leq \rho_{\pi_2}(r_2).$$
\end{definition}

For a $\rho$-monotonic cobordism, one of the most mysterious terms in the index formula is nonnegative, so we may focus on the rest. 

\begin{lemma}For a $\rho$-monotonic cobordism $(W,\mathbf{E},\pi)$ with $b_1(W) = b^+(W) = 0$ and rational homology sphere ends, the only bad reducibles are $\{z_1, z_2\}$ where $z_1 - z_2$ is a torsion class on $H^2(W;\mathbb Z)$ that restricts trivially to the ends. In particular, bad reducibles can only exist if $H_1 Y_1 \oplus H_1 Y_2 \to H_1 W$ fails to be surjective.
\end{lemma}

\begin{proof}The assumption of $\rho$-monotonicity means 
$$\frac{\rho_{\pi_2}(r_2) - \rho_{\pi_1}(r_1)}{2} \geq 0,$$ 
so 
$$I(D^{\text{irred}}_{\mathbf{A},\pi}) \geq -2p_1 \mathbf{A} + 1 - \frac{S(r_1)+S(r_2)}{2}.$$ 
If both $r_i$ are $SO(2)$-reducible but not fully reducible, $$1 - \frac{S(r_1)+S(r_2)}{2} = 0;$$ if precisely one of the $r_i$ is fully reducible then that same term is $-1$, and if both of the $r_i$ are fully reducible then the final term is $-2$. Becase $b^+ = 0$, we have $-2p_1 \mathbf{A} = -2(z_1 - z_2)^2 \geq 0$. Therefore for a $\rho$-monotonic cobordism, $I(D^{\text{irred}}_{\mathbf{A},\pi}) \geq -2$, with equality if and only if both $r_i$ are fully reducible and $-2 p_1 \mathbf{A} = 0$. 

Because $H^2_c(W;\mathbb Q) \cong H^2(W;\mathbb Q)$, and the intersection form is nondegenerate negative definite on $H^2_c$, we see that 
$$-2p_1 \mathbf{A} = -2(z_1 - z_2)^2 = 0$$ 
iff $z_1 - z_2 = 0 \in H^2(W;\mathbb Q)$; this is the same as saying that $z_1 - z_2$ is a torsion class. The assumption that the restriction to the ends is fully reducible is precisely the same as saying that $z_1 - z_2$ restricts trivially to the ends. Applying the universal coefficient theorem, we obtain a class in $\text{Ext}(H_1 W, \mathbb Z)$ which restricts trivially to $\text{Ext}(H_1 Y_1, \mathbb Z) \oplus \text{Ext}(H_1 Y_2, \mathbb Z)$, and so 
$$\text{Ext}(H_1 W) \to \text{Ext}(H_1 Y_1) \oplus \text{Ext}(H_1 Y_2)$$ 
is not injective. The natural isomorphism $$\text{Ext}(A, \mathbb Z) \cong \text{Hom}(A, S^1)$$ for finite abelian groups $A$ implies that if the map on $\text{Ext}$ is not injective, then the map $H_1 Y_1 \oplus H_1 Y_2 \to H_1 W$ fails to be surjective.

Because $I(D^{\text{irred}}_{\mathbf{A},\pi})$ is an even integer, if it is larger than $-2$, it is nonnegative, so $r$ is a good reducible. 
\end{proof}

\begin{definition}\label{admiss-cob}Let $(W, \mathbf{E})$ be a cobordism between 3-manifolds with $SO(3)$-bundles and signature data $(Y_1, E_1, \sigma_1)$ and $(Y_2, E_2, \sigma_2)$. We say that $(W,\mathbf{E})$ is \textit{weakly admissible} if one of the following holds.
\begin{itemize}
    \item The negative end $(Y_1, E_1)$ is admissible, meaning that $w_2(E_1)$ only lifts to non-torsion classes in $H^2(Y_1;\mathbb Z)$.
    \item $\beta w_2(\mathbf{E}) \neq 0 \in H^3(W;\mathbb Z)$, where $\beta$ is the integral Bockstein homomorphism. 
    \item $b_1(W) = b^+(W) = 0$, for every $\beta \in \textup{Red}(W, \mathbf E)$ restricting to $\beta_i$ on the ends, $$\rho(\beta_2) - \rho(\beta_1) + \sigma_{\pi_2}(\beta_2) - \sigma_{\pi_1}(\beta_1) \geq 0,$$ and $H_1(Y_1) \oplus H_1(Y_2) \to H_1(W)$ is surjective. That is, $(W,\mathbf{E})$ is $\rho$-monotonic and supports no bad reducibles. 
    \item $\mathbf{E}$ is non-trivial, $b_1(W) = 0$ and $b^+(W) > 0$. 
\end{itemize}
Later we will need a similar notion for $U(2)$-bundles. If $\widetilde{\mathbf{E}}$ is a $U(2)$-bundle on the cobordism, we say that it is weakly admissible if its reduction to an $SO(3)$-bundle $\mathbf{E}$ is weakly admissible. 
\end{definition}

The weakly admissible $U(2)$-bundles correspond to the first, third, and fourth cases above; the second case precisely \textit{means} that $\mathbf{E}$ admits no lift to a $U(2)$-bundle, which implies that it admits no reducible connections.

We will soon see that every item on this list admits a regular perturbation. It should be noted that in fact we can achieve regular perturbations when $(Y_2, E_2)$ is admissible but the negative boundary component is not. We do not include these in our definition of weakly admissible bundles as they do not glue together well in general. The composite of a cobordism from a rational homology sphere to an admissible bundle, and then back to a rational homology sphere, need not be weakly admissible.

For the definition above, we have the following.

\begin{lemma}\label{composable-cob}The composite of two weakly admissible cobordisms remains weakly admissible.
\end{lemma}
\begin{proof}First, it is clear that the composite of the first type of weakly admissible cobordism with any other remains weakly admissible. That the composite of any cobordism and one with $\beta w_2 \mathbf{E} \neq 0$ still has $\beta w_2 \mathbf{E} \neq 0$ follows immediately from the naturality of the Bockstein and the fact that $w_2$ is natural under restriction. Suppose we have weakly admissible cobordisms of the third or fourth type that are not of the first or second. Then their boundaries are rational homology spheres, and the composite also has $b_1(W) = 0$; the term $b^+(W)$ is additive under gluing along rational homology spheres, so the composite of any of the third or fourth type with the fourth type is again of the fourth type.

So what remains to check is that the composite of any cobordism of the third type with another of that type remains of that same type. It is clear that the $\rho$-monotonicity condition is additive; the interesting thing is to ask about is the homological condition. At first glance, it is mysterious why composites of cobordisms satisfying this condition should still satisfy this condition; it is made more clear by remembering the point of that condition.

Let $\mathbf{A}$ be a reducible on the cobordism $W$ with $b_1(W) = 0$, corresponding to cohomology classes $\{z_1, z_2\}$; write $r_i$ for the restriction of $\mathbf{A}$ to the corresponding boundary component.
The $\rho$-monotonicity condition and homological condition, combined, are equivalent to the following: the index $$I(D^\text{irred}_{\mathbf{A},\pi}) = -2(z_1-z_2)^2 - 2b^+(W) + \frac{\rho_{\pi_2}(r_2) - \rho_{\pi_1}(r_1)}{2} + 1 - \frac{S(r_1) + S(r_2)}{2}$$ is nonnegative for all $(\mathbf{A},\pi)$. (This is the index computed above Definition \ref{badred}.)

Now let $W_1$ and $W_2$ be weakly admissible cobordisms between rational homology spheres equipped with reducibles $\mathbf{A}_i$; write $\mathbf{A}_{12}$ for the reducible on the composite cobordism $W$, and $r_1, r_2, r_3$ for the restrictions of $s$ to the successive 3-manifolds that serve as boundary components of $W_1$ and $W_2$. It is easy to see that $$I(D^{\text{irred}}_{\mathbf{A}_{12},\pi}) = I(D^{\text{irred}}_{\mathbf{A}_1,\pi}) + I(D^{\text{irred}}_{\mathbf{A}_2,\pi}) + S(r_2) - 1.$$

Now the index $I(D^{\text{irred}}_{\mathbf{A}_i,\pi})$ is non-negative for all reducibles; so for the composite, we find $I(D^{\text{irred}}_{\mathbf{A}_{12},\pi}) \geq S(r_2) - 1.$ But $r_2$ is a reducible, and so $S(r_2)$ is either $1$ or $3$; thus we obtain $I(D^{\text{irred}}_{\mathbf{A}_{12},\pi}) \geq 0$, as desired.
\end{proof}

Consider the situation that $W_1$ is a cobordism from a homology sphere to a rational homology sphere, and $W_2$ is a cobordism from a rational homology sphere back to a homology sphere, so that $H_1 W_i$ are both isomorphic to the middle boundary component. Then clearly the composite cobordism fails the homological condition, even though it tautologically satisfies $\rho$-monotonicity (all reducibles are full reducibles, where the $\rho$-invariant and signature data are both zero). So even though they satisfy the homological condition, one of $W_1$ and $W_2$ must fail to be $\rho$-monotonic. 

For every weakly admissible cobordism, all reducible configurations that exist for generic perturbation $\pi$ have non-negative normal index. (In the admissible case, we may choose the perturbation $\pi$ so that there are no $\pi$-ASD connections on the cobordism.)

In the first three cases, this remains true when considering paths of perturbations. In the final case, reducibles only appear generically in families of perturbations of dimension $b^+(W)$, but do not necessarily have non-negative index; when they arise, we may not achieve transversality. Strengthening the $\rho$-monotonicity requirement does not solve this problem.

We conclude this section with some remarks on gradings.

\begin{lemma}\label{cob-grading-rel}
Let $(W,\mathbf{E}): (Y_1,E_1,\pi_1) \to (Y_2, E_2,\pi_2)$ be a cobordism between 3-manifolds equipped with weakly admissible bundles and regular perturbations. If $\alpha_i$ and $\beta_i$ are choices of critical orbits, then $$\textup{gr}^W(\alpha_1, \alpha_2) - \textup{gr}^W(\beta_1, \beta_2) =  \textup{gr}(\alpha_1, \beta_1) - \textup{gr}(\alpha_2, \beta_2).$$
\end{lemma}
This follows immediately from the additivity property of the grading and the fact that $\text{gr}(\beta, \beta) = 0$.

\begin{corollary}\label{cob-grading-abs}Let $(W,\mathbf{E}): (Y_1, \pi_1) \to (Y_2, \pi_2)$ be a cobordism between rational homology spheres equipped with the trivial bundle and small regular perturbations. Write $\theta_i$ for the corresponding trivial connections. We have 
$$w_2(\mathbf{E}) \in H^2(W, \partial W;\mathbb Z/2),$$
and so we may use the Pontryagin square to write 
$$w_2(\mathbf{E})^2 \in H^4(W,\partial W;\mathbb Z/4) = \mathbb Z/4,$$ 
and $2w_2(\mathbf{E})^2 \in 2\mathbb Z/8$.

Then $$\textup{gr}^W(\alpha_1, \alpha_2) = -2w_2(\mathbf{E})^2  +3(b^1 - b^+) + \textup{gr}(\theta_2, \alpha_2) - \textup{gr}(\theta_1, \alpha_1).$$

Suppose $(W,\mathbf{E},\mathbf{A})$ is a cobordism from $(Y,E,\pi,\alpha)$ to itself, where $Y$ is a rational homology sphere equipped with a regular perturbation $\pi$ and $\pi$-flat connection $\alpha$. Then $$\text{gr}^W(\alpha, \alpha) = -2w_2(\mathbf{E})^2 + 3(b^1(W)-b^+(W)).$$
\end{corollary}

\begin{proof}Recall Definition \ref{grading} that $\overline{\text{gr}}(\alpha, \beta)$ is given by $\text{ind}(Q^\nu_{\mathbf{A},\pi})$ for any connection $\mathbf{A}$ connecting $\alpha$ and $\beta$ with perturbation $\pi$ limiting to the fixed perturbations on the ends, and $$\text{gr}(\alpha, \beta) = \overline{\text{gr}}(\alpha, \beta) + 3 - \dim \alpha.$$
\noindent
Proposition \ref{index-calculation} then provides us with $$\overline{\text{gr}}^W(\theta_1, \theta_2) = -2w_2(\mathbf{E})^2  - 3(1-b^1+b^+),$$ using the fact that $p_1 \equiv w_2^2 \pmod 4$ and $S(\theta) = 3$, as well as the vanishing of the invariants $\rho(\theta) = \sigma_{\pi}(\theta) = 0$. Therefore $$\text{gr}^W(\theta_1, \theta_2) = -2w_2(\mathbf{E})^2 + 3(b^1 -b^+).$$

Then the conclusion is simply a special case of the previous lemma. 

The final claim follows from the same index calculation we used to calculate $\text{gr}^W(\theta_1, \theta_2)$; the key points are that the $\sigma_\pi(\alpha)$ terms cancel, and that 
$S(\alpha)$ is $3 - \dim \alpha.$
\end{proof}

We may say something about the gradings of reducibles. We begin with the fully reducible connections.

\begin{proposition}\label{gr-fullred}Let $Y$ be a rational homology 3-sphere equipped with a trivial bundle and small regular perturbation $\pi$. If $\Theta$ and $\Theta'$ are fully reducible connections, then $\textup{gr}(\Theta, \Theta') \in 4\mathbb Z/8$.
\end{proposition}
\begin{proof}It suffices to prove this when $\pi = 0$; the spectral flow description above, as well as the assumption that $\pi$ is sufficiently small that no eigenvalues cross the weight $\delta$, implies the grading is the same for arbitrary perturbation.

Now recall that there is an action of $\mathcal G_E/\mathcal G^e_E = H^1(Y;\mathbb Z/2)$ on both the 3-dimensional configuration space $\widetilde{\mathcal B}^e_{E,k}$. The unperturbed Chern-Simons functional is invariant under the full gauge group, and so the critical set is preserved by this action. Therefore, we also have an action on the 4-dimensional configuration space $\sqcup_{\alpha, \beta}\widetilde{\mathcal B}^e_{\pE,k,\delta}(\alpha,\beta)$; we take the disjoint union here because the action of $x$ takes $\widetilde{\mathcal B}^e_{\pE,k,\delta}(\alpha,\beta)$ to $\widetilde{\mathcal B}^e_{\pE,k,\delta}(x \cdot \alpha,x \cdot \beta)$.

The identification in Proposition \ref{red-enumerate} implies that $H^1(Y;\mathbb Z/2)$ acts freely and transitively on the set of fully reducible critical points (or what is the same, the set of fully reducible points). So there is a unique $x \in H^1(Y;\mathbb Z/2)$ with $x \cdot \Theta = \Theta'$. Now pick a connection $\mathbf{A} \in \widetilde{\mathcal B}^e_{\pE, k,\delta}(\Theta, \Theta')$. Because the action of $H^1(Y;\mathbb Z/2)$ preserves the ASD equations, the index of $Q_\mathbf{A}$ agrees with the index of $Q_{x \cdot \mathbf{A}}$.

Therefore, $\text{gr}(\Theta, \Theta') = \text{gr}(\Theta', \Theta)$; here we use that $2x = 0$, so 
$$x \cdot \Theta' = (x + x) \cdot \Theta = \Theta.$$ 
But we know $$2\text{gr}(\Theta, \Theta') = \text{gr}(\Theta, \Theta') + \text{gr}(\Theta', \Theta) = \text{gr}(\Theta, \Theta) = 0.$$ Because these take values in $\mathbb Z/8$, we see that $\text{gr}(\Theta, \Theta')$ is a multiple of $4$.
\end{proof}

We may make a similar observation about general reducibles (both $SO(2)$- and fully reducible connections), with a completely different proof.

\begin{proposition}\label{gr-red}Let $Y$ be a rational homology 3-sphere equipped with a weakly admissible bundle $E$ and a small regular perturbation $\pi$. If $\alpha$ and $\beta$ are reducible critical orbits, then $\textup{gr}(\alpha, \beta)$ is even.
\end{proposition}
\begin{proof}First, we show that there exists a cobordism $(W,\mathbf{E}): (Y, E) \to (Y, E)$ so that there is a reducible connection on $(W,\mathbf{E})$ restricting to $\alpha$ and $\beta$ on the corresponding ends. To see this, recall from Proposition \ref{action2} that components of reducibles are classified by pairs of cohomology classes $\{z_1, z_2\} \subset H^2(W;\mathbb Z)$ so that $z_1 + z_2$ is a fixed integral lift of $w_2 \mathbf{E}$ (and, in particular, if there are reducible components such an integral lift exists). We see, therefore, that it suffices to show that every oriented closed $3$-manifold equipped with a pair of cohomology classes $(z_1, z_2) \subset H^2(Y;\mathbb Z)$ is null-bordant through an oriented 4-manifold equipped with a pair of cohomology classes $(z_1^W, z_2^W)$ that restrict to the $z_i$ on the boundary. Because pairs of cohomology classes are classified by maps to $\mathbb{CP}^\infty \times \mathbb{CP}^\infty$, we are asking that $\Omega_3^{\text{or}}(\mathbb{CP}^\infty \times \mathbb{CP}^\infty) = 0$. This follows immediately from the existence of the Atiyah-Hirzebruch spectral sequence, the fact that $\mathbb{CP}^\infty \times \mathbb{CP}^\infty$ has cohomology only in even degrees, and $\Omega_i^{\text{or}} = 0$ for $1 \leq i \leq 3$, as this implies the terms $E^2_{k, 3-k} = 0$ for all $k$, and hence the same is true of $E^\infty$.

Now if $\eta_i^W$ are the associated complex line bundles to $(z_1^W, z_2^W)$, the $SO(3)$-bundle is 
$$\mathbf{E} \cong \mathbb R \oplus \eta_1^W \otimes (\eta_2^W)^{-1}$$
and the reducible component is labelled by $\{z_1^W, z_2^W\}$; choosing such a bounding manifold for both $(Y,\alpha)$ and $(\overline Y, \beta)$, we may simply take the connected sum to obtain the desired cobordism $(Y,E,\alpha) \to (Y,E,\beta)$.

Now recall from Lemma \ref{cob-grading-rel} that $$\text{gr}^W(\beta, \alpha) - \text{gr}^W(\beta, \beta) =  \text{gr}(\alpha, \beta).$$
To show that $\text{gr}(\alpha, \beta)$ is even, our goal is to show that $\text{gr}^W(\beta, \alpha) \equiv \text{gr}^W(\beta, \beta) \pmod 2.$ Pick a reducible connection $\mathbf{A}$, asymptotic to $\beta$ at $-\infty$ and $\alpha$ at $+\infty$; the previous discussion amounts to saying we may do so.

First, we remark that the $\mathbf{A}$ enjoys a splitting $\theta \oplus A$, for a connection $A$ on a complex line bundle and $\theta$ the trivial connection on the trivial real line bundle. We may thus write $I(Q_{\mathbf{A}}) = I(Q_{\theta}) + I(Q_{A})$. Because the index of a complex linear operator is even, we see that $I(Q^\nu_{\mathbf{A}}) \equiv I(Q_\theta) \pmod 2$ and $$I(Q_\theta) = b^1(W) - 1 - b^+(W).$$ Therefore, because $\beta$ is reducible and hence $3 - \dim \beta$ is odd, we find that $$\text{gr}^W(\beta, \alpha) \equiv b^1(W) - b^+(W) \pmod 2.$$

As for $\text{gr}^W(\beta, \beta)$, the final part of Corollary \ref{cob-grading-abs} says that $$\text{gr}^W(\beta, \beta) = -2w_2(\mathbf{E})^2 + 3(b^1(W)-b^+(W)).$$ Reducing modulo $2$, we find that $\text{gr}^W(\alpha, \alpha) \equiv b^1(W) - b^+(W) \pmod 2,$ and so $$\text{gr}^W(\alpha, \alpha) \equiv \text{gr}^W(\beta, \alpha) \pmod 2,$$ as desired. 
\end{proof}

\section{Transversality for the cylinder and cobordisms}\label{sec:4d-trans}
\begin{theorem}\label{trans2}Let $E$ be a weakly admissible bundle over a 3-manifold $Y$. Suppose a perturbation $\pi_0\in \mathcal P_{E,\delta}$ has been chosen. Then by Definition \ref{pertdef} of this space of perturbations, the enumeration of reducibles in Proposition \ref{red1} holds, and the critical set $\mathfrak C_{\pi_0}$ is a finite set of nondegenerate $SO(3)$-orbits. If $\mathcal O$ is a small $SO(3)$-invariant neighborhood of these orbits, and $\mathcal P_{E,\mathcal O, \delta}$ the space of perturbations $\pi \in \mathcal P_{E,\delta}$ with $f_\pi\big|_{\mathcal O}= f_{\pi_0} \big|_{\mathcal O}$, then for a residual set of $\pi \in \mathcal P_{E,\mathcal O,\delta}$, the $\mathbb R$-reduced moduli spaces of framed instantons $\widetilde{\mathcal M}^0_{E,z,\pi}(\alpha, \beta)$ between any two critical orbits are cut out nondegenerately. Hence these are smooth $SO(3)$-manifolds of dimension $\textup{gr}_z(\alpha, \beta) + \dim \alpha - 1$ unless $\alpha = \beta$ and $z$ is trivial, in which case $\widetilde{\mathcal M}^0_{E,0,\pi}(\alpha, \alpha) = \alpha$.
\end{theorem}

\begin{proof}First, we remark on the reducibles: by Proposition \ref{red-flowlines}, the only reducible solutions are constant. We want to verify that the solutions are nondegenerate at constant trajectories; this amounts to saying that for the operator $$\frac{d}{dt} + \widehat{\text{Hess}}_{\mathbf{A},\pi}: \Omega^0 \oplus \Omega^1 \to \Omega^0 \oplus \Omega^1$$ has no cokernel other than the constant trajectories at $\text{ker}(\Delta_\alpha) \subset \Omega^0(\mathfrak g_E)$. This is true by applying a standard separation of variables argument to the adjoint operator.

Now the only points to worry about are irreducibles, for which this theorem is standard: see \cite[Section~5.5.1]{Don} or \cite{KM1}. (This is precisely where we use the perturbations in $\mathcal P_{E, \mathcal O}$ which agree with our original perturbations in a neighborhood of the $\pi$-flat connections.)
\end{proof}

Before we continue the proof of transversality for cobordisms, we will need the following lemma.

The key assumption we started with was that \emph{$I$ is nonnegative}, as the following lemma makes clear. The non-equivariant index $0$ case is written in \cite[Appendix~A.3]{Salamon}.

\begin{lemma}\label{CokernelManifold}Fix a compact Lie group $G$ and two separable $G$-Hilbert spaces $X$ and $Y$. The Banach manifold $\mathcal F_G(X, Y)$ of $G$-equivariant Fredholm maps $X \to Y$ decomposes into disjoint open sets $\mathcal F_G^I(X, Y)$, where $I \in RO(G)$ is an element of the Grothendieck group on finite-dimensional $G$-representations, labelling the index of an operator.

Each $\mathcal F^I_G(X, Y)$ is stratified by locally closed subsets $\mathcal F^{V,W}_G(X, Y)$, where $V$ and $W$ are finite-dimensional $G$-representations with $[V] - [W] = I$. These are defined to consist of those operators $T \in \mathcal F_G(X, Y)$ with $\textup{ker}(T) \cong V$ and $\textup{coker}(T) \cong W$.

Each $\mathcal F^{V, W}_G(X, Y)$ is a smooth submanifold of $\mathcal F_G(X,Y)$; the normal space at an operator $T \in \mathcal F^{V, W}_G(X, Y)$ is isomorphic to $\textup{Hom}_G(\textup{ker} \;T, \textup{coker} \;T)$.
\end{lemma}

\begin{proof}Let $T \in \mathcal F^{V,W}_G(X, Y)$, and split $X = X_0 \oplus X_1$ and $Y = Y_0 \oplus Y_1$, where $X_0 = \text{ker}(T), X_1 = X_0^\perp$, and $Y_1 = \text{Im}(T)$, while $Y_0 = Y_1^\perp$; then by assumption $T_{11}$ is an isomorphism. A neighborhood of $T$ in $\mathcal F^{V,W}_G(X, Y)$ consists of equivariant maps $T' = \begin{pmatrix}A & B\\ C & D\end{pmatrix}$, written in block-matrix form, where
$$A: X_0 \to X_0, \; B: X_1 \to X_0, \; C: X_0 \to X_1, \; D: X_1 \to X_1,$$
where $D$ is an isomorphism. 

For $T'$ sufficiently close to $T$, we have $T' \in \mathcal F^{V,W}_G(X,Y)$ if and only if the projection map $\ker(T') \to \ker(T) = X_0$ is an isomorphism. For $T' \in \mathcal F^{V,W}_G(X,Y)$, denote the equivariant inverse of the above projection by $v: X_0 \to \text{ker}(T')$ and write it in components as $(x, v_1(x))$. The defining property of $v_1$ is that $$(A + Bv_1, C + Dv_1) = 0.$$ Because $D$ is an isomorphism, we may write $v_1 = - D^{-1} C$, and then we see that $A= BD^{-1} C.$ Conversely, it is easy to see that if $A = BD^{-1} C$ then the projection map $\text{ker}(T') \to X_0$ is an isomorphism with inverse $(1, -D^{-1} C)$. So $A = BD^{-1} C$ is a defining equation for $\mathcal F^{V,W}_G(X, Y)$ near zero. The map 
$$\mathcal F^{V,W}_G(X,Y) \to \text{Hom}_G(X_0, Y_0), \quad T' \mapsto A - BD^{-1} C$$
has derivative $T' \mapsto A$ at zero, which is a surjective linear map, and hence we see that $\mathcal F^{V, W}_G(X, Y)$ is a smooth manifold near $T$.
\end{proof}

For the most difficult part of the argument --- transversality normal to the reducible locus --- we will need the following lemma. 

\begin{lemma}\label{ker-end}Let $\pi$ be a regular perturbation on $(Y, E)$, and consider the $\pi$-perturbed ASD operator on $[0,\infty) \times Y$ with Coulomb condition. 
Writing a 1-form on $[0,\infty) \times Y$ as $\psi(t) + dt \wedge \phi(t)$ for time-varying 0- and 1-forms $\phi$ and $\psi$, and a self-dual 2-form as $(dt \wedge *\omega(t))^+$ for a time-varying 1-form, this operator is
$$Q_{\mathbf{A},\pi}(\phi_t, \psi_t) = \left(\phi'_t - d_{\mathbf{A}(t)}^* \psi_t, \psi'_t - d_{\mathbf{A}(t)} \phi_t + \beta_0(t) \cdot (D_{\mathbf{A}(t)})\left(\textup{cs}+f_\pi\right)\psi_t\right);$$ 
its adjoint is the same, with the signs on the time-derivatives negated. Henceforth, we denote $$L(t) = \beta_0(t) D_{\mathbf{A}(t)}\left(\textup{cs}+f_\pi\right).$$
Then the following are true.
\begin{enumerate}
\item If $\psi$ is a time-dependent 1-form with $Q^*_{\mathbf{A},\pi}\psi = 0$, and furthermore $\psi(0) = 0$, then in fact $\psi = 0$ on the whole of $[0,\infty) \times Y$. 
\item If $Q_{\mathbf{A},\pi}(\phi, \psi) = 0$, and furthermore $\psi(0) = -d_{\mathbf{A}(0)} \sigma$ for some 0-form $\sigma$, then one may extend $\sigma$ to a time-dependent $0$-form $\sigma_t$ with 
$$(\phi_t, \psi_t) = (\sigma'_t, -d_{\mathbf{A}(t)} \sigma_t);$$ 
that is, if an element of the kernel is exact on the boundary, it is globally exact.
\end{enumerate}
\end{lemma}

\begin{proof}
The first statement follows immediately from Lemma \ref{Unique}; the second is more subtle. 

First, set $\sigma_t = \Delta_{\mathbf{A}(t)}^{-1}d_{\mathbf{A}(t)}^*\psi_t$. Replacing $(\phi, \psi)$ with $(\phi, \psi) - d_{\mathbf{A}} \sigma$ to the given element of the kernel, we have $d_{\mathbf{A}(t)}^* \psi_t = 0$ for all $t$. We aim to show that this $(\phi, \psi)$ is itself the differential of some 0-form on $[0,\infty) \times Y$. 

The linearized ASD equation now takes the form 
$$\psi'_t + L(t)\psi_t = d_{A(t)}\phi_t.$$ 

Fix an isomorphism, $C^1$ in $t$, $$e_t: \text{ker}(d_{\mathbf{A}(0)}^*) \to \text{ker}(d_{\mathbf{A}(t)}^*).$$ We use this to write $\psi_t = e_t(f_t)$, where $f_t$ is a time-dependent element of $\text{ker}(d_{\mathbf{A}(0)}^*)$. Then we may write $\psi'_t = e_t(f'_t) + e_t'(f_t)$, where $e'_t: \text{ker}(d_{\mathbf{A}(0)}^*) \to \Omega^1$ is a continuous linear map, varying continuously in time. 

Projecting to $\text{ker}(d_{\mathbf{A}(t)}^*)$ (and calling the projection operator $\Pi_t$), we obtain the following equation on $\text{ker}(d_{\mathbf{A}(0)}^*)$:

$$f'_t + (e_t^{-1}\Pi_t L(t) e_t)(f_t) = -(e_t)^{-1}(\Pi_t e'_t)(f_t).$$

The operator $e_t^{-1} (\Pi_t L(t) e_t)$ differs from a symmetric and densely defined operator by a bounded operator, and the operator $e_t^{-1}(\Pi_t e'_t)$ is bounded, so this satisfies the conditions of Lemma \ref{Unique}; therefore, $f_t = 0$ for all time, and so $\psi = 0$. The ASD equation becomes 
$$d_{\mathbf{A}(t)} \phi_t = 0.$$ 

We need to show that $(\phi_t, 0)$ is exact. But if $\sigma_t = \int_0^t \phi_s ds$, we have 
$$d_{\mathbf{A}} \sigma_t = (\phi_t, \int_0^t d_{\mathbf{A}(s)} \sigma_s ds) = (\phi_t, 0).$$ 

So any element of the kernel which is exact on the boundary is exact on the whole of $[0,\infty) \times Y$. 
\end{proof}

\begin{theorem}\label{trans3}Let $(W,\mathbf{E},\pi): (Y_1, E_1,\pi_1) \to (Y_2, E_2,\pi_2)$ be a weakly admissible cobordism, where $\pi \in \mathcal P^{(4)}_\delta$; this means, in particular, that the $\pi_i$ are regular perturbations, and the enumeration of reducibles Proposition \ref{red1} applies. Write, as usual, $\mathcal P^{(4)}_c$ for the space of perturbations which agree with the $\pi_i$ on the ends.

For any fixed constant $C > 0$, there is an open neighborhood $(\pi_-, \pi_+) \in U_C \subset \mathcal P^{(4)}_c$ with the following significance.

There is a residual set of $\pi \in U_C$ so that for any $\pi$ in this residual set, every moduli space $\widetilde{\mathcal M}_{\mathbf{E},z,\pi'}(\alpha, \beta)$ with energy at most $C$ is cut out transversely, and hence is a smooth $SO(3)$-manifold of dimension $\textup{gr}_z(\alpha, \beta) + \dim \alpha$ if it contains irreducible connections. In particular, regular perturbations\footnote{Here we use regular in the sense that moduli spaces up to some energy bound are cut out regularly; we will only need those moduli spaces with $\textup{gr}_z \leq 10$, which corresponds to an energy bound.} exist arbitrarily close to $(\pi_1, 0, \pi_2)$.
\end{theorem}

\begin{proof}As before, we argue inductively on reducibility type. We have already seen in Lemma \ref{Trivial-transverse} that the set of perturbations for which the fully reducible connections are cut out transversely in a neighborhood $(\pi_-,\pi_+) \subset U$ of any end perturbation $(\pi_-, \pi_+)$ on a weakly admissible cobordism. Furthermore, we saw in the proof of Proposition \ref{red-solutions-nobplus} and Proposition \ref{red-bplus-pos} that for a residual set of $\pi \in U$, we achieve transversality internally to the reducible locus. The same argument applies just as well to irreducible connections. The only remaining difficulty is achieving transversality \emph{normal} to the reducible locus.

Recall that above the neighborhood $U$, we have a Hilbert manifold of parameterized reducible instantons $U\mathcal M^{\text{red}}_{k,\delta}$, equipped with a projection to $U$; this projection is proper so long as you restrict to $\pi$-instantons with energy bounded by a fixed constant $C$. A fixed component of this space corresponds to instantons which respect a splitting of topological type $\mathbf{E} \cong \mathbb R \oplus \eta$.

Over this Hilbert manifold there is a Banach bundle, whose fiber above $(\mathbf{A},\pi)$ is the space of $SO(2)$-equivariant Fredholm maps

$$\left(\text{ker}(d_{\mathbf{A}}^*)_{k,\delta}\right) \to \Omega^{2,+}_{k-1,\delta}(W;\eta);$$
a more invariant way of specifying the domain is $N_{\mathbf{A}} \widetilde{\mathcal B}^{\text{red}}_{\mathbf{A},k,\delta}$, the normal bundle to the reducible locus.

Write $U\mathcal F$ for this Banach bundle with projection $$U \mathcal F \to U \mathcal M^{\text{red}}_{k,\delta}.$$

There is a stratified sub-bundle of this bundle, which is a locally closed union of manifolds: the subspace $U \mathcal F^{\geq k}$ of those Fredholm operators with cokernel isomorphic to $\mathbb C^j$, for some $j \geq k$, with $SO(2)$ acting with weight one. By Lemma \ref{CokernelManifold}, the normal space to $U \mathcal F^j$ at $(\mathbf{A},\pi, D^{\text{irred}}_{\mathbf{A},\pi})$ is
$$\text{Hom}_{SO(2)}\big(\text{ker}(D^{\text{irred}}_{\mathbf{A},\pi}), \text{coker}(D^{\text{irred}}_{\mathbf{A},\pi})\big).$$
Here $D^{\text{irred}}$ is the ASD operator restricted to the bundle $\eta$; it has no gauge fixing condition, or rather has the gauge fixing condition built into the domain of the map.

Our next goal is to show that the section $D^{\text{irred}}_{\mathbf{A},\pi}: U\mathcal M^{\text{red}}_{k,\delta} \to U\mathcal F$ is transverse to $U\mathcal F^{\geq 1}$ \emph{along $\mathcal M^{\text{red}}_{(\pi_-,\pi_+), k, \delta}$} --- that is, we are only showing that this map $D^{\text{irred}}$ is transverse to $U\mathcal F^{\geq 1}$ along the fiber above $0$. This is true for all fully reducible connections by the definition of $U$, so we may as well assume $\mathbf{A}$ is $SO(2)$-reducible.

Now suppose $(\mathbf{A},(\pi_-, \pi_+))$ is a reducible perturbed instanton whose perturbation has no interior part. Recall that if $W$ has cylindrical ends $[0,\infty) \times Y$, we define the interior holonomy perturbations on the interior of the submanifold
$$W' = W \setminus \big((1,\infty) \times Y\big).$$

Our first claim is that for any point $p$ in the interior of $W'$, there is a collection of loops $\gamma_1, \cdots, \gamma_N$ based at $p$ so that the equivariant map

$$d_{\mathbf{A}}\text{Hol}^{\vec \gamma}: \text{ker}(D^{\text{irred}}_{\mathbf{A},\pi}) \to (\mathfrak g_{E_p})^N$$ is an embedding.

This follows because $\omega$ must be exact over the interior of $W'$ by Lemma \ref{dense}, and therefore by the second point of Lemma \ref{ker-end} $\omega$ is exact along $(0,\infty) \times Y$ as well. The issue is that if $\sigma_1, \sigma_2$ are the 0-forms with $d_{\mathbf{A}} \sigma_i = \omega$ over the interior and ends, respectively, it's not clear that $\sigma_1 = \sigma_2$ over $(0,1) \times Y$; however, their difference is a closed 0-form on $(0,1) \times Y$.

If $\mathbf{A}$ is reducible over $(0,1) \times Y$, then a unique continuation argument shows that it is reducible of the same class over all of $(0,\infty) \times Y$; in particular, the closed $0$-form $\sigma_1 - \sigma_2$ extends to a closed 0-form on the whole of $(0,\infty) \times Y$. Adding this to $\sigma_2$, we find that indeed there is a globally defined 0-form $\sigma \in \Omega^0_{k+1,\delta,\text{ext}}(W;\mathfrak g_{\mathbf{E}})$ so that $\omega = d_{\mathbf{A}} \sigma$, and thus $\omega$ is zero as a tangent vector in $\widetilde{\mathcal B}_{\mathbf{E},k,\delta}.$

In fact, because this derivative map is $SO(2)$-equivariant, its image must lie in some complement to the $SO(2)$-fixed subspace $(i\mathbb R)^N$.

The second claim is that one may choose a finite set of points $p_1, \cdots, p_m$ in the interior of $W'$ so that the evaluation map

$$\text{coker}(D^{\text{irred}}_{\mathbf{A},\pi}) \to \bigoplus_{i=1}^m \Lambda^{2,+}(\mathfrak g_{\mathbf{E}_{p_i}})$$
is injective; again as above we may restrict to $\Lambda^{2,+}(\mathfrak g/i\mathbb R)$.

This is where we use, in an essential way, that $\pi$ does not have a component which is an interior perturbation --- to run the elliptic regularity argument that identifies $\psi \in L^2_{k,\delta}$ requires that the perturbations $D_{\mathbf{A}}\widehat{\nabla}_\pi$ extend to Sobolev spaces $L^2_{j,\text{loc}}$ for all $-k \leq j \leq k$, whereas we only have that extension through $0 \leq j$ when the perturbation includes an interior part. When there is no interior part, however, we have enough regularity that we may point-evaluate as above.

So suppose that we have enough regularity to point-evaluate, and that the above map was \emph{not} injective, for any choice of points $p_i$ in the interior of $W'$. This would imply that there was some $\psi \in \text{coker}(D^{\text{irred}}_{\mathbf{A},\pi})$ which was identically zero on the whole of $W'$. But then, by the first point of Lemma \ref{ker-end}, $\psi$ is zero on all of $W$. Thus some collection of points $p_i$ suffices.

Now $\mathbf{A}$, being $SO(2)$-reducible but not fully reducible, we may assume that the curves above (at each point $p_i$) are chosen so that
$$\text{Hol}^{\vec \gamma_j}(\mathbf{A}) \in U(1)^{N_j} \setminus (\pm 1)^{N_j}:$$
here we identify $\Gamma_{\mathbf{A}} \subset \widetilde{\text{Aut}}(\mathbf{E}_b)$ with $U(1)$. Call this element $H_j$ for convenience.

This means that when we choose the equivariant map
$$r_j: \widetilde{\text{Aut}}(\mathbf{E}_{p_j})^N \to \Lambda^{2,+}\left(\mathfrak g_{\mathbf{E}_{p_j}}\right)$$
so that the derivative $d_{H_j} r_j$ is an arbitrary $SO(2)$-equivariant map
$$\mathfrak{g}_{\mathbf{E}_{p_j}}^{N_j} \to \Lambda^{2,+}\left(\mathfrak g_{\mathbf{E}_{p_j}}\right),$$
and in particular can restrict to an arbitrary chosen $SO(2)$-equivariant map

$$\left(\mathfrak g_{\mathbf{E}_{p_j}}/i\mathbb R\right)^{N_j} \to \Lambda^{2,+}\left(\mathfrak g_{\mathbf{E}_{p_j}}/i\mathbb R\right).$$

Thus, finally, we may choose these maps $r_j$ so that, upon restricting to the image of $\text{ker}(D^{\text{irred}}_{\mathbf{A},\pi})$ and projecting to the image of $\text{coker}(D^{\text{irred}}_{\mathbf{A},\pi})$, the maps $r_j$ give an arbitrary chosen $SO(2)$-equivariant map
$$\text{ker}(D^{\text{irred}}_{\mathbf{A},\pi})\to \text{coker}(D^{\text{irred}}_{\mathbf{A},\pi}).$$

As usual, now, we may approximate these collection of loops and maps $r_j$ by elements of the fixed countable dense set of interior holonomy perturbations $\pi_i$, choosing the self-dual 2-forms supported in a neighborhood of the points to very well approximate $\delta$-masses at those points.

So we may approximate an arbitrary $SO(2)$-equivariant map from the kernel to the cokernel as a map induced by $D^{\text{irred}}_{\mathbf{A}}\widehat{\nabla}_{\pi_i}$.

Thus, the map
$$D^{\text{irred}}: U \mathcal M^{\text{red}}_{k,\delta} \to U\mathcal F$$ is indeed transverse to $U \mathcal F^{\geq 1}$ above $(\pi_-, \pi_+) \in U$. Now recall from Lemma \ref{C-proper} that, so long as you restrict to those instantons of energy at most $C$, the map

$U \mathcal M^{\text{red},\leq C}_{k,\delta} \to U$ is proper, and in particular closed, so the tube lemma applies: there is a neighborhood $(\pi_-, \pi_+) \in U_C \subset U$ so that the restriction of $D^{\text{irred}}$ to $U_C\mathcal M^{\text{red}}$ is transverse to $U_C\mathcal F$ for those $\pi$-instantons with energy at most $C$. It is safe to assume that $C$ is not the energy of any reducible $\pi$-instanton above $U_C$, so that this is still a smooth manifold; we drop the $C$ from the notation from $\mathcal M^{\text{red}}$ for convenience.

Write $$(D^{\text{irred}})^{-1}(U_C\mathcal F^{\geq 1}) = U_C \mathcal M^{\text{red},\geq 1}_{k,\delta};$$ this is a union of a countable collection of locally closed submanifolds of positive codimension. By Sard's theorem, the regular values of the projection $$U_C \mathcal M^{\text{red},\geq 1}_{k,\delta} \to U_C$$ form a residual set; for $\pi$ to be regular simply means that $\text{coker}(D^{\text{irred}}_{\mathbf{A},\pi}) = 0$ for all reducible $\pi$-instantons $\mathbf{A}$.

This is what we set out to prove.
\end{proof}

It is worth recalling in the above that we may choose a constant $C$ so that for all $\pi \in U_C$, any moduli space of instantons with $\text{gr}_z(\alpha, \beta) \leq 10$ consists of instantons with energy at most $C$. In particular, there is no need to take $C$ arbitrarily large, and no risk of the open sets $U_C$ shrinking to zero as one attempts to do so. 

\section{Gluing}\label{sec:4d-gluing}
We follow the approach to gluing given in \cite[Chapter~19]{KMSW}, and in particular we need to briefly discuss weighted Sobolev spaces and the perturbed ASD equations on compact cylinders.

We write $Z^T = [-T, T] \times Y$ and $Z^\infty = [0, \infty) \times Y \sqcup (-\infty, 0] \times Y$; we view $Z^\infty$ as a limit of the $Z^T$, stretching until what used to be $0$ becomes the point at $\infty$. We uniformly have $\partial Z^T = Y \sqcup \overline Y$. The function spaces of interest to us are the Sobolev spaces of sections of vector bundles over $Z^T$. The $L^2_k$ Sobolev space is the completion of the space of compactly supported smooth sections on $Z^T$; notice that there are no boundary conditions on these sections. As a remark before continuing, the restriction map to the boundary takes value in the $L^2_{k-1/2}$ Sobolev space.

We define the \emph{weighted} Sobolev spaces on $Z^T$ as in \cite[Page~72]{Lin}: let $\sigma_\delta: \mathbb R \to [-\delta, \delta]$ be an odd smooth function such that $\sigma_\delta(t) = -\delta$ for $t \geq 1$ and $\sigma_\delta(t) = \delta$ for $t \leq -1$. For each finite $T \geq 2$ we let $g_{T,\delta}$ be the positive smooth even function on $[-T, T]$ which is equal to $1$ on the boundary and has $\sigma_\delta = \log(g_{T,\delta})'$. For $t \in [1, T]$, we have $g_{T,\delta} = e^{-\delta(t-T)}$, an exponentially decreasing function with final value $1$. Then for sections on the finite cylinder $Z^T$ we set 
$$\|f\|_{L^2_{k,\delta}} = \|g_{T, \delta} \cdot f\|_{L^2_k}.$$ 
For $T = \infty$ then as before we will use $g_{\infty, \delta} = e^{|t|\delta}$. For finite $T$, this is equivalent to the usual $L^2_k$ norm, with implicit constants growing exponentially in $T$. 

Now as before we may introduce the moduli space of instantons on the finite cylinder $Z^T$ in the same component as the constant solution $\gamma_A$ in the usual way: we consider the configuration space $\gamma_A + \Omega^1_{k,\delta}(Z^T;\mathfrak g_E)$ and quotient by the space $\mathcal G^{e,h}_{E,k+1,\delta}$ of gauge transformations. In the case $T = \infty$, the $h$ (for harmonic) denotes that these gauge transformations should be asymptotic on the noncompact end to elements of $\Gamma_A$, possibly different on each component. 
These moduli spaces with no boundary conditions are infinite-dimensional and depend on the Sobolev index $k$, as the restriction of an $L^2_k$ ASD connection to the boundary may only be of Sobolev class $L^2_{k-1/2}$, in any $L^2_{k+1/2}$ gauge equivalence class.

When defining \emph{framed} moduli spaces on the finite cylinder $Z^T$, we set 
$$\widetilde{\mathcal M}_{A,k} := \mathcal A_{A,k}(Z^T) \times E_{(-T,b)}\big/ \mathcal G^{e}_{\pE,k+1,\delta}.$$ 
Here $A$ denotes that we are based at the constant trajectory $\gamma_A$. Because these moduli spaces depend on the Sobolev index $k$, we include it in our notation.

In the infinite case $Z^\infty$, there should be \emph{two} framings, one for each component, which have the same value in the orbit through $A$ in $\widetilde{\mathcal B}^e_E$ when we take the holonomy to $\pm \infty$. To mirror the case of finite cylinders, we write this as $$\widetilde{\mathcal A}_{A,k,\delta}(Z^\infty) := \mathcal A_{A,k,\delta}(Z^\infty) \times E_{(0^+,b)} \times E_{(0^-, b)};$$ the notation $(0^+, b)$ singifies that this lies in the component $[0, \infty) \times Y$.

Write $\widetilde{\mathcal M}'_{A,k,\delta}(Z^\infty)$ for the quotient. Now we may take $\mathbf{A}$-holonomy to $(\infty, b)$ or $(-\infty, b)$, respectively, and project the framing factors to $E_b/\Gamma_\alpha$; and we write $\widetilde{\mathcal M}_{A,k,\delta}(Z^\infty)$ for the subset on which the two framings project to the same element of $E_b/\Gamma_\alpha$.

For all $T$, there is a restriction map $$R^T: \widetilde{\mathcal M}_{A,k,\delta}(Z^T) \to \widetilde{\mathcal B}^e_{E,k-1/2}(Y \sqcup \overline Y).$$ When $T = \infty$, an element of $\widetilde{\mathcal M}_{A,k,\delta}(Z^\infty)$ is an equivalence class of connection with a framing at each basepoint on the two boundary components. Therefore, $R^\infty$ may be given by literally restricting the connection and framing to the boundary (up to gauge equivalence). 

For finite $T$ this is a little more complicated. On the left boundary component $Y$, this is the usual restriction map (the framed basepoint is on $Y$), but on $\overline Y$, it is given by restriction in the connection coordinate and parallel transport from $-T$ to $+ T$ along $\mathbb R \times \{b\}$ in the framing coordinate.

Now we turn to slices for the gauge group action so that we may compute neighborhoods of $\widetilde{\mathcal M}_{A,k,\delta}$ as an equation on a linear space.

As before, we may consider the Coulomb slice $\text{ker}(d_{\gamma_A}^*) \subset \Omega^1_{k,\delta}(Z^T;\mathfrak g_E)$ where $\gamma_A$ denotes the constant trajectory at $A$; in particular $$d_{\gamma_A}^*(dt \wedge \sigma(t) + \omega(t)) = -\sigma'(t) + d_A^* \omega(t).$$ The Coulomb slice is no longer a slice for the gauge group action: if we try to solve the equation $d^*_{\gamma_A}(d_{\gamma_A} \sigma + \omega) = 0$, to write an arbitrary element of $\Omega^1$ as a sum of an element of $\text{Im}(d_{\gamma_A})$ and an element of $\text{ker}(d_{\gamma_A}^*)$, we find that 
$$d_{\gamma_A}^* d_{\gamma_A} \sigma = - d^*_{\gamma_A}\omega$$ 
has a unique solution $d_{\gamma_A} \sigma$ for each $\omega$ and \emph{fixed boundary values} $d_{A} \sigma \big|_{Y \sqcup \overline{Y}}$.

We obtain the \emph{Coulomb-Neumann slice} around $\gamma_A$, written $\mathcal{CN}_{A,k,\delta}$: the subset of $\omega \in \Omega^1_{k,\delta}(Z^T;\mathfrak g_E)$ on which $d_{\gamma_A}^* \omega = 0$ and, writing $\omega = \phi + dt \wedge \psi,$ where $\psi$ is a time-dependent $0$-form, we demand that 
$$\psi\big|_{\partial Z^T} = 0.$$
Every connection on $Z^T$ sufficiently close to $\gamma_A$ is gauge equivalent to one on the Coulomb-Neumann slice, and the only remaining ambiguity is that $\mathcal{CN}_{A,k,\delta}$ carries the action of the stabilizer $\Gamma_{\gamma_A}$ in the gauge group, and $\gamma_A + a$ is gauge equivalent to $\gamma_A + u(A)$. In particular, observe that $\mathcal{CN}_{A,k,\delta} \times_{\Gamma_{\gamma_A}} SO(3)$ gives an $SO(3)$-invariant neighborhood of $\gamma_A$ in $\widetilde{\mathcal B}(Z^T)$.

Before moving on, we recall Lin's abstract Morse-Bott gluing theorem. Let $E$ be a vector bundle over a closed oriented 3-manifold $Y$, and let $L$ be a self-adjoint elliptic operator acting on $L^2_k(Y;E)$ with kernel $H_0$. In what follows, we consider the operator $D := d/dt + L$ acting on finite cylinders $Z^T = [-T,T] \times Y$ and infinite cylinders $Z^\infty = [0,\infty) \times Y \sqcup (-\infty, 0] \times Y$.

On the infinite cylinder, we write 
$$L^2_{k,\delta,\text{ext}}(Z^\infty;\pE) := H_0 + L^2_{k,\delta}(Z^\infty;\pE),$$ 
where $H_0$ indicates sections $s_{h_0}$ in $\text{ker}(D)$ on $Z^\infty$ which are constant in time at some element $h_0$ in $H_0 = \text{ker}(L)$. In fact, we write 
$$\mathcal E^T_\delta = L^2_{k,\delta}(Z^T;\pE), \;\;\; \mathcal E^\infty_\delta = L^2_{k,\delta,\text{ext}}(Z^\infty, \pE)$$ 
and 
$$\mathcal F^T_\delta = L^2_{k-1,\delta}(Z^T;\pE), \;\;\; \mathcal F^\infty_\delta = L^2_{k-1,\delta}(Z^\infty, \pE).$$ 

We have a projection map $\Pi^\infty_0: \mathcal E^\infty_\delta \to H_0$ given by taking the asymptotic value at $\pm \infty$ and projection maps $\Pi^T_0: \mathcal E^T_\delta \to H_0$ given by projection onto the subspace of constant sections in the $L^2_{k,\delta}$ norm (equivalently, in the $L^2_{0,\delta}$ norm).

Suppose we are given a bounded linear operator 
$$\Pi: L^2_{k-1/2}(Y \sqcup \overline Y; E) \to H$$ 
for some Hilbert space $H$; by restriction to the boundary this induces maps 
\begin{align*}\Pi: L^2_{k,\delta}(Z^T;\pE) &\to H\\
\Pi: L^2_{k,\delta,\text{ext}}(Z^\infty;\pE) &\to H.
\end{align*}

We now assume, crucially, that the map $$(D, \Pi^\infty_0, \Pi): \mathcal E^\infty_\delta \to \mathcal F^\infty_\delta \oplus H_0 \oplus H$$ is an isomorphism, which implies the same for $\delta'$ sufficiently close to $\delta$. Now suppose we are given a \emph{non-linear} map 
$$\alpha: C^\infty(Z^T;\pE) \to L^2_{\text{loc}}(Z^T;\pE),$$ 
obtained by restriction to slices from a map $\alpha_0: C^\infty(Y;\pE) \to L^2_{\text{loc}}(Y;\pE)$. We assume that $\alpha$ extends to a smooth map 
$$L^2_k([-1,1] \times Y; \pE) \to L^2_{k-1}([-1,1] \times Y; \pE)$$ 
with $\alpha(h_0)$ for every $h_0 \in H_0$ and we assume that $\alpha$ is purely nonlinear, in the sense that its derivative at $0 \in L^2_k([-1,1] \times Y;\pE)$ is zero. This implies that $\alpha$ defines smooth maps $\mathcal E_\delta^T \to \mathcal F_\delta^T$ with the same property.

Now write $$F^T = D + \alpha: \mathcal E^T_\delta \to \mathcal F^T_\delta$$ and $$M(T) = \left(F^T\right)^{-1}(0) \subset \mathcal E_\delta^T.$$ The following is \cite[Proposition~3.5.15]{Lin}.

\begin{proposition}For $T \in [T_0, \infty]$, the sets $M(T)$ are Hilbert submanifolds of $\mathcal E^\infty_\delta$ in a neighborhood of $0$. There exist $\eta > 0$ and smooth maps $$u(T, -): B_\eta(H_0 \oplus H) \to M(T)$$ which are diffeomorphisms onto their image and have $$(\Pi_0, \Pi) u(T, (h_0, h)) = (h_0, h) \; \text{ and } \; u(T, (h_0, 0)) = s_{h_0}.$$ For $T \in [T_0, \infty]$, we have a map $$\mu_T: B_\eta(H_0 \oplus H) \to L^2_{k-1/2}(Y \sqcup \overline Y;E)$$ obtained as a composition of $u(T, -)$ with restriction to the boundary. For each $T$, $\mu_T$ is a smooth embedding; $\mu_T$ is smooth as a function on $[T_0, \infty) \times B_\eta(H_0 \oplus H)$, and $\mu_T \to \mu_\infty$ as $T \to \infty$ in the $C^\infty_{\text{loc}}$ topology on maps $B_\eta(H_0 \oplus H) \to L^2_{k-1/2}(Y \sqcup \overline Y;E)$. Finally, there is an $\eta' > 0$ independent of $T$ so that the images of the maps $u(T, -)$ contains all solutions $u \in M(T)$ with $\|u\|_{L^2_{k,\delta}} \leq \eta'.$
\end{proposition}

We will now state our main gluing lemma. In what follows, we will use Lin's abstract gluing theorem to prove our lemma in three cases, depending on the reducibility type of the orbit $\alpha$, along with some mild changes (the flavor of the changes depending on how reducible $\alpha$ is).

\begin{lemma}\label{maingluinglemma}Let $(Y,E,\pi)$ be a closed oriented 3-manifold equipped with a weakly admissible $SO(3)$-bundle and regular perturbation $\pi$. We write $\alpha$ for a $\pi$-critical orbit; if $A$ is a connection in the gauge equivalence class of $\alpha$, we write the corresponding Coulomb slice 
$$\mathcal K_\alpha = \textup{ker}(d^*_A) \cap \Omega^1_{k-1/2}(Y;\mathfrak g_E),$$ 
and $B(\mathcal K_\alpha)$ for its unit ball. The Hilbert space $\mathcal K_\alpha$ carries the action of the stabilizer of $A$ in the gauge group; thinking of $\mathcal K_\alpha$ as the normal space to a point in the orbit $\alpha$ in $\widetilde{\mathcal B}^e_{k-1/2}(Y)$, this is the same as the action of $\Gamma_\alpha$, the stabilizer of a point of $\alpha$ in $SO(3)$.

We may thus extend $\mathcal K_\alpha$ to a vector bundle 
$$\widetilde{\mathcal K}_\alpha := SO(3) \times_H \mathcal K_\alpha$$ 
over the orbit $\alpha$, equipped with an $SO(3)$ action (acting on the factor of $SO(3)$ on the left). We write the associated unit disc bundle as $B(\widetilde{\mathcal K}_\alpha)$.

There is a $T_0$ so that for all $T \in [T_0, \infty]$, we may find smooth, $SO(3)$-equivariant maps $$\widetilde u(T, -): B(\widetilde{\mathcal K}_\alpha) \to \widetilde{\mathcal M}_{\alpha, k,\delta}(Z^T)$$ which are diffeomorphisms onto neighborhoods of the constant solution $\gamma_\alpha$, and such that the map $$\widetilde{\mu}_T: B(\widetilde{\mathcal K}_\alpha) \to \widetilde{\mathcal B}^e_{E,k-1/2}(Y \sqcup \overline Y),$$ given by composing $u(T, -)$ with the restriction map $$R^T: \widetilde{\mathcal M}_{\alpha, k,\delta}(Z^T) \to \widetilde{\mathcal B}^e_{E,k-1/2}(Y \sqcup \overline Y)$$ described above using holonomy in the framing coordinate for finite $T$, has the following properties.

First, $\widetilde{\mu}_T$, being the composition of equivariant maps, is $SO(3)$-equivariant; $\widetilde{\mu}_T$ is a smooth embedding of $B(\widetilde{\mathcal K}_\alpha)$ for all $T \in [T_0, \infty]$, and $\widetilde{\mu}_T$ is smooth as a function of $[T, \infty) \times B(\widetilde{\mathcal K}_\alpha)$. Though not smooth as a function on $[T_0, \infty]$, we at least have $\widetilde{\mu}_T \to \widetilde{\mu}_\infty$ in the $C^\infty_{\textup{loc}}$ topology on $B(\widetilde{\mathcal K}_\alpha)$ as $T \to \infty$.

Finally, there is an $\eta > 0$, independent of $T$, so that the images of $\widetilde u(T, -)$ contain all solutions $[\gamma] \in \widetilde{\mathcal M}_{\alpha, k, \delta}(Z^T)$ such that $$\|\gamma - \gamma_\alpha\|_{L^2_{k,\delta}(Z^T)} \leq \eta.$$

The vector space $\mathcal K_\alpha$ has an eigenspace decomposition $\mathcal K^-_\alpha \oplus \mathcal K^+_\alpha = \mathcal K_\alpha$; correspondingly we have an $SO(3)$-equivariant fiber product decomposition $$B(\widetilde{\mathcal K}_\alpha) \cong B(\widetilde{\mathcal K}^-_\alpha) \times_{\alpha} B(\widetilde{\mathcal K}^+_\alpha).$$ We also have the decomposition $$\widetilde{\mathcal M}_{\alpha, k, \delta}(Z^\infty) = \widetilde{\mathcal M}_{\alpha, k, \delta}(\mathbb R^{\geq 0} \times Y) \times_\alpha \widetilde{\mathcal M}_{\alpha, k, \delta}(\mathbb R^{\leq 0} \times Y).$$ The map $\widetilde u(\infty, -)$ respects this decomposition: it is a fiber bundle map over $\alpha$, and on each fiber it respects the product structure.
\end{lemma}

The easiest case is when $\alpha$ is irreducible, as then we may solve the gluing problem in a Morse (not Morse-Bott) setting.

\begin{proof}[Proof of Lemma \ref{maingluinglemma} in the irreducible case]We will be studying the ASD operator with a gauge fixing condition; the linear operator $L$ on $Y$ is the perturbed operator $\widehat{\text{Hess}}_{A, \pi}$. In the irreducible case, this operator has no kernel, and Lin's gluing theorem reduces to the abstract gluing theorem \cite[Theorem~18.3.5]{KMSW}.

First we will apply the abstract gluing lemma to see the corresponding statement about \emph{unframed} moduli spaces: we will see that there are maps
$$u(T, -): B(\mathcal K_\alpha) \to \mathcal M_{\alpha, k,\delta}(Z^T)$$ 
parameterizing solutions of the perturbed ASD equations in a neighborhood of $\gamma_A$ satisfying the same conditions (without the equivariance); the corresponding restriction maps $R^T$ are given by the actual restriction of connections to the boundary (as opposed to before, where it involved a global quantity, the holonomy).

The argument in \cite[Section~18.4]{KMSW} applying the abstract gluing theorem to the Seiberg-Witten equations readily applies all the same to the ASD equations.

We sketch the argument, as we will use essentially the same argument in the reducible case. Fix a connection $A$ on $Y$ in the gauge equivalence class of $\alpha$. We identify a neighborhood of $\alpha$ in $\mathcal B^e_{E,k-1/2}$ as $\text{ker}(d_A^*)$ via the Coulomb slice.

On the cylinder $Z^T$, we identify the spaces $$\Omega^1_{k,\delta}(Z^T;\mathfrak g_E)$$ and $$\left(\Omega^0 \oplus \Omega^{2,+}\right)_{k,\delta}(Z^T;\mathfrak g_E)$$ with $$L^2_{k,\delta}(Z^T;\pi^* \! (\mathbb R \oplus T^*Y) \otimes \mathfrak g_E);$$ the first by writing every $1$-form as $\psi + dt \wedge \sigma$ for $\psi$ a time-dependentent $1$-form on $Y$, and the second by writing every self-dual 2-form as $(dt \wedge \psi)^+$.

Now note from this isomorphism that when we have $a \in \Omega^1_{k,\delta}(Z^T;\mathfrak g_E)$, we may keep track of more information at the boundary than just $a \big|_{Y \sqcup \overline Y}$; this restriction kills any term of the form $dt \wedge \sigma$. So we write $$r: \Omega^1_{k,\delta}(Z^T;\mathfrak g_E) \to \Omega^1_{k-1/2}(Y \sqcup \overline Y; \mathfrak g_E) \oplus \Omega^0_{k-1/2}(Y \sqcup \overline Y; \mathfrak g_E),$$ to record both the restriction of $a$ and its normal value $\sigma(\pm T,b)$, or $\sigma(0^\pm,b)$ in the case of the infinite cylinder.

On $\Omega^1_{k-1/2}(Y;\mathfrak g_E)$, we may project to $\text{ker}(d_A^*)$, which has an eigenspace decomposition $\mathcal K^-_\alpha \oplus \mathcal K^+_\alpha$ for the action of $D_{A,\pi}$. (The positive and negative eigenspaces swap upon orientation-reversal.) Write $$\Pi: \Omega^1_{k-1/2}(Y \sqcup \overline Y; \mathfrak g_E) \oplus \Omega^0_{k-1/2}(Y \sqcup \overline Y; \mathfrak g_E) \to \mathcal K^-_\alpha \oplus \Omega^0_{k-1/2}(Y;\mathfrak g_E) \oplus \mathcal K^+_\alpha \oplus \Omega^0_{k-1/2}(\overline Y; \mathfrak g_E);$$ this is the spectral projection on the $\text{ker}(d_A^*)$ term (projecting to negative eigenvalues on $Y$ and positive eigenvalues on $\overline Y$), and records the data of $\sigma \big|_{\partial Z^T}$.

Now, on the linear space $\Omega^1_{k,\delta}(Z^T;\mathfrak g_E)$, consider the equations

\begin{align*} d_{\gamma_A}^+ a + (a \wedge a)^+ + \left(dt \wedge \nabla_\pi(\gamma_A + a)\right)^+ = 0 \\
d_{\gamma_A}^* a = 0 \\
\Pi (a \big|_{\partial Z^T}) = c \end{align*}

Here $$c \in \mathcal K^-_\alpha \oplus \Omega^0_{k-1/2}(Y;\mathfrak g_E) \oplus \mathcal K^+_\alpha \oplus \Omega^0_{k-1/2}(\overline Y; \mathfrak g_E);$$ we henceforth write this space as $H = H^- \oplus H^+$ to simplify notation. The above equations are simply the ASD equations with boundary conditions corresponding to the spectral projection of the $\Omega^1(Y)$ component and restriction of the $\Omega^0(Y)$ component.

The linearization of these equations is \begin{align*} Q_{\gamma_A, \pi} a = 0 \\
\Pi(a \big|_{\partial Z^T}) = 0,\end{align*}
so the most important thing to check is that \emph{these linearized equations determine an isomorphism} 
$$\Omega^1_{k}(Z^\infty;\mathfrak g_E) \to (\Omega^{2,+} \oplus \Omega^0)_{k-1}(Z^\infty; \mathfrak g_E) \oplus H;$$ 
equivalently it suffices to show that 
$$\Omega^1_k([0,\infty) \times Y; \mathfrak g_E) \to (\Omega^{2,+} \oplus \Omega^0)_{k-1}([0,\infty) \times Y; \mathfrak g_E) \oplus H^-$$ 
is an isomorphism. To do this, we write $Q_{\gamma_A, \pi} = \frac d{dt} + \widehat{\text{Hess}}_{A,\pi}$, and split 
$$\Omega^0(Y) \oplus \Omega^1(Y) = \Omega^0 \oplus \text{Im}(d_A) \oplus \text{ker}(d_A^*);$$ 
rewriting $\text{Im}(d_A) \cong \Omega^0_{k+1/2}(Y;\mathfrak g_E)$, we have 
$$\widehat{\text{Hess}}_{A,\pi} = \begin{pmatrix}0 & - \Delta_A \\ -1 & 0\end{pmatrix},$$
which is invertible because $A$ is irreducible. On the other component, our operator takes the form $D_{A,\pi}$, which is invertible because $\pi$ is a regular perturbation.

First let us see that $(Q_{\gamma_A,\pi}, \Pi)$ is surjective. We may write an arbitrary element of the domain as $$(\sigma_s, \psi_s, h, c) \in (\Omega^{2,+} \oplus\Omega^0)_{k-1}([0,\infty) \times Y; \mathfrak g_E) \oplus H^-,$$ where $\sigma_s$ is a time-dependent $0$-form and $\psi_s$ is a time-dependent 1-form, thought of as the self-dual 2-form via $(dt \wedge \psi_s)^+$; $h \in \mathcal K^-_\alpha$ and $c \in \Omega^0_{k-1/2}(Y;\mathfrak g_E)$.

We may decompose $\Omega^0(Y)$ into the eigenspaces of $\Delta_A$ and $\text{ker}(d_A^*)$ into the eigenspaces of $D_{A,\pi}$. For $\phi \in \Omega^1_k([0, \infty) \times Y;\mathfrak g_E)$ we write $\phi_s = \omega_s + d_A \eta_s$, where $\omega_s \in \text{ker}(d_A^*)_{k-1/2}$. We will explicitly write out $(Q_{\gamma_A,\pi}, \Pi)$ in terms of these. 

Using the eigenspace decompositions, we may write 
\begin{align*}\sigma_s &= \sum b_\lambda(s) \sigma_\lambda \\
\eta_s &= \sum c_\lambda(s) \sigma_\lambda \\
\psi_s &= \sum d_\lambda(s) \omega_\lambda,
\end{align*}
where the coefficients range over a basis of eigenfunctions.

Now we find that these must satisfy $\sigma'(s) = \Delta_A \eta(s)$ and $\eta'(s) = \sigma(s)$. These immediately give us
\begin{align*}b_\lambda &= c'_\lambda, \\
b'_\lambda(s) &= \lambda c_\lambda(s).\end{align*}
Combining these we solve to find 
$$c_\lambda(t) = c_\lambda(0) \cdot e^{\pm \sqrt{\lambda} t}.$$ 
However, under the assumption that $\phi \in L^2_k$ and hence $\eta \in L^2_{k+1}$, the sign must be negative. By specifying $c$ above, we specify the values of $c_\lambda(0)$ for all $\lambda$, and hence both $b$ and $c$ for all time. Because $\sum c_\lambda(0)$ defines an $L^2_{k-1/2}$ function on $Y$ and the sum $\sum e^{-\sqrt{\lambda}} c_\lambda(0)$ satisfies an elliptic equation, regularity implies that this defines an $L^2_k$ function on $[0,\infty) \times Y$ (and similarly for $b$).

Now for $\psi_s$ we get the formula $\omega_\lambda'(s) = \lambda \omega_\lambda(s)$. This can only contribute to the solution if $\lambda < 0$, as otherwise $\omega$ is not $L^2$. We specified the values of $\omega_\lambda(0)$ for $\lambda < 0$ in the choice of $h \in H^-$, and thus have specified $\omega$ for all time; as above it has the desired regularity. We have constructed a solution to the equations of the desired regularity and with desired boundary values, so $(Q,\Pi)$ is surjective. 

To see that $(Q, \Pi)$ is injective, we remark that the above argument always produced \emph{unique} solutions.

This in hand, the abstract gluing theorem gives us a map $$u(T, -): B_\eta\left(\mathcal K^-_\alpha \oplus \mathcal K^+_\alpha \oplus \Omega^0_{k-1/2}(Y \sqcup \overline Y)\right) \to \Omega^1_{k,\delta}(Z^T;\mathfrak g_E)$$ parameterizing solutions in a neighborhood of $0$. To use this to get a parameterization of $\gamma_A$ in $\mathcal M_{A,k,\delta}(Z^T)$, we observe that a neighborhood of $\gamma_A$ in the latter is given by the Coulomb-Neumann slice $\mathcal{CN}_{A,k,\delta}$, and this parameterization is compatible with the restriction maps; so to obtain the desired parameterizations, we simply need to restrict the domain of $u$ to $B_\eta(\mathcal K^-_\alpha \oplus \mathcal K^+_\alpha)$, imposing the Neumann gauge condition that the $dt \wedge \sigma$ component is zero on the boundary.

Recall that projection $\widetilde{\mathcal M}(Z^T) \to \mathcal M(Z^T)$ forms a principal $SO(3)$-bundle over the irreducible instantons. Our goal is to choose a section of this over $B(\mathcal K_\alpha)$, sitting inside $\widetilde{\mathcal M}(Z^T)$ via $u(T, -)$, which will automatically give us an equivariant map $SO(3) \times B(\mathcal K_\alpha) \to \widetilde{\mathcal M}(Z^T)$, as desired. But we should be careful in how we choose this lift so that the maps $\mu_T$ have the desired properties.

We identified a neighborhood of $\gamma_A$ in $\mathcal M_{A,k,\delta}(Z^T)$ with a subset of $\mathcal{CN}_{A,k,\delta}(Z^T)$. We may thus define the section above $\mathcal{CN}_{A,k,\delta}(Z^T)$; choose the section to be 
\begin{align*}\mathcal{CN}_{A,k,\delta}(Z^T) &\to SO(3) \times \mathcal{CN}_{A,k,\delta}(Z^T) \\
\mathbf{A} &\mapsto (p, \mathbf{A})
\end{align*}
for $p$ a fixed framing.

Then writing $B(\widetilde{\mathcal K}_\alpha) = SO(3) \times \widetilde{\mathcal K}_\alpha$, our parameterization $$\widetilde u(T,-): SO(3) \times B(\mathcal K_\alpha) \to \widetilde{\mathcal M}(Z^T)$$ is given as $(\tilde u)(g,h) = \tilde s(u(h)) \cdot g$ in the `framed Coulomb-Neumann slice' $SO(3) \times \mathcal{CN}_{A,k,\delta}$. That this is a smooth embedding for all $T$ follows from the corresponding fact for unframed moduli spaces, and because $\mu_T$ is assumed smooth in finite $T$, so is $\widetilde\mu_T$; similarly $\widetilde{\mu}_T \to \widetilde{\mu}_\infty$ in the $C^\infty_{\text{loc}}$ topology.

Finally, because the unframed map $u(\infty, -)$ respects the product structure and we defined $\tilde u$ to be a fiber bundle map, the last statement is true.
\end{proof}

We quickly discuss linear models for the configuration space $\widetilde{\mathcal B}_{A,k,\delta}(Z^T)$ when $A$ is fully reducible. When $T$ is finite, choosing a base framing $p$, a neighborhood of $0$ in the Coulomb slice $$\mathcal{CN}_{A,k,\delta}(Z^T) \subset \mathcal A_{A,k,\delta}(Z^T) \times \{p\} \subset \widetilde{\mathcal A}_{A,k,\delta}(Z^T)$$ projects to a chart in $\widetilde{\mathcal B}_{\pE,k,\delta}(Z^T)$ around $\gamma_A$, so we may effectively study neighborhoods of $\gamma_A$ in the moduli space on finite cylinders by studying a subset of the equations in Coulomb-Neumann gauge. We will consider the Coulomb-Neumann slice as a subset of $\Omega^1_{k,\delta}(Z^T)$ and work there. In this case, the space $H_0 = \text{ker}(L)$ above is the Lie algebra $\mathfrak g_b$; when considered as a subset of $\Omega^1(Z^T)$, these consist of the 1-forms $dt \wedge \sigma$, where $\sigma$ is a fixed $A$-parallel section of $\mathfrak g_b$.

The discussion of the Coulomb-Neumann slice remains correct on the infinite cylinder: everything in a neighborhood of $\gamma_A$ in $\widetilde{\mathcal B}_{A,k,\delta}(Z^\infty)$ is gauge equivalent to an asymptotically decaying connection in Coulomb-Neumann gauge. Recall here that our definition of $\widetilde{\mathcal B}_{A,k,\delta}(Z^\infty)$ involves \emph{two framings}: one on each boundary component. Applying $A$-parallel gauge transformations on each component, we may change those framings arbitrarily while obtaining another connection in Coulomb-Neumann gauge.

In fact, while we will consider the space $\mathcal E^\infty_\delta = \Omega^1_{k,\delta,\text{ext}}(Z^\infty;\mathfrak g_E)$ of 1-forms which asymptotically decay to $dt \wedge \sigma$, in the case of the infinite cylinder the extra kernel arising as $dt \wedge \sigma$ is somehow illusory: applying a gauge transformation we may change any connection in this extended space to a connection in Coulomb-Neumann gauge.

Similarly, for any $a \in \mathcal{CN}_{A,k,\delta}$, we have by definition of Coulomb-Neumann gauge the integration by parts formula 
$$0 = \langle d_{\gamma_A}^* a, b \rangle_{L^2} = \langle a, d_{\gamma_A} b\rangle_{L^2}.$$ 
Above, the projection $\Pi_0$ was projection onto the constant sections $dt \wedge \sigma$ in the $L^2_{0,\delta}$ inner product; then this is defined to be the same as $\langle a, g_{T, \delta}^2 dt \wedge \sigma\rangle_{L^2},$ where $g_{T,\delta}$ is the weight function in the definition of Sobolev spaces on finite cylinders.

But we may choose an antiderivative $\frac{d}{dt} G_{T,\delta} = g_{T,\delta}^2$, and then the right-hand side of this inner product is $d_{\gamma_A} (G_{T,\delta} \sigma)$, because $\sigma$ is $A$-parallel. Therefore, for any solution in Coulomb-Neumann gauge, $\Pi^T_0 a = 0$. With this, we begin the proof.

\begin{proof}[Proof of Lemma \ref{maingluinglemma} in the fully reducible case]The equations we look to study with Lin's abstract gluing theorem are nearly the same as last time: now they are 
\begin{align*}d_{\gamma_A}^+ a + (a \wedge a)^+ + \left(dt \wedge \nabla_\pi(\gamma_A + a)\right)^+ &= 0 \\
d_{\gamma_A}^* a &= 0 \\
\Pi (a \big|_{\partial Z^T}) &= c\\
\Pi_0^T a &= h_0.
\end{align*}

That our operator $(D, \Pi_0, \Pi)$ satisfies the invertibility assumption follows by the same separation of variables argument as in the irreducible case.

Lin's theorem, therefore, gives us a parameterization 
$$u(T, h_0, h): B_\eta(H_0 \oplus H) \to \Omega^1_{k,\delta}(Z^T)$$ 
of solutions to the above equations on $Z^T$ with $(h_0, h) \in B_\eta(H_0 \oplus H)$; recall that here $H_0$ consists of 1-forms of the form $dt \wedge \sigma$ for $\sigma \in \text{ker}(\Delta_A)$ and $$H = \mathcal K_A^- \oplus \Omega^0(Y; \mathfrak g_E) \oplus \mathcal K_A^+ \oplus \Omega^0(\overline Y; \mathfrak g_E).$$

For his application, Lin restricts to the subdomain $H_0 \oplus \mathcal K_A^- \oplus \mathcal K_A^+$; based on the discussion above, we see that in our case we should actually restrict to the subdomain $\mathcal K_A^- \oplus \mathcal K_A^+$ and ignore the terms coming from $H_0$. We abuse notation and write $u(T, h) = u(T, 0, h)$. As long as $T$ is finite, all solutions thus obtained are in Coulomb-Neumann gauge, and this defines a parameterization of a small neighborhood of $\gamma_A$ in $\widetilde{\mathcal M}_{A,k}(Z^T)$.

Before going on, we should explain why this map $u$ is equivariant; this follows as a consequence of uniqueness of solutions. That is, if $g$ is a gauge transformation preserving $\gamma_A$, then $u(gh_0, gh)$ is the unique solution $a$ to the equations

\begin{align*}d_{\gamma_A}^+ a + (a \wedge a)^+ + \left(dt \wedge \nabla_\pi(\gamma_A + a)\right)^+ &= 0 \\
d_{\gamma_A}^* a &= 0 \\
\Pi (a \big|_{\partial Z^T}) &= gc\\
\Pi_0^T a &= 0
\end{align*}

Here $c = (h, 0, 0)$ and $\Pi$ is the projection operator to 
$$\mathcal K^- \oplus \mathcal K^+ \oplus \Omega^0(Y \sqcup \overline{Y};\mathfrak g_E) \oplus H_0.$$ 
But this equation is also satisfied by $gu(T, h)$, because the projections are equivariant and the ASD equations are invariant under gauge transformations.

So $u(T,h)$ defines a map $$B_\eta(\mathcal K_\alpha) \to \mathcal{CN}_{A,k,\delta}(Z^T) \times \{p\} \subset \widetilde{\mathcal A}_{A,k,\delta}(Z^T)$$ for $T$ finite, where $p$ is a fixed framing; projecting to $\widetilde{\mathcal B}_{A,k,\delta}(Z^T)$, this gives a neighborhood of $\gamma_A$ in $\widetilde{\mathcal M}_{A,k,\delta}(Z^T)$.

So we have our parameterization of a neighborhood of $\gamma_A$ for all $T \in [T_0, \infty]$; we should check the desired properties of the restriction maps. The abstract gluing theorem guarantees certain properties of the maps $\mu^\pm_T: B_\eta(H) \to \text{ker}(d_A^*)_{k-1/2}$ obtained as the restriction of $\widetilde u(T, h)$ to the corresponding boundary component. However, our gluing lemma asks for properties of a restriction map \emph{whose definition includes the framing}, written $\widetilde \mu_T$. In the case of restriction to the leftmost boundary component, there is no change, as the framing is already fixed on the left boundary component; there is no need to perform any holonomy. But the right enpoint map is more complicated.

Let's be precise. Let $(\mathbf{A},p)$ be a framed connection on $Z^T$, where the framing is chosen at $(-T, b)$. Then the definition of the right endpoint map is $$(\mathbf{A},p) \mapsto \left(\mathbf{A}\big|_{\overline Y}, \text{Hol}^{(-T,b) \to (T,b)}_{\mathbf{A}} p\right);$$ we then pass to the quotient by the gauge group action (as this is an equivariant map).

Now, we identified a neighborhood of $A$ in $\mathcal B_{E,k-1/2}(\overline Y)$ with the projection of the Coulomb slice $\text{ker}(d_A^*) \times \{p\} \subset \widetilde{\mathcal A}_{E,k-1/2}(\overline Y)$, where $p$ is the same fixed framing as above. Therefore, identifying a neighborhood of $\gamma_A$ in $\widetilde{\mathcal M}_{A,k,\delta}(Z^T)$ with the corresponding subset of the Coulomb-Neumann slice (and a neighborhood of $A$ with the corresponding Coulomb slices), the restriction $\widetilde{\mu}_T^-$ to the left boundary component is still given by the same 
$$\mu_T^-: \mathcal{CN}_{A,k,\delta}(Z^T) \to \text{ker}(d_A^*)_{k-1/2}(Y)$$ 
as in the abstract gluing theorem, because our framing is chosen fixed at the left boundary component to begin wit. However, the restriction to the right boundary component is now given by 
$$\widetilde{\mu}_T^+(a) = \sigma^T(a) \cdot \mu^T(a) \in \text{ker}(d^*_A)_{k-1/2}(\overline Y),$$ 
where we have identified $\sigma^T(a) \in SO(3)$ with an $A$-parallel gauge transformation.

Precisely, the map $\sigma^T$ is given by $$\sigma^T(a) = \left(\text{Hol}_{\gamma_A + a}^{(-T,b) \to (T,b)}\right)^{-1},$$ using the same framing $p$ at $\pm T$ to identify this isomorphism $E_{(-T,b)} \to E_{(T, b)}$ with an element of $SO(3)$. We have $\sigma(0) = 1$: the connection $\gamma_A$ is in temporal gauge and so there is no holonomy across $\mathbb R \times \{b\}$.

What we need to do first is see that $\sigma^T(u(T, h))$ converges in the $C^\infty_{\text{loc}}$ topology. We will then use this to define the map $\widetilde u(\infty, -)$, as an appropriate modification of $u(\infty, 0, -)$.

First we recall how holonomy is defined. First, given $a \in \gamma_A + \Omega^1_{k,\delta}(Z^T)$, there is a natural restriction map to $$d/dt + \Omega^1_{k-3/2,\delta}([-T, T], \mathfrak g_b) = d/dt + dt \wedge \Omega^0_{k-3/2, \delta}([-T, T], \mathfrak g_b)$$ to the line $\mathbb R \times \{b\}$. The holonomy along the path only depends on the values of $a$ on this path, and more precisely, if we decompose $a = dt \wedge \eta(t) + \psi(t)$, it only depends on $\eta(b,t)$. Precisely, let $\gamma(t): [-T, T] \to SO(3)$ be the unique solution to the differential equation 
$$\gamma'(t) \gamma^{-1}(t) = \eta(b,t) \;\; \text{and} \;\; \gamma(-T) = 1,$$ 
guaranteed as long as $k-3/2 \geq 1$, so that $\eta$ is $C^1$. Solving this differential equation gives a smooth map 
$$\text{Lift}: L^2_{k-3/2,\delta}([-T, T], \mathfrak g_b) \to L^2_{k-5/2,\delta}([-T, T], SO(3)),$$ 
and then evaluation at $T$ is clearly smooth.

If we write $u_+(h)$ as the component of $u(\infty, h)$ on the cylinder $[0,\infty) \times Y$, and similarly for $u_-(h)$, we write 
$$U(T, h) = \tau_T^*u_+(h) + \tau_{-T}^* u_-(h).$$ 
We see from \cite[Lemma~3.5.19]{Lin} that $U(T, h) - \widetilde u(T, h)$ goes to $0$ in the $C^\infty_{\text{loc}}$ topology. So it suffices to check that $\sigma^T(U(T, h))$ has a limit in $C^\infty_{\text{loc}}.$

Now we may be very explicit. By translation, consider instead the interval $[0, 2T]$; we may see that $\sigma^T(U(T, h))^{-1}$ is the solution to the differential equation $$\gamma'(t) \gamma^{-1}(t) = u_+(h)(t,b) + u_-(h)(t-2T,b)\;\; \text{and} \;\; \gamma(0) = 1,$$ evaluated at $2T$. We may instead consider this as the product $g_2 g_1$ of two elements of $SO(3)$: first, $g_1$ is the time-$T$ value of that same equation; second, $g_2$ is the value at $T$ of the solution to the differential equation 
$$\gamma'(t) \gamma^{-1}(t) = u_-(h)(-t,b) + u_+(h)(2T-t,b) \;\; \text{and} \;\; \gamma(0) = 1.$$
That $g = g_2 g_1$ is just the statement that the holonomy of two paths, traversed in succession, is the product of their holonomies and that the holonomy of a path traversed in reverse is the inverse of the holonomy.

We examine $g_1$; the analysis of the other term is very similar. Notice that the stated differential equation only needs the values of the given functions on $[0,T]$, so we introduce a cutoff function $\beta$ which has $\beta(t) = 1$ for $t \leq 0$ and $\beta(t) = 0$ for $t \geq 1/2$, and $\beta_T(t) = \beta(t-T)$; write $$q_T(h) = u_+(h)(t,b) + \beta_T(t)\left(u_-(h)(t - 2T, b)\right).$$ It is clear that $U(T,h)$ and $q_T(h)$ give the same time-$T$ holonomy map, as they agree on $[0,T]$. Our first claim, to be proved below, is that the time-dependent map $B_\eta(\mathcal K_\alpha) \to L^2_{k-3/2, \delta}([0,\infty), \mathfrak g_b)$ given by sending $h \mapsto \beta_T(t) u_-(h)(t-2T, b)$ converges to $0$ in $C^\infty_{\text{loc}};$ given that this is true, our goal becomes finding the $C^\infty_{\text{loc}}$ limit as $T \to \infty$ of $$\text{ev}_{T}\text{Lift}\left(u_+(h)(t,b)\right).$$ Now there is nothing left to do here, as $u_+$ is time-independent! We see that the output converges in the $C^\infty_{\text{loc}}$ topology to $$\text{Hol}^{(0^+, b) \to (\infty, b)}_{\gamma_A + u_+(h)}.$$

Similarly we may identify the limit of $g_2$ as the inverse of the corresponding $\text{Hol}^{(0^-,b) \to (-\infty,b)}_{u_-(h)}$.

We conclude, therefore, that $\sigma^T(u(T,h))$ converges in $C^\infty_{\text{loc}}$ to $\sigma^\infty(h)$, which is defined to be the inverse of $$\left(\text{Hol}^{(0^-,b) \to (-\infty,b)}_{u_-(h)}\right) \left(\text{Hol}^{(0^+, b) \to (\infty, b)}_{\gamma_A + u_+(h)}\right).$$

This makes some intuitive sense; the appearance of $\sigma^T$ corresponds to the twist in framing as we move the left endpoint framing to the right endpoint. So on the infinite cylinder, we should be forced to take holonomy from left endpoint to right endpoint.

Based on the above, we define the map $$\widetilde u(\infty, -): B_\eta(\mathcal K_\alpha) \to \mathcal{CN}_{A,k,\delta}(Z^\infty) \times SO(3) \times SO(3) \subset \widetilde{\mathcal A}_{A,k,\delta}(Z^\infty)$$ to be $$\widetilde u(\infty, h) = (u(\infty, h), p, \sigma^\infty(h) \cdot p);$$ we project from this framed Coulomb slice to the configuration space $\widetilde{\mathcal M}(Z^\infty)$. Then the above discussion shows that indeed, $\widetilde{\mu}_T \to \widetilde{\mu}_\infty$ as $T \to \infty$.

Now to prove the claim that $\beta_T(t)\left(u_+(h)(t-2T, b)\right)$ converges to $0$ in $C^\infty_{\text{loc}}$. The point is that the application of the inverse function theorem defining $u_+$ can be factored through a Sobolev space with smaller exponent, so that one may find pointwise bounds on all of the derivatives of $u_+$: $$\|\mathcal D_m u_+(h)(t-2T, b)\| \leq C_m(h) e^{-\delta'(t-2T)},$$ where the $C_m$ are continuous functions of $h$ depending on $m$; the same is true when adding in the bump function (whose derivatives all have compact support).

Because the support of $\beta_T u_+(t-2T)$ is contained in $[0, T+1]$, we obtain bounds by $C_m(h) e^{- \delta'(1-T)}$; taking $T \to \infty$ we see that this goes to $0$ in $C^\infty_{\text{loc}}$.
\end{proof}

What remains is the case of $SO(2)$-reducible critical orbits. We use a mix of the above techniques; we work in the Coulomb-Neumann slice as a model for the normal space to $\gamma_A$, and then extend this to an $SO(3)$-invariant parameterization.

\begin{proof}[Proof of Lemma \ref{maingluinglemma} in the $SO(2)$-reducible case]
We study the same equations as above, where $\gamma_A$ is an $SO(2)$-reducible constant trajectory. We start by choosing $p \in S^2 = E_{(0,b)}/\Gamma_A$, a choice of framing \emph{modulo} the action of $\Gamma_A$; restrict to the subspace $\widetilde{\mathcal M}'(Z^T)$ given by those pairs $(\mathbf{A},q)$ with framings that project to $p$ in the quotient. (This is the space that has a neighborhood modelled by the solutions to the equations in the Coulomb-Neumann slice.) We find just as above that Lin's abstract gluing theorem again provides us with an $S^1$-equivariant parameterization $B_\eta(\mathcal K_\alpha) \to \widetilde{\mathcal M}(Z^T)$. The issue is in verifying that the restriction maps to the right end, which depend on holonomy, converge as $T \to \infty$ in $C^\infty_{\text{loc}}$, as before.

We may choose an arbitrary lift of $p$ to an actual framing $\tilde p$. Note that because $\gamma_A$ is in temporal gauge, the holonomy from $-T$ to $T$ is the identity. Now choosing the domain of the parameterization $u$ small enough, we may demand that $\text{Hol}^{(-T,b) \to (T, b)}_{\gamma_A + u(T, 0, h)} \tilde p$ lies in an $S^1$-invariant neighborhood of $\Gamma_A \cdot \tilde p$. Exponentiating the normal bundle, we may choose an $S^1$-equivariant diffeomorphism $U \cong S^1 \times D^2$. This gives us an $S^1$-equivariant projection map $P: U \to S^1$. We write $$\overline{\text{Hol}}^{(-T,b) \to (T,b)}_{\gamma_A + u(T,0,h)} = P\text{Hol}^{(-T,b) \to (T,b)}_{\gamma_A + u(T,0,h)}\tilde p;$$ this is precisely the amount of holonomy in the `$S^1$-direction', and the analogue of $\sigma^T(u)^{-1}$ for the $SO(2)$-reducible case (that is, the inverse of this term is precisely the rotation appearing in the right-most restriction map.)

We saw above that $\text{Hol}^{(-T,b) \to (T,b)}_{u(T, 0, h)}$ does converge in $C^\infty_{\text{loc}}$ as $T \to \infty$, and to $$\left(\text{Hol}^{(0^-, b) \to (-\infty, b)}_{u(\infty, 0, h)}\right)^{-1} \text{Hol}^{(0^+, b) \to (\infty, b)}_{u(\infty, 0, h)}.$$ Because the projection $P: U \to S^1$ is smooth, the same is true of $\overline{\text{Hol}}^{(-T,b) \to (T,b)}_{u(T, 0, h)}$; we use this to define the map $u(\infty, 0, h)$ on $B_\eta(\mathcal K_\alpha) \to \widetilde{\mathcal M}_{A,k,\delta}'(Z^\infty)$ as before, and see that it is $SO(2)$-equivariant.

All that is left is to extend this to an $SO(3)$-equivariant parameterization of the whole of $\widetilde{\mathcal M}_{A,k,\delta}(Z^T)$. We write the neighborhood of $(\gamma_A, \tilde p)$ as 
$$(p \cdot S^1) \times_{S^1} \mathcal{CN}_{\gamma_A, k, \delta}(Z^T).$$
That is, we allow ourselves to vary over all framings that lie above the fixed $p \in S^2$, and then quotient by the natural action of $\Gamma_A$. Then there is a canonical $S^1$-equivariant map 
$$S^1 \times_{S^1} \mathcal{CN}_{\gamma_A, k, \delta}(Z^T) \to SO(3) \times_{S^1} \mathcal{CN}_{\gamma_A, k, \delta}(Z^T).$$ 
This is our desired section. We thus extend our parameterization to an $SO(3)$-equivariant parameterization $SO(3) \times_{S^1} B_\eta(\mathcal K_\alpha) \to \widetilde{\mathcal M}(Z^T)$ of a neighborhood of the orbit through $\gamma_A$, as desired.
\end{proof}

With this, the proof of the gluing theorem follows essentially the same lines as in \cite[Chapters~19~and~24.7]{KMSW}. We outline the procedure and show what mild modifications are necessary. For the rest of this section, if $(W,\mathbf{E},\pi)$ is a manifold with cylindrical ends and perturbation $\pi$, we demand that $\pi$ is regular.

First, the moduli spaces $\widetilde{\mathcal M}_{E,z,\pi,k}(Z^T)$ of $\pi$-perturbed instantons on the cylinder, finite or infinite, are cut out transversely inside $\widetilde{\mathcal B}^e_{E,z,\pi,k}(Z^T)$, and therefore form a smooth Hilbert manifold no matter the perturbation $\pi$. The proof of this is no different than \cite[Theorem~17.3.1]{KMSW}; the point is that we have a unique continuation result for zeroes of the adjoint operator $Q^*_{\mathbf{A}, \pi}$, as the equations are of gradient-flow type.

Fix a compact manifold $W$ with boundary $Y_1 \sqcup \overline{Y_2}$, with cylindrical metric near the boundary, and $SO(3)$-bundle $\mathbf{E}$ which restricts to bundles $E_i$ over the respective ends $Y_i$. Suppose we have fixed small perturbations $\pi_i$ on $(Y_i, E_i)$ so that there are finitely many nondegenerate critical orbits of $\text{cs}_{Y_i}+f_{\pi_i}$. Further, we assume $(W,\mathbf{E})$ is weakly admissible, so it admits \emph{some} perturbation achieving transversality on all moduli spaces of bounded dimension on the infinite cobordism $\hat W$.

On the cylindrical ends of the infinite manifold $\hat W$, we fix the constant perturbation $(dt \wedge \nabla_{\pi_\pm}(\mathbf{A}))^+$. This is dampened by a cutoff function on the end 
$$[0, \infty) \times (Y_1 \sqcup \overline{Y_2}),$$ equal to $0$ for $x \leq 1$ and equal to $1$ for $x \geq 2$. There is a further interior holonomy perturbation supported in the complement of $[1, \infty) \times (Y_1 \sqcup \overline{Y_2})$. 

For $L \geq 2$, if $W^L$ is the complement of $(L,\infty) \times (Y_1 \sqcup \overline{Y_2})$, then the moduli space of instantons $\widetilde{\mathcal M}_{\mathbf{E},\pi,k}(W^L)$ on the compact manifold $W^L$ is a smooth Hilbert manifold in a neighborhood of $\widetilde{\mathcal M}_{\mathbf{E},k,\delta}(\hat W)$. This is a matter of showing that if $\mathbf{A}$ is a $\pi$-instanton on $\hat W$, then any nonzero element of the kernel of $Q^*_{\mathbf{A},\pi}$ on $W^L$ must restrict nontrivially to the boundary.\\

For suppose to the contrary that we had $\psi \in \Omega^{2,+}_{-k}(W^L;\mathfrak g_{\mathbf{E}})$ with $\psi|_{\partial W^L} = 0$. Elliptic regularity allows one to obtain $\psi \in L^2_k$ on $W^L \setminus W$, where the perturbation arises from the gradient of cylinder functions on slices. For the same reason, we can argue that $\psi|_{\partial W^L} = 0$ means by unique continuation that $\psi|_{\partial W^L \setminus W} = 0$. In particular, $\psi$ extends by $0$ to an element 
$$\hat \psi \in \Omega^{2,+}_{-k,\delta}(\hat W;\mathfrak g_{\mathbf{E}}), \;\; \text{with } Q^*_{\mathbf{A},\pi}\hat \psi = 0.$$

The assumption that $\pi$ is a regular perturbation then means that $\hat \psi = 0$ globally, and so $\psi = 0$, as desired. 

We write 
$$\widetilde{\mathcal M}_{\mathbf{E},k,\delta} \subset \widetilde V^L_{\mathbf{E},\pi,k} \subset \widetilde{\mathcal M}_{\mathbf{E},\pi,k}(W^L),$$
the middle term denoting an open set of $\pi$-instantons $\mathbf{A}$ on $W^L$ for which $Q_{\mathbf{A},\pi}: \Omega^1_k(W^L;\mathfrak g_{\mathbf{E}}) \to \Omega^2_{k-1}(W^L;\mathfrak g_{\mathbf{E}})$ is surjective. 

From here, we may describe moduli spaces on the manifold $W$ with infinite cylindrical ends as a fiber product. Let $I^+_j$ and $I^-_\ell$ be a finite sequence of intervals ($0 \leq j \leq n, 0 \leq \ell \leq m$), where the first negative interval is $I^-_0 = (-\infty, 0]$, the first positive interval is $I^+_0 = [0,\infty)$, and all other intervals are finite. Write $\widetilde{\mathcal B^\pm} = \widetilde{\mathcal B}^e_{E,k-1/2}(Y^\pm)$. There are evaluation maps $$\text{ev}_\pm: \widetilde{V}^L_{\mathbf{E},\pi} \to \widetilde{\mathcal B}^- \times \widetilde{\mathcal B}^+,$$ and similarly for $\widetilde{\mathcal M}_{\mathbf{E},\pi,k}(Z^{I^\pm_j}) \to \widetilde{\mathcal B}^\pm \times \widetilde{\mathcal B}^\pm$. For the infinite ends, we only have one evaluation map. We assemble all of this into a map $$(R_-, R_+): \widetilde{V}^L_{\mathbf{E},\pi} \times_{j=0}^n \widetilde{\mathcal M}_{E,\pi,k}(Z^{I^+_j}) \times_{\ell=0}^m \widetilde{\mathcal M}_{E,\pi,k}(Z^{I^-_{\ell}}) \to \left(\times_{j=1}^n\widetilde{\mathcal B}^+ \times_{\ell=1}^m \widetilde{\mathcal B}^-\right)^2.$$

\begin{lemma}\label{fibertest}Let $W$ be the cobordism with cylindrical ends attached. There is a natural map from the fiber product $\textup{Fib}(R_-, R_+)$ to $\widetilde{\mathcal M}_{\mathbf{E},\pi,k}(\hat W)$, which is a homeomorphism. Because any instanton $\mathbf{A}$ in the image is assumed to be cut out transversely, the corresponding element of $\widetilde{V}^L_{\mathbf{E},\pi} \times \prod \widetilde{\mathcal M}(Z^I)$, the map $R_- \times R_+$ is transverse to the diagonal.
\end{lemma}
We do not repeat the proof, which follows essentially as in \cite[Theorem~19.1.4]{KMSW}.

We now want a canonical way to cut up the real line (or multiple copies of the real line, in the case of broken trajectories) so that every instanton is written as a fiber product, as above, where most intervals are of a fixed large length $2L$ (where most of the energy is supported), and otherwise are of large variable length $T \geq T_0$. The mechanism for this is given in \cite[Section~19.2]{KMSW}. For $A \in \widetilde{\mathcal B}^e_{E_i,k-1/2}(Y_i)$, we write $c_i(A) = \|\nabla_A(\text{cs}+f_\pi)\|_{L^2}$.
Fix $\epsilon > 0$ so that any nontrivial instanton on $\mathbb R \times Y_i$ has $c_i(\mathbf{A}(t)) \leq \epsilon$ for some $t \in \mathbb R$. Let $\beta$ be a cutoff function equal to $1$ for $x \geq \epsilon$ and $0$ for $x \leq \epsilon/2$.

\begin{definition}Let $I = [-L, L]$. An instanton on $I \times Y$ is said to be \emph{centered} if
\begin{enumerate}
\item $c(\mathbf{A}(t)) \leq \epsilon/2$ for $t \in [-L, -L+1] \cup [L-1, L]$, \\
\item $c(\mathbf{A}(t)) \geq \epsilon$ for some $t \in [-L, L]$, and\\
\item The center of mass $$\int t \beta\big(c(\mathbf{A}(t))\big) dt \Big/ \int \beta\big(c(\mathbf{A}(t))\big) dt$$ is zero.
\end{enumerate}
\end{definition}

The space of centered instantons forms a smooth Hilbert manifold. Then we have the following analogue of \cite[Proposition~24.7.3]{KMSW}, with the same proof.

\begin{lemma}\label{finalbit}For any compact subset \begin{align*}K \subset \widetilde{\mathcal M}_{E_1, k, \delta}&(\alpha_0, \alpha_1) \times_{\alpha_1} \cdots \times \widetilde{\mathcal M}_{E_1, k, \delta}(\alpha_{n-1}, \alpha_n)\\
&\times_{\alpha_n} \widetilde{\mathcal M}_{\mathbf{E},k,\delta}(\alpha_n, \beta_0) \times_{\beta_0} \cdots \times_{\beta_{m-1}} \widetilde{\mathcal M}_{E_2, k, \delta}(\beta_{m-1}, \beta_m),\end{align*} there is an $L_0$ so that for all $\infty > L \geq L_0$, we have a neighborhood $K \subset V(K,L) \subset \overline{\mathcal M}_{\mathbf{E},k,\delta}(\alpha_0, \beta_m)$ for which \begin{enumerate}
\item For any $\mathbf{A} \in V(K,L)$, the restriction of $\mathbf{A}$ to $\partial W^L$ has $c(\mathbf{A}\big|_{\partial W^L}) \leq \epsilon/2$, \\
\item Any $\mathbf{A} \in V(K,L)$ admits a unique collection of cylinders of length $2L$ in the complement of $W^L$ so that $\mathbf{A}\big|_{I}$ is a centered instanton, while \\
\item The complement of both $W^L$ and these intervals consists of cylinders of length at least $T_0$, where $T_0$ is larger than the least $T$ for which Lemma \ref{maingluinglemma} applies to each of $Y_1$ and $Y_2$.
\end{enumerate}

As a corollary, we have a map $d: V(K,L) \to (0,\infty]^{n+m}$ measuring the distance between the centers of successive centered intervals; for the ends closest to $W$, we consider $\partial W$ to be the location of that corresponding center. It is easy to see that the values of $d$ do not depend on $K$ or $L$; in particular, taking the union over all $V(K, L)$, we obtain a map $V \to (0,\infty]^{n+m}$ from a neighborhood of the stratum in $\overline{\mathcal M}$, the neighborhood $V$ consisting of all instantons that may be partitioned as above.
\end{lemma}

At the same time, we have a map \begin{align*}&\text{ev}_T: (0,\infty]^{n+m} \times \widetilde{\mathcal M}_{E_1, k,\delta}^{(-\infty, 0] \times Y_1}(\alpha_0) \times_{i=1}^n B_\eta(\widetilde{\mathcal K}_{\alpha_i})\times^n \widetilde{\mathcal M}_{\text{cen}}(Z^L_1) \times \widetilde V^L_{\mathbf{E},\pi,k}\\
&\times_{j=0}^{m-1} B_\eta(\widetilde{\mathcal K}_{\beta_j}) \times^m \widetilde{\mathcal M}_{\text{cen}}(Z^L_2) \times \widetilde{\mathcal M}_{E_2,k,\delta}^{[0,\infty) \times Y_2}(\beta_m) \to \left(\times^n\widetilde{\mathcal B}_{E_1,k-1/2} \times^m \widetilde{\mathcal B}_{E_2, k-1/2}\right)^2,\end{align*} the map in every case given by restriction; perhaps most notably here, for the $B_\eta(\mathcal K_A)$ factors we are using the maps $\widetilde u(T, -)$ of Lemma \ref{maingluinglemma} to parameterize a neighborhood of $\gamma_A$ in $\widetilde{\mathcal M}_{\gamma_{\alpha_i},k,\delta}(Z^T)$, where $T \in [T_0, \infty]$, and then restricting to the boundary; Lemma \ref{maingluinglemma} included the fact that these restriction maps converge in $C^\infty_{\text{loc}}$ as $T \to \infty$.

We call the domain $(0,\infty]^{n+m} \times \mathcal M$ and the codomain $\mathcal N$. This map is smooth on each stratum of $(0, \infty]^{n+m}$, and for a convergent sequence $T_n \in (0,\infty]^{n+m}$, the maps $\text{ev}_{T_n}$ converge to the appropriate limit in the $C^\infty_{\text{loc}}$ topology. We may identify $\text{ev}_T^{-1}(\Delta)$ in the above fiber product as the neighborhood $V(L)$ in $\overline{\mathcal M}_{\mathbf{E},k,\delta}(\alpha_0, \beta_m)$. Now we may apply \cite[Lemma~19.3.3]{KMSW}.

\begin{theorem}\label{gluing-end}Let $W$ be a compact 4-manifold equipped with two boundary components, cylindrical metric near the boundary, an $SO(3)$-bundle $\mathbf{E}$ restricting to the pullback of fixed bundles over the ends, and a curve $\gamma$ which is cylindrical at the basepoints on the ends. Further let $W$ be equipped with a perturbation $\pi$ so that $\overline{\mathcal M}_{\mathbf{E}, \pi, z}(\alpha, \beta)$ is cut out transversely whenever the expected dimension is at most 10.

Then suppose we are given any open stratum $\sigma \subset \overline{\mathcal M}_{\mathbf{E}, \pi, z}(\alpha, \beta)$, where $\text{gr}_z(\alpha, \beta) \leq 10 - \dim \alpha$. In that case, we may find a neighborhood $V(\sigma)$ and an $SO(3)$-invariant map $V(\sigma) \to (0,\infty]^{n+m}$ which is stratum-preserving, a submersion on each stratum, and a topologically trivial fiber bundle for the inverse image of a neighborhood of $\infty$. In particular, $\overline{\mathcal M}_{\mathbf{E}, \pi, z}(\alpha, \beta)$ is a topological manifold with corners and a smooth structure on each stratum.
\end{theorem}

\begin{proof}Above we constructed a map $\text{ev}_T:  (0,\infty]^{n+m} \times \mathcal N \to \mathcal R$, where $\mathcal N$ and $\mathcal R$ are Hilbert manifolds; here $T \in (0,\infty]^{n+m}$. For fixed $T$, the map $\text{ev}_T$ is smooth, and $T \mapsto \text{ev}_T$ is convergent in the $C^\infty_\text{loc}$ topology because the same is true of the restriction maps $\mu_T$ on $B(\widetilde{\mathcal K}_{\alpha_i})$ and $B(\widetilde{\mathcal K}_{\beta_j})$. Thus the aforementioned lemma tells us that $\text{ev}^{-1}(\Delta) \to (0,\infty]^{n+m}$ is a topological submersion, which means precisely that $V(\sigma) \to (0,\infty]^{n+m}$ is a trivial fiber bundle near $\infty$. That the projection is smooth and a submersion on each stratum follows because $\mu_T$ is actually smooth in $T$ on each stratum.
\end{proof}

We summarize the results of this section as follows.

\begin{proposition}\label{gluing}For a regular perturbation $\pi$, and $\textup{gr}_z(\alpha, \beta) \leq 10 - \dim \alpha$ as in Corollary \ref{nobubbling}, $\overline{\mathcal M}_{E,z,\pi}(\alpha,\beta)$ can be given the natural structure of a compact topological $SO(3)$-manifold with corners with a smooth structure on each open stratum. We have the following decomposition of the boundary: $$\partial \overline{\mathcal M}_{E,z,\pi}(\alpha, \beta) = \bigsqcup_{\gamma; z_1 \ast z_2 = z} \overline{\mathcal M}_{E,z_1,\pi}(\alpha, \gamma) \times_\gamma \overline{\mathcal M}_{E,z_2,\pi}(\gamma, \beta).$$

The same is true on a cobordism $W$ as long as $\textup{gr}^W_z(\alpha, \beta) \leq 10 - \dim \alpha$. In that case, we have the decomposition
\begin{align*}\partial \overline{\mathcal M}^W_{\mathbf E,z,\pi}(\alpha, \beta) &= \bigsqcup_{\gamma \in \mathfrak C_{\pi_1}; z_1 \ast z_2 = z} \overline{\mathcal M}_{E_1,z_1,\pi_1}(\alpha, \gamma) \times_\gamma \overline{\mathcal M}^W_{\mathbf{E},z_2,\pi}(\gamma, \beta)\\
&\bigsqcup_{\zeta \in \mathfrak C_{\pi_2}; z_1 \ast z_2 = z} \overline{\mathcal M}^W_{\mathbf{E},z_1,\pi} (\alpha, \zeta) \times_\zeta \overline{\mathcal M}_{E_2,z_2,\pi_2}(\zeta, \beta).
\end{align*}
\end{proposition}

\section{Families of metrics and perturbations}\label{famsec}\label{sec:4d-families}
Let $S$ be a compact smooth manifold with corners. Let $(W,\mathbf{E})$ be a cylindrical end cobordism $(Y_1, E_1) \to (Y_2, E_2)$, and suppose each $(Y_i, E_i)$ is equipped with a regular metric and perturbation $(g_i,\pi_i)$. A family of metrics and perturbations parameterized by $S$, written $\pi_S$, is the data of:
\begin{itemize}
\item a smooth metric on the bundle $p_W^*TW$ on $S \times W$, restricting to the product metric $dt^2 + g_i$ on $\mathbb R \times Y_i$ on each end, and
\item a smooth map $S \to \mathcal P^{(4)}_c$, the latter being the space of perturbations on $W$ which agree with $\pi_i$ on the corresponding ends. We usually restrict to a small neighborhood $U_C$ of the perturbation $(\pi_-,\pi_+)$ with no interior part, as in Theorem \ref{trans3}.
\end{itemize}

We may then define the parameterized moduli spaces $\widetilde{\mathcal M}_{\mathbf{E},z,S}(\alpha, \beta)$ of pairs $(s,\mathbf{A})$, where $\mathbf{A}$ is a $\pi_s$-perturbed $L^2_{k,\delta}$ instanton going between critical orbits $\alpha$ and $\beta$ in a fixed component of trajectories $z$. Fiberwise compactifying by ideal instantons and broken trajectories we obtain $\overline{\mathcal M}_{\mathbf{E},z,S}(\alpha, \beta)$; the result is compact by a version of Proposition \ref{propertrajectories} which allows for variations of metric on the interior (which requires no change in the argument).

There is a bundle over $S \times \widetilde{\mathcal B}^e_{\mathbf{E}, z,k,\delta}(\alpha, \beta)$, written $\mathcal S_S$, whose fiber above $(s, \mathbf{A})$ is $\Omega^{2,+,s}_{k-1,\delta}(W;\mathfrak g_E)$; the notation $+,s$ indicates that we are taking the self-dual 2-forms with respect to the metric $g_s$.

\begin{definition}\label{regfamily}At each perturbed instanton $(s, \mathbf{A}) \in \widetilde{\mathcal M}_{\mathbf{E},z,S}(\alpha, \beta)$, taking the derivative of the $SO(3)$-equivariant map $\sigma: S \times \widetilde{\mathcal B}_{\mathbf{E}, z,k,\delta}(\alpha, \beta) \to \mathcal S$ defined by the instanton equation induces a map on the normal space to each orbit, $$(d\sigma)_{\mathbf{A}, s}: T_s S \times \textup{ker}(d_{\mathbf{A}}^{*,s})_{k,\delta} \to \Omega^{2,{+,s}}_{k-1,\delta}(W;\mathfrak g_E).$$

We say the parameterized moduli spaces $\widetilde{\mathcal M}_{\mathbf{E},z,S}(\alpha, \beta)$ are cut out regularly if $(d\sigma)_{\mathbf{A},s}$ is surjective for all perturbed instantons $(s,\mathbf{A})$ in the homotopy class $z$. In this case the parameterized moduli space $\widetilde{\mathcal M}_{\mathbf{E},z,S}(\alpha, \beta)$ is a smooth $SO(3)$-manifold of dimension $\textup{gr}_z(\alpha, \beta) + \dim S$.

We say the compactified parameterized moduli spaces $\overline{\mathcal M}_{\mathbf{E},z,S}(\alpha, \beta)$ are cut out regularly if every moduli space $\widetilde{\mathcal M}_{\mathbf{E},w,S}(\gamma, \gamma')$ appearing in the compactification is cut out regularly.

If $\overline{\mathcal M}_{\mathbf{E},z,S}(\alpha, \beta)$ is cut out regularly for all $(\alpha, \beta, z)$ with
$$\textup{gr}_z(\alpha, \beta) + \dim S \leq 10 - \dim \alpha,$$
we say that $\pi_S$ is regular.
\end{definition}

Recall Definition \ref{admiss-cob} of weakly admissible cobordism. In the first two cases, no reducibles arise for any perturbation, and still do not arise in families. In the third case, there is a unique reducible of each type for every perturbation, and in particular for the moduli spaces parameterized by a family of perturbations there is a $\dim S$-dimensional space of reducibles of each type.

In the final case, we have $b^+(W) > 0$ and $\mathbf{E}$ is nontrivial. In this case, for a \emph{generic} perturbation, there are no reducibles; the expected dimension of the space of reducibles is $-b^+(W)$. This means that they \emph{do} arise in generic families when $\dim S \geq b^+(W)$. Something else needs to be asserted to guarantee that in families, reducibles are still cut out transversely. Instead of adding more homological conditions (akin to $\rho$-monotonicity, but even more strict), we will simply endeavor to avoid the reducibles.

\begin{proposition}\label{famtrans}Suppose $(W,\mathbf{E})$ is equipped with a metric, a weakly admissible $SO(3)$-bundle and regular perturbations $\pi_\pm$ on the ends. There is an open set $(\pi_-, 0, \pi_+) \in U_C \subset \mathcal P^{(4)}_c$ of the space of perturbations which equal $\pi_\pm$ on the ends. The set $U_C$ is the small neighborhood around $(\pi_-, 0, \pi_+)$ given in Theorem \ref{trans3}, in which all fully reducible connections are cut out transversely and the parameterized space of reducibles $U_C \mathcal M^{\text{red}}$ is transverse to the space of Fredholm operators with non-trivial cokernel at those reducibles with energy $\leq C$. We may demand that $U_C$ is contractible by passing to a small ball around $(\pi_-, 0, \pi_+)$.

Suppose $\pi_0, \pi_1 \in U_C$ are regular perturbations.

If either $b^+(W) \neq 1$ or that $\mathbf{E}$ admits no $SO(2)$-reducible connections of class $L^2_{k,\delta}$ with $\pi_\pm$-flat limits, then there is a regular family $\pi_t$ of metrics and perturbations $\pi_t: [0,1] \to U_C$, all with the same values on the ends.

If $b^+(W) = 1$ and $\mathbf{E}$ supports $SO(2)$-reducibles, suppose $\pi_t: [0,1] \to U_C$ is a path of metrics and perturbations such that no $\pi_t$ supports a reducible instanton. Then we may modify the perturbations $\pi_t$ by an arbitrarily small amount on the interior of $[0,1]$ so that $\pi_t$ forms a regular family of metrics and perturbations.

If $K$ is any compact family of metrics on $W$, constant on the ends, then replacing the open set $U_C$ above with a smaller open set $U_{C,K}$, everything above still applies while allowing us to vary the metric in the family $K$.
\end{proposition}

\begin{proof}This follows from the same strategy as in Theorem \ref{trans3}. In the first case (that $b^+(W) \neq 1$ or $\mathbf{E}$ supports no reducibles), pick an arbitrary path $\pi_t: [0,1] \to U_C$; in the second case we are supplied such a path. By the assumption that $\pi_t \in U_C$, all fully reducible connections are cut out transversely for all $t$. We have a parameterized moduli space of reducibles, and the map
$$U_C \mathcal M^{\text{red}} \to U_C$$
has index $-b^+(W)$; we have already assumed that $\pi_0, \pi_1$ are regular values, so by a small perturbation we may assume that $\pi_t$ is transverse to this projection. That is to say that for the family $\pi_t$, the reducibles are cut out transversely inside the reducible locus. When $b^+(W) > 1$, this means there are no reducible $\pi$-ASD connections; when $b^+(W) = 1$, we already assumed there were no reducible $\pi$-ASD connections; when $b^+(W) = 0$, there is a finite set of $\pi_t$-ASD connections for each $t$.

Now as in Theorem \ref{trans3} we have a countable union of submanifolds of even codimension
$$U_C\mathcal M^{\text{red},\geq 1} \subset U_C \mathcal M^{\text{red}},$$
restricting to those reducible connections of energy $\leq C$. This is the space of reducible connections supporting a splitting $\mathbf{E} \cong \mathbb R \oplus \eta$ so that the normal ASD operator
$$D^{\text{irred}}_{\mathbf{A},\pi}: \Omega^1_{k,\delta}(W;\eta) \to \Omega^{2,+}_{k-1,\delta}(W;\eta)$$
has non-trivial cokernel. The codimension of each of these is even because the normal space is isomorphic to $\text{Hom}_{\mathbb C}(\mathbb C^{i+j}, \mathbb C^i)$ for $i \geq 0, j > 0$, which is a non-trivial complex vector space; the fact that the index $i$ is positive for each reducible $\mathbf{A}$ is part of the assumption that $W$ is weakly admissible.

By assumption, $\pi_0$ and $\pi_1$ are regular values of the projection
$$U_C\mathcal M^{\text{red},\geq 1} \to U_C,$$
so by a small perturbation we may assume that $\pi_t$ is transverse to this projection as well. Because the codimension of each manifold in $U_C\mathcal M^{\text{red},\geq 1}$ is at least 2, this means that for any reducible $\pi_t$-instanton $\mathbf{A}$, the normal ASD operator is surjective.

That is to say, for the path $\pi_t$, any reducible $\pi_t$-instanton $\mathbf{A}$ of energy $\leq C$ is cut out transversely, and not just internally to its own locus.

Now we conclude by playing the same game with the irreducible connections: we have a parameterized moduli space $U_C \mathcal M^{\text{irred}}$, and we demand that the path $\pi_t$ is transverse to the projection
$$U_C \mathcal M^{\text{irred}} \to U_C,$$
which we can make true by a small perturbation on the interior of $[0,1]$.

The set $U_C$ depends on the chosen metric on $W$, but as in the statement of the proposition, we may vary over a compact set of metrics by passing to a smaller open set $U_{C,K}$.
\end{proof}

\begin{remark}Except in the case where $\mathbf{E}$ admits no reducible connections whatsoever, there are obstructions to extending the above theorem to arbitrary families of metrics and perturbations, even when $b^+ = 0$. Generically, one expects the set of reducibles for which the operator $Q^{\text{irred}}_{\mathbf{A},\pi}$ has nontrivial cokernel to form a codimension $I+2$ family, where $I$ is the index of this operator. The definition of regular family of metrics and perturbations does nothing to help this: the map we want to be surjective, $$\widehat{\nabla}_{s'(\pi)} \mathbf{A} : T_s S \to \text{coker}(Q^{\text{irred}}),$$ is an \emph{equivariant} map. The former space has trivial circle action, and the latter space has circle action of weight 1: the only equivariant map is zero. So if our family of metrics and perturbations supports a reducible which is not cut out transversely in the usual sense, it won't be cut out transversely in the family sense, either.

When $\mathbf{E}$ admits no reducible connections, there is no difficulty applying the usual ideas: there is a smooth map $\pi: \mathcal P \mathcal M \to \mathcal P$, and we choose a map $S \to \mathcal P$ transverse to $\pi$ extending a given map $\partial S \to P$ transverse to $\pi$.
\end{remark}

\begin{remark}\label{interpolate}Suppose we are given metrics $g_0, g_1$ on $W$, which give open sets $U_{C,i}$ depending on the metric; and suppose we are given regular $\pi_i \in U_{C,i}$. Choosing a path $g_t$ between these metrics, we see that we may choose a generic path between the two perturbations, so long as they are chosen in a smaller set $U_{C,g_t} \subset U_{C,0} \cap U_{C,1}$.

To actually construct a generic path between the perturbations, then, first choose a perturbation $\pi_i' \in U_{C,g_t}$, and choose a generic path from $\pi_i$ to $\pi_i'$ inside of $U_{C,i}$. Then choose a generic path from $\pi_0'$ to $\pi_1'$ inside of $U_{C,g_t}$.
\end{remark}

Assuming $S$ is a regular family of metrics and perturbations, we have smooth moduli spaces $\widetilde{\mathcal M}_{\mathbf{E},z,S}(\alpha, \beta)$ with smooth projection maps $\widetilde{\mathcal M}_{\mathbf{E},z,S} \to S$. The space $\widetilde{\mathcal M}_{\mathbf{E},z,S}$ is not compact but the fiberwise Uhlenbeck compactification $\overline{\mathcal M}_{\mathbf{E},z,S}$ is. As long as
$$\text{gr}^W_z(\alpha, \beta) + \dim \alpha + \dim S \leq 10,$$
this being at least the dimension of $\widetilde{\mathcal M}_S/SO(3)$, there is no Uhlenbeck bubbling arising in this compactification. Note that here it is possible for $\widetilde{\mathcal M}_{\mathbf{E},z,S}$ to be nonempty even if $\overline{\text{gr}}_z(\alpha,\beta)$ is negative. In this case, the set of $s \in S$ so that $\pi_s$ is nonempty is a compact submanifold with corners of codimension $-\overline{\text{gr}}^W_z(\alpha, \beta)$.

We then have the following analogue of Theorem \ref{gluing-end}.

\begin{proposition}\label{famgluing}If $\pi_S$ is a regular family of metrics and perturbations indexed by $S$, and $\textup{gr}^W_z(\alpha, \beta) + \dim S \leq 10 - \dim \alpha$, then $\overline{\mathcal M}_{\mathbf{E},z,S}$ has the natural structure of a compact topological $SO(3)$-manifold with corners and a smooth structure on each stratum. Its boundary has the decomposition
\begin{align*}\partial \overline{\mathcal M}_{\mathbf{E},z,S}(\alpha, \beta) = &\bigcup_{\substack{\gamma \in \mathfrak C_{Y_1}\\ z_1 * z_2 = z}}\overline{\mathcal M}_{E_1,z_1,\pi}(\alpha, \gamma) \times_\gamma \overline{\mathcal M}_{\mathbf{E},z_2,S}(\gamma, \beta)\\ &\bigcup_{\substack{\zeta \in \mathfrak C_{Y_2} \\ z_1 * z_2 = z}} \overline{\mathcal M}_{\mathbf{E},z_1,S}(\alpha, \zeta) \times_\zeta \overline{\mathcal M}_{E_2,z_2,\pi}(\zeta, \beta) \;\cup \; \overline{\mathcal M}_{\mathbf{E},z,\partial S}(\alpha, \beta).
\end{align*}

In particular, suppose $S = [0,1]$, and write $W$ as the composite of cobordisms with cylindrical ends $W_1$ and $W_2$, obtained by spliting together the positive end of $W_1$ with the negative end of $W_2$. Write the broken perturbation $\pi(1)$ as $\pi_i$ on each of the $W_i$, where $\pi_-$ is the perturbation on the negative end of $W_1$ and $\pi_+$ is the perturbation on the positive end of $W_2$; write the common perturbation in the center as $\pi_b$. So long as $\textup{gr}_z^W(\alpha, \beta) \leq 9-\dim \alpha$, we have a decomposition
\begin{align*}\partial \overline{\mathcal M}_{\mathbf E,S,z}(\alpha, \beta) = \bigcup_{\substack{\gamma \in \mathfrak C_{\pi_-}\\ z_1 \ast z_2 = z}} &\overline{\mathcal M}_{E_-,z_1,\pi_-}(\alpha, \gamma) \times_\gamma \overline{\mathcal M}_{\mathbf{E},z_2,S}(\gamma, \beta)\\
\bigcup_{\substack{\eta \in \mathfrak C_{\pi_b}\\ z_1 \ast z_2 = z}} &\overline{\mathcal M}_{\mathbf{E}_1,z_1,\pi_1}(\alpha, \eta) \times_\eta \overline{\mathcal M}_{\mathbf{E}_2,z_2,\pi_2}(\eta, \beta)\\
\bigcup_{\zeta \in \mathfrak C_{\pi_+}; z_1 \ast z_2 = z} &\overline{\mathcal M}_{\mathbf{E},z_1,S}(\alpha, \zeta) \times_\zeta \overline{\mathcal M}_{E_+,z_2,\pi_+}(\zeta, \beta)\\
\cup \; \; \; \; &\overline{\mathcal M}_{\mathbf{E},z,\pi(0)}(\alpha, \beta).\\
\end{align*}
\end{proposition}

\begin{proof}Fix an open stratum $\sigma \subset \overline{\mathcal M}_{\mathbf{E},z,S}(\alpha, \beta)$ corresponding to $n$ breakings along $Y_1$ and $m$ breakings along $Y_2$. The result will follow essentially as in Theorem \ref{gluing-end}. There we described a gluing map $$\text{ev}_T: (0,\infty]^{n+m} \times \mathcal N \to \mathcal R$$ between Hilbert manifolds. The manifold $\mathcal N$ consisted of products of spaces parameterizing $n$ long (length $L$) broken pieces of framed instantons on $[-L, L] \times Y_1$, and similarly $m$ on $(Y_2, E_2)$, a space parameterizing instantons on cylinders $[-T, T] \times Y_i$ for all $T_0 \leq T \leq \infty$ sufficiently large, and a space parameterizing instantons on the compact manifold $W^L$ (which is $W$ with the cylindrical ends $[0, \infty) \times Y$ truncated to $[0, L] \times Y$). The codomain $\mathcal R$ consisted of restrictions of these to $2n$ copies of the configuration space of framed connections on $Y_1$ and $2m$ copies of the configuration space of framed connections on $Y_2$. This is smooth on each stratum of $(0,\infty]^{n+m}$, and the individual maps $\text{ev}_T$ converge in $C^\infty_{\text{loc}}$ as
$$T_k \to T \in (0,\infty]^{n+m}.$$
The moduli space of instantons on $W$, compactified by broken trajectories, is identified with $\text{ev}^{-1}(\Delta)$, the subset of instantons on each piece so that the restrictions to corresponding boundary pieces agree.

When we allow the metric and perturbation to vary in some family $S$, this variation occurs in a compact part of the manifold with cylindrical ends $\hat W$. As long as $L$ above is taken large enough, this variation only takes place in $W^L$, and hence only affects that term of the product; we may replace the term
$$\widetilde{\mathcal M}_{\mathbf{E},\pi,k,\delta}(\hat W) \subset \widetilde V^L_{\mathbf{E},\pi} \subset \widetilde{\mathcal M}_{\mathbf{E},\pi}(W^L)$$
instead with a parameterized version $\widetilde{V}^L_{\mathbf{E},S}$ of instantons on $W^L$ for which
\begin{align*}\Omega^1(W^L) \oplus T_s S &\to \Omega^{2,+}(W^L)\\
(\omega, \pi') &\mapsto D_{\mathbf{A},\pi} + \widehat{\nabla}_{\pi'}(\mathbf{A})\end{align*}
is surjective.

The modified evaluation map above is written
$$(0,\infty]^{n+m} \times \mathcal N_S \to \mathcal R,$$
and is again smooth on each stratum of $(0,\infty]^{n+m}$, with convergent sequences of $T_k \in (0,\infty]^{n+m}$ giving sequences of maps $\text{ev}_{T_k}$ which converge in $C^\infty_{\text{loc}}$. Again we have that $\overline{\mathcal M}_{\mathbf{E},z,S} = \text{ev}^{-1}(\Delta)$.

Therefore we may apply \cite[Lemma~19.3.3]{KMSW}, which asserts that
$$\text{ev}^{-1}(\Delta) \to (0,\infty]^{n+m}$$
is a topological submersion, meaning that it is a topologically trivial fiber bundle near $\infty$. In particular, this gives a chart for a neighborhood of $\sigma \subset \overline{\mathcal M}_{S,z}$, an open embedding $(0,\infty]^{n+m} \times \sigma \to \overline{\mathcal M}_{S,z}$ that is the identity on $\sigma$, is stratum-preserving, and is a diffeomorphism on each open stratum.
\end{proof}

We will need a version of these results for more general families of metrics and perturbations termed `families of broken metrics' (see \cite{Dae} and \cite{KM3}, where this notion was used to study spectral sequences from Khovanov homology).

Suppose are given a 4-manifold $W$ with cylindrical ends $(-\infty, 0] \times Y_1$ and $[0, \infty) \times Y_2$, and a sequence of Riemannian 3-manifolds $M_1, \cdots, M_\ell$, equipped with an isometric embedding $$\sqcup \varphi_i: (-1-\epsilon, 1+\epsilon) \times M_i \to W.$$ Each interval $[-1, 1] \times M_i$ will be equipped with a perturbation $\widehat{\nabla}_{\pi_i}$ which restricts to some neighborhood of $(-c_i, c_i) \times M_i$ as $(dt \wedge \nabla_{\pi_i})^+$ for some regular perturbation $\nabla_{\pi_i}$ on the pair $(M_i, \mathbf{E}_i)$ (that is, so that the critical orbits are cut out nondegenerately and all moduli spaces up to some sufficiently large dimension are also cut out nondegenerately). We say that the perturbation is adapted to the cut.

To define the restriction map for the framed instanton moduli spaces, we need to assume that the path $\gamma: \mathbb R \to W$ is adapted to the cuts, in the sense that there is a sequence $t_1, \cdots, t_\ell$ of real numbers with $t_i +2 < t_{i+1}$ so that for the interval $t \in [t_i-1, t_i+1]$, we have $\gamma(t_i+t) = \varphi_i(t-t_i, x_i)$ for some fixed $x_i \in M_i$. In particular, $\gamma$ passes through the $M_i$ in the order listed.

The canonical family of broken metrics on $W$ is parameterized by $[0, \infty]^\ell$; for a fixed element $(T_1, \cdots, T_\ell)$, one replaces the isometric copies of $[-1, 1] \times M_i$ with isometric copies of
\begin{align*}[-1-T_i, 1+T_i] \times M_i \;\; &\text{if} \;\; T_i < \infty, \text{ or}\\
\Big([-1,\infty) \sqcup (-\infty, 1]\Big) \times M_i \;\; &\text{if} \;\; T_i = \infty,
\end{align*}
and replaces the image of $\gamma$ appropriately. The perturbation $\widehat{\nabla}_{\pi_i}$ associated to $(T_1, \cdots, T_\ell)$ is defined by assuming that its value on the interval $(-c_i - T_i, c_i+T_i)$ is equal to $(dt \wedge \nabla_{\pi_i})^+$ as above; this remains true in the infinite case. We say that the metric parameterized by $(\infty, \cdots, \infty)$ is an \emph{$\ell$-times broken metric on $W$, with cuts along $M_1, \cdots, M_\ell$}.

A \emph{family of broken metrics and perturbations} is given by a compact smooth manifold with corners $S$ equipped with a family of metrics and perturbations on $W$, possibly containing broken metrics and perturbations. If $s \in S$ corresponds to one of these broken metrics, cut along $M_1, \cdots, M_\ell$, then we demand there is a chart $\varphi: (T, \infty]^\ell \times U \to S$, where $U$ is some topological manifold with corners and $\varphi$ is a stratum-preserving open embedding, so that $\varphi(T_1, \cdots, T_\ell, x)$ is the element of the canonical family of broken metrics and perturbations associated to $(T_1, \cdots, T_\ell)$ described above. The variance in metric and perturbation parameterized by $x \in U$ occurs in the complement of the canonical intervals $I \times M_i$. The union over the images of $\varphi(U)$ in $S$, as the $\varphi$ vary over charts for neighborhoods of a specific sequence of cuts $M_1, \cdots, M_\ell$, is called the \emph{cut stratum for $M_1, \cdots, M_\ell$}.

If $S$ parameterizes a family of broken metrics and perturbations, we define the parameterized moduli space $\widetilde{\mathcal M}_S$ as follows. For a fixed broken metric and perturbation $\pi$, cut along $M_1, \cdots, M_\ell$, we are left with a sequence $W_0, \cdots, W_\ell$ of manifolds with two cylindrical ends; other than the outer two, we have $W_i$ a manifold with incoming cylindrical end $(-\infty, 0] \times M_i$ and outgoing cylindrical end $[0,\infty)\times M_{i+1}$. Each of these has a corresponding collection of moduli spaces $$\bigsqcup_{\substack{\alpha_i, \alpha_{i+1} \in \mathfrak C_{\pi_i}\\ z_i \in \pi_1(\widetilde{\mathcal B}^e_E(\alpha_i, \alpha_{i+1}))}}\widetilde{\mathcal M}^{W_i}_{\mathbf{E},z_i,\pi}(\alpha_i, \alpha_{i+1});$$ we say that an element of $\widetilde{\mathcal M}^W_{\mathbf{E},z,\pi}$ is a sequence of framed instantons
$$\mathbf{A} \in \widetilde{\mathcal M}^{W_i}_{\mathbf{E},z_i,\pi}(\alpha_i, \alpha_{i+1})$$
so that the concatenation $z_1 * \cdots * z_\ell$ is $z$. That is, writing the space given by the disjoint union of $\pi_i$-critical orbits on $M_i$ as $\widetilde R_{\pi_i}(M_i)$, we have

$$\widetilde{\mathcal M}^W_{\mathbf{E},\pi} := \widetilde{\mathcal M}^{W_0}_{\mathbf{E},\pi} \times_{\widetilde R_{\pi_1}(M_1)} \cdots \times_{\widetilde R_{\pi_\ell}(M_\ell)} \widetilde{\mathcal M}^{W_\ell}_{\mathbf{E},\pi},$$
which then may be decomposed as a disjoint union according to the paths $(z_1, \cdots, z_\ell)$.

The moduli space of the family is then defined to be $$\widetilde{\mathcal M}_{\mathbf{E},z,S}(\alpha_0, \alpha_n) = \cup_{s \in S} \widetilde{\mathcal M}_{\mathbf{E},z,\pi_s}(\alpha_0, \alpha_n),$$ topologized so that $(s_n, \mathbf{A}_n)$ converges to $(s,\mathbf{A})$ if $s_n \to s$ and $\mathbf{A}_n \to \mathbf{A}$ in $L^2_k$ on compact sets.

We say that the individual moduli space $\widetilde{\mathcal M}_{\mathbf{E},z,\pi}(\alpha_0, \alpha_n)$ is \emph{regular} if each component moduli space $\widetilde{\mathcal M}_{\mathbf{E}_i,z_i,\pi}(\alpha_i, \alpha_{i+1})$ is cut out regularly in dimensions at most $10$.

We say that a family of broken metrics and perturbations parameterized by $S$ is \emph{regular} if, for each cut stratum $\sigma$ of $S$ parameterizing $\ell$-broken metrics with cuts along $M_1, \cdots, M_\ell$, but varies metric and perturbation elsewhere on $W$, we have that each of the component moduli spaces $\widetilde{\mathcal M}_{\mathbf{E}_i,z_i,\sigma}(\alpha_i, \alpha_{i+1})$ is regular. For $\ell = 0$, this just says that the open submanifold of $S$ parameterizing unbroken metrics should be regular in the sense of Definition \ref{regfamily}.

One may then define the fiberwise Uhlenbeck compactification $\overline{\mathcal M}_S$ as before, by setting $\overline{\mathcal M}_s$ to be the fiber products $$\overline{\mathcal M}_{\mathbf{E}_0,\pi_0} \times_{\widetilde R_{\pi_1}(M_1)} \cdots \times_{\widetilde R_{\pi_\ell}(M_\ell)} \overline{\mathcal M}_{\mathbf{E}_\ell,\pi_\ell}$$ of Uhlenbeck compactifications (by ideal instantons and broken trajectories); $\overline{\mathcal M}_{\mathbf{E},\pi,S}$ is the union of these over $s \in S$, topologized as above. One easily verifies that $\overline{\mathcal M}_S$ is compact by the same argument as the case of families of unbroken metrics; now if $s_n \to s$ is a sequence that limits to a metric and perturbation with more breakings, energy slides off the newly-infinite ends. If $\overline{\mathcal M}_S$ is regular and of dimension at most $10$, then this compactification only includes broken trajectories, not ideal instantons.

On the composite cobordism $W$, we will want to consider the composite perturbation, but this involves a term along the neck which does not come from what we call `interior perturbations'. We will need to dampen this perturbation out if we hope to compare the moduli spaces on the composite, with generic perturbation, to those of the individual manifolds $W_i$.

We resolve this by considering a 2-parameter family of perturbations. Suppose $S$ is a family of broken metrics and perturbations; let $(\pi, t) \in S \times [0,1]$ denote the perturbation which is the same as $\pi$, but instead of being $\beta_0(s) \widehat{\nabla}_{\pi_i}(\mathbf{A})$ on the neck, we instead use $t\beta_0(s) \widehat{\nabla}_{\pi_i}(\mathbf{A})$; for $t = 1$, this is what we already had, and for $t = 0$, this is identically zero.

To see why this is necessary or desirable, recall where the open set $U_C$ of Theorem \ref{trans3} comes from: it is a set of perturbations $\pi$ so that the map
$$D^{\text{irred}}: U_C\mathcal M^{\text{red}} \to U_C \mathcal F$$
sending a reducible $\pi$-instanton to its normal ASD operator, is transverse to those operators with non-trivial cokernel . This open set is guaranteed because we see that $D^{\text{irred}}$ is transverse to $U_C \mathcal F$ along the perturbations $(\pi_-, 0, \pi_+)$ with no part coming from `interior holonomy perturbations'.

The point is the following. Suppose $W$ is a manifold with two cylindrical ends, arising from gluing the $W_i$ together. Suppose further that $\pi_-, \pi_1, \cdots, \pi_+$ is a sequence of regular perturbations on the `cutting' 3-manifolds, and we consider the perturbation $(\pi_-, 0, t\pi_1, 0, \cdots, t\pi_n, 0, \pi_+)$ on the composite, where we dampen the perturbation on the necks as above. Then the same argument that shows $D^{\text{irred}}$ is transverse along $(\pi_-, 0, \pi_+)$ shows that $D^{\text{irred}}$ is transverse just as well along $(\pi_-, 0, \cdots, 0, \pi_+)$ - this is guaranteed so long as there is no `interior holonomy perturbation' component.

Given this, suppose we have a family of broken metrics and perturbations $S$; we will demand that the perturbation on the necks takes a `dampened' form such as the above. Then if $\mathcal P^{(4)}_{c,i}$ are the spaces of interior perturbations on $W_i$, there is an open set $U_{C,S,i} \subset \mathcal P^{(4)}_{c,i}$ so that $D^{\text{irred}}$ is transverse to the operators with non-trivial cokernel along any
$$(\pi_-, \pi_0^{\text{int}}, \pi_1, \cdots, \pi_n^{\text{int}}, \pi_+)$$ with
$$\pi_i^{\text{int}} \in U_{C,S,i}$$ for all $i$.

The reason that this is an open condition is the following. Suppose $\mathbf{A}_n$ is a sequence of instantons on $W^L$, with gluing factors $L$ finite but going to infinity, with $\mathbf{A}_n$ converging to a broken instanton $\mathbf{A}$. If all of the components of $\mathbf{A}$ are regular, then so is $\mathbf{A}_n$, for large $n$.

\begin{proposition}\label{famcut}Let $S$ parameterize a family of broken metrics, and $\pi: S \to \prod_i U_{C,S,i}$ parameterize a family of perturbations on $(W,\mathbf{E})$ so that $\partial S$ is a regular family.

Then the family of broken metrics on $S$ only consists of unbroken metrics and perturbations on the interior of $S$, and one may modify $S$ by an arbitrarily small amount on the interior to make it into a regular family. In particular, given a broken metric and perturbation $\pi_0 \in \prod_i U_{C,S,i}$, there is a regular path $\pi_t$ of metrics and perturbations, with $\pi_t$ unbroken for $t > 0$, and so that the perturbation $\pi_1$ a regular perturbation on the composite cobordism $W$ with zero perturbation on the `neck'.

Now suppose that $S$ parameterizes a regular family of broken metrics and perturbations on $(W,\mathbf{E})$. Let $\sigma \subset S$ be a cut stratum of codimension $\ell$, so parameterizing metrics with a sequence of cuts along $M_1, \cdots, M_\ell$. For any stratum $\widetilde \sigma \subset \overline{\mathcal M}_{\mathbf{E},\sigma,z}(\alpha, \beta)$ of codimension $n + k_1 + \cdots + k_\ell + m$, parameterizing broken $s$-perturbed instantons for $s \in \sigma$ with $n$ breakings along the incoming end $Y_1$, with $k_i$ breakings along the internal end $M_i$, and with $m$ breakings along the outgoing end $Y_2$. Then there is a neighborhood $$V(\sigma) \subset \overline{\mathcal M}_{\mathbf{E},z,S}(\alpha, \beta)$$ and an $SO(3)$-invariant map $V(\sigma) \to (0,\infty]^{n+k_1 + \cdots + k_\ell +m}$ which is stratum-preserving, a submersion on each stratum, and a topologically trivial fiber bundle for the inverse image of a neighborhood of $\infty$. In particular, $\overline{\mathcal M}_{\mathbf{E},S,z}(\alpha, \beta)$ is a compact topological manifold with corners and a smooth structure on each stratum, with the same decomposition as in Proposition \ref{famgluing}.
\end{proposition}

\begin{proof}
The transversality claim follows without real change as in Proposition \ref{famtrans}. By the discussion before this proposition, if we have a broken perturbation $\pi_0$, we may choose a path $\pi_t$ so that $\pi_1$ is of the form $(\pi_-, \pi, \pi_+)$, where the compactly supported perturbation $\pi$ is an interior holonomy perturbation --- it has nothing corresponding to the holonomy perturbation on the neck.

The rest follows the same lines as the proof of the gluing theorem: first, we may use Lemma \ref{maingluinglemma} to build charts $\mu^i_T: B_{\eta_i}(\widetilde H_i) \to \widetilde{\mathcal M}_{\alpha_i, k, \delta}(Z_i^T)$ near the constant solution $\gamma_{\alpha_i}$, where $Z_i^T = [-T, T] \times M_i$ for finite $T$ and $$\big([0, \infty) \cup (-\infty, 0]\big) \times M_0$$ for $T = \infty$. For uniformity of notation, set $M_0 = Y_1$ and $M_{\ell+1} = Y_2$. Then as in the proof of Lemma \ref{finalbit}, for sufficiently large fixed $L$ may assemble a map with domain the product

\begin{align*}(0,\infty]^{n+k_1 + \cdots + k_\ell + m} &\times \widetilde{\mathcal M}^{(-\infty \times Y_1}_{E_1, k, \delta}(\alpha_0) \times^n B_{\eta_0}(\widetilde H_0) \times^n \widetilde{\mathcal M}_{\text{cen}}(Z^L_0)\\
&\times \widetilde{V}^{W'_0}_{\sigma,\mathbf{E},k} \times^{k_1} B_{\eta_1}(\widetilde H_1) \times^{k_1} \widetilde{\mathcal M}_{\text{cen}}(Z^L_1)\\
&\;\;\;\;\;\cdots\\
&\times_\sigma \widetilde{V}^{W'_\ell}_{\sigma,\mathbf{E},k} \times^{k_\ell} B_{\eta_\ell}(\widetilde H_\ell) \times^{k_\ell} \widetilde{\mathcal M}_{\text{cen}}(Z^L_\ell)\\
&\times^m B_{\eta_{\ell+1}}(\widetilde H_{\ell+1}) \times^m \widetilde{\mathcal M}_{\text{cen}}(Z^L_{\ell+1}) \times \widetilde{\mathcal M}^{[0,\infty) \times Y_2}_{E_2,k,\delta}(\beta).
\end{align*}

Recall that $\widetilde{V}^W_{\mathbf{E}, k} \subset \widetilde{\mathcal M}^{W'}_{\mathbf{E}, k}$ is an open neighborhood of those instantons which extend to the whole of $W$.

Here the $\times_\sigma$ indicates that the moduli spaces $$\widetilde{V}^{W'_0} \times \widetilde{V}^{W'_1} \times \cdots \times \widetilde{V}^{W'_\ell}$$ should then be considered as those moduli spaces altogether parameterized by $\sigma$, but not that we parameterize each by $S$ separately and take the product. Secondly, here we mean by $W_i$ the compact manifolds given by trunacting the cylindrical ends of $W'_i$, written $\text{End}(W_i) = [0, \infty) \times \partial W'_i$, at some large finite $[0, L]$.

The codomain of the map is restriction to the boundary of each piece (this uses parallel transport along the path $\gamma$):

$$\left(\times^n \widetilde{\mathcal B}^e_{Y_1, k-1/2} \times^{k_1} \widetilde{\mathcal B}^e_{M_1, k-1/2} \times \cdots \times^{k_\ell} \widetilde{\mathcal B}^e_{M_\ell, k-1/2} \times^m \widetilde{\mathcal B}^e_{Y_2,k-1/2}\right)^2.$$

Call the  domain $(0, \infty]^{n+k_1 + \cdots + k_\ell + m} \times \mathcal N$ and the codomain $\mathcal R$. Mapping to an element of the diagonal in the codomain means that the corresponding instantons on each piece may be glued to a (possibly broken) instanton on the sequence of manifolds with cylindrical ends $W_0, \cdots, W_\ell$. The assumption that $\partial S$ is regular means, in particular, that $\sigma$ is a regular family of metrics and perturbations, and so the map $\text{ev}_\infty: \mathcal N \to \mathcal R$ is transverse to the diagonal. By convergence in $C^\infty_{\text{loc}}$, the same is true for other nearby $\text{ev}_T$, and in particular the families $(T, \sigma)$ are regular families of metrics and perturbations for $T$ sufficiently large. This is the first part of the proposition.

Now \cite[Lemma~19.3.3]{KMSW} gives the second part: the fact that this map $$(0, \infty]^{n+k_1 + \cdots + k_\ell + m} \times \mathcal N \to \mathcal R$$ is transverse to the diagonal at $\infty$ means that the projection $$\text{ev}^{-1}(\Delta) \to (0, \infty]^{n+k_1 + \cdots + k_\ell + m}$$ is a topological submersion at $\infty$, which is smooth on each stratum, and hence we may obtain a corresponding chart for a neighborhood of $\text{ev}_\infty^{-1}(\Delta) \subset \text{ev}^{-1}(\Delta)$, which is smooth on each stratum. This space $\text{ev}_\infty^{-1}(\Delta)$ contains a neighborhood $V(\sigma)$ of the stratum of $\overline{\mathcal M}_S$ described in the statement, and so we may conclude.
\end{proof}


\chapter{Orientations}\label{chap:5}
\section{Orientability and determinant line bundles}\label{sec:or-detlines}
We begin by briefly recalling some generalities on determinant lines. Let $H_1$ and $H_2$ be Hilbert spaces, and let $\mathcal F(H_1, H_2)$ be the space of Fredholm maps. We would like to claim that there is a natural line bundle over $\mathcal F(H_1, H_2)$, called the \textit{determinant line bundle}, whose fiber over $T$ is isomorphic to 
$$\Lambda^*(\text{ker}(T)) \otimes \Lambda^*(\text{coker}(T))^*);$$ 
we intuit this as being the determinant line of the virtual vector space $\text{ker}(T) - \text{coker}(T)$.

If one attempts to do this naively (with literally those fibers), the result does not have continuous transition functions due to jumps in the kernel. Instead, given any Fredholm operator $T_0$, choose a finite-dimensional subspace $J \subset H_2$ so that 
$$T_0 \oplus \text{Id}_J: H_1 \oplus J \to H_2$$ 
is surjective, and instead define the determinant line bundle to be 
$$\Lambda^* \text{ker}(T \oplus \text{Id}_J) \otimes \Lambda^*(J^*)$$ 
for $T$ near $T_0$. Because $T + \text{Id}_J$ is surjective for $T$ close enough to $T_0$, this gives a well-defined line bundle over this chart; one should argue that the line bundles are identical for different choices of $J$, and that doing this on an open cover defines a legitimate line bundle on $\mathcal F(H_1, H_2)$.

The idea of the determinant line bundle originates in \cite{quillen} in the context of Riemann surfaces; the operators $\overline \partial_J$ define a family of Fredholm operators over the moduli space of all Riemann surfaces. Our application is similar, as with most applications of this construction. Also see \cite[Section~20.2]{KMSW}, in which the authors describe the construction of the determinant line bundle, and explain how to construct a canonical associative isomorphism $$q: \det(\mathcal F(H_1, H_2)) \otimes \det(\mathcal F(K_1, K_2)) \to \det(\mathcal F(H_1 \oplus K_1, H_2 \oplus K_2)),$$ its definition involves signs that depend on the choice of $J$ in a local description as above. All in all, we have smooth real line bundles over each $\mathcal F(H_1, H_2)$ which have canonical direct-sum isomorphisms.

Now let $(W,\mathbf{E})$ be a Riemannian 4-manifold equipped with weakly admissible $SO(3)$-bundle and regular perturbation $\pi$, with two cylindrical ends: one incoming $(Y_1, E_1,\pi_1)$ and one outgoing $(Y_2, E_2,\pi_2)$. If $\alpha$ is a $\pi_1$-flat connection on $Y_1$, and similarly $\beta$ is $\pi_2$-flat, then the space of connections $\mathcal A_{\mathbf{E}, k, \delta}(\alpha, \beta)$ admits a smooth map $$\mathcal A_{\mathbf{E}, k, \delta}(\alpha, \beta) \to \mathcal F\left(\Omega^1_{k,\delta}, \Omega^{2,+}_{k-1,\delta} \oplus \Omega^0_{k-1,\delta}\right),$$ sending a connection $\mathbf{A}$ to the perturbed ASD operator $$Q_{\mathbf{A},\pi} = (d^+_{\mathbf{A}} + D_{\mathbf{A}} \widehat{\nabla}_\pi, d^*_{\mathbf{A}});$$ pulling back the determinant line bundle on the space of Fredholm operators, we find that we have a line bundle over $\mathcal A^{(4)}_{\mathbf{E}, k, \delta}$. 

Furthermore, because $Q_{\mathbf{A}}$ is invariant under the gauge group, we see that $Q$ descends to an $SO(3)$-invariant map $\widetilde{\mathcal B}^e_{\mathbf{E},k,\delta} \to \mathcal F$, and hence that there is an $SO(3)$-equivariant line bundle $\det(Q)$ over this space.

\begin{lemma}\label{global-orientability}The line bundle $\det(Q)$ over $\widetilde{\mathcal B}^e_{\mathbf{E},k,\delta}$ is trivializable.
\end{lemma}
\begin{proof}For convenience of notation we drop the super/subscripts on $\mathcal B$; everything here is modulo even gauge.

We begin with Donaldson's stabilization trick, adapted for $SO(3)$-bundles (as opposed to $SU(2)$-bundles). There is a natural map $\mathcal A_{\mathbf{E},k,\delta} \to \mathcal A_{\mathbf{E} \oplus \mathbb C, k, \delta}$ given by taking the direct sum with the trivial connection; this is equivariant for the even gauge group (where for an $SO(n)$-bundle, `even' means that the gauge transformation lifts to a section of $\text{Aut}(E) \times_{SO(n)} \text{Spin}(n)$, the action by conjugation). Therefore, this map descends to a map on the quotients $\widetilde{\mathcal B}^e_{\mathbf{E}} \to \widetilde{\mathcal B}^e_{\mathbf{E} \oplus \mathbb C}.$ Write $\det(Q')$ for the determinant line bundle for connections over $\mathbf{E} \oplus \mathbb C$; for a connection of the form $\mathbf{A} \oplus \theta$, the operator $Q'$ splits as a direct sum of the ASD deformation complex for the three bundles $\mathfrak g_E, \;\mathbb C^* \otimes E,$ and  $\mathbb R$, where on the first term we use the connection induced by $\mathbf{A}$, on the second term the connection $\theta \otimes \mathbf{A}$, and on the final term the trivial connection. Note that the second term is a \emph{complex linear} operator, and hence has canonically trivial determinant.

Thus this stabilization map pulls back $\det(Q')$ to $\det(Q) \oplus \det(Q_\theta)$, where the final term is a trivial line bundle. In particular, if $\det(Q')$ is trivializable, then $\det(Q)$ is as well. Write $\mathcal B_{(\ell)}$ for the configuration space of connections on $\mathbf{E} \oplus \mathbb C^{\ell}$, modulo even gauge, and write $\widetilde{\mathcal B}_{(\ell)}$ for the smooth $SO(2\ell + 3)$-manifold of framed configurations on $\mathbf{E} \oplus \mathbb C^\ell$, whose quotient is $\widetilde{\mathcal B}_{(\ell)}$.

Next, because the reducible subspace is of infinite codimension, transversality implies that the inclusion of the irreducible subspace $\widetilde{\mathcal B}_{(\ell)}^* \hookrightarrow \widetilde{\mathcal B}_{(\ell)}$ is a weak homotopy equivalence. In particular the restriction map $$[\widetilde{\mathcal B}_{(\ell)}, \mathbb{RP}^\infty] \to [\widetilde{\mathcal B}_{(\ell)}^*, \mathbb{RP}^\infty]$$ is a bijection, so any line bundle which is trivializable over $\widetilde{\mathcal B}_{(\ell)}^*$ is necessarily globally trivializable. Now our bundle is crucially an $SO(3)$-equivariant line bundle over $\widetilde{\mathcal B}_{(\ell)}^*$, so it is the pullback of some line bundle over $\mathcal B^{*}_{(\ell)}$, the space of irreducible connections modulo gauge (with no framing). Our goal now is to show that this line bundle is trivial, which will start from a calculation of $\pi_1 \mathcal B^*_{(\ell)}$. By definition, this space is the quotient of $\mathcal A^*$ by the even gauge group $\mathcal G^{e,h}_{\mathbf{E} \oplus \mathbb C^\ell}$; because $\mathcal A^*$ is again the complement of a union of infinite codimension submanifolds, we see that it too is contractible, and so $\pi_1 \mathcal B^*_{(\ell)} = \pi_0 \mathcal G^{e,h}_{\mathbf{E} \oplus \mathbb C^\ell}$. Next, we must calculate this group of components. This is where stabilization is first useful.

Fix a principal $G$-bundle $P$ over a CW complex $Y$. It is not difficult to show that the pointed mapping space $\text{Map}^P_*(Y, BG)$ gives a model for $B\mathcal G(P)$, where the superscript indicates we are only interested in the component which classifies the bundle $P$. Now, because $BSO$ is an $H$-space (with product structure given by direct sum of vector bundles), in fact every component of $\text{Map}_*(Y, BSO)$ is homotopy equivalent to any another. For $n$ sufficiently large, 
$$\pi_i \text{Map}^P_*(Y, BSO) \cong \text{Map}^P_*(Y, BSO(n))$$ 
(this is true for $i=1$ already as soon as $n = 7$). To compute the fundamental group of the quotient by all gauge transformations, it suffices to compute the fundamental group of the identity component of $\text{Map}_*(Y, BSO)$. As before, because we are only interested in the quotient by even gauge transformations, we actually need only compute $\pi_1 \text{Map}_*(Y, B\text{Spin})$. Because there is a fibration 
$$F \to B\text{Spin} \to K(\mathbb Z, 4)$$ 
with $F$ 6-connected, we see that 
$$\pi_1 \text{Map}_*(Y, B\text{Spin}) \cong \pi_1 \text{Map}_*(Y, K(\mathbb Z, 4)).$$ 
Now, it is a theorem of Thom \cite{ThomMappingSpaces} that we have a splitting $$\text{Map}_*(Y, K(G,n)) = \prod_{i=1}^n K(H^{n-i}(Y;G), i),$$ and in particular, $\pi_1 \text{Map}_*(Y, K(\mathbb Z, 4)) \cong H^3(Y;\mathbb Z)$.

In our case, the gauge group $\mathcal G^{e,h}$ consists of even gauge transformations which are asymptotically parallel; this fits into an exact sequence $\mathcal G^e \to \mathcal G^{e,h} \to \Gamma_\alpha \times \Gamma_\beta$, where the groups $\Gamma$ are the stabilizers of the corresponding connections and $\mathcal G^e$ is the gauge group of asymptotically trivial gauge transformations. in particular, we see that $\pi_0 \mathcal G^e \to \pi_0 \mathcal G^{e,h}$ is a surjection. The homotopy type of the group of gauge transformations which are asymptotically trivial is the same as the relative gauge group over $(W, \partial W)$; the discussion of the above paragraph goes through without difficulty to this relative setting, replacing cohomology with relative cohomology. We conclude that 
$$\pi_0 \mathcal G^e = H^3(W, \partial W; \mathbb Z) \cong H_1(W;\mathbb Z)$$ 
by Poincar\'e duality, and in particular $H_1(W;\mathbb Z) \to \pi_1 \widetilde{\mathcal B}_{\mathbf{E} \oplus \mathbb C^\ell}$ is surjective.

Now that we know what the loops in this space are, the goal is to show that $w_1(\det(Q'))$ pairs trivially against them. Now we are in standard territory, and the proof follows exactly as in \cite[Lemma~3.23]{donaldson1987orientation} (modifying the Poincar\'e duality argument to account for the boundary): if one chooses a base connection $\mathbf{A}$, a concentrated charge-one instanton $J$ on $S^4$, and a loop $[\gamma] \in H_1(W)$ in the interior of $W$, then one performs a family connected-sum construction along the loop to construct a loop of connections. One may identify both that the corresponding loop in $\pi_1 \widetilde{\mathcal B}_{\mathbf{E} \oplus \mathbb C^\ell}$ is the image of $[\gamma]$, and that the determinant line bundle over this loop of connections is trivial (if one so desires, this can be carried out using the index-gluing constructions along the boundary of a 4-manifold, which will be described later in this section).

Therefore, the pairing of $w_1(\det(Q'))$ against $H_1(W;\mathbb Z)$ is identically zero, and so $\det(Q')$ is trivializable. Being the pullback of the sum of this bundle and a trivial line bundle, the same is true of $\det(Q)$, as desired.
\end{proof}

Now choose a framed $\pi$-instanton $[\mathbf{A}, p] \in \widetilde{\mathcal M}_{\mathbf{E},z,\pi}$, where $\pi$ is a regular perturbation. The tangent space $T_{[\mathbf{A}, p]} \widetilde{\mathcal M}_{\mathbf{E},z,\pi}$ is identified with the first cohomology of the deformation complex $$\Omega^0_{k+1,\delta}(\mathfrak g_E) \xrightarrow{d_\mathbf{A}, \text{ev}_b} \Omega^1_{k,\delta}(\mathfrak g_E) \oplus \mathfrak g \xrightarrow{D_{\mathbf{A},\pi}} \Omega^{2,+}_{k-1,\delta}(\mathfrak g_E),$$ where $\mathfrak g$ means the fiber of $\mathfrak g_E$ above $b$, and $\text{ev}_b \sigma = \sigma(b)$. We denote this complex $C_{\mathbf{A}}$. Here we require $k$ to be sufficiently large that every $L^2_{k+1}$ function is continuous, so that evaluation is a continuous operation. Write $\mathfrak g_{\mathbf{A}}$ to denote the subspace of $\mathfrak g$ that extends to $\mathbf{A}$-parallel gauge transformations (so that $\mathfrak g_{\mathbf{A}}$ is the Lie algebra of $\Gamma_{\mathbf{A}}$), and write $e: \mathfrak g_{\mathbf{A}} \to \Omega^0_{k+1,\delta,\text{ext}}$ for the operator that takes a point-value and sends it to the corresponding parallel section.

Take a map $p_{\mathbf{A}}: \mathfrak g_{\mathbf{A}} \to \Omega^0(\mathfrak g_E)$ which has $\langle p_{\mathbf{A}}(h), e(h) \rangle > 0$ for all $h \neq 0 \in \mathfrak g_{\mathbf{A}}$. We write the operator $$\widetilde Q_{\mathbf{A},\pi}: \Omega^1_{k,\delta} \oplus \mathfrak g \to \Omega^{2,+}_{k-1,\delta} \oplus \Omega^0_{k-1,\delta}$$ for $\widetilde Q = Q + p_{\mathbf{A}}$ (suppressing the $p_{\mathbf{A}}$ from notation; as we shall see shortly, this is not so depraved).

\begin{lemma}The composite map $\textup{ker}(\widetilde Q) \hookrightarrow \textup{ker}(D_{\mathbf{A},\pi}) \twoheadrightarrow H^1(C_{\mathbf{A}})$ is an isomorphism if $p_{\mathbf{A}}$ is sufficiently small.
\end{lemma}
\begin{proof}By assumption that $\pi$ is a regular perturbation, the complex $C_{\mathbf{A}}$ has cohomology concentrated in degree one, and $\text{coker}(Q_{\mathbf{A},\pi}) = \mathfrak g_{\mathbf{A}}$. By assumption that $p_{\mathbf{A}}(h)$ pairs nontrivially with $e(h)$ for any $h \in \mathfrak g_{\mathbf{A}}$, we see that $\widetilde Q_{\mathbf{A},\pi}$ is surjective. Because the index of the operator $\widetilde Q$ agrees with the index of the complex $C$, we see that $\text{ker}(\widetilde Q)$ and $H^1(C_{\mathbf{A}})$ have the same dimension; so it suffices to show surjectivity.

So suppose $D_{\mathbf{A},\pi} \omega = 0$; our goal is to show there is some $\sigma \in \Omega^0_{\text{ext}}$ so that 
$$\Delta_{\mathbf A} \sigma + p_{\mathbf{A}}(\sigma(b)) = - d_{\mathbf A}^*\omega.$$ 
The map $\Delta_{\mathbf{A}}: \Omega^0_{k+1,\text{ext},\delta} \to \Omega^0_{k-1,\delta}$ is a self-adjoint operator with kernel $\mathfrak g_{\mathbf{A}}$.

In particular, because $p_{\mathbf{A}}: \text{ker}(\Delta_{\mathbf{A}}) \to \text{coker}(\Delta_{\mathbf{A}})$ is assumed to be an isomorphism, the map $\Delta_{\mathbf{A}} + tp_{\mathbf{A}}$ is surjective for sufficiently small $t$. So (replacing $p_{\mathbf{A}}$ with $tp_{\mathbf{A}}$ if necessary, $t$ small) we see that $\Delta_{\mathbf{A}} + p_{\mathbf{A}}: \Omega^0_{k+1,\delta,\text{ext}} \to \Omega^0_{k-1,\delta}$ is an isomorphism. Therefore we may indeed solve the desired equation, so that we see the map $\text{ker}(\widetilde Q) \to H^1(C_{\mathbf{A}})$ is surjective (and hence an isomorphism), as desired.
\end{proof}

In particular, we have a canonical isomorphism \emph{in a neighborhood of $[\mathbf{A},p]$} between $\det(T\widetilde{\mathcal M})$ and $\det(\widetilde Q)$. The sign of this isomorphism is independent of the auxiliary choice of $p_{\mathbf{A}}$. What we look to investigate is, then, the latter operator.

Recall now that the index of Fredholm operators is invariant under homotopy through Fredholm operators; because the map $Q_{\mathbf{A},\pi} + tp_{\mathbf{A}}$, for $t \in [0,1]$, gives a homotopy through Fredholm operators (as $\mathfrak g$ is finite-dimensional), we see that $I(\widetilde{Q}_{\mathbf{A},\pi}) \cong I(Q_{\mathbf{A},\pi}) \oplus \mathfrak g$, and in particular 
$$\det(\widetilde{Q}_{\mathbf{A},\pi}) \cong \det(Q_{\mathbf{A},\pi}) \otimes \det(\mathfrak g).$$ 
Therefore, a choice of orientation of the Lie group $SO(3)$ and an orientation of $Q$ canonically gives us an orientation of $\widetilde Q$. This is completely independent of the choice of $r$, and we may as well remove $r$ from the discussion of the index theory, and instead study $Q_{\mathbf{A},\pi}$ itself. 

In particular, a choice of orientation of the line bundle $\det(Q_{\mathbf{A},\pi})$ for $\mathbf{A}$ in some component of $\widetilde{\mathcal B}^e_{\mathbf{E}}$, and a fixed orientation of $\mathfrak so(3)$, therefore gives a global orientation of the framed moduli space $\widetilde{\mathcal M}$ sitting inside that component.

If $\mathbf{A} \in \widetilde{\mathcal B}^e_{\mathbf{E},z,k,\delta}(\alpha, \beta)$, where here $z \in \pi_0 \widetilde{\mathcal B}^e_{\mathbf{E}, k, \delta}$ labels a path-component, we write $\Lambda^W_z(\alpha, \beta)$ for the two-element set of orientations of $\det(Q_{\mathbf{A}}) \otimes \det(\mathfrak g)$. As above, a choice of element of $\Lambda^W_z(\alpha, \beta)$ induces an orientation on the moduli space $\widetilde{\mathcal M}_{\mathbf{E},z,\pi}(\alpha, \beta)$.

\section{A gluing operation for determinant lines}\label{sec:or-gluing1}
Let $(W_i, \mathbf{E}_i)$ be two cobordisms (considered as manifolds with two cylindrical ends), so that the positive end of $(W_1,\mathbf{E}_1)$ is isometric to the negative end of $(W_2,\mathbf{E}_2)$. Here we include the possibility that either or both of the $W_i$ are cylinders; write $(W_{12}^T, \mathbf{E}_{12})$ for the composite cobordism, glued by cutting off the positive end of $W_1$ (say, $[0, \infty) \times Y_2$) at $\{T\} \times Y_2$ and similarly for the negative end of $W_2$. If it is not necessary for clarity, we drop the $T$ from notation. If $z_1$ and $z_2$ are components of $\widetilde{\mathcal B}_{\mathbf{E}_1, k, \delta}(\alpha, \beta)$ and $\widetilde{\mathcal B}_{\mathbf{E}_2, k, \delta}(\beta, \gamma)$, respectively, then denote the apparent `composite component' as $z_{12}$. 

What we would like is to take an element of each of $\Lambda^{W_1}_{z_1}(\alpha, \beta)$ and $\Lambda^{W_2}_{z_2}(\beta, \gamma)$, and output an element of $\Lambda^{W_{12}}_{z_{12}}(\alpha, \gamma)$ (possibly depending on some extra input coming from $\beta$ alone).

To do this, we will use the mechanism of gluing indices along boundary components. We begin by setting up some notation.

Let $X_1$ and $X_2$ be 4-manifolds equipped with $SO(3)$-bundles $\mathbf{E}_i$, possibly with cylindrical ends, so that $\partial X_i$ is nonempty and the metric is of product type near the boundary. We decompose each $\partial X_i$ into a union $\partial_- X_i \cup \partial_+ X_i$ of connected components (either term possibly empty) which we call `positive' or `negative', and orient them so that a neighborhood of $\partial_+ X_i$ is isometric to $(-t, 0] \times \partial_+ X_i$ and a neigbhorhood of $\partial_- X_i$ is isometric to $[0,t) \times \partial_- X_i$. Suppose furthermore that $(\partial_+ X_1, \mathbf{E}_1)$ is oriented isometric to $(\partial_- X_2, \mathbf{E}_2)$. Write $$Y_1 = \partial_- X_1, \;\;\; Y_2 = \partial_+ X_1 = \partial_- X_2 \;\;\; Y_3 = \partial_+ X_2,$$ and write $Y = Y_1 \cup Y_2 \cup Y_3$.

Then we may form the composite $X_{12} = X_1 \cup_{Y_2} X_2$, equipped with $SO(3)$-bundle $\mathbf{E}_{12}$.

Suppose each $(X_i, \mathbf{E}_i)$ is equipped with a choice of 4-manifold perturbation, regular on the ends, which is also constant in time near the boundary. Suppose each is equipped with a connection $\mathbf{A}_i$ of regularity $L^2_{k,\delta}$, decaying towards $\pi$-flat connections on the ends, and suppose that $\mathbf{A}_{12}$ is a connection on $X_{12}$ of the same regularity which agrees with each $\mathbf{A}_i$ when restricted to each $X_i$. Write $A$ for the restriction of $\mathbf{A}_{12}$ to $Y$; it is of regularity $L^2_{k-1/2}$.

We are interested in studying the perturbed ASD operator on these compact manifolds with boundary. We should introduce a restriction on the boundary values to make this a Fredholm operator, as in Chapter \ref{sec:4d-gluing}. However, the choice there (the Coulomb-Neumann gauge) is not appropriate for us here; that boundary value problem is not continuously varying as we change the reducibility type of $A$. Instead, we use the standard Atiyah-Patodi-Singer spectral boundary value problem. On the composite $X_{12}$, the ASD operator takes the form $d/dt + L_{A,\pi}$ on $(-t, t) \times Y$ for the signature operator $$L_{A,\pi}: \Omega^0_{k-1/2}(Y) \oplus \Omega^1_{k-1/2}(Y) \to \Omega^0_{k-3/2}(Y) \oplus \Omega^1_{k-3/2}(Y),$$ written in matrix form as $$L_{A,\pi} = \begin{pmatrix}0 & -d_A^* \\ -d_A & D_{A,\pi}\end{pmatrix},$$ where $D_{A,\pi} = *d_A + D_A \nabla_\pi$. 

Henceforth we write $\Omega^0(Y) \oplus \Omega^1(Y)$ as $H$, only writing the Sobolev indices if necessary. When there are multiple $Y_i$, we write $$H_i := \Omega^0(Y_i) \oplus \Omega^1(Y_i).$$

Suppose $\lambda \in \mathbb R$ is not an eigenvalue of $L_{A,\pi}$ for any $A$ arising as the restriction of $\mathbf{A}$ to a boundary component. Then we may split $H$ as $H^{>\lambda} \oplus H^{< \lambda}$ as the closure of the linear span of eigenvectors with eigenvalue greater than, or less than, $\lambda$. Then we define the $\lambda$-weighted ASD operator on $X_1$ to be 
$$Q^\lambda_{\mathbf{A},\pi} = (D_{\mathbf{A},\pi}, d^*_{\mathbf{A}}, \Pi^{< \lambda}_1r, \Pi^{> \lambda}_2r): \Omega^1_{k,\delta}(X_1) \to \Omega^{2,+}(X_1) \oplus \Omega^0(X_1) \oplus H_1^{< \lambda} \oplus H_2^{> \lambda},$$ 
where $r$ is restriction to the appropriate boundary components and the operators $\Pi$ are the $L^2_{k-1/2}$ orthogonal projections onto the corresponding eigenspaces. The operator written on $X_2$ is the same, with the boundary value projections changed to $\Pi^{<\lambda}_2$ and $\Pi^{> \lambda}_3$, and similarly with $X_{12}$.

With this notation in hand, the following is essentially the content of the beginning of \cite[Section~20.3]{KMSW}, especially pages 383-384.

\begin{lemma}\label{index-gluing}Let $K$ be any compact family of connections $\mathbf{A}_{12}$ and perturbations as above, where $\lambda$ is never an eigenvalue for any $L_{A,\pi}$ for any $A$ given by restriction of $\mathbf{A} \in K$ to the boundary. Then the operators $Q^\lambda_{\mathbf{A}_i,\pi}$ form a family of Fredholm operators, and there is a canonical isomorphism of line bundles over $K$
$$\det(Q^\lambda_{\mathbf{A}_1, \pi}) \otimes \det(Q^\lambda_{\mathbf{A}_2, \pi}) \to \det(Q^\lambda_{\mathbf{A}_{12},\pi}),$$ given by a homotopy between the operators $Q^\lambda_{\mathbf{A}_1, \pi} \oplus Q^\lambda_{\mathbf{A}_2, \pi}$ and an operator whose kernel and cokernel are canonically identified with that of $Q^\lambda_{\mathbf{A}_{12},\pi}$.
\end{lemma}

What we really care about are connections $\mathbf{A} \in \widetilde{\mathcal B}_{\mathbf{E},z,k,\delta}(\alpha, \beta)$ on a manifold $W$ with two cylindrical ends and no boundary components. To make use of the above lemma, we should break these into pieces. For $T$ large, write $\mathbf{N}$ for the restriction of $\mathbf{A}$ to $(-\infty, -T] \times Y_1$, $\mathbf{P}$ for the restriction to $[T, \infty) \times Y_2$, and $\mathbf{C}$ for the restriction to the compact piece lying between these.

We fix a small constant $\epsilon > 0$. What we demand of our choice of $T$, writing $\mathbf{N} = p^*\alpha + a$ for $p^*\alpha$ the pullback connection $d/dt + d_\alpha$ and $a \in \Omega^1_{k,\delta}\left( (-\infty, -T] \times Y_1\right)$, is the following. First, that the signature operator $L_{\alpha + t a(-T)}$ never has an eigenvalue of absolute value $\epsilon$ for any $t \in [0,1]$. Second, that the only eigenvalue of $L_{\alpha,\pi}$ in $[-\epsilon, \epsilon]$ is $0$, which corresponds to $\text{ker}(L_{\alpha, \pi}) = \mathfrak g_\alpha$, the space of $\alpha$-parallel sections of $\mathfrak g_E$; this is the entire kernel by the assumption that $\alpha$ is a nondegenerate critical orbit. We demand similarly for $\mathbf{P}$. This is a smallness condition on $\mathbf{N}$ and $\mathbf{P}$, or rather the size of the differences between those and constant connections. In particular, the operators $Q^{-\epsilon}_{\mathbf{N},\pi}$ and $Q^{-\epsilon}_{\mathbf{P},\pi}$ are both Fredholm operators, and there are canonical isomorphisms $$\det(Q_{\mathbf{N},\pi}^{-\epsilon}) \cong \det(Q_{p^*\alpha,\pi}^{-\epsilon}), \;\;\;\; \det(Q_{\mathbf{P},\pi}^{-\epsilon}) \cong \det(Q_{\mathbf{\beta},\pi}^{-\epsilon}).$$ This remains true for families of such connections, so long as $T$ is chosen to satisfy the given properties for all connections in the family; this is always possible if the family is compact.

Then by Lemma \ref{index-gluing}, we have an isomorphism $$\det(Q_{p^*\alpha,\pi}^{-\epsilon}) \otimes \det(Q_{\mathbf{C},\pi}^{-\epsilon}) \otimes \det(Q_{p^*\beta,\pi}^{-\epsilon}) \cong \det(Q_{\mathbf{A},\pi}).$$ To reduce this to studying $\mathbf{C}$ alone, we should understand explicitly the operators corresponding to each end.

\begin{lemma}Write $Z^- = (-\infty, 0] \times Y_1$ and $Z^+ = [0, \infty) \times Y_2$.

The operator $$Q_{p^*\alpha, \pi}^{-\epsilon}: \Omega^1_{k,\delta}(Z^-) \to \Omega^{2,+}_{k-1,\delta}(Z^-) \oplus \Omega^0_{k-1,\delta}(Z^-) \oplus H_1^{< -\epsilon}$$ is an isomorphism, and so there is a canonical trivialization $\det(Q_{p^*\alpha, \pi}^{-\epsilon}) \cong \mathbb R$. The operator $$Q_{p^*\beta, \pi}^{-\epsilon}: \Omega^1_{k,\delta}(Z^+) \to \Omega^{2,+}_{k-1,\delta}(Z^+) \oplus \Omega^0_{k-1,\delta}(Z^+) \oplus H_1^{> -\epsilon}$$ is injective with cokernel canonically identified with $\mathfrak g_\beta$, the space of $\beta$-parallel sections of $\mathfrak g_E$, and so we have a canonical isomorphism $\det(Q_{p^*\beta, \pi}^{-\epsilon}) \cong \det(\mathfrak g_\beta)^*$, giving a trivialization as soon as we orient the vector space $\mathfrak g_\beta$.
\end{lemma}
\begin{proof}Because the operator $$(D_{p^*\alpha, \pi}, d^*_{p^*\alpha}) = (d^+_{p^*\alpha} + D_\alpha \widehat{\nabla}_\pi, d^*_{p^*\alpha})$$ takes the form $d/dt + L_{\alpha, \pi}$, and identically with $\beta$, we may solve the equation $(D_{p^*\alpha, \pi}, d^*_{p^*\alpha})\omega = 0$ by separation of variables, decomposing $\Omega^0(Y_1) \oplus \Omega^0(Y_2)$ into the eigenspaces of $L_{\alpha,\pi}$. A solution with boundary value $\phi \in L^2_{k-1/2}$ exists if and only if $\phi$ is in the closure of the span of eigenvectors with negative eigenvalues. Because $H_1^{< -\epsilon} = H_1^{< 0}$ by assumption on $\epsilon$, we see that the operator $Q_{p^*\alpha, \pi}^{-\epsilon}$ is invertible, as desired. However, for the other cylindrical end, we have $H_2^{> -\epsilon} = H_2^{\geq 0}$, we see that we have cokernel equal to $H_2^0 = \text{ker}(L_{\alpha,\pi}) = \mathfrak g_\alpha$, by the assumption that $\alpha$ is nondegenerate.
\end{proof}

We use these facts about the index to construct a comparison map between certain determinant line bundles.

\begin{definition}\label{analytic-gluing}Let $\mathbf{A}_1$ and $\mathbf{A}_2$ be connection on cobordisms $W_1$ and $W_2$ as above, going between $\alpha$ and $\beta$ or $\beta$ and $\gamma$, respectively. If $\mathbf{A}_{12}$ is a connection on the composite $W_{12}^T$ which is uniformly close to $\mathbf{A}_1$ on $W_1^{\leq T}$ and $\mathbf{A}_2$ on $W_2^{\geq -T}$, the there is a canonical isomorphism $$\rho_{\textup{an}}: \det\left({Q_{\mathbf{A}_1, \pi}}\right) \otimes \det(\mathfrak g_\beta) \otimes \det\left(Q_{\mathbf{A}_2, \pi}\right) \cong \det\left(\widetilde{Q}_{\mathbf{A}_{12}, \pi}\right),$$ given by decomposing $$\det\left({Q_{\mathbf{A}_1, \pi}}\right) \cong \det\left({Q^{-\epsilon}_{\mathbf{C}_1, \pi}}\right) \otimes \det(\mathfrak g_\beta)^*$$ and $$\det\left({Q^{-\epsilon}_{\mathbf{A}_2, \pi}}\right) \cong \det\left({Q_{\mathbf{C}_2, \pi}}\right) \otimes \det(\mathfrak g_\gamma)^*,$$ as above, and applying the canonical isomorphism $\det(\mathfrak g_\beta)^* \det(\mathfrak g_\beta) \cong \mathbb R$, as well as the isomorphism $$\det\left({Q^{-\epsilon}_{\mathbf{C}_1, \pi}}\right) \otimes \det\left({Q^{-\epsilon}_{\mathbf{C}_2, \pi}}\right) \cong \det\left({Q^{-\epsilon}_{\mathbf{C}_{12}, \pi}}\right)$$ given by Lemma \ref{index-gluing}, by first applying a small homotopy taking $\mathbf{C}_i$ to the restriction of $\mathbf{C}_2$.

We call $\rho_{\textup{an}}$ the \emph{analytic gluing map}. There is a corresponding version of this map for the operators $\widetilde Q$ which have a factor of $\mathfrak g$ to account for the $SO(3)$-action; this is written
$$\widetilde \rho_{\textup{an}}: \det\left(\widetilde{Q}_{\mathbf{A}_1, \pi}\right) \otimes \det(\mathfrak g_\beta^\perp)^* \otimes \det\left(\widetilde{Q}_{\mathbf{A}_2, \pi}\right) \cong \det\left(\widetilde{Q_{\mathbf{A}_{12}, \pi}}\right),$$ given by decomposing $$\det\left(\widetilde{Q}_{\mathbf{A}_i, \pi}\right) \cong \det\left(Q_{\mathbf{A}_i, \pi}\right) \otimes \det(\mathfrak g)$$ and using the canonical isomorphism $\det(\mathfrak g) \cong \det(\mathfrak g_\beta) \otimes \det(\mathfrak g_\beta^\perp)$, where here we implicilty demand that $\mathfrak g = \mathfrak g_\beta \oplus \mathfrak g_\beta^\perp$ is an oriented splitting if we choose an orientation on any two of these.

In particular, this gives an isomorphism $$\Lambda^{W_1}_{z_1}(\alpha, \beta) \times_{\mathbb Z/2} \Lambda(\mathfrak g_\beta) \times_{\mathbb Z/2} \Lambda^{W_2}_{z_2}(\beta, \gamma) \cong \Lambda^{W_{12}}_{z_{12}}(\alpha, \gamma),$$ where $\Lambda(\mathfrak g_\beta)$ is the set of orientations of $\mathfrak g_\beta$.

Given any compact family of $\mathbf{A}_i$ as above, the analytic gluing map $\rho_{\text{an}}$ gives an isomorphism of line bundles over the parameter space, and identically with $\widetilde{\rho}_{\text{an}}$.
\end{definition}

Note that this is usually relevant to us in the case where one of $W_i$ is a cylinder. If neither are cylinders, the orientations that actually appear in practice are of the parameterized ASD operator, with parameter given by some family of broken metrics. 

\section{The gluing diffeomorphism is orientation-preserving}\label{sec:or-gluing2}
In the previous section, we defined a gluing isomorphism $\rho_{\text{an}}$ between certain determinant lines. We also have a similar isomorphism coming from a different source: the inverse function theorem. Recall from Chapter \ref{sec:4d-gluing} that, for certain connected $SO(3)$-invariant open subsets $U_{\alpha \beta} \subset \widetilde{\mathcal M}_{\mathbf{E}_1,z_1}(\alpha, \beta)$ and $U_{\beta \gamma} \subset \widetilde{\mathcal M}_{\mathbf{E}_2, z_2}(\beta, \gamma)$, Proposition \ref{gluing} gives a diffeomorphism 
$$\text{gl}: U_{\alpha \beta} \times_\beta U_{\beta \gamma} \to U_{\alpha \gamma}.$$ 
In the compactification by broken trajectories, the closure of $\text{gl}(U_{\alpha \beta} \times_\beta U_{\beta \gamma})$ includes the open face 
$$\widetilde{\mathcal M}_{\alpha \beta} \times_\beta \widetilde{\mathcal M}_{\beta \gamma} \subset \partial \overline{\mathcal M}_{\alpha \gamma}.$$ 
We remind the reader here that when considering $\widetilde{\mathcal M}_{\mathbf{E}, k, \delta}$ in the case where $W$ is a cylinder, this is the moduli space of \emph{parameterized} trajectories: we do not reduce by the $\mathbb R$-action.

In particular, suppose we are given a pair of framed instantons $([\mathbf{A}_1,p_1], [\mathbf{A}_2,p_2])$ which project to the same framing at $\beta$. Taking the derivative of $\text{gl}$, we obtain an isomorphism $$\widetilde{\rho}_{\text{gm}}: \det\left(\widetilde{Q}_{\mathbf{A}_1, \pi}\right) \otimes \det(\mathfrak g_\beta)^* \otimes \det\left(\widetilde{Q}_{\mathbf{A}_2, \pi}\right) \cong \det\left(\widetilde{Q}_{\mathbf{A}_{12}, \pi}\right);$$ note that the $SO(3)$-action on a point gives an isomorphism $T_p \beta \cong \mathfrak g_\beta^\perp$, which explains this factor. Because the gluing map is equivariant, this descends to an isomorphism $$\rho_{\text{gm}}: \det\left({Q_{\mathbf{A}_1, \pi}}\right) \otimes \det(\mathfrak g^\perp_\beta) \otimes \det\left({Q_{\mathbf{A}_2, \pi}}\right) \cong \det\left({Q_{\mathbf{A}_{12}, \pi}}\right).$$
We call $\rho_{\text{gm}}$ the \emph{geometric gluing map}. If the open sets $U_{\alpha \beta}$ and $U_{\beta \gamma}$ are chosen so that the restriction of $\mathbf{A} \in U_{\alpha \beta}$ to each end is close enough to the constant trajectory, and similarly for the restriction of $\mathbf{A} \in U_{\beta \gamma}$ to each end, that we may apply Lemma \ref{index-gluing}. Then we obtain \emph{two} isomorphisms
$$\rho: \det\left({Q_{\mathbf{A}_1, \pi}}\right) \otimes \det(\mathfrak g_\beta) \otimes \det\left({Q_{\mathbf{A}_2, \pi}}\right) \cong \det\left({Q_{\mathbf{A}_{12}, \pi}}\right).$$ 
The next part of this section is dedicated to showing that these two isomorphisms are positive scalar multiples of one another, and hence both induce the same orientation on $\det\left({Q_{\mathbf{A}_{12}, \pi}}\right)$.

The map $\rho_{\text{gm}}$ (or rather, $\widetilde{\rho}_{\text{gm}}$) is not difficult to understand at the level of kernels and cokernels; it is less clear how to understand $\widetilde{\rho}_{\text{an}}$ at this level, so this is our next immediate goal. 

\begin{lemma}\label{ker-coker}Let $V, W, $ and $W_b$ be Hilbert spaces; suppose we have a continuous family of Fredholm maps $$A_t = (A, r_t): V \to W \oplus W_b,$$ with $A: V \to W$ surjective. Suppose that $A_0$ is surjective, and that the map $$r_t - r_0: \textup{ker}(A) \to W_b$$ has sufficiently small norm, uniformly in $t$. Let $J \subset W_b$ be a finite-dimensional subspace so that $\widetilde A_t: V \oplus J \to W \oplus W_b,$ given by $\begin{pmatrix}A & 0\\ r_t & 1\end{pmatrix}$, is surjective for all $t$. Then the composite $$\textup{ker}(\widetilde{A}_0) \hookrightarrow V \oplus J \twoheadrightarrow \textup{ker}(\widetilde A_t)$$ is an isomorphism, where the first map is inclusion and the second map is orthogonal projection. Similarly, the map $\textup{ker}(A_0) \to \textup{ker}(A_1)$ is injective, and the above gives us an isomorphism $\textup{ker}(A_0) \oplus \textup{coker}(A_1) \cong \textup{ker}(A_1)$. This gives an isomorphism between orientation sets $\Lambda(A_0) \cong \Lambda(A_1)$, the same as that induced by the homotopy through Fredholm operators $A_t$.
\end{lemma}

\begin{proof}First, we argue in slightly more generality. Suppose we are given two Fredholm maps $T_i: X \to Y$ of Hilbert spaces, and that $\|T_2\big|_{\text{ker}(T_1)}\| < \epsilon$. Write $\pi_2: X \to \text{ker}(T_2)$ for the orthogonal projection, and $C_2 = \text{ker}(\pi_2)$, the orthogonal complement to $\text{ker}(T_2)$. Then the map $T_2: C_2 \to Y$ is an isomorphism onto its (closed) image, and in particular enjoys an inequality of the form $\|T_2 v\|_Y \geq c\|v\|_{C_2}$ for some $c > 0$. For any $v \in \text{ker}(T_1)$, we have $$\|(1-\pi_2)v\| \leq \frac 1c \|T_2 (1-\pi_2)v\| = \frac 1c \|T_2 v\|\leq \frac{\epsilon}{c} \|v\|;$$ the equality follows because $T_2 \pi_2 v = 0$, by definition of $\pi_2$. As soon as $\epsilon < c$, we see that $\pi_2 v$ must be nonzero. In particular, so long as $\|T_2\big|_{\text{ker}(T_1)}\|$ is sufficiently small, the composite $\text{ker}(T_1) \hookrightarrow X \twoheadrightarrow \text{ker}(T_2)$ is injective.

Given an element $(v, j)$ of $\text{ker}(\widetilde A_0)$, we have $$\widetilde A_t(v,j) = (Av, r_t v + j) = (0, (r_t - r_0) v);$$ because $v \in \text{ker}(A)$, we see that $$\|\widetilde A_t\big|_{\text{ker}(\widetilde A_0)} \| \leq \|(r_t - r_0)\big|_{\text{ker}(A)}\|.$$ So if $\|(r_t - r_0)\big|_{\text{ker}(A)}\|$ is sufficiently small, we see that $\text{ker}(\widetilde A_0) \to \text{ker}(\widetilde A_t)$ as defined above is injective, for all $t$; because the operators $\widetilde A_t$ are all homotopic through Fredholm operators and assumed to be surjective, these kernels all have the same dimension, and hence this injection is an isomorphism.

By the assumption that $A_0$ is surjective, we may identify $\text{ker}(\widetilde A_0)$ as $\text{ker}(A_0) \oplus J$. If we write $\text{coker}(\widetilde A_t) \oplus J'_t \cong J$, then choosing a section of the projection $\text{ker}(\widetilde A_t) \to J'_t$, we may similarly identify $\text{ker}(\widetilde A_t) \cong \text{ker}(A_t) \oplus J'_t$. The isomorphism $$\text{ker}(A_0) \oplus J \cong \text{ker}(A_t) \oplus J'_t$$ is isotopic to the direct sum of an isomorphism $\text{ker}(A_0) \oplus \text{coker}(A_t) \to \text{ker}(A_t)$ with the identity map on $J'_t$. Here, the map $\text{coker}(A_t) \to \text{ker}(A_t)$ is given (up to a small homotopy) by first choosing a section $s$ of the surjective map 
$$\text{ker}(A) \xrightarrow{-r_0} W_b \to \text{coker}(A_t),$$ 
which supplies us with a map 
\begin{align*}s': \text{coker}(A_t) &\to \text{ker}(A_0) \oplus J\\ 
s'(j) &= (s(j), r_0 s(j));\end{align*}
because the space of sections is contractible this map only depends on the choice of section up to an isotopy. Then as usual we apply orthogonal projection to $s(j)$ to obtain a map $\text{ker}(A_0) \oplus J \to \text{ker}(A_t)$.

This, then, is identified as an isomorphism of virtual vector spaces 
$$(\text{ker}(A_0), 0) \to (\text{ker}(A_t), \text{coker}(A_t)),$$
the map on kernels just being the composite $\text{ker}(A_0) \to V \to \text{ker}(A_t)$ as above.

Given that we have a family of such isomorphisms, starting with the identity, identifies this map with the isomorphism $\Lambda(A_0) \to \Lambda(A_1)$ induced by the homotopy through Fredholm operators.
\end{proof}

Let $W_1$ and $W_2$ be cobordisms, with the positive end of $W_1$ and the negative end of $W_2$ modeled on $(Y, E)$. Suppose we have regular framed instantons $(\mathbf{A}_1,p_1), (\mathbf{A}_2,p_2)$ on the respective cobordisms (possibly cylinders) which limit to the critical orbit $\beta$; the gluing map provides a framed instanton $(\mathbf{A}_{12}^T,p)$ on the glued-up cobordism 
$$W_{12}^T = W_1^{\leq T} \cup_Y W_2^{\geq -T},$$ 
obtained by truncating the two cobordisms at some large $\{\pm T\} \times Y$ and pasting them together. The defining property of $(\mathbf{A}_{12}^T,p)$, arising from the inverse function theorem, is as follows. First, for large enough $T$ the restriction to $W_1^{\leq T/2}$ and $W_2^{\geq -T/2}$ is uniformly close to $(\mathbf{A}_1, p_1)$ and $(\mathbf{A}_2, p_2)$, respectively (here we use parallel transport to move the framing between various basepoints). Second, the restriction to the region in between these two pieces, isometric to $[-T/2, T/2] \times Y$, is uniformly small. We write $[-T/2, T/2] \times Y =: Z^{T/2}$. If $\pi: E_b \to E_b/\Gamma_{\beta}$ is the projection and $q$ the framing at some basepoint $b$ on $Z^{T/2}$, then $\pi q$ is uniformly close to the common limiting value $q$ of $p_1$ and $p_2$.

In this context, we consider 
$$V = \Omega^1_{k,\delta}(W_1^{\leq T/2}) \oplus \Omega^1_{k,\delta}(Z^{T/2}) \oplus \Omega^1_{k,\delta}(W_2^{\geq -T/2}),$$ 
while 
$$W = (\Omega^{2,+} \oplus \Omega^0)_{k-1,\delta}(W_1^{\leq T/2}) \oplus (\Omega^{2,+} \oplus \Omega^0)_{k-1,\delta}(Z^{T/2}) \oplus (\Omega^{2,+} \oplus \Omega^0)_{k-1,\delta}(W_2^{\geq -T/2}),$$ 
and 
$$W_b = (\Omega^0 \oplus \Omega^1)_{k-1/2}(\{-T/2, T/2\} \times Y).$$ 

The operator $A: V \to W$ is given by $(d^*, d^+)$, while there are two restriction operators $V \to W_b$ of interest to us: one is $r_0 = (\text{res}, -\text{res})$, taking the restriction on each positive boundary component and negative the restriction on each negative boundary component. The other is $r_1 = (\Pi^{>-\epsilon}, \Pi^{<-\epsilon}),$ the spectral projection onto the $>-\epsilon$ eigenspaces on positive ends, and the $<\epsilon$ eigenspaces on negative ends. The operator $A_1$ has the property that it splits as a direct sum of operators corresponding to the pieces $W_1^{\leq T/2}, Z^{T/2}, W_2^{\geq -T/2}$, while the operator $A_0$ has kernel and cokernel naturally isomorphic to $Q_{\mathbf{A}_{12},\pi}$ itself. Here the appropriate choice of $L^2_{k,\delta}$ norm on the direct sum in fact has a weight on the $Z^{T/2}$ component, as in Chapter \ref{sec:4d-gluing}, essentially weighted by a symmetric function equal to $e^{\delta(t+T/2)}$ on $[-T/2, -t)$ for some small $t >0$.

We may interpolate through these projections via 
$$r_t = (\Pi^{>-\epsilon} + t\Pi^{<-\epsilon}, -t\Pi^{>-\epsilon} - \Pi^{< -\epsilon});$$ 
as discussed in \cite[Section~20.3]{KMSW}, this gives a homotopy $A_t = (A, r_t)$ through Fredholm operators. Furthermore, if $T$ is chosen large enough, then $\|(r_t - r_0)\big|_{\text{ker}(A)}\|$ is arbitrarily small, uniformly in $t$. Observe here that 
$$r_t - r_0 = t(\Pi^{<-\epsilon},-\Pi^{>-\epsilon}),$$ 
so we may focus attention on the case $t = 1$. For convenience, we only pay attention to the factor $\Pi^{<-\epsilon}$ as we approach $\{T/2\} \times Y$ from the left.

Write $\omega_{T/4} = \omega\big|_{\{T/4\} \times Y}$ and $\omega_{T/2} = \omega\big|_{\{T/2\} \times Y}$. For $T$ very large, we have that $\mathbf{A}_{12} \big|_{W_1^{[T/4, T/2]}}$ is uniformly close to the constant trajectory at the $\pi$-flat connection $\beta$, as follows because $\mathbf{A}_1$ decays exponentially. For some constant $C$, independent of $T$, we have $\|\omega_{T/4}\|_{k-1/2} \leq C\|\omega\|_{k,\delta}$ for all $\omega \in L^2_{k,\delta}$; this is just the claim that the trace along a hypersurface is a continuous map (and the fact that a neighborhood of $\{T/4\} \times Y$ is isometric to $(-t, t) \times Y$, independent of $T$). Writing the ASD operator for the constant trajectory at $\beta$ as an ODE in the eigenvalues of 
$$L_{\beta, \pi}: (\Omega^0\oplus \Omega^1)_{k-1/2}(Y) \to (\Omega^0 \oplus \Omega^1)_{k-3/2}(Y),$$ 
and writing $\omega_{T/2} = \omega^+_{T/2} \oplus \omega^-_{T/2}$ for the spectral decomposition, we see that $$\|\omega_{T/2}^-\| \leq \frac{e^{-\delta T/4}}{2} \|\omega_{T/4}^-\| \leq \frac{e^{-\delta T/4}}{C} \|\omega\|$$ for any $\omega \in \text{ker}(A)$, where $\delta$ is the absolute value of the least eigenvalue of $L_{\beta, \pi}$, and the factor of $2$ is simply a fudge factor to account for the fact that $\mathbf{A}_1$ is not \emph{literally} the constant trajectory at $\beta$ in the relevant portion of the cobordism. Here we use that $-\epsilon < 0$ to conclude that this component of $r_t - r_0$ is uniformly small.

For the negative boundary components, a similar argument applies, but now one must exploit the exponential weights in the definition of the $L^2_{k,\delta}$ norm on the compact cylinder and the fact that $\epsilon$ is chosen less than $\delta$ to get the desired bound (now instead there is a factor of $e^{(\epsilon - \delta)T/4}$).

This is almost sufficient to apply the above lemma, except for the condition that $A_0$ is surjective; rather, we have a canonical identification $\text{coker}(A_0) \cong \mathfrak g_{\mathbf{A}_{12}}$, the tangent space to the stabilizer $\Gamma_{\mathbf{A}_{12}}$. This is easy enough to dispatch; recall the definitions of the extended operators $\widetilde{Q}_{\mathbf{A},\pi}$ from earlier in this section, depending on a certain choice of map $p: \mathfrak g \to \Omega^0(W)$. These were chosen precisely so that $\widetilde Q_{\mathbf{A},\pi}$ is surjective for a regular instanton $\mathbf{A}$. In this context, take $\widetilde V = V \oplus \mathfrak g \oplus \mathfrak g$ (thinking of each $\mathfrak g$ as being a choice of framing on one of the components $W_1^{\leq T/2}$ or $W_2^{\geq -T/2}$, based at say $\gamma(-T)$ and $\gamma(T)$ for the base curve $\gamma$), while $\widetilde W = W \oplus \mathfrak g$ and $\widetilde W_b = W_b$. The map $\mathfrak g \oplus \mathfrak g \to \widetilde W$ is the expected maps $p_i$ on the $\Omega^0(W_i)$ factors, while the map $\mathfrak g \oplus \mathfrak g \to \mathfrak g$ is the identity in the first factor, and $-\text{Hol}^{\gamma(T) \to \gamma(-T)}_{\mathbf{A}_{12}}$ in the second factor; if the points $\gamma(\pm T)$ are chosen to lie on the boundary of the two pieces, then the composite of this map with the projection $\mathfrak g \to \mathfrak g_\beta^\perp$ is very close to the identity, because $\mathbf{A}_{12}$ is sufficiently close to the constant trajectory at $\beta$.

Then the map $\widetilde A: \widetilde V \to \widetilde W$ is surjective, as is $\widetilde A_0: \widetilde V \to \widetilde W \oplus \widetilde W_b$, and the maps $p_i$ contribute minimally to the boundary-evaluation, so these satisfy the assumptions of the lemma.

\begin{corollary}If $T$ is sufficiently large, the index-theoretic isomorphism $\widetilde{\Lambda}(\mathbf{A}_{12}) \cong \widetilde{\Lambda}(\mathbf{A}_1^{\leq T/2}) \Lambda(\mathfrak g)^* \widetilde{\Lambda}(\mathbf{A}_2^{\geq -T/2})$ is given, at the level of kernels and cokernels, by an injection $$\textup{ker}(\widetilde{Q}_{\mathbf{A}_{12}}) \to \textup{ker}(\widetilde{Q}^{\leq T/2}_{\mathbf{A}_1}) \oplus_{\mathfrak g} \textup{ker}(\widetilde Q^{\geq -T/2}_{\mathbf{A}_2}),$$ obtained first by restriction to $\Omega^1(W_1^{\leq T/2} \sqcup W_2^{\geq -T/2})$ and second by projection to the kernel. (Here $\oplus_{\mathfrak g}$ means we take the kernel of the natural projection of the direct sum to $\mathfrak g$.)

Similarly, for sufficiently large cutoffs $T$, the index-thereotic isomorphism $$\widetilde{\Lambda}(\mathbf{A}_1) \Lambda(\mathfrak g_\beta^\perp)^* \widetilde{\Lambda}(\mathbf{A}_2) \to \widetilde{\Lambda}(\mathbf{A}_1^{\leq T/2}) \Lambda(\mathfrak g)^*  \widetilde{\Lambda}(\mathbf{A}_2^{\geq -T/2})$$ may be described at the level of kernels and cokernels as a map $$\textup{ker}(\widetilde Q_{\mathbf{A}_1}) \oplus_{\mathfrak g_\beta^\perp} \textup{ker}(\widetilde Q_{\mathbf{A}_2}) \to \textup{ker}(\widetilde Q^{\leq T/2}_{\mathbf{A}_1}) \oplus \textup{ker}(\widetilde Q^{\geq -T/2}_{\mathbf{A}_2}),$$ which may be described explicitly again as restriction to the corresponding manifolds with boundary, and then projection to the latter kernel.
\end{corollary}
\begin{proof}The only point not outlined above is a discussion of the middle piece, isometric to $[-T/2, T/2] \times Y$. Given that $\mathbf{A}_{12}$ is uniformly close to the constant trajectory at $\beta$ for large enough $T$, and that for large enough $T$ the ASD map on the cylinder (with boundary conditions) is an isomorphism for the constant trajectory at $\beta$, we see that its determinant line is canonically trivial, and does not enter into the discussion.
\end{proof}

Now consider the diagram
\[
\begin{tikzcd}\widetilde{\Lambda}(\mathbf{A}_1) \Lambda(\mathfrak g_\beta^\perp)^* \widetilde{\Lambda}(\mathbf{A}_2) \arrow{r} \arrow{d} & \widetilde{\Lambda}(\mathbf{A}_1^{\leq T/2}) \Lambda(\mathfrak g)^* \widetilde{\Lambda}(\mathbf{A}_2^{\geq -T/2}) \\
\widetilde{\Lambda}(\mathbf{A}_{12}) \arrow{ur}
\end{tikzcd}
\]

The vertical map may be taken to be either $\Lambda(\widetilde{\rho}_{\text{gm}})$ or $\Lambda(\widetilde{\rho}_{\text{an}})$. The horizontal map and upper-right map are as discussed above: they maybe understood either as the inverses of index-theoretic gluing isomorphisms, or via projection maps between various kernels.

If the vertical map is taken to be $\widetilde{\rho}_{\text{an}}$, then this diagram commutes (in fact, $\rho_{\text{an}}$ is \textit{defined} so that this diagram commutes). If we take the vertical map to be $\widetilde{\rho}_{\text{gm}}$, this diagram still commutes: we are claiming that the composite map 
$$(\mathbf{A}_1, p_1, T, \mathbf{A}_2, p_2) \to (\mathbf{A}_{12}^T, p) \to (\mathbf{A}_{12}^{\leq T/2} p'_1, \mathbf{A}_{12}^{\geq -T/2}, p'_2)$$ 
is uniformly close to restriction, and hence this map commutes (up to a small homotopy) at the level of kernels. Because the diagram commutes with either choice, $\Lambda(\widetilde \rho_{\text{gm}}) = \Lambda(\widetilde \rho_{\text{an}})$: the two gluing maps are the same.

If we have fixed an orientation of $\mathfrak g_\beta$, then via the recipe given by the analytic gluing map, an orientation of $\Lambda_z^{W_1}(\alpha, \beta)$ and an orientation of $\Lambda_w^{W_2}(\beta, \gamma)$ induces an orientation of $\Lambda_{zw}^{W_{12}}(\alpha, \gamma)$; it also naturally an induces an orientation of the fiber product 
$$\widetilde{\mathcal M}_{\mathbf{E}_1, z, k, \delta}(\alpha, \beta) \times_\beta \widetilde{\mathcal M}_{\mathbf{E}_2, w, k, \delta}(\beta, \gamma),$$ 
and what we learned above is that the geometric gluing map, which gives rise to a diffeomorphism between an open subset of this fiber product and $\widetilde{\mathcal M}_{\mathbf{E}_{12}, zw, k, \delta}(\alpha, \gamma)$ is \emph{orientation-preserving}, having oriented the latter via the recipe given by the analytic gluing map.\\

Suppose that $W_1$ is a cylinder. Then if $\overline{\mathcal M}'_{\mathbf{E}_1, z}(\alpha, \beta)$ is the `compactification' of $\widetilde{\mathcal M}_{\mathbf{E}_1, z}(\alpha, \beta)$ by broken trajectories,\footnote{To make sense of this, the broken trajectories are not parameterized on each component, but rather the stratum corresponding to a $k$-broken trajectory is quotiented by the action of $\mathbb R^{k-1} \subset \mathbb R^k$, sitting inside as the subset with zero sum.} quotienting by $\mathbb R$ gives the usual $\overline{\mathcal M}_{\mathbf{E}_1, z}(\alpha, \beta)$ of unparameterized broken trajectories. We always orient this so that the diffeomorphism $$\mathbb R \times \overline{\mathcal M}_{\mathbf{E}_1, z, k,\delta}(\alpha, \beta) \cong \overline{\mathcal M}'_{\mathbf{E}_1, z, k,\delta}(\alpha, \beta)$$ changes orientation by a factor of $(-1)^{\dim \alpha}$. This convention is chosen so that the corresponding diffeomorphism between \emph{fibers above $\alpha$} are orientation preserving. When later discussing Floer homology, this is the more natural convention; it is the convention which reduces to the standard one for unframed moduli spaces of irreducible instantons.\\

The notation in the following is slightly different than in Proposition \ref{gluing}, to allow for a more uniform discussion.

\begin{proposition}\label{boundary-orientation}If $(Y,E)$ is a 3-manifold equipped with regular perturbation, then if $\overline{\mathcal M}_{E,k,\delta}(\alpha, \beta)$ is the compactified moduli space of unparameterized flowlines, then Proposition \ref{gluing} gives a decomposition of the boundary $$\partial \overline{\mathcal M}_{E,z,\pi}(\alpha, \beta) = \bigcup_{\substack{\gamma \in \mathfrak C_\pi\\ z_1 \ast z_2 = z}} \overline{\mathcal M}_{E,z_1,\pi}(\alpha, \gamma) \times_\gamma \overline{\mathcal M}_{E,z_2,\pi}(\gamma, \beta).$$ Suppose we orient $\Lambda_z(\alpha, \beta)$ using the analytic gluing map and fixed orientations of the three of $\Lambda_{z_1}(\alpha, \gamma), \Lambda_{z_2}(\gamma, \beta)$, and $\mathfrak g_\gamma$. If $d$ is the dimension of $\overline{\mathcal M}_{E,z_1,\pi}(\alpha, \gamma)$, then the boundary orientation on $$\overline{\mathcal M}_{E,z_1,\pi}(\alpha, \gamma) \times_\gamma \overline{\mathcal M}_{E,z_2,\pi}(\gamma, \beta)$$ differs from the fiber product orientation by a factor of $(-1)^{d+1}$.

Now let $(W,\mathbf{E})$ be a cobordism equipped with regular perturbation. By Proposition \ref{gluing}, each compactified moduli space $\overline{\mathcal M}_{\mathbf{E},z,\pi}$ has boundary given as the union \begin{align*}\partial \overline{\mathcal M}_{\mathbf E,z,\pi}(\alpha, \beta) &= \bigcup_{\substack{\gamma \in \mathfrak C_{\pi_1}\\ z_1 \ast z_2 = z}} \overline{\mathcal M}_{E_1,z_1,\pi_1}(\alpha, \gamma) \times_\gamma \overline{\mathcal M}_{\mathbf{E},z_2,\pi}(\gamma, \beta)\\
&\bigcup_{\substack{\zeta \in \mathfrak C_{\pi_2}\\ z_1 \ast z_2 = z}} \overline{\mathcal M}_{\mathbf{E},z_1,\pi} (\alpha, \zeta) \times_\zeta \overline{\mathcal M}_{E_2,z_2,\pi_2}(\zeta, \beta).\end{align*} This decomposes the boundary into two types of components: whether breaking occurs at the negative end or at the positive end.

Then if we orient $\Lambda_z(\alpha, \beta)$ using the analytic gluing map and fixed orientations of the three of $\Lambda_{z_1}(\alpha, \gamma), \Lambda_{z_2}(\gamma, \beta)$, and $\mathfrak g_\gamma$, the orientation of $$\overline{\mathcal M}_{E_1,z_1,\pi_1}(\alpha, \gamma) \times_\gamma \overline{\mathcal M}_{\mathbf{E},z_2,\pi}(\gamma, \beta)$$ as a boundary stratum disagrees with the orientation induced by the fiber product by a factor of $(-1)^{\dim \alpha}$. If the dimension of $\overline{\mathcal M}_{\mathbf{E},z_1,\pi}$ is $d$, then the orientation of $$\overline{\mathcal M}_{\mathbf{E},z_1,\pi} (\alpha, \zeta) \times_\zeta \overline{\mathcal M}_{E_2,z_2,\pi_2}(\zeta, \beta)$$ as a boundary stratum is $(-1)^{d+1}$ that of the orientation given as a fiber product.
\end{proposition}

\begin{proof}The argument is no different from \cite[Proposition~20.5.2]{KMSW} and \cite[Proposition~25.1.1]{KMSW}, respectively. We write $\widetilde{\mathcal M}^0_{E, z, \pi}(\alpha, \beta)$ for the quotient of the space of parameterized flowlines by the $\mathbb R$ action, not compactified. This is (locally!) oriented diffeomorphic to $$\alpha \times \Bbb R \times M^{0,\text{fib}}_{\alpha \beta},$$ where we suppress much of the notation for convenience; this is the fiber above a point in $\alpha$ of the $\Bbb R$-reduced moduli space. 

In the first case, the sign arises from the orientation-preserving local homeomorphism $$\widetilde{\mathcal M}_{E,z_1,\pi}(\alpha, \gamma) \times_\gamma \times \widetilde{\mathcal M}_{E,z_2,\pi}(\gamma, \beta) \to \widetilde{\mathcal M}_{E,z,\pi}(\alpha, \beta).$$ In terms of the decomposition above, this (locally) provides an oriented diffeomorphism $$\alpha \times \Bbb R \times \widetilde{\mathcal M}^{0,\text{fib}}_{\alpha \gamma} \times_\gamma \gamma \times \Bbb R \times \widetilde{\mathcal M}^{0,\text{fib}}_{\gamma \beta} \to \alpha \times \Bbb R \times \widetilde{\mathcal M}^{0,\text{fib}}_{\alpha \beta}.$$ Passing to the $\Bbb R$-reduced moduli spaces (and recalling that this incurs a sign change of $(-1)^{\dim \alpha}$ on both sides) provides an oriented local diffeomorphism $$\alpha \times \widetilde{\mathcal M}^{0,\text{fib}}_{\alpha \gamma} \times_\gamma \gamma \times \Bbb R \times \widetilde{\mathcal M}^{0,\text{fib}}_{\gamma, \beta} \to \alpha \times \widetilde{\mathcal M}^{0,\text{fib}}_{\alpha \beta} = \widetilde{\mathcal M}^0_{\alpha \beta}.$$ 

Passing to the boundary component corresponding to broken trajectories corresponds to translating the second instanton to $+\infty$, which corresponds to taking $t \to -\infty$ in the $\Bbb R$ factor above by Remark \ref{rmk:t-action}. Therefore, the boundary orientation corresponds to passing to the end of this moduli space corresponding to the negative end of $\Bbb R$. To compute this orientation, shuffle the $\Bbb R$ factor to the front; this introduces a sign change of $(-1)^d$. Finally, $-\infty$ is the negative end of $\Bbb R$, so introduces an extra $-1$ to the orientation comparison. 

The argument in the case of a cobordism is essentially identical: shuffle the relevant $\Bbb R$-factor to the front, and pass to the corresponding end. In the first case, we want to translate the first instanton all the way to $-\infty$, so this involves shuffling $\Bbb R$ across the factor of $\alpha$ and restricting to the $\{+\infty\}$ end, giving the sign. In the second case, we want to translate the second instanton to $+\infty$, so this involves shuffling $\Bbb R$ across the moduli space of $W$ and passing to $\{-\infty\}$, giving the stated sign.
\end{proof}

\section{Canonical orientations for determinant lines}\label{sec:or-canonical}
What remains is to find a recipe so that, given some choices for each $\alpha$ and $\beta$, and a choice depending only on the underlying cobordism (in the case that it is not just a cylinder), we are given natural orientations of each $\Lambda^W_z(\alpha, \beta)$. These must compose approproiately under gluing. If we change one of the choices for $\alpha$ then the sign of the orientation on each $\Lambda_z(\alpha, \beta)$ should change, and similarly with swapping the choice for $\beta$. If one changes the choice for $W$, then every orientation should change uniformly. This assignment of orientations is what we call in the title of this section `canonical orientations'.

To do this, we will ultimately require that the cobordism (or possibly some larger cobordism) admits a reducible connection. Already this places the demand that $\beta w_2(E) = 0$, but this is not enough to pin down orientations.

Supposing that $\mathbf{A}$ is an $SO(2)$-reducible, we have a splitting $Q_{\mathbf{A},\pi} \cong Q_{\eta, \pi} \oplus Q_{\theta, \pi}$. If we chose a complex structure on $\eta$, then the term $Q_{\eta, \pi}$ is complex linear with respect to this complex structure, which would give an isomorphism 
$$\det(Q_{\mathbf{A},\pi}) \cong \det(Q_{\theta, \pi}),$$ 
thus reducing discussion to the case of the trivial connection (possibly perturbed), but a different choice of complex structure will give a different isomorphism. We need to ensure that a choice of such a complex structure is part of our data.

For the first time in this text we need to choose (for each 3-manifold $Y$ and each cobordism $W$) not an $SO(3)$-bundle over the manifold, but rather a $U(2)$-bundle. An $SO(3)$-bundle has a lift to a $U(2)$-bundle if and only if $\beta w_2(E) = 0 \in H^3(-;\mathbb Z)$. If $\widetilde{\mathbf{E}}$ is a rank 2 complex vector bundle with fixed connection on its determinant line, then a reducible connection $\mathbf{A}$ induces a splitting $\widetilde{\mathbf{E}} \cong \eta \oplus \zeta$, where $\eta \otimes \zeta \cong \det(\widetilde{\mathbf{E}})$ and $\eta, \zeta$ are complex line bundles. Thus for the trace-$0$ adjoint bundle we obtain an isomorphism $\mathfrak g_{\widetilde{\mathbf{E}}} \cong \mathbb R \oplus (\eta \otimes \zeta^{-1})$, respecting the splitting induced by $\mathbf{A}$; thus we have fixed a complex structure on the second component of the above splitting.

We begin this with some preliminaries.

\begin{lemma}Let $(W,\mathbf{E})$ be a complete Riemannian manifold equipped with a $U(2)$-bundle with no boundary components and some number of cylindrical ends, modelled on either $(-\infty, 0] \times Y_i$ or $[0, \infty) \times Y_i$; in the former case we say $Y_i$ is a negative end, in the latter a postive end. Suppose $W$ is equipped with a perturbation which is regular on each end, and a choice of $\alpha_i \in \mathfrak C_\pi(Y_i)$ for each 3-manifold $Y_i$ the ends are modelled on. Write $\widetilde{\mathcal B}^e_{\mathbf{E},k,\delta}(\alpha)$ for the configuration space of $L^2_{k,\delta}$ connections, asymptotic to the $\alpha_i$ on the corresponding ends.

Then $\pi_0 \mathcal B^e_{\mathbf{E},k,\delta}(\alpha) \cong \mathbb Z$, this isomorphism affine over $\pi_1 \mathcal B^e_{\mathbf{E}_i} \cong \mathbb Z$ for any end $(Y_i, \mathbf{E}_i)$.

If all of the $\alpha_i$ are trivial, then there is a unique component $z \in \pi_0 \mathcal B^e_{\mathbf{E},k,\delta}(\alpha)$ which supports a reducible connection.
\end{lemma}

\begin{proof}The statement about components in the space of connections is little more than the classification of $U(2)$-bundles over a compact 4-manifold with boundary with fixed isomorphism class on the boundary in terms of their first Chern class $c_1 \in H^2(W;\mathbb Z)$ and Pontryagin class $p_1 \in \Lambda \subset H^4(W, \partial W; \mathbb R)$. Here $\Lambda$ is a subset affine over $8\pi^2 \mathbb Z$, the latter defined by a curvature integral with respect to some connection with fixed boundary components. Gluing in the nontrivial positive (determined by the same curvature integral) generator of $\pi_1 \mathcal B^e_{\mathbf{E}_i}$ increases this by $8\pi^2$, as expected. We denote this operation as $z \mapsto z+1$.

The fact that one and only one component supports a reducible connection follows from the enumeration of reducible components in Proposition \ref{action3} (or rather, a version allowing more ends).
\end{proof}

\begin{lemma}In the situation above, there is a canonical isomorphism $\Lambda_z^W(\alpha) \cong \Lambda_{z+1}^W(\alpha)$, compatible with the gluing maps.

Suppose all of the $\alpha_i$ are trivial, and $z$ is the corresponding component supporting a reducible connection, and the perturbation is zero on the ends of $W$. Then if we write $$Q^W_\theta: \Omega^1_{k,\delta}(W;\mathbb R) \to \Omega^{2,+}_{k,\delta}(W, \mathbb R) \oplus \Omega^0_{k,\delta}(W; \mathbb R),$$ for $(d^+, d^*)$, we have an isomorphism $\Lambda_z^W(\alpha) \cong \Lambda^W(\theta)$, where the latter is the set of orientations of $Q^W_\theta$. This is true for any choice of Sobolev indices on the various ends, so long as they never coincide with eigenvalues of $\widehat{\textup{Hess}}_{\alpha_i}$ or $\widehat{\textup{Hess}}_\theta$.
\end{lemma}
\begin{proof}The first isomorphism follows by gluing charge-$1$ instantons on $S^4$ into the cobordism. By the assumption that the perturbation is compactly supported in $W$, we may modify $Q^W_{\mathbf{A},\pi}$ by a homotopy through Fredholm operators to $Q^W_{\mathbf{A}}$.

We discussed the second claim above for $SO(2)$-reducible connections on $SO(3)$-bundles. When we have a lift of the $SO(3)$-bundle to a $U(2)$-bundle, if the reducible corresponds to some splitting $\mathbf{E} \cong \mathbb R \oplus \xi$, then we obtain a canonical complex orientation on the complement of the $\mathbb R$-factor, giving the desired claim. For $SO(3)$-reducible connections, instead we have $Q^W_{\mathbf{A}} = Q^W_\theta \otimes \mathfrak g$, and so an orientation of one canonically induces an orientation of the other (supposing an orientation of $\mathfrak g$ is given).
\end{proof}

Thus we drop $z$ from the notation $\Lambda^W(\alpha, \beta)$.

If one of $\alpha$ or $\beta$ is irreducible, we cannot appeal to reducible connections on the cobordism as above. So we need to 'cap off' the two ends, as appropriate. Being precise about this is a little intricate, as described in the following definition.

\begin{definition}Let $(W,\widetilde{\mathbf{E}})$ be a cobordism $(Y_1, E_1) \to (Y_2, E_2)$ equipped with $U(2)$-bundle, a choice of perturbation on the ends (possibly zero, allowing for it to fail to be regular), and a choice of critical orbits $\alpha, \beta$ on the negative and positive end.

The symbol $\Lambda^W(\alpha, \beta)$ always means the usual two-element set of orientations of $\det(Q_{\mathbf{A},\pi})$, where $\mathbf{A}$ is a connection asymptotic to $\alpha$ and $\beta$ on the appropriate ends, and $Q_{\mathbf{A},\pi}$ has domain $\Omega^1_{k,\delta}$ and codomain $\Omega^{2,+}_{k-1,\delta} \oplus \Omega^0_{k-1,\delta}$. (If one of these terms $\alpha, \beta$ is trivial, we sometimes do not write the corresponding $\theta$.) When unadorned, the symbol $\Lambda^W(\theta)$ denotes the same for the trivial connection on the trivial bundle.

If $(W,\widetilde{\mathbf{E}})$ is a cobordism from some $(Y_-, \textup{triv})$ to $(Y_1, E_1)$, for which the perturbation is zero on the incoming end, we say $W$ is an \emph{incoming cap}; if $(Y_2, E_2)$ is the incoming end and $(Y_+, \textup{triv})$ is the outgoing end, with zero perturbation on the outgoing end, we call $(W,\widetilde{\mathbf{E}})$ an \emph{outgoing cap}, and otherwise call $W$ an \emph{intermediate cobordism}.

If $W$ is an incoming, intermediate, or outgoing cobordism respectively, we write $\Lambda^W_-(\theta), \Lambda_i^W(\theta), $ and $\Lambda^W_+(\theta)$ for the two-element sets of orientations of the following three operators:

\begin{align*}Q^-_{\theta}: \Omega^1_{\delta, -\delta} &\to (\Omega^{2,+} \oplus \Omega^{0})_{\delta}\\
Q^i_\theta: \Omega^1_{-\delta} &\to (\Omega^{2,+} \oplus \Omega^{0})_{-\delta, \delta}\\
Q^+_\theta: \Omega^1_{-\delta,\delta} &\to (\Omega^{2,+} \oplus \Omega^0)_{-\delta, \delta}\end{align*}

Here the subscripts indicate the Sobolev weights on the ends; if only one subscript appears, it is the Sobolev weight on both ends, and if two appear, they are the Sobolev weights on the negative and positive ends, in that order.

Whenever we have two-element sets $\Lambda$ and $\Lambda'$, we write $\Lambda \Lambda'$ to denote $\Lambda \times_{\mathbb Z/2} \Lambda'$, the $\mathbb Z/2$ action the canonical free involution on both sets.
\end{definition}

The point of the choices of Sobolev indices is that these operators (for the unperturbed trivial connection) enjoy an additivity property, immediate from spectral flow arguments.

\begin{lemma}Let $W_0, W, W_1$ denote incoming, intermediate, and outgoing cobordisms equipped with $U(2)$-bundles, as above. If $\hat W$ denotes the composite of these, then we have canonical isomorphisms $$\Lambda^{W_0}(\alpha) \Lambda(\mathfrak g_\alpha) \Lambda^W(\alpha, \beta) \Lambda(\mathfrak g_\beta) \Lambda^{W_1}(\beta) \cong \Lambda^{\hat W}(\theta)$$ and $$\Lambda^{W_0}_-(\theta) \Lambda^W_i(\theta) \Lambda^{W_1}_+(\theta) \cong \Lambda^{\hat W}(\theta).$$

Given a 3-manifold with $U(2)$-bundle $(Y, \tilde E)$, given any two incoming caps $(W_0, \widetilde{\mathbf{E}})$ and $(W_0', \widetilde{\mathbf{E}}')$, the sets $\Lambda(\mathfrak g_\alpha) \Lambda^W(\alpha) \Lambda^W_-(\theta)$ and $\Lambda(\mathfrak g_\alpha) \Lambda^{W'}(\alpha) \Lambda^{W'}_-(\theta)$ are canonically isomorphic. The appropriate modification is true for the outgoing caps, as well.
\end{lemma}
\begin{proof}These follow from the usual gluing of operators on manifolds with cylindrical ends (for which the Sobolev weights on the ends match up appropriately, or there is a corresponding term to account for the spectral flow). The composite of the terms in the first displayed equation are naturally isomorphic, rather, to the orientation set of a nontrivial connection over $\hat W$, but as above we may reduce this to the case of a reducible connection and then to the trivial connection, as above.

Once we have this, using twice the orientation-reversal $\overline{W_0}$ as a positive cap, using these canonical isomorphisms we further get canonical isomorphisms $$\Lambda(\mathfrak g_\alpha) \Lambda^{W_0}(\alpha) \Lambda^{W_0}_-(\theta) \cong \Lambda^{\overline{W_0}}(\alpha)  \Lambda^{\overline{W_0}}_+(\theta) \cong \Lambda(\mathfrak g_\alpha) \Lambda^{W'_0}(\alpha) \Lambda^{W'_0}_-(\theta),$$ as desired.
\end{proof}

The collection of isomorphisms above between the different possible 2-element sets $\Lambda(\mathfrak g_\alpha) \Lambda^W(\alpha) \Lambda^W_-(\theta)$ determine an equivalence relation on the set $$\bigsqcup_{\substack{(W, \widetilde{\mathbf{E}})\\ $W$ \text{ incoming cap}}} \Lambda(\mathfrak g_\alpha) \Lambda^{W_0}(\alpha) \Lambda^{W_0}_-(\theta),$$ whose quotient is a 2-element set we write $\Lambda_-(\alpha)$. Similarly the 2-element set $\Lambda_+(\alpha)$ is a quotient of various copies of $\Lambda^{W_1}_+(\theta)\Lambda^{W_1}(\alpha)\Lambda(\mathfrak g_\alpha)$ as $W_1$ varies over positive caps.

We thus have from the first part of the above lemma a natural isomorphism $$\Lambda_-(\alpha) \Lambda^W_i(\theta) \Lambda_+(\beta) \cong \Lambda^W(\alpha, \beta).$$ This was essentially our goal, though we are not quite done.

\begin{lemma}There is a canonical isomorphism $\Lambda_+(\beta) \cong \Lambda(\mathfrak g_\beta)\Lambda_-(\beta)$.
\end{lemma}
\begin{proof}Writing this out explicitly, we have chosen incoming and outgoing caps $W_0, W_1$, we are trying to construct a trivialization of $$\Lambda^{W_1}_+(\theta)\Lambda^{W_1}(\beta) \Lambda(\mathfrak g_\beta) \Lambda^{W_0}(\beta) \Lambda^{W_0}_-(\theta),$$ here exploting the isomorphism between orientations on a vector space and their dual to reverse the order of the first few terms. By gluing the middle three terms, we obtain an isomorphism to $\Lambda^{W_1}_+(\theta) \Lambda^{\hat W}(\theta) \Lambda^{W_0}_-(\theta)$; splitting up the middle term, this is isomorphic to $\Lambda^{W_1}_+(\theta) \Lambda^{W_0}_-(\theta) \Lambda^{W_1}_+(\theta) \Lambda^{W_0}_-(\theta)$. As this takes the form $\Lambda^2$ for some 2-element set $\Lambda$, it is canonically trivial.
\end{proof}

Therefore, if we simply write $\Lambda(\alpha) := \Lambda_+(\alpha)$, we have above found a canonical isomorphism $\Lambda(\mathfrak g_\alpha)\Lambda(\alpha) \Lambda^W_i(\theta) \Lambda(\beta) \cong \Lambda^W(\alpha, \beta)$; then we obtain the gluing isomorphism 
$$\Lambda^W(\alpha, \beta) \Lambda(\mathfrak g_\beta) \Lambda^{W'}(\beta, \gamma) \cong \Lambda^{W \circ W'}(\alpha, \gamma)$$ 
simply by paing off adjacent like terms.

The last thing to be clear about is precisely what the middle orientation set $\Lambda^W_i(\theta)$ is.

\begin{definition}\label{homology-or}A \emph{homology orientation} of $W$, a cobordism with cylindrical ends and incoming end $(-\infty, 0] \times Y_1$, is an orientation of the real vector space $H^1(W) \oplus H^{2,+}(W) \oplus H^1(Y_1)$.
\end{definition}
Because $W$ and $Y$ are connected, note that we have a canonical isomorphism $H^0(W) \oplus H^0(Y_1) \cong \mathbb R^2$, and in particular carries a canonical orientation. Because we have 
\begin{align*}\text{ker}(Q^i_\theta) &\cong H^1(W)\\
\text{coker}(Q^i_\theta) \cong H^{2,+}(W) \oplus &H^1(Y_1) \oplus H^0(W) \oplus H^0(Y_1),\end{align*}
a homology orientation canonically induces an orientation of $\Lambda^W_i(\theta)$.

\begin{remark}Homology orientations have a natural associative composition law, given by the index-gluing described above. However, the composition law may be described explicitly, without passing to a discussion of Fredholm operators. An explicit formula for this law was found in \cite{scaduto2015instantons}, where an explicit understanding was crucial to discuss the differentials in a spectral sequence to Khovanov homology. We do not need this here, and so will not discuss Scaduto's results in any more detail.
\end{remark}

We assemble the content of this section into a proposition.

\begin{proposition}\label{orientations}For any 3-manifold $Y$ equipped with $U(2)$-bundle $\tilde E$ and regular perturbation $\pi$, there are canonical 2-element sets $\Lambda(\alpha)$ for each critical orbit $\alpha \in \mathfrak C_\pi$, and we have the canonical isomorphism $\Lambda(\alpha, \beta) \cong \Lambda(\mathfrak g_\alpha)\Lambda(\alpha) \Lambda(\beta)$. Therefore, the moduli space $\widetilde{\mathcal M}_z(\alpha, \beta)$ may be given an orientation if we choose an element of $\Lambda(\alpha)$ and $\Lambda(\beta)$, and an orientation of the orbit $\alpha$. Choosing the other element of either set will negate the orientation of $\widetilde{\mathcal M}(\alpha, \beta)$.

If one fixes a choice of element of each of $\Lambda(\alpha), \Lambda(\beta), $ $\Lambda(\gamma)$, and an orientation of $\alpha$, then the orientation these induce on $\overline{\mathcal M}_{E,\pi}(\alpha, \gamma) \times_{\gamma} \overline{\mathcal M}_{E,\pi}(\gamma, \beta)$ via the fiber product differs from the orientation induced as a component of $\partial \overline{\mathcal M}_{E,\pi}(\alpha, \beta)$ by a sign of $(-1)^{d+1}$, where $d$ is the dimension of $\widetilde{\mathcal M}_{E,\pi}(\alpha, \gamma)$. 

Let $(W,\widetilde{\mathbf{E}})$ be a cobordism from $(Y_1, E_1)$ to $(Y_2, E_2)$, equipped with a regular perturbation $\pi$ that restricts to $\pi_i$ on the ends. Then a choice of element of each $\Lambda(\alpha)$ and $\Lambda(\beta)$, as well as a homology orientation of $W$ and an orientation of $\alpha$, give rise to an orientation of $\overline{\mathcal M}_{\mathbf{E},\pi}(\alpha, \beta)$; swapping any one of these elements will negate this orientation.

The boundary components of $\overline{\mathcal M}_{\mathbf{E},\pi}(\alpha, \beta)$ arise in two pieces: those of the form 
$$\overline{\mathcal M}_{E_1,z_1,\pi_1}(\alpha, \gamma) \times_\gamma \overline{\mathcal M}_{\mathbf{E},z_2,\pi}(\gamma, \beta),$$ 
where $\gamma \in \mathfrak C_{\pi_1}$, and those of the form 
$$\overline{\mathcal M}_{\mathbf{E},z_1,\pi}(\alpha, \zeta) \times_\zeta \overline{\mathcal M}_{E_2, z_2, \pi_2}(\zeta, \beta),$$ 
where $\zeta \in \mathfrak C_{\pi_2}$.

Given a choice of element of each $\Lambda(\alpha), \Lambda(\beta), \Lambda(\gamma),$ and $\Lambda(\zeta)$ as above, as well as a choice of homology orientation on $W$ and orientation of $\alpha$, the orientations on $\partial \overline{\mathcal M}$ and on the fiber products differ by a factor of $(-1)^{\dim \alpha}$ in the first case; they differ by a sign of $(-1)^{d_W+1}$ in the second case, where $d_W$ is the dimension of $\widetilde{\mathcal M}_{\mathbf{E},z_1,\pi}(\alpha, \zeta)$.

Similarly, if $S$ is a 1-parameter family of metrics on $W$ abutting to a broken metric, the decomposition of the boundary of $\overline{\mathcal M}_{\mathbf{E},z,S}(\alpha, \beta)$ of Proposition \ref{famgluing} gives the boundary components orientations which differ by the orientation induced by analytic gluing by a factor of $(-1)^{\dim \alpha}, (-1)^{d_1+1}, (-1)^{d_{12}}, $ and $(-1)^{\dim \alpha - 1}$, respectively (the last term coming from the boundary orientation of $\{0\} \subset [0,1]$), where $d_1$ is the dimension of $\overline{\mathcal M}_{E_-, z_1}(\alpha, \gamma)$, while $d_{12}$ is the dimension of $\overline{\mathcal M}_{\mathbf{E},z_1,S}(\alpha, \eta)$. 
\end{proposition}

\part{Equivariant homology from instantons}\label{part2}


\chapter{Floer homology}\label{chap:6}
\section{Geometric homology}\label{sec:Floer-gchain}
It will be useful in what follows to have a chain complex computing the singular homology of a smooth manifold $X$, whose generators are smooth maps from oriented manifolds with corners (as opposed to continuous maps from simplices). Fix a principal ideal domain $R$. 

The definitions below are modeled on \cite[Chapter~4.1]{Lin}, which in turn follow the simpler \cite{Lip}. For technical reasons, we must introduce the notion of strong $\delta$-chain. The reader will not be led astray in what follows by pretending every occurrence of ``strong $\delta$-chain" means ``compact smooth manifold with corners"; if we could achieve that level of smoothness on the instanton moduli spaces, this sequence of definitions would not be necessary to set up our homology theory. 

\begin{definition}A $d$-dimensional strong $\delta$-chain is a compact topological space $P$ with a stratification $$P^d \supset P^{d-1} \supset \cdots \supset P^0 \supset \varnothing$$ by closed subsets, so that $P^e \setminus P^{e-1}$ is decomposed as a finite disjoint union of \emph{smooth} manifolds of dimension $e$, written $\sqcup_{i=1}^{m_e} M_i^e$; the top stratum $P \setminus P^{d-1}$ has only one open face $M_1^d = P \setminus P^{d-1}$ in its decomposition. (Note that $M_i^e$ need not be connected!) We denote the closure of any one of these manifolds an $e$-dimensional \emph{face}, and write it as $\Delta$. We write $\Delta^\circ$ for the interior of a face (one of the manifolds $M_i^e$ in the disjoint union $P^e \setminus P^{e-1}$) and call it an \textup{open face}.

We demand that whenever a codimension $e$ face $\Delta_0$ is contained in a codimension $(e-2)$ face $\Delta_2$, there are exactly two codimension $(e-1)$ faces $\Delta'$ with $\Delta_0 \subset \Delta' \subset \Delta_2.$

Whenever $\Delta_1 \subset \Delta_2$, we assign a set $n(\Delta_1, \Delta_2)$ (which we will write as $n$, or $n_{12}$ when the faces are not implicit) with cardinality $\dim \Delta_1 - \dim \Delta_2$, an open neighborhood $\Delta_1^\circ \subset W(\Delta_1, \Delta_2) \subset \Delta_2$, and a map $r: W \to [0,\varepsilon)^n$ so that $\Delta_1^\circ = r^{-1}(0)$.

We also assign a space $W \hookrightarrow EW(\Delta_1, \Delta_2)$, where $EW$ is a \emph{topological manifold with corners and a smooth structure on each stratum} of dimension $d+m(\Delta_1, \Delta_2)$, for which the map from $W$ is smooth on each stratum, and equipped with a map $\tilde r: EW(\Delta_1, \Delta_2) \to [0,\varepsilon)^n$ extending $r: W \to [0,\varepsilon)^n$. We demand the following.
\begin{itemize}
\item There is a vector bundle $V(\Delta_1, \Delta_2) \to EW$ of rank $m(\Delta_1, \Delta_2)$, and a section $\sigma$ of $V$, smooth and transverse to the zero section on each stratum, so that $\sigma^{-1}(0) = W$.
\item Write $EW^k$ for the inverse image of $R_k \subset [0,\varepsilon)^n$, the set of points for which exactly $k$ coordinates are zero; we demand that the map $EW^k \to R_k$ is a smooth submersion.
\item $\tilde r$ is a fiber bundle projection. The restriction of $\sigma$ to $\tilde r^{-1}(0)$ is transverse to the zero section, whose fiber above zero is $\Delta_1^\circ$. In particular, $EW$ is diffeomorphic to $[0,\varepsilon)^{n_{12}} \times \tilde r_{12}^{-1}(0)$, and $\tilde r_{12}^{-1}(0)$ is diffeomorphic to a neighborhood of the zero section in the restriction of $V$ to $\Delta_1^\circ$.
\end{itemize}

These are compatible in the following sense. Associated to a sequence of inclusions $\Delta_1 \subset \Delta_2 \subset \Delta_3$ we have an inclusion of sets $n(\Delta_1, \Delta_2) \hookrightarrow n(\Delta_1, \Delta_3)$, an embedding $EW(\Delta_1, \Delta_2) \hookrightarrow EW(\Delta_1, \Delta_3),$ as well as an embedding of vector bundles $V(\Delta_1, \Delta_2) \hookrightarrow V(\Delta_1, \Delta_3)$ covering this; the section $\sigma_{12}$ is the restriction of $\sigma_{13}$. The map $n_{12} \to n_{13}$ of sets induces a stratum-preserving embedding $[0,\varepsilon)^{n_{12}} \hookrightarrow [0,\varepsilon)^{n_{13}}$ so that, with respect to this inclusion, $\tilde r_{12}$ is the restriction of $\tilde r_{13}$. Finally, we demand that with respect to these embeddings, $$EW(\Delta_1, \Delta_2) = \tilde r_{13}^{-1}\left([0,\varepsilon)^{n_{12}}\right).$$
\end{definition}

This complicated definition is in fact more or less forced on us by a few simple requirements. First, our chains should include compact topological manifolds with corners with a smooth structure on each stratum (satisfying the combinatorial condition). Second, they should be closed under transverse intersections (in particular, the inverse image of a regular value should be a chain). This is the property we need to ensure that excision holds in our coming homology theory based on maps from strong $\delta$-chains.

Already we run into trouble: because the smooth structure on the strata do not interact, all we see is that the transverse intersection of two topological manifolds with smooth structures on each stratum is a space \textit{stratified} by smooth manifolds. For instance, it is not hard to construct a $\theta$-shaped graph as the zero set of a continuous function $f$ on $[0,1] \times \mathbb R$ which is smooth on each stratum and for which $0$ is a regular value on each stratum. 

\begin{figure}[h]
    \centering
    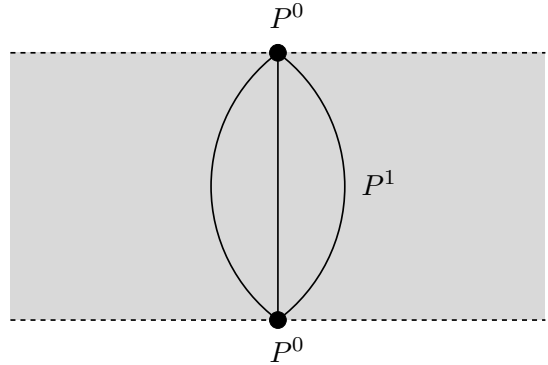
    \caption{The $\delta$-chain $P \subset [0,1] \times \mathbb R$ in this figure is the zero set of a continuous function $f: [0,1] \times \mathbb R \to \mathbb R$ which is smooth and transverse to zero on each stratum. $P^0 = f^{-1}(0) \cap \{0,1\} \times \mathbb R$ is two points, while $P^1 \setminus P^0$ is three open arcs.}
    \label{fig:1}
\end{figure}

This is what forces us to think of spaces equipped with ``local thickenings" which actually \emph{are} manifolds. In the case of the $\theta$-shaped graph, the three open arcs of the graph form the unique top-dimensional open face, whose local thickening is a small open neighborhood of these in the strip; the vector bundle is the trivial line bundle, and $\sigma$ is the function $f$. The stratum $P^0$ consists of two points, for which the neighborhood $W(\Delta_0, \Delta_1)$ is a small neighborhood of one of these points (which looks like a chicken foot). The local thickenings are small 2-dimensional neighborhoods $[0,\varepsilon) \times (-t,t)$ and the map $\tilde r$ is projection to $[0,\varepsilon)$. Again, the vector bundle is the trivial line bundle and $\sigma = f$.

\begin{figure}[h]
    \centering
    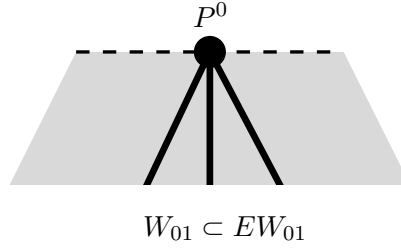
    \caption{The stratum $P^0 \subset P^1$ has a neighborhood $P^0 \subset W(P^0, P^1) \subset P^1$ (the `chicken feet'), along with a further neighborhood $W(P^0, P^1) \subset EW(P^0, P^1)$ so that $EW(P^0, P^1)$ is a smooth surface with corners and $Ef_{01}: EW \to \mathbb R$ is a smooth map with $W_{01} = Ef_{01}^{-1}$. Here, we may take $EW(P^0, P^1)$ to a neighborhood of $W_{01}$ inside of $[0,1] \times \mathbb R$, and $Ef_{01}$ to simply be $f|_{EW}$. In this figure we depict the top half of $EW_{01}$ (which should include a neighborhood of both the top and bottom point).}
    \label{fig:2}
\end{figure}

The third requirement is that faces of a chain should also be chains, and thus we are forced to define local thickenings not just for the whole space, but also for inclusions $\Delta_1 \subset \Delta_2$. 
Finally, it is crucial that the combinatorial boundary operator satisfies $\partial^2 = 0$; this is why we demand the combinatorial condition on faces. At the level of manifolds with corners, it says that a bigon is a manifold with corners, but a teardrop (unigon) is not.

Essentially, the definition of strong $\delta$-chain is one possible approach to formalize a space with local charts modelling a neighborhood of each point as the zero set of a smooth function (smooth section of a vector bundle) on a compact topological manifold with smooth structures on each stratum, and compatibility relations between the charts.

The notion of strong $\delta$-chain here is a slight modification of the notion of $\delta$-chain given in \cite[Definition~24.7.1]{KMSW} and \cite[Definition~4.1.1]{Lin}. 

There are two differences. First, for a strong $\delta$-chain the set $W \subset EW$ is cut out as the zero set of a transverse (on each stratum) section of a vector bundle over $EW$; for a $\delta$-chain, instead the map is simply to some vector space $\mathbb R^k$. The definition of $\delta$-chain is attempting to formalize the kind of spaces that appear as the monopole moduli spaces, while the definition of strong $\delta$-chain is meant to capture all transverse intersections of compact topological manifolds with corners and smooth structures on each stratum. These transverse intersections may be written as the inverse image of a diagonal $\Delta \subset Y \times Y$; the vector bundle in question is $N\Delta$, and if $Y$ has interesting topology this bundle may be nontrivial, so we need to include zero sets of sections of vector bundles. Locally near any point, these are the same notion, but they may not be the same in a neighborhood of an entire stratum.

Second, the map $EW \to \mathbb R^k$ in the definition of $\delta$-chain is required to satisfy certain positivity properties, depending on the values of $\tilde r(EW) \subset [0,\varepsilon)^n$; we do not make any such demand for strong $\delta$-chains. This requirement originates from the boundary-obstructedness phenomenon in monopole moduli spaces, which we do not encounter.

\begin{definition}Suppose we are given the following data. 
\begin{itemize}
\item For each pair of faces $\Delta_1 \subset \Delta_2$, open subsets $U_{12} \subset W_{12}$ and $U'_{12} \subset W'_{12}$, as well as open subsets $EU_{12}, EU_{12}'$ of the extensions, 
\item homeomorphisms $\varphi_{12}: U_{12} \to U_{12}'$, $\tilde \varphi_{12}: EU_{12} \to EU_{12}'$, which are diffeomorphisms over each stratum, 
\item an isomorphism of vector bundles $V_{12} \cong V'_{12}$ over these open subsets $EU$ and $EU'$, covering the diffeomorphism $\tilde \varphi_{12}$,
\end{itemize} 
all of which respect all the associated structure and are compatible for triples $\Delta_0 \subset \Delta_1 \subset \Delta_2$. Then we say the two strong $\delta$-chains are \emph{germ equivalent}.

\end{definition}

This definition fits with the above intuition for strong $\delta$-chains as spaces equipped with a special kind of chart modeling them as a zero set: if we pass to smaller open subsets containing the zero set, this should still present a perfectly good, equivalent chart at a point. 

Next, we will show the important property that strong $\delta$-chains are, in a sense, closed under fiber products. Germ equivalence was partly introduced as a means to this end. To start, we should define what a map from a strong $\delta$-chain is. 

\begin{definition}Suppose $X$ is a smooth manifold, and $P$ is a strong $\delta$-chain. We say a map $f: P \to X$ is a continuous map $f: P \to X$ on the underlying topological space which is smooth on each stratum, and for each inclusion of faces $\Delta_1 \subset \Delta_2$, an extension $Ef_{12}: EW_{12} \to X$ which is smooth on each stratum. Given a sequence of faces $\Delta_0 \subset \Delta_1 \subset \Delta_2$, recall that $EW_{02} \hookrightarrow EW_{12}$; we demand that $Ef_{02}$ is the restriction of $Ef_{12}$.

If $f_1: P_1 \to X$ and $f_2: P_2 \to X$ are two maps from strong $\delta$-chains to $X$, we say that they are \emph{transverse} if they are transverse on each stratum, and for any pair of faces in $P_1$ and any pair of faces in $P_2$, the extensions $Ef_1$ and $Ef_2$ are transverse on each stratum in a neighborhood of $W_1 \times W_2$ in $EW_1 \times EW_2$.

\end{definition}



\begin{proposition}Suppose $f_i: P_i \to X$ are transverse maps from germ equivalence classes of strong $\delta$-chains $P_1, P_2$ to a smooth manifold $X$. Then $P_1 \times_X P_2$ has the natural structure of a germ equivalence class of strong $\delta$-chain.
\end{proposition}


\begin{proof}
The compact space $P_1 \times_X P_2$ is still stratified by smooth manifolds (the fiber products of the original open faces). For a pair of inclusions $\Delta_1^i \subset \Delta_2^i$ of faces of $P_i$, the open set $W$ corresponding to $$\Delta_1^1 \times_X \Delta_1^2 \subset \Delta_2^1 \times_X \Delta_2^2$$ is the fiber product of the corresponding open sets. 

Take $EW$ to be a small open neighborhood of $EW_1 \times_X EW_2 \subset EW_1 \times EW_2$ (well-defined after passing to germs); this is an open subset of a manifold with corners, hence also a manifold with corners. Take this neighborhood to be so small so that $(Ef_1 \times Ef_2)(EW)$ lies in a neighborhood of the diagonal $\Delta_X \subset X \times X$. In particular, $EW$ comes equipped with a map $Eg: EW \to N\Delta_X$ to the normal bundle of the diagonal, so that $(Eg)^{-1}(0) = EW_1 \times_X EW_2$. 

The vector bundle $V$ should be taken to be $V_1 \oplus V_2 \oplus (Ef_1 \times Ef_2)^*(N \Delta_X)$, where $N\Delta_X$ is the normal bundle to the diagonal. The section is given by $\sigma = (\psi_1, \psi_2, Eg)$, whose zero set is canonically identified with $P_1 \times_X P_2$. 

The restriction maps are $r = r_1 \times r_2$ and $\widetilde r = \widetilde r_1 \times \widetilde r_2$. Verifying the axioms is tedious but straightforward, as is verifying that this construction is well-behaved after passing to germ equivalence classes. 
\end{proof}

\begin{remark}\label{rmk:assoc}
It is straightforward but tedious to verify that the above construction is associative when the triple fiber product is defined. 
\end{remark}

To use coefficients other than $\mathbb Z/2$, we must introduce a notion of orientation.
\begin{definition}Let $P$ be a strong $\delta$-chain.

Suppose we have an orientation on the top stratum of $P$, and for a codimension 1 face $\Delta_1 \subset P$, suppose we have oriented the open manifold $\Delta_1^\circ$. We may find an isomorphism $\text{det }(TEW_{12}) \cong \det{V_{12}}$ above $W_{12}^\circ$ using the fiber bundle isomorphism $EW \cong [0,\varepsilon) \times \tilde r_{12}^{-1}(0),$ using that $\tilde r_{12}^{-1}(0)$ is isomorphic to a neighborhood of the zero section of the restriction of $V_{12}$ to $\Delta_1^\circ$, and that we have oriented $\Delta_1^\circ$. Then
$$W_{12}^\circ := r_{12}^{-1}(0, \varepsilon)^n = \tilde r_{12}^{-1}(0, \varepsilon)^n \cap \sigma^{-1}(0)$$
may be oriented as the zero set of the section $EW_{12}^\circ \to V_{12}$ using the above isomorphism of determinant bundles. This gives an orientation on an open subset of $P^d \setminus P^{d-1}$. If the same as the orientation already given, we say that $\Delta_1$ has the \emph{boundary orientation}; we express that this is possible for some orientation on $\Delta_1$ by saying that $\Delta_1$ is \emph{consistently orientable with respect to $P$}.

For any sequence of codimension 1 faces $\Delta_k \subset \cdots \subset \Delta_1 \subset P$, we may thus inductively define the notion of any sequence being consistently orientable with respect to $P$. 

We say that $P$ is an \emph{oriented $\delta$-chain} if $P$ is equipped with an orientation on its top stratum $P^d \setminus P^{d-1}$, satisfying the follwing conditions:
\begin{itemize}
\item Any sequence of codimension 1 faces of $P$ is consistently orientable with respect to $P$,
\item Given a codimension $2$ inclusion $\Delta_0 \subset \Delta_2$, there are two intermediate faces $\Delta_1$ and $\Delta_1'$; we demand that if we fix an orientation on $\Delta_2$ (induced from $P$ by some sequence of codimension 1 faces as above), then the boundary orientation that arises on $\Delta_0$ from the sequence $\Delta_0 \subset \Delta_1 \subset \Delta_2$ is negative that from the sequence $\Delta_0 \subset \Delta_1' \subset \Delta_2$. 
\end{itemize}
\end{definition}

The definition of boundary orientation above is relevant for the following calculation for strong $\delta$-chains. This is proved for $\delta$-chains as \cite[Theorem~21.3.2]{KMSW}; instead of subdividing our strong $\delta$-chains into pieces that may be given the structure of $\delta$-chains in the sense of Kronheimer-Mrowka, instead we spell out their proof in simpler language as an exercise in understanding the definition of strong $\delta$-chains. The proof follows similar lines as a discussion in \cite[Section~2.3]{SS-Invol}.

\begin{lemma}\label{Stokes}If $P$ is a 1-dimensional oriented strong $\delta$-chain, then $P^0$ consists of a finite number of oriented points whose count, signed by the boundary orientation, is equal to zero.
\end{lemma}
\begin{proof}If $x \in P^0$, there is an open subset $W_x \subset P$ with a map to $[0,\varepsilon)$, as well as a larger space $EW_x \cong [0,\varepsilon) \times V_x$ so that $W_x \hookrightarrow EW_x$, with $x$ mapping to $(0,0)$, and $W_x^\circ$ mapping into $(0,\varepsilon) \times V_x$; the neighborhood $W_x \subset P$ is the zero set of a map $\sigma: [0,\varepsilon) \times V_x \to V_x$ which is smooth and transverse to zero on each stratum. Write $\sigma_t(x) = \sigma(t,x)$. That $\sigma$ is transverse to zero on the boundary implies that $\sigma^{-1}(0) \cap \left(\{0\} \times V_x\right)$ is discrete (and in particular discrete inside $W$). In particular, $P^0$ is a discrete set in the compact space $P$, so it is finite. For convenience, after passing to a possibly smaller neighborhood and reparameterizing, we assume $\sigma_0 = \text{Id}$ in our coordinates.

We also see that $\sigma$ does not vanish on a sufficiently small sphere $\{0\} \times S_\delta(V_x)$ in the boundary of the local thickening near $x$; the inclusion map from this sphere to $V_x \setminus \{0\}$ is degree 1 (regardless of how we orient $V_x$, as long as we orient these compatibly). Extending this to a disc $D_x$ in $[0,\varepsilon) \times V_x$ which only intersects $\{0\} \times V_x$ in $S_\delta(V_x)$ and is transverse to $\sigma^{-1}(0)$, we see that $D_x \cap \sigma^{-1}(0)$ is a finite set of points whose oriented sum is 1, orienting them as $\text{det }d(\sigma_t)(y)$ at an intersection point $(t,y)$.

$D_x$ bounds a ball one dimension larger; deleting the part of $W_x$ contained in the interior of this ball for each $x$, what is left of $P^1$ is a compact oriented 1-manifold $L$ with boundary $\bigsqcup_{x \in P^0} \left(D_x \cap \sigma^{-1}(0)\right)$. The boundary orientation convention is that a point $(t,y)$ in this intersection is oriented positively as the boundary of $L$ if the sign of $\text{det}(d\sigma_t)(y)$ agrees with the sign of $x$. In particular, the signed sum of points in $\partial L$ is $$\sum_{x \in P^0} \text{sgn}(x) \#\left(D_x \cap \sigma^{-1}(0)\right) = \sum_{x \in P^0} \text{sgn}(x),$$ where the signed count is as in the previous paragraph. The signed count of points in the boundary of an oriented compact 1-manifold is zero, and thus the desired count is zero.
\end{proof}

Now that these objects have been introduced, we can introduce the \emph{geometric chain complex} computing singular homology. First we will write down the generating set.

\begin{definition}\label{basicchain}Let $X$ be a smooth manifold. Fix once and for all a Hilbert space $H$ equipped with a Hilbert bundle $V_H$.

A \emph{basic chain} of degree $d$ in $X$ is a germ equivalence class of maps $\sigma: P \to X$, where $P\subset H$ is a compact connected subspace equipped with the structure of an oriented strong $\delta$-chain of dimension $d$, whose thickenings $EW_{12}$ are compatibly embedded subspaces of $H$, and whose vector bundles $V_{12}$ are subbundles of $V_H$. 

Two basic chains $\sigma_i: P_i \to X$ are \emph{isomorphic} if there is an orientation-preserving homeomorphism $f: P_1 \to P_2$ which is a diffeomorphism on each stratum, with $\sigma_2f = \sigma_1$, and an extension for each $\Delta_1 \subset \Delta_2$ to $Ef_{12}: EW_{12}^1 \to EW_{12}^2$ which are diffeomorphisms on each stratum compatible with restriction. 
A basic chain $\sigma: P \to X$ is \emph{achiral} if it is isomorphic to itself with the opposite orientation, and \emph{chiral} if it is not achiral. A basic chain $\sigma: P \to X$ has \emph{small image} if there is map $g: Q \to X$ from a $\delta$-chain $Q$ of strictly smaller dimension.
\end{definition}

Write $\widetilde C_*(X;R)$ for the graded vector space freely generated in degree $d$ by the basic chains of degree $d$. One step closer to the desired chain complex, there is a quotient $\overline C_*(X;R)$ obtained from the relations
$$- [\sigma: P \to X] = [\sigma: \overline P \to X],$$
and if $2 \neq 0 \in R$ we further set any achiral basic chains $\sigma$ equal to zero.

This is still a free graded $R$-module, freely generated by isomorphism classes of chiral basic chains, with exactly one choice of orientation for each nontrivial basic chain appearing as a generator. (The achiral chains are precisely the basis elements of $\widetilde C_*(X;R)$ which are sent to zero in $\overline C_*(X;R)$; if $2 = 0 \in R$ then no basis elements are set to zero.)

There is a geometric boundary operator $\widetilde\partial: \widetilde C_*(X;R) \to \widetilde C_{*-1}(X;R)$, sending each (isomorphism class of) basic chain $\sigma: P \to X$ to the sum of the faces of $P^{d-1}$ equipped with their boundary orientation. Because this is compatible with orientation-reversal, it descends to an operator $\overline \partial: \overline C_*(X;R) \to \overline C_{*-1}(X;R)$. 

\begin{lemma}We have $\overline \partial^2 = 0$.
\end{lemma}

\begin{proof}This immediately follows from the combinatorial conditions on the definition of $\delta$-chain and orientation thereof.

Given any face $\Delta_2 \subset P$ of codimension 2, there are two intermediate faces of codimension 1. The boundary of $\partial P$ then contains two copies of the face $\Delta_2$: once from each of the intermediate faces. The assumption on orientations were that these two copies were oriented negatively to one another, so they sum to zero. Therefore, $\overline\partial^2 \sigma = 0 \in \overline C_*(X;R)$. 
\end{proof}

As is familar already from cubical definitions of singular homology, the homology of $\overline C_*(\text{pt};R)$ is not obviously concentrated in degree zero unless we impose some further degeneracy requirements. (These degeneracy requirements will later turn out to be essential in our definition of the instanton Floer complex.) This is furnished by Lipyanskiy's notion of small image. However, the basic chains of small image do not span a subcomplex of $\overline C_*(X;R)$: for instance, a basic chain of dimension larger than $X$ is automatically of small image, but if the chain is of dimension $\dim X + 1$, there is no such guarantee that its boundary has small image. This inspires us to make the following definition.

\begin{definition}\label{degeneracy}A basic chain $\sigma: P \to X$ is \emph{degenerate} if both $\sigma$ has small image and $\partial \sigma$ is a disjoint union of basic chains of small image and an achiral chain. The span of the degenerate basic chains forms a subcomplex $D_*(X;R) \subset \overline C_*(X;R)$; this is precisely the submodule spanned by basic chains with small image and whose boundary in $\overline C_*(X;R)$ is a sum of elements of small image. We define the geometric chain complex of $X$ (with coefficients in $R$) as $$C_*^{\textup{gm}}(X;R) := \overline C_*(X;R)/D_*(X;R).$$
\end{definition}

Again, $C_*^{\textup{gm}}(X;R)$ is degreewise $R$-free. It is functorial under smooth maps, and is supported in degrees $[0, \dim X + 1]$; the map $\partial: C^{\text{gm}}_{\dim X + 1}(X;R) \to C^{\text{gm}}_{\dim X}(X;R)$ is an isomorphism onto its image. It is also well-behaved with respect to transverse intersections.

\begin{definition}\label{trans-chains}Suppose $f_i: P_i \to X$ is a countable family $\mathcal F$ of maps from $\delta$-chains to $X$. The subcomplex of $C_*^{\textup{gm}}(X;R)$ spanned by nondegenerate chains transverse to all of the $f_i$ is written $C_*^{\textup{gm},\mathcal F}(X;R)$.
\end{definition}

The following lemma can be proved as an inductive application of transversality theorems to each stratum, using that $P^e \setminus P^{e-1}$ is a manifold and $P^e$ is compact. A proof in the only slightly different setting of $\delta$-chains (not strong) is given in \cite[Lemma~4.1.14]{Lin}. 

\begin{lemma}\label{qi-trans}The inclusion of $C_*^{\textup{gm},\mathcal F}(X) \hookrightarrow C_*^{\textup{gm}}(X)$ is a quasi-isomorphism.
\end{lemma}

Most importantly for us, these chain complexes have chain-level fiber product maps, as long as we're potentially willing to pass to a quasi-isomorphic subcomplex. First we need an additional definition.

\begin{definition}\label{FOrient}Let $f: Z \to X$ be a map from a strong $\delta$-chain to a smooth manifold $X$ which is a submersion on each stratum. We say that a \emph{fiber orientation} of $f$ is an orientation of the $\delta$-chain $f^{-1}(x)$ for any $x \in X$. 
\end{definition}

\begin{lemma}\label{fibprod}Suppose $X$ is a connected smooth manifold. Let $Z$ be a strong $\delta$-chain, and suppose we have maps $e_-: Z \to X, e_+: Z \to Y$, for which $e_-$ is a submersion on each stratum of $Z$. Further, suppose $e_-$ is equipped with a fiber orientation. Then there is an induced map $C^{\textup{gm}}_*(X;R) \to C^{\textup{gm}}_*(Y;R)$ of degree $(\textup{dim } Z - \textup{dim } X)$ given by sending a basic chain $\sigma: P \to X$ to $e_+: P \times_{e_-} Z\to Y$, where $$P \times_{e_-} Z = \{(p, z) \in P \times Z \mid \sigma(p) = e_-(z)\}$$ and $e_+(p,z) = e_+(z)$. We denote this chain $\sigma \times_{e_-} Z$. Here we orient the fiber product so that $$T_p P \oplus T^{\textup{fib}}_zZ \cong T_{p,z} (P \times_X Z)$$ is an oriented splitting.

If $\dim Z > \dim X + \dim Y + 1$, this chain map is identically zero. These maps satisfy the formula $$\partial(\sigma \times_{e_-} Z) =  \partial \sigma \times_{e_-} Z + (-1)^{\dim P + \dim X}\sigma \times_{e_-} \partial Z.$$
\end{lemma}

The second statement follows because $C^{\textup{gm}}_*(Y;R)$ vanishes in degrees larger than $\dim Y + 1$, and the final statement because taking the boundary in the final case involves commuting the desired outward-pointing normal across the first two factors (of $P$ and the base $X$). 

\begin{remark}\label{rmk:assoc2}
The construction of Lemma \ref{fibprod} is associative, because the fiber product construction is associative. 
\end{remark}

Now we should justify the claim that this is just a fancy way to write down singular homology with desirable chain-level properties. To make sense of the Eilenberg-Steenrod axioms for \emph{smooth manifolds}, we work with the notion of \emph{admissible pairs} $(X,A)$, where $X$ is a smooth manifold without boundary, and $A$ is a closed (in the sense of point-set topology) submanifold of $X$ of codimension zero. These were introduced in \cite{Schwarz} to prove the equivalence between Morse homology and singular homology, and used in \cite[Theorem~4.1.13]{Lin} to prove that a very similar homology theory to ours (using $\delta$-chains, instead of strong $\delta$-chains) agrees with singular homology; our proof follows similar lines as his.

\begin{theorem}\label{eilenberg}The functor $H_*^{\textup{gm}}(X;R)$ from smooth manifolds and maps to graded $R$-modules satisfies the Eilenberg-Steenrod axioms for a homology theory: \begin{enumerate}
    \item The induced map $H_*^{\textup{gm}}(X;R) \to H_*^{\textup{gm}}(Y;R)$ is a homotopy invariant, 
    \item There is a natural relative long exact sequence relating 
    $$H_*^{\textup{gm}}(X,Y) := H\left(C_*^{\textup{gm}}(X)/C_*^{\textup{gm}}(Y)\right)$$ 
    to $H_*^{\textup{gm}}(X)$ and $H_*^{\textup{gm}}(Y)$, 
    \item If $\textup{pt}$ denotes the one-point space, then $H_*^{\textup{gm}}(\textup{pt};R)$ is a copy of $R$ concentrated in degree zero, and
    \item If $A$ is codimension zero, and $Y = \partial A$ is a codimension 1 submanifold of $X$, bounding $A$ on one side and $B$ on the other, then the inclusion induces an isomorphism $H_*^{\textup{gm}}(B, Y) \cong H_*^{\textup{gm}}(X, A)$.
\end{enumerate}
As a result, there is an isomorphism $H_*^{\textup{sing}}(X;R) \cong H_*^{\textup{gm}}(X;R)$ for all $X$, natural under smooth maps.
\end{theorem}
\begin{proof}
We prove each property separately, and then explain why this restricted class of Eilenberg-Steenrod axioms is enough.

Homotopy invariance follows the expected strategy: write down a chain homotopy of the induced map, given by sending a basic chain $P$ to $P \times I$, with map given by composing the map $\sigma: P \to X$ with the homotopy $f_t$.

The relative long exact sequence is a matter of homological algebra (it is induced by a short exact sequence of chain complexes).

To check the first nontrivial axiom, observe that $C_*^{\textup{gm}}(\textup{pt};R)$ is concentrated in degrees 0 and 1. We saw in Lemma \ref{Stokes} that the boundary map $$C_1^{\textup{gm}}(\textup{pt};R) \to C_0^{\textup{gm}}(\textup{pt};R) = R$$ is identically zero; but a chain of small image with zero boundary has boundary of small image, and is in particular degenerate, so $C_1^{\textup{gm}}(\textup{pt};R) = 0$. Therefore 
$$H_*^{\textup{gm}}(\textup{pt}) = C_*^{\textup{gm}}(\textup{pt}) = R.$$
Excision is the hardest property to verify. If we write $C_*^{\textup{gm}}(A \cup B)$ for the image of $C_*^{\textup{gm}}(A) \oplus C_*^{\textup{gm}}(B)$ in $C_*^{\textup{gm}}(X)$, then there is a diagram of exact sequences
$$\begin{CD}0 @>>> C_*^{\textup{gm}}(A) @>>> C_*^{\textup{gm}}(A \cup B) @>>> C_*^{\textup{gm}}(B,Y) @>>> 0 \\
@. @VVV @VVV @VVV @. \\
0 @>>> C_*^{\textup{gm}}(A) @>>> C_*^{\textup{gm}}(X) @>>> C_*^{\textup{gm}}(X,A) @>>> 0.
\end{CD}$$
By the five lemma and the induced map on homology long exact sequences, if we can show the middle vertical arrow is a quasi-isomorphism, so too will be the last arrow, as desired. Now, if we write $C_*^{\textup{gm},Y}(X)$ for the quasi-isomorphic subcomplex of strong $\delta$-chains transverse to the submanifold $Y$, there is a map 
$$\rho: C_*^{\textup{gm},Y}(X) \to C_*^{\textup{gm}}(A \cup B)$$ 
given by sending each basic chain $\sigma: P \to X$ to the sum of $P \cap A$ and $P \cap B$; these are again $\delta$-chains because of the transversality hypothesis, and their boundary is a subdivision of the original boundary into $(\partial P) \cap A$ and $(\partial P) \cap B$, respectively, and so this is a chain map. If $P$ is a cycle, then $\rho(P)$ is homologous to $P$ itself, the bounding chain given by $P \times I$, with the usual $\delta$-structure on one end, and the `broken' $\delta$-structure $P = (P \cap A) \cup (P \cap B)$ on the other. In particular, we see that the inclusion $C_*^{\textup{gm}}(A \cup B) \to C_*^{\textup{gm}}(X)$ is surjective on homology. Injectivity is similar: given a chain whose boundary is transverse to $P$, we can represent it instead by a chain transvere to $P$ with homologous boundary still transverse to $P$; cutting it into two pieces, we get a chain in $C_*^{\textup{gm}}(A \cup B)$ whose boundary is homologous to the original.

Now we should explain why the homotopy category of admissible pairs of smooth manifolds (with the homotopy type of a finite CW complex) is equivalent to the homotopy category of pairs of finite CW complexes. Putting a relative CW structure on each admissible pair via Morse theory, we have a natural inclusion from the category of admissible pairs to the category of finite CW complexes; by smooth and (relative) smooth approximation, this is a fully faithful functor, and it suffices to show that it's surjective on objects. Given a pair $(X,A)$, find an embedding of this in some large Euclidean space; then a small open neighborhood $U_X$ of the image of $X$ (called a regular neighborhood) is homotopy equivalent to $X$ itself and so that $\partial U_X$ is a manifold, and we may choose an open neighborhood $U_A$ of $A$ satisfying the same property, with $\overline{U}_A \subset U_X$. Then the pair $(U_X, \overline{U}_A)$ is the desired admissible pair.
\end{proof}

If $R$ has characteristic 2, we can simplify definitions by removing all reference to orientations; this is even desirable, as it allows us to refer to the identity map of a non-oriented compact manifold as a chain.

\section{Equivariant instanton homology}\label{sec:Floer-eqinst}
We are now ready to define the framed instanton chain complex $\widetilde{CI}(Y,E,\pi;R)$, using the instanton moduli spaces. The material of sections 1-4 amounts to the following package of theorems about these instanton moduli spaces. Recall from Definition \ref{admiss-3m} that for a $U(2)$-bundle $E$ over a 3-manifold $Y$ to be weakly admissible, either $b_1(Y) = 0$ or $c_1(E)$ is not twice some element of $H^2(Y;\mathbb Z)$.

In the following, the regular perturbation $\pi$ is chosen from an open subset $\mathcal P_{E,\delta}$ of the Banach space $\mathcal P_E$ of perturbations, defined in Definition \ref{pertdef}. 

\begin{theorem}\label{ModuliPackage}Let $Y$ be a 3-manifold equipped with a weakly admissible $U(2)$-bundle $E$, a basepoint $b$, and a choice of metric and regular perturbation $\pi$ (which always exists). Then:\begin{enumerate}
\item The collection of \emph{critical orbits} of the perturbed Chern-Simons function $\text{cs}+\pi$ on the configuration space $\widetilde{\mathcal B}_E$ is a finite set of $SO(3)$-orbits. We write this set of orbits as $\mathfrak C_\pi$.\\
\item Write $c_1(E) = \lambda \in H^2(Y; \mathbb Z)$. Define $\textup{Pair}(H^2)$ to be the set of unordered pairs of integral cohomology classes on $Y$; that is, the set $H^2(Y;\mathbb Z) \times H^2(Y;\mathbb Z)$ modulo the relation $(z,w) \sim (w,z)$. The set of reducible critical orbits may be identified with the set $$\mathcal R_\lambda(Y,E) \subset\textup{Pair}(H^2)$$ given by those $(z,w)$ with $z+w = \lambda$. The fully reducible critical orbits are those for which $z = w$, and there are hence $H^1(Y;\mathbb Z/2)$ of them if $Y$ is a rational homology sphere with $c_1(E)$ divisible by $2$, and none otherwise. \\
\item There is a number $\text{\em{gr}}_z(\alpha, \beta) \in \mathbb Z$, assigned to each pair of critical orbits $\alpha, \beta$ and homotopy class of path $z$ between them in $\widetilde{\mathcal B}_E$. If $w$ is a path from $\beta$ to $\gamma$ and $z \ast w$ is the concatenation of paths, we have $$\text{\em{gr}}_{z \ast w}(\alpha, \gamma) = \text{\em{gr}}_{z}(\alpha, \beta) + \text{\em{gr}}_{w}(\beta, \gamma).$$ For different paths $z,w$, we have $\text{\em{gr}}_z(\alpha, \beta) - \text{\em{gr}}_w(\alpha, \beta) \in 8\mathbb Z$, so that $\text{\em{gr}}(\alpha, \beta) \in \mathbb Z/8$ is well-defined.\\
\item If $\alpha$ and $\beta$ are reducible, then $\textup{gr}(\alpha, \beta)$ is even. If $\alpha$ and $\beta$ are fully reducible, then $\textup{gr}(\alpha, \beta)$ is divisible by $4$.\\
\item Associated to each pair $(\alpha, \beta)$ of critical orbits and homotopy class $z$ is a smooth $SO(3)$-manifold (possibly empty) $\widetilde{\mathcal M}^0_{E,z,\pi}(\alpha, \beta)$ of dimension $\text{\em{gr}}_{z}(\alpha, \beta)+\dim \alpha-1$. It comes equipped with equivariant smooth maps $$\alpha \xleftarrow{e_-} \widetilde{\mathcal M}^0_{E,z,\pi}(\alpha, \beta) \xrightarrow{e_+} \beta.$$
\item For each critical orbit $\alpha$, there is an associated 2-element set $\Lambda(\alpha)$. A choice of element of each of $\Lambda(\alpha)$ and $\Lambda(\beta)$ induces an orientation on the fiber of $$e_-: \widetilde{\mathcal M}^0_{E,z,\pi}(\alpha, \beta)$$ for all $z$; in the language of Definition \ref{FOrient}, this is a fiber orientation of $e_-$. Negating either one of these choices negates the resulting fiber orientation.\\
\item If $\textup{gr}_z(\alpha, \beta) \leq 10 - \dim \alpha$, there is a natural compactification $$\widetilde{\mathcal M}^0_{E,z,\pi}(\alpha, \beta) \subset \overline{\mathcal M}_{E,z,\pi}(\alpha, \beta)$$ into a compact topological $SO(3)$-manifold with corners and a smooth structure on each stratum. The endpoint maps extend to equivariant maps from $\overline{\mathcal M}^0$ which are smooth on each stratum; we will use the same notation $e_\pm$ for these extended maps. In the trivial case $\alpha = \beta$ with homotopy class $z = 0$, we consider $\overline{\mathcal M}_{E,0,\pi}(\alpha, \alpha)$ to be empty.\\
\item The action of $SO(3)$ on $\overline{\mathcal M}_{E,z,\pi}(\alpha, \beta)$ is free.\\
\item Given any choice of element of each of $\Lambda(\alpha), \Lambda(\beta)$, and $\Lambda(\gamma)$, so that all of the relevant moduli spaces are fiber oriented, there is an oriented decomposition $$(-1)^{\dim \alpha}\partial \overline{\mathcal M}_{E,z,\pi}(\alpha, \beta) \cong \bigsqcup_{\substack{\gamma, w_1, w_2 \\ w_1 \ast w_2 = z}} (-1)^{\textup{gr}_{w_1}(\alpha, \gamma)} \overline{\mathcal M}_{E,w_1,\pi}(\alpha, \gamma) \times_\gamma \overline{\mathcal M}_{E,w_2,\pi}(\gamma, \beta).$$
\end{enumerate}
\end{theorem}

Observe that the relation $\text{gr}_z(\alpha, \beta) = \dim \overline{\mathcal M}_{E,z,\pi}(\alpha, \beta) - \dim \alpha + 1$ is compatible with the dimensions of the manifolds involved in the gluing formula.

\begin{proof}For a fixed metric $g$, the existence of an admissible perturbation $\pi$ is the combination of Theorem \ref{trans1} and Theorem \ref{trans2}; these guarantee that the critical orbits are isolated, and that the moduli spaces of trajectories between them are smooth manifolds of the appropriate dimension, respectively. That there are only finitely many critical orbits follows from Lemma \ref{prop1} (that the derivative of our perturbed Chern-Simons functional is a proper map). The enumeration of reducible critical orbits is a combination of Corollary \ref{orbits1} and Proposition \ref{redcrit}; as mentioned after Definition \ref{pertdef}, one of the reasons we demand $\pi \in \mathcal P_{E,\delta}$ is so that this enumeration remains true. The integer $\text{gr}_z(\alpha, \beta) = \overline{\text{gr}}_z(\alpha, \beta) - \dim \alpha$ is defined in Definition \ref{grading}, where we also explain why $\overline{\text{gr}}_z$ counts the dimension of the fiber of $\widetilde{\mathcal M} \cong \mathbb R \times \widetilde{\mathcal M}^0$ above a point in $\alpha$. That the grading is well-defined modulo $8$ is Corollary \ref{gr-mod-8}. The calculation of relative gradings between full reducibles is Proposition \ref{gr-fullred} and arbitrary reducibles is Proposition \ref{gr-red}. The existence of a compactification only in terms of fiber products of lower-dimensional moduli spaces is Corollary \ref{nobubbling}; that the resulting object is a topological manifold with corners with a smooth structure on each stratum is Proposition \ref{gluing}. That the $SO(3)$-action is free means that there are no reducible solutions; this is guaranteed by Proposition \ref{red-flowlines}. The moduli spaces carry compatible orientations by Proposition \ref{orientations}, writing the sign in terms of the grading function $\text{gr}_z$.
\end{proof}

We will also need a similar package in the case of cobordisms $W$. We include this here as well. Recall the definition of weakly admissible bundle over a cobordism from Definition \ref{admiss-cob}, as well as Definition \ref{homology-or} of homology orientations of a cobordism $W$.

In the following, we choose the 4-manifold perturbation $\pi$ from a contractible open set $U \subset \mathcal P^{(4)}_c$ of perturbations which restrict to fixed elements of $\mathcal P_{E_i, \delta}$ on the ends; this open set was defined in the process of proving Theorem \ref{trans3} (when comparing, take the constant $C$ large enough that all instantons $\mathbf{A}$ in a trajectory $z$ with $\text{gr}_z(\alpha, \beta) \leq 10$ have energy less than $C$). It contains the perturbation $(\pi_1, 0, \pi_2)$ supported on the ends, with no interior part. 

If $\pi \in U$, we say that $\pi$ is \emph{$W$-small}.

\begin{theorem}\label{CobordismPackage}Suppose $(W,\mathbf{E})$ is an oriented 4-manifold with two cylindrical ends, with incoming end modeled on $(Y_1,E_1)$ and outgoing end modeled on $(Y_2,E_2)$; suppose the $E_i$ are weakly admissible, and that we have chosen regular perturbations $\pi_i \in \mathcal P_{E_i, \delta}$. Furthermore suppose that $\mathbf{E}$ is a weakly admissible $U(2)$-bundle over $W$. Then we have the following.
\begin{enumerate}
\item For every pair of critical orbits $\alpha \subset \mathfrak C_{\pi_1}(Y_1)$ and $\beta \subset \mathfrak C_{\pi_2}(Y_2)$, there is a set of configurations of framed connections on $W$ from $\alpha$ to $\beta$, denoted $\widetilde{\mathcal B}_{\mathbf{E}}(\alpha, \beta)$. The set of components is written $\pi_0 \widetilde{\mathcal B}_{\mathbf{E}}(\alpha, \beta)$, and carries a free and transitive action of both $\pi_1 \widetilde{\mathcal B}_{E_1}(\alpha)$ and similarly with $E_2$; in particular, it carries an affine identification $$\pi_0 \widetilde{\mathcal B}_{\mathbf{E}}(\alpha, \beta) \cong \mathbb Z.$$ For each $z$ in this set, there is an integer $\textup{gr}^W_z(\alpha, \beta)$, satisfying the same additivity formula as in the previous theorem, which is independent of $z$ after reducing modulo 8. \\
\item Let $\pi$ be a $W$-small regular perturbation on $W$, restricting to regular perturbations $\pi_i$ on the cylindrical ends corresponding to $Y_i$. For each pair of critical orbits $\alpha, \beta$, with $z$ as above, we associate a smooth $SO(3)$-manifold $\widetilde{\mathcal M}_{\mathbf{E},z,\pi}(\alpha, \beta)$, which is of dimension $\text{\em{gr}}^W_z(\alpha, \beta)+\dim \alpha$ so long as $\widetilde{\mathcal M}_z(\alpha, \beta)$ contains an irreducible connection. The term $\textup{gr}^W_z(\alpha, \beta)$ is independent of $z$ modulo $8$. When $W$ is equipped with a smoothly embedded path $\gamma: \mathbb R \to W$ which agrees with $(-t,b_1) \in (-\infty, 0] \times Y_1$ and $(t, b_2) \in [0,\infty) \times Y_2$ for large enough $|t|$, then $\widetilde{\mathcal M}$ is imbued with smooth equivariant endpoint maps to $\alpha$ and $\beta$. A residual set of $\pi \in U$ are regular.\\
\item If $\textup{gr}^W_z(\alpha, \beta) \leq 10 - \dim \alpha$, then there is a natural compactification of this $SO(3)$-manifold to a compact topological $SO(3)$-manifold with corners and a smooth structure on each stratum, $\overline{\mathcal M}^W_{\mathbf{E},z,\pi}(\alpha, \beta)$. The endpoint maps extend to $\overline{\mathcal M}^W$.\\
\item If $W$ is given a homology orientation as in Definition \ref{homology-or}, and we choose generators of the 2-element sets $\Lambda(\alpha)$ and $\Lambda(\beta)$ from the previous proposition, these induce a fiber orientation on $e_-: \overline{\mathcal M}_{\mathbf{E},z,\pi}(\alpha, \beta)\to \alpha$. Negating any one of these choices negates the corresponding orientation.\\
\item Suppose $b_1(W) = b^+(W) = 0$. So long as all three of the perturbations $\pi_1, \pi, \pi_2$ are taken sufficiently small, the set of reducible orbits in $\overline{\mathcal M}_{\mathbf{E},\pi}(\alpha_1, \alpha_2)$ may be identified as the following set, written $\textup{Red}(W,\mathbf{E})$. Write $c_1(\mathbf{E}) = \lambda \in H^2(W;\mathbb Z)$, and write $\lambda_i = \lambda\big|_{Y_i}$. There is an induced map $\textup{Pair}(H^2 W) \to \textup{Pair}(H^2 Y_i)$ given by restriction; the orbit $\alpha_i$ corresponds to some $$r_i \in \mathcal R_{\lambda_i}(Y_i, E_i) \subset \textup{Pair}(H^2 Y_i).$$ Then $\textup{Red}(W, \mathbf{E})$ is the subset of $\textup{Pair}(H^2 W)$ consisting of pairs $\{z, w\}$ restricting to $r_i$ on the corresponding ends, and with $z+w = \lambda$. The set of full reducibles is taken to the subset with $z = w$. If $b^+(W) \neq 0$ and $\pi$ is a regular perturbation, then no reducibles arise in the moduli space $\overline{\mathcal M}$.\\
\item Assuming we have chosen a homology orientation of $W$ and an element of all relevant orientation sets, there is a fiber oriented decomposition \begin{align*}(-1)^{\dim \alpha}\partial \overline{\mathcal M}^W_{\mathbf{E},z,\pi}(\alpha, \beta) \cong &\bigcup_{\substack{\gamma \subset \mathfrak C_{\pi_2}, w_1, w_2 \\ w_1 \ast w_2 = z}} (-1)^{\textup{gr}^W_{w_1}(\alpha, \gamma)+1}\;\; \overline{\mathcal M}^W_{\mathbf{E},w_1,\pi}(\alpha, \gamma) \times_\gamma \overline{\mathcal M}_{E_2,w_2,\pi_2}(\gamma, \beta)\\ &\bigcup_{\substack{\gamma \subset \mathfrak C_{\pi_1}, w_1, w_2 \\ w_1 \ast w_2 = z}} \overline{\mathcal M}_{E_2,w_1,\pi_1}(\alpha, \gamma) \times_\gamma \overline{\mathcal M}^W_{\mathbf{E},w_2,\pi}(\gamma, \beta).\end{align*}
\end{enumerate}
\end{theorem}

\begin{proof}That such a regular perturbation exists, and that they are generic in $U$, is Theorem \ref{trans3}. The grading function was defined in Definition \ref{grading} and seen to be well-defined mod 8 in Corollary \ref{gr-mod-8}, the compactification with no bubbling is given in Corollary \ref{nobubbling}, and the manifold-with corners structure (and decomposition of the boundary) given by Proposition \ref{gluing}. The orientations are provided by Proposition \ref{orientations}. The identification of reducible orbits was given by Proposition \ref{action3} and Proposition \ref{red-solutions-nobplus}.

The path $\gamma$ is necessary to choose a basepoint in $W$ (as $\gamma(1/2)$, say) and then to parallel transport the chosen framing to the basepoints $\gamma(0), \gamma(1)$ on the boundary components.

The fiber orientations and oriented decomposition were provided by Proposition \ref{orientations}.
\end{proof}

Lastly, to prove invariance of these induced maps, we will want a version of this package for \emph{families} of perturbations $\pi$. Recall both the definition of \emph{broken metric} and \emph{regular family of metrics and perturbations} from Chapter \ref{famsec}. We only use this notation for the case of families indexed by the interval $[0,1]$. 

Also recall that Definition \ref{admiss-cob} of weakly admissible cobordisms partitions these into three classes: first, the admissible cobordisms for which one end is an admissible bundle (so supports no reducible flat connections). Second, there are the weakly admissible cobordisms with rational homology sphere ends, $b_1(W) = 0$, and $b^+(W) > 0$ and $\mathbf{E}$ non-trivial. Third, there are those with $b_1(W) = b^+(W) = 0$ and which satisfy the `$\rho$-monotonicity inequality' 
$$\rho_{\pi_1}(r_1) \leq \rho_{\pi_2}(r_2)$$ and support no bad reducibles, meaning that $H_1 Y_1 \oplus H_1 Y_2 \to H_1 W$ is surjective.

(The second condition Definition \ref{admiss-cob} is never satisfied for $U(2)$-bundles.)

\begin{theorem}\label{FamilyPackage}Suppose $(W,\mathbf{E}): (Y_1, E_1, \pi_-) \to (Y_2, E_2, \pi_+)$ is a weakly admissible cobordism between 3-manifolds equipped with regular perturbations $\pi_\pm \in \mathcal P_{E_i, \delta}$.

Suppose $I = [0,1]$ parameterizes a family of metrics on $W$ (with fixed cylindrical ends); the metrics $g_0$ and $g_1$ may be cut along a separating 3-dimensional submanifold of $W$, as in Chapter \ref{sec:4d-families}, itself equipped with a regular perturbation. Then there is a contractible open set $U_g$ with
$$(\pi_-, 0, \pi_+) \subset U_g \subset U \subset \mathcal P^{(4)}_c$$ such that the following are true. Suppose $(W,\mathbf{E})$ is equipped with a pair of regular perturbations $\pi(0), \pi(1) \in U_g \subset \mathcal P^{(4)}_c$ which equal the $\pi_\pm$ on the ends.

\begin{enumerate}
\item If $b^+(W) \neq 1$ or $\mathbf{E}$ supports no reducibles, there is a path $\pi: [0,1] \to U_g$ so that $\pi(t)$ and $g_t$ form a regular family of metrics and perturbations.
\item If $b^+(W) = 1$, then for $\pi_0, \pi_1$ as above and any path $\pi: [0,1] \to U_g$ that suppots no reducible instantons, we may perturb $\pi(t)$ on the interior to form a regular family of metrics and perturbations.\\
\item If $\pi(t)$ is a regular family of perturbations on $[0,1] = I$, there are smooth $SO(3)$-manifolds $\widetilde{\mathcal M}^W_{\mathbf{E},z,\pi(t)}(\alpha, \beta)$ of dimension $\text{\em{gr}}_z(\alpha, \beta) + \dim \alpha + 1$. \\
\item A homology orientation of $W$ and a choice of element of each of $\Lambda(\alpha)$ and $\Lambda(\beta)$ induces a fiber orientation on $$e_-: \widetilde{\mathcal M}^W_{\mathbf{E},z,\pi(t)}(\alpha, \beta) \to \alpha,$$ as before. These orientations negate under orientation-reversal of any one of these choices.\\
\item As long as $\text{\em{gr}}^W_z(\alpha, \beta) \leq 9 - \dim \alpha$, the moduli space $\widetilde{\mathcal M}_{\mathbf{E},z,\pi(t)}(\alpha, \beta)$ has a natural compactification $\overline{\mathcal M}_{\mathbf{E},z,\pi(t)}(\alpha, \beta)$ satisfying the same properties as before.\\
\item If an element is chosen from each relevant orientation set as above, then there is a fiber oriented decomposition \begin{align*}\partial \overline{\mathcal M}_{\mathbf E,z,\pi(t)}(\alpha, \beta) = \bigcup_{\gamma \in \mathfrak C_{\pi_-}; z_1 \ast z_2 = z} &(-1)^{\textup{gr}_{z_1}(\alpha, \gamma)}\overline{\mathcal M}_{E_1,z_1,\pi_-}(\alpha, \gamma) \times_\gamma \overline{\mathcal M}^W_{\mathbf{E},z_2,\pi(t)}(\gamma, \beta)\\
\bigcup_{\zeta \in \mathfrak C_{\pi_+}; z_1 \ast z_2 = z} &(-1)^{\textup{gr}^{W}_{z_1}(\alpha, \zeta)} \overline{\mathcal M}^W_{\mathbf{E}, z_1,\pi(t)}(\alpha, \zeta) \times_\zeta \overline{\mathcal M}_{z_2}(\zeta, \gamma)\\
\cup \; \; \; \; &\overline{\mathcal M}^W_{\mathbf{E},z,\pi(1)}(\alpha, \beta) \cup -\overline{\mathcal M}^W_{\mathbf{E},z,\pi(0)}(\alpha, \beta).\\
\end{align*}

\end{enumerate}
\end{theorem}

\begin{proof}
This is precisely the content of Chapter \ref{famsec}; in particular, the transversality and gluing results in the absence of broken metrics are given by Proposition \ref{famtrans} and Proposition \ref{famgluing}, respectively; the extension to the case of broken metrics, including the third item above, is given by Proposition \ref{famcut}.
\end{proof}

We will use the fiber product maps of Lemma \ref{fibprod} associated to the $SO(3)$-equivariant endpoint maps $$e_-: \overline{\mathcal M}_{E,z,\pi}(\alpha, \beta) \to \alpha, \;\;\; e_+: \overline{\mathcal M}_{E,z,\pi}(\alpha, \beta) \to \beta$$ to define the framed instanton differential, as opposed to the usual counting of points in 0-dimensional moduli spaces. (Note that the degree of the fiber product map is precisely $\text{gr}_z(\alpha, \beta) - 1$ for moduli spaces on the cylinder, and precisely $\text{gr}_z^W(\alpha, \beta)$ on a cobordism.) Lemma \ref{fibprod} guarantees that when
$$\dim \overline{\mathcal M}(\alpha, \beta) \geq \dim \alpha + \dim \beta + 2$$
the map $\sigma \mapsto \sigma \times_{e_-} \overline{\mathcal M}(\alpha, \beta)$ is identically zero on the chain level. Because all orbits have dimension at most $3$, and
$$\dim \overline{\mathcal M}(\alpha, \beta) = \text{gr}_z(\alpha, \beta) + \dim \alpha+1,$$
we see that these fiber product maps vanish when $\text{gr}_z(\alpha, \beta) \geq 6$. As a result, we need not panic about the Uhlenbeck bubbling arising in large-dimensional moduli spaces. This same observation arose in the original definition of equivariant instanton homology given in \cite{AB}. In their case, the differential was defined via pullback and integration-over-the-fiber of differential forms, and whenever the differential form is in degree lower than the dimension of the fiber, the integral is zero.

\begin{definition}\label{relgrading}A \emph{relative $\mathbb Z/8$ grading} on a set $S$ is defined by a function 
$$r: S \times S \to \mathbb Z/8,$$ 
with $r(x,y) + r(y,z) = r(x,z)$. Such an additive function is equivalent to a function $i: S \to \mathbb Z/8$ considered up to translation; that is, $i \sim i'$ if $i(x) = i'(x) + c$, for some fixed $c$. Given a function $i$, one obtains an additive function from its differences: $$r_i(x,y) = i(y) - i(x).$$ This gives the same function for any equivalent $i, i'$. Conversely, evaluation defines functions $i_s(x) = r(x,s)$, with $i_s(x) = i_t(x) + i(t,s)$, so that $i_s \sim i_t$. 

A relatively graded chain complex splits as a direct sum $C = \oplus_{i \in \mathbb Z/8} C_i$, where the differential is the direct sum over $d_i: C_i \to C_{i-1}$. 
\end{definition}

The relevant relative grading to us is $\text{gr}: \mathfrak C_\pi \times \mathfrak C_\pi \to \mathbb Z/8$. Picking a critical orbit $\rho$ \emph{arbitrarily}, we write $i(\alpha) = \text{gr}(\rho, \alpha)$. If $\mathbf{E}$ is trivial, we set $i(\alpha) = i(\theta, \alpha)$, where $\theta$ is the trivial connection. This gives us an absolute grading when $\mathbf{E}$ is trivial (though changing the trivialization may change $i$ by a multiple of $4$). 

Just as there is a relative grading on critical orbits, the framed instanton chain complex $\widetilde{CI}$ has a relative grading. As a relatively graded $R$-module, $\widetilde{CI}$ is defined by

$$\widetilde{CI}_*(Y,E,\pi;R) := \bigoplus_{\alpha \subset \mathfrak C_\pi} C_{*}^{\textup{gm}}(\alpha;R)[i(\alpha)] \otimes_{R[\mathbb Z/2]} R[\Lambda(\alpha)].$$ Here $\Lambda(\alpha)$ is the 2-element orientation set discussed in Theorem \ref{ModuliPackage} (6). The group $\mathbb Z/2$ acts on $\Lambda(\alpha)$ by swapping the two elements, and on $C_*^{\text{\textup{gm}}}(\alpha;R)$ by negation.

If $\sigma_1: P \to \alpha$ and $\sigma_2: Q \to \beta$ are basic chains, the relative grading $|\sigma_1| - |\sigma_2| \in \mathbb Z/8$ is given as $$|\sigma_1| - |\sigma_2| =  (\dim P + i(\alpha)) - (\dim Q - i(\beta))\mod{8}.$$ Because $i(\alpha) - i(\beta) = \text{gr}(\rho, \alpha) - \text{gr}(\rho, \beta)$, this simplifes using the additivity formula of Theorem \ref{ModuliPackage} (2) to $$|\sigma_1| - |\sigma_2| = \dim P - \dim Q + \text{gr}(\alpha, \beta).$$ This further implies $(|\sigma_1| - |\sigma_2|) + (|\sigma_2| - |\sigma_3|) = |\sigma_1| - |\sigma_3|$, as expected.

The differential is given (termwise on basic chains $\sigma: P \to \alpha$) by

$$\partial_{CI} \sigma = \partial \sigma + \sum_{\substack{\beta \subset \mathfrak C_\pi\\z, \text{gr}_z(\alpha, \beta) \leq 5}} (-1)^{\dim \sigma} \sigma \times_{e_-} \overline{\mathcal M}_{E,z,\pi}(\alpha, \beta).$$ 
Here $\partial \sigma$ denotes the differential inside $C_*^{\textup{gm}}(\alpha)$ and the chain $\sigma \times_{e_-} \mathcal M$ is defined as in Lemma \ref{fibprod}; we orient the moduli spaces using choices of elements of $\Lambda(\alpha)$ and $\Lambda(\beta)$. Because swapping elements of these orientation sets negates the orientation on the fiber product, the map $C_*^{\text{gm}} \otimes R[\Lambda(\alpha)] \to C_*^{\text{gm}} \otimes R[\Lambda(\beta)]$ descends to the quotient under the two $\mathbb Z/2$ actions. We drop the orientation sets from notation as much as is reasonably possible.

The index demand on the sum ensures that all of the moduli spaces $\overline{\mathcal M}$ appearing in the sum are compact oriented topological manifolds with corners and a smooth structure on each stratum (as a consequence of Theorem \ref{ModuliPackage} (7)). The reason we can do this without concern is that the fiber product with any larger-dimensional moduli spaces is \textit{identically zero} --- we can consider our sum as being a formal truncation of what ``should be" the instanton differential $\partial_{CI}$, where we throw out moduli spaces of dimension too large to contribute. This happens precisely when the degree of the fiber product map is larger than $\dim \beta + 1$; because $SO(3)$-orbits have dimension at most 3, this is true when $\text{gr}_z(\alpha, \beta) -1 > 4$.

First observe that the differential decreases the relative grading by one: if $\sigma: P \to \alpha$ is a basic chain, taking the fiber product gives a basic chain $\sigma': P \times_{e_-} \overline{\mathcal M}_z \to \beta$. The relative grading between these is $$|\sigma| - |\sigma'| = \dim P - \dim \left(P \times_{e_-} \overline{\mathcal M}_z\right) + \text{gr}_z(\alpha, \beta).$$ Lemma \ref{fibprod} tells us that the dimension of the fiber product is $\dim P + \dim \overline{\mathcal M}_z - \dim \alpha$. Combining these with $\dim \overline{\mathcal M}_z = \text{gr}_z(\alpha, \beta) + \dim \alpha - 1$, we see that the relative grading $$|\sigma| - |\sigma'| = \dim P - (\dim P + \text{gr}_z(\alpha, \beta) - 1) + \text{gr}_z(\alpha, \beta) = 1.$$

\begin{lemma}\label{complex}$\widetilde{CI}_*(Y,E,\pi;R)$ is a chain complex. That is, $\partial_{CI}^2 = 0$. Furthermore, the right action of $C_*^{\textup{gm}}(SO(3);R)$ on $\widetilde{CI}$, acting on each $\oplus C_*^{\textup{gm}}(\alpha;R)$ on the right (induced by the right action of $SO(3)$ on $\alpha$) gives $\widetilde{CI}$ the structure of a dg-module.
\end{lemma}

\begin{proof}It is clear that, for a basic chain $\sigma: P \to \alpha$,
\begin{align*}\partial_{CI}^2 \sigma = \partial^2 \sigma &+ \sum_{\substack{\beta \in \mathfrak C_\pi\\z, \text{gr}_z(\alpha, \beta) \leq 5}} (-1)^{\dim \sigma - 1} (\partial \sigma) \times_{e_-} \overline{\mathcal M}_{E,z,\pi}(\alpha, \beta) \\
&+ \sum_{\substack{\beta \in \mathfrak C_\pi\\z, \text{gr}_z(\alpha, \beta) \leq 5}} (-1)^{\dim \sigma} \partial\left(\sigma \times_{e_-} \overline{\mathcal M}_{E,z,\pi}(\alpha, \beta)\right)\\
&+ \sum_{\substack{\gamma, \beta \in \mathfrak C_\pi\\z, \text{gr}_z(\alpha, \gamma) \leq 5\\ w, \text{gr}_w(\gamma,\beta) \leq 5}}(-1)^{2\dim \sigma + \text{gr}(\alpha, \gamma)-1}\sigma \times_{e_-} \left(\overline{\mathcal M}_{E,z,\pi}(\alpha, \gamma) \times_{\gamma} \overline{\mathcal M}_{E,z,\pi}(\gamma, \beta)\right).\end{align*} 

First, $\partial^2 \sigma = 0$ because $C_*^{\textup{gm}}(\alpha;R)$ is a chain complex. Using the decomposition $$\partial(\sigma \times_{e_-} \overline{\mathcal M})= \partial \sigma \times_{e_-} \overline{\mathcal M} + (-1)^{\dim \sigma + \dim \alpha} \sigma \times_{e_-} \overline{\mathcal M}$$ of Lemma \ref{fibprod}, and cancelling the $(\partial \sigma) \times_{e_-} \overline{\mathcal M}$ terms, this reduces to

\begin{align*}\partial_{CI}^2 \sigma &= \sum_{\substack{\beta \in \mathfrak C_\pi\\z, \text{gr}_z(\alpha, \beta) \leq 5}} (-1)^{\dim \alpha}\sigma \times_{e_-} \left(\partial \overline{\mathcal M}_{E,z,\pi}(\alpha, \beta)\right) \\ 
&+ \sum_{\substack{\gamma, \beta \subset \mathfrak C_\pi\\z, \text{gr}_z(\alpha, \gamma) \leq 5\\ w, \text{gr}_w(\gamma,\beta) \leq 5}} (-1)^{\text{gr}(\alpha, \gamma) - 1} \sigma \times_{e_-} \left(\overline{\mathcal M}_{E,z,\pi}(\alpha, \gamma) \times_{\gamma} \overline{\mathcal M}_{E,z,\pi}(\gamma, \beta)\right).\end{align*} 

The terms in the second sum can only be nonzero when $\text{gr}_z(\alpha, \gamma) + \text{gr}_w(\gamma, \beta) \leq 5$, which is to say that $\text{gr}_{z \ast w}(\alpha, \beta) \leq 5$. After eliminating terms in the sum which vanish for dimension reasons, we're left with 

\begin{align*}(-1)^{\dim \alpha}\partial_{CI}^2 \sigma &= \sum_{\substack{\beta \in \mathfrak C_\pi\\z, \text{gr}_z(\alpha, \beta) \leq 5}}\sigma \times_{e_-} \left(\partial \overline{\mathcal M}_{E,z,\pi}(\alpha, \beta)\right) \\ 
&+ \sum_{\substack{\gamma, \beta \in \mathfrak C_\pi\\z, w, \text{gr}_{z \ast w}(\alpha,\beta) \leq 5}} (-1)^{\dim \alpha + \text{gr}(\alpha, \gamma) - 1} \sigma \times_{e_-} \left(\overline{\mathcal M}_{E,z,\pi}(\alpha, \gamma) \times_{\gamma} \overline{\mathcal M}_{E,z,\pi}(\gamma, \beta)\right).\end{align*}

This is zero by the decomposition of the boundary given in Theorem \ref{ModuliPackage} (9).

That the action of $C_*^{\textup{gm}}(SO(3);R)$ makes $\widetilde{CI}$ into a $C_*^{\textup{gm}}(SO(3);R)$-module is clear from the fact that each summand $C_*^{\textup{gm}}(\alpha;R)$ is, and that the fiber product map of Lemma \ref{fibprod} is a $C_*^{\textup{gm}}(SO(3);R)$-module homomorphism.
\end{proof}

Suppose we have a weakly admissible bundle $\mathbf{E}$ over a cobordism $W$ from $(Y_1, E_1)$ to $(Y_2, E_2)$. To fix the orientations of these moduli spaces, we need to choose a homology orientation on $W$; we will suppress this from notation and just refer to ``a cobordism $W$".

\begin{lemma}\label{inducedmap}Suppose $$(W,\mathbf{E},\pi): (Y_1, E_1,\pi_1) \to (Y_2, E_2,\pi_2)$$ is a cobordism where $\mathbf{E}$ is a weakly admissible bundle and $\pi$ is a choice of metric and regular perturbation restricting to the $\pi_i$ on the ends. Furthermore suppose $W$ is equipped with an embedded path $\gamma$ between the basepoints $b_1$ and $b_2$. Then there is an induced $C_*^{\textup{gm}}(SO(3);R)$-equivariant chain map $$\widetilde F_{W,\mathbf{E},\pi,\gamma}: \widetilde{CI}(Y_1, E_1, \pi_1;R) \to \widetilde{CI}(Y_2, E_2, \pi_2; R).$$ 
\end{lemma}

\begin{proof}
The map is defined analogously to the differential itself; we only need to define its value on a basic chain $\sigma: P \to \alpha$, where $\alpha \subset \mathfrak C_{\pi_1}$ is a critical orbit on $Y_1$. Here the value is $$\widetilde F_{W,\mathbf{E},\pi,\gamma}\left(\sigma\right) = \sum_{\substack{\beta \subset \mathfrak C_{\pi_2}\\ z, \text{ gr}^W_z(\alpha, \beta) \leq 4}} (-1)^{\dim \sigma} \sigma \times_{e_-} \overline{\mathcal M}_{\mathbf{E}, z, \pi},$$ where the map to $\beta$ is defined using the positive endpoint map $e_+: \overline{\mathcal M}_{\mathbf{E}, z, \pi} \to \beta$. Defining this endpoint map is the essential place the path $\gamma: \mathbb R \to W$ is used. Here the sum is up to $\text{gr}_z(\alpha, \beta) = 4$ instead of $5$, because for cobordisms $\dim \overline{\mathcal M} = \text{gr} + \dim \alpha$, whereas there is an extra factor of $+1$ on the cylinder, account for the translation action.

That this is an equivariant chain map follows as in the previous lemma from the decomposition of the moduli space $\overline{\mathcal M}_{\mathbf{E}, z, \pi}$ of Theorem \ref{CobordismPackage} (6).
\end{proof}

These maps are essentially independent of the perturbation chosen on the cobordism, except in the case that $W$ is weakly admissible and $b^+(W) = 1$. (In the statement of the following, recall the definition of weakly admissible from Definition \ref{admiss-cob}.)

\begin{lemma}\label{invariance}Given two $W$-small regular perturbations $\pi_0, \pi_1$ on a weakly admissible cobordism $(W, \mathbf{E})$, there is a regular family of $W$-small perturbations $\pi(t)$ parameterized by $[0,1]$ interpolating between them, restricting to the same fixed perturbations on the ends for all $t$, or possibly some chain of them, so long as either $b^+(W) \neq 1$ or $\mathbf{E}$ supports no reducible connections. This induces a $C_*^{\textup{gm}}(SO(3);R)$-equivariant chain homotopy between
$$\widetilde F_{W,\mathbf{E}, \pi_0, \gamma} \simeq \widetilde F_{W,\mathbf{E}, \pi_1, \gamma}.$$
If $(W,\mathbf{E})$ is weakly admissible and $b^+(W) = 1$, then consider the set of metrics and regular perturbations on $W$, equipped with the equivalence relation $(g_0, \pi_0) \sim (g_1, \pi_1)$ if there is a path $(g_t\pi_t)$ between them so that no $(g_t,\pi_t)$ suppots a reducible $\pi_t$-ASD connection. If $\pi_0 \sim \pi_1$, then we have a chain homotopy as above.
\end{lemma}

\begin{proof}When $W$ supports no reducible connections, or has $b^+(W) \neq 1$, then the existence of such a regular family of perturbations is Theorem \ref{FamilyPackage} (1), and when $b^+(W) = 1$ it is point (2). Some care should be taken when varying the metric; first we should modify the perturbations $\pi_i$ to be small relative to the path of metrics $g_t$, and then we should construct a regular path between them.

That the moduli spaces of relevant dimension have compactifications to manifolds with corners and a smooth structure on each stratum follows from Theorem \ref{FamilyPackage} (5). The chain map is given on a basic chain $\sigma: P \to \alpha$ by sending $$\sigma \mapsto \sum_{\substack{\beta \in \mathfrak C_{Y_2}\\ z, \text{ gr}^W_z(\alpha, \beta) + 1 \leq 4}} \sigma \times_{e_-} \overline{\mathcal M}_{\mathbf{E}, z,\pi(t)}.$$ That this is a chain homotopy from $\widetilde F_{W,\mathbf{E}, \pi_0, \gamma}$ to $\widetilde F_{W,\mathbf{E}, \pi_1, \gamma}$ follows from the decomposition of the moduli space given by Theorem \ref{FamilyPackage} (6).

The main interesting point to make here is that these moduli spaces may be nonempty even though $\text{gr}_z(\alpha, \beta)$ is negative. This corresponds to the existence of $\pi(t)$-instantons for $t \in (0,1)$ that are regularly cut out in the \textit{family} --- but because of the index they cannot possibly be cut out regularly (considered as an instanton for the single perturbation $\pi(t)$).
\end{proof}

We need to relate maps arising from the geometric composition of cobordisms and the composition of chain maps.

\begin{lemma}\label{composition}Suppose we have two weakly admissible cobordisms $(W_i, \mathbf{E}_i, \pi_i)$ from $Y_i$ to $Y_{i+1}$ with regular perturbations $\pi_i \in \mathcal P^{(4)}_{W_i, L_i, \delta}$. Further suppose each $W_i$ is equipped with a path $\gamma_i$, so that the positive end of $W_1$ agrees with the negative end of $W_2$, and the paths concatenate to form a smooth path $\gamma$ in the composite cobordism. Denote the composite $(W, \mathbf{E})$. Then there is a regular family of metrics and perturbations $\pi(t)$ so that $\pi(0)$ is a regular perturbation and metric on $W$, while $\pi(1)$ is the broken metric (and perturbation) corresponding to the obvious way to glue the two cobordisms. This induces a $C_*^{\textup{gm}}(SO(3);R)$-equivariant chain homotopy $$\widetilde F_{W_2, \mathbf{E}_2, \pi_2, \gamma_2} \circ \widetilde F_{W_1, \mathbf{E}_1, \pi_1, \gamma_1} \simeq \widetilde F_{W, \mathbf{E}, \pi(0), \gamma}.$$
\end{lemma}
\begin{proof}Recall that the moduli spaces for the broken metric are by definition given by $$\overline{\mathcal M}_{\mathbf{E},z,\pi(1)}(\alpha, \gamma) = \bigcup_{\substack{\beta \in \mathfrak C_{Y_2}\\ z_1 * z_2 = z}} \overline{\mathcal M}_{\mathbf{E}_1,z_1,\pi_1} \times_{\beta} \overline{\mathcal M}_{\mathbf{E}_2, z_2, \pi_2}.$$ Thus the chain map given by a broken metric is precisely the composite of its component cobordisms.

Lemma \ref{composable-cob} guarantees that the composite of weakly admissible cobordisms remains weakly admissible. The conditions in Theorem \ref{FamilyPackage} (2)-(3) were that $\pi(t)$ admits $SO(2)$-reducible instantons if and only if $b^+(W) = 0$, and all fully reducible instantons are cut out transversely. These are assumed to be true for the two pieces of the broken metric. If $\pi(t)$ is the path from a broken metric to a non-broken metric, the same will hold for sufficiently small $t$. Some care should be taken to ensure that the perturbation $\pi(0)$ is $W$-small and has no part on the neck corresponding to the perturbation that used to lie `at infinity', as in the discussion preceding Proposition \ref{famcut}. Then by a small perturbation, we may ensure that $\pi(t)$ is a regular family of perturbations.

Given that, we may define the induced map as 
$$\sigma \mapsto \sum_\beta (-1)^{\dim \sigma} \sigma \times_{e_-}\overline{\mathcal M}_{\mathbf{E},z,\pi(t)}(\alpha, \beta).$$ 
There are no additional difficulties in verifying this is a chain homotopy between the map induced by the perturbation $\pi(0)$ and the map induced by the broken metric/perturbation $\pi(1)$.
\end{proof}

The chain maps $\widetilde F$ are invariant under diffeomorphisms of the cobordisms.

\begin{lemma}\label{DiffInvar}Suppose $W_1$ and $W_2$ are cobordisms equipped with weakly admissible bundles $\mathbf{E}_i$, regular perturbations $\pi_i$, and paths $\gamma_i$ between the basepoints of the ends. Furthermore suppose that the ends of the two cobordisms $(W_i, \mathbf{E}_i, \pi_i,\gamma_i)$ agree (so that we think of them as cobordisms between the same manifolds). Suppose $\varphi: W_1 \to W_2$ is a diffeomorphism, equal to the identity on the ends and with $\varphi(\gamma_1) = \gamma_2$, and $\Psi: \varphi^*\mathbf{E}_2 \cong \mathbf{E}_1$ an isomorphism of $U(2)$-bundles so that $\Psi$ takes $\varphi^*\pi_2$ to $\pi_1$. Furthermore suppose $\varphi$ preserves the homology orientations on the $W_i$. Then the chain maps $\widetilde F_{W_i, \mathbf{E}_i, \pi_i, \gamma_i}$ are identical.
\end{lemma}

\begin{proof}This data induces a diffeomorphism $$\Psi': \overline{\mathcal M}_{\mathbf{E}_1,\pi_1}(\alpha, \beta) \to \overline{\mathcal M}_{\mathbf{E}_2,\pi_2}(\alpha,\beta).$$ That this diffeomorphism preserves the endpoint maps follows because $\varphi(\gamma_1) = \gamma_2$; it preserves orientation because $\varphi$ preserves the homology orientation of the cobordisms. Given any basic chain $\sigma: P \to \alpha,$ the map $$P \times_{e_-} \overline{\mathcal M}_{\mathbf{E}_1} \to P \times_{e_-} \overline{\mathcal M}_{\mathbf{E}_2}$$ induced on the fiber product is a diffeomorphism preserving $e_+: P \times_{e_-} \overline{\mathcal M}_{\mathbf{E}_i} \to \beta$. By definition, this means these two basic chains are isomorphic, hence equal in $C_*^{\textup{gm}}(\beta;R)$.
\end{proof}

\begin{remark}In dimension 4, any two embedded paths which agree near the ends and are homotopic relative to their boundary are in fact isotopic relative to their boundary. The isotopy extension theorem then provides a diffeomorphism of $W$ fixing the boundary and taking $\varphi(\gamma_2) = \gamma_1$. Hence from the previous lemma $\widetilde F_{W, \mathbf{E}, \pi, \gamma_1} = \widetilde F_{W, \mathbf{E}, \varphi^*\pi, \gamma_2}$. So homotopic paths induce the same map in homology.
\end{remark}

Combining all of these, we see that we have a functor from a sort of cobordism category to a homotopy category of chain complexes. We make this precise in the following definition. Recall Definition \ref{sigdata} of the finite set $\sigma(Y,E)$ of \textit{signature data} on $(Y, E)$.

There are two categories of cobordisms relevant to us. The first includes the data of the perturbations; the second removes it as much as possible.

\begin{definition}\label{cobcat}
Let $\mathsf{Cob}_{3,b}^{U(2),w,\pi}$ denote the following category.\footnote{The notation $b$ is meant to stand for based, and $w$ for weakly admissible; the $\pi$ signifies that the perturbation data are included in an essential way.} The objects are closed oriented Riemannian 3-manifolds $Y$, equipped with a choice of basepoint $b \in Y$, a weakly admissible $U(2)$-bundle $E \to Y$, and a regular perturbation $\pi \in \mathcal P_{E, \delta}$ for some $\delta > 0$. 

The morphisms in $\mathsf{Cob}_{3,b}^{U(2),w,\pi}$ will be formal concatenations of morphisms of the following type (called `unbroken cobordisms').

An unbroken cobordism $(Y_1, E_1, \pi, b_1) \to (Y_2, E_2, \pi_2, b_2)$ is given by the data of $(W, \widetilde{\mathbf{E}}, \pi,\gamma, \varphi,\psi)$. Here $W$ is an oriented, homology oriented, Riemannian 4-manifold with two cylindrical ends, with chosen oriented isometries of these ends to $(-\infty, 0] \times Y_1$ and $[0,\infty) \times Y_2$, respectively, while $\mathbf{E}$ is a weakly admissible $U(2)$-bundle with specified isomorphisms to the pullback of $E_i$ on the cylindrical ends. The data $\pi$ is a $W$-small regular perturbation restricting to the perturbations $\pi_i$ on the ends, and $\gamma: \mathbb R \to W$ is an embedded path which agrees with $(t, b_i)$ on the ends, following the specified isometries above.

A morphism in $\mathsf{Cob}_{3,b}^{U(2),w,\pi}$ is simply a finite sequence of morphisms $W_i$ (for $1 \leq i \leq n$) so that the target of $W_i$ is the source of $W_{i+1}$ when $1 \leq i < n$. These morphisms are thought of as Riemannian manifolds with `broken metric and perturbation'.

Let $\mathsf{Cob}_{3,b}^{U(2),w}$ denote the following category. Its objects are closed oriented 3-manifolds $Y$ equipped with a basepoint $b \in Y$, a weakly admissible $U(2)$-bundle $E \to Y$, and a \emph{signature datum} $\sigma \in \sigma(Y,E)$. A morphism $(Y_1, \tilde E_1, \sigma_1, b_1) \to (Y_2, \tilde E_2, \sigma_2, b_2)$ is given by the data of $(W, \widetilde{\mathbf{E}}, [g,\pi],\gamma, \varphi,\psi)$. Here $W$ is a compact oriented, homology oriented, 4-manifold; $\mathbf{E}$ is a weakly admissible $U(2)$-bundle; $[g,\pi]$ is an equivalence class as in Lemma \ref{invariance} of metric $g$ and $(W,g)$-small regular perturbation $\pi$\footnote{By definition, unless $b^+(W) = 1$, all perturbations are equivalent, and so this additional data is vacuous!} restricting to perturbations $\pi_i$ on the ends with associated signature data $\sigma_{\pi_i} = \sigma_i$; and $\gamma: [0,1] \to W$ is an embedded path which is cylindrical in a collar of the boundary; $\varphi: \partial W \to Y_1 \sqcup \overline{Y_2}$ is a diffeomorphism sending $\gamma(i)$ to $b_{i+1}$; $\psi$ is an isomorphism $E_1 \sqcup E_2 \to \varphi^*\mathbf{E}$.

Two such cobordisms $(W_i, \widetilde{\mathbf{E}}_i, [\pi_i], \gamma_i, \varphi_i)$ are to be considered the same morphism in $\mathsf{Cob}_{3,b}^{U(2),w}$ if there is a diffeomorphism $\phi: W_1 \to W_2$ with $\varphi_2 \phi = \varphi_1$, a bundle isomorphism $\Psi: \widetilde{\mathbf{E}}_1 \to \phi^*\widetilde{\mathbf{E}}_2$ taking $\Psi^*[\pi_2] = [\pi_1]$ and a homotopy relative to the endpoints between $\gamma_1$ and $\phi^{-1}(\gamma_2)$.
\end{definition}
\noindent
The target category is the following.

\begin{definition}Let $\mathsf{Kom}_{C_*SO(3);R}^{r,\mathbb Z/8}$ denote the category with objects dg-modules over $C_*^{\textup{gm}}(SO(3);R)$ equipped with a relative $\mathbb Z/8$ grading, and whose morphisms are relatively graded $C_*^{\textup{gm}}(SO(3);R)$-equivariant chain maps. There is an equivalence relation on the morphisms in this category --- $C_*^{\textup{gm}}(SO(3);R)$-equivariant chain homotopy --- and there is a category $$\mathsf{Ho}\left(\mathsf{Kom}_{C_*SO(3);R}^{r,\mathbb Z/8}\right),$$ the homotopy category of right $C_*^{\textup{gm}}(SO(3);R)$-modules, whose morphisms are equivalence classes of equivariant chain maps. If we denote the category of relatively $\mathbb Z/8$-graded $R$-modules with a graded action of $H_*(SO(3);R)$ as $\mathsf{Mod}^{r,\mathbb Z/8}_{H_*SO(3);R}$, then taking homology gives a functor $$\mathrm{Ho}\left(\mathsf{Kom}_{C_*SO(3);R}^{r,\mathbb Z/8}\right) \to \mathsf{Mod}^{r,\mathbb Z/8}_{H_*SO(3);R}.$$
\end{definition}

The following is essentially immediate from definitions, and we record it as a lemma. 

\begin{lemma}\label{framed-functor-broken}Sending $(Y,E,b,\pi)$ to the framed instanton complex $\widetilde{CI}(Y, E, \pi; R)$ defines a functor $$\mathsf{Cob}_{3,b}^{U(2),w,\pi} \to \mathsf{Kom}_{C_*SO(3);R}^{r,\mathbb Z/8}.$$
\end{lemma}
\begin{proof}Because morphisms $\mathsf{Cob}_{3,b}^{U(2),w,\pi}$ are given by formal concatenations of a given generating set of morphisms, this merely asserts that we have defined induced maps for that generating set. This was done in Lemma \ref{inducedmap}. 
\end{proof}

Finally, we see that framed instanton homology is functorial on the category of cobordisms \emph{without} perturbation data.

\begin{theorem}\label{framed-functor}There is a functor $$\widetilde I: \mathsf{Cob}_{3,b}^{U(2),w} \to \mathsf{Mod}_{H_*SO(3);R}^{r,\mathbb Z/8},$$ so that for any regular perturbation $\pi$ on $(Y,E)$, the group $\widetilde I(Y,E,\sigma;R)$ is canonically isomorphic to the homology groups $\widetilde I(Y, E, \pi;R)$ of the chain complexes defined above. Following this canonical isomorphism, the induced map of a cobordism $(W,\mathbf{E},[\pi])$ is given by the induced map on homology of the map of chain complexes defined above.
\end{theorem}

\begin{proof}The previous lemmas showed that $\widetilde{CI}$ is a complex of the appropriate type, and showed that cobordisms $W$ in this category --- when equipped with a perturbation $\pi$ --- induce chain maps, as long as $\pi$ is sufficiently small. An isomorphism between $(W,\widetilde{\mathbf{E}}, \gamma, \varphi, \psi)$ (here, taking the path $\gamma_1$ diffeomorphically to $\gamma_2$) takes the moduli spaces (and their endpoint maps) of the first to those of the second, and hence the chain maps are the same on the nose. In particular, given regular perturbations $\pi_1, \pi_2$ on $(Y,E)$ with signature data $\sigma$, any path $\pi_t$ between them gives a sufficiently small perturbation on $\mathbb R \times Y$ (possibly with time-varying metric), and hence gives an induced map $\widetilde F_{\mathbb R \times Y,\pE,\pi_t,b}$; by Lemma \ref{invariance}, any two such chain maps are homotopic. 

We will use this to pin down the group $I(Y, E, \sigma)$ as something independent of the choice of metric and perturbation. 

Write $\rho_{12}: \widetilde I(Y,E,\pi_1) \to \widetilde I(Y,E,\pi_2)$ for the unique isomorphism between the groups associated to the perturbations $\pi_i$ (abbreviated for the moment $\widetilde I(\pi_i)$) defined by choosing a perturbation along $\mathbb R \times Y$ as above. Set $$\widetilde I(Y,E,\sigma;R) := \left(\bigoplus_{\pi; \; \sigma_\pi = \sigma} \widetilde I(Y,E,\pi;R)\right)\big/\left((x \in \widetilde I(\pi_1) \sim \rho_{12} x \in \widetilde I(\pi_2)\right).$$ The map $\widetilde I(Y, E,\pi;R) \to \widetilde I(Y, E;R)$ induced by inclusion into the direct sum is an isomorphism for any regular perturbation $\pi$. 

A cobordism with cylindrical ends $(W,\mathbf{E},\pi)$ restring to $\pi_i$ on the appropriate ends induces a map on each summand $\widetilde F_{W,\widetilde{\mathbf{E}},\gamma}: \widetilde I(Y_1,E_1,\pi_1) \to \widetilde I(Y_2,E_2,\pi_2)$. The invariance result of Lemma \ref{invariance} implies that this map is the same, independent of the choice of $\pi$ extending the $\pi_i$ on the ends; furthermore, Lemma \ref{composition} implies that these compose as expected with the canonical isomorphisms $\rho_{12}$ above. Therefore, this defines an induced map $\widetilde F_{W,\mathbf{E},\gamma}: \widetilde I(Y_1, E_1, \sigma_1) \to \widetilde I(Y_2, E_2, \sigma_2)$ by attaching cylindrical ends on each boundary component and choosing some sufficiently small regular perturbation.

We saw that the chain maps were independent of the perturbation $\pi$ and homotopy class $\gamma$ up to chain homotopy in Lemma \ref{invariance} (unless $b^+(W) = 1$, in which case there is an equivalence relation on perturbations), and that they compose as a functor should up to chain homotopy in Lemma \ref{composition}.
\end{proof}

The choice of $U(2)$ bundle is necessary to pin down signs, both in the differential on the chain complex and in the cobordism maps. But at least at the level of Floer homology groups, the $U(2)$-bundle is mostly inessential, as we will see in the following.

\begin{lemma}\label{Cocycle}Let $J$ be a nonempty finite set and $A$ be an abelian group. Write $a\delta_{ij}$ for the function $J^2 \to A$ which equals $a$ if $i = j$ and is zero otherwise; then denote
$a \delta_i: J^2 \to A$ by $(a\delta_i)(x,y) = (a\delta_{xi}) - (a\delta_{yi})$.

Suppose $s: J^2 \to A$ be a set map satisfying $$s(x,y) + s(y,z) = s(x,z).$$ 

Then there is a choice $a_j \in A$ for each $j \in J$ so that $$s = \sum_j a_j \delta_j;$$ the choice of $a_j$ is unique up to a uniform change $a'_j = a_j + b$, where $b$ is independent of $j$. 
\end{lemma}
\begin{proof}Choose a basepoint $p \in J$, and set $a_j = s(p,j)$; then $s' = s - \sum_j a_j \delta_j$ still satisfies $s'(x,y) + s'(y,z) = s'(x,z)$, and now further satisfies $s'(p,j) = 0$ for all $j$. Therefore for any $i,j \in J$, we have $$s'(p,i) + s'(i,j) = s'(p,j),$$ and so $s'$ is identically zero, as desired.

The essential uniqueness of $a_j$ follows from knowing the kernel of the map $(a_j) \mapsto \sum_j a_j \delta_j$. If $\sum_j a_j \delta_j(i,k) = 0$ for all $i,k$, then we obtain from the definition of $\delta_j$ that $a_i - a_k = 0$ for all $i,k$, as desired.
\end{proof}

\begin{corollary}\label{no-u2}Let $(Y,E,\pi)$ be a 3-manifold equipped with a weakly admissible $SO(3)$-bundle and regular perturbation $\pi$. Suppose we are given any two $U(2)$-bundles $\tilde E_i$ with an isomorphism $\tilde E_i \times_{U(2)} SO(3) \cong E$; we may define the four flavors of instanton Floer homology groups $I^\bullet(Y,\tilde E_i,\pi)$. 

Then there is an isomorphism $f: I^\bullet(Y,\tilde E_1, \pi) \to I^\bullet(Y, \tilde E_2, \pi)$ of relatively graded $R$-modules with an action of either $H_*(SO(3);R)$ or $H^{-*}(BSO(3);R)$, as appropropriate to the flavor $\bullet$. While this isomorphism is not canonically defined, it is defined up to a sign; that is, $\{f, -f\}$ is a canonically defined pair of isomorphisms between these two a priori different Floer homology groups.
\end{corollary}

\begin{proof}Note that the underlying graded $R$-module (or graded $C_*(SO(3);R)$-module, as appropriate) is independent of the choice of $U(2)$-bundle; the only contribution the choice of $U(2)$-bundle makes is an orientation of the moduli spaces $\overline{\mathcal M}(\alpha, \beta)$, or rather the fiber above a point in $\alpha$, assuming we choose an element of two orientation sets $\Lambda_\alpha$ and $\Lambda_\beta$. If we change the $U(2)$-bundle, the orientation sets change. Choosing an isomorphism $\Lambda_\alpha \to \Lambda'_\alpha$ for each $\alpha$ gives an identification $\overline{\mathcal M}_{\alpha \beta} \cong \overline{\mathcal M}'_{\alpha \beta}$, which is either oriented or unoriented, depending on a sign $s(\alpha, \beta) \in \pm 1$ (note that this sign is the same for any connected component of the moduli space); because the orientations are compatible with gluing, we find that $s(\alpha, \beta) s(\beta, \gamma) = s(\alpha, \gamma)$. 

We may abstract somewhat. The situation we have (for any of the instanton homology complexes) is that there is a decomposition 
$$CI^\bullet(Y, \tilde E_i, \pi) \cong \bigoplus_{j \in J} C_j \times_{\pm 1} \Lambda^i_j$$ as graded $R$-modules; here $J = \mathfrak C_\pi$, the set of $\pi$-flat orbits, and $\Lambda^i_j$ is a 2-element orientation set for each $j$, which depends on the choice of bundle $\tilde E_i$. The differential may be written in matrix form as $(d^i_{jk})$; pinning down an actual map $C_j \to C_k$ requires choosing an element of $\Lambda^i_j$ and $\Lambda^i_k$.

Choose an arbitrary isomorphism $r_j: \Lambda^1_j \cong \Lambda^2_j$ for each $j$; this gives a map 
$$f_r: CI^\bullet(Y, \tilde E_1, \pi) \to CI^\bullet(Y, \tilde E_2, \pi).$$
This map satisfies 
$$f_r d^1_{jk} = s(j,k)f_r d^2_{jk},$$
with as above $s(j,k) s(k,\ell) = s(j,\ell).$ Following Lemma \ref{Cocycle}, there are exactly two functions $t: J \to \{\pm 1\}$ with $t(j)t(k) = s(j,k).$ If we change the chosen identifications $r_j: \Lambda^1_j \to \Lambda^2_j$ to 
$$(tr)_j := t(j)r_j: \Lambda^1_j \to \Lambda^2_j,$$ 
correspondingly we find that the new map $f_{tr}$ is in fact a chain map, and hence an isomorphism of chain complexes (with action of whatever dg-algebra is appropriate). Because there are only two choices of $r$ so that $f_r$ is a chain map (and the two choices are negatives of one another), we have a canonical pair of isomorphisms $\{f, -f\}$ between $I^\bullet(Y,\tilde E_1, \pi)$ and $I^\bullet(Y, \tilde E_2, \pi)$, where $\tilde E_i$ are two $U(2)$-lifts of a given $SO(3)$-bundle.
\end{proof}

\begin{remark}There is a natural partial ordering on $\sigma(Y,E)$: an element $\sigma$ defines a function $f_\sigma: \text{Red}(Y, E) \to \mathbb Z_{\geq 0}$, and we say that $\sigma \leq \sigma'$ if $f_\sigma(\alpha) \leq f_{\sigma'}(\alpha)$ for all $\alpha$. The problem with attempting to show that the instanton homology groups are invariant of the signature datum is that the cylinder is only a weakly admissible cobordism $(Y, E, \sigma) \to (Y, E, \sigma')$ if $\sigma \leq \sigma'$. Any attempt to resolve this needs somehow to cope with reducible connections on the cylinder which cannot be made to be cut out transversely; perhaps the obstruction bundle technique, as described briefly at the end of \cite{Don} and used to great effect in \cite{taubes1984self}, is one such tool.
\end{remark}

\begin{remark}\label{infty-stuff}Suppose we work instead in the category $\mathsf{Cob}_{3,b}^{SO(3),a}$ of admissible cobordisms. Then in fact, a much stronger invariance property of the maps $\widetilde F_{W,\mathbf{E}}$ is true: the homotopies between induced maps are homotopic, and so on. 

A precise formulation of this statement uses the language of \emph{quasicategories} (sometimes called $(\infty,1)$-\emph{categories}); a good introduction to the language is \cite{HLS}. There is a quasicategory $\pi\mathcal{C}ob_{3,b}^{U(2),a}$ of 3-manifolds equipped with perturbations, whose morphism spaces are simplicial sets whose $n$-simplices are $\Delta^n$-indexed admissible families of (possibly broken) perturbations on cobordisms. The framed instanton chain complex lifts to a functor from this quasicategory to the quasicategory of chain complexes, and furthermore the forgetful functor $\pi\mathcal{C}ob_{3,b}^{U(2),a} \to \mathsf{Cob}_{3,b}^{U(2),a}$ to the homotopy category is an equivalence of quasicategories. This is essentially the statement that the space of admissible perturbations on a given cobordism is contractible, and a version of this would be provided by a natural extension of Theorem \ref{FamilyPackage} for simplices; the only obstructions to such an extension lies at reducible connections, which the admissibility assumption is used to avoid. Using this, it's possible to find a chain-level version of $\widetilde{CI}$ that is functorial under cobordisms on the nose, as opposed to up to homotopy, and well-defined up to essentially unique natural equivalence. In fact, this is done using the quasicategorical Kan extension; the definition of $I(Y,E,\sigma)$ above is a very special case of a Kan extension. We don't see any need for this structure, and so leave the details to the interested reader.

Any extension of such a result to all weakly admissible $(W,\widetilde{\mathbf{E}})$ would require a substantially different notion of ``regular family", as there are obstructions to achieving transversality in the standard sense at reducible connections for families. 

In this framework, the space of perturbations on a fixed 3-manifold form a quasi-category (actually, a Kan complex): $0$-simplices are pairs $(Y,\tilde E)$ of a weakly admissible 3-manifold and regular perturbation $\pi$, and an $n$-simplex starting at $(Y,\tilde E,\pi_1)$ and ending at $(Y,\tilde E,\pi_2)$ is a regular family of perturbations on $\mathbb R \times Y$ parameterized by the $(n-1)$-dimensional simplex, so that for all perturbations in this family, the restriction of $\pi$ to the the corresponding end is the fixed perturbation $\pi_i$. Then \ref{FamilyPackage} (1) guarantees that the space of perturbations \emph{with fixed signature perturbation} is contractible, and that this quasi-category is equivalent to the poset whose elements are signature data on $Y$, and there is a morphism $\sigma \to \sigma'$ iff $\sigma \leq \sigma'$.
\end{remark}

\section{The index filtration}\label{sec:Floer-filt}
Recall that the relatively $\mathbb Z/8$-graded complex $$\widetilde{CI}(Y,E,\pi;R)$$ is defined as a graded $C_*^{\textup{gm}}(SO(3);R)$-module to be $$\bigoplus_{\alpha \in \mathfrak C_\pi} C_{*}^{\textup{gm}}(\alpha;R)[i(\alpha)],$$ where $i(\alpha) = \text{gr}(\rho, \alpha) \in \mathbb Z/8$; this grading is well-defined up to a translation (arising from choosing a different base orbit $\rho$), and thus defines a relative grading on $\widetilde{CI}$. Here we drop the orientation sets $\Lambda(\alpha)$ from notation; if so desired, there is no harm in making a choice of orientation for each critical orbit.

Now if $z$ is a homotopy class of path from $\rho$ to $\alpha$, write $$i(\alpha,z) = \text{gr}_z(\rho, \alpha) \in \mathbb Z/8.$$ There is a $\mathbb Z$-graded complex $\widetilde{CI}_\textup{unr}(Y,E,\pi;R)$ given as $$\bigoplus_{\substack{\alpha \subset \mathfrak C_\pi\\z \in \pi_1\widetilde{\mathcal B}_E, \rho, \alpha}} C_{*}^{\textup{gm}}(\alpha;R)[i(\alpha,z)].$$

If $(\alpha_i, z_i)$ are two such labelled critical orbits, there is a unique homotopy class of path $w$ from $\alpha_1$ to $\alpha_2$ so that $z_2 \simeq w \ast z_1$. Then using the additivity of $\text{gr}$, we see that the degree difference between $\sigma_i: P_i \to \alpha_i$ is given as \begin{align*}|\sigma_1| - |\sigma_2| &= \left(\dim P_1 - i(\alpha_1, z_1)\right) - \left(\dim P_2 - i(\alpha_2, z_2)\right)\\
&= \dim P_1 - \dim P_2 + \left(\text{gr}_{z_2}(\rho, \alpha_2) - \text{gr}_{z_1}(\rho, \alpha_1)\right)\\
&= \dim P_1 - \dim P_2 + \text{gr}_w(\alpha_1, \alpha_2).\end{align*}

We thus see that the grading induced by choosing $\rho$ as a basepoint induces the expected relative grading.

We write a basic chain $\sigma: P \to \alpha$ corresponding to a critical orbit $\alpha$ labelled by $z$ as $(\sigma, z)$. The differential of $\sigma: P \to (\alpha, z)$ is given by $$\partial_{CI}(\sigma, z) = (\partial \sigma, z) + \sum_{\substack{\beta \in \mathfrak C_\pi\\w \in \pi_1(\widetilde{\mathcal B}_E, \alpha, \beta)\\ \text{gr}_w(\alpha, \beta) \leq 5}}\!\!\! \left(\sigma \times_{e_-} \overline{\mathcal M}_{E,w,\pi}(\alpha, \beta),z \ast w\right).$$

Each critical orbit $\alpha$ has a unique homotopy class of path $1 \in \pi_1(\widetilde{\mathcal B}_E, \alpha)$ with $\text{gr}_1(\alpha, \alpha) = 8$; the element $1$ is a cyclic generator of the fundamental group, which is isomorphic to $\mathbb Z$. The complex $\widetilde{CI}_\textup{unr}(Y,E,\pi;R)$ has a periodicity isomorphism, sending the component $C_*^{\textup{gm}}(\alpha;R)$ labelled by $z$ identically to the component labelled by $z + 1$. That this commutes with the differential above is only the statement $w + 1 = 1 + w$ for any path $w$.

The $8\mathbb Z$-periodic, $\mathbb Z$-graded complex $\widetilde{CI}_\textup{unr}(Y,E,\pi;R)$ is called the \emph{unrolled complex} of $\widetilde{CI}(Y, E, \pi; R)$ in Chapter \ref{PeriodicMachine}; the quotient by the above periodicity isomorphism just forgets about the labelling by $z$.

$\widetilde{CI}_\textup{unr}(Y,E,\pi;R)$ carries an honest filtration by \text{index}. We write $$F_s\widetilde{CI}_\textup{unr}(Y,E,\pi;R) = \bigoplus_{\substack{\alpha \in \mathfrak C_\pi\\z \in \pi_1(\widetilde{\mathcal B}_E, \rho, \alpha)\\i(\alpha, z) \leq s}} C_{*}^{\textup{gm}}(\alpha;R)[i(\alpha,z)].$$ Because a nonempty moduli space $\overline{\mathcal M}_{E,w,\pi}(\alpha, \beta)$ can only exist if $$\text{gr}_z(\rho, \alpha) > \text{gr}_{z \ast w}(\rho, \beta),$$ the differential decreases index and hence preserves the filtration. Furthermore, because the filtration is defined by taking a direct sum of $C_*^{\textup{gm}}(\alpha;R)$ for $(\alpha, z)$ satisfying the index bound, this is a filtration by $C_*^{\textup{gm}}(SO(3);R)$-modules.

Because $\text{gr}_{z+1}(\rho, \alpha) = \text{gr}_z(\rho, \alpha) + 8$, this is a periodic filtration in the sense of Definition \ref{pfilt}. The associated graded module is $$\text{gr}_p \widetilde{CI}_\textup{unr}(Y, E, \pi; R) \cong \bigoplus_{\substack{\alpha \in \mathfrak C_\pi \\ i(\alpha) = p \bmod 8}}\!\! C^{\text{gm}}_*(\alpha; R)[p].$$

In particular, because $C_*^{\text{gm}}(X;R)$ is supported in degrees at most $\dim X + 1$, and all of our components $\alpha$ are $SO(3)$-orbits, each associated graded piece is bounded, supported in degrees $[0, 4]$.

Now note that the instanton differential $\partial_{CI}: \widetilde{CI}_\textup{unr} \to \widetilde{CI}_\textup{unr}$ decomposes as $\partial_{CI} = \partial_0 + \partial_1 + \cdots + \partial_5$, where $\partial_0(\sigma, z) = (\partial \sigma, z)$ and for $k>0$, $$\partial_{k}(\sigma, z) = \sum_{\substack{\beta \in \mathfrak C_\pi \\ w \in \pi_1(\widetilde{\mathcal B}_E, \alpha, \beta)\\ \text{gr}_w(\alpha, \beta) = k}} \left(\sigma \times_{e_-} \overline{\mathcal M}_{E,w,\pi}(\alpha, \beta),z*w\right).$$ This is the component of the differential that decreases filtration by $k$, and is given by those fiber product maps with moduli spaces which increase dimension by $k-1$. The decomposition of $\partial_{CI}$ into the $\partial_k$ ends at $\partial_5$ because fiber products with moduli spaces with $\text{gr}_w(\alpha, \beta) > 5$ are identically zero. 

Furthermore, observe that $\dim \overline{\mathcal M}_z(\alpha, \beta) = \dim \alpha + \text{gr}_z(\alpha, \beta) - 1$, while for any nonempty moduli space we have $\dim \overline{\mathcal M}_z(\alpha, \beta) \geq 3$ because by Theorem \ref{ModuliPackage} (8), the $SO(3)$-action is free. Therefore, if $\overline{\mathcal M}_z(\alpha, \beta)$ is nonempty and $\alpha$ is an $SO(2)$-reducible, necessarily $\text{gr}_z(\alpha, \beta) > 1$; if $\alpha$ is fully reducible, then necessarily $\text{gr}_z(\alpha, \beta) > 3$. In particular, $\partial_1$ is identically zero on reducible orbits; in addition both $\partial_2$ and $\partial_3$ vanish on fully reducible orbits.

Therefore, we have the following, which a special case of Theorem \ref{4flavors-periodic}, item (3).

\begin{theorem}\label{IndexSS}Let $(Y,E)$ be a weakly admissible bundle and $\pi$ a regular perturbation. There is a $(\mathbb Z/8,\mathbb Z)$-bigraded spectral sequence of $H_*(SO(3);R)$-modules so that the $\mathbb Z/8$-grading is relative, whose $E^1$ page is $$E^1_{p,q} = \bigoplus_{\substack{\alpha \subset \mathfrak C_{\pi} \\ \textup{gr}(\alpha) = p}} H_q(\alpha;R)$$ this spectral sequence converges to $E^\infty_{p,q} \cong \textup{gr}_p \widetilde{I}_{p+q}(Y,E,\pi;R)$, the pth component of the associated graded vector space of $\widetilde I(Y,E,\pi;R)_{p+q}$.

On the unrolled $(\mathbb Z, \mathbb Z)$-graded spectral sequence, we have $$E^\infty_{p,q} \cong \textup{gr}_p \widetilde I^\textup{unr}_{p+q}.$$

For a class $[x] \in E^r$ which has a representative $x \in \widetilde{CI}$ with $d_i x = 0$ for $i <r$, the spectral sequence differential $d_r[x]$ may be identified with $[\partial_r x]$. In particular, the differential $d_1$ on $E^1$ is only nonzero on the irreducible components, where it is given by counting points in 0-dimensional moduli spaces $\overline{\mathcal M}_{E,z,\pi}(\alpha, \beta)/SO(3)$ (note that $\beta$ need not be irreducible).

The index spectral sequence for $\widetilde{I}$ degenerates on the $E^5$ page for dimension reasons.
\end{theorem}

\begin{proof}The only conditions we need to check to apply Theorem \ref{4flavors-periodic}(3) is that each associated graded module $\text{gr}_p \widetilde{CI}_{\text{unr}}$ (for each fixed $p$) is bounded in degree. These follow by definition: $C_*^\text{gm}(X;R)$ is supported in degrees $[0, \dim X + 1]$, and $\text{gr}_p \widetilde{CI}_{\text{unr}}$ is a finite direct sum of $SO(3)$-orbits with $i(\alpha) = p$, so is supported in degrees $[p, p+4]$. This gives us the existence of the spectral sequence, as well as the fact that it computes the associated graded of $\widetilde I$. \\

Degeneration takes only slightly more work. That the differentials on elements $[x]$ with lifts that have $d_ix = 0$ for $i < r$ are induced by $\partial_r$ is an elementary diagram chase. In fact, that this is true holds more generally: Wall's notion of a \emph{multicomplex} (see Definition \ref{multicomplex}) is a complex so that the differential splits nicely into components $d_r$ that decrease the filtration by $r$, and $C$ is equipped with the induced differential $$d = d_0 + d_1 + d_2 + \cdots$$ Here, $d_r = \partial_r$ are the terms which decrease index by precisely $r$. Because $\partial_r$ increases the dimension of a geometric chain by $r-1$ and the homology groups are supported in degrees at most three, these summands $d_r$ vanish for $r \ge 5$.
\end{proof}

As an immediate corollary, we may define and calculate the Euler characteristic $\chi(\widetilde{I}(Y,E,\sigma))$, where $\sigma$ is a signature datum.

\begin{corollary}Let $(Y,E,\sigma)$ be a 3-manifold equipped with a weakly admissible bundle and signature datum $\sigma$. Then $\widetilde{I}(Y,E, \sigma)$ is finitely generated, and 
$$\chi(\widetilde{I}(Y,E, \sigma)) = |H^2(Y;\mathbb Z)| = |H_1(Y;\mathbb Z)|.$$ 
If $b_1(Y) > 0$, we interpret the right-hand-side to be zero.
\end{corollary}

This follows because $$\chi(\widetilde{I}(Y,E, \sigma)) = \chi(\text{gr} \;\widetilde{I}(Y,E, \sigma));$$ because $\chi(SO(3)) = 0$, the only contribution to this Euler characteristic is from the reducibles. Then the result follows from the enumeration of reducibles given by Proposition \ref{red-enumerate}, which have even index by the calculation of Proposition \ref{gr-red} (and so all contribute positively to the Euler characteristic sum). The isomorphism $H^2(Y;\mathbb Z) = H_1(Y;\mathbb Z)$ is Poincar\'e duality.

\begin{remark}In general it's rare for a multicomplex that a class $[x] \in E^r$ has a representative with all $d_i x = 0$ for $i < r$, but is relatively common in the $\widetilde I$ spectral sequence. The discussion simplifies substantially if $\frac 12 \in R$, so we make this assumption; then we know that $H_*(SO(3);R) = R \oplus R[3]$ and $$H_*(SO(3)/SO(2);R) = H_*(S^2;R) = R \oplus R[2].$$
\noindent
Choosing an arbitrary class $\rho \in \mathfrak C_\pi$, we write the absolute grading $$i(\alpha) = \text{gr}(\rho, \alpha).$$ We define

\begin{align*}C^{\text{irr}}_p &= \bigoplus_{\substack{\alpha \in \mathfrak C_\pi\\ \Gamma_\alpha = 1 \\ i(\alpha) =p}} R\\ C^{SO(2)}_p &= \bigoplus_{\substack{\alpha \in \mathfrak C_\pi\\ \Gamma_\alpha \cong SO(2) \\ i(\alpha) =p}} R \\ C^{\theta}_p &= \bigoplus_{\substack{\alpha \in \mathfrak C_\pi\\ \Gamma_\alpha = SO(3) \\ i(\alpha) =p}} R\end{align*}
\noindent
With this notation, the $E^1$ page is given as
\[ E^1_{p,q} = \begin{cases}
      C^{\text{irr}}_{p} \oplus C^{SO(2)}_p \oplus C^\theta_p & q = 0 \\
      C^{SO(2)}_p & q=2 \\
      C^{\text{irr}}_p & q=3
   \end{cases}
\]
The $E^1$ differential on $C^{\text{irr}}_p$ counts $SO(3)$-orbits in 3-dimensional moduli spaces. For $\alpha \in E^1_{p,0}$ the first page differential is given by counting points in unframed moduli spaces with $\text{gr}(\alpha, \beta) = 1$; $$d_1\alpha = \sum_{\beta, \text{gr}(\alpha, \beta) = 1} \left(\# \overline{\mathcal M}(\alpha, \beta)/SO(3)\right) \beta,$$ and $d_1 \alpha$ is zero if $\alpha$ is reducible.

Because the $E^1$ differential vanishes on the $SO(2)$-irreducibles, we see that the differential $d_2$ on $[\alpha] \in E^2_{p,2}$, for $\alpha \in C^{SO(2)}_p \subset E^1_{p,2}$, is given by a similar formula as above: counting, for each $\beta$ with $\text{gr}(\alpha, \beta) = 2$, the (signed) number of $SO(3)$-orbits of trajectories between them: $n(\alpha, \beta) = \# \overline{\mathcal M}(\alpha, \beta)/SO(3)$ for $\beta$ irreducible; so $$d_2[\alpha] = \sum_{\beta, \text{gr}(\alpha, \beta) = 2} n(\alpha, \beta) [\beta]_{p-2,3}.$$ $d_2$ vanishes everywhere else for dimension reasons. Finally, $d_4$ on $[\theta] \in E^4_{0,0}$ (when $Y$ is a rational homology sphere and so has trivial connections, and $0 \in V$ makes sense) counts $SO(3)$-orbits $\overline{\mathcal M}(\theta, \beta)/SO(3)$ between an irreducible and an irreducible with $\text{gr}(\theta, \beta) = 4$. (This is the map Donaldson calls $D_2$ in \cite[Section~7.1]{Don}.)
\end{remark}

Once we develop the equivalence to the Donaldson model $DCI$ in Chapter \ref{secDCI}, the terms in this spectral sequence become much easier to calculate, and visible at the chain level.

\section{Four flavors of instanton homology}\label{sec:fourflavors}
We will soon apply Theorem \ref{4flavors-periodic} to define three additional flavors of equivariant instanton homology, $I^+(Y), I^-(Y)$, and $I^\infty(Y)$. To explain what these different flavors represent, it is useful to analogize to the situation of finite $G$-CW complexes. To correctly state certain stabilization pheneomena, we must ensure everything in sight has a basepoint. Let $X_+$ denote the disjoint union of a $G$-space $X$ with a disjoint fixed basepoint. The smash product of pointed spaces is 
$$X \wedge Y = (X \times Y)\big/(X \times *_Y, *_X \times Y).$$ 
A pointed $G$-CW complex is a $G$-space assembled from cells of the form 
$$D^i \wedge (G/H)_+ = (D^i \times G/H)\big/(* \times G/H),$$ 
where the basepoint is in the boundary sphere and attaching maps preserve the basepoint. Given a pointed $G$-CW complex $X$, its most well-known homological invariants are its nonequivariant homology $H_*(X)$ and its Borel equivariant homology $$H_*^G(X) = H_*\big((X \wedge EG_+)/G\big),$$ both taken relative to a basepoint (equivalently, taking reduced homology). These are both invariants of $X$ up to equivariant homotopy equivalence. In fact, write $\Sigma X$ denotes the reduced suspension 
$$S^1 \wedge X = ([0,1] \times X)/([0,1] \times *, \{0\} \times X, \{1\} \times X)$$ 
of a pointed $G$-space. Because $H_{*+1}(\Sigma X) \cong H_*(X)$ and $\Sigma X \wedge_G EG_+ = \Sigma (X \wedge_G EG_+)$, the equivariant homology groups $H^G_*$ are invariants of $X$ up to \emph{stable} equivariant homotopy equivalence.

In the stable homotopy category of finite $G$-CW complexes, where the suspension operator $\Sigma$ has an inverse and there is a sphere $S^n$ for every $n \in \mathbb Z$, there is a contravariantly functorial Spanier-Whitehead duality operator $D_G$ sending $$D_G\big(S^n \wedge (G/H)_+\big) \cong S^{-n} \wedge (G/H)_+.$$ Explicitly, if $V \cong \mathbb R^N$ is a $G$-representation, and $X \hookrightarrow S^{V+1}$ is an equivariant embedding into the one-point compactification of $V \oplus \mathbb R$, then $D_GX$ is the desuspension $\Sigma^{-V}(S^{V+1} \setminus X)$.

In the nonequivariant case, the remarkable Alexander duality theorem identifies $H^{-*}(DX) \cong H_*(X)$, by $$H^{-*}(DX) = H^{-*}\left(\Sigma^{-n}(S^{n+1} \setminus X)\right) = H^{n-*}(S^{n+1} \setminus X) \cong H_*(X).$$
Though the double dual involved in writing $H^{-*}(DX)$ outputs a homology theory (a composition of two contravariant functors is a covariant functor), Alexander duality identifies it as the homology theory we started with, so this operation doesn't produce anything new. Note that we used negatively graded cohomology, $H^{-*}$, so that the gradings in this formula would work out.

However, in the equivariant case, something altogether new happens: applying the above formula to orbits, $$H^{-*}_G\Big(D_G\big((G/H)_+\big)\Big) = H^{-*}_G\big((G/H)_+\big) = H^{-*}(BH),$$ which is not $H_*^G\big((G/H)_+\big) = H_*(BH)$, nor some degree shift of it: $H^{-*}(BH)$ is concentrated in nonpositive degrees, and $H_*(BH)$ is concentrated in nonnegative degrees. Thus $H^{-*}_G(D_GX)$ is an altogether different homology theory, which deserves to be called \emph{coBorel homology}. \cite{GreM} and \cite{Man} write this $cH^G_*(X)$. Following \cite{Jones}, and because $cH^G_*(X)$ is usually supported in the negative direction, we prefer to denote it $H^-_G(X)$. Correspondingly, we write $H^+_G(X) = H^G_*X$, as this is usually supported in the positive direction.

Using the cap product of cohomology and homology, and pulling back cohomology classes from $BG$, we find that $H^+_*(X)$ is a module over $H^{-*}(BG)$ (note the negative grading, as cohomology classes contract against homology, decreasing degree). Using the cup product instead, $H^-_*(X)$ is also a module over $H^{-*}(BG)$.

If $G$ is finite or connected, is a homomorphism relating these two homology theories, the \emph{norm map}: 
$$N_G: H^+_G(X)[\dim G] \to H^-_G(X).$$
Note the degree shift by $\dim G$; we think of the norm map as a sort of averaging operator. For a general compact Lie group there is a twist involving the character $\pi_0 G \to \pm 1$ given by the determinant of the adjoint representation. This norm map is defined by topological means, using spectra; we will have to use a purely algebraic approach to constructing it.

There is a final equivariant homology theory, \emph{Tate homology}, written $H^\infty_*(X)$ and fitting into a long exact sequence $$\cdots \to H^+_{*- \dim G}(X) \xrightarrow{N_G} H^-_{*}(X) \to H^\infty_{*}(X) \to H^+_{*-\dim G -1}(X) \to \cdots$$ Remarkably, the norm map is an $H^{-*}(BG)$-module homomorphism, and Tate homology also has the structure of an $H^{-*}(BG)$-module for which all maps in the above exact sequence are module homomorphisms. The relevant claim in the purely algebraic setting is somewhat more subtle.

The other essential property of Tate homology is that $H^\infty_G(G) = 0$, so that Tate homology of a finite $G$-CW complex can be calculated (essentially) from its subcomplex of points with nontrivial stabilizer. This is extremely useful for calculation, as we will see later in the context of $I^\infty$.

The following table summarizes the analogies between equivariant homology of $G$-spaces and the various Floer homologies of 3-manifolds. The monopole Floer homology is defined to be a sort of $S^1$-equivariant Floer homology, but the $S^1$ symmetry is less visible in Heegaard Floer theory. In any case, the monopole Floer and Heegaard Floer homology groups are modules over $H^{-*}(BS^1;\mathbb Z) = \mathbb Z[U]$, where $|U| = -2$.

\renewcommand{\arraystretch}{1.4}
\begin{center}
\begin{tabular}{|c|c|c|c|}
\hline
Heegaard Floer & Monopole Floer & Instanton Floer & Equivariant homology\\
\hline
\hline
$HF^+_*(Y)$ & $\widecheck{HM}_*(Y)$ & $I^+_*(Y)$ & $H^+_*(X) = H_*(EG_+ \wedge_G X)$ \\
\hline
$HF^-_*(Y)$ & $\widehat{HM}_*(Y)$ & $I^-_*(Y)$ & $H^-_*(X) = H^{-*}_G(D_GX)$ \\
\hline
$HF^\infty_*(Y)$ & $\overline{HM}_*(Y)$ & $I^\infty_*(Y)$ & $H^\infty_*(X)$ \\
\hline
$\widehat{HF}_*(Y)$ & $\widetilde{HM}_*(Y)$ & $\widetilde I_*(Y)$ & $H_*(X)$ \\
\hline
\end{tabular}
\newline
\end{center}

\noindent In our setting, we have a functorial dg-module with dg-algebra action $$\widetilde{CI}(Y,\tilde E,\pi;R) \curvearrowleft C_*^{\textup{gm}}(SO(3);R),$$ acting from the right because $\widetilde{\mathcal B}_{E}$ carries a right $SO(3)$-action; these are right orbits. This is a reasonable notion of `chain complex with action of $SO(3)$', and we would like to obtain invariants that behave like equivariant homology of a $G$-space $X$. The majority of the appendix develops chain complexes $C^\bullet(A;M)$ whose homology gives us equivariant homology theories for $\mathbb Z$-graded dg-modules over a dg-algebra.

The definitions for $\mathbb Z/8$-graded complexes are more delicate, and carried out in Chapter \ref{PeriodicMachine}. Here we take the point of view that the periodic filtration on the module $M$ should give rise to a periodic filtration on $C^\bullet_{SO(3)}(M)$, and to do so we require that $\text{gr}_p \tilde M$ is bounded. Our complexes $\widetilde{CI}$ are equipped with a periodic filtration, and each $\text{gr}_p \widetilde{CI}_\textup{unr}(Y, \tilde E, \pi; R)$ is bounded. In particular, Theorem \ref{4flavors-periodic} is readymade for us to apply to $\widetilde{CI}$.

\begin{theorem}\label{4flavors-instanton}Suppose $R$ is a PID. There is a dg-algebra $C^-\Big(C_*^{\textup{gm}}(SO(3);R)\Big)$ whose homology algebra is graded isomorphic to $H^{-*}(BSO(3);R)$. We denote this dg-algebra by $\mathcal C^-$.

There are functors $$CI^\bullet: \mathsf{Cob}_{3,b}^{U(2),w,\pi} \to \mathsf{Ho}\left(\mathcal C^-\text{-}\mathsf{Mod}_R^{r,\mathbb Z/8}\right),$$ where $\bullet \in \{+, -, \infty\}$, given by applying the constructions $$C^\bullet\left(C_*^{\textup{gm}}(SO(3);R), \widetilde{CI}_\textup{unr}(Y,\tilde E,\pi;R)\right)$$ of Theorem \ref{4flavors-periodic} for periodically filtered modules with finite support on each associated graded piece. After taking homology, we then have functors $$I^\bullet: \mathsf{Cob}_{3,b}^{U(2),w} \to \mathsf{Mod}^{r,\mathbb Z/8}_{H^{-*}(BSO(3);R)}.$$ These fit into an exact triangle $$\cdots \to I^+(Y,E;R) \xrightarrow{[3]} I^-(Y,E;R) \to I^\infty(Y,E;R) \xrightarrow{[-4]} I^+(Y,E;R) \to \cdots$$ where all arrows are $R$-module maps, the arrow from $I^- \to I^\infty$ is an $H^{-*}(BSO(3);R)$-module map, and the other two connecting maps are $H^{-*}(BSO(3);R)$-module maps if $2$ is either invertible or zero in $R$.

Given any other bounded, graded dg-algebra $A$ and right $A$-module $M$, equipped with a periodic filtration each of whose associated graded groups $\textup{gr}_p \tilde M$ is bounded, there are $A$-homology groups $H^\bullet_A(M)$. If there is an algebra map which is a quasi-isomorphism $$f: C_*^{\textup{gm}}(SO(3);R) \to A,$$ and an $f$-equivariant filtered module map $g: \widetilde{CI}(Y, \tilde E, \pi; R) \to M$ whose induced map on associated graded complexes is a quasi-isomorphism, then there is an induced canonical isomorphism $I^\bullet \cong H^\bullet_A(M)$, equivariant under the actions of $H^-_A(R)$ under the induced isomorphism $H^{-*}(BSO(3);R) \cong H^-_A(R)$.
\end{theorem}

\begin{proof}
One thing has been simplified in the statement above: we have removed the reference to a `twisted homology group' $I^{+,\text{tw}}(Y,E;R)$ which appears in the appendix, applying Theoerm \ref{4flavors-periodic} (8) to replace it with $I^+(Y,E;R)[3]$. As stated in Theorem \ref{4flavors-periodic}, this requires that $A = C_*^{\text{gm}}(SO(3);R)$ satisfies `weak Poincar\'e duality of degree $3$' to do this replacement at the level of $R$-modules, and `strong Poincar\'e duality of degree $3$' to do this replacement at the level of $H^{-*}(BSO(3);R)$-modules. The weak duality result --- for any compact Lie group $G$ --- is provided in Proposition \ref{prop:G-is-weak-PD} (one has to pass between $C_*^{\text{gm}}$ and the singular chain algebra using Proposition \ref{sm-gm-comparison}). We only establish the strong duality result for $SO(3)$ when $2 = 0$ or $2$ is invertible in $R$, stated as Corollary \ref{cor:iso-for-SO3}.
\end{proof}

We also have spectral sequence for the equivariant homology groups. As with the non-equivariant homology, we state this in slightly more detail than in Theorem \ref{4flavors-periodic}, point (4).

\begin{theorem}\label{EqIndexSS}For each of $\bullet \in  \{+, -, \infty\}$, if $(Y, E)$ is a weakly admissible bundle equipped with regular perturbation $\pi$, there is a $(\mathbb Z/8, \mathbb Z)$-bigraded spectral sequence of $H^{-*}(BSO(3);R)$-modules from $$E^1_{p,q} = \bigoplus_{\substack{\alpha \in \mathfrak C_\pi\\ \textup{gr}(\alpha) = p}} H^\bullet_{SO(3)}(\alpha; R)$$ whose target is the associated equivariant instanton homology $I^\bullet(Y, E, \pi; R)$.

For $I^-$, we have an isomorphism of the unrolled $(\mathbb Z, \mathbb Z)$-bigraded spectral sequence $E^\infty_{p,q} \cong \textup{gr}_p I^{-, \textup{unr}}_{p+q}$. For any of the theories, a filtered chain map $CI^\bullet(Y, E,\pi) \to CI^\bullet(Y', E',\pi')$ which induces an isomorphism on any finite page $E^r$ induces an isomorphism on $I^\bullet$.
\end{theorem}

\begin{proof}The existence of this spectral sequence, and the fact that we may detect quasi-isomorphisms by isomorphisms on a finite page $E^r$, is precisely Theorem \ref{4flavors-periodic}, point (4), together with the fact that each complex $\text{gr}_p \widetilde{CI}_{\textup{unr}}$ is bounded above (in fact, bounded). As for the specific calculation of the $E^1$ page, our associated graded complex $$\text{gr}_p \widetilde{CI}_\textup{unr}(Y, E, \pi; R) = \bigoplus_{\substack{\alpha \in \mathfrak C_\pi\\ \textup{gr}(\alpha) = p}} C_*^{\text{gm}}(\alpha; R)$$ is a bounded $C_*^{\text{gm}}(SO(3);R)$-module. By Proposition \ref{sm-gm-comparison}, there is a chain of quasi-isomorphisms connecting $C_*^{\text{gm}}(SO(3);R)$ to the algebra $C_*(SO(3);R)$, and similarly a chain of equivariant quasi-isomorphisms between $C_*^{\text{gm}}(\alpha; R)$ to $C_*(\alpha; R)$. In particular, we may identify $$H^\bullet\left(C_*^{\text{gm}}(SO(3)), C_*^{\text{gm}}(\alpha; R)\right) \cong H^\bullet_{SO(3)}(\alpha; R)$$ using Theorem \ref{orbitcalc}. 
\end{proof}

As before, these spectral sequences will become significantly more computable after the introduction of the Donaldson models $\overline{DCI}^\pm$ for $CI^-$ and $CI^+$ in Chapter \ref{secDCI}.

\begin{remark}It is important to note that while the homology theory for $\mathbb Z$-graded complexes $H_A^+(M)$ described at first in the appendix sends $A$-equivariant quasi-isomorphisms $M \to M'$ to isomorphisms on $H_A^+$, this is \emph{not true} for arbitrary filtered $A$-equivariant quasi-isomorphisms between $\mathbb Z/8$-graded complexes. An acyclic $\mathbb Z/8$-graded complex $M$ whose associated graded complex is not acyclic need not have $H_A^+(M) = 0$. 

What remains true, and will be used extensively, is that if $M, M'$ are $\mathbb Z/8$-graded $A$-modules equipped with a periodic filtration, and $f: M \to M'$ is a filtered dg-module map inducing an isomorphism on the associated graded homology groups (the $E^1$ page of the spectral sequence), then $f$ induces an isomorphism on $H_A^+(M)$. The assumption that $M$ is acyclic does not imply that any finite page of the spectral sequence is zero, even if $M$ is a finitely generated $\Bbb Z/8$-graded module, so there is no contradiction to the previous paragraph! 

The reason this happens is that these objects are defined via a completion with respect to the given filtration, which ensures that the associated spectral sequences converge. However, quasi-isomorphisms between the uncompleted complexes which are not quasi-isomorphisms on the associated graded (or some later page of the spectral sequence) do not necessarily induce quasi-isomorphisms between the completions.\end{remark}


\chapter{Examples, calculations, and comparisons}\label{chap:7}
For the entirety of this chapter, we choose elements in each orientation set $\Lambda(\alpha)$ arbitrarily, and suppress these orientation sets from notation.

\section{Equivariant instanton homology for admissible bundles}\label{sec:ex-admissible}
When $E$ is a nontrivial admissible $U(2)$-bundle over an oriented 3-manifold $Y$, all critical orbits are free. The complex $\widetilde{CI}(Y,E;R)$ resembles the equivariant Morse complex of a finite-dimensional free $SO(3)$-manifold $M$, so our heuristic is that the equivariant homology groups $I^\bullet(Y,E;R)$ should behave like the equivariant homology $H^\bullet_G(M)$: 
\begin{align*}H^-_G(M;R)[\dim M] &\cong H^+_G(M;R) \cong H(M/G;R)\\
H^\infty_G(M;R) &= 0.\end{align*}

\begin{theorem}\label{admissible-iso}Suppose $E$ is a nontrivial admissible bundle over a 3-manifold $Y$. Then $I^\infty(Y, E) = 0$, and if $CI(Y, E)$ is Floer's instanton chain complex for admissible bundles, there is a natural quasi-isomorphism $CI^+(Y, E;R) \to CI(Y, E;R)$. In particular, the equivariant Floer homology $I^+(Y, E)$ of a nontrivial admissible bundle is Floer's original instanton homology group $I(Y, E)$. Furthermore, there is also a functorial quasi-isomorphism $CI(Y, E;R)[3] \to CI^-(Y, E;R)$. 
\end{theorem}
\begin{proof}We first prove the statement for $CI^\infty$. The $E^1$ page of the index spectral sequence is $$\bigoplus_{\substack{\alpha \in \mathfrak C, \;n \in \mathbb Z \\ n = i(\alpha) \!\!\!\mod 8}} H^\infty_{*}(\alpha;R)[i(\alpha)].$$ Now all orbits are free orbits $\alpha \cong SO(3)$, and so we can calculate this as $$\bigoplus H^\infty_{*}(SO(3);R)[i(\alpha)].$$ That $H^\infty_{SO(3)}(SO(3))$ vanishes is one of the defining features of Tate homology: see Proposition \ref{eq-package} (4). So the $E^1$ page is identically zero, and by the $I^\infty$ part of Theorem \ref{EqIndexSS} we see that the map $0 \to CI^\infty(Y, E;R)$ is a quasi-isomorphism, as desired.

Floer's instanton complex $CI(Y, E)$ is defined as the free $R$-module generated by critical orbits $\alpha \subset \mathfrak C$, the relative grading defined still as $\text{gr}_z(\alpha, \beta)$; the differential counts points in zero-dimensional moduli spaces of \textit{unframed} instantons between $\alpha$ and $\beta$. As a graded $\mathcal C^-$-module, 
\begin{align*}CI^+ &= \bigoplus_{\alpha \in \mathfrak C} C^+(\alpha;R)[i(\alpha)],\\
CI &= \bigoplus_{\alpha \in \mathfrak C} R[i(\alpha)].\end{align*}
The differential is a sum of the induced differential on each factor $C^+(\alpha;R)$ and the maps $C^+(\alpha;R) \to C^+(\beta;R)$ induced by the $C^\text{gm}_*(SO(3))$-equivariant (anti)chain maps $$\sigma \mapsto (-1)^{|\sigma|} \sigma \times_{e_-} \overline{\mathcal M}(\alpha, \beta).$$ 

Recall from Lemma \ref{groupcoho} that $H_{SO(3)}^+(SO(3);R) = R$, concentrated in degree zero. There is a chain map $\varphi: CI^+(Y, E;R) \to CI(Y, E;R)$ given on the summand for each critical orbit by the augmentation map $C^+_{SO(3)}(\alpha;R) \to R$ induced by the point-counting map $C_*(\alpha;R) \to R$; this map kills everything in degree larger than $i(\alpha)$ and sends $C^+_0(\alpha;R) = C_0(\alpha;R) \to H^+_{i(\alpha)} = R$ by the augmentation. 

To see that $\varphi$ is a chain map, note that if $\sigma \in CI_k^+(Y, E;R)$ is not sent to zero, it can be written as a sum of points in 
$$\bigoplus_{\substack{\alpha \in \mathfrak C_\pi \\ i(\alpha) = k}}C_0(\alpha;R).$$ 
So it suffices to check that $\varphi(\tilde d\sigma) = d\varphi(\sigma)$ for $\sigma = p$, where $p$ is a point in $\alpha$. The only components of $d\sigma \in CI^+$ not automatically sent to zero are those in some $C_0(\beta;R) \subset C^+(\beta;R)$, which arise from taking fiber products with moduli spaces $\overline{\mathcal M}(\alpha, \beta)$ of dimension 3. Because the $SO(3)$ action on $\overline{\mathcal M}$ is free, the count of points of $p \times_\alpha \overline{\mathcal M}(\alpha, \beta)$ is the same as that of $\overline{\mathcal M}(\alpha, \beta)/SO(3)$, Floer's unframed moduli space of instantons. In particular, we find that $\varphi$ is a chain map, as desired.

Now, on the unrolled complexes, this map is a filtered chain map for the tautological periodic filtration of $CI(Y, E;R)$ by degree, and Theorem \ref{EqIndexSS} says that the above map induces an isomorphism $$E^1\left(CI^+(Y, E;R)\right) \to E^1\left(CI(Y, E;R)\right) = CI(Y, E;R).$$ Therefore, by the comparison theorem, the induced map $$I^+(Y, E;R) \to I(Y, E;R)$$ is an isomorphism. That these isomorphisms fit into a commutative square with the corresponding cobordism maps is proved the same way. 

The statement about $I^-$ follows because of the exact triangle relating $I^\infty$, $I^+$, and $I^-$ and the fact that $I^\infty$ vanishes; one may alternately define a chain map $CI[3] \to CI^-$ explicitly, almost precisely as above: instead of projecting to $R$ using point-counts, we include $R$ via the inclusion of the fundamental class.
\end{proof}

Suppose $\frac 12 \in R$; then $H^{-*}(BSO(3);R) = R[U]$, where $|U| = -4$. Floer's $I(Y,E;R)$ also carries an action by this ring (when $\frac 12 \in R$), and we should check that the above map preserves the $U$-action. Because the $U$-action carries things down vertically in the spectral sequence, but $E^2(CI^+)$ is here concentrated on a single horizontal line, this is not a theorem well-suited to a spectral sequence proof. Instead, we must get our hands dirty at the chain level. This will be delayed until the following section, where we construct a simpler chain-level model of $\widetilde{CI}$.

\section{Comparison with Donaldson's theory}\label{secDCI}
In this section, we describe a finite-dimensional complex computing framed instanton homology $\widetilde{I}$, and then explain how to use it to calculate the equivariant instanton groups. While it seems likely that this is possible over any PID $R$, the situation is drastically simpler when $\frac 12 \in R$. So for the rest of this section, $R$ is a PID in which $2$ is invertible; the ring $R$ will be dropped from the notation whenever possible. 

If $2$ is invertible in $R$, then $H_*(SO(3);R)$ is isomorphic to a single copy of $R$ in degrees $0$ and $3$ and zero otherwise. We denote this $R$-algebra $\Lambda(u)$, the exterior algebra on a generator $u$ in degree 3. There is a dg-algebra homomorphism $$i: \Lambda(u) \to C_*^{\textup{gm}}(SO(3);R)$$ given by picking out the identity in degree 0 and the fundamental class in degree 3; this map induces the identity on homology. By Proposition \ref{dg-invar}, because $i$ is a quasi-isomorphism, we have 
$$H^\bullet(C_*^{\textup{gm}}(SO(3);R), \widetilde{CI}) \cong H^\bullet(\Lambda(u), \widetilde{CI}).$$ 
For the purposes of computing the equivariant instanton homology groups $I^+, I^-$, and $I^\infty$, it therefore suffices to consider $\widetilde{CI}$ as a module over $\Lambda(u)$.

Our goal is to write down a differential on the $\mathbb Z/8$-graded $R$-module $$DCI(Y, E,\pi;R) = \oplus_{\alpha} H_{*}(\alpha;R)[i(\alpha)],$$ which has a periodic filtration by index of the orbit $\alpha$. This differential should \textit{decrease} the filtration, and there should be a map $\widetilde{CI}(Y, E,\pi) \to DCI$ which is the identity on the $E^1$ page, or something like it. To produce this, we use what is called the \emph{homological perturbation lemma}. Recall its statement: 

\begin{lemma}[Homological perturbation lemma]\label{HPL}Suppose $(C,d)$ and $(C',d')$ are chain complexes, equipped with an inclusion $i: C \hookrightarrow C'$ and a projection $p: C' \to C$ so that $pi = 1$, both quasi-isomorphisms, as well as a degree 1 map $h: C' \to C'$ serving as a homotopy witnessing this. That is, $$ip = 1 + d'h + h d'.$$ This data is called deformation retract data, and is depicted $$(C,d) \substack{\xrightarrow{i}\\ \xleftarrow[p]{}} (C',d') \! \rcirclearrowleft \! h.$$

Suppose that $C'$ is equipped with a \emph{deformation}: an additional map $\delta: C \to C$ so that $(d' + \delta)^2 = 0$. Suppose that $h^n = 0$ for sufficiently large $n$. Write $A = \sum_{n=0}^\infty (\delta h)^n\delta$. 

Then $(C, d + pAi)$ is a chain complex, equipped with deformation retract data $$(C,d+pAi) \substack{\xrightarrow{i+hAi}\\ \xleftarrow[p+pAh]{}} (C',d'+\delta) \! \rcirclearrowleft \! h+hAh.$$

If $C$ and $C'$ are dg-modules over a dga, and all of $i, p,h, $ and $\delta$ are dg-module homomorphisms, then the same is true of the perturbed deformation retract data. In particular $(C,d+pAi)$ is homotopy equivalent to $(C',d'+\delta)$ as a dg-module. 

If $C$ and $C'$ are filtered, and $i, p, h,$ and $\delta$ all preserve the filtration (in the sense that $f(F_k) \subset F'_k$ for all $k$), then the same is true for the perturbed deformation retract data. In particular, $(C,d+pAi)$ is a filtered complex, which is filtered homotopy equivalent to $(C',d'+\delta)$.  
\end{lemma}

A very explicit reference is \cite{crainic2004perturbation}. 

To apply the homological perturbation lemma, first recall that we may write the differential on $\widetilde{CI}$ as $d + d_M$. The first term is the usual boundary operator on geometric chains, and the second term is the contribution from fiber products with moduli spaces. We want to apply the lemma to $(DCI, 0)$ and $(\widetilde{CI}, d)$, with $\delta = d_M$. 

Choose a basepoint $b_\alpha \in \alpha$ for each critical orbit $\alpha$; if $\alpha$ is an $SO(2)$-reducible, we demand that the stabilizer of $b_\alpha$ is the standard $SO(2) \subset SO(3)$. 

Our inclusion map $$i: DCI \to \widetilde{CI}$$ is given on $H_0(\alpha; R) \to C_*^{\text{gm}}(\alpha;R)$ by sending the generator to the basepoint $b_\alpha$, and on $H_{\dim \alpha}(\alpha;R) \to C_*^{\text{gm}}(\alpha;R)$ by sending a generator to the corresponding fundamental class of $\alpha$. (Note that the choice of fundamental class depends on an orientation of $\alpha$; we are making this choice by choosing a generator of top homology.)

To define the map $p$, choose basepoints $q_\alpha \in \alpha$. Recall from Definition \ref{trans-chains} and Lemma \ref{qi-trans} that, given a countable family $\mathcal F$ of maps from $\delta$-chains to a space $X$, the geometric chain complex has a quasi-isomorphic subcomplex $C_*^{\textup{gm},\mathcal F}(X;R)$, spanned by those chains which are transverse on every stratum to the countable family of maps. If $X$ is an $SO(3)$-orbit, this is a $\Lambda(u)$-invariant subcomplex, because the evaluation map $SO(3) \to X$ is a submersion. This means that the product of any chain with the fundamental class is transverse to any element of $\mathcal F$. 

To define the map $$p: C_*^{\text{gm},\mathcal F}(\alpha;R) \to H_*(\alpha;R),$$ take $\mathcal F$ to include the basepoint $q_\alpha$. In degree zero, any 0-dimensional $\delta$-chain has an underlying oriented point, and the projection is given by counting these with sign. In top degree, we send a $\delta$-chain $\sigma: P \to \alpha$ to the signed count $\# \sigma^{-1}(q_\alpha)$. Note that to define the sign of this signed count, we need to choose an underlying orientation of $\alpha$, which determines the generator of $H_{\dim \alpha}(\alpha; R)$ we should send $\sigma$ to. Changing the orientation of $\alpha$ changes this sign twice, so the resulting map is independent of any choice of orientation of $\alpha$.

Unfortunately, the fiber products with moduli spaces do not send chains transverse to $q_\alpha$ to chains transverse to $q_\beta$. So this submodule of $\widetilde{CI}$ is \emph{not} a subcomplex. To fix this, we enlarge the collection $\mathcal F$ that our chains must be transverse to. 

Denote $\overline{\mathcal M}^q_{\alpha \beta}$ for the fiber above $q_\beta$ of the endpoint map $$\overline{\mathcal M}_{\alpha \beta} \to \beta.$$ More generally, given a sequence of orbits $\gamma_1, \cdots, \gamma_n$, denote $\overline{\mathcal M}^q_{\alpha \vec{\gamma} \beta}$ denote the fiber above $q_\beta$ of the map $$\overline{\mathcal M}_{\alpha \gamma_1} \times_{\gamma_1} \cdots \times_{\gamma_n} \overline{\mathcal M}_{\gamma_n \beta} \to \beta.$$ 

For the orbit $\alpha$, we let the collection $\mathcal F$ consist of $q_\alpha$ and all of the endpoint maps of moduli spaces $$\overline{\mathcal M}^q_{\alpha \vec{\gamma} \beta} \to \alpha.$$ Because this collection is closed under fiber products with moduli spaces on the left, the condition of being transverse to $\mathcal F$ is preserved by $d_M$. So for this collection $\mathcal F$, we have a $\Lambda(u)$-subcomplex $$\widetilde{CI}^{\mathcal F}(Y,E,\pi;R) \subset \widetilde{CI}(Y,E,\pi;R)$$ for which the projection map $p$ above is defined and a $\Lambda(u)$-module map. Running the index spectral sequence, and using that the inclusion $C_*^{\textup{gm},\mathcal F}(\alpha;R) \hookrightarrow C_*^{\textup{gm}}(\alpha; R)$ is a quasi-isomorphism for each orbit, we see that the inclusion map $i: \widetilde{CI}^{\mathcal F} \to\widetilde{CI}$ is an equivariant quasi-isomorphism. Furthermore, it is a free $R$-submodule; its generators are precisely those generators of $\widetilde{CI}$ for which all strata are transverse to the given family of $\delta$-chains.

To ensure that the image of $i: DCI \to \widetilde{CI}$ lies in $\widetilde{CI}^{\mathcal F}$, we need to demand that $b_\alpha$ is a regular value of $\overline{\mathcal M}^{q}_{\alpha \vec{\gamma} \beta} \to \alpha$; one may symmetrically demand that $q_\beta$ is the regular value of the endpoint map from the fiber above $b_\alpha$. It is not difficult to ensure that the $q_\alpha$ satisfy this, by an application of Sard's theorem. 

An $R$-linear chain homotopy $h_\alpha$ is guaranteed to exist because $H_*(\alpha; R)$ is a free $R$-module and $ip = 1$ on homology. Any such map $h$ is $\Lambda(u)$-equivariant for $\alpha$ reducible for degree reasons. We will show that any such $h$ is also $\Lambda(u)$-equivariant for $\alpha$ irreducible. The equivariance condition is automatic except for elements $x \in C_0(\alpha)$. If $x$ is a singleton in $\alpha$, then $h(x)$ is a $1$-manifold with boundary $x - \{b_\alpha\}$; then $I \times SO(3) \to \alpha$ has boundary given by the difference $[\alpha] - [\alpha]$, hence represents zero in $C_4^{\text{gm}}(\alpha; R)$. Then equivariance is equivalent to the statement that $h([\alpha]) = 0$, which is automatic: because $(dh-hd)([\alpha]) = (ip-1)([\alpha]) = 0$, we have $d h([\alpha]) = 0$, but $d: C_4^{\text{gm}}(\alpha; R) \to C_3^{\text{gm}}(\alpha; R)$ is injective.\\

This is enough to apply the homological perturbation lemma and get \emph{some} result, a differential on $DCI(Y,E,\pi;R)$ which is $\Lambda(u)$-equivariantly homotopy equivalent to $\widetilde{CI}^{\mathcal F}(Y,E,\pi;R)$. 

Take the complex $\widetilde{CI}_1(Y,E,\pi;R)$ equipped with the differential $d+d_M$ induced by the above fiber products with moduli spaces; we described the construction of deformation retract data $$(DCI,0) \substack{\xrightarrow{i}\\ \xleftarrow[p]{}} (\widetilde{CI}_1,d+d_M) \! \rcirclearrowleft \! h$$ above. Applying the homological perturbation lemma, we arrive at a chain complex $(DCI, \partial_{DCI} = pd_M i + p d_M h d_M i + \cdots)$ with $\Lambda(u)$-action, and an equivariant homotopy equivalence to $(\widetilde{CI}_1, d+d_M)$. We now describe the differential $\partial_{DCI}$. 

Denote by $C_*^\text{irr}$ the free $R$-module on the irreducible orbits $\alpha \subset \mathfrak C_\pi$ (graded the same), $C_*^{U(1)}$ the free $R$-module on $U(1)$-reducible orbits, and $C_*^\theta$ the free $R$-module on the full reducibles. Set $$DCI(Y,E,\pi) = C_*^\text{irr} \oplus C_*^\text{irr}[3] \oplus C_*^{U(1)} \oplus C_*^{U(1)}[2] \oplus C_*^\theta.$$ As graded $R$-modules, $DCI(Y,E,\pi) = \bigoplus H_{*-i(\alpha)}(\alpha)$. 

The differential $\partial_{DCI}$ that $DCI$ inherits is given by a matrix 
$$\widetilde{\partial}_{CI} := \begin{pmatrix}\partial_1 & 0 & 0 & 0 & 0\\
U_{\text{Fl}} & -\partial_1 & V_4 & V_2 & D_2 \\
V_1 & 0 & 0 & 0 & 0\\ 
V_3 & 0 & 0 & 0 & 0\\
D_1 & 0 & 0 & 0 & 0\end{pmatrix}.$$ 
(We will momentarily define the terms in this matrix.) Here we write ``$\dim \sigma$" to mean $0$ for the component $C_*^{\text{irr}}$ and $3$ for the component $C_*^{\text{irr}}[3]$. 

The action of $u$ on $\widetilde{CI}$ is defined by the matrix whose only nonzero term is $\text{id}: C_*^\text{irr} \to C_*^\text{irr}[3]$. It is easy to see that $u\widetilde{\partial_{CI}} = \widetilde{\partial_{CI}}u = \partial_1$, considered as a map $C_*^\text{irr} \to C_{*-1}^\text{irr}[3]$.

Unfortunately, the differential is still more complicated than we would like: it takes the form $$p\sum_{n=0}^\infty (d_M h)^n d_M i,$$ where $d_M$ is the part of the differential coming from fiber products with moduli spaces. It is not impossible that many components of this may be nonzero: one takes a fiber product of the basepoint $b_\alpha$ with some moduli space, pushes it forward to $\beta$, applies the homotopy $h$, takes another fiber product, and then counts intersections with $b_\gamma$. One could easily imagine that, starting and ending at irreducibles $\alpha, \gamma$, one takes the fiber product with a moduli space that increases the dimension by $2$, cones towards $b_\beta$, and then takes a fiber product with a $0$-dimensional moduli space and gets a result that has nonzero degree above $b_\gamma$.

Nevertheless, these higher components of $\langle \partial_{DCI} \alpha, \beta\rangle$ can only be nonzero if $\overline{\mathcal M}_{\alpha \beta}$ has nonempty boundary, and we can describe the component $p d_M i$ completely.

As before, we denote $$\mathcal M_{\alpha\beta} = \overline{\mathcal M}_{\alpha\beta}/SO(3);$$ these are the spaces of unframed flowlines between the underlying connections $[\alpha], [\beta]$, having forgotten the framings. Define $X_{\alpha \beta} = (\alpha \times \beta)/SO(3)$. The operators in this matrix all arise of the form $T(\alpha) = \sum_\beta n(\alpha, \beta) \beta$, where the sum is taken over $\beta$ in the right degree, and $n(\alpha, \beta)$ arises either by counting points (with orientation) in $\mathcal M_{\alpha\beta}$, or degree of the map $\mathcal M_{\alpha\beta} \to X_{\alpha \beta}$, as measured by the number of points lying above $[b_\alpha, q_\beta]$. We've chosen $b_\alpha, q_\beta$ to ensure that this is a regular value (and in the case that both $\alpha$ and $\beta$ are $U(1)$-reducibles, to have different stabilizer). 

\begin{itemize}
    \item $\partial_1$ counts points in $\mathcal M_{\alpha\beta}$, when $\alpha, \beta$ are both irreducible and $\text{gr}(\alpha, \beta) = 1$.
    \item $D_1$ counts points in $\mathcal M_{\alpha\beta}$ when $\alpha$ is irreducible, $\beta$ is fully reducible, and $\text{gr}(\alpha, \beta) = 1$.
    \item $D_2$ counts points in $\mathcal M_{\alpha\beta}$ when $\alpha$ is fully reducible, $\beta$ is irreducible, and $\text{gr}(\alpha, \beta) = 4$. \footnote{Recall that this means that $\overline{\mathcal M}$ is of dimension $\dim \alpha + \text{gr}(\alpha, \beta) - 1$; because $\alpha$ is fully reducible, this means that $\overline{\mathcal M}$ is 3-dimensional. Because the endpoint map to $\beta$ is equivariant, and $\beta$ is irreducible, $\overline{\mathcal M}$ is a finite set of orbits with no stabilizer, and $\mathcal M$ is 0-dimensional.} 
    \item The $p d_M i$ component of $U_\text{Fl}$ counts the number of points in the fiber above $[b_\alpha, q_\beta]$ of the map $\mathcal M_{\alpha\beta} \to X_{\alpha \beta}$, when $\alpha$ and $\beta$ are both irreducible and $\text{gr}(\alpha, \beta) = 3$. 
    \item $V_1$ counts points in $\mathcal M_{\alpha\beta}$ when $\alpha$ is irreducible, $\beta$ is $U(1)$-reducible, and $\text{gr}(\alpha, \beta) = 1$.
    \item $V_2$ counts points in $\mathcal M_{\alpha\beta}$ when $\alpha$ is $U(1)$-reducible, $\beta$ is irreducible, and $\text{gr}(\alpha, \beta) = 2$.
    \item The $p d_M i$ component of $V_3$ counts points above $[b_\alpha, q_\beta]$ when $\alpha$ is irreducible, $\beta$ is $U(1)$-reducible, and $\text{gr}(\alpha, \beta) = 3$.
    \item The $p d_M i$ component of $V_4$ counts points above $[b_\alpha, q_\beta]$ when $\alpha$ is $U(1)$-reducible, $\beta$ is irreducible, and $\text{gr}(\alpha, \beta) = 4$.
\end{itemize}

The derivation of these terms is relatively self-explanatory; the sign arises from the factor $(-1)^{\dim \sigma}$ in $d_M$. The one small subtlety is why there is no matrix element corresponding to counting points above $[b_\alpha, q_\beta]$ when both $\alpha, \beta$ are $U(1)$-reducible and $\text{gr}(\alpha, \beta) = 3$. By Proposition \ref{ModuliPackage} (4), this is impossible: the relative grading of reducible connections is always even!\\

In some cases, it is possible to eliminate the higher terms in the differential $\partial_{DCI}$ by a homotopy, an approach suggested by \cite[Section~7.3.2]{Don}. This can be done by modifying the endpoint maps, as in the following lemma; we include this for comparison to \cite{Don}, as we will not need to carry out any homotopy to prove our main results.

\begin{lemma}\label{EndpointHomotopy}Suppose $e_{\alpha \beta}: \overline{\mathcal M}_{\alpha \beta} \to \alpha \times \beta$ are the equivariant endpoint maps associated to moduli spaces of instantons on a weakly admissible bundle $(Y,E)$ with respect a fixed perturbation $\pi$.

Suppose we choose a collection of $e^t_{\alpha \beta}: \overline{\mathcal M}_{\alpha \beta} \to \alpha \times \beta,$ homotopies through equivariant maps that are smooth on each stratum, one homotopy for each pair $(\alpha,\beta)$, so that $e^0_{\alpha \beta} = e_{\alpha \beta}$ above. Furthermore demand that these are compatible in the sense that $$e^t_{\alpha \beta} \times_\beta e^t_{\beta \gamma}: \overline{\mathcal M}_{\alpha \beta} \times_\beta \overline{\mathcal M}_{\beta \gamma} \to \alpha \times \gamma$$ agrees with the restriction of $e^t_{\alpha \gamma}$ to that stratum of $\overline{\mathcal M}_{\alpha \gamma}$.

We may define a chain complex $\widetilde{CI}_t(Y,E,\pi)$, identical as a $C_*^{\textup{gm}}(SO(3))$-module to $\widetilde{CI}$, but whose differential uses the fiber product with the same moduli spaces but endpoint maps $e^t_{\alpha \beta}$. Then $\widetilde{CI}(Y,E,\pi)$ is $\Lambda(u)$-equivariantly quasi-isomorphic to $\widetilde{CI}_1(Y,E,\pi)$.
\end{lemma}

\begin{proof}
We identify the simplex $\Delta^n$ with the subset $\Delta^n \subset [0,1]^n$ of nondecreasing sequences, with the induced orientation. We write $\Delta^n_i$ for the subset of $(t_1, \cdots, t_n)$ with $t_i = t_{i+1}$, equipped with the corresponding boundary orientation; here one interprets the cases $i = 0, n+1$ as $t_0 = 0$ and $t_{n+1} = 1$, respectively. There is a transparent order-preserving identification $\Delta^{n-1} \cong \Delta^n_i$, and this identification is oriented up to a sign of $(-1)^{i+1}$.

Let $\alpha, \gamma_1, \cdots, \gamma_n, \beta$ be a sequence of critical orbits. Write 

\begin{align*}\overline{\mathcal M}^{\Delta^{n+1}}_{\alpha \vec{\gamma} \beta} \subset &\Delta^n \times \mathcal M_{\alpha \gamma_1} \cdots \mathcal M_{\gamma_n \beta}\\
\overline{\mathcal M}^{\Delta^{n+1}}_{\alpha \vec{\gamma} \beta} = \{(t_1, \cdots, &t_{n+1}), \mathbf{A}_1, \cdots, \mathbf{A}_{n+1} \mid e_+^{t_i}\mathbf{A}_i = e_-^{t_{i+1}} \mathbf{A}_{i+1}\},
\end{align*}

whose fiber is oriented as $$\Delta^{n+1} \times \text{Fib}(\overline{\mathcal M}_{\alpha \gamma_1}) \times \cdots \times \text{Fib}(\overline{\mathcal M}_{\gamma_n \beta}),$$ where $\text{Fib}(\overline{\mathcal M}_{\alpha \beta})$ indicates a generic fiber above a point in $\alpha$. Given any oriented submanifold $S \subset \Delta^{n+1}$, restriction to $S$ gives another fiber-oriented manifold.

A straightforward but tedious computation shows that, as fiber-oriented manifolds, we have 

\begin{align*}(-1)^{\dim \alpha}\partial \overline{\mathcal M}^{\Delta^{n+1}}_{\alpha \vec{\gamma} \beta} \cong &\bigcup_{i=0}^{n+1} \overline{\mathcal M}^{\Delta^{n+1}_i}_{\alpha \vec{\gamma} \beta} \\
&\bigcup_{j=0}^n \bigcup_\eta (-1)^{\text{gr}(\alpha, \eta) + n+1} \overline{\mathcal M}^{\Delta^{n+2}_{j+1}}_{\alpha, \vec{\gamma}_{\leq j}, \eta, \vec{\gamma}_{>j}, \beta}.\end{align*}

Write $\widetilde{CI}_t$ for the chain complex given by using the moduli spaces with endpoint maps $e^t_\pm$; the underlying $R$-modules are all the same, it is only the differential that changes. Consider the maps $H_n: \widetilde{CI}_0 \to \widetilde{CI}_1$, given for each $n>0$ by

$$H_n(\sigma) = \sum_{\gamma_1, \cdots, \gamma_n, \beta} (-1)^{\sum_{i=1}^n (\text{gr}(\alpha, \gamma_i) + i-1)} \sigma \times_\alpha \overline{\mathcal M}^{\Delta^{n+1}}_{\alpha \vec{\gamma} \beta}.$$ 

For convenience, we write the sign here\footnote{We could remove this sign by orienting $\overline{\mathcal M}^{\Delta^{n+1}}$ by spreading the orientation of $\Delta^{n+1}$ throughout the iterated fiber product, placing the orientation corresponding to the $t_i$ direction before the $i$'th fiber product. While this may be conceptually clearer, the author hopes the current format is easier to read.} as $s(n, \alpha, \vec{\gamma}, \beta)$; it satisfies the recurrence 
$$s(n-1, \alpha, \vec{\gamma}_{<j}, \vec{\gamma}_{>j}, \beta) \cdot (-1)^{\text{gr}(\alpha, \gamma_j) + n-1} = s(n, \alpha, \vec{\gamma}, \beta).$$

Our desired chain map will be $\text{Id} + \sum_{n>0} H_n$. The map $H_n$ is identically zero for degree reasons when $n>3$, so this sum is finite. Because $H_n$ increases the filtration depth, we see that the induced map on associated graded complexes is the identity; so as long as this is a chain map, the usual spectral sequence arguments will imply it is a quasi-isomorphism. 

First we explicitly calculate $\partial_1 H_n - H_n \partial_0$. This is given by 

\begin{align*}(-1)^{\dim \sigma}(\partial_1 H_n - H_n \partial_0)(\sigma) &= \sum_{\vec{\gamma}, \beta} (-1)^{\dim \alpha}s(n,\alpha, \vec{\gamma}, \beta) \sigma \times_\alpha \partial \mathcal M^{\Delta^{n+1}}_{\alpha, \vec{\gamma}, \beta}\\
&- \sum_{\eta, \vec{\gamma}, \beta} s(n, \eta, \vec{\gamma}, \beta) (-1)^{(\text{gr}(\alpha, \eta)+1)(n+1)} \mathcal M^{\Delta^{n+2}_0}_{\alpha, \eta, \vec{\gamma}, \beta}\\
&+ \sum_{\vec{\gamma}, \eta, \beta} s(n, \alpha, \vec{\gamma}, \eta) (-1)^{\text{gr}(\alpha, \eta) + n+1} \mathcal M^{\Delta^{n+2}_{n+2}}_{\alpha, \vec{\gamma}, \eta, \beta}.
\end{align*}

The claim, then, is that the sum $\partial_1 - \partial_0 + \sum_{n >0} \partial_1 H_n - H_n \partial_0$ is a telescoping sum. Expanding out the decomposition of $\partial \mathcal M^{\Delta^{n+1}}_{\alpha, \vec{\gamma}, \beta}$ above, we see that this telescoping amounts to the claims that 

\begin{align*}s(n, \alpha, \vec{\gamma}, \beta) &= s(n-1, \alpha, \vec{\gamma}_{<i}, \vec{\gamma}_{>i}, \beta) \cdot (-1)^{\text{gr}(\alpha, \gamma_i)+n-1} \;\;\;\;\; 0 < i < n\\
s(n, \alpha, \vec{\gamma}, \beta) &= s(n-1, \gamma_1, \vec{\gamma}_{>1}, \beta) \cdot (-1)^{n\text{gr}(\alpha, \gamma_1) + n-1}\\
s(n, \alpha, \vec{\gamma}, \beta) &= s(n-1, \alpha, \vec{\gamma}_{<n}, \gamma_n) \cdot (-1)^{\text{gr}(\alpha, \gamma_n)+n-1}.
\end{align*}

The first is just the recurrence stated above, and the last is essentially the same. The middle term in addition uses that $\text{gr}(\alpha, \gamma_1) + \text{gr}(\gamma_1, \eta) = \text{gr}(\alpha, \eta)$.

Thus $\text{Id} + \sum_{n>0} H_n: \widetilde{CI}_0 \to \widetilde{CI}_1$ is an $SO(3)$-equivariant quasi-isomorphism. 
\end{proof}

We summarize the result of the homological perturbation lemma discussion as follows.

\begin{corollary}\label{DCI}The $R$-module with differential $(DCI,\partial_{DCI})$ described above is a dg-module over $\Lambda(u)$, and is $\Lambda(u)$-equivariantly homotopy equivalent to $\widetilde{CI}$. 
\end{corollary}

\begin{remark}In \cite[Section~7.3.3]{Don}, Donaldson sketches the definition of a complex $\widetilde{CF}$ associated with an integer homology sphere. Donaldson's complex includes the data of the filtration\footnote{Thinking of $DCI$ as a $\mathbb Z/8 \times \mathbb Z$-graded multicomplex with $d_0 = 0$, this is the filtration corresponding to the $\mathbb Z$ factor, while ours is the periodic filtration corresponding to the $\mathbb Z/8$ factor.} by degree of an element in $H_* \alpha$, as well as the action of $u \in \Lambda(u)$. 

When the holonomy maps $\mathcal M_{\alpha \beta} \to SO(3)$ lift to $SU(2)$ in a multiplicative vashion for $\dim \mathcal M_{\alpha \beta} \le 3$, Donaldson introduces in Section 7.3.2 a method for reducing the contribution of small-dimensional moduli spaces, so that the higher terms $p (d_M h)^n d_M i$ are zero. It is not typically possible to construct such a lift; it is impossible, for instance, for the Poincar\'e homology sphere $Y = \Sigma(2,3,5)$. In this general case, it is not entirely clear how to define the complex $\widetilde{CF}$.

When $\widetilde{CF}$ can be defined, it follows that $\widetilde{CI}$ is equivariantly homotopy equivalent to Donaldson's $\widetilde{CF}$.
\end{remark}

In what follows we use the complex $DCI$ to produce small models which we write as $\overline{DCI}^\pm$ for $CI^+$ and $CI^-$. In the case of integer homology spheres, these are essentially the same as what Donaldson writes $\overline{\overline{CF}}$ and $\uline{\uline{CF}}$.

As described in the appendix, the \emph{reduced bar construction} of an augmented dg-algebra $A$ and a right dg-module $M$ is given, as a graded $R$-module, as $$B(M,A,R) := \bigoplus_{n=0}^\infty M \otimes \overline A[1]^n,$$ where $\overline A = \ker \varepsilon$ and $\varepsilon$ is the augmentation. For us, $A = \Lambda(u)$, and $\overline A$ is just a copy of $R$ concentrated in degree $3$. Therefore, we may write $$B(M, A, R) = \bigoplus_{n=0}^\infty  M[4n].$$ We write this as $M \otimes_R R[U^*]$, where $|U^*| = 4$.\footnote{We choose the notation $U^*$ here because $U$ is reserved for the degree $-4$ operation coming from $H^{-4}(BSO(3);R)$.} The differential $d^+$ is seen to be 
$$d^+(m \otimes (U^*)^n) = um \otimes (U^*)^{n-1} + (-1)^n dm \otimes (U^*)^n,$$ 
where $(U^*)^{-1}$ should be interpreted as zero. Write $U_{\text{alg}}^+$ for the operator with $$U_{\text{alg}}^+\left(m \otimes (U^*)^n\right) = (-1)^n um \otimes (U^*)^{n-1};$$ again, this means that $U_{\text{alg}}^+(1) = 0$.

Now, when we apply this to a $\mathbb Z/8$-graded module $M$ (with finite basis of finite-dimensional $\Lambda(u)$-modules degreewise isomorphic to $R$) as in Chapter \ref{PeriodicMachine}, we need to \emph{full-complete}. This is a two-step process. 

First, we unroll $CI$ to $CI^{unr}$, a $\Bbb Z$-graded and $8\Bbb Z$-periodic complex. We then apply $$B(CI^{unr}, A, R) = CI^{unr} \otimes R[U],$$ the isomorphism one of $R$-modules. This is filtered by $B(CI^{unr}_{\le p}, A, R) = CI^{unr}_{\le p} \otimes R[U]$, where $CI^{unr}_{\le p}$ is the part of the complex given by summands with $i(\alpha) \le p$. The full-completed complex is obtained in two steps. First, we pass to a colimit to ensure that everything lies in the union of this filtration; we set $$\colim_p B(CI^{unr}_{\le p}, A, R) = \colim_p CI^{unr}_{\le p} \otimes R[U] = CI^{unr} \otimes R[U].$$ The last equality follows because for any module $M$ with an exhaustive filtration, we have $\cup_p B(F_p M, A, R) = B(M, A, R)$. 

Next, we take a quotient and a limit to ensure the filtration is Hausdorff and complete. It is Hausdorff to start with, but not complete --- this is where something changes. Here we take $$\lim_{q \to -\infty} B(CI^{unr}/CI^{unr}_{\le q}, A, R) = CI/CI^{unr}_{\le q} \otimes R[U].$$ In degree $d$, an element of this tensor product is a sequence $$(\dots, x_{d-8}, x_d)$$ of elements of $CI_d$ so that $x_i = 0$ if $i \le q$ (because of the quotient) and so that $x_i = 0$ for $i$ large enough (because we are taking finite sums). Taking the limit as $q \to -\infty$ we obtain sequences $(\cdots, x_{d-8}, x_d, x_{d+8}, \cdots)$ with $x_i = 0$ for $i$ large enough, but no similar condition for $i$ small. (Another way to think about this is as the completion of a multicomplex, as in Definition \ref{multicomplex}.) Quotienting by $\Bbb Z$, we may identify $$CI^+ = \prod_{n \geq 0} \widetilde{CI}[4n] \cong \widetilde{CI} \otimes_R R\llbracket U^*\rrbracket.$$ The action of $R[U] = C^-_{\Lambda(u)}$ is such that $U \cdot (U^*)^k = (U^*)^{k-1}$; we must restrict to this subalgebra as the action of the entire algebra does not extend to the completion. We denote $DCI^+ = C^+_{\Lambda(u)}(DCI),$ equipped with its action of $R[U^*]$.

By convention, if $A$ is an operator from one summand of $DCI$ to another, then we use $A$ to also refer to the corresponding map between summands of $DCI^+$ with 
$$A(x \otimes (U^*)^n) = (-1)^n Ax \otimes (U^*)^n.$$ 
Then the differential of $$DCI^+ \cong \left(C_*^\text{irr} \oplus C_*^\text{irr}[3] \oplus C_*^{U(1)} \oplus C_*^{U(1)}[2] \oplus C_*^\theta\right) \widehat{\otimes}_R R[U^*]$$ is given by the following matrix.
 
$$\partial_{DCI}^+ := \begin{pmatrix}\partial_1 & 0 & 0 & 0 & 0\\
U_{\text{alg}}^+ + U_{\text{Fl}} & -\partial_1 & V_4 & V_2 & D_2 \\ 
V_1 & 0 & 0 & 0 & 0\\ 
V_3 & 0 & 0 & 0 & 0\\ 
D_1 & 0 & 0 & 0 & 0\end{pmatrix}.$$

Now we apply the following analogue to \cite[Lemma~5]{SS-Invol} to our very similar situation. First, observe that the map $U_\text{Fl} + U_{\text{alg}}$ defines an isomorphism $$U^* C_*^{\text{irr}}\llbracket U^*\rrbracket \to C_*^{\text{irr}}[3] \llbracket U^*\rrbracket;$$ here we have used the completeness in $U^*$. The inverse is explicitly given by 

$$\sum_{n \geq 0} c_n (U^*)^n \mapsto \sum_{n \geq 0} \left(\sum_{j=0}^n (-1)^{j+1} U_{\text{Fl}}^{n-j} c_j\right)(U^*)^{n+1}.$$

We write $P^+$ for this inverse.

Consider the projection map $\epsilon_1: DCI^+ \to C_*^{\text{irr}}\llbracket U^*\rrbracket [3];$ we let $$\epsilon = (\epsilon_1, \epsilon_1 \partial^+): DCI^+ \to C_*^{\text{irr}}\llbracket U^*\rrbracket[3] \oplus C_*^{\text{irr}}\llbracket U^*\rrbracket[2].$$ The codomain is given the differential $\begin{pmatrix}0 & 1 \\ 0 & 0 \end{pmatrix}$, and $\epsilon$ is a chain map; in fact, it is surjective, as $(x,y)$ is in the image of $(P^+(y - \partial_1 x), x) \in DCI^+$. As the codomain is acyclic, the inclusion $\text{ker}(\varepsilon) \to DCI^+$ is a quasi-isomorphism. 

Write $$\overline{DCI}^+ = C_*^{\text{irr}} \oplus \left(C_*^{U(1)} \oplus C_*^{U(1)}[2] \oplus C_*^\theta\right) \otimes_R R\llbracket U^*\rrbracket$$ as a graded $R$-module. There is also a projection $\pi: DCI^+ \to \overline{DCI}^+$, and the composite $\pi i: \ker \epsilon \to \overline{DCI}^+$ is an isomorphism of $R$-modules. We wish to describe the differential (and $U$-action) on $\overline{DCI}^+$, pulling back via this isomorphism. The differential is written $(\pi i) \partial^+ (\pi i)^{-1}$, and similarly the $U$-map (whose degree is $-4$) is $(\pi i) U (\pi i)^{-1}$. 

Explicitly, $(\pi i)^{-1}x$ is the unique element $y \in DCI^+$ so that neither $y$ nor $\partial^+ y$ have components in $C_*^{\text{irr}}\llbracket U^*\rrbracket [3]$ with $\pi(y) = x$. If $$x = (x_0, x_1, x_2, x_3) \in C_*^{\text{irr}} \oplus \left(C_*^{U(1)} \oplus C_*^{U(1)}[2] \oplus C_*^\theta\right) \otimes_R R\llbracket U^*\rrbracket,$$ then 
\begin{align*}(i\pi)^{-1}(x) = (x_0+P^+(U_{\text{Fl}}x_0 + V_4 x_1 &+ V_2 x_2 + D_2 x_3), 0, x_1, x_2, x_3) \in \\ &\left(C_*^{\text{irr}} \oplus C_*^{\text{irr}}[3] \oplus C_*^{U(1)} \oplus C_*^{U(1)}[2] \oplus C_*^\theta\right)R\llbracket U^*\rrbracket.\end{align*} 
Applying $\partial^+$, projecting onto $\overline{DCI}^+$, and using that $\partial_1 P^+x$ is in $\text{ker}(\epsilon)$, we find that the desired matrix is

$$\overline\partial^+_{DCI} := \begin{pmatrix}\partial_1 & 0 & 0 & 0\\ 
V_1 + V_1P^+U_{\text{Fl}} & V_1 P^+ V_4 & V_1 P^+ V_2 & V_1 P^+ D_2\\ 
V_3 + V_3P^+U_{\text{Fl}} & V_3 P^+ V_4 & V_3 P^+ V_2 & V_3 P^+ D_2\\ 
D_1 + D_1P^+U_{\text{Fl}} & D_1 P^+ V_4 & D_1 P^+ V_2 & D_1 P^+ D_2\end{pmatrix}$$

A remarkable number of these terms are zero for degree reasons. Recall from Proposition \ref{gr-red} that the relative grading of any pair of reducible orbits is even. Now, $D_i$ and $V_i$ have degree of the same parity as $i$, and all terms in the infinite sum $P^+x$ have the same grading as $x$ modulo $4$; therefore the entries not in the left column are all of odd degree, and hence identically zero! Therefore, we have 

$$\overline{\partial}_{DCI}^+ = \begin{pmatrix}\partial_1 & 0 & 0 & 0 \\ 
V_1 + V_1P^+U_{\text{Fl}} & 0 & 0 & 0 \\ 
V_3 + V_3P^+U_{\text{Fl}} & 0 & 0 & 0 \\ 
D_1 + D_1P^+U_{\text{Fl}} & 0 & 0 & 0\end{pmatrix}$$

Similarly, the action of $U \in H^{-4}(BSO(3);R)$ was previously contraction against $U^*$; using $U \frown$ to denote contraction by $U$, the action $\overline U$ on $\overline{DCI}^+$ is given by $$\begin{pmatrix}-U_{\text{Fl}} & -V_4 & -V_2 & -D_2 \\ 0 & U \frown & 0 & 0\\ 0 & 0 & U \frown & 0 \\ 0 & 0 & 0 & U \frown \end{pmatrix}.$$

In the case that $Y$ is an integer homology sphere and $\mathbb Q \subset R$, up to a change of basis (just a scaling of each coordinate) the chain complex $\overline{DCI}^+(Y;R)$ is identical as a $U$-module to Donaldson's $\overline{\overline{CF}}(Y;R)$.

The full completion involved in defining $CI^-$ does not require passing to power series (now the interesting part is the colimit instead of the limit, but $CI^{unr}_{\le p}$ stabilizes in each degree, so the colimit is degreewise the same). As a result, we may write $$DCI^- = \left(C_*^\text{irr} \oplus C_*^\text{irr}[3] \oplus C_*^{U(1)} \oplus C_*^{U(1)}[2] \oplus C_*^\theta\right)[U],$$ with differential 
$$\partial_{DCI}^- := \begin{pmatrix}\partial_1 & 0 & 0 & 0 & 0\\
U^-_{\text{alg}} + U_{\text{Fl}} & -\partial_1 & V_4 & V_2 & D_2 \\ 
V_1 & 0 & 0 & 0 & 0\\ 
V_3 & 0 & 0 & 0 & 0\\ 
D_1 & 0 & 0 & 0 & 0\end{pmatrix}.$$ 

Unlike in the case of $DCI^+$, the operators $\partial_1, V_i, D_i$ do not pick up an extra sign when acting on $C[U]$. For $x \in C_*^{\text{irr}}$, we define the operator $$U^-_\text{alg}(x \otimes U^k) = (-1)^{|x| + k} x \otimes U^{k+1} \in C_*^{\text{irr}}[3][U];$$ this is the only sign in the above differential hidden in the notation. 

Now $$U_{\text{Fl}} + U_{\text{alg}}^-: C_*^{\text{irr}}[U] \to \left(C_*^{\text{irr}}[3]\right)[U]\big/(1 \cdot C_*^{\text{irr}}[3])$$ is an isomorphism; it is crucial here that we are working with a polynomial ring and not a power series ring. The inverse, $P^-$, is given on a basis element by $$P^-(x \otimes U^{n+1}) = \sum_{i=0}^n (-1)^{(i+1)(|x|+n)+i(i+1)/2} U_{\text{Fl}}^i x \otimes U^{n-i}.$$

Now we follow \cite{SS-Invol} more closely, and instead of taking the kernel of a map to an acyclic complex, we quotient by an acyclic subcomplex. The subcomplex $Z$ is spanned by $C_*^\text{irr}\llbracket U\rrbracket$ and its image under $\partial_{DCI}^-$, which is an injective map. The quotient $DCI^-/Z$ remains a $U$-module. If we set $$\overline{DCI}^- = C_*^{\text{irr}}[3] \oplus \left(C_*^{U(1)} \oplus C_*^{U(1)}[2] \oplus C_*^\theta\right) \otimes_R R[U]$$ we have a natural inclusion map $i: \overline{DCI}^- \hookrightarrow DCI^-$ so that the composite $$\pi i: \overline{DCI}^- \to DCI^-/Z$$ is an isomorphism. As above, we may compute the induced differential as

$$\overline{\partial}^-_{DCI} = \begin{pmatrix}-\partial_1 & V_4 - U_\text{Fl} \left(P^- V_4\right)_{U = 0} & V_2 - U_\text{Fl} \left(P^- V_2\right)_{U = 0} & D_2 - U_\text{Fl} \left(P^- D_2\right)_{U = 0}\\ 
0 & -V_1P^-V_4 & -V_1P^-V_2 & -V_1P^-D_2\\ 
0 & -V_3P^-V_4 & -V_3P^-V_2 & -V_3P^-D_2\\ 
0 & -D_1P^-V_4 & -D_1P^-V_2 & -D_1P^-D_2\end{pmatrix}.$$ 
Here the term $V_4$ in the topmost row, and similarly with $V_2$ and $D_2$ in that same row, denotes the composite $$C_*^{U(1)}[U] \xrightarrow{U = 0} C_*^{U(1)} \xrightarrow{V_4} C_*^{\text{irr}}[3];$$ in the bottom right $3 \times 3$ block matrix, there is no projection component to the terms $V_4, V_2$, or $D_2$. When we write for instance $(P^- V_4)_{U = 0}$, we mean the composite $$C_*^{U(1)}[U] \xrightarrow{P^- V_4} C_*^{\text{irr}}[3][U] \xrightarrow{U = 0} C_*^{\text{irr}}[3]:$$ that is, apply $P^-V_4$ then project to the constant term. 

Again, most of these terms are zero for degree reasons, and we may in fact write 
$$\overline{\partial}^-_{DCI} = \begin{pmatrix}-\partial_1 & V_4 - U_\text{Fl} \left(P^- V_4\right)_{U = 0} & V_2 - U_\text{Fl} \left(P^- V_2\right)_{U = 0} & D_2 -
U_\text{Fl} \left(P^- D_2\right)_{U = 0}\\ 0 & 0 & 0 & 0 \\ 0 & 0 & 0 & 0\\ 0 & 0 & 0 & 0\end{pmatrix}.$$

Finally, the action of $U$ is given on $\overline{DCI}^-$ by 

$$\begin{pmatrix}(-1)^{|x|+1}U_{\text{Fl}} & 0 & 0 & 0 \\ (-1)^{|x|+1} V_1 & U \smile & 0 & 0\\ (-1)^{|x|+1} V_3 & 0 & U\smile & 0 \\ (-1)^{|x|+1} D_1 & 0 & 0 & U \smile \end{pmatrix}.$$

Again, when $Y$ is an integer homology sphere, $\mathbb Q \subset R$, and the holonomy maps lift to $SU(2)$, the chain complex $\overline{DCI}^+(Y;R)$ is identical as a $U$-module to Donaldson's $\uline{\uline{CF}}(Y;R)$ up to scaling of basis.

When $(Y, E)$ is equipped with an admissible bundle, there are no reducible connections; in this case, the equivalence of $CI^+(Y,E;R)$ with $\overline{DCI}^+(Y,E;R)$ as an $R[U]$-module (and the same for the minus flavor) immediately gives is the following.

\begin{corollary}\label{admiss-U}When $\frac 12 \in R$, the equivalence $I^\pm(Y,E) \cong I(Y,E)$ for an admissible bundle $(Y,E)$ takes the $U$-map to the $U$-map up to sign.
\end{corollary}

\section{Instanton Tate homology}\label{sec:ex-tate}
We have the following theorems for Tate homology. Fix a ground ring $R$ with $\frac 12 \in R$; then we have a canonical isomorphism $H^-_{SO(3)}(R) \cong R[U]$ and a quasi-isomorphism $C_*(SO(3);R) \simeq \Lambda(u) := \Lambda$, where $\Lambda(u)$ denotes the exterior algebra on a degree 3 generator $u$. 

\begin{proposition}\label{InstTateLoc}The action of $U^*$ is an isomorphism on $I^\infty(Y,E;R)$. In fact, the map $I^-(Y,E;R) \to I^\infty(Y,E;R)$ may be identified with the localization $$I^-(Y,E;R) \to I^-(Y,E;R) \otimes_{R[U]} R[U, U^{-1}\rrbracket;$$ even more explicitly, we may identify the map $CI^-(Y,E;R) \to CI^\infty(Y,E;R)$ on the chain level using the reduced Donaldson model as $$\overline{DCI}^-(Y,E;R) \to \overline{DCI}^-(Y,E;R) \otimes_{R[U]} R[U, U^{-1}\rrbracket.$$
\end{proposition}
\begin{proof}We follow much the same lines as in Proposition \ref{TateLoc} in the $\mathbb Z$-graded case, and work with the Donaldson model $\overline{DCI}^-$. First, to see that the action of $U$ on $I^\infty$ is invertible, we look at the spectral sequence corresponding to the index filtration of $CI^\infty$. The $E^1$ page is identified with a direct sum of copies of 
$$H^\infty_{\Lambda}(\Lambda), \; H^\infty_{\Lambda}(R \oplus R[2]), \;\text{and } H^\infty_{SO(3)}(R);$$ 
the first corresponds to irreducible orbits, the second corresponds to $SO(2)$-reducibles and is just a sum of two copies of $H^\infty_\Lambda(R)$, and the last corresponds to full reducibles. In all cases, the action of $U$ is an isomorphism on each factor: that the first group is zero is one of the axioms of Tate homology, and that the second two have invertible $U$ action is precisely the calculation of Lemma \ref{PeriodicComputation}. Because the index filtration is complete, a map which is an isomorphism on the $E^1$ page is an isomorphism on homology, and so the action of $U$ is invertible in $I^\infty(Y,E;R)$.

Given a map $f: M \to M$ of degree $k$ of a dg-module, there is a definition given before Proposition \ref{TateLoc} of a chain complex $M[f^{-1}]$; it is the mapping cone of $M[t] \xrightarrow{1 - tf} M[t]$, where $|t| = -k$. This is defined so that 
$$H(M[f^{-1}]) = H(M)[f^{-1}],$$ the latter notation meaning the usual module-theoretic sense of inverting an element or map.

If one applies this to define $CI^\infty(Y,E;R)[U^{-1}]$, one finds that the index filtration is \emph{no longer complete}, even though the index filtration is complete on $CI^\infty$. One may then pass to the completion of $CI^\infty(Y,E;R)[U^{-1}]$; this is the filtered complex that the spectral sequence actually computes the homology of (and the $E^1$ pages of the spectral sequence for a filtered complex and its completion are the same). Precisely, if $F_s M$ is a filtered complex, the completion is $$\widehat M = \operatorname{lim}_{p \to \infty} M/F_{-p} M.$$ We write this completion as $CI^\infty(Y,E;R)\widehat{[U^{-1}]}$. There is a natural map $$CI^\infty(Y,E;R) \to CI^\infty(Y,E;R)\widehat{[U^{-1}]}$$ which we know to be an isomorphism on the $E^1$ page, because the action of $U^*$ is an isomorphism on the $E^1$ page of $CI^\infty(Y,E;R)$.

Because completion is natural, we have maps $$CI^-(Y,E;R)\widehat{[U^{-1}]} \to CI^\infty(Y,E;R)\widehat{[U^{-1}]};$$ again, the computations of Lemma \ref{PeriodicComputation} imply that this is an isomorphism on the $E^1$ page, and because these filtrations are complete, the same is true at the level of homology.

To conclude we need to find a quasi-isomorphism 
$$CI^-(Y,E;R)\widehat{[U^{-1}]} \simeq \overline{DCI}^-(Y,E;R) \otimes_{R[U]} R[U, U^{-1}\rrbracket.$$ 
This follows again via naturality of completion: there is a canonical map $$\overline{DCI}^-(Y,E;R)[U^{-1}] \to \overline{DCI}^-(Y,E;R) \otimes_{R[U]} R[U, U^{-1}]$$ of filtered complexes for which the filtration is \textit{not} complete (remember that the first term is defined to be the mapping cone of a certain map on $DCI^-[t]$). The same argument as before shows this map is an isomorphism on the $E^1$ page, so the corresponding map on completions 
$$\overline{DCI}^-(Y,E;R)\widehat{[U^{-1}]} \to \overline{DCI}^-(Y,E;R) \widehat{\otimes}_{R[U]} R[U, U^{-1}]$$ is a quasi-isomorphism. We conclude by observing that we may explicitly identify the completion as $$\overline{DCI}^-(Y,E;R) \widehat{\otimes}_{R[U]} R[U, U^{-1}] = \overline{DCI}^-(Y,E;R) \otimes_{R[U]} R[U, U^{-1}\rrbracket,$$ using the fact that $\overline{DCI}^-$ is a finite direct sum of copies of $R$ and $R[U]$.
\end{proof}

\begin{corollary}\label{tate-iso}Let $(Y,E)$ be a rational homology 3-sphere with $U(2)$-bundle and signature data, and write $c = c_1(E) \in H^2(Y;\mathbb Z)$. Then there is a natural $\mathbb Z/2$-graded isomorphism $$I^\infty(Y,E) \cong R[U^{1/2}, U^{-1/2}\rrbracket \otimes_{R[\mathbb Z/2]} R[H^2 Y].$$ Here the first term is Laurent series in $U^{-1/2}$, and the action of $\mathbb Z/2$ is $(-1) \cdot U^{k/2} = (-1)^k U^{k/2}$, while the action on $R[H^2 Y]$ is given by $(-1) \cdot e^{[x]} = e^{[c - x]}$.
\end{corollary}

\begin{proof}This is immediate from the computation of the differential on $\overline{DCI}^-$ in the previous section: there is no component of the differential going from the reducible part of $\overline{DCI}^-$ back to itself. At the same time, the reducible part is the only thing that survives to $\overline{DCI}^- \otimes_{R[U]} R[U, U^{-1}\rrbracket.$ Therefore the differential on this is identically zero, and the graded group $$\overline{DCI}^- \otimes_{R[U]} R[U, U^{-1}\rrbracket$$ itself calculates the Tate homology.

What we see, then, is that $I^\infty(Y,E;R)$ is a direct sum of copies of $R[U, U^{-1}\rrbracket$, one for each full reducible and two for each $SO(2)$-reducible (one of the copies shifted up by $2$). Now, the content of Proposition \ref{ModuliPackage} (2) is that pairs $\{z_1, z_2\}$ with $z_1 + z_2 = c$ and $z_1 \neq z_2$ are in bijection with the $SO(2)$-reducible critical orbits, and that the pairs with $z_1 = z_2$ are in bijection with the full reducibles.

We choose the perturbation to be small enough so that this enumeration of reducible flat connections holds.

Contributing two copies of $R[U, U^{-1}\rrbracket$, one shifted up by degree $2$, is equivalent to contributing a single copy of 
$$R[U^{1/2}, U^{-1/2}\rrbracket =R[U, U^{-1}\rrbracket \oplus U^{-1/2} R[U, U^{-1}\rrbracket,$$ 
recalling here that $|U^{-1/2}| = 2$.

The description in the corollary arises from the observation that we may give this a more succinct description, only in terms of the cohomology group $H^2(Y)$ and the class $c_1 E$. The point is that $$R[U^{1/2}, U^{-1/2}\rrbracket \otimes_{R[\mathbb Z/2]} R[H^2 Y]$$ may be described explicitly: for a pair $\{z_1, z_2\}$ with $z_1 + z_2 = c$ and $z_1 \neq z_2$, the involution swaps the two corresponding copies of $R[U^{1/2}, U^{-1/2}\rrbracket$, and therefore in the quotient we have identified these two towers; so we have a copy of $R[U^{1/2}, U^{-1/2}\rrbracket$ for each such pair. If, on the other hand, $z_1 = z_2$, the action identifies 
$$U^{k/2} \otimes e^{[z_1]} \sim (-1)^k U^{k/2} \otimes e^{[c - z_1]} = (-1)^k U^{k/2} \otimes e^{[z_1]}.$$ 
Because $2$ is invertible in $R$, we see that this kills precisely the terms of the form $U^{k/2}$; so such a pair contributes a copy of $R[U, U^{-1}\rrbracket.$

This is precisely the description given above as $$\overline{DCI}^-(Y,E,\pi;R) \otimes_{R[U]} R[U, U^{-1}\rrbracket.\qedhere$$
\end{proof}

\section{Examples of the $I^\bullet$ and the index spectral sequence}\label{examples}
\begin{example}Let $Y = S^3$ equipped with the trivial bundle. Then $\widetilde{CI}(Y;R) = R$ concentrated in degree $0$. Therefore, $CI^+(S^3;R)$ is given as the completion $\hat C_*(BSO(3);R)$. Recall that this means that for $i \in \mathbb Z/8$, we have $$\hat C_i(BSO(3);R) = \prod_{\substack{\mathbb Z \ni j \geq 0\\ j \equiv i \bmod 8}} C_j(BSO(3);R).$$ 
Therefore, we have
$$I^+_*(S^3;R) \cong \hat H_*(BSO(3);R).$$ 
Similarly, 
$$I^-_*(S^3;R) \cong H^{-*}_{\text{fin}}(BSO(3);R),$$ 
meaning that we take cohomology classes which correspond to functionals of finite support; that is, we do not complete in the direction of negative degrees:
$$H^{-i}_{\text{fin}}(BSO(3);R) = \bigoplus_{\substack{\mathbb Z \ni j \geq 0\\ j \equiv i \bmod 8}} H^{-j}(BSO(3);R).$$
Finally, as graded $R$-modules we have
$$I^\infty(S^3;R) = \hat H_*(BSO(3);R)[3] \oplus H^{-*}_{\text{fin}}(BSO(3);R).$$ 
The calculations of these groups, and their module structures, is given in Example \ref{SO3calc} for $\frac 12 \in R$ and $R = \mathbb Z/2$.

We are presented with a dichotomy. When $\frac 12 \in R$, we see that

\begin{align*}
I^+(S^3;R) &\cong R\llbracket U^{-1}\rrbracket, \\
I^-(S^3;R) &\cong R[U]\\
I^\infty(S^3;R) &\cong R[U, U^{-1}\rrbracket,
\end{align*} all as $R[U]$-modules; this periodicity calculation in Tate homology is the content of Lemma \ref{PeriodicComputation}, and the other two are straightforward calculations using the definitions of bar and cobar constructions.

However, when $2$ is not invertible, we do not have such a periodicity. The most dramatic case is when $R = \mathbb Z/2$; in that case, $$I^+(S^3;R) = (\mathbb Z/2)\llbracket w_2^*, w_3^*\rrbracket,$$ and $$I^-(S^3;R) = (\mathbb Z/2)[w_2, w_3].$$ Then $I^\infty(S^3;R)$ looks like a `bi-infinite cone': the Tate homology $H^\infty_{SO(3)}(\mathbb Z/2)$ has rank growing roughly linearly in the degree $|k|$. In fact, the action of $H^{-*}_{\text{fin}}(BSO(3);\mathbb Z/2)$ on $$\hat H_*(BSO(3);\mathbb Z/2) \subset H^\infty_{SO(3)}(\mathbb Z/2)$$ is nilpotent on each element, as opposed to having an element that induces a periodicity isomorphism $I^\infty(S^3) \to I^\infty(S^3)$!

This suggests that in the case of $R = \mathbb Z/2$, it would be more appropriate to study $\widetilde{CI}(Y;\mathbb Z/2)$ as a dg-module over some specific exterior algebra $\Lambda(u_1)$ or $\Lambda(u_2)$ in $C_*(SO(3);\mathbb F_2)$. Mostly we will content ourselves with the case $\frac 12 \in R$.
\end{example}

\begin{example}If $p > q$ are coprime integers (not necessarily prime), let $L(p,q)$ be the lens space, which by our convention is the quotient of $S^3$ by the $\mathbb Z/p$ action generated by $[1] \cdot (z, w) = (e^{2\pi/p}z, e^{2q\pi/p} w).$

First, we work with the trivial $U(2)$-bundle; the set of critical orbits now correspond to 
$$\text{Hom}(\pi_1, SU(2))/\!\sim \; = \text{Hom}(\mathbb Z/p, SU(2))/\!\sim,$$ 
where the equivalence elation is conjugacy in $SU(2)$. Because $\mathbb Z/p$ is abelian, and simultaneously commuting matrices are simultaneously diagonalizable, this is the same as $\text{Hom}(\mathbb Z/p, S^1)/\text{conj}$; here $\text{conj}$ is complex conjugation. Identify $\text{Hom}(\mathbb Z/p, S^1)$ with the $p$th roots of unity, and hence with $\mathbb Z/p$ again after fixing the generator $e^{2\pi i/p}$; thus reducibles correspond to $(\mathbb Z/p)/\pm 1 = [0, p/2]$, and full reducibles --- those fixed by the conjugation action, which is $\pm 1$ on $\mathbb Z/p$ --- correspond to $[0]$ for any $p$ and $[p/2]$ when $p$ is even. The $SO(2)$-reducibles correspond to $0 < i < p/2$. Finally, $0$ corresponds to the trivial connection.

Austin calculates in \cite{LensSpaceModuli} the expected dimension of the different components of the moduli space of unframed instantons (before quotienting by the translation action) between two flat connections on $L(p,q)$, and in particular the expected dimension of $\mathcal M(L(p,q), 0, i)$; the expected dimension of our $\overline{\mathcal M}(L(p,q), \theta, \alpha_i)$ is $\dim SO(3) - \dim \mathbb R = 2$ dimensions larger; the index $\text{gr}(\alpha_i) = \text{gr}(\theta, \alpha_i)$ differs from this dimension by subtracting $\dim \alpha_i$ and adding $1$. Thus we should add either $3$ or $1$ to Austin's result; we add $3$ in the case that $i$ is either $[0]$ or $[p/2]$, and we add $1$ otherwise. Set $\varepsilon(i) = 1$ if $0 < i < p/2$ and $\varepsilon(0) = \varepsilon(p/2) = 0$, and write $0 < q' < p$ for the unique integer with $qq' = ap + 1$.

Then the grading function $\delta(p,q,i) = \text{gr}(L(p,q), i) \in \mathbb Z/8$ is given by $$\delta(p,q,i) = \frac{8i^2q'}{p} - \varepsilon(i) + \frac{2}{p} \sum_{j=1}^{p-1} \cot\left(\frac{j \pi}{p}\right)\cot\left(\frac{jq \pi}{p}\right)\sin^2\left(\frac{2ij\pi}{p}\right) \pmod{8}.$$

Sasahira \cite[Corollary~5.3]{SasahiraLens} has given formulas for the $\delta(p,q,i)$ in terms of counting solutions to congruences.

A few observations about this complicated-looking function are in order: \begin{itemize}
\item We have $\delta(p,q,0) = 0$, as we should. \\
\item As Austin observes, $\delta(p,q,i)$ is even for all $i$.\\
\item If $p$ is even, $\delta(p,q,p/2) = 4q \cdot \frac{p}{2}$; because $q'$ must be odd, $\delta(4k, q, 2k) = 0$ and $\delta(4k+2, q, 2k+1) = 4$.\\
\item Because the sum only depends on $q$'s value modulo $p$, and $\cot$ is odd, when $0 < i < p/2$ we have $\delta(p,q,i) = - \delta(p,p-q,i) - 2$. The factor of $-2$ comes from $-2\varepsilon(i)$. It is trivially true by the above calculations that $\delta(p,q,i) = -\delta(p,p-q,i)$ when $i = 0$ or $p/2$.\\
\item If $qq' = ap + 1$, then $\delta(p, q, qi) = \delta(p, q', i).$ Here we are considering $d$ as a function from the integers, but notice that $\delta(p,q,i) = -\delta(p,q,i)$ and $\delta(p,q,i+p) = \delta(p,q,i)$. Then the claimed equality follows because the summation in $\delta(p,q',i)$ is just the summation in $\delta(p,q,qi)$, but with index $jq'$ instead of $j$; as $\times q': \left((\mathbb Z/p) \setminus 0\right) \to \left((\mathbb Z/p) \setminus 0\right)$ is a bijection, and $j$'s value modulo $p$ is all that's relevant, the sum is therefore the same.
\end{itemize}

It seems plausible that $\delta(p,q,ki) = \delta(p,q',i)$ for all $i$ iff $k = 1$ and $q' = q$ or $k = q$ and $q' = q^{-1} \pmod p$; calculation shows that this is true at least for $p \leq 15$. This would be consistent with the classification of lens spaces up to oriented homeomorphism.\\

We now investigate the differential. Because $\overline{\mathcal M}(\alpha, \beta)$ is $SO(3)$-free, if $\alpha$ or $\beta$ is fully reducible, the component $d_M: C_*(\alpha; R) \to C_*(\beta; R)$ is zero for degree reasons. Similar considerations apply unless $\alpha, \beta$ are $U(1)$-reducible and $\dim\overline{\mathcal M}(\alpha, \beta) = 3$, as this is necessarily an odd number at least $3$. In the three-dimensional case, we may homotope $e_+: \overline{\mathcal M}(\alpha, \beta) \to \beta \cong S^2$ so that the joint endpoint map \[\overline{\mathcal M}(\alpha, \beta) \to \alpha \times \beta \cong S^2 \times S^2\] has image in the diagonal; for this modified endpoint map, the fiber-product map is zero on the chain level. Extending this homotopy coherently to higher-dimensional moduli spaces and applying Lemma \ref{EndpointHomotopy}, we obtain an equivariant homotopy equivalent complex for which $d_M$ is identically zero.

Thus, as a dg $C^\text{gm}_*(SO(3))$-module, and with $L(p,q)$ equipped with the trivial $SO(3)$-bundle, we see that we have a quasi-isomorphism
$$\widetilde{CI}(L(p,q); R) \simeq \begin{cases}
R \oplus_{i=1}^{\lfloor p/2 \rfloor} C_*(S^2;R)[\delta(p,q,i)] & p \text{ odd} \\
R \oplus_{i=1}^{p/2 - 1} C_*(S^2;R)[\delta(p,q,i)] \oplus R[\delta(p,q,p/2)] & \text{else}
\end{cases}$$

Thus to write down $I^\bullet$ (where $\bullet = +, -, \infty$) we use Theorem \ref{orbitcalc}, which says that the equivariant homology groups of an orbit are given by the homology groups of the stabilizer (with a degree shift) as long as the stablizer is connected. So to take $I^\bullet$, we replace every appearance of $C_*(S^2;R)$ with $H^\bullet_{SO(2)}(R)$, with a degree shift of $2$ if $\bullet \in \{-,\infty\}$, and every appearance of $R[d]$ with $H^\bullet_{SO(3)}(R)[d]$. When $\frac 12 \in R$, we have isomorphisms as $R[U] = H^-_{SO(3)}(R)$-modules
\begin{align*}H_{SO(2)}^+(R) &\cong R[U^{1/2},U^{-1/2}\rrbracket\big/U^{1/2}R[U^{1/2}]\\
H_{SO(2)}^-(R) &\cong R[U^{1/2}] \\
H_{SO(2)}^\infty(R) &\cong R[U^{1/2},U^{-1/2}\rrbracket
\end{align*}

If $p$ is even, $L(p,q)$ carries a unique nontrivial $SO(3)$-bundle; choose a lift of this to a $U(2)$-bundle. We can identify the classes now with $\{a + \xi\}/\pm 1$, where $\xi \in \mathbb Z/p$ is a choice of odd number (the choice of $U(2)$ lift). This, then, may be identified with the odd numbers in $\mathbb Z/2p$ modulo $\pm 1$; no points are fixed, so all are $SO(2)$-reducibles. So we label the reducibles by $i = 1, 3, \cdots, p-1$; there are $p/2$ of them. Their relative grading is also given in \cite{LensSpaceModuli}, now as
\begin{align*}
\delta(p,q,i,i') &= \frac{2(i^2 - (i')^2)q'}{p} - 2  \\
+\frac 2p \sum_{i=1}^{p-1} &\cot\left(\frac{j \pi}{p}\right)\cot\left(\frac{jq \pi}{p}\right)\left(\sin^2\left(\frac{i'j\pi}{p}\right) -  \sin^2\left(\frac{ij\pi}{p}\right)\right) \pmod 8.
\end{align*}
Again, this is always even.

The result is here independent of the choice of $U(2)$-lift of the underlying $SO(3)$-bundle because there are no differentials, and thus we may ignore the orientation of the moduli spaces; in fact, we expect this in general (up to noncanonical isomorphism, the noncanonicity coming from sign choices).
\end{example}

\begin{example} Suppose $Y$ is an integer homology sphere and $E$ is the trivial bundle. We may use the results of Chapter \ref{secDCI} to determine the index spectral sequences explicitly.

We write $CI(Y)$ for the chain complex $(C_*^{\text{irr}}, \partial_1)$, Floer's original chain complex for integer homology spheres. Donaldson introduced in \cite[Section~7.1]{Don} the complex $\overline{CI}(Y)$, given as $C_*^{\text{irr}} \oplus R$ with differential $$\overline \partial = \begin{pmatrix}\partial_1 & 0 \\ D_1 & 0\end{pmatrix}.$$

Corollary \ref{DCI} shows that there is an equivariant filtered homotopy equivalence $\widetilde{CI}(Y) \simeq DCI(Y)$, and so to investigate the index spectral sequence computing $\tilde I$, we may do the same for the index spectral sequence on the finite-dimensional complex $DCI(Y)$. 

We consider $DCI(Y)$; this is, as an $R$-module, $CI(Y) \oplus CI(Y)[3] \oplus R$, with differential $$\begin{pmatrix}\partial_1 & 0 & 0 \\ U_\text{Fl} & -\partial_1 & D_2 \\ D_1 & 0 & 0\end{pmatrix}.$$ 

The differential splits as $d_1 + d_4$ into pieces which decrease the filtration the corresponding amount; 
$$d_1 = \begin{pmatrix}\partial_1 & 0 & 0 \\ 
0 & -\partial_1 & 0\\ 
D_1 & 0 & 0\end{pmatrix} \text{ and } d_4 =  \begin{pmatrix}0 & 0 & 0 \\ 
U_\text{Fl} & 0 & D_2 \\ 
0 & 0 & 0\end{pmatrix}.$$ 
Then the $E^2$ page of the index spectral sequence for $DCI(Y)$ is $$(DCI(Y), d_1) = \overline{CI}(Y) \oplus CI(Y)[3] \oplus R,$$ and so the $E^3$ page is $\overline{I}(Y) \oplus I(Y)[3]$. The matrix $d_4$ defines a chain map $$(U, D_2): \overline{CI}(Y) \to CI(Y)[3];$$ writing the induced map $f: \overline{I}(Y) \to I(Y)$, we see that the $E^5$ page is $$(\overline{I}(Y) \oplus I(Y)[3], f).$$ Therefore, the $E^6$ page (and all successive pages, for degree reasons) is $\text{ker}(f) \oplus \text{coker}(f)$.

If $R$ is a field, there is nothing left to say; there are no extension problems to resolve, and $\widetilde I(Y) \cong \text{ker}(f) \oplus \text{coker}(f)$.

Now consider the $CI^+$ spectral sequence by passing through $\overline{DCI}^+$. The complex $\overline{DCI}^+$ is, as a graded $R$-module, given by $CI(Y) \oplus R\llbracket U^*\rrbracket$, where $|U^*| = 4$, and the differential is given by 
$$\begin{pmatrix}\partial_1 & 0 \\ D_1 + D_1P^+U_\text{Fl} & 0\end{pmatrix},$$ 
where recall that by definition when $x \in CI(Y)$, we have $$(D_1 + D_1 P^+U_{\text{Fl}})x = \sum_{i \geq 0} (-1)^i D_1 U_{\text{Fl}}^{i+1} x (U^*)^{i+1}.$$

We decompose this differential as $d_1 + d_5 + d_9 + \cdots$, where $d_1 = \begin{pmatrix}\partial_1 & 0\\ D_1 & 0\end{pmatrix}$ and
$$d_{4i+1}x = \begin{pmatrix}0 & 0\\ (-1)^{i+1} D_1 U_{\text{Fl}}^ix \otimes (U^*)^i & 0 \end{pmatrix}.$$
Now this is a multicomplex in the sense of Wall, and we may identify the $E^2$ page as $\overline{I} \oplus U^*R\llbracket U^*\rrbracket$.

The map $d_5$ defines a homomorphism $\overline{I} \to U^*R$; write $\overline{I}_2$ for its kernel. Inductively, there is a homomorphism $$d_{4k+1}: \overline I_k \to (U^*)^k R,$$ and we write its kernel as $\overline I_{k+1}$. We may identify the $q = 0$ line in $E^{4k-2}$ page with $\overline I_{k}$ (identifying the successive differentials with zig-zags of the differentials $d_{4k+1}$, and using that any compositions of the differentials are zero). In particular, we identify the $q=0$ line in the $E^\infty$ page with $$\overline I_\infty = \bigcap \text{ker}(D_1U_{\text{Fl}}^k) \subset \overline I.$$

A similar discussion presents itself for the $CI^-$ spectral sequence, passing through $\overline{DCI}^-$. As a graded $R$-module, $$\overline{DCI}^-(Y) = CI(Y)[3] \oplus R[U];$$ thinking of an element of $R[U]$ as a polynomial $f$, we write $f_k$ for the coefficient of $U^k$; this has differential 
$$\overline{\partial}^- = \begin{pmatrix}-\partial_1 & \sum_{k \geq 0} (-1)^{k(k-1)/2} U_{\text{Fl}}^k D_2 f_k \\0 & 0\end{pmatrix}.$$ 
The top-right term means that $$\overline{\partial}^-(U^k) = (-1)^{k(k-1)/2} U_{\text{Fl}}^k D_2(1).$$

As before, we may split this into pieces $d_1 + d_5 + d_9 + \cdots$. The $E^1$ page is given as $CI[3] \oplus R$ with differential $$\begin{pmatrix}-\partial_1 & D_2\\0 & 0\end{pmatrix},$$ which (shifted by 3 and with a sign change) Donaldson calls $\uline{CI}(Y)$. The homology, written $\uline I$, is the $q = 0$ page of this spectral sequence. The differential $d_5$ (which is, up to sign, $U_{\text{Fl}} D_2$) determines a homomorphism $R \to \uline I$; its cokernel might be written $\uline I^{(2)}$, and may be identified with the $E^6$ page. Inductively, the differential $d_{4k+1}$ defines a homomorphism $R \to \uline I^{(k)}$, and its cokernel is written $\uline I^{(k+1)}$, and may be identified with the $E^{4k+2}$ page.

If $\Lambda$ is the exterior algebra on a generator of degree $3$, the norm map is an isomorphism $$H^+_\Lambda(\Lambda)[3] \to H^-_\Lambda(\Lambda)$$ (as a corollary of item (4) of Theorem \ref{eq-package}), and is zero for $H^+_\Lambda(R)[3] \to H^-_\Lambda(R)$ for degree reasons: the former is supported in degrees no less than $3$, and the latter is supported in nonpositive degrees.

Therefore, the induced map of the norm map on the $E^1$ page of the $I^+$ and $I^-$ spectral sequences is the identity on $C_*^{\text{irr}}$, and zero on $R\llbracket U^*\rrbracket$. We thus identify the image $\overline I_{(\infty)} \to \uline I^{(\infty)}$ as the quotient group $\overline I_{(\infty)} / \langle U_{\text{Fl}}^i D_2, 1\rangle,$ the second term in the quotient generating the reducible piece of $\overline I_{(\infty)}$. This\footnote{Despite the choice of notation, these groups are not related to the Tate homology groups.} is Fr\o yshov's reduced Floer group $\widehat I$; the first hints of this were introduced in the $d^\alpha$ homomorphisms of \cite{Fr2} and used to introduce the $h$-invariant in \cite{Fr1}.
\end{example}

\begin{example}For the sake of compactness of notation, we write $R_{n}$ in this example to mean $R[n]$ (a copy of $R$ in degree $n$), and similarly for bigradings $R_{p,q}$. Let us see the above in practice in the explicit case of $\Sigma(2,3,5)$. In this case, we have a reducible generator of $\widetilde{CI}$ in degree $0$, as well as an irreducible generator in each pair of degrees $(1,4)$ and $(5,0)$, with the action of $u$ sending the lower degree generator to the generator 3 degrees larger. So as a $u$-module our complex is given as $$R_0 \oplus (R_1 \oplus R_4) \oplus (R_5 \oplus R_0) = R_0 \oplus (R_1 \oplus R_5) \otimes \Lambda(u).$$

For degree reasons, $\partial_1 = 0$, as is $D_2: \widetilde{CI}_0 \to \widetilde{CI}_{7}$. Through the beautiful and explicit calculations of moduli spaces in \cite{austin1995equivariant}, we may identify that the map $D_1: \widetilde{CI}_1 \to \widetilde{CI}_0$ is the identity map $R \to R$, while $U_{\text{Fl}}: \widetilde{CI}_1 \to \widetilde{CI}_0$ and $U_{\text{Fl}}: \widetilde{CI}_5 \to \widetilde{CI}_4$ are multiplication by one of $4$ and $3$, respectively. Note that in this case, the three-dimensional moduli spaces used in the definition of $U_{\text{Fl}}$ have no boundary, so that $\partial_{DCI}$ has no `higher terms'.

The $E^1$ page of the above spectral sequence is $R_{0,0} \oplus R_{1,0} \oplus R_{1,3} \oplus R_{5,0} \oplus R_{5,3}$, with multiplication-by-$u$ map taking us up $3$ vertical degrees. The differential $d_1$ comprises only the map $D_1$, and therefore is the identity $R_{1,0} \to R_{0,0}$.

The $E^2$ page is therefore given by $R_{1,3} \oplus R_{5,0} \oplus R_{5,3}$. There are no differentials until the $E^5$ page, where the next differential is $U_{\text{Fl}}: R_{5,0} \to R_{1,3}$, which we know to be multiplication-by-$8$. Since we are working over a ground ring with $\frac 12 \in R$, we see that $$E^\infty(\widetilde{CI}(\Sigma(2,3,5)) = R_{5,3};$$ because there can be no extension problems, $$\widetilde I(\Sigma(2,3,5)) = R_0,$$ with trivial multiplication-by-$u$ map.

For the $CI^+$ spectral sequence, instead we start with $E^1$ page equal to $$R_{0,0}\llbracket U^*\rrbracket \oplus R_{1,0} \oplus R_{5,0},$$ where $|U^*| = (0,4)$. Then our first differential $d_1$ is given as the identity map $R_{1,0} \to R_{0,0}$; the $E^5$ page is then identified with $R_{0,4}\llbracket U^*\rrbracket \oplus R_{5,0}$. We may identify the differential $d_5$ with the multiplication-by-$8$ map $R_{5,0} \to R_{0,4}$. The $E^6$ page, then, is simply $R_{0,8}\llbracket U^*\rrbracket$, and there are no further differentials. As this is concentrated in a single vertical line, there are no extension problems, and we see $$I^+(\Sigma(2,3,5);R) \cong (U^*)^2 R\llbracket U^*\rrbracket,$$ with $U$-action given by contraction against $U^*$. It is clear from this that the above `reduced group' is $$\widehat{I}(\Sigma(2,3,5)) = 0.$$

In the $CI^-$ spectral sequence, on the other hand, all differentials past the $E^2$ page are identically zero, as they factor through the map $D_2$. In fact, $$I^-(\Sigma(2,3,5);R) \cong \overline{DCI}^-(\Sigma(2,3,5)) = R_0[U] \oplus (R_4 \oplus R_0)$$ as a chain complex; by the formula for the $U$-map from Chapter \ref{secDCI}, which includes a term corresponding to $D_1$, we see that the $U$-action is given by $$U \cdot (U^k, x, y) = (U^{k+1} + D_1x, 0, 0).$$ Therefore we may write $I^-(\Sigma(2,3,5);R) \cong R_4[U] \oplus R_0$ as a $U$-module, where $|U| = -4$, identifying $1 \in R_0[U]$ with $D_1 x$, where $x \in R_4$ is a generator.
\end{example}

\section{Orientation reversal and equivariant cohomology}\label{sec:ex-orrev}
For a dg-$A$-module $M$ equipped with a periodic filtration whose associated graded complex is bounded, we defined equivariant $\bullet$-homology chain complexes $C^\bullet_G(M)$, with an action of the ring $C^-_G(R) = \text{Hom}_{\text{fin}}(C^+_G(R), R)$, the negative chains of finite support.

We define the corresponding $\bullet$-cohomology, written $C_\bullet^G(M)$, to be $$\text{Hom}(C^\bullet_G(M), R),$$ the literal dual of $C^\bullet_G$. We then define the instanton $\bullet$-cohomology, $I_\bullet^*$, to be $H_\bullet^{SO(3)}(\widetilde{CI}),$ where precisely by $SO(3)$ we mean the dg-algebra $C_*^{\text{gm}}(SO(3);R)$. We have the following duality theorem.

\begin{theorem}Let $\frac 12 \in R$. There is a canonical isomorphism of $\Lambda(u)$-modules $\widetilde I(Y, E) \cong \widetilde I^*(\overline Y, E)$, and furthermore canonical isomorphisms of\newline $H^{-*}(BSO(3);R) = R[U]$-modules
\begin{align*}I^+_*(Y,E) &\cong I_-^*(\overline Y,E)\\
I^-_*(Y,E) &\cong I_+^*(\overline Y,E)\\
I^\infty_*(Y,E) &\cong I_\infty^*(\overline Y,E),\end{align*}
\end{theorem}
\begin{proof}We will show that $\text{DCI}(Y)^\vee$ is isomorphic, as a $\Lambda(u)$-module, to $\text{DCI}(\overline Y)$. This implies the rest of the results, as one may check that at the level of $\mathcal C^-_A$-modules, we have $\left(C^+_A(D)\right)^\vee = C^-_A(D^\vee)$ as well as $\left(C^-_A(D)\right)^\vee = C^+_A(D^\vee)$. These equalities preserve the norm map, so we get the same result for Tate homology.

Before doing this calculation, recall that if $D$ is a chain complex, its dual $D^\vee$ has the differential $\delta f = (-1)^{|f|+1} fd$. This is isomorphic as an $A$-module to the chain complex with $\delta f = f d$, with isomorphism given by $f \mapsto (-1)^{|f|(|f|+1)/2} f$. (Observe that this makes sense for any chain complex graded over $\mathbb Z/4N$ for any integer $N \geq 0$.) In what follows, we will use the second differential on the dual.

Recall the definition of the chain complex $$DCI(Y) = C_*^{\text{irr}}(Y) \oplus C_*^{\text{irr}}(Y)[3] \oplus C_*^{U(1)}(Y) \oplus C_*^{U(1)}(Y)[2] \oplus C_*^\theta(Y)$$ stated before Corollary \ref{DCI}. Passing to the dual (and slightly rewriting), we obtain $$DCI(Y)^\vee = C_{-*-3}^{\text{irr}}(Y) \oplus C_{-*-3}^{\text{irr}}(Y)[3] \oplus C_{-*-2}^{U(1)}(Y) \oplus C_{-*-2}^{U(1)}(Y)[2] \oplus C_{-*}^\theta(Y),$$ with differential $$\widetilde{\partial_{CI}}^\vee := \begin{pmatrix}\partial_1 & 0 & 0 & 0 & 0\\U_{\text{Fl}} & -\partial_1 & V_3 & V_1 & D_1 \\ V_2 & 0 & 0 & 0 & 0\\ V_4 & 0 & 0 & 0 & 0\\ D_2 & 0 & 0 & 0 & 0\end{pmatrix}.$$

Now if $\pi$ is a regular perturbation on $Y$, then the same perturbation is regular on $\overline Y$. Using the time-reversal symmetry $\mathbb R \times \overline Y \cong \mathbb R \times Y$, if $\mathbf{A}$ is a connection on the latter going from critical orbits $\alpha$ to $\beta$, then it is sent to a connection $\overline{\mathbf{A}}$ going from $\beta$ to $\alpha$ on $\mathbb R \times \overline Y$; in particular, the indices of the corresponding deformation operators are equal: $\text{ind}(Q^\nu_{\mathbf{A},\pi}) = \text{ind}(Q^\nu_{\overline{\mathbf{A}},\pi})$. Recalling the definition of grading from Definition \ref{grading}, we see that $$\text{gr}_Y(\alpha, \beta) + (\dim \alpha - \dim \beta) = \text{gr}_{\overline Y}(\beta, \alpha).$$ In particular, if $\beta$ is the trivial connection, we see that $i_{\overline Y}(\alpha) = -i_Y(\alpha) - \dim \alpha$. Thus the graded basis for $DCI(Y)^\vee$ described above is precisely a graded basis for $DCI(\overline Y)$.

Furthermore, this time-reversal symmetry gives isomorphisms $$\overline{\mathcal M}_{Y,z}(\alpha, \beta) \cong \overline{\mathcal M}_{\overline Y, -z}(\beta, \alpha).$$ Following this isomorphism, we may use the same basepoints $b_\alpha, q_\alpha$ and chain homotopies $h_\alpha$ to define $DCI(\overline Y)$ as for $DCI(Y)$. Once we do this, we have equalities of the operators
\begin{align*}D_1(Y)^\vee &= D_2(\overline Y),\\
V_1(Y)^\vee &= V_2(\overline Y),\\
V_3(Y)^\vee &= V_4(\overline Y),\\
U_{\text{Fl}}(Y)^\vee &= U_{\text{Fl}}(\overline Y).
\end{align*}

Therefore the isomorphism of graded vector spaces $DCI(Y)^\vee \cong DCI(\overline Y)$ described above, which also preserves the $\Lambda(u)$-action, is in fact an isomorphism of chain complexes. 
\end{proof}

\begin{remark}While this theorem is surely true for arbitrary principal ideal domains $R$, there are technical obstructions in proving it. One natural idea is to try to construct an equivariant pairing $\widetilde{CI}_*(Y) \otimes \widetilde{CI}_*(\overline Y) \to R$, or at least a pairing on a quasi-isomorphism $C_*^{\textup{gm}}(SO(3);R)$-submodule of that tensor product (one defined by the demand that the chains intersect transversely, as in \cite{McClure}). To do this, one needs to be able to put a module structure on the tensor product. If $A$ is a dg-algebra and we wish to endow the category of $A$-modules with a tensor product structure, we require a comultiplication $\Delta: A \to A \otimes A$. While this is easy in the model of $C_*(SO(3);R)$ with simplicial chains (given by the Alexander-Whitney map), it seems unlikely this is possible for $C_*^{\text{gm}}(SO(3))$; we have no way to cut a $\delta$-chain into canonical pieces. Perhaps a modification of this model exists that admits the structure of a bialgebra (or even a Hopf algebra), but this is not clear to the author.

A more straightforward approach is to extend the Donaldson model to work with arbitrary principal ideal domains $R$, where instead of $C_*(SO(3);R) \simeq \Lambda(u)$, one would instead use $C_*(SO(3);R) \simeq \mathcal A$, where $\mathcal A$ is a dg-algebra with a generator $a_i$ in each degree $0 \leq i \leq 3$, with the relations $$d(a_2) = 2a_1,\; a_1^2 = 0,\;\text{and } a_1a_2 = a_3 = -a_2 a_1.$$ The above proof then generalizes easily. The construction of such a Donaldson model, together with the complete determination of the differential $\partial_{DCI}$, will be carried out in forthcoming work with Daemi and Scaduto.
\end{remark}

\appendix

\chapter{Equivariant homology of dg-modules}
Let $G$ be a finite group acting on a topological space $X$. The first algebraic object studied in \emph{equivariant algebraic topology} is its \textit{equivariant (co)homology}: if $EG$ is a contractible $G$-space with a free and proper action of $G$, we may define the equivariant homology $H_*^G(X)$ as the homology of the \textit{Borel space} $(X \times EG)/G$.

As discussed in Chapter \ref{sec:fourflavors}, an investigation of equivariant Poincar\'e duality --- or rather, its failure --- reveals the existence of a second algebraic invariant of the $G$-space $X$: its `coBorel homology'. When $X$ is a manifold, this corresponds to $H^{\dim G - *}_G(X)$, and no Poincar\'e duality result relates this to the usual equivariant homology. There is, however, a natural comparison map between these --- and that map sits in an exact triangle with a third group which measures the failure of Poincar\'e duality, the `Tate homology' $H^\infty_G(X)$. 

The definitions explored in \cite{GreM} are defined for \textbf{spectra}. Our main interest is applying these constructions to define equivariant instanton homology groups; unfortunately, the necessary `Floer homotopy theory' needed to define an `$SO_3$-equivariant instanton Floer spectrum' seems quite far out of reach --- all that is available to us is a chain complex with an action of $SO_3$. Instead, we will want to give purely algebraic definitions of equivariant cohomology of chain complexes.

An algebraic approach to the Borel construction, which works more generally for any $\mathbb Z[G]$-module, replaces $C_*(X)$ with a resolution in the category of chain complexes over $\mathbb Z[G]$, takes the quotient by $\mathbb Z[G]$, and then computes the homology (or cohomology) of its totalization. We will want a generalization of this process that will work for compact Lie groups $G$ in some sense ``acting on a chain complex" $C$, which gives equivariant (co)homology in the case of $G$ acting on $C_*(X)$. 

We model this by considering dg-modules over the dg algebra $C_*(G)$. After that, we will describe the dual homology theory, called \textit{coBorel homology}, and a homology theory called \textit{Tate homology} that compares the two. Our approach to the Borel/bar and coBorel/cobar constructions are strongly inspired by \cite{GM} and \cite{barthel2014six}, while the approach to Tate homology is essentially that of \cite{Klein}. 

The instanton `chain complex' $\widetilde C(Y; R)$ is $\Bbb Z/8$-graded, not $\Bbb Z$-graded. To apply our constructions to this level of generality, we conclude with a discussion of extensions to the case of complexes graded over $\mathbb Z/2N$, equipped with an appropriate object called a \emph{periodic filtration}.\\

\noindent \textbf{Convention.} All chain complexes and dg-algebras are $\Bbb Z$-graded until the final section, where $\Bbb Z/N$-graded chain complexes are introduced; even there, dg-algebras remain $\Bbb Z$-graded.

\section{Bar constructions}
We begin by fixing conventions. From here onwards, $R$ is a principal ideal domain (PID) which will serve as the ground ring of all of our chain complexes.

A chain complex of $R$-modules is a collection of $R$-modules $C_i$, indexed by $i \in \mathbb Z$, equipped with an $R$-linear map $d_i: C_i \to C_{i-1}$ so that the composite $d_{i-1} d_i: C_i \to C_{i-2}$ is zero. The tensor product of $\mathbb Z$-graded chain complexes $C$ and $D$ is $$(C \otimes D)_n = \bigoplus_{i+j = n} C_i \otimes D_j,$$ with differential acting on $a \otimes b \in C_i \otimes D_j$ as $$d(a \otimes b) = (d a) \otimes b + (-1)^{i} a \otimes (db)$$ and swap isomorphism $$\tau(a \otimes b) = (-1)^{|a||b|} (b \otimes a).$$

We will introduce the basic constructions for modules over a differential graded algebra that we use to define the various flavors of equivariant homology. 

Let $A$ be a homologically graded unital dg-algebra (over the ground ring $R$); so $A$ is a $\mathbb Z$-graded chain complex of $R$-modules whose product satisfies the graded Leibniz rule $d(ab) = d(a)b + (-1)^{|a|} ad(b)$. Furthermore, assume $A$ has an augmentation $\epsilon: A \to R$ with $\epsilon(1) = 1$; then we may identify 
$$\text{ker}(\epsilon) = A/\langle 1\rangle =: \overline A,$$
and give $R$ the natural structure of an $A$ (bi)-module. In a differential graded left $A$-module $N$, the same Leibniz rule must hold and the unit should act by the identity; for right $A$-modules $M$, the sign in the Leibniz rule uses the grading of $N$ instead of the grading of $A$.

\begin{definition}\label{bardef}The (two-sided) bar construction $B(M,A,N)$ is the totalization of the resolution 
$$M \otimes N \leftarrow M \otimes \overline A \otimes N \leftarrow  M \otimes \overline A^{\otimes 2} \otimes N \leftarrow \cdots$$ 
So as a graded module, 
$$B(M,A,N) = \oplus_{n=0}^\infty M \otimes \overline A[1]^{\otimes n} \otimes N.$$ 
where $\overline A_k = \overline A[1]_{k+1}$. Writing a generic tensor product as $m \left[a_1 \mid \cdots \mid a_k\right] n$ where possibly $k = 0$, its differential is given as
\begin{align*}(-1)^k &\Bigg(dm [a_1, \mid \cdots \mid a_k ] n + \sum_{i=1}^k (-1)^{|m|+\varepsilon_{i-1}} m\left[a_1\mid \cdots \mid da_i \mid \cdots \mid a_k\right]n \\
&+ (-1)^{|m|+\varepsilon_k} m\left[a_1\mid \cdots \mid a_k\right] dn \Bigg)+ \Bigg(ma_1 [a_2 \mid \cdots \mid a_k] n \\
+&\sum_{i=1}^{k-1} (-1)^i m \left[a_1 \mid \cdots \mid a_i a_{i+1} \mid a_k\right] n + (-1)^k m \left[a_1 \mid \cdots \mid a_{k-1}\right] a_k n\Bigg)\
\end{align*}
Here $\varepsilon_i = |a_1| + \cdots + |a_i|$.
\end{definition} 

\noindent These sign conventions are those of \cite{GM}. Gugenheim and May define $\text{Tor}^A(M,N)$ for arbitrary pairs of dg-modules over a dga; their Corollary A.9 shows that $$\text{Tor}^A_*(M,N) = H_*B(M,A,N)$$ as long as $M$, $A$, and $N$ satisfy appropriate flatness hypotheses. We will state, and then assume, the relevant hypotheses shortly; they are there to ensure that $B(M,A,N)$ takes quasi-isomorphisms to quasi-isomorphisms.

Most of the time we are interested in the special case $B(M,A,R)$. If one thinks of $A = C_*(G)$ as the canonical dg-algebra of interest, acting on the module $M = C_*(X)$ for some right $G$-space $X$, then $B(M,A,R)$ models the Borel construction, also known as the `homotopy quotient', $$M_{hG} = M \times_G EG.$$ Note that whenever $a_k n$ appears in the differential above, this term is zero, as $a_k \in \overline A = \text{ker}(\epsilon)$, acting on $n \in R$ via the augmentation $\epsilon$.

The bar construction is functorial (as an $A$-module) under maps 
$$f: A \to A', \;g: M \to M', \; h: N \to N'$$ 
with $g(ma) = g(m)f(a)$ and similarly for $h$. It is also functorial under homotopies of such maps. Most commonly one either fixes the modules or $A$ when using this functoriality; for a map of triples 
$$(M,A,N) \to (M',A,'N'),$$ 
the corresponding map of bar constructions factors as 
$$B(M,A,N) \to B(M',A,N') \to B(M',A',N').$$
There is a canonical map, natural in $(M,A,N)$, 
\begin{equation}\label{comodule}B(M, A, N) \to B(M,A,R) \otimes B(R,A,N)\end{equation} 
given as 
$$m [a_1 \mid \cdots \mid a_k] n \mapsto \sum_{i=0}^k (-1)^{(k-i)(|m|+\varepsilon_i)}m[a_1 \mid \cdots \mid a_i] \otimes [a_{i+1}\mid \cdots \mid a_k] n.$$ 
This endows $BA := B(R,A,R)$ with the structure of a coalgebra, and $B(M,A,R)$ the structure of a right $BA$-comodule. 

\begin{definition}The cobar construction of $(N,A,M)$, where $A$ is a dg-algebra and both $N$ and $M$ are right $A$-modules, is the \emph{chain} complex $$cB(N,A,M) = \textup{Hom}_A(B(N,A,A),M).$$
\end{definition} 

This, too, is a special case Gugenheim and May's constructions; now we have $H_* cB(N,A,M) = \text{Ext}_A^{-*}(N,M)$ under suitable projectivity hypotheses. We will primarily be interested in the case $N = R$, as this represents the homotopy fixed points of $M$.

For concreteness, we note that $cB(N,A,M)$ is isomorphic as a graded $R$-module to 
$$\prod_{i=0}^\infty \text{Hom}_R(N \otimes \overline A[1]^{\otimes i}, M) = \text{Hom}_R(N \otimes BA, M),$$ 
graded so that $|a| + |\eta| = |\eta(a)|$. This is the negative of the usual grading (which would make $\text{Hom}_R(C,D)$ into a \textit{cochain} complex); for $M$ bounded and $A$ nonnegatively graded, the \textit{chain} complex we define here is bounded above but unbounded in the negative direction! Passing through the above isomorphism, the differential of an element $\eta: B(N,A,R) \to M$ of degree $d$ is given as

\begin{align*}(d\eta)(n[a_1 | \cdots | a_p]) &= d_M\big(\eta(n[a_1 | \cdots | a_p])\big) - (-1)^d \eta\left(d_{B(N,A,R)}n[a_1 | \cdots | a_p]\right) \\
&- (-1)^{d+p} \eta(n[a_1 | \cdots | a_{p-1}])a_p.
\end{align*}

The cobar construction is functorial for triples $(M,A,N)$; it is covariant in $M$, but contravariant in $A$ and $N$. The equivariance condition on morphisms is that we demand $f: A' \to A$ and $g: M \to M'$ satisfy $g(mf(a')) = g(m)a'$; the condition on $A$ and $N$ is the usual, that $h(n'a') = h(n')f(a')$.

Dual to equation (\ref{comodule}), $$cB(R,A,R) = \text{Hom}_A(B(R,A,A), R) \cong \text{Hom}_R(BA, R)$$ is naturally an algebra. The chain complex $cB(R,A,M)$ is naturally a left module over this algebra, following the diagram \begin{align*}&\text{Hom}_R(B(R,A,R),R) \otimes \text{Hom}_A(B(R,A,A), M) \\
&\to \text{Hom}_A(B(R,A,R) \otimes B(R,A,A), M)  \to \text{Hom}_A(B(R,A,A), M) = cB(R, A, M).\end{align*} 
Here the final map is dual to the $BA$-comodule structure on $B(R,A,A)$; note that the map $B(R,A,A) \to B(R,A,R) \otimes B(R,A,A)$ is a map of right $A$-modules, so it makes sense to apply $\text{Hom}_A$ to this map.

Finally, observe that $B(M, A, R)$ is a \emph{left} module over $cB(R,A,R)$, by applying the comodule structure and then the pairing. In the following map, we suggestively rewrite $cB(R,A,R) = BA^\vee$: 
\begin{align*} BA^\vee \otimes B(M,A,R) &\xrightarrow{1 \otimes \Delta} BA^\vee \otimes B(M,A,R) \otimes BA \\
\xrightarrow{\tau \otimes 1} B(M,A,R) &\otimes BA^\vee \otimes BA \xrightarrow{1 \otimes \operatorname{eval}} B(M,A,R).\end{align*}

The map $\Delta$ is the comodule structure and $\tau$ is the swap map $x \otimes y \mapsto (-1)^{|x||y|} y \otimes x$. Here, and elsewhere in this text, we write $M^\vee$ for an $R$-module $M$ to be its dual as an $R$-module, $\text{Hom}_R(M, R)$; if $M$ was a right $A$-module, then $M^\vee$ carries the structure of a left $A$-module by $(a \eta)(m) = \eta(ma)$. 

\section{Invariance}
Suppose we have a map of pairs $(M,A,N) \to (M',A',N')$ as above inducing a map $B(M,A,N) \to B(M',A',N')$. When is the induced map an isomorphism on homology? Because the bar construction is functorial under homotopies, this is true if the map of triples is a homotopy equivalence (of triples). More generally, we have the following theorem. These flatness restrictions are harmless for most purposes in topology and for use in this paper, but indicate that the approach taken here is too naive for general modules over a dga. This result, and much of \cite{GM}, is put into the powerful general framework of model categories in \cite{barthel2014six}, but we will not need this language here.

\begin{theorem}\label{bar-inv}Suppose we have a map of triples 
$$(g,f,h): (M,A,N) \to (M',A',N')$$ 
so that $f$, $g$ and $h$ all induce isomorphisms on homology (from here on we will say ``are quasi-isomorphisms"). Suppose further that $N, N', \overline A$ and $\overline A'$ are flat (as graded $R$-modules). Then the induced map $B(M,A,N) \to B(M',A',N')$ is a quasi-isomorphism.
\end{theorem}

Note in particular that we may take $N = R$ in the above theorem. A similar result is true of $cB(N,A,M)$ (making slightly stronger assumptions). 

\begin{theorem}\label{cb-inv}Suppose $M,N$ are $A$-modules and $N, N'$ are $A'$-modules. Let 
$$f: A' \to A, \;g: M \to M' \; h: N' \to N$$ 
be quasi-isomorphisms that are equivariant in the sense that 
$$g(mf(a')) = g(m)a' \text{ and } h(n'a') = h(n')f(a').$$
If all of $A, \; A', \; N$ and $N'$ are $R$-free, then $cB(N,A,M) \to cB(N',A',M')$ is a quasi-isomorphism.
\end{theorem}

These theorems are proved by appealing to natural spectral sequences associated to the bar and cobar constructions. All of our spectral sequences arise from filtered complexes (as opposed to more intricate constructions, such as exact couples), so let us introduce the concept. 

\begin{definition}\label{filtrations}
A \emph{filtered complex} is a chain complex $C$ equipped with an increasing (possibly unbounded) sequence of subcomplexes $$\cdots \subset F_{p-1} C \subset F_p C \subset F_{p+1} C \subset \cdots$$ 

Four special properties of filtrations $F_\bullet C$ deserve a name: \begin{enumerate}
\item A filtration is \emph{exhaustive} if $\bigcup_p F_p C = C$. 
\item A filtration is \emph{Hausdorff} if $\bigcap_p F_p C = 0$.
\item A filtration is \emph{complete} if the natural map $C \to \lim_{p \to -\infty} C/F_p C$ is an isomorphism.
\item A filtration is \emph{regular} if, for all $n \in \Bbb Z$, there exists $u(n) \in \Bbb Z$ so $H_n(F_p C) = 0$ for $p < u(n)$. 
\end{enumerate}

Given a filtered complex, the complex obtained by taking the minimal quotients $\textup{gr}_p C = F_p C/F_{p-1} C$ is called the \emph{associated graded complex}.

Finally, a \textup{filtered chain map} (of degree $m$ and filtration index $n$) is a chain map $f: C \to C'$ of degree $m$ between filtered complexes for which $f(F_p C) \subset F_{p+n} C'$ for all $p \in \Bbb Z$. 
\end{definition}

The associated graded complex can be understood as a bigraded group $E^0_{p,q} = \text{gr}_p C_{p+q}$, equipped with a differential $d^0 = \text{gr}(d)$ of bidegree $(p,q) = (0, -1)$. Taking the homology with respect to this differential gives $E^1_{p,q}$; careful inspection shows that $d: C_p \to C_{p-1}$ gives rise to a differential $d^1$ on $E^1$ of bidegree $(-1, 0)$. Continuing in this way, we obtain a \emph{spectral sequence}, a sequence of bigraded complexes $(E^r, d^r)$ where $d^r$ has bidegree $(-r, r-1)$, and $E^{r+1}_{p,q} = H(E^r_{p,q}, d^r)$. Intuitively, the successive $E^r$ give Taylor polynomial approximations of $H(C)$, considering only the part of the differential which decreases filtration by at most $r-1$.

Notice that a filtered chain map $f: C \to C'$ of degree $m$ and index $n$ induces a chain map of bidegree $(n,m-n)$ on the associated graded complex, and hence a map on the $E^1$ page of the spectral sequence, and so on. See \cite[Chapter XV]{CE} for a concise reference; their conventions are cohomological, as compared to our preference for homological gradings. 

We will be interested in two questions: 

\begin{itemize}
\item When can this spectral sequence be used to compute the homology $H(C)$? 
\item If that's not possible, when can this spectral sequence be used to detect isomorphisms?
\end{itemize}

The first statement should be qualified. The filtration on $C$ gives rise to a filtration $F_p H(C)$ on the homology, and our goal here is to compute the associated graded $\text{gr}_p H_{p+q}(C)$. Secondly, what we mean\footnote{The sense in which $E^r$ `converges to' $\text{gr} H(C)$ discussed here is weaker than that discussed in \cite[Chapter XV.2]{CE}, but is perhaps the notion more familiar to most topologists, and the one which can more easily be understood by simple inspection. Readers interested in a more general notion of `convergence' are directed to the given reference.} by \emph{compute} is that the $E^r_{p,q}$ should limit to this group $\text{gr}_p H_{p+q}(C)$. But there is not usually a canonical map $E^r_{p,q} \to E^{r+1}_{p,q}$, as there is rarely a canonical map from a complex to its homology; what should this mean? 

For regular filtrations, for each $(p,q)$ there exists $N$ so that for $r \ge N$ we have $d^r_{p,q} = 0$. As a result, the passage from $E^r_{p,q}$ to $E^{r+1}_{p,q}$ is given by quotienting by a subgroup, and in particular we have a sequence $$E^N_{p,q} \to E^{N+1}_{p,q} \to \cdots$$ and we may set $E^\infty_{p,q} = \colim_r E^r_{p,q}$. 

The following result is stated at the conclusion of \cite[Chapter XV.4]{CE}, and follows immediately from the discussion thusfar in that chapter. (They impose the exhaustive and complete Hausdorff assumptions at the beginning of Chapter XV.)

\begin{proposition}\label{ssCompute}
If $C$ is an exhaustive, complete, Hausdorff, regular complex, there is a canonical isomorphism $$\textup{gr}_p(H_{p+q} C) \cong E^\infty_{p,q} = \colim_r E^r_{p,q}.$$ 
\end{proposition}

We will only ever care about filtrations that are exhaustive and complete Hausdorff (but it is still important to check these conditions). Completeness in particular is often subtle and delicate. On the other hand, we will care about both regular and non-regular filtrations. The regularity assumption is easy to check in practice; for instance, if for each $n$ the sequence $F_p C_n$ is bounded below in $p$, we automatically have regularity.\\

If $F_\bullet C$ is not regular (but still exhaustive, and complete Hausdorff) then we cannot use the filtration to compute the groups themselves, but we can still \emph{check isomorphisms}, which is for our purposes usually more important. 

\begin{proposition}[Eilenberg--Moore comparison theorem]\label{ssComparison}
Suppose $C$ and $C'$ are filtered complexes whose filtrations are exhaustive and complete Hausdorff. If $f: C \to C'$ is a filtered map which induces an isomorphism $f^r: E^r(C) \to E^r(C')$ for $r \ge 0$ on some finite page of the spectral sequence, then $f$ induces an isomorphism $H(C) \to H(C')$. 
\end{proposition}

This is proved as \cite[Theorem 5.5.11]{weibel1995introduction}, and the proof is short: the mapping cone of $f$ can be equipped with a filtration for which $E^r(\text{Cone}(f)) = 0$. This implies that the filtration on $\text{Cone}(f)$ is regular, and it follows from Proposition \ref{ssCompute} that $H(\text{Cone}(f)) = 0$, which therefore implies that $f$ is a quasi-isomorphism.

\begin{remark}\label{rmk:completion}
Suppose $C$ is equipped with an arbitrary filtration. Consider the filtered complex, called the \emph{full completion} in \cite{GOH}, $$\hat C = \lim_{q \to -\infty} \colim_{p \to \infty} F_p C/F_q C,\quad  \text{with filtration} \quad F_p \hat C = \lim_{q < p} F_p C/F_q C.$$ Then the filtration on $\hat C$ is exhaustive and complete Hausdorff by construction. If the filtration on $C$ is regular, the associated filtration on $\hat C$ is regular, and the spectral sequence $E^r(C)$ actually converges to (the associated graded of) the homology of $\hat C$. Similarly, if $f: C \to C'$ is a filtered chain map, the associated spectral sequences can be used to detect if the associated map $\hat f: \hat C \to \hat C'$ is a quasi-isomorphism; it says little about $f$ itself.
\end{remark}

Notice that if a filtration is bounded below, it is automatically both complete and Hausdorff; if it is bounded above, it is automatically exhaustive. Another situation that frequently arises is that the filtration is the totalization of a \textit{multicomplex}; in this case, the full completion is easy to describe. We state the definition here (first appearing in \cite{wall_1961}).

\begin{definition}\label{multicomplex}A multicomplex is a bigraded $R$-module $M_{s,t}$ with differentials $d_r: M_{s,t} \to M_{s-r,t+r-1}$ for $r \geq 0$ so that $$\sum_{i+j = n} d_i d_j = 0.$$ The associated (completed) filtered complex is the subcomplex $$\hat C_n = \prod_{s \to -\infty} M_{s, n-s} \subset \prod_{s \in \mathbb Z} M_{n, n-s}$$ consisting of those sequences $(x_s)$ with $x_s = 0$ for sufficiently large $s$. (That is to say, the product is only taken in the negative direction.) The differential is given as $\sum_{r \geq 0} d_r$.
\end{definition}

\begin{proposition}The completed complex $\hat C$ of a multicomplex $M$ is complete Hausdorff, and we can identify the $E^1$ page of the associated (conditionally converging) spectral sequence with $H(M,d_0)$ equipped with the differential $H(d_1)$.
\end{proposition}

These tools in hand, we can prove the two invariance results we need.

\begin{proof}[Proof of Theorem A.1]
Filter $B(M,A,N)$ by $$F_p B(M,A,N) = M \otimes \left(\oplus_{i=0}^p \overline A^{\otimes i} \otimes N\right).$$ Then we calculate the $E^1$ page of the associated spectral sequence as the (bigraded) complex $$H(M \otimes N) \leftarrow H(M \otimes \overline A \otimes N) \leftarrow \cdots$$ Because the map $B(M,A,N) \to B(M',A',N')$ preserves the filtration, it induces a map of spectral sequences; if the maps $$H(M \otimes \overline A^{\otimes i} \otimes N) \to H(M' \otimes \overline A'^{\otimes i} \otimes N')$$ are isomorphisms, then we would have proved that $B(g,f,h)$ induces an isomorphism on the $E^1$ page and hence all pages of the spectral sequence.

We prove this fact inductively on the number of tensor factors. 

\begin{lemma}If $F: X \to X'$ is a quasi-isomorphism of degreewise $R$-flat complexes, and $G: Y \to Y'$ is a quasi-isomorphism of arbitrary complexes, then $F \otimes G: X \otimes Y \to X' \otimes Y'$ is a quasi-isomorphism.
\end{lemma}

This follows immediately from comparing the short exact sequences (which exist because $X$, $X'$ are degreewise flat)
$$0 \to \bigoplus_{p+q=n} H_p(X) \otimes H_q(Y) \to H_n(X \otimes Y) \to \bigoplus_{p+q=n-1} \text{Tor}(H_p X, H_qY) \to 0,$$
because the outside terms only depend on $HX$ and $HY$. This is where we need $R$ to be a PID; otherwise we would need to apply the Kunneth spectral sequence, and would need boundedness assumptions. One may use the spectral sequence argument to extend this result to Dedekind domains, as in \cite[Lemma~2.2]{kriz1995operads}. 

The filtration of $B(M,A,N)$ is trivially complete and Hausdorff, as $F_{-1} = 0$. It is exhaustive because the infinite direct sum is the union of its finite direct summands.

Because $N, N', \overline A$, and $\overline A'$ are degreewise $R$-flat, we have $B(g,f,h)$ is an isomorphism on the $E^2$ pages. It follows from Proposition \ref{ssComparison} that $B(g,f,h)$ is a quasi-isomorphism.
\end{proof}

\begin{proof}[Proof of Theorem A.2]
Now the filtration is $$F_{-p} cB(N,A,M) = \text{Hom}_A(B(N,A,A)/F_p B(N,A,A), M);$$ that is, it consists of those functionals that vanish on $F_p B(N,A,A)$. Let us abbreviate $B(N,A,A) =: B$. This filtration is now bounded above, and thus is automatically exhaustive. Note that the intersection $$\lim_{p} \text{Hom}_A(B/F_p, M) = \text{Hom}_A(\colim_p B/F_p, M).$$ The colimit is an exact functor, so $$\colim_p B/F_p = B/\colim_p F_p = B/B = 0,$$ because the filtration of $B$ was exhaustive, and thus the filtration is Hausdorff. Now, the short exact sequence $0 \to F_p \to B \to B/F_p \to 0$ is \emph{split} as $A$-modules (not dg), identifying $B/F_p$ with the summand $\oplus_{i > p} N \otimes \overline{A}^{\otimes i} \otimes A$ (but ignoring the differential). This implies that 
$$\text{Hom}_A(B, M)\big/\text{Hom}_A(B/F_p, M) \cong \text{Hom}_A(F_p, M)$$ 
as $R$-modules ($\text{Hom}$ takes split exact sequences to split exact sequences). Applying 
\begin{align*}\lim \text{Hom}_A(B, M)\big/\text{Hom}_A(B/F_p, M) &\cong \lim \text{Hom}_A(F_p, M) \\
&= \text{Hom}_A(\colim F_p, M) = \text{Hom}_A(B, M),
\end{align*}
we see that the filtration is complete (the last equality because the filtration on $B$ was exhaustive).

Using the isomorphism of $cB(N,A,M)$ as a graded $R$-module to $\text{Hom}(N \otimes BA, M)$, we identify the $E^1$ page as the totalization of the double complex $$H\big(\text{Hom}_R(N,M)\big) \to H\big(\text{Hom}_R(N \otimes \overline A, M)\big) \to \cdots$$
We need to see that $\text{Hom}(N \otimes \overline A^{\otimes i}, M) \to \text{Hom}(N' \otimes \overline A'^{\otimes i}, M')$ is a quasi-isomorphism for all $i$; then the theorem will be proved. We write this as a lemma: the only property we will use is that $N \otimes \overline A^{\otimes i}$ is free for any $i$ (a tensor product of free modules is free).
\end{proof}

\begin{lemma}If $X$ and $X'$ are $R$-free, and $F: X' \to X$ is a quasi-isomorphism, then for any quasi-isomorphism $G: Y \to Y'$, the induced map $\textup{Hom}_R(X, Y) \to \textup{Hom}_R(X', Y')$ is a quasi-isomorphism.\end{lemma}
\begin{proof}
When both $X$ and $d(X) \subset X$ are complexes of projective modules, we have a natural Kunneth short exact sequence $$0 \to \prod_{p+q = n-1} \text{Ext}^1(H_p X, H_{-q} Y) \to H_{-n} \text{Hom}(X, Y) \to \prod_{p+q = n} \text{Hom}(H_p X, H_{-q} X) \to 0;$$ from this the theorem is clear, as long as every submodule $d(X)$ of a free module of arbitrary rank is free; this is true for PIDs. For the statement of the Kunneth theorem and that submodules of free modules are projective (and hence free, as projective modules over a PID are free), see \cite[Exercises~3.6.1-3.6.2]{weibel1995introduction}.
\end{proof}

We will often be interested in thinking of one-variable versions of the bar and cobar constructions as providing homology theories for $A$-modules; we will introduce new notation for the sake of compactness. Our notation is chosen to fit with both \cite{Jones} and \cite{OS}.

\begin{definition}Suppose that $A$ is an augmented dg-algebra with each $A_n$ free over $R$, and $M$ a right $A$-module. We define the positive and negative $A$-chains on $M$ to be 
\begin{align*}C_A^+(M) &:= B(M,A,R)\\
C_A^-(M) &:= cB(R,A,M).\end{align*} 
The homology of $C_A^+(M)$ is denoted $H_A^+(M)$, the (Borel, or positive) $A$-homology of $M$, while the homology of $C_A^-(M)$ is denoted $H_A^-(M)$, the coBorel (or negative) $A$-homology of $M$. 
\end{definition}
Note that both $C^+_A(M)$ and $C^-_A(M)$ are \emph{covariant} in $M$! When $A$ is clear and we are not interested in varying it, we omit it from the notation.

\section{The dualizing complex and Tate homology}
Tate homology, constructed in this section for dg-modules over a dg-algebra and written as $H^\infty_A(M)$, has appeared in the literature in many forms. Our approach here is essentially a chain-level interpretation of \cite{Klein} (which was written in the context of Tate homology of $G$-spectra). Tate homology may be viewed, in some sense, as the homology theory that arises when you kill off free objects. The classic reference to Tate homology of spectra is \cite{GreM}; we warn that if $X$ is a $G$-space, our $H^\infty_{C_*(G)}(C_*(X);R)$ is more analogous to what they would call $t_{HR \wedge G}(HR)(X)$, not $t_G(HR)(X)$. (We are `chainifying' the group from the start as well.) The idea of Tate homology is beautifully developed in an abstract homotopical setting in \cite{Green}, which surely includes as a special case the content of this section. A recent approach to Tate homology of G-spectra via localization of $(\infty,1)$-categories appears in \cite{NS}; this level of abstraction has the advantage of giving multiplicativity and uniqueness results that are not as easily available otherwise (in some cases, not available at all).

We take our current approach as it seems to minimize input energy, at the cost perhaps of some conceptual clarity and multiplicativity results (we do not construct a product on $C^\infty_A(R)$, for instance).

Because $A$ is a left $A$-module, the complex $cB(R, A, A) = C^-_A(A)$ inherits the structure of a left $A$-module, by $(a\eta)([b_1 \mid \cdots \mid b_n]b) = a\cdot \eta([b_1 \mid \cdots \mid b_n]b)$.

\begin{definition}The dualizing complex $D_A$ of a degreewise $R$-free dg-algebra $A$ is the left $A$-module $C^-_A(A)$.\end{definition}

\begin{definition}\label{tate-def}Let $A$ be a degreewise $R$-free dg-algebra and $M$ a right $A$-module. The map $$N_M: B(M, A, D_A) \to C^-_A(M)$$ given by 
\[ N_M\left(m[a_1 \mid \cdots \mid a_k] \psi\right) = \begin{cases} 
      0 & k > 0 \\
      m \cdot \psi & k = 0
   \end{cases}
\] 
is called the \textit{norm map} of $M$. The product $m \cdot \psi$ makes $\psi$ into an $M$-valued functional by using the fact that $\psi$ is $A$-valued and $M$ is a right $A$-module. The mapping cone of $N_M$ is denoted $C^\infty_A(M)$, the \textit{Tate complex} of $M$. Its homology is the Tate homology $H^\infty_A(M)$ of $M$.
\end{definition}

The following is important enough to record as a lemma.
\begin{lemma}\label{flatcheck}$D_A$ is degreewise $R$-flat.\end{lemma}
\begin{proof}
The algebra $A$ is degreewise $R$-free, and in particular torsion-free. Because $R$ is a PID, $R$-flatness is equivalent to ($R$-) torsion-freeness. The complex $\text{Hom}_A(B(R,A,A), A)$ is $R$-torsion free: if $\psi: B(R,A,A) \to A$ has $(r \psi)(x) = 0$ for all $x$, then $\psi(x) = 0$ for all $x$ as $A$ is torsion-free. Therefore $D_A$ is torsion-free, hence degreewise flat. 
\end{proof}

To justify the name ``Tate homology", we show that this satisfies part of the corresponding versions of Klein's axioms defining Tate cohomology \cite{Klein}, skipping the complete verification that $H^\infty_A$ and $H\left(B(-, A, D_A)\right)$ are homology theories: the remaining axioms state that these preserve homotopy pullbacks and filtered homotopy colimits. These are not hard to verify, but we will not use them. 

\begin{theorem}The functor $C^\infty_A(M)$ from right $A$-modules to chain complexes satisfies Klein's axioms specifying the Tate homology of $M$: \begin{enumerate}
    \item $C^\infty_A(M)$ preserves weak equivalences;
    \item $H^\infty_A(X \otimes A) = 0$, where $X$ is a finite-dimensional chain complex of free $R$-modules and the right $A$-module structure is given by acting on $A$;
    \item There is a map $C^-_A(M) \to C^\infty_A(M)$, natural in $M$, whose homotopy fiber preserves weak equivalences. 
    \end{enumerate}
\end{theorem}

\begin{proof}
It follows from Theorem \ref{bar-inv} that $B(M, A, D_A)$ preserves weak equivalences in $M$, because by Lemma \ref{flatcheck} above, $D_A$ is flat. The Tate complex is the mapping cone (homotopy cofiber) of $N_M: B(M, A, D_A) \to C_*^-(M)$, and both the domain and the codomain of the norm map preserve weak equivalences, so $C_*^\infty(M)$ does as well. This is Axiom 1.

For Axiom 2, note that 
\begin{align*}B(X \otimes A, A, D_A) &\cong X \otimes B(A, A, D_A) \\
cB(R, A, X \otimes A) &= \text{Hom}_A(B(R,A,A), X \otimes A) \cong X \otimes D_A.\end{align*} 
The final identification uses the assumption that $X$ is finitely-generated and free. 
 
Then under these identifications, the norm map is identified with the projection $X \otimes B(A, A, D_A) \to X \otimes D_A$; but $B(A,A,N) \to N$ is a homotopy equivalence for any left $A$-module $N$, so the mapping cone of the norm map is contractible.

The first part of Axiom 3 is obvious (it is the inclusion into a mapping cone), and the second part almost so: the homotopy fiber = mapping cocone is naturally equivalent to $B(M, A, D_A)$.\end{proof}

We henceforth write $H^{+,\textup{tw}}_A(M)$ for $H(B(M, A, D_A))$, and call it the \emph{twisted Borel homology}. We will investigate its relationship to $H^+_A(M)$ later. 

We won't prove Klein's uniqueness theorem that these axioms do uniquely characterize Tate cohomology, but rather use it as motivation that we have the correct definition. (It seems likely that some variation of his argument works in this context.)

In addition to the above, it is important to observe that there is a natural left action of $C^-_A(R)$ on $D_A$. This action commutes with the left action of $A$ described at the beginning of this section, and hence the formula $$\varphi \cdot m[a_1 | \cdots | a_p] \psi = (-1)^{|\varphi|(|m| + \epsilon_p+p)} m[a_1 | \cdots | a_p] (\varphi \cdot \psi)$$ defines a dg-module structure of $C^-_A(R)$ on $B(M, A, D_A).$ The norm map $N_M$ is easily seen to be $C^-_A(R)$-equivariant. Therefore, the mapping cone inherits the structure of a left $C^-_A(R)$-module, and the natural map $C^-_A(M) \to C^\infty_A(M)$ is a module homomorphism.

The following theorem summarizes everything we have assembled about the three $A$-homology functors $H^{+,\textup{tw}}, H^-, H^\infty$.

\begin{theorem}\label{eq-package}
Let $A$ be a dg-algebra over a commutative PID $R$ which is ($R-$) flat in each degree $A_n$. There are functors $$H^{+,\textup{tw}}_A(M), \; \; H^-_A(M) \; \; H^\infty_A(M),$$
from dg $A$-modules to graded $R$-modules, satisfying the following properties.
\begin{enumerate}
    \item The functors send short exact sequences of $A$-modules to exact triangles and preserve weak equivalences.
    \item $H_A^-(R)$ is a ring, and each of these homology theories carry a natural left module structure over $H_A^-(R)$.
    \item There is an exact triangle of $H_A^-(R)$-modules $$H^{+,\textup{tw}}_A(M) \to H_A^-(M) \to H_A^\infty(M) \xrightarrow{[-1]} H^{+,\textup{tw}}_A(M) \to \cdots$$
    \item $H_A^\infty(X \otimes A) = 0$ when $X$ is a finitely-generated chain complex of free $R$-modules and $X \otimes A$ is given the canonical right action.
    \end{enumerate}
\end{theorem}

We state an invariance theorem for equivariant homology with respect to quasi-isomorphic dgas.

\begin{proposition}\label{dg-invar}Suppose $f: A\to A'$ is a quasi-isomorphism of algebras, each degreewise free over $R$. This induces a functor $F: \mathsf{Mod}_{A'} \to \mathsf{Mod}_A$ via restriction of scalars, and there are natural isomorphisms $H^\bullet_A(FM) \xrightarrow{\cong} H^\bullet_{A'}(M),$ for the homology theories $\bullet \in \{(+,\textup{tw}), -, \infty\}$. These natural isomorphisms induce an $H^-_A(R)$-equivariant natural isomorphism of exact triangles.
\end{proposition}

\begin{proof}
The fact that $cB(R,A,M) \leftarrow cB(R,A', M)$ is a quasi-isomorphism is Theorem \ref{cb-inv}; for that reason, the maps $$D_A = cB(R,A,A) \to cB(R,A,A') \leftarrow cB(R,A',A') = D_{A'}$$ are quasi-isomorphisms. This is the crucial place we need $A$ and $A'$ to be degreewise free. 

To make the diagrams smaller, let us write 
\begin{align*}B_A &= B(D_A, A, FM), \\
cB_A &= C^-_A(FM) = \text{Hom}_A(B(A,A,R), FM), \\
B' &= B(C^-_A(A'), A', M).\end{align*}
Then the induced maps on bar constructions $B_A \to B' \leftarrow B_{A'}$ are quasi-isomorphisms by Theorem \ref{bar-inv}, using that $D_A$, $C^-_A(A')$, and $D_{A'}$ are all flat $R$-modules, which is Lemma \ref{flatcheck}. 

Comparing the Tate homology groups is more complicated. The following diagram commutes, which is an easy check left to the reader. 

$$\begin{tikzcd}
B_A \arrow{r} \arrow{d}{N_A} & B' \arrow{ld}{N'} \arrow[leftarrow]{r} & B_{A'} \arrow{d}{N_{A'}} \\
cB_A \arrow{d} \arrow{rd} \arrow[leftarrow]{rr} & {} & cB_B \arrow{d} \\
C(N_A) \arrow[dashrightarrow]{r} &C(N') \arrow[dashleftarrow]{r} &C(N_{A'})
\end{tikzcd}$$
\newline
Here $N': \text{Hom}_A(B(A,A,R), A') \otimes_{A'} B(A',A',M) \to \text{Hom}_A(B(A,A,R), M')$ is obtained from the projection map $B(A', A', M) \to M$. The maps in the bottom row are the natural maps induced on a cone, and we denote the norm maps $N_A$ and $N_{A'}$ instead of the usual $N_M$ to reflect the different choice of algebra $M$ is a module over.

Because all of the maps in the top two rows are quasi-isomorphisms, the maps on the bottom row are also quasi-isomorphisms by an application of the five lemma. In particular, inverting the bottom-right quasi-isomorphism on homology, we have an isomorphism $H^\infty_A(FM) \to H^\infty_{A'}(M)$.

All of the maps in this diagram are module homomorphisms with respect to either $C^-_A(R)$ or $C^-_{A'}(R)$, as appropriate; this implies that the isomorphisms on homology are $H^-_A(R)$-equivariant. 
\end{proof}

\begin{remark}\label{rmk:dg-invar}
Suppose $f: A \to A'$ is a quasi-isomorphism and $M$ is a right $A$-module. Then there is a zig-zag of quasi-isomorphisms $M \leftarrow B(M, A, A) \to B(M, A, A') = M'$ of right $A$-modules; further, the $A$-module structure on $B(M,A,A')$ can be understood as the restriction of an $A'$-module structure. The above arguments then imply that this zigzag induces a natural isomorphism $H^\bullet_A(M) \cong H^\bullet_{A'}(M')$. As a result, we may fairly freely pass between different quasi-isomorphic models for a given dga. 
\end{remark}

\section{Periodicity in Tate homology}\label{sec:periodic}
Tate homology, in some cases, is \emph{periodic}: there is some class in $H^-_A(R)$ so that its action on $H^\infty_A(M)$ is an isomorphism when $M$ is a finite $A$-module (in an appropriate sense). A beautiful reference for this phenomenon is \cite[Section~III.16]{GreM} in the setting of genuine $G$-spectra, where $G$ is a compact Lie group \emph{that acts freely on some sphere}. Without this assumption, periodicity phenomena often fail; see for instance \cite{NegProducts}.

Here we endeavor only to prove an analogue of it in a simple case which will suffice for our purposes, and in a computational manner. We write $A = \Lambda := \Lambda(u_n)$, where $|u_n| = n$. We have $H^-_\Lambda(R) = R\llbracket U\rrbracket$, where $|U| = -n-1$, and $H^+_\Lambda(R) = R[U^*]$ with $|U| = n+1$, where the action of $U$ is contraction against $U^*$; these are immediately clear from the definition of bar and cobar construction here, which have no nonzero differentials.

\begin{remark}
Notice that $H^-_\Lambda(R)$ is only graded-commutative if $|n|$ is odd or $R$ has characteristic two. This corresponds to the fact that $\Lambda$ is a graded bialgebra precisely when one of those conditions holds; more generally, $H^-_A(R)$ is graded-commutative whenever $A$ is a dg-bialgebra.
\end{remark}

We see from the tautological exact triangle that, \emph{as a graded $R$-module}, we have $H^\infty_\Lambda(R) = R\llbracket U, U^{-1}]$ (where here we have suggestively rewritten $U^* = U^{-1}$); we know from the fact that the exact triangle are maps of $H^-_\Lambda(R)$-modules that the action of $U$ on these increases the power of $U$ by $1$, \emph{except possibly for the action on $U^{-1}$}. The content of the following crucial lemma is that $U \cdot U^{-1} = 1$.

\begin{lemma}\label{PeriodicComputation}Let $\Lambda$ be the exterior algebra over $R$ on an element $u$ in degree $n$. We have an $H^-_\Lambda(R)$-module isomorphism $H^\infty_\Lambda(R) \cong R \llbracket U, U^{-1}]$, with $U \cdot U^i = U^{i+1}$ for all $i \in \mathbb Z$.
\end{lemma}
\begin{proof}We write this out very explicitly; all of what follows is a transcription of definitions. First, we may write $D_\Lambda = R\llbracket \eta\rrbracket \otimes \Lambda$, with differential 
$$d(\eta^k \otimes x) = (-1)^{nk} \eta^{k+1} \otimes ux.$$ 
To avoid confusion later, we have used the notation $\eta$ where previously we wrote $U$: it is the degree $(-n-1)$ functional on $B\Lambda$ which sends $U^*$ to $1 \in R$. As graded $R$-modules we have
$$B(R, \Lambda, D_\Lambda) = R[U^*] \otimes R\llbracket \eta\rrbracket \otimes \Lambda(u);$$ 
under this isomorphism, the differential of $B(R, \Lambda, D_\Lambda)$ is taken to
$$d((U^*)^k \otimes \eta^j \otimes 1) = (-1)^{n(j+k)+k}(U^*)^k \otimes \eta^{j+1} \otimes u + (-1)^k (U^*)^{k-1} \eta^j \otimes u,$$ and is zero on $(U^*)^k \otimes \eta^j \otimes u$. Combining both the definition of the action of $C_\Lambda^-$ on $D_\Lambda$ and the extension to an action on $B(R, \Lambda, D_\Lambda)$, the action of $U$ on this complex is 
$$U \cdot ((U^*)^k \otimes \eta^j \otimes p) = (-1)^{(n+1)k+j} (U^*)^k \otimes \eta^{j+1} \otimes p$$ for any $p \in \Lambda(u)$. Finally, the norm map $N_\Lambda: B(R, \Lambda, D_\Lambda) \to C^-_\Lambda(R)$ is given by $\eta^k \mapsto U^k$ for any $k \geq 0$ and is otherwise zero.

By definition, the Tate complex $$C^\infty_\Lambda(R) = B(R, \Lambda, D_\Lambda)[1] \oplus C^-_\Lambda(R)$$ is the mapping cone of the norm map; this means that its differential is 
$$d_\infty = \begin{pmatrix}-d & 0\\N_\Lambda & d\end{pmatrix}.$$ 
The action of $U^*$ is the same as before on each component. Now note that 
$$H_* B(R, \Lambda, D_\Lambda) \cong H^+_\Lambda[n]$$ 
as an $H^-_\Lambda$-module, and a chain in $B(R, \Lambda, D_\Lambda)$ representing the degree $n$ element of $H^+_\Lambda[n]$ is given by $(U^*)^0 \otimes \eta^0 \otimes u$. Our goal, then, is to show that 
$$U \cdot ((U^*)^0 \otimes \eta^0 \otimes u) = (U^*)^0 \otimes \eta \otimes u$$ 
is homologous in $C^\infty_\Lambda(R)$ to $U^0 \in C^-_\Lambda(R)$. But from the formula for the differential, we have 
$$d_\infty((U^*)^0 \otimes \eta^0 \otimes 1) = -(U^*)^0 \otimes \eta \otimes u + U^0;$$ this is precisely what we wanted.
\end{proof}
We can use this to prove the following localization theorem.

\begin{proposition}\label{TateLoc}Let $M$ be a dg-module over $\Lambda$, degreewise free over $R$; suppose $M$ has a finite filtration $0 = F_{-1}M \subset F_0 M \subset \cdots \subset F_n M = M$ so that each piece $F_k M/F_{k-1}M$ of the associated graded dg-module is quasi-isomorphic to a finite direct sum of copies of $R$ and $\Lambda$. Then the action of $U \in H^-_\Lambda$ on $H^\infty_\Lambda(M)$ is an isomorphism. Therefore the natural map $H^-_\Lambda(M) \to H^\infty_\Lambda(M)$ factors through $H^-_\Lambda(M)[U^{-1}]$, and the map $$H^-_\Lambda(M)[U^{-1}] \to H^\infty_\Lambda(M)$$ is an isomorphism, natural for dg-module homomorphisms of such $M$.
\end{proposition}

To prove this, we first define the notion of \emph{inverting an endomorphism} of a dg-module with respect to a dg-module homomorphism $f: M \to M$ of degree $k$, written $M[f^{-1}]$. This is \textbf{NOT} the strict notion of inversion of an element familiar in module theory, $$\operatorname{colim} \left(M \xrightarrow{f} M \xrightarrow{f} M \xrightarrow{f} \cdots\right),$$ because this rarely plays well with taking homology. Instead it is a \emph{homotopy colimit}. The simplest way to phrase this (taken from \cite[Definition~24.5]{HRanicki}) is that $M[f^{-1}]$ is the mapping cone of the map $1-tf: M[t] \to M[t]$, where $t$ is a polynomial generator in degree $-k$. Then immediate from the exact triangle on homology, and the fact that $(1-tf)_*: H(M)[t] \to H(M)[t]$ is injective, we see that $H(M[t^{-1}]) = H(M)[t^{-1}]$, where on the right side we are inverting an element of a module in the usual sense. (Here we are using that the module-theoretic notion may be defined perfectly well as $M[t]/(1-tf)$.) 

\begin{proof}[Proof of Proposition \ref{TateLoc}]First we show that the action of $U^*$ on $H^\infty_\Lambda(M)$ is an isomorphism. This is equivalent to showing that the natural map $C^\infty_\Lambda(M) \to C^\infty_\Lambda(M)[U^{-1}]$ is a quasi-isomorphism. To see this, we use that the filtration $F_k M$ induces a filtration $F_k C^\infty_\Lambda(M) = C^\infty_\Lambda(F_k M)$, and similarly a filtration on $C^\infty_\Lambda(M)[U^{-1}]$; these are complete because the filtration is finite. Now the $E^1$ page of the corresponding spectral sequence is given as a direct sum of copies of $H^\infty_\Lambda(R)$ and $H^\infty_\Lambda(\Lambda)$ (or, respectively, the results of inverting $U$), and the action of $U$ on the $E^1$ page is the direct sum of the corresponding actions. But $H^\infty_\Lambda(\Lambda) = 0$, and the fact that $U^*$ is an isomorphism on $H^\infty_\Lambda(R)$ was the content of Lemma \ref{PeriodicComputation}. So the natural (filtered) map $C^\infty_\Lambda(M) \to C^\infty_\Lambda(M)[U^{-1}]$ is an isomorphism on the $E^1$ page, and therefore a quasi-isomorphism. This proves the first part of the theorem.

The rest follows similar lines: there is a natural map $$C^-_\Lambda(M)[U^{-1}] \to C^\infty_\Lambda(M)[U^{-1}],$$ and our goal is to show that this map is a quasi-isomorphism; because the map $C^-_\Lambda(M) \to C^-_\Lambda(M)[U^{-1}]$ is identified on homology with the result of inverting $U$, the desired result follows. But using the same filtration as the above, the $E^1$ page of the first spectral sequence is a direct sum of copies of $H^-_\Lambda(R)[U^{-1}]$ and $H^-_\Lambda(\Lambda)[U^{-1}]$, and the $E^1$ page of the second spectral sequence is a direct sum of corresponding copies of $H^\infty_\Lambda(R)$ and $H^\infty_\Lambda(\Lambda)$. Then the only thing to observe is that what we already know: $H^-_\Lambda(R)[U^{-1}] \to H^\infty_\Lambda(R)$ is an isomorphism, and that $H^-_\Lambda(\Lambda)$ is a copy of $R$ concentrated in degree zero, so that $H^-_\Lambda(\Lambda)[U^{-1}] = 0$. Therefore the above map is an isomorphism on the $E^1$ page of the corresponding spectral sequence, and therefore a quasi-isomorphism.\end{proof}

\section{Simplifying the twisted Borel homology}\label{borel}
In this section, we give conditions under which there is a natural isomorphism $H^{+,\textup{tw}}_A(M) \cong H^+_A(M)[n]$ for some degree shift $n$. We would, furthermore, like this isomorphism to preserve the action of $H^-_A(R)$. The desired conditions for us amount to a sort of \emph{Poincar\'e duality} assumptions, which will be true for the groups $C_*(G)$ of most interest to us. 

\begin{definition}
Suppose $A$ is an augmented dg-algebra. Consider $A^\vee = \text{Hom}_R(A, R)$ as an $A$-bimodule, via $$(a \phi b)(x) = (-1)^{|a|(|\phi|+|b|+|x|)}\phi(bxa).$$ 
\begin{itemize}
\item We say that $A$ satisfies \emph{strong Poincar\'e duality of degree $n$} if there is a zig-zag of $A$-bimodule quasi-isomorphisms $A \simeq A^\vee[n]$. 
\item We say that $A$ satisfies \emph{weak Poincar\'e duality of degree $n$} if there is a zig-zag of augmented dg-algebra quasi-isomorphisms relating $A$ to an augmented dg-algebra $\tilde A$ with the following property: there exists a quasi-isomorphism $q: \tilde A \to \tilde A^\vee[n]$ of right $\tilde A$-modules, and $\tilde A$ contains a cycle $F \in \tilde A_n$ satisfying $\tilde a \cdot F = \epsilon(\tilde a)F$ for all $\tilde a \in \tilde A$, and for which $q(F)$ is homologous to the augmentation $\epsilon \in \tilde A^\vee$. 
\end{itemize}
\end{definition}

Despite the names, note that \textbf{strong Poincar\'e duality does not necessarily imply weak Poincar\'e duality}, because of the new hypothesis on the existence of an `invariant fundamental class' $F$. Rather, strong Poincar\'e duality essentially supposes the existence of a \textit{bimodule} equivalence between $A$ and $A^\vee[n]$, while weak Poincar\'e duality supposes slightly more than the existence of a \textit{right module} equivalence between $A$ and $A^\vee[n]$. A mere right module equivalence is insufficient for the argument we want to make below, and we need the added hypothesis about the cycle $F$; nonetheless, this second condition remains weaker, in the sense that far more examples in nature satisfy weak duality than strong duality.

These hypotheses will be relevant to simplifying the dualizing module $D_A$, because we have equivalences $$R[n] \simeq \text{Hom}_R(B(R, A, A), R)[n] \cong \text{Hom}_A(B(R,A,A), A^\vee)[n]$$ of left $A$-modules; the first map is the quasi-isomorphism given by sending $1$ to the augmentation map $\epsilon$, and the second the isomorphism given by sending $\eta$ to the map $\eta': EA \to A^\vee$ defined by $\eta'(x)(a) = \eta(xa)$. 

Therefore, if we can replace $A$ with $A^\vee$ as a right module, we can at least conclude that $H_*(D_A) = R[n]$; and if we can replace $A$ with $A^\vee$ as a bimodule, we see that $D_A \simeq R[n]$ as left $A$-modules, which is essentially enough for the desired result.

\begin{remark}This is a very special case of what is called a \emph{Gorenstein condition} in the literature, which usually amounts to something like an assumption that there exists a bimodule isomorphism $\text{Hom}_A(R, A) \cong A[n]$, interpreted in a derived sense. Our condition implies this because of the chain of equivalences 
\begin{align*}\text{Hom}_A(R,A) \simeq \text{Hom}_A(R, \text{Hom}_R(A, R))[n] &= \text{Hom}_R(A \otimes_R R, \text{Hom}_A(A, R))[n] \\
&= \text{Hom}_R(A, R)[n] \simeq A.\end{align*}
The second equality uses the tensor-hom adjunction for (bi)modules over two rings, and all appearances of ``$\text{Hom}$" are derived. We will not use the Gorenstein condition except through the above Poincar\'e duality condition; this remark is purely motivational.

A particularly nice reference, which applies to the dg-algebra case, is \cite{dwyer2006duality}.
\end{remark}

\begin{theorem}\label{PDiso-weak}
If $A$ is $R$-flat and satisfies weak Poincar\'e duality of degree $n$, then there is a natural isomorphism of graded $R$-modules $H^{+,\textup{tw}}_A(M) \cong H^+_A(M)[n]$. 
\end{theorem}
\begin{proof}
If $f: A' \to A$ is a dg-algebra quasi-isomorphism, there is a corresponding isomorphism $H^+_{A'}(M) \cong H^+_A(M)$ as $R$-modules, and similarly with the twisted Borel homology. Therefore, it suffices to prove the stated claim for the algebra $\tilde A$ related to $A$ by a zig-zag of quasi-isomorphisms. To simplify notation, in what follows we assume $A = \tilde A$.

As discussed above, this implies that $D_A \simeq R[n]$ as $R$-modules; but to replace $B(M, A, D_A)$ with $B(M, A, R)$ requires an equivalence between $D_A$ and $R$ as left $A$-modules. To do so, we will construct a specific map $R[n] \to D_A$ which is $A$-invariant and induces a quasi-isomorphism, using our knowledge of $H_*(D_A)$. 

Let $\psi: EA \to A$ be the $A$-equivariant map given by $$\psi(a[a_1 | \cdots | a_n]) = \begin{cases} \epsilon(a) F & n = 0 \\ 0 & n > 0 \end{cases}$$ where $F$ is the fundamental class assumed in our definition of weak Poincar\'e duality. It is straightforward to verify that $\psi$ defines a cycle of degree $n$ in $D_A$, and we let $f_\psi: R[n] \to D_A$ be the chain map $f_\psi(r) = r\psi$. Note that this map is $A$-equivariant: we have $$f_\psi(a \cdot r) = f_\psi(\epsilon(a) r) = \epsilon(a) r \psi,$$ while $(a \cdot f_\psi)(r) = a r\psi$; but using that $a \cdot F = \epsilon(a) F$ and that $R$ is central in $A$, we have $$(a f_\psi)(r)(b) = (ar\psi)(b) = (ar)\epsilon(b) F = \epsilon(a) r \epsilon(b) F = \epsilon(a) r \psi(b) = f_\psi(ar)(b).$$

It remains to verify that $f_\psi$ is a quasi-isomorphism. To see this, we trace our equivalences backwards. The quasi-isomorphism $$q_*: \text{Hom}_A(B(R, A, A), A) \to \text{Hom}_A(B(R, A, A), A^\vee)$$ sends $\psi$ to $$\psi'(a[a_1 | \cdots | a_n]) = \begin{cases} \epsilon(a) q(F) & n = 0 \\ 0 & n > 0 \end{cases}.$$ Because $q(F)$ is homologous to the augmentation, up to adding a coboundary we may suppose $\psi'(1) = \epsilon$. Next, the isomorphism $\text{Hom}_A(B(R,A,A), A^\vee) \cong \text{Hom}_R(B(R,A,A), R)$ sends $\psi'$ to $$\psi''(a[a_1 | \cdots | a_n]) = \psi'(a[a_1 | \cdots | a_n])(1) = \begin{cases} \epsilon(a) & n = 0 \\ 0 & n > 0 \end{cases}$$ This is precisely the map obtained as the image of $1$ under the quasi-isomorphism $R \to \text{Hom}_R(B(R,A,A), R)$, and therefore $\psi$ generates the homology of $D_A$, so that $f_\psi$ is a quasi-isomorphism.

The natural isomorphism $H^{+, \textup{tw}}_A(M) \cong H^+_A(M)[n]$ then follows from the invariance theorem for bar constructions, Theorem \ref{bar-inv}.
\end{proof}

\begin{theorem}\label{PDiso}If $A$ is $R$-flat and satisfies strong Poincar\'e duality of degree $n$, then there is a natural isomorphism $H^{+,\textup{tw}}_A(M) \cong H^+_A(M)[n]$ of $H^-_A$-modules.
\end{theorem}
\begin{proof}
The argument is essentially unchanged if the zig-zag of equivalences between $A$ and $A^\vee$ has length larger than one, so we suppose we have a single bimodule quasi-isomorphism $A \to A^\vee[n]$. 

As discussed above, this induces a chain of equivalences \begin{align*}D_A = \text{Hom}_A(B(R,A,A), A) &\simeq \text{Hom}_A(B(R,A,A), A^\vee)[n] \\ &\cong \text{Hom}_R(B(R,A,A), R)[n] \simeq R[n]\end{align*} as left $A$-modules. It is better to stop in this chain before the end, as this equivalence $D_A \simeq \text{Hom}_R(B(R,A,A), R)$ is equivariant under the actions of $C^-_A(R)$ essentially by definition: the product structure uses the left $BA$-comodule structure of $B(R, A, A)$, and all of the maps above only involved the rightmost factor. We write $B(R, A, A) = EA$ for convenience; this equivalence says $D_A \simeq (EA)^\vee$, this equivalence equivariant under both the left $A$-action and the left $C^-_A(R) = (BA)^\vee$-action.

Now we have a map $$B(M, A, (EA)^\vee) \to B(M, A, R) \otimes_R B(R, A, (EA)^\vee) = B(M, A, R) \otimes_R EA \otimes_A (EA)^\vee.$$ The pairing $(EA) \otimes (EA)^\vee \to R$ is $A$-equivariant, so it factors through $EA \otimes_A (EA)^\vee$, and thus we have a composite $\eta: B(M, A, (EA)^\vee) \to B(M, A, R)$. The claim is that $\eta$ is a quasi-isomorphism, and is equivariant under the $C^-_A(R) = (BA)^\vee$-action.

To see that this map is a quasi-isomorphism, consider the map $\tau: B(M,A,R) \to B(M,A,(EA)^\vee)$ induced by the dual of the augmentation $EA \to R$. Because the map $R \to (EA)^\vee$ is a quasi-isomorphism, our invariance result Theorem \ref{bar-inv} establishes that $\tau$ is also a quasi-isomorphism. Further, the composite $\eta \tau$ is equal to the identity. It follows that $\eta$ is also a quasi-isomorphism.

As for equivariance, we should check this at the chain level following the sign conventions in \cite{TLSigns}. Write $\psi$ for an element of $(EA)^\vee$ and $\beta$ for an element of $(BA)^\vee$; then a generic element of $B(M, A, (EA)^\vee) = B(M, A, A) \otimes_A (EA)^\vee$ is written as $m[a_1 | \cdots | a_n ] \psi$. Following the given map and then applying $\beta$, we have
\begin{align*}&\beta \otimes m[a_1 | \cdots | a_n] \psi \mapsto \beta \otimes \sum_{i=0}^n (-1)^{s_1} m[a_1 | \cdots | a_i] \otimes [a_{i+1} | \cdots | a_n]\psi \\
&\mapsto \beta \otimes \sum_{i=0}^n (-1)^{s_2} m[a_1 | \cdots | a_i] \psi([a_{i+1} | \cdots | a_n]) \\
&\mapsto \beta \otimes \sum_{0 \leq j \leq i \leq n} (-1)^{s_3} m[a_1 | \cdots | a_j] \otimes [a_{j+1} | \cdots | a_i] \psi([a_{i+1} | \cdots | a_n])\\
&\mapsto \sum_{0 \leq j \leq i \leq n} (-1)^{s_4} m[a_1 | \cdots | a_j] \otimes \beta([a_{j+1} | \cdots | a_i]) \psi([a_{i+1} | \cdots | a_n]).
\end{align*}

Here the exponents are
\begin{align*}
s_1 &= (n-i)(|m|+\epsilon_i)\\
s_2 &= (n-i)(|m|+\epsilon_i) + |\psi|(\epsilon_n - \epsilon_i)\\
s_3 &= (i-j)(|m|+\epsilon_j)+(n-i)(|m|+\epsilon_i) + |\psi|(\epsilon_n - \epsilon_i)\\
s_4 &= (|\beta|+i-j)(|m|+\epsilon_j)+(n-i)(|m|+\epsilon_i) + |\psi|(\epsilon_n - \epsilon_i).
\end{align*}

Here recall that $\epsilon_i := |a_1| + \cdots + |a_i| + i$ is the degree of $[a_1| \cdots | a_j]$. We have used the rule that the swap map is $a \otimes b \mapsto (-1)^{|a||b|} b \otimes a$.

On the other hand, if we apply $\beta$ first to $\psi$ and then follow the given map, we have
\begin{align*}&\beta \otimes m[a_1 | \cdots | a_n] \psi \mapsto (-1)^{t_1} m[a_1| \cdots | a_n] (\beta \cdot \psi) \\
&\mapsto \sum_{0 \leq j \leq n} (-1)^{t_2} m[a_1 |\cdots |a_j] \otimes [a_{j+1} | \cdots | a_n] (\beta \cdot \psi) \\
&\mapsto \sum_{0 \leq j \leq n} (-1)^{t_3}m[a_1|\cdots |a_j] (\beta \cdot \psi)([a_{j+1} | \cdots | a_n])\\
&= \sum_{0 \leq j \leq i \leq n} (-1)^{t_4}m[a_1 |\cdots |a_j](\beta \otimes \psi)([a_{j+1} | \cdots | a_i] \otimes [a_{i+1} | \cdots | a_n])\\
&= \sum_{0 \leq j \leq i \leq n} (-1)^{t_5} m[a_1 | \cdots | a_j] \beta([a_{j+1}| \cdots | a_i]) \psi([a_{i+1}| \cdots | a_n]),
\end{align*}

where the exponents are 
\begin{align*}
t_1 &= |\beta|(m+\epsilon_n)\\
t_2 &= (n-j)(|m|+\epsilon_j)+|\beta|(|m|+\epsilon_n)\\
t_3 &= (n-j)(|m|+\epsilon_j)+|\beta|(|m|+\epsilon_n) + (|\beta|+|\psi|)(\epsilon_n - \epsilon_j) \\
t_4 &= (n-j)(|m|+\epsilon_j)+|\beta|(|m|+\epsilon_n) + (|\beta|+|\psi|)(\epsilon_n - \epsilon_j)+ (n-i)(\epsilon_i -\epsilon_j)\\
t_5 &= (n-j)(|m|+\epsilon_j)+|\beta|(|m|+\epsilon_n) + (|\beta|+|\psi|)(\epsilon_n - \epsilon_j)+ (|\psi| + n-i)(\epsilon_i -\epsilon_j).
\end{align*}

\noindent The first and second expressions clearly agree, at least up to the signs on each factor. To check that $(-1)^{t_4} = (-1)^{s_5}$, observe the congruences mod 2
\begin{align*}&(n-j)(|m|+\epsilon_j) + |\beta|(|m|+\epsilon_n) + (|\beta| + |\psi|)(\epsilon_n - \epsilon_j) + (|\psi|+n-i)(\epsilon_i - \epsilon_j)\\
&\equiv (n-j+|\beta|)|m| + |\psi| \epsilon_n + (|\psi| + n-i) \epsilon_i + (|\beta|+i-j) \epsilon_j \\
&\equiv (|\beta| + i -j)(|m| + \epsilon_j) + (n-i)(|m|+\epsilon_i) + |\psi|(\epsilon_n - \epsilon_i).
\end{align*}

\noindent These are easily seen by breaking the top and bottom formulas (which are the relevant exponents of $-1$) into the components labelled by $|m|, \epsilon_n, \epsilon_i,$ and $\epsilon_j$ in the center formula.

The natural isomorphism $H^{+, \textup{tw}}_A(M) \cong H^+_A(M)[n]$ then follows from the invariance theorem for bar constructions, Theorem \ref{bar-inv}, and we have seen this is an isomorphism of $H^-_A(R)$-modules.
\end{proof}

\section{Spectral sequences}
Making some further mild assumptions on the algebras, we have useful spectral sequences for calculating the various flavors of $A$-homology.

\begin{proposition}If $A$ is a non-negatively graded dg-algebra ($A_n = 0$ for $n < 0$), there is a projection of dg-algebras $\pi: A \to H_0 A$, through which we can have $A$ act on $H_qM$. Then for any $\bullet \in \{+, -,\infty\}$ there exists an exhaustive and complete Hausdorff filtration of $C^\bullet_A(M)$ by $C^-_A(R)$-submodules for which the $E^1$ page of the associated spectral sequence is $$E^2_{p,q} = H_p^\bullet(A, H_q M).$$ 
If $H_0 A = R$, the action of $A$ on $H(M)$ is trivial; if $H(M)$ is flat over $R$, we may identify $E^2_{p,q} \cong H^\bullet_p(A) \otimes H_q(M)$.

The filtration of $C^\bullet_A(M)$ is also regular if $\bullet = +$ or if $H_*(M)$ is bounded above. 
\end{proposition}

In particular, we may use this spectral sequence for computations of $H^+_A(M)$, or for all of these theories when $H_*(M)$ is bounded above. In all cases, we may use these spectral sequences to detect quasi-isomorphisms. 

\begin{proof}
Recalling that $$C^+_A(M) = B(M,A,R) \cong M \otimes BA$$ as graded modules, we filter $$F_p B(M, A, R)_n = \oplus_{i \leq p} (M)_{n-i} \otimes BA_{i};$$ that is, $F_p B(M,A,R)$ consists of elements with total $BA$-degree at most $p$. Notice that this filtration is transparently exhaustive and Hausdorff; completeness and regularity follow from the fact that $BA$ is concentrated in non-negative degrees.

The differential on $B(M,A,R)$ can be written (ignoring signs) as three terms:

\begin{align*}d(m[a_1 \mid \cdots \mid a_n]) =& \big ((dm) [a_1 \mid \cdots \mid a_n]\big )\\ &+ \big (m a_1[a_2 \mid \cdots \mid  a_n]\big )\\ &+ \big (m d([a_1 \mid \cdots \mid a_n])\big ).\end{align*}

The first and last term clearly preserve the filtration, and the second does because $A_n = 0$ for $n < -1$. The first differential of the associated spectral sequence is the differential on the associated graded complex, which is $d_0 = d_M \otimes 1_{BA}$. We can thus identify the $E^1$ page with $BA \otimes HM$ (remember $BA$ is degreewise $R$-flat). The piece that decreases filtration by exactly $1$ is $d_1$, given by the differential in $BA$ and multiplication by elements of $A_0$. Therefore we identify that the $E^1$ page is given as $B(HM, A, R)$, where a positive degree element of $A$ acts trivially (these components of the differential decrease the filtration by at least $2$) and the action of $A_0$ factors through $A_0 \to H_0 A$. This filtration is a pleasant example of the filtration on the totalization of a \textit{multicomplex}, as in Definition \ref{multicomplex}. 

Finally, observe that the action of $C^-_A(R)$ is filtered, as the action is defined by contraction against $BA$ and by assumption $A$ is non-negatively graded, so $C^-_A$ has no elements of positive degree.\\

The proof for the $H^-$ and $H^\infty$ spectral sequences follow mostly as before. As graded modules, we may write $$cB(R,A, M)_n = \prod_{i \to -\infty} cB(R,A, R)_i \otimes M_{n-i},$$ where the symbol $\prod_{i \to -\infty}$ is as in Definition \ref{multicomplex}, and filter as above $$F_p cB(R,A,M)_n = \prod_{i \leq p} cB(R,A,R)_i \otimes M_{n-i}.$$ Again, this filtration is transparently exhaustive and Hausdorff, and is complete because we take a product as $i \to -\infty$; if $H_*(M)$ is bounded below (say, $H_k(M) = 0$ for $k > d$) then  for $p < n - d$ we have $H_n(F_p C^-_A(M)) = 0$: the same filtration remains exhaustive and complete Hausdorff on $F_p C^-_A(M)$, but now the degree-$n$ part of the $E^2$ page is $$E^2(F_p C^-_A(M))_{i,j} = \begin{cases} H_i^-(A) \otimes H_j(M) & i \le p \\ 0 & i > p \end{cases},$$ so if $p < n - d$ then because $i + j = n$ and $i \le p$ we have $j > d$, so that the given group always vanishes with the given restriction on $p$. (A general form of this argument is spelled out more explicitly in the language of convergence in \cite[Theorem 5.5.10(2)]{weibel1995introduction}.)

We needed $A$ to be nonnegatively graded for this to be a filtration of complexes, and identify the $E^2$ page as before, and similarly observe that the filtered subcomplexes are $C^-_A(R)$-submodules.\\

As for $C^\infty_A(M)$, as a graded $R$-module we have $$C^\infty_A(M)_n = \prod_{i \to -\infty} C^\infty_A(R)_i  \otimes M_{n-i}.$$ The arguments for $C^-$ apply again without change.
\end{proof}

\section{Group algebras}\label{sec:gpalg}
The most important application of the $A$-homology functors is for $A = C_*(G;R)$, where $G$ is a topological group (but most importantly a compact Lie group). The product structure is given as the composite $$C_*(G) \otimes C_*(G) \xrightarrow{EZ} C_*(G \times G) \xrightarrow{\times} C_*(G),$$ where $EZ$ is the Eilenberg-Zilber map which sends $\sigma \otimes \tau$ to a sum of simplices in a standard triangulation of the product $\Delta^i \times \Delta^j$; see \cite[Section 26(b)]{FHT} for more details on this product and other structure on $C_*(G)$. The unit is $[e] \in C_0(G)$ and the augmentation $C_0(G) \to R$ is the natural augmentation (add up points). These makes $C_*(G;R)$ into an associative augmented algebra. For the rest of this section, we frequently abuse notation and write $G$ for $C_*(G;R)$ unless there is danger of confusion.

We begin with the following well-known fact to justify our definition of group homology; a more detailed proof may be found in \cite[Theorem~3.9]{GM}. Here $EG$ denotes a contractible left $G$-space that $G$ acts freely on, and for a right $G$-space $X$ we denote $X_{hG} = (X \times EG)/G$ for the Borel construction on $X$.

\begin{lemma}\label{groupcoho}
Let $G$ be a compact Lie group and $X$ a right $G$-space. Then there is an algebra isomorphism $H^-_G(R) \cong H^{-*}(BG; R)$, as well as compatible $H^-_G(R)$-module isomorphisms \begin{align*}
H_G^+(C_*(X;R)) &\cong H_*(X_{hG}; R) \\
H_G^-(C^{-*}(X; R)) &\cong H^{-*}(X_{hG}; R) 
\end{align*}
\end{lemma}

Notice that $C^{-*}(X;R)$ is naturally a \textbf{left} dg $A$-module, not a right module; so what we really mean here by $H_G^-(C^{-*}(X; R))$ is the homology of the complex $\text{Hom}_A(B(A,A,R), C^{-*}(X; R))$. This is harmless; because $C_*(G)$ enjoys a dg-algebra anti-automorphism $\iota: C_*(G) \to C_*(G)$ given by pointwise inversion of simplices with $\iota^2 = 1$, we may use $\iota$ to pass back and forth between left- and right- modules at will. 

\begin{proof}
Write $M = C_*(X; R)$. We first observe there is a quasi-isomorphism natural in $X$ between $\psi_X: C_A^+(M) \to C_*(X_{hG})$. This follows because $(EG \times X)/G$ can be given as the bar construction $B(X,G,\text{pt})$ in topological spaces, and the natural map $$B(C_*(X),C_*(G),R) \to C_*\big(B(X,G,\text{pt});R\big)$$ given levelwise by Eilenberg-Zilber maps is a quasi-isomorphism (filter the former by number of tensor factors, and the latter by number of join factors; apply the K\"unneth theorem levelwise). Taking $X = \text{pt}$, in particular we have a natural quasi-isomorphism $\psi_\ast: C_A^+(R) \to C_*(BG;R)$. It is straightforward to see that $\psi_\ast$ preserves the coalgebra structures, and that $\psi_X$ is a comodule map with respect to $\psi_\ast$. Taking duals gives us an algebra quasi-isomorphism $C^{-*}(BG;R) \to C^-_A(R)$, establishing the first claim. 

More generally, observe that for any right $A$-module $M$, we have an isomorphism $C_A^+(M)^\vee \cong C_A^-(M^\vee),$ preserving the natural $C_A^-$ action. Take $M = C_*(X;R)$ to see that there is an isomorphism of $R$-modules as stated in the lemma; because $\psi_X$ respects comodule structures, its dual respects module structures, which establishes the second claim.
\end{proof}

Our goal is to calculate what we need to make Theorem \ref{eq-package} practically useful. There are two pieces to this goal: first, simplifying the twisted Borel homology into a more calculable object; second, computing the flavors $H^\bullet_G(G/H)$ for the various (right) $G$-orbits.\\ 

\begin{proposition}\label{prop:G-is-weak-PD}
Suppose $G$ is a compact Lie group of dimension $n$. Then $C_*(G;R)$ satisfies weak Poincar\'e duality of degree $n$, and there is a natural isomorphism $H^{+,\text{tw}}_{G}(M) \cong H^+_G(M)[n]$ as graded $R$-modules. 
\end{proposition}

\begin{proof}
First, replace $C_*(G)$ with the quasi-isomorphic algebra $$\bar C_d(G) = \begin{cases} C_d(G) & d < n \\ C_n(G)/B_n(G) & d = n \\ 0 & d > n \end{cases}$$ where we truncate to degrees $[0, n]$; that this again produces a dg-algebra follows because $C_* G$ is concentrated in non-negative degrees (otherwise $B_n(G)$ might not be closed under multiplication by $0$-chains). Now $\bar Z_n(G) = H_n(G;R) \cong R^{\pi_0 G}$ and if one orients $\mathfrak g$ there is a canonical left-invariant orientation on $G$, hence a canonical generator of $\bar Z_n(G) \in \bar C_n(G)$. This is transparently invariant under left-multiplication by any $g \in G$ (so $C_0(G)$ acts via the augmentation), and $C_d(G)$ acts trivially for degree reasons for $d > 0$. So $\bar C_*(G)$ has a fundamental class.

Write $S^{\mathfrak g}$ for the one-point compactification of the \emph{right} adjoint representation, which can be identified with $B_\epsilon(\mathfrak g)/S_\epsilon(\mathfrak g)$ for any $\epsilon > 0$. There is an adjoint-equivariant map $\log: G \to S^{\mathfrak g}$ given by collapsing the complement of a neighborhood of the identity. More precisely, the exponential map $\exp: B_\epsilon(\mathfrak g) \to G$ is an embedding, and the image of the interior of $B_\epsilon$ is an open subset $G_{< \varepsilon} \subset G$. The exponential map thus gives an equivariant homeomorphism $\exp: S^{\mathfrak g} \cong G/G_{\ge \varepsilon}.$ The logarithm map is the composite $$G \to G/G_{\ge \varepsilon} \xrightarrow{\exp^{-1}} S^{\mathfrak g}.$$ 

We may define a map $$\Phi: C_*(G) \to \text{Hom}(C_* G, \tilde C_* S^{\mathfrak g}), \quad g \mapsto \Phi_g(h) = \log_*(gh).$$ The Hom complex has a natural left action via $(\alpha \cdot g')(h) = \alpha(g'h)$, and this map is clearly left-equivariant. One may compose $\Phi$ with any cocycle $\psi: \tilde C_*(S^{\mathfrak g}) \to R[n]$ which takes the fundamental class to $1$ to obtain a left-equivariant map $\hat \Phi: C_*(G) \to \text{Hom}(C_* G, R)[n]$. 

Finally, $\hat \Phi$ is a quasi-isomorphism. To see this, observe that the slightly modified map $\Phi'_g(h) = \psi\left(\log_*(g\iota(h))\right)$ is a model for Poincar\'e duality: working with field coefficients and passing to homology classes, the corresponding pairing corresponds to taking the intersection product of cycles, which is well-known to be a perfect pairing. It follows that $\hat \Phi$ is a quasi-isomorphism for any field coefficients, and thus for $\Bbb Z$ coefficients, and thus with coefficients in any PID $R$. 

Finally, to apply Theorem \ref{PDiso-weak}, one needs to know that we can pass between quasi-isomorphic dgas with reckless abandon; this follows from Proposition \ref{dg-invar} and the discussion in Remark \ref{rmk:dg-invar}.
\end{proof}

\begin{remark}
On the one hand, it seems as though this argument \emph{should} be able to prove more than it does, because the logarithm map is adjoint-equivariant. In fact, \cite[Theorem 10.1]{Klein} establishes a space-level duality result, identifying $G$ as having the same stable bi-equivariant homotopy type as $\text{Map}(G, S^{\mathfrak g})$ by an Atiyah duality argument, where the right action mixes between the actions on $G$ and $S^{\mathfrak g}$: if $f$ is in the mapping space, we have $(fh)(g) = f(hg) \text{ad}_h$. One might hope to use this same strategy above, but one runs into the difficulty that while we have a well-behaved map $C_*(G) \otimes C_*(G) \to \tilde C_*(S^{\mathfrak g})$, we cannot curry it to an equivariant map $C_*(G) \to \text{Hom}(C_* G, \tilde C_* S^{\mathfrak g})$: $C_*(G)$ carries the structure of a bialgebra --- so one can put $G$-module structures on tensor products and Hom-spaces --- but the antipode does not satisfy the Hopf algebra identities on the nose (only up to homotopy), which is the structure used to get an equivalence between $\text{Hom}_{C_* G}(L \otimes M, N)$ and $\text{Hom}_{C_* G}(L, \text{Hom}_R(M, N))$, so we do not obtain a curried map. 

It is plausible that one can improve this to a map which is $A_\infty$-equivariant on the right by arguing that $C_*(G)$ satisfies the Hopf identities up to `coherent homotopy', which would be enough to establish most of the desired duality claims. However, one would expect a slightly different result depending on whether the adjoint representation is orientation-preserving or not (if not, then the right module structure on the dual should be twisted by the adjoint determinant character). This should be expected, as the adjoint determinant is the obstruction to choosing an \emph{bi-invariant} orientation, hence a bi-invariant fundamental class. 

This is of no consequence for the cases we are interested in, where the relevant groups are connected.
\end{remark}

With this, we can assemble the last of basic groups $H_G^+$, $H_G^-$, and $H_G^\infty$.

\begin{corollary}\label{cor:calc} Let $G$ be any compact Lie group. Then as $R$-modules, we have $H_* B(D_G, G, R) = H_{*-d}(BG;R),$ Thus, by Theorem \ref{eq-package} (4), $$H_{G}^\infty(\text{pt};R)_k = \begin{cases} H_{k-d-1}(BG; R) & k > d \\
H^{-k}(BG;R) & k \leq 0\end{cases}$$
\end{corollary}

Still, we would like to be able to keep track of the $H^-_G$-module structure, so the above is somewhat disappointing. What we need, more generally, is that $C_*(G)$ has a better-behaved model: 

\begin{definition}
We say a compact Lie group $G$ is \emph{strongly dualizable} over $R$ if  $C_*(G;R)$ satisfies strong Poincar\'e duality. 
\end{definition}

It seems plausible that this holds for all compact Lie groups $G$ \emph{whose adjoint representation is orientable}; in particular, for connected $G$. However, as outlined above, this seems difficult to prove and is not relevant to our main results. Instead, we will establish this slightly more generally than the cases of most interest: $G = SO(3)$ and for which either $\frac 12 \in R$ or $2 = 0$ in $R$. 

\begin{proposition}\label{prop:strongduals}
The exterior algebra $\Lambda_R(u)$ satisfies strong Poincar\'e duality for any ring $R$. Finite groups are strongly dualizable over any ring. The group $G = SO(n)$ is strongly dualizable over $\Bbb F_2$-algebras $R$.
\end{proposition}

\begin{proof}
For the exterior algebra on a generator with $|u| = d$, the shifted dual $A^\vee[d] = \text{Hom}(A, R[d])$ is isomorphic as $R$-modules to $R[0] \oplus R[d]$, with generators $\phi_0, \phi_d: \Lambda(u) \to R$ given by $\phi_0(u) = 1$ and $\phi_d(1) = 1$. The bimodule action is given by $$(\phi u)(y) = \phi(uy), \quad (u \phi)(y) = (-1)^{d(|\phi|+|y|}) \phi(yu).$$ Notice that $$(\phi_0 \cdot u)(1) = \phi_0(u) = 1, \quad (u \cdot \phi_0)(1) = \phi_0(u) = 1,$$ so that $u \cdot \phi_0 = \phi_0 \cdot u = \phi_d$. Similarly, $u \cdot \phi_d = \phi_d \cdot u = 0$. Thus the map $\Lambda \to \Lambda^\vee[d]$ given by sending $1 \mapsto \phi_0$ and $u \mapsto \phi_d$ is a bimodule isomorphism. \\

For finite groups, one may replace $C_*(G)$ with the reduced chains (setting degenerate chains to zero), say $K_*(G)$. Then $K_*(G) = R[G]$ concentrated in degree zero, and inversion gives a bimodule equivalence between $K_*(G)$ and $\text{Hom}(K_* G, R)$. 

For $G = SO(n)$, one may use a cellular model. Equip $G$ with the cell structure described in \cite[Chapter 3.D]{Hatcher}. This cell structure has the property that it is \emph{multiplicative}, meaning that $\mu: G \times G \to G$ is a cellular map (with product cell structure on the domain); this implies that $C_*^{CW}(G)$ is a dg-algebra. It comes with a collection of basic cells $e^0, e^1, \cdots, e^{n-1}$ which satisfy $$d_{CW}(e^i) = \begin{cases} 2e^{i-1} & i \text{even} \\ 0 & i \text{ odd}\end{cases},$$ and for which all other cells can be expressed as $e^I = e^{i_1} \cdots e^{i_m},$ where $I = (i_1, \cdots, i_m)$ is some decreasing sequence of integers $n > i_1 > \cdots > i_m > 0$. Finally, one may check by hand that $$e^i \cdot e^i = 0 \text{ if } i > 0, \quad e^i \cdot e^j = (-1)^{ij+1} e^j e^i \text{ if } 0 < i < j < n.$$ The Leibniz rule and associativity now determine the full dg-algebra structure of $C_*^{CW}(SO_n;R) = A$. It follows that if $1 = -1$ the differential is trivial and the product is strictly commutative. Henceforth, we work over rings $R$ with $1 = -1$. Write $d = n(n-1)/2 = \dim SO_n$. 

Now let us verify that we have a bimodule \emph{isomorphism} $\Phi: A \cong A^\vee[d]$. Fix the isomorphism $\text{ev}: A_d \cong R[d]$ to be the map sending the single $d$-cell $e^{n-1} \cdot e^{n-2} \cdots e^1$ to $R$. Because the differential of this cell is zero, $\text{ev}$ is a cocycle. It is then straightforward to see that the map $\Phi$ is given by sending $x \mapsto \Phi_x$ defined by $$\Phi_x(y) = \text{ev}(xy).$$ and further it's clear that $$\Phi_{axb}(y) = \Phi(axby) = \Phi(xbya) = (a \Phi_x b)(y),$$ so this indeed gives a bimodule map. Finally, $\Phi$ is an isomorphism: if $x = e^I = e^{i_1} \cdots e^{i_k}$, the only cell 

Lastly, let's verify that there is a zig-zag of quasi-isomorphisms between $C_*(SO_n; R)$ and $C_*^{CW}(SO_n; R)$. Say a simplex $\sigma: \Delta^k \to SO_n$ is \emph{cellular} if all $i$-dimensional faces map into the $i$-dimensional subcomplex of $SO_n$. Because this cell structure is multiplicative, the cellular simplices define a subalgebra $C_*^{c\Delta}(SO_n; R) \to C_*(SO_n; R)$, and a cellular approximation argument implies this map is a quasi-isomorphism. There is also a natural projection map $C_*^{c\Delta}(SO_n;R) \to C_*^{CW}(SO_n; R)$ given by sending a simplex $\sigma: \Delta^k \to SO_n^{(k)}$ to the degrees of $\sigma/\partial \sigma: S^k \to SO_n^{(k)}/SO_n^{(k-1)}$ above each $k$-cell. It is again straightforward to verify that this is also an algebra quasi-isomorphism, hence $SO_n$ is strongly dualizable.
\end{proof}

\begin{remark}
Interestingly, it does not appear true that we have bimodule isomorphisms $A \cong A^\vee[d]$ over the integers: working this out for $n = 3$ is enlightening. (The difficulty arises from the strange sign in the commutation relation $e^i \cdot e^j = (-1)^{ij+1} e^j \cdot e^i$.) This does not rule out hope that $SO_3$ is strongly dualizable over $\Bbb Z$, only that this cannot be seen using the rigid cellular model described above. 
\end{remark}

When $\frac 12 \in R$, there is a quasi-isomorphism $\Lambda_R(u) \to C_*(SO_3; R)$. Applying Proposition \ref{prop:strongduals} and Theorem \ref{PDiso}, we obtain the following corollary.

\begin{corollary}\label{cor:iso-for-SO3}
For $G = SO(3)$, if $2 \in R$ is either invertible or zero, then for $A = C_*(SO(3); R)$ and $M$ any right $A$-module, there is a natural isomorphism $H^{+,\text{tw}}_A(M) \cong H^+_A(M)[3]$ as graded $H^-_A$-modules.
\end{corollary}

We can generalize the previous calculations to that of the equivariant homologies of orbits. Let us begin by setting up notation for the relevant Poincar\'e duality lemma. 

Recall for the moment that our actions are \emph{right} actions, and hence the appropriate orbit space is the set of \emph{right} cosets; as we never work with spaces of left cosets, we prefer to write this quotient manifold as $G/H$. Our first ingredient is an $H$-equivariant Poincar\'e duality result. First, observe that $C_*(G/H)$ is a right $G$-module, via the action $G/H \times G \to G/H$ given by $Hg \cdot g' = Hgg'$. Also notice that $\mathfrak g/\mathfrak h = T_{He}(G/H)$ carries a right action of $H$, the derivative of the translation action. 

We write $\lambda_{\mathfrak g/\mathfrak h}$ for the $C_*(H)$-module which is a copy of $R$ concentrated in degree zero, where the action factors through the map $C_*(H) \to H_0(H) \cong R[\pi_0 H]$; the action of $\pi_0 H$ is via the determinant of the representation on $\mathfrak g/\mathfrak h$. (This is the local system corresponding to $w_1(G/H)$.) 

Finally, consider the complex $\text{Hom}_{C_* H}(C_* G, \lambda)$ of \emph{right} $C_*(H)$-equivariant maps. This carries a \emph{right} module structure over $C_*(G)$ via $(\alpha \cdot g)(g') = \alpha(gg')$.

\begin{lemma}\label{PD-for-H}
Suppose that $H \subset G$ is a closed subgroup, and that either $2 = 0$ in $R$ or that the isotropy representation of $H$ on $\mathfrak g/\mathfrak h$ is orientable. Then there is a zigzag of right module quasi-isomorphisms relating $C_*(G/H)$ and $\text{Hom}_{C_* H}(C_* G, \lambda_{\mathfrak g/\mathfrak h})[\dim G - \dim H]$. 
\end{lemma}
\begin{proof}
The argument of Proposition \ref{prop:G-is-weak-PD} follows with minimal change. 

The tangent space $\mathfrak g/\mathfrak h = T_{He}(G/H)$ carries a natural right action of $H$ (the stabilizer of $He \in G/H$) by taking the derivative of the right action of $H$ on $G/H$ at the identity coset. We write $S^{\mathfrak g/\mathfrak h}$ for the one-point compactification; there is an $H$-equivariant collapse map $\log_H: G/H \to S^{\mathfrak g/\mathfrak h}$. The desired map is given by $$x \mapsto \alpha_x, \quad \text{where}\quad \alpha_x(y) = \log_H(x \cdot y).$$ Each $\alpha_x$ is $C_* H$ equivariant as $$\alpha_x(yh) = \log_H(x \cdot yh) = \log_H(x \cdot y) \cdot h = \alpha_x(y) \cdot h,$$ and $\alpha$ is $C_*(G)$-equivariant as $$\alpha_{xg}(y) = \log_H(xg \cdot y) = \log_H(x \cdot gy) = \alpha_x(gy) = (\alpha_x \cdot g)(y).$$ Now work with the reduced version $\overline C_*(S^{\mathfrak g/\mathfrak h})$, killing all chains in degrees larger than $\dim G - \dim H$ and killing boundaries in that degree. Finally, there is an $H$-equivariant quasi-isomorphism $$\lambda_{\mathfrak g/\mathfrak h}[\dim G - \dim H] \to \overline C_*(S^{\mathfrak g/\mathfrak h})$$ given by sending $1$ to a chosen fundamental class; $\pi_0(H)$ acts on this fundamental class by the determinant of the isotropy representation, so this map is $H$-equivariant. This gives a zig-zag of $C_*(G)$-equivariant maps $$C_*(G/H) \xrightarrow{\alpha} \text{Hom}_{C_* H}(C_* G, \overline C_*(S^{\mathfrak g/\mathfrak h})) \leftarrow \text{Hom}_{C_* H}(C_* G, \lambda_{\mathfrak g/\mathfrak h}).$$

That $\alpha$ is a quasi-isomorphism once again rests on the fact that it provides a model for the Poincar\'e duality equivalence. In the case that $\lambda_{\mathfrak g/\mathfrak h}$ is non-trivial as an $H$-module, this uses non-orientable Poincar\'e duality, which gives an equivalence between $C_*(G/H)$ and $C^*(G/H; \lambda_{\mathfrak g/\mathfrak h})$, where the latter is given by the local system with monodromy given by the orientation character of $G/H$, which factors as $\pi_1(G/H) \to \pi_0(H) \xrightarrow{\det \mathfrak g/\mathfrak h} \{\pm 1\}$. 
\end{proof}

In particular, when $H$ is connected or at least has orientable isotropy representation, we have $C_*(G/H) \simeq \text{Hom}_{C_* H}(C_* G, R)$. 

This is enough for us to calculate the equivariant cohomology of orbits.

\begin{theorem}\label{orbitcalc}If $G$ is a compact Lie group and $H$ is a closed subgroup, let $G/H$ be the orbit of right cosets of $H$ (which is thus a right $G$-space). Then we have the following isomorphisms of $H_G^-$-modules: \begin{enumerate}
    \item $H_G^+(G/H) \cong H_H^+(R)\cong H_*(BH;R)$\\
    \item $H_G^-(G/H) \cong H_H^-(\lambda_{\mathfrak g/\mathfrak h})[\dim G - \dim H] \cong H^{\dim G - \dim H-*}(BH;\lambda_{\mathfrak g/\mathfrak h})$\\
    \item $H_G^\infty(G/H) \cong H_H^\infty(\lambda_{\mathfrak g/\mathfrak h})[\dim G - \dim H]$.
\end{enumerate}

Here, $H_G^-$ acts via the restriction map $H_G^- \to H_H^-$. 
\end{theorem}

\begin{proof}The map $C_*(G) \to C_*(G/H)$ is $C_*(H)$-invariant, inducing a quasi-isomorphism of right $G$-modules $B(R, H, G) \to C_*(G/H)$ (that this is a quasi-isomorphism follows from \cite[Theorem~3.9]{GM}). There is further a canonical isomorphism of chain complexes $$B(B(R,H,G),G,R) \cong B(R,H,B(G,G,R)).$$ Notice that all of these equivalences have been $C^-_G$-equivariant. Finally, the inclusion quasi-isomorphism $B(H,H,R) \to B(G,G,R)$ induces a $C^-_G$-equivariant quasi-isomorphism $$B(R,H,H) = B(R,H,B(H,H,R)) \simeq B(R, H, B(G, G, R)),$$ where the action of $C^-_G$ on $B(H,H,R)$ factors through the restriction $C^-_G \to C^-_H$.

Thus, (1) follows. The same argument in general identifies $H_G^+(G \otimes_H M)$ with $H_H^+(M)$ as $H^-_G$-modules. The isomorphism to group homology follows from Lemma \ref{groupcoho}.\\

\noindent (2) Applying Lemma \ref{PD-for-H} and the invariance result Theorem \ref{cb-inv}, we have the following sequence of equivalences (suppressing both $C_*$ and the dimension shift from notation)
\begin{align*}cB(R, G, G/H) &= \text{Map}_G\left(B(R,G,G), G/H\right) \simeq \text{Map}_G\left(B(R,G,G), \text{Map}_H(G, \lambda_{\mathfrak g/\mathfrak h})\right)\\
\cong \text{Map}_H&(B(R,G,G), \lambda_{\mathfrak g/\mathfrak h}) \simeq \text{Map}_H(B(H,H,R), \lambda_{\mathfrak g/\mathfrak h}) = cB(R,H, \lambda_{\mathfrak g/\mathfrak h}).\end{align*}

The inverse of the second equivalence is given via $\eta(b)(g) \mapsto \eta(bg)$. The last equivalence follows from \cite[Page~11]{GM} and verifying that, in that language, $B(G,G,R)$ is a proper split Kunneth resolution of $R$ as an $H$-module, which follows similar lines as verifying that $C_*(G)$ is split as an $H$-module, implicit in the proof of (1); it is induced by the map $B(H, H, R) \to B(G, G, R)$. Finally, all of these maps are transparently $C^-_G$-equivariant, in the final case acting on $B(H, H, R)$ by restriction. The isomorphism to group cohomology follows from Lemma \ref{groupcoho}, adapted to twisted coefficient systems.

We can combine the previous two parts and the long exact sequence of Theorem \ref{eq-package} (4) to verify the isomorphism on Tate homology. This follows from the following homotopy commutative diagram, where every vertical arrow is a quasi-isomorphism and horizontal arrow is an appropriate modification of the norm map.

\[\begin{tikzcd}
	{B(G/H, G, cB(R, G, G)) } & {cB(R, G, G/H)} \\
	{B(R, H, cB(R,G,G))} & {cB(R, G, G/H)} \\
	{B(R, H, cB(R, H, \text{Th}(\mathfrak g/\mathfrak h)))} & {cB(R, H, \tilde S^{\mathfrak g/\mathfrak h})} \\
	{B(R,H, cB(R,H,H\otimes \lambda_{\mathfrak g/\mathfrak h}))[\dim G - \dim H]} & {cB(R, H, \lambda_{\mathfrak g/\mathfrak h})[\dim G - \dim H]}
	\arrow[from=1-1, to=1-2]
	\arrow[from=1-1, to=2-1]
	\arrow[from=1-2, to=2-2]
	\arrow[from=2-1, to=2-2]
	\arrow[from=2-1, to=3-1]
	\arrow[from=4-1, to=3-1]
	\arrow[from=3-1, to=3-2]
	\arrow[from=2-2, to=3-2]
	\arrow[from=4-2, to=3-2]
	\arrow[from=4-1, to=4-2]
\end{tikzcd}\]

The second vertical maps are induced by the map of pairs (recall that $cB$ is contravariant in the algebra and covariant in the module) $$(G,G) \to \left(H, \text{Th}(\mathfrak g/\mathfrak h)\right),$$ the first map inclusion $H \hookrightarrow G$ and the second a collapse map to the Thom space of the normal bundle to $H$ inside $G$. That this defines a quasi-isomorphism on coBorel constructions follows by an argument using the Poincar\'e duality equivalences, together with the Thom isomorphism. The bottom-left vertical map is induced by a chain-level Thom isomorphism map, and the bottom-right map is given by picking a fundamental class. 

Finally, one may identify $$\text{Hom}_H(EH, H \otimes \lambda) \cong \text{Hom}_H(EH, H) \otimes \lambda \quad \text{and} \quad B(R, H, M \otimes \lambda) \cong B(\lambda, H, M)$$ for any character $\lambda$ by explicit isomorphisms:

If $\pi: C_*(H;R) \to H_0 H \to R$ is the composite of the projection and the map which sends $h \in \pi_0 H$ to $\rho(h) \in \{\pm 1\}$ (where $\rho$ is the representation defining the character $\lambda$), then the latter map is given by $$[a_1 \mid \cdots \mid a_k] m \mapsto \pi(a_1 \cdots a_k) [a_1 \mid \cdots \mid a_k] m,$$ and the former map is given by $$\eta \mapsto \eta', \quad \eta'([a_1 \mid \cdots \mid a_k]a) = \pi(a_1 \cdots a_k a)\eta([a_1 \mid \cdots \mid a_k] a).$$ Chasing the norm map through these equivalences, one finds that we have identified the mapping cone of $N^G_{G/H}$ with the mapping cone of $N^H_{\lambda_{\mathfrak g/\mathfrak h}}$, and that all maps are $C^-_G$-equivariant, giving the desired result. 
\end{proof}

We conclude this section with some related results which we use in the main text.

\begin{example}\label{SO3calc}We may apply this to calculate the three cases relevant to us in this text: $G = SO(3)$ and $H$ one of the three subgroups $\{e\}, SO(2),$ and $SO(3)$. In every case $H$ is connected, and so the representation $\lambda_{\mathfrak g/\mathfrak h}$ of $\pi_0 H$ is trivial and is precisely $R[3 - \dim H]$. Thus $H^\bullet_{SO(3)}(SO(3)/H; R) = H^\bullet_H(R)$, with a dimension shift if appropriate.

When $H = \{e\}$, the Tate homology is trivial, and $H^+_{\{e\}}(R) = H^-_{\{e\}}(R) = R$ concentrated in degree zero. (This is just the axiomatic property of Tate homology.) 

When $H = SO(2)$, we are left with $H^+_{SO(2)}(R) = H_*(\mathbb{CP}^\infty)$ and $H^-_{SO(2)}(R) = H^{-*}(\mathbb{CP}^\infty)$ by Lemma \ref{groupcoho}, and the Tate homology is a splicing of these. By Proposition \ref{TateLoc}, as modules over $H^-_{SO(2)}(R)$, we may write the Borel homology, coBorel homology, and Tate homology respectively as $$R[V], R\llbracket V^{-1}\rrbracket, R[V, V^{-1}\rrbracket,$$ where $|V| = 2$. This isomorphism is $H^-_{SO(3)}(R)$-equivariant for Borel and coBorel homology, and is $H^-$ equivariant for Tate homology if $2$ is either invertible or zero in $R$. 

Because $H^-_{SO(3)}(R) \to H^-_{SO(2)}(R)$ sends $$p_1 \mapsto V^2 \in H^{-*}(BSO(2);R),$$ we learn the module structure of these over $\langle p_1\rangle \in H^{-*}(BSO(3);R)$. When $\frac 12 \in R$, we have 
$$H^{-*}(BSO(3);R) \cong H^{-*}(BSU(2);R) \cong R\llbracket p_1\rrbracket,$$ 
so this determines the module structure over $H^-_{SO(3)}(R)$. 

On the other hand, we have that $$H^{-*}(BSO(3);\mathbb Z/2) = (\mathbb Z/2)\llbracket w_2, w_3\rrbracket$$ is a power series ring in two generators (where here one takes $H^{-*} = \prod H^{-k}$ as opposed to a direct sum; the latter would give a polynomial ring). This is standard, and proved in \cite[Proposition~3.12]{hatchervector}. A calculation with integral coefficients is given in \cite{brown1982cohomology}, and in particular shows that $p_1$ restricts to $e^2$ on oriented 2-plane bundles, where $e$ is the Euler class in $H^-_{SO(2)}$. 
\end{example}

Lastly, to apply these calculations to geometric chains, we should know that we can pass between the singular chain dga and the geometric chain dga.

\begin{proposition}\label{sm-gm-comparison} There is a chain of algebra quasi-isomorphisms $$C_*^{\textup{gm}}(SO(3);R) \simeq C_*(SO(3);R).$$
\end{proposition}
\begin{proof}
We neglect to mention the coefficients $R$ throughout, as they play no major role.

First we should replace singular chains with something more easily comparable with the degeneracy relations involved in $C_*^{\textup{gm}}$. We define the chain complex of smooth singular chains $C_*^\textup{sm}(M)$ to be the set of smooth maps from $\Delta^n \to M$; then it has a quotient $C_*^\textup{smd}(M;R)$ after we carry out the previous identifications under orientation-preserving diffeomorphisms and quotienting by the same degeneracies. There are, in this situation, two forgetful maps $C_*^\textup{smd}(M) \leftarrow C_*^\textup{sm}(M) \to C_*(M)$. The easiest way to see that all of these maps are quasi-isomorphisms is to prove that all of these theories are in fact homology theories on smooth manifolds $M$, and that the induced map on $M = \text{pt}$ is an isomorphism on homology. In fact the rightmost chain complexes are identical for $M = \text{pt}$. On the other hand, $C_*^\textup{smd}(\text{pt}) = R$, a copy of the ground ring concentrated in degree zero. All higher-dimensional chains are degenerate. Because $H_*(\text{pt}) = R_{(0)}$, and the given maps do the obvious things to points (which are cycles generating $H_*(\text{pt})$), they are all homology isomorphisms. A proof that these are homology theories follows similar lines as Theorem \ref{eilenberg}.

It is also possible to show these quasi-isomorphisms extremely explicitly, using e.g. smooth approximation to find cycles in $C_*^\textup{sm}$ homologous to any given cycle in $C_*$.

So we now need to relate the chain complex of smooth simplices modulo degeneracy, $C_*^\textup{smd}(M)$, to $C_*^{\textup{gm}}(M)$.

The theory interpolating between these is the chain complex of \emph{triangulated geometric chains} on a smooth manifold $M$. This chain complex $C_*^\textup{tgm}(M)$ is functorial under smooth maps, and its homology groups define a homology theory (as before).

Precisely, a \emph{triangulated basic chain} on $M$ is a compact smooth oriented manifold with corners $P$ equipped with a smooth triangulation (that is, a homeomorphism $f: |X| \to P$, where $|X|$ is the realization of a simplicial complex, so that $f$ is a diffeomorphism from each closed simplex onto its image), and a smooth map $\sigma: P \to M$. Two triangulated basic chains are isomorphic if there is an orientation-preserving diffeomorphism $\varphi: P \to P'$ and an isomorphism of simplicial complexes $\psi: X \to X'$ so that $\varphi' f = f' |\psi|$ and $\sigma' f = f' \sigma$. The triangulated geometric chain complex $C_*^\textup{tgm}(M)$ is defined following the same procedure as for $C_*^{\textup{gm}}(M)$: identify orientation-reversals with their negative, and quotient by the subcomplex of degenerate chains (triangulated basic chains for which the images of $\sigma|_{\Delta^k}$ and $\partial \sigma|_{\Delta^k}$ are both contained in the image of some smooth manifold of smaller dimension than $k$, resp $k-1$, for each component simplex $\Delta^k$ ).

Now observe that there are natural chain maps $$C_*^{\textup{gm}}(M) \leftarrow C_*^\textup{tgm}(M) \to C_*^\textup{smd}(M).$$ The left map is given by forgetting the triangulation and thinking of a smooth manifold with corners as a very special kind of $\delta$-chain, and the right map is given by sending a triangulated $n$-chain to the sum of its component simplices: $\sigma: P \to M$ to $$\sum_{\Delta^n \subset X} \sigma\big|_{\Delta^n}.$$ As before, because these are homology theories, to show that these are quasi-isomorphisms in general it suffices to check that these maps are isomorphisms on $H_*(\text{pt})$. Each of these chain complexes (because of the nondegeneracy requirements) are simply a copy of $R$ concentrated in degree zero, and the induced map between them is the identity.

Now, if $G$ is a Lie group, $C_*^\textup{smd}(G)$ is a dg-algebra using the same product operation as simplicial chains, the Eilenberg-Zilber product. One must check that this respects the added relations in $C_*^{\textup{smd}}$. This also defines a product on $C_*^\textup{tgm}(M)$: now if $\sigma: P \to M$ and $\eta: Q \to M$ are triangulated basic chains, their product is $\sigma \times \eta: P \times Q \to M$ equipped with the product triangulation (again, triangulate each component $\Delta^k \times \Delta^n$ in a standard way). With these algebra structures, the maps $$C_*^\textup{smd}(M) \leftarrow C_*^\textup{tgm}(M) \to C_*^{\textup{gm}}(M)$$ are in fact dg-algebra homomorphisms. It is even easier to see that the map $C_*^\textup{smd}(M) \leftarrow C_*^\textup{sm}(M) \to C_*(M)$ are dg-algebra homomorphisms, giving us the desired zig-zag of dg-algebra equivalences.
\end{proof}

\section{Periodic homological algebra}\label{PeriodicMachine}
In this section, we set up a version of the machinery before, built to work for complexes graded over $\mathbb Z/2N$ which are finite in a suitable sense. Here, a \textbf{$\Bbb Z/2N$-graded chain complex} is a collection of $R$-modules $C_i$, indexed by $i \in \Bbb Z/2N$, together with $R$-linear maps $d_i: C_i \to C_{i-1}$ which satisfy $d_{i-1} d_i = 0$. 

These objects are essentially equivalent to a more standard, $\Bbb Z$-graded object: dg-modules over $R[y, y^{-1}]$, where $|y| = 2N$, which we will call `$2N$-periodic complexes'.\footnote{This equivalence is one of many reasons to demand our complex is graded by $\Bbb Z/k$ for an even integer $k$: if $|y|$ is odd and $R$ has characteristic other than two, then $R[y, y^{-1}]$ is not a graded algebra, as $y \cdot y \ne - y \cdot y$! It is also difficult, for related reasons, to define a tensor product of $\Bbb Z/k$-graded complexes for $k$ odd.} Given a $\Bbb Z/2N$-graded complex $C$, we may pass to its \textit{unrolled complex} $\widetilde C_*$ of $C$, defined as $\widetilde C_k = C_{k \!\!\! \mod 2N}$, with the obvious differential. The identity map gives an isomorphism $\varphi: \widetilde C \to \widetilde C[2N]$, and this isomorphism gives $\widetilde C$ the structure of an $R[y, y^{-1}]$-module. Chain maps between $\Bbb Z/2N$-graded modules lift to $R[y, y^{-1}]$-linear chain maps $\widetilde C \to \widetilde C'$, and the lift is well-defined up to multiplication by some $y^k$. On the other hand, given an $R[y, y^{-1}]$-module $\widetilde C$,  we recover the $\mathbb Z/2N$-graded complex by picking an interval $[i,i+2N-1] \subset \mathbb Z$, and for $k \in [i, i+2N-1]$, define $$C_{k \!\!\! \mod 2N} = \widetilde C_k;$$ the differential is defined in the obvious way except on $C_{i\!\!\!\mod 2N}$, where it is defined as $dx = \varphi(\widetilde d x)$; this makes sense as $\varphi(\widetilde dx)$ is in degree $i-1+2N = i+2N-1$. The definitions immediately imply that the $\mathbb Z/2N$-graded complex $C$ did not depend on the choice of representative interval. An $R[y,y^{-1}]$-linear chain map then descends to a $\Bbb Z/2N$-graded chain map; the map associated to $f$ and $y^k \cdot f$ agree, for any integer $k$.\\

In the main text, we will frequently want to use a filtration resembling the index filtration on the Morse-Bott complex of a finite-dimensional compact smooth manifold, but because the index is only defined in $\mathbb Z/8$, this doesn't make sense! So if we want to use spectral sequences to check that maps are quasi-isomorphisms, we must do something at least slightly more subtle.

Because $H_k(C) = H_k(\widetilde C)$ for all $k$, if we're trying to show that a map $C \to C'$ of $\mathbb Z/2N$-graded complexes is a quasi-isomorphism, this is equivalent to showing that the same is true for the lift map $\tilde f: \widetilde C \to \widetilde C'$ of $2N\mathbb Z$-periodic complexes. We may more or less pass freely between these notions. To avoid notational irritation, we ignore the periodicity isomorphism $\varphi$ --- up to isomorphism of $2N\mathbb Z$-periodic complexes, we can take $\widetilde C_k = \widetilde C_{k+2N}$ on the nose, and $\varphi = \text{Id}[2N]$.

Instead of attempting to filter a $\mathbb Z/2N$-graded complex, we find that the appropriate notion seems to the following.

\begin{definition}\label{pfilt}Let $\widetilde C$ be a $2N\mathbb Z$-periodic complex. A \emph{periodic filtration} on $\widetilde C$ is a filtration $\cdots \subset F_s \widetilde C \subset F_{s+1} \widetilde C \subset \cdots$ with $F_{s+2N}\widetilde C = (F_s \widetilde C)[2N]$; that is, $$F_{s+2N} \widetilde C_{t+2N} = F_s \widetilde C_t.$$ 
\end{definition}

So while the filtration on each of the individual abelian groups $F_s C_k$ may stabilize (so $F_r C_k = C_k$ for large $r$) --- as indeed is frequently the case for us --- the filtration itself is infinite in both directions (i.e., neither is the complex $F_{-r}\widetilde C$ equal to zero for any $r \geq 0$, nor is $F_r \widetilde C = \widetilde C$ for any $r$).

We now pass to the spectral sequence $E^r_{s,t}$ of the filtered complex $\widetilde C$. The associated graded complex is $$E^0_{s,t}\widetilde C = F_s \widetilde C_{s+t}/F_{s-1} \widetilde C_{s+t},$$ and the fact that the filtration is periodic implies that $E^0_{s,t} = E^0_{s+2N,t}$. This identification is the map induced by the periodicity isomorphism $\widetilde C \to \widetilde C[2N]$, which is a filtered chain map; because it induces an isomorphism on the $E^0$ page, the same is true for all pages $E^r$ of the spectral sequence: there is a periodicity isomorphism $E^r_{s,t} \to E^r_{s+2N,t}$ preserving the differentials.

If so desired, we may thus view this filtration as inducing a $(\mathbb Z/2N,\mathbb Z)$-bigraded spectral sequence $E^r_{[s],t}$, with $[s] \in \mathbb Z/2N$. This may be pictured as a cylinder, where differentials wrap around.

If we have a map $f: \widetilde C \to \widetilde C'$ of unrolled complexes, compatible with a periodic filtration of each, it induces a map of the spectral sequences $E^r_{s,t}$. We would like to know when we can check that the map $f$ is a quasi-isomorphism from the corresponding fact about the $E^2$ page of this spectral sequence. Because the `unrolled' spectral sequence $E^r_{s,t}$ (considered with $(\mathbb Z, \mathbb Z)$-bigrading, $2N\mathbb Z$-periodic in the first grading) is a whole-plane spectral sequence, it is difficult to prove and often false that the $E^\infty$ page actually \emph{calculates} the homology groups. However, as long as the periodic filtration is \emph{complete, exhaustive, and Hausdorff}, we at least still know from Proposition \ref{ssComparison} that we may detect quasi-isomorphisms from isomorphisms on the $E^2$ page; further, if the spectral sequence degenerates on some finite page, we may indeed calculate the associated graded homology groups from the $E^\infty$ page.

Now let $A$ be a bounded, non-negatively graded dg-algebra satisfying Poincar\'e duality of degree $n$, as in Chapter \ref{borel}. Let $M$ be a right dg-$A$-module, graded over $\mathbb Z/2N$, equipped with a complete, exhaustive, and Hausdorff periodic filtration. 

We would like to apply the constructions of the previous sections to construct `equivariant homology' complexes $C^\bullet_A(M)$, with corresponding periodic filtrations. The clear thing to try is to pass to the unfolded complex $\tilde M$ and consider the corresponding complexes $C^\bullet_A(\tilde M)$. However, because $\tilde M$ is unbounded in both directions, these filtrations are no longer complete: $BA$ has elements in arbitrarily large degrees, which by pairing with elements of $\tilde M$ in arbitrarily low filtration may contribute a sequence of nonzero elements of $C^\bullet_A(\tilde M)$ of arbitrarily low filtration but fixed degree. The spectral sequence will not converge to its homology, but to its full completion's, as in Remark \ref{rmk:completion}.\\

As suggested by \cite[Definition 5.35]{GOH}, these spectral sequences are important enough that we simply pass to the full completion. We set $$\hat C^\bullet_A(\tilde M)_k = \text{lim}_q \colim_p C^\bullet_A(F_p \tilde M_k/F_q \tilde M_k),$$ equipped with the filtration $$F_p \hat C^\bullet_A(\tilde M) = \lim_{q < p} C^\bullet_A(F_p \tilde M/F_q \tilde M).$$ This satisfies the necessary periodicity, so forgetting gradings to $\Bbb Z/2N$ we obtain a $\Bbb Z/2N$-graded complex $\hat C^\bullet_A(M)$ equipped with a periodic filtration. Furthermore, defined in this way, the new filtration is tautologically complete, exhaustive, and Hausdorff. 

We may assemble this into the following homology theories for periodically graded $A$-modules.

\begin{theorem}\label{4flavors-periodic}Let $A$ be an $R$-free dg-algebra. If $M$ is a $\mathbb Z/2N$-graded right $R$-free $A$-module, equipped with a complete, exhaustive, Hausdorff periodic filtration, there are $\mathbb Z/2N$-graded complexes $\hat C^\bullet_A(M)$, for $\bullet \in \{+, \;(+, \textup{tw}), \; -,\; \infty\}$, satisfying the following properties.
\begin{enumerate}
\item Each $\hat C^\bullet_A(M)$ is equipped with a left action of $C^-_A(R) := C^-_A$.
\item The $\hat C^\bullet_A(M)$ are functorial as $C^-_A$-modules under $A$-module maps $f: M \to M'$ of $\mathbb Z/2N$ graded complexes, which are filtered in the sense that there is some $d$ for which $\tilde f(F_s \tilde M) \subset F_{s+d} \tilde M'$ for all $s$.
\item There is a $(\mathbb Z/2N,\mathbb Z)$-bigraded spectral sequence of $H(A)$-modules, $$H_{p+q}(\textup{gr}_p(M)) \rightarrow H_{p+q}(M).$$ If each $\textup{gr}_p M$ is bounded above, the spectral sequence computes the target.
\item There is a $(\mathbb Z/2N,\mathbb Z)$-bigraded spectral sequence of $H^-_A(R)$-modules $$H^\bullet_A(\textup{gr}(\tilde M)) \rightarrow H^\bullet_A(M),$$ which detects isomorphisms in the target. If each $\textup{gr}_p \tilde M$ is bounded above and $A$ is non-negatively graded, the $H^-$ spectral sequence computes the target.
\item A filtered $A$-module map $M \to M'$ which is a quasi-isomorphism on the associated graded complexes induces a quasi-isomorphism $C^\bullet_A(M) \to C^\bullet_A(M').$
\item If the associated graded is free over $A$, in the sense that $\textup{gr}_p \tilde M \cong X_p \otimes A$ for each $p$ and some finite-dimensional $R$-free complex $X_p$, then $H^\infty_A(M) = 0$.
\item There is a long exact sequence of $H^-_A(R)$-modules $$\cdots \to H^{+,\textup{tw}}_A(M) \to H_A^-(M) \to H_A^\infty(M) \xrightarrow{[-1]} H^{+,\textup{tw}}_A(M) \to \cdots$$ natural under filtered $A$-module maps.
\item If $A$ satisfies weak Poincar\'e duality of degree $n$, we have a natural $R$-module isomorphism $H^{+,\textup{tw}}_A(M) \cong H^+_A(M)[n]$. If $A$ satisfies strong Poincar\'e duality of degree $n$, this natural isomorphism is one of $H^-_A(R)$-modules.
\end{enumerate}
\end{theorem}

\begin{proof}
Each $C^\bullet_A(F_p \tilde M)$ is a $C^-_A$-module natural under $A$-module morphisms, so the defining limit and colimit $\lim_q \colim_p C^\bullet_A(F_p \tilde M/F_q \tilde M)$ are limits and colimits of $C^-_A$-modules; the resulting object is therefore also a $C^-_A$-module, natural in morphisms between the (periodically) filtered complexes. Thus, (1) and (2) both follow.

As for (3), the periodic filtration on $\tilde M$ is assumed to be a complete exhaustive Hausdorff filtration of $A$-modules, so as discussed above we have a spectral sequence of $H(A)$-modules to $H(M)$ which detects isomorphisms in the target, and similarly for the existence of the spectral sequences in (4) --- the complexes $\hat C^\bullet_A(M)$ are defined precisely so that these spectral sequences exist. What remains is to discuss when they compute the homology of the target. 

In general, suppose $\tilde M_k$ is a periodically filtered complex. By assumption, each $\text{gr}_p(\tilde M_k) \cong \text{gr}_{p+2N}(\tilde M_{k+2N})$ is bounded above in $k$, which has the same bound for $k \equiv k' \pmod 2N$. Because there are only finitely many equivalence classes modulo $2N$, it follows that there exists a uniform $C$ so that if $\text{gr}_p \tilde M_k \ne 0$, we have $k - p \le C$. It follows that $F_p \tilde M_k$ is only nonvanishing in degrees $k \le C + p$. In particular, $H_n(F_p \tilde M)$ vanishes for $p < n - C$, so the filtration is regular. The full statement of (3) now follows. 

As for (4), we want $F_p \hat C^\bullet = \lim_{q < p} C^\bullet(F_p \tilde M/F_q \tilde M)$ to be bounded above; for that it suffices by the above argument that each $C^\bullet(\text{gr}_p \tilde M)$ is bounded above. This follows for $C^-$ when $A$ is non-negatively graded (as $\text{Hom}_A(B(R,A,A), \text{gr}_p \tilde M)$ is supported in degrees no larger than $\tilde M$ is).

Now (5) follows thanks to the existence of these spectral sequences and the invariance results Theorem \ref{bar-inv} and Theorem \ref{eq-package}(1), which imply that an $A$-equivariant quasi-isomorphism on $\text{gr} \tilde M$ induces a quasi-isomorphism on $H^\bullet_A(\text{gr} \tilde M)$. 

Next, the vanishing theorem \ref{eq-package} (4) implies that if $\text{gr} \tilde M \cong X \otimes A$, we have $H^\infty_A(\text{gr} \tilde M) = 0,$ and so the spectral sequence $H^\infty_A(\text{gr} \tilde M) \to H^\infty_A(M)$ collapses at the $E^1$ page with $E^1 = 0$. Because this spectral sequence detects isomorphisms, the map $0 \to C^\infty_A(\text{gr} \tilde M)$ is a quasi-isomorphism, and so under the assumptions of (6) we indeed have $H^\infty_A(M) = 0$. 

Item (7), the existence of this long exact sequence, follows from the fact that there is an isomorphism $\text{Cone}(\widehat N_M) \cong \widehat C^\infty_A(M)$. This is explained in somewhat more detail in \cite[Lemma 5.38]{GOH}. 

Finally, item (8) on the Poincar\'e duality isomorphisms follows by applying the arguments of Chapter \label{borel}: in the case of weak Poincar\'e duality, there is a natural map $B(N, A, R)[n] \to B(N, A, D_A)$ which induces a quasi-isomorphism of $R$-modules. Applying this to $N = F_p \tilde M$, we obtain a natural transformation $\hat C^+_A(M)[n] \to \hat C^{+,\text{tw}}(M)$ which induces an isomorphism on the $E^1$ page of the associated spectral sequence, and therefore a quasi-isomorphism on the target. Similarly in the case of strong Poincar\'e duality, only now the map is $B(N, A, D_A) \to B(N, A, R)[n]$ and it is $C^-_A$-equivariant. 
\end{proof}

\backmatter
\bibliographystyle{amsalpha}
\bibliography{biblio}
\printindex

\end{document}